\documentclass[a4]{article}

\font\es=eufm10

\def\ga{\mbox{\es {a}}} 
\def\gb{\mbox{\es {b}}} 
 
\def\gd{\mbox{\es {d}}} 
\def\gge{\mbox{\es {e}}} 
\def\gf{\mbox{\es {f}}} 
\def\gg{\mbox{\es {g}}} 
\def\gh{\mbox{\es {h}}}

\def\gq{\mbox{\es {q}}}

\def\gu{\mbox{\es {u}}}

\def\gA{\mbox{\es {A}}} 
 
\def\gC{\mbox{\es {C}}} 
\def\gD{\mbox{\es {D}}} 
\def\gE{\mbox{\es {E}}}

\def\gJ{\mbox{\es {J}}} 
\def\gK{\mbox{\es {K}}} 
 
\def\gM{\mbox{\es {M}}} 
\def\gN{\mbox{\es {N}}} 
 
\def\gP{\mbox{\es {P}}}

\def\gS{\mbox{\es {S}}}

\def\gW{\mbox{\es {W}}} 
\def\gX{\mbox{\es {X}}}

\def\sl{\mbox{\es {sl}}}
\def\so{\mbox{\es {so}}}
\def\sp{\mbox{\es {sp}}}
\def\su{\mbox{\es {su}}}
\def\tr{\mbox{\rm {tr}}}
\def\ad{\mbox{\rm {ad}}}
\def\Ad{\mbox{\rm {Ad}}}
\def\det{\mbox{\rm {det}}}
\def\diag{\mbox{\rm {diag}}}
\def\Iso{\mbox{\rm {Iso}}}
\def\Hom{\mbox{\rm {Hom}}}
\def\Ker{\mbox{\rm {Ker}}}
\def\ov{\overline}
\def\wti{\widetilde}
\def\dsum{\displaystyle \sum}

\def\sR{\mbox{\boldmath $\scriptstyle{R}$}} 
\def\sC{\mbox{\boldmath $\scriptstyle{C}$}} 
\def\sH{\mbox{\boldmath $\scriptstyle{H}$}} 
\def\dfrac#1#2{\displaystyle \frac{#1}{#2}}

\def\C{\mbox{\boldmath $C$}}

\def\H{\mbox{\boldmath $H$}}

\def\R{\mbox{\boldmath $R$}}

\def\X{\mbox{\boldmath $X$}}                     
\def\Y{\mbox{\boldmath $Y$}}      
\def\Z{\mbox{\boldmath $Z$}} 
\def\sR{\mbox{\boldmath $\scriptstyle{R}$}} 
 
\def\0{\mbox{\boldmath {0}}}    
\def\1{\mbox{\boldmath {1}}}      
\def\2{\mbox{\boldmath {2}}}      
\def\3{\mbox{\boldmath {3}}}      
\def\4{\mbox{\boldmath {4}}}      
\def\5{\mbox{\boldmath {5}}}      
\def\6{\mbox{\boldmath {6}}}      
\def\7{\mbox{\boldmath {7}}}      
\def\8{\mbox{\boldmath {8}}}      
\def\9{\mbox{\boldmath {9}}}      
\def\a{\mbox{\boldmath $a$}}
\def\b{\mbox{\boldmath $b$}}
\def\c{\mbox{\boldmath $c$}}
\def\dd{\mbox{\boldmath $d$}}
\def\e{\mbox{\boldmath $e$}}

\def\m{\mbox{\boldmath $m$}}
\def\n{\mbox{\boldmath $n$}}

\def\u{\mbox{\boldmath $u$}}
\def\v{\mbox{\boldmath $v$}}
\def\w{\mbox{\boldmath $w$}}
\def\x{\mbox{\boldmath $x$}}
\def\y{\mbox{\boldmath $y$}}
\def\z{\mbox{\boldmath $z$}}

\usepackage{epic,eepic}
\begin{document}
\baselineskip=14pt    

\vspace{20cm}
\begin{center}
{\Large{Exceptional Lie groups}}
\end{center}
\vspace{5mm}
\begin{center}
Ichiro Yokota
\end{center}

\vspace{5mm}

\begin{center}
{\Large{ Preface}}
\end{center}

In the end of 19 century, W. Killing and E. Cartan classified the complex simple Lie algebras, called $A_n, B_n, C_n, D_n$ (classical type) and $G_2, F_4, E_6, E_7, E_8$ (exceptional type). These simple Lie algebras and the corresponding compact simple Lie groups have offered many subjects in mathematicians. Especially, exceptional Lie groups are very wonderful and interesting miracle in Lie group theory.  
\vspace{1mm}

Now, in the present book, we describe simply connected compact exceptional simple Lie groups $G_2, F_4, E_6, E_7, E_8$, in very elementary way. The contents are given as follows. We 
first construct all simply connected compact exceptional Lie groups $G$ 
concretely. Next, we find all involutive automorphisms $\sigma$ of $G$, and 
determine the group structures of the fixed points subgroup $G^\sigma$ by $\sigma$. Note that they correspond to classification of all irreducible compact symmetric spaces $G/G^\sigma$ of exceptional type, and that they also correspond to classification of all non-compact exceptional simple Lie groups. Finally, we determined the group structures of the maximal subgroups of maximal rank. At any rate, we would like this book to be used in mathematics and physics.
\vspace{1mm}

The author thanks K. Abe, K. Mituishi, T. Miyasaka, T. Miyashita, T. Sato, O. Shukuzawa, K. Takeuchi and O. Yasukura for their advices and encouragements.

\begin{flushright}
\begin{tabular}{l}
Ichiro Yokota 
\end{tabular}
\end{flushright} 

\newpage

\vspace{3mm}

\begin{center}
{\Large{\bf Contents}}
\end{center}
\vspace{3mm}

\begin{center}
{\large{I. \quad Exceptional Lie group $G_2$}}
\end{center}

\begin{center}
\begin{tabular}{lr}
1.1. Cayley algebra $\gC$ \qquad \qquad \qquad & $\cdots \;\; 1$
\vspace{1mm}\\

1.2. Compact exceptional Lie group $G_2$ \qquad & $\cdots \;\; 2$
\vspace{1mm}\\

1.3. Outer automorphism of the Lie algebra $\gd_4$ \qquad & $\cdots \;\; 3$
\vspace{1mm}\\

1.4. Lie algebra $\gg_2$ of $G_2$ \qquad & $\cdots \;\; 8$
\vspace{1mm}\\

1.5. Lie subalgebra $\su(3)$ of $\gg_2$ \qquad & $\cdots \;\; 9$
\vspace{1mm}\\

1.6. Simplicity of ${\gg_2}^C$ \qquad & $\cdots \;\; 11$
\vspace{1mm}\\

1.7. Killing form of ${\gg_2}^C$ \qquad & $\cdots \;\; 13$
\vspace{1mm}\\

1.8. Roots of ${\gg_2}^C$ \qquad & $\cdots \;\; 14$
\vspace{1mm}\\

1.9. Automorphism $w$ of order $3$ and subgroup $SU(3)$ of $G_2$ & \qquad $\cdots \;\; 16$
\vspace{1mm}\\

1.10. Involution $\gamma$ and subgroup $(Sp(1) \times Sp(1))/\Z_2$ of $G_2$ & $\cdots \;\; 19$
\vspace{1mm}\\

1.11. Center $z(G_2)$ of $G_2$ \qquad & $\cdots \;\; 21$
\vspace{1mm}\\

1.12. Complex exceptional Lie group ${G_2}^C$ \qquad & $\cdots \;\; 21$
\vspace{1mm}\\

1.13. Non-compact exceptional Lie group $G_{2(2)}$ of type $G_2$ \qquad \qquad & $\cdots \;\; 22$
\vspace{1mm}\\

1.14. Principle of triality in $SO(8)$ \qquad & $\cdots \;\; 23$
\vspace{1mm}\\

1.15. Spinor group $Spin(7)$ \qquad & $\cdots \;\; 26$
\vspace{1mm}\\

1.16. Spinor group $Spin(8)$ \qquad & $\cdots \;\; 28$
\vspace{3mm}
\end{tabular}
\end{center}

\begin{center}
{\large{II. \quad Exceptional Lie group $F_4$}}
\end{center}

\begin{center}
\begin{tabular}{lr}
2.1. Exceptional Jordan algebra $\gJ$ \qquad & $\cdots \;\; 31$
\vspace{1mm}\\

2.2. Compact exceptional Lie group $F_4$ \qquad & $\cdots \;\; 33$
\vspace{1mm}\\

2.3. Lie algebra $\gf_4$ of $F_4$ \qquad & $\cdots \;\; 36$
\vspace{1mm}\\

2.4. Simplicity of ${\gf_4}^C$ \qquad & $\cdots \;\; 41$
\vspace{1mm}\\

2.5. Killing form of ${\gf_4}^C$ \qquad & $\cdots \;\; 45$
\vspace{1mm}\\

2.6. Roots of ${\gf_4}^C$ \qquad & $\cdots \;\; 47$
\vspace{1mm}\\

2.7. Subgroup $Spin(9)$ of $F_4$ \qquad & $\cdots \;\; 51$
\vspace{1mm}\\

2.8. Connectedness of $F_4$ \qquad & $\cdots \;\; 55$
\vspace{1mm}\\

2.9. Involution $\sigma$ and subgroup $Spin(9)$ of $F_4$ \qquad & $\cdots \;\; 57$
\vspace{1mm}\\

2.10. Center $z(F_4)$ of $F_4$ \qquad & $\cdots \;\; 58$
\vspace{1mm}\\

2.11. Involution $\gamma$ and subgroup $(Sp(1) \times Sp(3))/\Z_2$ of $F_4$ & $\cdots \;\; 59$
\vspace{1mm} \\

2.12. Automorphism $w$ of order $3$ and subgroup $(SU(3) \times SU(3))/\Z_3$\\
\quad of $F_4 $ & $\cdots 63$
\vspace{1mm}\\

2.13. Complex exceptional Lie group ${F_4}^C$ \qquad & $\cdots \;\; 65$
\vspace{1mm}\\

2.14. Non-compact exceptional Lie groups $F_{4(4)}$ and $F_{4(-20)}$ of type $F_4$ & $\cdots \;\; 66\;$
\end{tabular}
\end{center}
\vspace{5mm}

\begin{center}
{\large{III. \quad Exceptional Lie group $E_6$}}
\end{center}

\begin{center}
\begin{tabular}{lr}
3.1. Compact exceptional Lie group $E_6$ \qquad & $\cdots \;\; 68$
\vspace{1mm}\\

3.2. Lie algebra $\gge_6$ of $E_6$ \qquad & $\cdots \;\; 68$
\vspace{1mm}\\

3.3. Simplicity of ${\gge_6}^C$ \qquad & $\cdots \;\; 70$
\vspace{1mm}\\

3.4. Element $A \vee B$ of ${\gge_6}^C$ \qquad & $\cdots \;\; 71$
\vspace{1mm}\\

3.5. Killing form of ${\gge_6}^C$ \qquad & $\cdots \;\; 73$
\vspace{1mm}\\

3.6. Roots of ${\gge_6}^C$ \qquad & $\cdots \;\; 75$
\vspace{1mm}\\

3.7. Involution $\tau$ and subgroup $F_4$ of $E_6$ \qquad & $\cdots \;\; 82$
\vspace{1mm}\\

3.8. Connectedness of $E_6$ \qquad & $\cdots \;\; 82$
\vspace{1mm}\\

3.9. Center $z(E_6)$ of $E_6$ \qquad & $\cdots \;\; 85$
\vspace{1mm}\\

3.10. Involution $\sigma$ and subgroup $(U(1) \times Spin(10))/\Z_4$ of $E_6$ & $\cdots \;\; 86$
\vspace{1mm}\\

3.11. Involution $\gamma$ and subgroup $(Sp(1) \times SU(6))/\Z_2$ of $E_6$  & $\cdots \;\; 90$
\vspace{1mm}\\

3.12. Involution $\tau\gamma$ and subgroup $Sp(4)/\Z_2$ of $E_6$ \qquad & $\cdots \;\; 95$
\vspace{1mm}\\

3.13. Automorphism $w$ of order $3$ and subgroup $(SU(3) \times SU(3) \times$ \\\quad $SU(3))/\Z_3$ of $E_6$ \qquad $\;\;$& $\cdots \;\; 100$
\vspace{1mm}\\

3.14. Complex exceptional Lie group ${E_6}^C$ \qquad & $\cdots \;\; 105$
\vspace{1mm}\\

3.15. Non-compact exceptional Lie groups $E_{6(6)}, E_{6(2)}, E_{6(-14)}$ and \\\quad $E_{6(-26)}$ of type $E_6$ \quad & $\cdots \;\; 105 \;$ 
\end{tabular}
\end{center}
\vspace{3mm}

\begin{center}
{\large{IV. \quad Exceptional Lie group $E_7$}}
\end{center}

\begin{center}
\begin{tabular}{lr}
4.1. Freudenthal vector space $\gP^C$ \qquad & $\cdots \;\; 106$
\vspace{1mm}\\

4.2. Compact exceptional Lie group $E_7$ \qquad & $\cdots \;\; 108$
\vspace{1mm}\\

4.3. Lie algebra $\gge_7$ of $E_7$ \qquad & $\cdots \;\; 108$
\vspace{1mm}\\

4.4. Simplicity of ${\gge_7}^C$ \qquad & $\cdots \;\; 114$
\vspace{1mm}\\

4.5. Killing form of ${\gge_7}^C$  \qquad & $\cdots \;\; 116$
\vspace{1mm}\\

4.6. Roots of ${\gge_7}^C$ \qquad & $\cdots \;\; 118$
\vspace{1mm}\\

4.7. Subgroups $E_6$ and $U(1)$ of $E_7$  \qquad & $\cdots \;\;124$
\vspace{1mm}\\

4.8. Connectedness of $E_7$  \qquad & $\cdots \;\; 126$
\vspace{1mm}\\

4.9. Center $z(E_7)$ of $E_7$  \qquad & $\cdots \;\;129$
\vspace{1mm}\\

4.10. Involution $\iota$ and subgroup $(U(1) \times E_6)/\Z_3$ of $E_7$ \qquad & $\cdots \;\;130$
\vspace{1mm}\\

4.11. Involution $\sigma$ and subgroup $(SU(2) \times Spin(12))/\Z_2$ of $E_7$ \qquad & $\cdots \;\;133$
\vspace{1mm}\\

4.12. Involution $\tau\gamma$ and subgroup $SU(8)/\Z_2$ of $E_7$  & $\cdots \;\;142$
\vspace{1mm}\\

4.13. Automorphism $w$ of order $3$ and subgroup $(SU(3) \times SU(6))/\Z_3$ \\
\quad of $E_7$ \quad & $\cdots \;\;146$
\vspace{1mm}\\

4.14. Complex exceptional Lie group ${E_7}^C$ \qquad & $\cdots \;\;155$
\vspace{1mm}\\

4.15. Non-compact exceptional Lie groups $E_{7(7)}, E_{7(-5)}$ and $E_{7(-25)}$ \\ \quad of type $E_7$ \qquad & $\cdots \;\;155 \;$
\end{tabular}
\end{center}
\vspace{3mm}

\begin{center}
{\large{V. \quad Exceptional Lie group $E_8$}}
\end{center}
\begin{center}
\begin{tabular}{lr}

5.1. Lie algebra ${\gge_8}^C$  \qquad & $\cdots \;\; 157$
\vspace{1mm}\\

5.2. Simplicity of ${\gge_8}^C$ \qquad & $\cdots \;\; 157$
\vspace{1mm}\\

5.3. Killing form of ${\gge_8}^C$ \qquad & $\cdots \;\; 159$
\vspace{1mm}\\

5.4. Complex exceptional Lie group ${E_8}^C$ \qquad & $\cdots \;\;160$
\vspace{1mm}\\

5.5. Compact exceptional Lie group $E_8$ \qquad & $\cdots \;\; 161$
\vspace{1mm}\\

5.6. Roots of ${\gge_8}^C$  \qquad & $\cdots \;\;163$
\vspace{1mm}\\

5.7. Involution $\upsilon$ and subgroup $(SU(2) \times E_7)/\Z_2$ of $E_8$ \qquad & $\cdots \;\; 173$
\vspace{1mm}\\

5.8. Involution $\wti{\lambda}\gamma$ and subgroup $Ss(16)$ of $E_8$ \qquad & $\cdots \;\; 176$
\vspace{1mm}\\

5.9. Center $z(E_8)$ of $E_8$ \qquad & $\cdots \;\; 184$
\vspace{1mm}\\

5.10. Automorphism $w$ of order $3$ and subgroup$(SU(3) \times E_6)/\Z_3$ of \\
\quad $E_8$  \qquad & $\cdots \;\;184$
\vspace{1mm}\\

5.11. Automorphism $w_3$ of order $3$ and subgroup $SU(9)/\Z_3$ of $E_8$  \qquad & $\cdots \;\;186$
\vspace{1mm}\\

5.12. Automorphism $z_5$ of order $5$ and subgroup$(SU(5) \times SU(5))/\Z_5$\\ \quad of $E_8$  \qquad & $\cdots \;\;192$
\vspace{1mm}\\

5.13. Non-compact exceptional Lie groups $E_{8(8)}$ and $E_{8(-24)}$ of $E_8$ \qquad & $\cdots \;\; 199 \;$ 
\end{tabular}
\end{center}

\newpage
\vspace{5mm}

\begin{center}
{\Large{\bf Notation}}
\end{center}
\vspace{4mm}

$\R$, $\C = \R \oplus \R e_1$, $\H = \R \oplus \R e_1 \oplus \R e_2 \oplus \R e_3$ denote the fields of real, complex, quaternion numbers, respectively. 
\vspace{2mm}

For $\R$-vector space $V$, its complexification $\{u + iv \, | \, u, v \in V\}$ is denoted by $V^C$. The complex conjugation in $V^C$ is denoted by $\tau$:
$$
      \tau(u + iv) = u - iv. $$
$\R^C$ is briefly denoted by $C$.
\vspace{2mm}

For $K$-vector space $V$ ($K = \R, \C, C$), $\Iso_K(V)$ denotes all $K$-linear isomorphisms of $V$. For a $K$-linear mapping $f$ of $V$, $V_f$ denotes $\{v \in V \,| \, f(v) = v \}$. 
\vspace{2mm}

For $K$-vector spaces $V, W$ ($K = \R, C$), $\Hom_K(V, W)$ denotes all $K$-homomorph-\\ism $f : V \to W$. $\Hom_K(V, V)$ is briefly denoted by $\Hom_K(V)$.
\vspace{2mm}

Let $G$ be a group and $\sigma$ an automorphism of $G$. Then $G^\sigma$ denotes $\{g \in G \, | \, \sigma(g) = g\}$. For $s \in G$, $G^s$ denotes $\{g \in G \, | \, sgs^{-1} = g\}$. 
\vspace{2mm}

For topological spaces $X, Y$, $X \simeq Y$ denotes that $X$ and $Y$ are homeomorphic.
\vspace{2mm}

For groups $G, G'$, $G \cong G'$ denotes that $G$ and $G'$ are isomorphic as groups. Isomorphic two groups $G, G'$ are often identified: $G = G'$.
\vspace{2mm}

$M(n, K)$ denotes all $n \times n$ matrices with entries in $K$.
\vspace{2mm}

$E = \diag(1, \cdots, 1) \in M(n, K)$ is the unit matrix.
\vspace{2mm}

For $A \in M(n, K)$, ${}^tA $ denotes the transposed matrix of $A$ and $A^*$ denotes the conjugate transposed matrix of $A$: $A^* = {}^t\ov{A}$. 
\vspace{2mm}

$O(n) = \{A \in M(n, \R) \, | \, {}^tAA = E\}$ (orthogonal group),
\vspace{1mm}

$SO(n) = \{A \in O(n) \, | \, \det A = 1 \}$ (spcecial orthogonal group),
\vspace{1mm}

$U(n) = \{A \in M(n, \C) \, | \, A^*A = E\}$ or $\{A \in M(n, C) \, | \, \tau({}^tA)A = E\}$ (unitary group),
\vspace{1mm}

$SU(n) = \{A \in U(n) \, | \, \det A = 1 \}$ (speccial unitary group),
\vspace{1mm}

$Sp(n) = \{A \in M(n, \H) \, | \, A^*A = E\}$ (symplectic group).
\vspace{2mm}

For a Lie group $G$, its Lie algebra is denoted by the correspoding small Germann letter $\gg$. For example, $\su(n)$ is the Lie algebra of the special unitary group $SU(n)$.
\vspace{4mm}

{\large{Special notations}}
$$
   \gamma, \quad \sigma, \quad \iota, \quad \kappa, \quad \lambda, \quad \mu, \quad \tau, \quad \upsilon, \quad \chi, \quad g, \quad w. $$

\newpage

\setcounter{page}{1}

\begin{center}
\large{\bf Exceptional Lie group $G_2$}
\end{center}
\vspace{4mm}

{\bf 1.1. Cayley algebra $\gC$}
\vspace{3mm}

We denote the division Cayley algebra by $\gC$. We now explain this algebra. 
Consider an 8 dimensional $\R$-vector space with basis $\{e_0 = 1,e_1,e_2,e_3,
e_4,e_5,e_6,e_7 \}$ and define a multiplication between them as follows. In the
 figure below, the multiplication between $e_1,e_2,e_3$ is defined as
$$
            e_1e_2 = e_3, \quad e_2e_3 = e_1, \quad e_3e_1 = e_2, $$
and defined similarly on the other lines. For example, $e_1e_6 = e_7$, $e_4e_7 
= e_3$ etc. We regard that $e_2, e_5, e_7$ are also collinear, for example, 
$e_5e_7 = e_2$. $e_0 = 1$ is the unit of the multiplication and assume
$$
      {e_i}^2 = -1 ,\;\; i \neq 0, \quad 
      e_ie_j = -e_je_i ,\;\; i \neq j , i \neq 0 , j \neq 0,$$
and the distributive law. Thus $\gC$ has a multiplication.  $x1, x \in \R$ is 
briefly denoted by $x$. In $\gC$, the conjugate $\ov{x}$, an inner product $(x,
 y)$, the length $|x|$ and the real part $R(x)$ are defined respectively by

\vspace{-22mm}                   
\begin{picture}(400,200)
 \path(124,1.4)(130,-0.2)(124,-1.8)
 \path(114,10.0)(120,11.55)(115,7.0)
 \path(226,10.0)(220,11.55)(225,7.0)
 \path(226.7,20.0)(225,26.0)(229.4,22.0)
 \path(207.0,28.0)(210,34.0)(210.0,27.0)
 \path(168.5,107.0)(170,101.0)(171.5,107.0)
 \path(154.6,97.5)(153,91.8)(157.5,96.3)
 \put(166,124){$e_1$}
 \put(125,62){$e_2$}
 \put(90,-5){$e_3$}
 \put(171,30){$e_4$}
 \put(166,-6){$e_5$}
 \put(242,-5){$e_6$}
 \put(206,62){$e_7$}
 \thicklines
 \path(100,0)(240,0)(170,121.2)(100,0)
 \put(170,40.4){\circle{80.8}}
 \path(100,0)(205,60.6)
 \path(170,0)(170,121.2)
 \path(135,60.6)(240,0)
\end{picture}
\vspace{5mm}

$$
\begin{array}{cc}
    \ov{x_0 + \dsum_{i=1}^7x_ie_i} = x_0 - \dsum_{i=1}^7x_ie_i, & 
    \Big(\dsum_{i=0}^7x_ie_i, \dsum_{i=0}^7y_ie_i \Big) = \dsum_{i=0}^7x_iy_i,
\vspace{1mm}\\
   |x| = \sqrt{(x, x)}, & R\Big(x_0 + \dsum_{i=1}^7x_ie_i \Big) = x_0.
\end{array}$$
For $x \in \gC$, $x \neq 0$, we denote $\dfrac{\ov{x}}{|x|^2}$ by $x^{-1}$, 
then we have $xx^{-1} = x^{-1}x = 1$, and $\gC$ satisfies all axioms of a field
 except the associative law $x(yz) = (xy)z$. Of course the commutative law $xy =
 yx$ does not hold. Since the associative law does not hold, calculations in 
$\gC$ are complicated, however the following relations hold. (See Freudenthal [6] or Yokota [58]).
 For $a, b, x, y \in \gC$, we have
\vspace{2mm}

$\begin{array}{lll}
1  \quad (xy, xy) = (x, x)(y, y),\quad |xy| = |x||y|.
\vspace{1mm}\\

2  \quad (ax, ay) = (a, a)(x, y) = (xa, ya).
\vspace{1mm}\\

3  \quad (ax, by) + (bx, ay) = 2(a, b)(x, y),
\vspace{1mm}\\

4  \quad (ax, y) = (x, \ov{a}y), \quad (xa, y) = (x, y\ov{a}).
\vspace{1mm}\\

5  \quad \ov{\ov{x}} = x, \quad \ov{x + y} = \ov{x} + \ov{y}, \quad \ov{xy} = 
\ov{y} \, \ov{x}.
\vspace{1mm}\\

6  \quad (x, y) = (y, x) = \dfrac{1}{2}(\ov{x}y + \ov{y}x) = 
	\dfrac{1}{2}(x\ov{y} + y\ov{x}), \quad \ov{x}x = x\ov{x} = |x|^2. 
\vspace{1mm}\\

7  \quad a(\ov{a}x) = (a\ov{a})x, \quad a(x\ov{a}) = (ax)\ov{a}, \quad 
x(a\ov{a}) = (xa)\ov{a}.
\vspace{1mm}\\

  \,\, \quad a(ax) = (aa)x, \quad a(xa) = (ax)a, \quad x(aa) = (xa)a,
\vspace{1mm}\\

8  \quad \ov{b}(ax) + \ov{a}(bx) = 2(a,b)x = (xa)\ov{b} + (xb)\ov{a}.
\vspace{1mm}\\

9  \quad \mbox{We use a notation}\;\; \{x, y, z \} = (xy)z - x(yz), \;\;
\mbox{called the associator of}\;\; x, y, z.
\end{array}$

\noindent Then, we have
\vspace{1mm}\\
$\begin{array}{lll}

   \qquad  \{x, y, a \} = \{y, a, x \} = \{a, x, y \} = -\{y, x, a \} = -\{x, 
a, y \} = -\{a, y, x \}.
\end{array}$
\vspace{1mm}

\noindent For example, we have
$$
\begin{array}{l}
          (ax)y + x(ya) = a(xy) + (xy)a,
\vspace{1mm}\\
          (xa)y + (xy)a = x(ay) + x(ya),
\vspace{1mm}\\
          (ax)y + (xa)y = a(xy) + x(ay).
\end{array}$$

$\begin{array}{ll}
10 \quad (ax)(ya) = a(xy)a \quad \mbox{(Moufang's formula)}.
\vspace{1mm}\\
11 \quad R(xy) = R(yx), \quad R(x(yz)) = R(y(zx)) = R(z(xy)) \; (= R(xyz)).
\end{array}$
\vspace{1mm}

\noindent For an orthonormal basis $\{ 1, a_1, a_2, \cdots, a_7 \}$ of $\gC$ 
with respect to the inner product $(x, y)$, the following 12.1 $\sim$ 12.3 hold.\vspace{1mm}

$\begin{array}{lll}
12.1  \quad a_i(a_jx) = -a_j(a_ix), \quad \mbox{in particular,} \quad 
	a_ia_j = -a_ja_i, \quad i \neq j.
\vspace{1mm}\\
12.2  \quad a_i(a_ix) = -x, \quad \mbox{in particular,} \quad {a_i}^2 = -1.
\vspace{1mm}\\
12.3  \quad a_i(a_ja_k) = a_j(a_ka_i) = a_k(a_ia_j), \quad i, j, k \;\;
\mbox{are distinct.}
\end{array}$
\vspace{4mm}

{\bf 1.2. Compact exceptional Lie group $G_2$}
\vspace{3mm}

{\bf Definition.} The group $G_2$ is defined to be the automorphism group of 
the Cayley algebra $\gC$:
$$
        G_2 = \{ \alpha \in \Iso_{\sR}(\gC) \, | \,  
           \alpha(xy) = (\alpha x)(\alpha y) \}.$$

{\bf Lemma 1.2.1.} {\it For $\alpha \in G_2$, we have}
$$
       (\alpha x, \alpha y) = (x,y), \quad x,y \in \gC.$$

{\bf Proof.} For $\alpha \in G_2$ and $x \in \gC$, we have
$$
          \ov{\alpha x} = \alpha\ov{x}. $$
To prove this, it is sufficient to show
$$
       \alpha 1 = 1,\;\; \ov{\alpha e_i} = - \alpha e_i,\quad i=1,2,\cdots,7.$$

\noindent From $(\alpha 1)(\alpha 1) = \alpha(1 \cdot 1) = \alpha 1$ and 
$\alpha 1 \neq 0$, we have $\alpha 1 = 1$, while from  $(\alpha e_i)
(\alpha e_i) = \alpha(e_ie_i) = \alpha(-1) = -\alpha 1 = -1$ for $i \neq 0$, we
 get $\ov{\alpha e_i} =  - \alpha e_i$. Now,
\begin{eqnarray*}
    (\alpha x,\alpha y) \!\!\! &=& \!\!\! \dfrac{1}{2}((\alpha x)(\ov{\alpha y}) + (\alpha y)(\ov{\alpha x})) = \frac{1}{2}((\alpha x)(\alpha \ov{y}) + 
(\alpha y)(\alpha \ov{x})) 
\vspace{1mm}\\
       \!\!\! &=& \!\!\! \alpha\Big(\frac{1}{2}(x\ov{y} + y\ov{x})\Big) = 
\alpha((x,y)) = (x,y).
\end{eqnarray*}

{\bf Theorem 1.2.2.}  $G_2$ {\it is a compact Lie group.}
\vspace{2mm}

{\bf Proof.} $G_2$ is a compact Lie group as a closed subgroup of the 
orthogonal group 
$$
     O(8) = O(\gC) = \{ \alpha \in \Iso_{\sR}(\gC) \, | \, (\alpha x, \alpha y)
 =(x, y) \}. $$

{\bf Remark.} Since $\alpha 1 = 1$ for $\alpha \in G_2$, $G_2$ is a subgroup of the orthogonal group $O(7) = \{ \alpha \in O(\gC) \, | \, \alpha 1 = 1 \}$, that is, $G_2 \subset O(7)$.
\vspace{4mm}

{\bf 1.3. Outer automorphisms of Lie algebra $\gD_4$}
\vspace{3mm}

In order to study the Lie algebra $\gg_2$ of the group $G_2$, we consider the 
Lie algebra
$$ 
        \gD_4 = \so(8) = \so(\gC) = \{ D \in \Hom_{\sR}(\gC) \, | \, (Dx, y) + (x, Dy) =
 0 \}$$
of the Lie group $SO(8)$.
\vspace{2mm}

We define $\R$-linear mappings $G_{ij} : \gC \to \gC,$ $i, j = 0, 1, \cdots , 
7, i \neq j$ satisfying
$$
  G_{ij}e_j = e_i,\;\;\; G_{ij}e_i = -e_j,\;\;\; G_{ij}e_k = 0, \;\, k \neq i,j.$$
Then $G_{ij} \in \gD_4$ and $\{ G_{ij} \, | \, 0 \leq i < j \leq 7 \}$ forms an
 $\R$-basis of $\gD_4$. Furthermore we define $\R$-linear mappings $F_{ij} : 
\gC \to \gC$, $i,j = 0,1,\cdots ,7$, $i \neq j$ by
$$
  F_{ij}x = \displaystyle{\frac{1}{2}}e_i(\overline{e}_jx), \quad x \in \gC. $$

{\bf Lemma 1.3.1.} {\it For $i, j = 0, 1, \cdots, 7, i \neq j$, we have $F_{ij}
 \in \gD_4$, and when $i < j,$ $F_{ij}$ is expressed in terms of $G_{ij}$ as 
follows.}
$$
\begin{array}{ll}
\left\{\begin{array}{l}
      2F_{01} = \, \, \, \,  G_{01} + G_{23} + G_{45} + G_{67} 
\vspace{1mm}\\
      2F_{23} = \, \, \, \,  G_{01} + G_{23} - G_{45} - G_{67} 
\vspace{1mm}\\
      2F_{45} = \, \, \, \,  G_{01} - G_{23} + G_{45} - G_{67} 
\vspace{1mm}\\
      2F_{67} = \, \, \, \,  G_{01} - G_{23} - G_{45} + G_{67},
\end{array}\right. 
&  
\left\{\begin{array}{l}
      2F_{02} = \, \, \, \, G_{02} - G_{13} - G_{46} + G_{57} 
\vspace{1mm}\\
      2F_{13} =      - G_{02} + G_{13} - G_{46} + G_{57} 
\vspace{1mm}\\
      2F_{46} =      - G_{02} - G_{13} + G_{46} + G_{57} 
\vspace{1mm}\\
      2F_{57} = \, \, \, \, G_{02} + G_{13} + G_{46} + G_{57},
\end{array}\right.
\end{array}$$
\vspace{1.5mm} 
$$
\begin{array}{ll}
\left\{\begin{array}{l}
      2F_{03} = \, \, \, \, G_{03} + G_{12} + G_{47} + G_{56} 
\vspace{1mm}\\
      2F_{12} = \, \, \, \, G_{03} + G_{12} - G_{47} - G_{56} 
\vspace{1mm}\\
      2F_{47} = \, \, \, \, G_{03} - G_{12} + G_{47} - G_{56} 
\vspace{1mm}\\
      2F_{56} = \, \, \, \, G_{03} -G _{12} - G_{47} + G_{56},
\end{array}\right. 
&  
\left\{\begin{array}{l}
      2F_{04} = \, \, \, \, G_{04} - G_{15} + G_{26} - G_{37} 
\vspace{1mm}\\
      2F_{15} =      - G_{04} + G_{15} + G_{26} - G_{37} 
\vspace{1mm}\\
      2F_{26} = \, \, \, \, G_{04} + G_{15} + G_{26} + G_{37} 
\vspace{1mm}\\
      2F_{37} =      - G_{04} - G_{15} + G_{26} + G_{37},
\end{array}\right.
\vspace{1.5mm} \\
\left\{\begin{array}{l}
      2F_{05} = \, \, \, \, G_{05} + G_{14} - G_{27} - G_{36} 
\vspace{1mm}\\
      2F_{14} = \, \, \, \,  G_{05} + G_{14} + G_{27} + G_{36} 
\vspace{1mm}\\
      2F_{27} =      - G_{05} + G_{14} + G_{27} - G_{36} 
\vspace{0.5mm}\\
      2F_{36} =      - G_{05} + G_{14} - G_{27} + G_{36},
\end{array}\right. 
&  
\left\{\begin{array}{l}
      2F_{06} = \, \, \, \, G_{06} - G_{17} - G_{24} + G_{35} 
\vspace{1mm}\\
      2F_{17} =      - G_{06} + G_{17} - G_{24} + G_{35} 
\vspace{0.5mm}\\
      2F_{24} =      - G_{06} - G_{17} + G_{24} + G_{35} 
\vspace{1mm}\\
      2F_{35} = \, \, \, \, G_{06} + G_{17} + G_{24} + G_{35},
\end{array}\right.
\end{array}$$
\vspace{-1.5mm}
$$
\left\{\begin{array}{l}
      2F_{07} = \, \, \, \, G_{07} + G_{16} + G_{25} + G_{34} 
\vspace{1mm}\\
      2F_{16} = \, \, \, \, G_{07} + G_{16} - G_{25} - G_{34} 
\vspace{1mm}\\
      2F_{25} = \, \, \, \, G_{07} - G_{16} + G_{25} - G_{34} 
\vspace{1mm}\\
      2F_{34} = \, \, \, \, G_{07} - G_{16} - G_{25} + G_{34}.
\end{array}\right.
$$ 
{\it In particular, $\{ F_{ij} \, | \, 0 \leq i < j \leq 7 \}$ forms an $\R$-basis 
of $\gD_4$.}
\vspace{3mm}

{\bf Definition.} We define $\R$-linear mappings $\kappa,\pi,\nu : \gD_4 \to \gD_4$ respectively by
\begin{eqnarray*}
     (\kappa D)x \!\! &=& \!\! \overline{D\overline{x}}, \quad x \in \gC, \\
      \pi(G_{ij}) \!\! &=& \!\! F_{ij},\quad i,j = 0,1,\cdots, 7, i \neq j, \\
       \nu \!\! &=& \!\! \pi\kappa.
\end{eqnarray*}

{\bf Lemma 1.3.2.} {\it The mappings $\kappa,\pi,\nu$ are automorphisms of the 
Lie algebra $\gD_4$}: 
$$ 
             \kappa,\pi,\nu \in \mbox{Aut}(\gD_4). $$

{\bf Proof.} Since $\kappa^2 = 1$ , $\kappa$ is an $\R$-linear isomorphism of 
$\gD_4$. We have
$$
\begin{array}{l}
    [\kappa D_1,\kappa D_2]x = (\kappa D_1)(\kappa D_2x) - (\kappa D_2)(\kappa 
D_1x) = \overline{D_1(\ov{\kappa D_2x})} - \overline{D_2(\ov{\kappa D_1 x})} 
\vspace{1mm} \\ \qquad \quad
     = \overline{D_1D_2\overline{x}}-\overline{D_2D_1\overline{x}} = 
\kappa(D_1D_2 - D_2D_1)x = \kappa[D_1,D_2]x,\quad x \in \gC. 
\end{array} $$
Hence $\kappa \in \mbox{Aut}(\gD_4)$. Since $\pi$ maps the $\R$-basis 
$\{ G_{ij} \, | \, 0 \le i < j \le 7\}$ to the $\R$-basis $\{ F_{ij} \, | \, 0 \le i < j \le 7\}$ of 
$\gD_4$, $\pi$ is an $\R$-linear isomorphism. To show that $\pi$ is an 
automorphism of $\gD_4$, it is sufficient to show that
$$
         [\pi G_{ij}, \pi G_{kl}] = \pi[G_{ij}, G_{kl}],$$
which in turn would follow from the relations
$$
\begin{array}{ll}
    [F_{ij}, F_{jk}] = F_{ik}, & \quad i,j,k \; \mbox{are distinct}, 
\vspace{1mm}\\
    {[}F_{ij}, F_{kl}{]} = 0,    & \quad i,j,k,l \; \mbox{are distinct}.
\end{array} $$
In the following calculations, let $i,j,k,l$ be all distinct and $i,j,k,l \neq 
0$. For $x \in \gC$,
$$ 
\begin{array}{l}
   [F_{i0}, F_{0k}]x = (F_{i0}F_{0k}- F_{0k}F_{i0})x = - \displaystyle{\frac{1}
{4}}e_i(e_k x) + \dfrac{1}{4}e_k(e_i x)
\vspace{1mm}\\
\qquad \qquad \;\;\;
     = \dfrac{1}{2}e_k (e_i x) = F_{ik}x, 
\end{array}$$
\begin{eqnarray*}
   {[}F_{i0}, F_{kl}{]}x \!\!\! &=& \!\!\! (F_{i0}F_{kl} - F_{kl}F_{i0})x = 
\frac{1}{4}e_i 
       (e_l (e_k x)) - \frac{1}{4}e_l(e_k(e_i x))\\
     \!\!\! &=& \!\!\! -\frac{1}{4}e_l (e_i(e_k x)) + \frac{1}{4}e_l(e_i(e_k x)) = 0,\\
       {[}F_{ij}, F_{jk}{]}x \!\!\! &=& \!\!\! (F_{ij}F_{jk} - F_{jk}F_{ij})x 
= \frac{1}{4}e_j          (e_i(e_k(e_jx))) - \frac{1}{4}e_k(e_j(e_j(e_ix)))\\
     \!\!\! &=& \!\!\! \frac{1}{4}e_i(e_k(e_j(e_jx))) + \frac{1}{4}e_k(e_ix) = 
          -\frac{1}{4}e_i(e_kx) + \frac{1}{4}e_k(e_ix)\\
     \!\!\! &=& \!\!\! \frac{1}{2}e_k(e_ix) = F_{ik}x,\\
       {[}F_{ij}, F_{kl}{]}x \!\!\! &=& \!\!\! (F_{ij}F_{kl} - F_{kl}F_{ij})x = \frac{1}{4}e_j           (e_i(e_l(e_kx))) - \frac{1}{4}e_l(e_k(e_j(e_ix)))\\
     \!\!\! &=& \!\!\! \frac{1}{4}e_j(e_k(e_i(e_lx))) - \frac{1}{4}e_j(e_l(e_k(e_ix)))= \cdots = 0.
\end{eqnarray*} 
Hence $\pi \in \mbox{Aut}(\gD_4)$. Finally, since $\nu = \pi\kappa$, we have $\nu \in \mbox{Aut}(\gD_4)$.
\vspace{2mm}

{\bf Definition.} For $a \in \gC$, we define $\R$-linear mappings $L_a,R_a,T_a 
: \gC \to \gC$ respectively by
$$
     L_ax = ax , \quad R_ax = xa , \quad T_ax = ax + xa = (L_a + R_a)x, \quad 
x \in \gC. $$
Hereafter, we denote by $\gC_0$ the subset $\{ a \in \gC \, | \, \overline{a} 
= -a \}$ of $\gC$.
\vspace{3mm}

{\bf Lemma 1.3.3.} {\it For $a \in \gC_0$, we have}
\vspace{1mm}

(1) \quad $L_a,R_a,T_a \in \gD_4$.
\vspace{1mm}

(2) \quad $\kappa L_a = -R_a, \;\; \kappa R_a = -L_a, \;\,\;\; \kappa T_a = 
-T_a.$ 
\vspace{1mm}

(3) \quad $\pi L_a = T_a, \quad \, \pi R_a = -R_a, \quad \pi T_a = L_a.$ 
\vspace{1mm}

(4) \quad $\nu L_a = R_a, \quad \, \nu R_a = -T_a, \quad \nu T_a = - L_a. $ 
\vspace{2mm}

{\bf Proof.} (1) $(L_a x, y) = (ax, y) = (x, \overline{a}y) = - (x,ay) = 
-(x,L_a y)$, $x,y \in \gC$. Hence $L_a \in \gD_4$. Similarly, $R_a \in \gD_4$ 
and $T_a = L_a + R _a \in \gD_4$.
\vspace{1mm}

(2) $(\kappa L_a)x = \overline{L_a\overline{x}} = \overline{a\overline{x}} = 
x\overline{a} = - xa = - R_ax$, $x \in \gC$. Hence $\kappa L_a = - R_a$. The 
others can be similarly obtained.
\vspace{1mm}

(3) It is sufficient to show that these relations hold for $a = e_i$, $i = 1,
\cdots,7$. We have 
$$
         L_{e_i} = 2F_{i0}, \quad T_{e_i} = 2G_{i0}.$$    
Indeed,  
$$
\begin{array}{l}
     L_{e_i}x = e_ix = 2F_{i0}x, \quad x \in \gC,
\vspace{1mm}\\
     T_{e_i}x = e_ix + xe_i = 
\left\{\begin{array}{ll}
                   2e_i, & \quad x = 1 \\
                   -2,   & \quad x = e_i \\
                   0, & \quad x = e_j,\;\; j \neq 0,i.
\end{array}\right.             
\end{array} $$
It follows that
\vspace{2mm}

\qquad \quad
  $\pi L_{e_i} = \pi(2F_{i0}) = -2 \pi F_{0i} = -2G_{0i}\;\mbox{(Lemma 1.3.1)}\; = 2G_{i0} = T_{e_i}$,
\vspace{1mm}

\qquad \quad
  $\pi T_{e_i} = \pi(2G_{i0}) = 2F_{i0} = L_{e_i}$,
\vspace{1mm}

\qquad \quad
  $\pi R_{a} = \pi(T_a - L_a) = L_a - T_a = - R_a$.
\vspace{2mm}

\noindent \hspace*{5mm}
(4) follows immediately from (2) and (3), since $\nu = \pi\kappa$. 
\vspace{3mm}
         
{\bf Lemma 1.3.4.} {\it The Lie algebra $\gD_4$ is generated by} $\{ L_a \, | 
\, a \in \gC_0 \}$ ({\it by taking at most one time the Lie bracket} 
$[\;\;\;,\;\;\;]$), {\it that is, any $D \in \gD_4$ is expressed by} 
$$
   D = L_a + \sum_{i}[L_{b_i}, L_{c_i}], \quad a,b_i,c_i \in \gC_0.$$

{\bf Proof.} Let $\gD'$ be the Lie subalgebra of $\gD_4$ generated by 
$\{ L_a \, | \, a \in \gC_0 \}$. Since $L_{e_i} = 2F_{i0}$ and $[L_{e_i},
 L_{e_j}] = 4[F_{i0}, F_{j0}] = - 4F_{ij}$, $i \neq 0, j \neq 0, i \neq j$, we 
see that $\gD'$ contains $F_{ij}$, $i,j = 0,1,\cdots,7$, $i \neq j$. Since 
$\{F_{ij}$, $i < j$ \} is an $\R$-basis of $\gD_4$ (Lemma 1.3.1), $\gD'$ 
coincides with $\gD_4$.
\vspace{3mm}

{\bf Theorem 1.3.5.} {\it In the automorphism group $\mbox{Aut}(\gD_4)$ of $\gD_4$}, {\it the subgroup $\gS_3$ 
generated by $\kappa$ and $\pi$ is isomorphic to the symmetric group $S_3$ of 
degree} 3. {\it Furthermore, $\kappa,\pi,\nu$ have the following relations.}
$$
    \kappa^2 = 1, \quad \pi^2 = 1, \quad \nu^3 = 1, \quad \nu =\pi \kappa.$$

{\bf Proof.} Since $\{ L_a \, | \,  a \in \gC_0 \}$ generates $\gD_4$ (Lemma 1.3.4), it is sufficient to check $\kappa^2 = 1$, $\pi^2 = 1$, $\nu^3 = 1$ etc. for $L_a$.  However, these relations follow from Lemma 1.3.3.(2),(3),(4). The 
mapping $f : \gS_3 \to S_3$ defined by the correspondence
$$
\begin{array}{lll}  
        1 \to \pmatrix{1 & 2 & 3 \cr 
                               1 & 2 & 3}, &
        \kappa \to \pmatrix{1 & 2 & 3 \cr
                            2 & 1 & 3}, &
        \pi \to \pmatrix{1 & 2 & 3 \cr
                         3 & 2 & 1}, 
\vspace{1mm} \\
        \nu \to \pmatrix{1 & 2 & 3 \cr
                         2 & 3 & 1}, &
        \nu^2 \to \pmatrix{1 & 2 & 3 \cr
                           3 & 1 & 2}, &
        \nu\pi \to \pmatrix{1 & 2 & 3 \cr
                            1 & 3 & 2}
\end{array} $$
gives an isomorphism as groups.
\vspace{3mm}

{\bf Theorem 1.3.6.} ({\bf Principle of infinitesimal triality in $\gD_4$}). 
{\it For any $D_1 \in \gD_4$, there exist $D_2, D_3 \in \gD_4$ such that
$$
         (D_1x)y + x(D_2y) = D_3(xy), \quad x,y \in \gC.$$
Also such $D_2$, $D_3$ are uniquely determined for $D_1$ and we have} 
$$
        D_2 = \nu D_1, \quad D_3 = \pi D_1.$$

{\bf Proof.} If $a \in \gC_0$, then $L_a,R_a,T_a \in \gD_4$ (Lemma 1.3.3.(1)) 
and the equality $(ax)y + x(ya) = a(xy) + (xy)a$ implies that
$$
    \displaylines{\hfill
         (L_ax)y + x(R_ay) = T_a(xy), \quad x,y \in \gC.
    \hfill\mbox{(i)}}$$
Similarly, for $b \in \gC_0$, we have
$$
    \displaylines{\hfill
         (L_bx)y + x(R_by) = T_b(xy), \quad x,y \in \gC.
    \hfill\mbox{(ii)}}$$
Applying $T_a$ on (ii) and using (i), we have
$$
     (L_aL_bx)y + (L_bx)(R_ay) + (L_ax)(R_by) + x(R_aR_by) =T _aT_b(xy).$$
Exchanging $a$ for $b$, and subtracting from the above, we get 
$$                    
     \displaylines{\hfill
          ([L_a,L_b]x)y + x([R_a,R_b]y) = [T_a,T_b](xy).
     \hfill\mbox{(iii)}}$$
Since $D_1 \in \gD_4$ is expressed as
$$
      D_1 = L_a + \sum[L_b,L_c], \quad a,b,c \in \gC_0 $$
(Lemma 1.3.4), putting $D_2 = R_a + \sum[R_b, R_c]$, $D_3 = T_a + \sum[T_b,T_c]$
 and using (i), (iii), we obtain
$$
         (D_1x)y + x(D_2y) = D_3(xy), \quad x,y \in \gC.$$
Next, we shall show that $D_2$ and $D_3$ are determined uniquely for $D_1$. To 
prove this, it is sufficient to show for $D_1 = 0$ that $D_2 = D_3 = 0$. Now, 
suppose
$$
         x(D_2y) = D_3(xy), \quad x,y \in \gC.$$
Putting $x = 1$, we have $D_2y = D_3y$, so that $D_2 = D_3$ ($= D$). Therefore
$$
      \displaylines{\hfill
          x(Dy) = D(xy), \quad x,y \in \gC.
      \hfill\mbox{(iv)}}$$
Putting $D1 = p$, we have $2(p, 1) = (p, 1) + (1, p) = (D1, 1) + (1, D1) = 0,$ 
which implies $p \in \gC_0$. Furthermore, putting $y = 1$ in (iv), we have $xp =
 Dx$, so (iv) becomes 
$$
        x(yp)=(xy)p,\quad \mbox{for all}\;\; x,y \in \gC.$$ 
We therefore see that $p \in \R$, so that $p = 0$ since $p \in \gC_0$. Hence 
$Dx = xp = 0$, and so $D = 0$. This proves the uniqueness. Finally, if we 
express $D_1 \in \gD_4$ as $D_1 = L_a + \sum[L_b, L_c]$, $a, b, c \in \gC_0$, 
then $D_2 = R_a + \sum[R_b, R_c]$, $D_3 = T_a + $ $\sum[T_b, T_c]$ from the 
arguments above and the uniqueness. Hence we have $D_2 = \nu D_1$, $D_3 = 
\pi D_1$ (Lemma 1.3.3).
\vspace{3mm}

{\bf Lemma 1.3.7.} {\it For $D_1,D_2,D_3 \in \gD_4$, the relation
$$
        (D_1x)y + x(D_2y) = (\kappa D_3)(xy), \quad x,y \in \gC $$ 
implies}
$$
\begin{array}{l}
         (D_2x)y + x(D_3y) = (\kappa D_1)(xy), \quad x,y \in \gC,
\vspace{1mm}\\
         (D_3x)y + x(D_1y) = (\kappa D_2)(xy), \quad x,y \in \gC. 
\end{array} $$

{\bf Proof.} For $D_2 \in \gD_4$, there exist ${D_3}', {D_1}' \in \gD_4$ such 
that$$
      (D_2x)y + x({D_3}'y) = (\kappa {D_1}')(xy), \quad x,y \in \gC, $$ 
and ${D_3}' = \nu D_2$, $\kappa {D_1}' = \pi D_2$ (Theorem 1.3.6). The 
assumption of the lemma is $D_2 = \nu D_1$, 
$\kappa D_3 = \pi D_1$ (Theorem 1.3.6). Hence, using Theorem 1.3.5, we have
$$
\begin{array}{l}
    {D_3}' = \nu D_2 = \nu\nu D_1 = \nu^{-1}D_1 = \kappa\pi D_1 = \kappa\kappa D_3 = D_3,
\vspace{1mm} \\
    {D_1}' = \kappa\pi D_2 = \nu^{-1}D_2 = D_1.
\end{array} $$
\vspace{2mm}

{\bf 1.4. Lie algebra $\gg_2$ of $G_2$}
\vspace{3mm}

{\bf Theorem 1.4.1.} {\it The Lie algebra $\gg_2$ of the Lie group $G_2$ is 
given by}
$$
      \gg_2 = \{ D \in \Hom_{\sR}(\gC) \, | \, D(xy) = (Dx)y + x(Dy) \}. $$

{\bf Proof.} If $D \in \Hom_{\sR}(\gC)$ satisfies $(\exp tD)(xy) = ((\exp tD)x)
((\exp tD)y), t \in \R$, then by differentiating with respect to $t$ and 
putting $t = 0$, we get $D(xy) = (Dx)y + x(Dy)$. Conversely, if $D \in 
\Hom_{\sR}(\gC)$ satisfies $D(xy) = (Dx)y + x(Dy)$, then it is not difficult to
 verify that $\alpha = \exp tD$ satisfies $\alpha(xy) = (\alpha x)(\alpha y)$.
\vspace{2mm}

In order to study the Lie algebra $\gg_2$, we need the following Lie algebra 
$\gb_3 = \so(7)$ of $SO(7)$:
$$
\begin{array}{lll}
      \gb_3 \!\!\! &=& \!\!\! \{ D \in \gD_4 \, | \, D1 = 0 \} 
\vspace{1mm}\\
         \!\!\! &=& \!\!\! \{ D \in \gD_4 \, | \, \kappa D = D \} = \{ D \in 
\gD_4\, | \, \nu D = \pi D \}.
\vspace{3mm}
\end{array}$$

{\bf Lemma 1.4.2.} {\it $\gg_2$ is a Lie subalgebra of $\gb_3$. Moreover we 
have}
\begin{eqnarray*}
    \gg_2 \!\!\! &=& \!\!\! \{ D \in \gD_4 \, | \, \nu D = D, \pi D = D \} \\
  \!\!\! &=& \!\!\! \{ D \in \gD_4 \, | \, \lambda D = D, \lambda \in \gS_3 \} 
\\
  \!\!\! &=& \!\!\! \{ D \in \gb_3 \, | \, \pi D = D \}.
\end{eqnarray*}

{\bf Proof.} Let $D \in \gg_2$. Putting $x = y = 1$ in $D(xy) = (Dx)y + x(Dy)$,
  we get $D1 = 0$. We next show that
$$
            Dx \in \gC_0, \quad x \in \gC.$$
If $i \neq 0$, then $(De_i)e_i + e_i(De_i) = D(e_ie_i) = D(-1) = 0$, and so $De_i 
\in \gC_0$. Together with $D1 = 0$, we have $Dx \in \gC_0, x \in \gC$. Now, note that
$$
     xy + yx = -(x\overline{y} + y\overline{x}) = -2(x,y), \quad x,y \in \gC_0.  $$
Applying $D$ on the relation above, we have $(Dx)y + x(Dy) + (Dy)x + y(Dx) = 0$. Since $Dx,Dy \in \gC_0$, we get
$$
        (Dx, y) + (x, Dy) = 0, \quad x,y \in \gC. $$
Hence $\gg_2 \subset \gb_3$. The relation in the lemma 
follows easily from Theorem 1.3.5 and Theorem 1.3.6.
\vspace{3mm}

{\bf Theorem 1.4.3.} {\it Any element of $\gg_2$ is expressed by the sum of 
elements of the following seven types.}
\vspace{2mm}

\qquad \qquad  
        $\lambda G_{23} + \mu G_{45} + \nu G_{67}, \quad 
       - \lambda G_{13} - \mu G_{46} + \nu G_{57}$, 
\vspace{1mm}

\qquad \qquad 
        $\lambda G_{12} + \mu G_{47} + \nu G_{56}, \quad
    - \lambda G_{15} + \mu G_{26} - \nu G_{37}, \quad \lambda,\mu,\nu \in \R$ 
\vspace{1mm}  

\qquad \qquad  
        $\lambda G_{14} - \mu G_{27} - \nu G_{36}, \quad
     - \lambda G_{17} - \mu G_{24} + \nu G_{35}, \quad \lambda + \mu + \nu = 0$
\vspace{1mm}

\qquad \qquad \qquad \qquad \qquad  
   $\lambda G_{16} + \mu G_{25} + \nu G_{34}$. 
\vspace{2mm}

\noindent {\it Conversely, the above seven elements belong to $\gg_2$. In 
particular, the dimension of $\gg_2$ is} 14:
$$
          \dim\gg_2 = 14. $$

{\bf Proof.} Let $D \in \gg_2$. Since $D \in \gb_3$ (Lemma 1.4.2), $D = \sum_{0<i<j}\lambda_{ij}G_{ij}$, $\lambda_{ij} \in \R$. The condition $\pi D = D$ (Lemma 1.4.2) implies that
$$
      \dsum_{0<i<j}\lambda_{ij} F_{ij} = \sum_{0<i<j}\lambda_{ij}G_{ij}.$$
Applying the above element on $1 \in \gC$, we have
$$
         \dsum_{0<i<j}\lambda_{ij}e_je_i = 0. $$
Replacing $e_je_i$ by $e_k$ and comparing the coefficient of each $e_1,e_2,\cdots,e_7$, we have
\vspace{2mm}

\qquad \qquad \qquad $\;\;$
        $\lambda_{23} + \lambda_{45} + \lambda_{67} = 0, \qquad 
       -\lambda_{13} - \lambda_{46} + \lambda_{57} = 0$, 
\vspace{1mm}

\qquad \qquad \qquad $\;\;$
        $\lambda_{12} + \lambda_{47} + \lambda_{56} = 0, \qquad
       -\lambda_{15} + \lambda_{26} - \lambda_{37} = 0$, 
\vspace{1mm}

\qquad \qquad \qquad $\;\;$
        $\lambda_{14} - \lambda_{27} - \lambda_{36} = 0, \qquad
       -\lambda_{17} - \lambda_{24} + \lambda_{35} = 0$, 
\vspace{1mm}

\qquad \qquad \qquad \qquad \qquad \qquad $\;\;$
        $\lambda_{16} + \lambda_{25} + \lambda_{34} = 0$, 
\vspace{2mm}

\noindent from which the first result follows. Conversely, any of these seven elements $D$ obviously belongs to $\gb_3$ and the condition $\pi D = D$ (Lemma 1.4.2) is verified from the table of Lemma 1.3.1.
\vspace{4mm}

{\bf 1.5. Lie subalgebra $\su(3)$ of $\gg_2$}
\vspace{3mm}

The Cayley algebra $\gC$ naturally contains the field $\C$ of complex numbers 
as$$
          \C = \{ x_0 + x_1e_1 \, | \, x_i \in \R \}.$$
Any element $x \in \gC$ is expressed by
\begin{eqnarray*}
    x \!\!\! &=& \!\!\! x_0 + x_1e_1 + x_2e_2 + x_3e_3 + x_4e_4 + x_5e_5 + 
x_6e_6 + x_7e_7 \quad (x_i \in \R) 
\vspace{1mm}\\
      \!\!\! &=& \!\!\! (x_0 + x _1e_1) + (x_2 + x_3e_1)e_2 + (x_4 + x_5e_1)e_4
 + (x_6 + x_7e_1)e_6,
\end{eqnarray*}
that is,
$$
           x = a + m_1e_2 + m_2e_4 + m_3e_6, \quad a, m_i \in \C.$$
We associate such element $x$ of $\gC$ with the element of $\C \oplus \C^3$ 
$$
           a + \pmatrix{m_1 \cr
                        m_2 \cr
                        m_3}. $$
In $\C \oplus \C^3$, we define a multiplication, an inner product $(\;\;\;,\;\;
\;)$ and a conjugation $\overline{{\;}^{\;}\;\;}$ respectively by
\begin{eqnarray*}
       (a + \m)(b + \n) \!\!\! &=& \!\!\! (ab - \langle \m, \n \rangle ) + 
(a\n + \ov{b}\m - \ov{\m \times \n}), 
\vspace{1mm}\\
        (a + \m, b + \n) \!\!\! &=& \!\!\! (a, b) + (\m, \n), 
\vspace{1mm}\\
        \ov{a + \m} \!\!\! &=& \!\!\! \ov{a} - \m, 
\end{eqnarray*}
where the real valued symmetric inner product $(\m, \n)$, the Hermitian inner 
product $\langle \m, \n \rangle$  and the exterior product $\m \times \n$ are 
usually defined respectively by
$$
     (\m, \n) = \frac{1}{2}(\m^{*}\n + \n^{*}\m) = \sum_{i=1}^3(m_i,n_i), 
\quad      \langle \m, \n \rangle = \sum_{i=1}^3m_i\ov{n}_i, $$
\vspace{-2mm}
$$
          \m \times \n = \pmatrix{m_2n_3 - n_2m_3 \cr
                                  m_3n_1 - n_3m_1 \cr 
                                  m_1n_2 - n_1m_2} $$
for $\m = \pmatrix{m_1 \cr m_2 \cr m_3}$, $\n = \pmatrix{n_1 \cr n_2 \cr n_3} 
\in \C^3$. Since these operations correspond to their respective operations in 
$\gC$, hereafter, we identify $\C \oplus \C^3$ with $\gC$ , that is,
$$
                   \C \oplus \C^3 = \gC. $$

We shall study the following subalgebra $(\gg_2)_{e_1}$ of $\gg_2$:
$$
         (\gg_2)_{e_1} = \{ D \in \gg_2 \, | \, De_1 = 0 \}. $$

{\bf Theorem 1.5.1.} \qquad \quad $(\gg_2)_{e_1} \cong \su(3)$. 
\vspace{2mm}

{\bf Proof.} We define a mapping $\varphi_* : \su(3) = 
\{D \in M(3, \C) \, | \, D^* = - D, \tr(D) = 0 \} \to (\gg_2)_{e_1}$ by
$$
      \varphi_*(D)(a + \m) = D\m, \quad a + \m \in \C \oplus \C^3 = \gC. $$
We first prove that $\varphi_*(D) \in (\gg_2)_{e_1}$. For elements
$$
\begin{array}{llll}
  e_1(E_{11} - E_{22}), & e_1(E_{22} - E_{33}), & E_{12} - E_{21}, & e_1(E_{12} + 
E_{21}), 
\vspace{1mm}\\
  E_{13} - E_{31}, & e_1(E_{13} + E_{31}), & E_{23} - E_{32}, & e_1(E_{23} + 
E_{32})
\end{array} $$
of an $\R$-basis of $\su(3)$ (where $E_{kl} \in M(3, \R)$ is the matrix with 
the $(k,l)$-entry is 1 and otherwise are 0), we have
$$
\begin{array}{ll}
  \varphi_*(e_1(E_{11} - E_{22})) = - G_{23} + G_{45}, & 
  \varphi_*(e_1(E_{22} - E_{33})) = - G_{45} + G_{67}, 
\vspace{1mm}\\
  \varphi_*(E_{12} - E_{21}) = G_{24} + G_{35}, & 
  \varphi_*(e_1(E_{12} + E_{21})) = - G_{25} + G_{34},  
\vspace{1mm}\\
  \varphi_*(E_{13} - E_{31}) = G_{26} + G_{37}, & 
  \varphi_*(e_1(E_{13} + E_{31})) = - G_{27} + G_{36},  
\vspace{1mm}\\
  \varphi_*(E_{23} - E_{32}) = G_{46} + G_{57}, & 
  \varphi_*(e_1(E_{23} + E_{32})) = - G_{47} + G_{56}.
\end{array} $$
Hence $\varphi_*(D) \subset \gg_2$ (Theorem 1.4.3). Clearly $\varphi_*(D)e_1 = 0$, so that $\varphi_*(D) \subset (\gg_2)_{e_1}$. Obviously $\varphi_* : \su(3) \to \gg_2$ is a homomorphism as Lie algebras and is injective, so that we identify $\su(3)$ and $\varphi_*(\su(3)$. We set
$$
\begin{array}{ll}
  S_1 = 2G_{12} - G_{47} - G_{56}, &  S_2 = 2G_{13} - G_{46} + G_{57}, 
\vspace{1mm}\\
  S_3 = 2G_{14} + G_{27} + G_{36}, &  S_4 = 2G_{15} + G_{26} - G_{37}, 
\vspace{1mm}\\
  S_5 = 2G_{16} - G_{25} - G_{34}, &  S_6 = 2G_{17} - G_{24} + G_{35}, 
\end{array}$$
and let $\gS$ be the $\R$-vector subspace of $\gg_2$ spanned by $S_1, \cdots, S_6$. Then we have the following decomposition of $\gg_2$.
$$
        \gg_2 = \su(3) \oplus \gS. $$
Now, we shall show that $\varphi_* : \su(3) \to (\gg_2)_{e_1}$ is onto. Let $B \in (\gg_2)_{e_1}$. Denote
$$
      B = D + \dsum_{i=1}^6x_iS_i, \quad D \in \su(3), x_i \in \R. $$
From the condition $Be_1 = 0$, we have
$$
   -2x_1e_2 - 2x_2e_3 - 2x_3e_4 - 2x_4e_5 - 2x_5e_6 - 2x_6e_7 = 0. $$
Hence $x_1 = \cdots = x_6 = 0$, so that $B = D \in \su(3)$. Thus the proof of Theorem 1.5.1 is completed. 
\vspace{4mm}

{\bf 1.6. Simplicity of ${\gg_2}^C$}
\vspace{3mm}

Let $\gC^C = \{ x_1 + ix_2 \, | \,  x_1, x_2 \in \gC \}$ be the complexification  of the Cayley algebra $\gC$.  In the same manner as in $\gC$, we can also define in $\gC^C$ the multiplication $xy$, the inner product $(x, y)$ such that they satisfy properties 1 $\sim$ 12.3 of Section 1.1 (except some formulas about the length $|x|$). $\gC^C$ is called the complex Cayley algebra. $\gC^C$ has two complex conjugations, namely,
$$
      \ov{x_1 + ix_2} = \ov{x}_1 + i\ov{x}_2,\;\; \tau(x_1 + ix_2) = x_1 - ix_2, \quad x_i \in \gC.$$
The complex conjugation $\tau$ is a complex-conjugate linear transformation of 
$\gC^C$ and satisfies
$$
        \tau(xy) = (\tau x)(\tau y), \quad x, y \in \gC^C. $$

For a while, in the Lie algebra $\su(3) \subset \gg_2$, we use the following notations.
$$
\begin{array}{c}
          H_1 = - G_{23}+ G_{45}, \quad H_2 = - G_{45} + G_{67}
\vspace{1mm}\\
        L_{12} = \;\; G_{24} + G_{35}, \quad L_{21} = - G_{25} + G_{34}, \quad L_{13} = \;\; G_{26} + G_{37},
\vspace{1mm}\\
        L_{21} = - G_{27} + G_{36}, \quad L_{23} = \;\; G_{46} + G_{57},
\vspace{3mm}
 \quad L_{32} = - G_{47} + G_{56}.
\end{array}$$

{\bf Lemma 1.6.1.} {\it The Lie brackets $[D, S], D \in \su(3), S \in \gS$ are given as follows.}
$$
  \matrix{
\begin{array}{rrrrrrrrrrrrrrr}
{} &\vline& S_1 & S_2 & S_3 & S_4 & S_5 & S_6
\\ \hline
          H_1 &\vline & S_2 & - S_1 & - S_4 & S_3 & 0 & 0 
\\
          H_2 &\vline & 0 & 0 & S_4 & - S_3 & - S_6 & S_5
\\ 
           L_{12} &\vline & -S_3 & - S_4 & S_1 & S_2 & 0 & 0 
\\
           L_{21} &\vline & S_4 & - S_3 & S_2 & -S_1 & 0 & 0 
\\
           L_{13} &\vline & - S_5 & - S_6 & 0 & 0 & S_1 & S_2 
\\
           L_{31} &\vline & S_6 & - S_5 & 0 & 0 & S_2 & - S_1 
\\
           L_{23} &\vline & 0 & 0 & - S_5 & - S_6 & S_3 & S_4
\\
           L_{32} &\vline & 0 & 0 & S_6 & - S_5 & S_4 & - S_3 
\end{array}}
$$ 

{\bf Lemma 1.6.2.} $\gS$ {\it is a $\su(3)$-irreducible $\R$-module, and hence we have}
$$
               [\su(3), \gS] = \gS. $$

{\bf Proof.} Evidently $\gS$ is a $\su(3)$-$\R$-module. Let $W$ be a non-zero $\su(3)$-invariant $\R$-snbmodule of $\gS$. If $W$ contains some $S_k$, then, from the table of Lemma 1.6.1, we can see that $W$ contains all $S_k, k = 1, 2, \cdots, 6$, and hence $W = \gS$. Now, let $S = \sum_{k=1}^6x_kS_k$, $x_k \in \R$  be a non-zero element of $W$ and assume that $x_1 \neq 0$ (without the loss of generality). Applying $H_1$ on $S$, we have $x_1S_2 - x_2S_1 - x_3S_4 + x_4S_3 \in W$. Next, applying $L_{13}$ and $L_{31}$ on it, we have
$$  
     - x_1S_6 + x_2S_5 \in W \quad \cdots \mbox{(i)} \quad \mbox{and} \quad 
     - x_1S_5 - x_2S_6 \in W \quad \cdots \mbox{(ii)} $$
Taking (ii)$\times x_2 -$ (i)$\times x_1$, we have $({x_1}^2 + {x_2}^2)S_5 \in W$. Since ${x_1}^2 + {x_2}^2 \ne 0$, we have $S_5 \in W$ and so $W = \gS$. Consequently the irreducibility of $\gS$ is proved. Finally, since $[\su(3), \gS]$ is a $\su(3)$-invariant $\R$-submodle of $\gS$, from the irreducibility of $\gS$, we have $[\su(3), \gS] = \gS$. 
\vspace{3mm}

{\bf Theorem 1.6.3.} {\it The Lie algebra ${\gg_2}^C$ is simple and so $\gg_2$ is also simple.}
\vspace{2mm}

{\bf Proof.}  We shall prove that $\gg_2$ is simple, because the proof of that ${\gg_2}^C$ is simple is the same as the case of $\gg_2$. We use the decomposition of $\gg_2$
$$
        \gg_2 = \su(3) \oplus \gS \quad \mbox{(Theorem 1.5.1)}. $$
Let $p : \gg_2 \to \su(3)$ and $q : \gg_2 \to \gS$ be projections of $\gg_2 = \su(3) \oplus \gS$. Now, let $\ga$ be a non-zero ideal of $\gg_2$. Then $p(\ga)$ is an ideal of $\gg_2$. Indeed, if $D \in p(\ga)$, then there exists $S \in \gS$ such that $D + S \in \ga$. For any $D' \in \su(3)$, we have 
$$
         \ga \ni [D', D + S] = [D',D] + [D', S], \quad [D', S] \in \gS $$
(Lemma 1.6.1), hence $[D', D] \in p(\ga)$. 
\vspace{1mm}

We show that either $\su(3) \cap \ga \neq \{ 0 \}$ or $\gS \cap \ga \neq \{ 0 \}$. Assume that $\su(3) \cap \ga = \{ 0 \}$ and $\gS \cap \ga = \{ 0 \}$. The 
mapping $p|\ga : \ga \to \su(3)$ is injective because $\gS \cap \ga = \{ 0 \}$. Since $p(\ga)$ is a non-zero ideal of $\su(3)$ and $\su(3)$ is simple, we have $p(\ga) = \su(3)$. Hence $\dim\ga = \dim p(\ga) = \dim\su(3) = 8$. On the other hand, since $\su(3) \cap \ga = \{ 0 \}$, $q|\ga : \ga \to \gS$ is also injective, we have $\dim\ga \le \dim\gS = 6$. This leads to a 
\vspace{1mm}
contradiction.

We now consider the following two cases.
\vspace{1mm}

(1) Case $\su(3) \cap \ga \neq \{ 0 \}$. From the simplicity of $\su(3)$, we have $\su(3) \cap \ga = \su(3)$, hence $\ga \supset \su(3)$. On the other hand, we have
$$
    \ga \supset [\ga, \gS] \supset [\su(3), \gS] = \gS \;\; \mbox{(Lemma 1.6.1)}. $$
Hence $\ga \supset \su(3) \oplus \gS = \gg_2$.
\vspace{1mm}

(2) Case $\gS \cap \ga \neq \{ 0 \}$. Choose a non-zero element $S \in \gS \cap \ga \subset \ga$. Under the actions of $\su(3)$, we can see that $S_1 \in \ga$ (Lemma 1.6.1). Hence $ 0 \neq 4H_1 + 2H_2 = [S_1, S_2] \in \ga$. So this case can be reduced to the case (1). Thus we have $\ga = \gg_2$, which proves the simplicity of $\gg_2$.   
\vspace{4mm}

{\bf 1.7. Killing form of ${\gg_2}^C$}
\vspace{3mm}

{\bf Lemma 1.7.1.} {\it In $\su(3) \subset \gg_2$, the Lie brackets between $H_1$ and $L_{ij}, L_{ji}$ of Section} 1.6 {\it are given by}
$$
\begin{array}{llll}
     [H_1, L_{12}] = 2L_{21}, &  [H_1, L_{21}] = - 2L_{12}, & [H_1, L_{13}] = L_{31},
\vspace{1mm}\\
     {[}H_1, L_{31}] = - L_{13}, &  [H_1, L_{23}] = - L_{32}, & [H_1, L_{32}] = L_{23}.
\end{array}$$ 

{\bf Theorem 1.7.2.} {\it The killing form $B_2$ of the Lie algebra ${\gg_2}^C$ is given by}
$$
         B_2(D_1, D_2) = 4\tr(D_1D_2), \quad D_i \in {\gg_2}^C.$$

{\bf Proof.} Since $\tr(D_1D_2)$ is a ${\gg_2}^C$-adjoint invariant bilinear form of ${\gg_2}^C$ and ${\gg_2}^C$ is simple (Theorem 1.6.3), there exists $k \in C$ such that
$$
         B_2(D_1, D_2) = k\tr(D_1D_2). $$
To determine this $k$, let $D_1 = D_2 = H_1$. Then from
$$
\begin{array}{ll}
   [H_1, [H_1, L_{12}]\,] = [H_1, 2L_{21}] = - 4L_{12}, &
   [H_1, [H_1, L_{21}]\,] = [H_1, -2L_{12}] = - 4L_{21},
\vspace{1mm}\\
   {[}H_1, [H_1, L_{13}]\,] = [H_1, L_{31}] = - L_{13}, &
   [H_1, [H_1, L_{31}]\,] = [H_1, - L_{13}] = - L_{31}, 
\vspace{1mm}\\
    {[}H_1, [H_1, L_{23}]\,] = [H_1, - L_{32}] = - L_{23}, &
    [H_1, [H_1, L_{32}]\,] = [H_1,  L_{23}] = - L_{32}, 
\end{array}$$
$$\begin{array}{ll}
     {[}H_1, [H_1, S_1]\,] = [H_1, S_2] = - S_1, &
     [H_1, [H_1, S_2]\,] = [H_1, - S_1] = - S_2,
\vspace{1mm}\\
     {[}H_1, [H_1, S_3]\,] = [H_1, - S_4] = - S_3, &
     [H_1, [H_1, S_4]\,] = [H_1, S_3] = - S_4,
\vspace{1mm}\\
     {[}H_1, [H_1, S_5]\,] = 0, \; [H_1, [H_1, S_6]\,] = 0
\end{array}$$
(Lemma 1.6.1), we have
$$
    B_2(H_1, H_1) = \tr((\mbox{ad}H_1)^2) = (-4) \times 2 + (-1) \times 8 = - 16. $$
On the other hand,
$$
\begin{array}{lll}
    H_1H_1e_2 = H_1e_3 = - e_2, &  H_1H_1e_3 = - H_1e_2 = - e_3, 
\vspace{1mm}\\
    H_1H_1e_4 = - H_1e_5 = - e_4, &  H_1H_1e_5 = H_1e_4 = - e_, 
\vspace{1mm}\\
   \, H_1H_1e_i = 0 \;\; \mbox{otherwise}.
\end{array}$$
Hence $\tr(H_1H_1) = (-1) \times 4 = - 4$. Therefore $k = 4$. 
\vspace{4mm}

{ \bf 1.8. Roots of ${\gg_2}^C$}
\vspace{3mm}

We recall the $C$-Lie isomorphism $f_* : \sl(3, C) \to \su(3)^C$ and the embedding $\varphi_* : \su(3)^C \to {\gg_2}^C$,
$$
\begin{array}{l}
        f_*(A) = \varepsilon A - \ov{\varepsilon}\,{}^t\!A, \qquad \;\; \varepsilon = \dfrac{1}{2}(1 + ie_1),
\vspace{1mm}\\
       \varphi_*(D)(a + \m) = D\m, \quad  a + \m \in \C^C \oplus (\C^3)^C = \gC^C,
\end{array} $$
and we regard $\sl(3, C)$ as a subalgebra of ${\gg_2}^C$ under the composition of $f_*$ and $\varphi_*$. Further we know that the Lie algebra $\sl(3, C)$ has roots $\pm(\lambda_k - \lambda_l), 1 \le k < l \le 3$ relative to the Cartan subalgebra 
$$\gh =
      \Big\{\pmatrix{\lambda_1 & 0 & 0 \cr
               0 & \lambda_2 & 0 \cr
               0 & 0 & \lambda_3} \Big| \, \lambda_i \in \C, \lambda_1 + \lambda_2 + \lambda_3 = 0 \Big\}$$
 of $\sl(3, \C)$, and $E_{kl}$ is a root vector associated with the root $\lambda_k - \lambda_l$.
\vspace{3mm}

{\bf Theorem 1.8.1.} {\it The rank of the Lie algebra ${\gg_2}^C$ is} 2. {\it The roots of ${\gg_2}^C$ relative to some Cartan subalgebra of ${\gg_2}^C$ are given by}
$$
   \pm(\lambda_1 - \lambda_2), \quad  \pm(\lambda_1 - \lambda_3), \quad \pm(\lambda_2 - \lambda_3), 
   \quad \pm\lambda_1, \quad \pm\lambda_2 \quad \pm\lambda_3 $$
with $\lambda_1 + \lambda_2 + \lambda_3 = 0$.
\vspace{2mm}
                                    
{\bf Proof.} \quad $\gh = \{ - i\lambda_1G_{23} - i\lambda_2G_{45} - i\lambda_3G_{67} \in {\gg_2}^C \, | \, \lambda_k \in C \} \subset \sl(3, C) \subset {\gg_2}^C$
\vspace{1mm}

\noindent is an abelian subalgebra of ${\gg_2}^C$ (it will be a Cartan subalgebra of ${\gg_2}^C$). The roots of $\sl(3, C)$ is also roots of ${\gg_2}^C$, so we have the table of roots and associated root vectors as follows.
\begin{eqnarray*}
    \pm(\lambda_1 - \lambda_2) &:& \pm(G_{24} + G_{35}) + i(- G_{25} + G_{34}),
\vspace{1mm}\\
    \pm(\lambda_1 - \lambda_3) &:& \pm(G_{26} + G_{37}) + i(- G_{27} + G_{36}),
\vspace{1mm}\\
    \pm(\lambda_2 - \lambda_3) &:& \pm(G_{46} + G_{57}) + i(- G_{47} + G_{56}).
\end{eqnarray*}
The remainder roots and associated root vectors are found as follows.
\begin{eqnarray*}
 \pm\lambda_1 &:& (2G_{12} - G_{47} - G_{56}) \pm i(2G_{13} - G_{46} + G_{57}),
\vspace{1mm}\\
 \pm\lambda_2 &:& (2G_{14} + G_{27} + G_{36}) \pm i(2G_{15} + G_{26} - G_{37}),
\vspace{1mm}\\
 \pm\lambda_3 &:& (2G_{16} - G_{25} - G_{34}) \pm i(2G_{17} - G_{24} + G_{35}).
\end{eqnarray*}

{\bf Theorem 1.8.2.} {\it In the root system of ${\gg_2}^C$ of Theorem} 1.8.1, 
$$
         \alpha_1 = \lambda_1 - \lambda_2, \quad \alpha_2 = \lambda_2 $$
{\it is a fundamental root system of the Lie algebra ${\gg_2}^C$ and
$$
      \mu = 2\alpha_1 + 3\alpha_3 $$
is the highest root. The Dynkin diagram and the extended Dynkin diagram of ${\gg_2}^C$ are respectively given by}

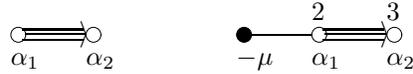
\setlength{\unitlength}{1mm}
\begin{picture}(100,20)
\put(30,10){\circle{2}} \put(29,6){$\alpha_1$}
\put(30.6,10.7){\line(1,0){8}}
\put(30.9,10.1){\line(1,0){8}}
\put(30.6,9.3){\line(1,0){8}}
\put(38.1,9.2){$\rangle$}
\put(40,10){\circle{2}} \put(39,6){$\alpha_2$}
\put(60,10){\circle*{2}} \put(59,6){$-\mu$}
\put(60,10){\line(1,0){9}}
\put(70,10){\circle{2}} \put(69,6){$\alpha_1$} \put(69,12){$2$}
\put(70.5,10.7){\line(1,0){8}}
\put(70.8,10.1){\line(1,0){8}}
\put(70.5,9.3){\line(1,0){8}}
\put(78.1,9.2){$\rangle$}
\put(80,10){\circle{2}} \put(79,6){$\alpha_2$} \put(79,12){$3$}
\end{picture}
\vspace{1mm}

{\bf Proof.} All positive roots of ${\gg_2}^C$ are expressed by
$$
\begin{array}{cccc}
     \lambda_1 - \lambda_2 = \alpha_1, &  \;\;\; \lambda_1 - \lambda_3 = 2\alpha_1 + 3\alpha_2, &  \lambda_2 - \lambda_3 = \alpha_1 + 3\alpha_2, 
\vspace{1mm}\\
   \qquad \qquad \;\;\lambda_1 = \alpha_1 + \alpha_2, & \lambda_2 = \alpha_2, &\;\;\;\;\; - \lambda_3 = \alpha_1 + 2\alpha_2.
\end{array}$$ 
Hence ${\mit\Pi} = \{ \alpha_1, \alpha_2 \}$ is a fundamental root system of ${\gg_2}^C$. The real part of $\gh_{\sR}$ of $\gh$ is
$$
    \gh_{\sR} = \{ -i\lambda_1G_{23} -i\lambda_2G_{45} -i\lambda_3G_{67} \, |\, \lambda_i \in \R, \lambda_1 + \lambda_2 + \lambda_3 = 0 \} $$
and the Killing form $B_2$ on $\gh_{\sR}$ is given by
$$
         B_2(H, H') = 8\dsum_{k=1}^3\lambda_k{\lambda_k}', $$
for $H = - i\lambda_1G_{23} - i\lambda _2G_{45} - i\lambda_3G_{67}, 
    H' = - i{\lambda_1}'G_{23} - i{\lambda_2}'G_{45} - i{\lambda_3}'G_{67} \in \gh_{\sR}$ (Theorem 1.7.2). Now, the canonical element $H_{\alpha_i} \in \gh_{\sR}$ associated with $\alpha_i \;(B_2(H_{\alpha}, H) = \alpha(H), H \in \gh_{\sR})$ are determined as follows.
$$
        H_{\alpha_1} = - \dfrac{1}{8}iG_{23} + \dfrac{1}{8}iG_{45},
\quad
       H_{\alpha_2} = \dfrac{1}{24}iG_{23} - \dfrac{1}{12}iG_{45} + \dfrac{1}{24}iG_{67}.$$
Hence we have
$$
\begin{array}{l}
    (\alpha_1, \alpha_1) = B_2(H_{\alpha_1}, H_{\alpha_1}) = 8\dfrac{1}{8}\dfrac{1}{8}(1 + 1) = \dfrac{1}{4},
\vspace{1.5mm}\\
    (\alpha_2, \alpha_2) = B_2(H_{\alpha_2}, H_{\alpha_2}) = 8\dfrac{1}{24}\dfrac{1}{24}(1 + 4 + 1) = \dfrac{1}{12},
\vspace{1.5mm}\\
    (\alpha_1, \alpha_2) = B_2(H_{\alpha_1}, H_{\alpha_2}) = 8\dfrac{1}{8}\dfrac{1}{24}(-1 - 2) = - \dfrac{1}{8},
\vspace{1mm}\\
     (-\mu, -\mu) = \dfrac{1}{4}, \quad (- \mu, \alpha_1) = - \dfrac{1}{8}, \quad (- \mu, \alpha_2) = 0.
\end{array}$$
Using them, we can draw the Dynkin diagram and the extended Dynkin diagram of 
\vspace{2mm}
${\gg_2}^C$.     

According to Borel-Siebenthal theory (Borel and Siebenthal [4]), the Lie algebra $\gg_2$ has two subalgebras as maximal subalgebras with the maximal rank 2.
\vspace{1mm}

(1) One is a subalgebra of type $C_1 \oplus C_1$ which is obtained as the fixed points by an involution $\gamma$ of $\gg_2$.
\vspace{1mm}

(2) The other is a subalgebra of type $A_2$ which is obtained as the fixed points by an automorphism $w$ of order 3 of $\gg_2$.
\vspace{1mm}

These subalgebras will be realized in the group $G_2$ in Theorem 1.10.1 and Theorem 1.9.4, respectively.
\vspace{4mm}

{\bf 1.9. Automorphism $w$ of order 3 and subgroup $SU(3)$ of $G_2$}
\vspace{3mm}

We shall study the following subgroup $(G_2)_{e_1}$ of $G_2$:
$$
        (G_2)_{e_1} = \{\alpha \in G_2 \, | \, \alpha e_1 = e_1 \}. $$

{\bf Theorem 1.9.1.} \qquad \qquad $(G_2)_{e_1} \cong SU(3).$
\vspace{2mm}

{\bf Proof} (cf. Theorem 1.5.1). We define a mapping $\varphi : SU(3) \to (G_2)_{e_1}$ by
$$
        \varphi(A)(a + \m) = a + A\m, \quad
         a + \m \in \C \oplus \C^3 = \gC. $$
We first prove that $\varphi(A) \in (G_2)_{e_1}$. For $\alpha = \varphi(A), A \in SU(3)$ and $x = a + \m, y = b + \n \in \C \oplus \C^3 = \gC,$  using that if $A \in SU(3)$ then $\wti{A}$ (which is the adjoint matrix of $A$) $= A^{-1} = A^*$, we have
\begin{eqnarray*}
     (\alpha x)(\alpha y)
     \!\!\! &=& \!\!\! (a + A\m)(b + A\n)
\vspace{1mm}\\
     \!\!\! &=& \!\!\!(ab - \langle A\m, A\n \rangle) + (aA\n + \ov{b}A\m - \overline{A\m \times A\n})
\vspace{1mm}\\
     \!\!\! &=& \!\!\! (ab - \langle \m, A^*A\n \rangle) + (aA\n + \ov{b}A\m - \overline{^t\wti{A}(\m \times \n}))
\vspace{1mm}\\
     \!\!\! &=& \!\!\! (ab - \langle \m, \n \rangle) + A(a\n + \bar{b}\m - \overline{\m \times \n}) 
\vspace{1mm}\\
     \!\!\! &=& \!\!\! \varphi(A)((a + \m)(b + \n)) = \alpha(xy).                  \end{eqnarray*}       
Hence $\varphi(A) \in G_2.$ Clearly $\varphi(A)e_1 = e_1,$ so that $\varphi(A) \in (G_2)_{e_1}.$  Evidently $\varphi$ is a homomorphism. We show that $\varphi$ is onto. Let $\alpha \in (G_2)_{e_1}$. Note that $\alpha$ induces a \C-linear transformation of $\C^3$. Now let 
$$
              \alpha e_2 = \a_1, \quad \alpha e_4 = \a_2, \quad
              \alpha e_6 = \a_3 $$
and construct a matrix $A = (\a_1, \a_2, \a_3) \in M(3, \C)$.  From $(\alpha e_2)(\alpha e_4) = \alpha(e_2e_4) = - \alpha e_6$, we have $\a_1\a_2 = - \a_3$, namely, $- \langle \a_1, \a_2 \rangle - \overline{\a_1 \times \a_2} = - \a_3$, hence we have
$$
       \langle \a_1, \a_2 \rangle = 0, \quad
        \a_3 = \ov{\a_1 \times \a_2}. $$
\noindent Similarly we have $\langle \a_2, \a_3 \rangle = \langle \a_3, \a_1 \rangle = 0$. Moerover from $(\alpha e_k)(\alpha e_k) = \alpha(e_ke_k) = \alpha(-1) = -1$, we have $\langle \a_k, \a_k \rangle = 1$, hence $A \in U(3)$. Finally $\det\,A = (\a_3, \a_1 \times \a_2) = (\a_3, \ov{\a}_3) = \langle \a_3, \a_3 \rangle = 1$ (where $(\a, \b)$ is the ordinary inner product in $\C^3$: $(\a, \b) = {}^t\a\b$). Hence we have $A \in SU(3)$ and $\varphi(A) = \alpha$, which shows that $\varphi$ is onto. $\Ker\,\varphi = \{ E \}$ is easily obtained. Thus we have the isomorphism $SU(3) \cong (G_2)_{e_1}$.    
\vspace{3mm}

{\bf Theorem 1.9.2.} \qquad \qquad $G_2/SU(3) \simeq S^6. $
\vspace{2mm}

{\bf Proof.} $S^6 = \{a \in \gC \, | \, \bar{a} = - a, |a| = 1 \}$ is a 6 dimensional sphere. Since the group $G_2$ is a subgroup of $O(7) = \{\alpha \in O(\gC) \, | \, \alpha 1 = 1 \}$ (Remark of Theorem 1.2.2), $G_2$ acts on $S^6$. We shall show that this action is transitive. To prove this, it is sufficient to show that any element $a \in S^6$ can be transformed to $e_1 \in S^6$ by some $\alpha \in G_2$. Now, for $a_1 \in S^6$, choose any element $a_2 \in S^6$ such that $(a_1, a_2) = 0$. Let
$$
                   a_3 = a_1a_2.  $$
Then $a_3 \in S^6$ and $a_3$ satisfies $(a_1, a_3) = (a_2, a_3) = 0$. Choose any element $a_4 \in S^6$ such that $(a_1, a_4) = (a_2, a_4) = (a_3, a_4) = 0$. Let$$
        a_5 = a_1a_4, \quad a_6 = a_4a_2, \quad a_7 = a_3a_4.$$
Then the set $\{a_0 = 1,a_1,a_2,\cdots,a_7 \}$ is an orthonormal $\R$-basis of $\gC$. Indeed, $|a_i| = 1, 0 \leq i \leq 7$ are trivial, and we need to verify that $(a_i,a_j) = 0, i \neq j$. However this can be checked by direct calculations such as
$$
\begin{array}{l} 
    (a_4, a_7) = (a_4, a_3a_4)= (1, a_3)(a_4, a_4) = 0, 
\vspace{1mm}\\
    (a_1, a_6) = (a_1, a_4a_2) = -(a_1a_2, a_4) = -(a_3, a_4) = 0,
\vspace{1mm}\\
    (a_3, a_6) = (a_3, a_4a_2)= -(a_3a_2, a_4) = (a_1, a_4) = 0,\, \; \mbox{etc.}
\end{array} $$
Now, since $\{e_0 = 1,e_1,e_2,\cdots,e_7 \}$ and $\{a_0 = 1, a_1,a_2,\cdots,a_7 \}$ are both orthonormal $\R$-bases, the $\R$-linear isomorphism $\alpha : \gC \to \gC$ satisfying 
$$
            \alpha e_i =a_i, \quad i = 0,1,\cdots,7$$
belongs to $O(7)$: $\alpha \in O(7)$. Moreover we claim that $\alpha \in G_2$, that is, $\alpha$ satisfies
$$
         \alpha(xy) = (\alpha x)(\alpha y), \quad x,y \in \gC.$$
To show this, it is sufficient to note that
$$
       \alpha(e_ie_j) = (\alpha e_i)(\alpha e_j), \quad i,j = 0,1,\cdots,7  $$
which can be also checked by direct calculations such as
\begin{eqnarray*}
    (\alpha e_4)(\alpha e_7) \!\!\! &=& \!\!\! a_4 a_7 = a_4(a_3 a_4) = - a_4(a_4a_3) = a_3 = \alpha e_3 = \alpha (e_4e_7),
\vspace{1mm}\\
    (\alpha e_1)(\alpha e_6) \!\!\! &=& \!\!\! a_1 a_6 = a_1 (a_4 a_2)= - a_4(a_1a_2) = - a_4 a_3 = a_3a_4 = a_7 
\vspace{1mm}\\
                             \!\!\! &=& \!\!\! \alpha e_7 = \alpha(e_1e_6),\\
    (\alpha e_3)(\alpha e_6) \!\!\! &=& \!\!\! a_3a_6 = (a_1a_2)(a_4a_2) = -(a_2 a_1)(a_4 a_2) = - a_2(a_1a_4)a_2 
\vspace{1mm}\\
   \!\!\! &=& \!\!\! - a_2a_5a_2 = a_2a_2a_5 = - a_5 = - \alpha e_5 = \alpha(e_3e_6), \,\; \mbox{etc.}
\end{eqnarray*}
Hence $\alpha \in G_2$ and $\alpha e_1 = a_1$, and so $\alpha^{-1}a_1 = e_1$.  This shows the transitivity. The isotropy subgroup $(G_2)_{e_1}$ of $G_2$ at $e_1$ is $SU(3)$ (Theorem 1.9.1). Thus we have the homeomorphism $G_2/SU(3) \simeq S^6$.
\vspace{2mm}

Since $S^6$ and $SU(3)$ are both simply connected, from $G_2/SU(3) \simeq S^6$ (Theorem 1.9.2), we see that $G_2$ is also simply connected. Hence we have the 
\vspace{3mm}
following theorem.

{\bf Theorem 1.9.3.} {\it $G_2 = \{ \alpha \in \Iso_{\sR}(\gC) \, | \, \alpha(xy)=(\alpha x)(\alpha y) \}$ is a simply connected compact simple Lie group}.
\vspace{3mm}

{\bf Remark 1.} Since $G_2$ is connected, $G_2$ is contained in $SO(7)=\{ \alpha \in O(7) \, | \, \det\,\alpha = 1 \}$: $G_2 \subset SO(7)$.
\vspace{1mm}

{\bf Remark 2.} Since we know that the dimension of the group $G_2$ as $\dim G_2 = \dim\gg_2 = 14$ (Theorem 1.4.3), $G_2/SU(3) \simeq S^6$ is proved as follows. The group $G_2$ acts on $S^6$. The isotropy subgroup $(G_2)_{e_1}$ of $G_2$ at $e_1$ is $SU(3)$ (Theorem 1.9.1) and $\dim(G_2/(G_2)_{e_1}) = \dim G_2 - \dim SU(3) = 14 - 8 = 6 = \dim S^6$. Therefore we have $G_2/SU(3) \simeq S^6.$ 
\vspace{2mm}

Using the mapping $\varphi : SU(3) \to G_2$, we define a mapping $w : \gC \to \gC$ by
$$
         w = \varphi(\diag(\omega_1, \omega_1, \omega_1)) $$
where $\omega_1 = - \displaystyle{\frac{1}{2}} + \displaystyle{\frac{\sqrt{3}}{2}}e_1 \in \C \subset \gC$. This $w$ is defined as
$$
     w(a + \m) = a + \omega_1\m, \quad a + \m \in \C \oplus \C^3 = \gC. $$
Then $w \in G_2$ and $w^3 = 1$.        
\vspace{2mm}

We shall study the following subgroup $(G_2)^w$ of $G_2$:
$$
     (G_2)^w = \{ \alpha \in G_2 \, | \, w\alpha = \alpha w \}. $$

{\bf Theorem 1.9.4.} \qquad \quad $(G_2)^w = (G_2)_{e_1} \cong SU(3).$
\vspace{2mm}

{\bf Proof.} Recall the mapping $\varphi : SU(3) \to G_2$ of Theorem 1.9.1. We first show that $\varphi(SU(3)) \in (G_2)^w$. Indeed, for $A \in SU(3)$ and $a + \m \in \C \oplus \C^3 = \gC$, we have
\begin{eqnarray*}
      w\varphi(A)(a + \m) \!\!\! &=& \!\!\! w(a + A\m) = a + \omega_1A\m \\
          \!\!\! &=& \!\!\! a + A\omega_1\m = \varphi(A)w(a + \m).
\end{eqnarray*}
Hence $w\varphi(A) = \varphi(A)w$, so that $\varphi(A) \in (G_2)^w$. Conversely, let $\alpha \in (G_2)^w$. We consider the $\R$-vector subspace $\gC_w = \{ x \in \gC \, | \, wx = x \}$ of $\gC$, then $\gC_w = \C$. Since $\alpha$ satisfies $w\alpha = \alpha w$, $\gC_w$ is invariant under $\alpha$. Since the restriction of $\alpha$ to $\gC_w$ induces an automorphism of $\C$, we have
$$
        \alpha e_1 = e_1 \quad \mbox{or} \quad \alpha e_1 = - e_1. $$
In the latter case, consider the mapping $\gamma_1 : \gC \to \gC$ defined by $\gamma_1(a + \m) = \overline{a} + \overline{\m}$. Then we have $\gamma_1 \in G_2$ and $\gamma_1e_1 = - e_1$. Let $\beta = \gamma_1\alpha$. Since $\beta e_1 = e_1$, we have $\beta \in SU(3) \subset (G_2)_{e_1}$ (Theorem 1.9.1) $\subset (G_2)^w$. Therefore, $\gamma_1 = \beta\alpha^{-1} \in (G_2)^w$, which is a contradiction. Indeed, 
$$
     \omega_1\overline{\m} = w(\gamma_1\m) = \gamma_1(w\m) = \overline{\omega_1\m} = \overline{\omega}_1\overline{\m} \quad \mbox{for all $\m \in \C^3$,} $$
which is false. Therefore $\alpha e_1 = e_1$, and so $\alpha \in SU(3)$ (Theorem 1.9.1). Thus we have $(G_2)^w = (G_2)_{e_1}$.
\vspace{4mm}

{ \bf 1.10. Involution $\gamma$ and subgroup $(Sp(1) \times Sp(1))/\Z_2$ of $G_2$}
\vspace{3mm}

The Cayley algebra $\gC$ naturally contains the field $\H$ of quaternions as
$$
      \H = \{ x_0 + x_1e_1 + x_2e_2 + x_3e_3 \, | \, x_i \in \R \}.$$
Any element $x \in \gC$ is expressed by
\begin{eqnarray*}
    x \!\!\! &=& \!\!\! x_0 + x_1e_1 + x_2e_2 + x_3e_3 + x_4e_4 + x_5e_5 + x_6e_6 + x_7e_7 \quad (x_i \in \R) 
\vspace{1mm}\\
      \!\!\! &=& \!\!\! (x_0 + x _1e_1 + x_2e_2 + x_3e_3) + (x_4 + x_5e_1 - x_6e_2 + x_7e_3)e_4,
\end{eqnarray*}
that is,
$$
           x = m + ae_4, \quad m,a \in \H.$$
In $\H \oplus \H e_4$, we define a multiplication, an inner product $(\;\;\;,\;\;\;)$, a conjugation $\overline{{\;}^{\;}\;\;}$ and an $\R$-linear transformation $\gamma$ respectively by
\begin{eqnarray*}
       (m + ae_4)(n + be_4) \!\!\! &=& \!\!\! (mn - \overline{b}a) + (a\overline{n} + bm)e_4,\\
        (m + ae_4, n + be_4) \!\!\! &=& \!\!\! (m, n)+(a, b),\\
        \overline{m + ae_4} \!\!\! &=& \!\!\! \overline{m} - ae_4, \\
        \gamma(m + ae_4) \!\!\! &=& \!\!\! m - ae_4.
\end{eqnarray*}
Since these operations correspond to their respective operations in $\gC$, hereafter, we identify $\H \oplus \H e_4$ with $\gC$, that is,
$$
             \H \oplus \H e_4 =  \gC. $$

We shall study the following subgroup $(G_2)^{\gamma}$ of $G_2$:
$$
  (G_2)^{\gamma} = \{ \alpha \in G_2 \, | \, \gamma \alpha =\alpha \gamma \}.$$

{\bf Theorem 1.10.1.} \quad $(G_2)^{\gamma} \cong (Sp(1) \times Sp(1))/\Z_2, \, \, \Z_2 = \{ (1, 1), (-1, -1) \}. $
\vspace{2mm}

{\bf Proof.} We define a mapping $\varphi : Sp(1) \times Sp(1) \to (G_2)^{\gamma}$ by
$$
      \varphi(p,q)(m + ae_4) = qm\overline{q} + (pa\overline{q})e_4, \quad m + ae_4 \in \H \oplus \H e_4 = \gC.$$
We first show that $\varphi(p,q) \in (G_2)^{\gamma}$. For $\alpha = \varphi(p,q),$  $p, q \in Sp(1)$ and $x = m + ae_4$, $y = n + be_4 \in \H \oplus \H e_4 = \gC$, we have
\begin{eqnarray*}
        (\alpha x)(\alpha y) \!\!\! &=& \!\!\! 
 (qm\overline{q} + (pa\overline{q})e_4)(qn\overline{q} + (pb\overline{q})e_4)
\vspace{1mm}\\
         \!\!\! &=& \!\!\! ((qm\overline{q})(qn\overline{q}) - 
                       (\overline{pb\overline{q}})(pa\overline{q})) +
       ((pa\overline{q})(\overline{qn\overline{q}})+(pb\overline{q})(qm\overline{q}))e_4 
\vspace{1mm}\\
         \!\!\! &=& \!\!\! q(mn - \overline{b}a)\overline{q} + (p(a\overline{n} + bm)\overline{q})e_4 
\vspace{1mm}\\
         \!\!\! &=& \!\!\! \varphi(p, q)((m + ae_4)(n + be_4)) = \alpha(xy).
\end{eqnarray*}
Hence $\varphi(p,q) \in G_2$. Clearly $\gamma\varphi(p,q) = \varphi(p,q)\gamma$, so that $\varphi(p,q) \in (G_2)^{\gamma}$. Evidently $\varphi$ is a homomorphism. We shall show that $\varphi$ is onto. Let $\alpha \in (G_2)^{\gamma}$. Since $\alpha$ satisfies $\gamma\alpha = \alpha\gamma$, $\gC_{\gamma} = \{ x \in \gC \, | \, \gamma x = x \} = \H$ is invariant under $\alpha$, so that the restriction of $\alpha$ to $\gC_{\gamma}$ induces an automorphism of $\H$. Hence there exists $q \in Sp(1)$ satisfying
$$
            \alpha m = qm\overline{q}, \quad m \in \H  $$
(Proposition 108 of Yokota [58]). By putting $\beta = \varphi(1,q)^{-1}\alpha$, we have $\beta \in (G_2)^{\gamma}$ and $\beta|\H = 1$. Therfore $\beta$ induces a transformation of $\gC_{-\gamma} = \{ x \in \gC \,| \, \gamma x = -x \} = \H e_4$. By putting $\beta e_4 = pe_4$, where $p \in \H$, we have $|p| = |pe_4| = |\beta e_4| = |e_4| = 1$, which implies $p \in Sp(1)$. From
$$
     \beta(m + ae_4) = \beta m + (\beta a)(\beta e_4) = m + a(pe_4) = m + (pa)e_4 = \varphi(p,1)(m + ae_4),$$
we have $\beta = \varphi(p,1)$. Therefore we have
$$
      \alpha =\varphi(1,q)\beta = \varphi(1,q)\varphi(p,1) = \varphi(p,q), \quad (p, q) \in Sp(1) \times Sp(1),$$
which shows that $\varphi$ is onto. $\Ker\,\varphi = \{(1,1),(-1,-1) \} = \Z_2$ is easily obtained. Thus we have the isomorphism $(Sp(1)\times Sp(1))/\Z_2 \cong (G_2)^{\gamma}.$ 
\vspace{2mm}

{\bf  Remark 1.} \qquad \qquad $(Sp(1) \times Sp(1))/\Z_2 \cong SO(4)$.
\vspace{1mm} 

\noindent Indeed, a mapping $f : Sp(1) \times Sp(1) \to SO(4) = SO(\H)$ defined by
$$
        f(p, q)x = px\ov{q}, \quad x \in \H $$
induces the isomorphism $(Sp(1) \times Sp(1))/\Z_2 \cong SO(4)$ as groups.
\vspace{1mm}

{\bf Remark 2.} Since $(G_2)^\gamma$ is connected as the fixed points subgroup by an automorphism $\gamma$ of the simply connected group $G_2$, the fact that $\varphi : Sp(1) \times Sp(1) \to (G_2)^\gamma$ is onto can be proved as follows. The elements 
$$
\begin{array}{ll}
      2G_{12} - G_{47} - G_{56}, & - G_{47} + G_{56}, 
\vspace{1mm}\\
      2G_{13} - G_{46} - G_{57}, & \;\;\; G_{46} + G_{57}, 
\vspace{1mm}\\
      2G_{23} - G_{45} - G_{67}, & - G_{45} + G_{67}
\end{array} $$
forms an $\R$-basis of $(\gg_2)^\gamma$. So $\dim(\gg_2)^\gamma = 6 = 3 + 3 = \dim(\sp(1) \oplus \sp(1))$. Hence $\varphi$ is onto. 
\vspace{4mm}

{\bf 1.11. Center $z(G_2)$ of $G_2$}
\vspace{3mm}

{\bf Theorem 1.11.1.} {\it The center $z(G_2)$ of the group $G_2$ is trivial}\,: 
$$
               z(G_2)=\{1\}. $$

{\bf Proof.} Let $\alpha \in z(G_2)$. From the commutativity with $\gamma$: 
$\gamma\alpha = \alpha\gamma$, we have $\alpha \in SO(4)$ (Remark 1 of Theorem 1.10.1) and so $\alpha \in z(SO(4))$. Since the center $z(SO(4))$ of $SO(4)$ 
is the group of order 2, we have
$$ 
      \alpha = 1 \quad \mbox{or} \quad \alpha =\gamma.$$
However $\gamma \notin z(G_2)$ (Theorem 1.10.1). Hence $\alpha = 1$.
\vspace{4mm}

{\bf 1.12. Complex exceptional Lie group ${G_2}^C$}
\vspace{4mm}

{\bf Definition.} The group ${G_2}^C$ is defined to be the automorphism group 
of the complex Cayley algebra $\gC^C$: 
$$
      {G_2}^C = \{ \alpha \in \Iso_C(\gC^C) \, | \, \alpha(xy) = (\alpha x)(\alpha y) \}.$$

{\bf Lemma 1.12.1.} {\it For $\alpha \in {G_2}^C$, we have}
$$
           (\alpha x, \alpha y) = (x, y), \quad x, y \in \gC^C. $$

{\bf Proof.} The equality $\alpha x = (\alpha 1)(\alpha x)$ holds for all $x \in \gC^C$, and so we have $\alpha 1 = 1$.  We shall show that
$$
         \ov{\alpha e_k} = - \alpha e_k, \quad k = 1, 2, \cdots, 7. $$
Note that $(\alpha e_k)(\alpha e_k) = \alpha(e_ke_k) = \alpha(-1) = -1$. Let $x = \alpha e_k$ and $N(x) = x\ov{x} \in C$. Then $xx = - 1$ and $N(x)N(x) = 1$, which shows that $N(x) = \pm 1$. If $N(x) = - 1$, then $\ov{x} = - xx\ov{x} = - xN(x) = x$, so $x \in C$. From $xx = -1$, we have $x = \pm i$. Then $\alpha(e_k) = \pm \alpha(i)$, and so $e_k = \pm i$, which is a contradiction. Hence $N(x) = 1$, that is, $x\ov{x} = 1$ and $\ov{x} = -xx\ov{x} = - xN(x) = - x$. Thus we have
$$
             \ov{\alpha x}=\alpha \ov{x}, \quad x \in \gC^C.$$
Now we have
\begin{eqnarray*}
      (\alpha x, \alpha y) \!\!\! &=& \!\!\! 
     \frac{1}{2}((\alpha x)(\ov{\alpha y}) + (\alpha y)(\ov{\alpha x}))=\frac{1}{2}((\alpha x)(\alpha \ov{y}) + (\alpha y)(\alpha \ov{x}))
\vspace{1mm}\\
     \!\!\! &=& \!\!\! \alpha\Big(\frac{1}{2}(\alpha(x\ov{y} + y\ov{x})\Big) = \alpha((x,y)) = (x,y).
\end{eqnarray*}

We define a positive definite Hermitian inner product $\langle x, y \rangle $ in $\gC^C$ by
$$
         \langle x, y \rangle = (\tau x, y), $$
For $\alpha \in \Hom_C(\gC^C)$, we denote the complex conjugate transpose of  $\alpha$ with respect to the inner product $\langle x, y \rangle $ by $\alpha^{*}$: $\langle \alpha^{*}x, y \rangle = \langle x, \alpha y \rangle$. 
\vspace{3mm}

{\bf Lemma 1.12.2.} (1) {\it For} $\alpha \in {G_2}^C$, {\it we have} $\alpha^{*} = \tau\alpha^{-1}\tau \in {G_2}^C$.
\vspace{1mm}

(2) {\it For any} $\alpha \in G_2$, {\it its complexificated mapping $\alpha^C : \gC^C \to \gC^C$ belongs to} ${G_2}^C$: $\alpha^C \in {G_2}^C$. {\it Identifying $\alpha$ with $\alpha^C$, we regard $G_2$ as a subgroup of ${G_2}^C$}: $G_2 \subset {G_2}^C$. {\it Now, for $\alpha \in {G_2}^C$, we have, $\alpha \in G_2$ if and only if $\tau\alpha = \alpha\tau$, that is, }
$$
            G_2 = \{ \alpha \in {G_2}^C \, | \, \tau\alpha = \alpha\tau \}. $$ 

{\bf Proof.} (1) $\langle \alpha^{*}x, y \rangle = \langle x, \alpha y \rangle =(\tau x, \alpha y) = (\alpha^{-1}\tau x, y) = \langle \tau\alpha^{-1}\tau x, y \rangle $ for all $x, y \in \gC^C$. Hence $\alpha^{*} = \tau\alpha^{-1}\tau \in {G_2}^C$. 
\vspace{1mm}

(2) Let $\alpha \in {G_2}^C$ satisfy $\tau\alpha = \alpha\tau$. Then since $\tau\alpha x = \alpha\tau x = \alpha x$, we have $\alpha x \in \gC$ for $x \in \gC$. Hence $\alpha$ induces an $\R$-transformation $\alpha'$ of $\gC$ and $\alpha' \in G_2$, further we have $\alpha = (\alpha')^C$.
\vspace{3mm}

{\bf Theorem 1.12.3.} {\it The polar decomposition of the group} ${G_2}^C$ {\it is given by}
$$
             {G_2}^C \simeq G_2 \times \R^{14}.$$
{\it In particular}, ${G_2}^C$ {\it is a simply connected complex Lie group 
of type $G_2$}.
\vspace{2mm}

{\bf Proof.} Evidently ${G_2}^C$ is an algebraic subgroup of $\Iso_C(\gC^C) = GL(8, C)$. If $\alpha \in {G_2}^C$, then $\alpha^{*} \in {G_2}^C$ (Lemma 1.12.2.(1)). Hence, from Chevalley's lemma (Chevalley [5]), we have
$$
     {G_2}^C \simeq ({G_2}^C \cap U(\gC^C)) \times \R^d  = G_2 \times \R^d, $$
where $U(\gC^C) = \{\alpha \in \Iso_C(\gC^C) \, | \, \langle \alpha X, \alpha Y \rangle = \langle X, Y \rangle\}$ and $d = \dim{G_2}^C - \dim G_2 = 2 \times 14 - 14 = 14.$ Since $G_2$ is simply connected (Theorem 1.9.3), ${G_2}^C$ is also simply connected. The Lie algebra of the group ${G_2}^C$ is ${\gg_2}^C$, so that ${G_2}^C$ is a complex simple Lie group of type $G_2$.
\vspace{4mm}

{\bf 1.13. Non-compact exceptional Lie group $G_{2(2)}$ of type $G_2$}
\vspace{3mm}

In $\gC' = \H \oplus \H{e_4}'$, we define a multiplication by
$$
        (m + a{e_4}')(n + b{e_4}') = (mn + \ov{b}a) + (a\ov{n} + bm){e_4}'. $$
This algebra $\gC'$ is called the split Cayley algebra, and is isomorphic to $(\gC^C)_{\tau\gamma}  = \{ x \in \gC^C \, | \, \tau\gamma x = x \}$ as algebras. Now, the group $G_{2(2)}$ is defined to be the automorphism group of the split Cayley algebras $\gC'$:
$$
   G_{2(2)} = \{ \alpha \in \Iso_{\sR}(\gC') \, | \, \alpha(xy) = (\alpha x)(\alpha y) \}. $$
and which can also be defined by

$$
       G_{2(2)} = ({G_2}^C)^{\tau\gamma} =
      \{ \alpha \in {G_2}^C \, |\, \tau\gamma\alpha = \alpha\gamma\tau \}. $$

{\bf Theorem 1.13.1.} {\it The polar decomposition of the Lie groups $G_{2(2)}$ is given by}
$$
       G_{2(2)} \simeq (Sp(1) \times Sp(1))/\Z_2 \times \R^{8}.$$  

{\bf Proof.} In the split Cayley algebra $\gC'$, the inner product $(x, y)'$ is defined as $(x, y)' = \dfrac{1}{2}(x\ov{y} + y\ov{x})$. If we define an inner product $(x, y)$ of $\gC'$ by $(\gamma x, y)'$, then $(x, y)$ is a positive definite inner product. For $\alpha \in G_{2(2)}$, the transpose ${}^t\alpha$ of $\alpha$ with respect to the inner product $(x, y)$ is $^t\alpha = \gamma\alpha^{-1}\gamma \in G_{2(2)}$. Since $G_{2(2)}$ is an algebraic subgroup of $\Iso_{\sR}(\gC') = GL(8, \R)$, from Chevalley's lemma, we have
\begin{eqnarray*}
  G_{2(2)} \!\!\! &\simeq& \!\!\! (G_{2(2)} \cap O(8)) \times \R^d
\vspace{1mm}\\
        \!\!\! &=& \!\!\! (({G_2}^C)^{\tau\gamma})^{\gamma} \times \R^d 
                = (({G_2}^C)^{\tau})^{\gamma} \times \R^d 
                = (G_2)^{\gamma} \times \R^d
\vspace{1mm}\\
        \!\!\! &=& \!\!\! (Sp(1) \times Sp(1))/\Z_2 \times \R^d \;\, \mbox{(Theorem 1.10.1)},\;\; d = 8.
\end{eqnarray*}
where $O(8) = \{ \alpha \in \Iso_{\sR}(\gC') \, | \, (\alpha x, \alpha y) = (x, y) \}$.
\vspace{3mm} 

{\bf Theorem 1.13.2.} {\it The center $z(G_{2(2)})$ of the group $G_{2(2)}$ is trivial}: 
$$
              z(G_{2(2)}) = 
\vspace{4mm}
\{1\}. $$

{\bf 1.14. Principle of triality in $SO(8)$}
\vspace{3mm}

For $a \in S^{n-1} = \{ a \in \R^n \, | \, (a, a) = 1 \}$, we define an element  $D_a \in O(n) = O(\R^n) = \{ A \in \Iso_{\sR}(\R^n) \, | \, (Ax, Ay) = (x, y) \}$ by
$$
        D_ax = x - 2(x, a)a, \quad x \in \R^n. $$
$D_a$ is called the reflection with respect to the hyperplane orthogonal to $a$. Its determinant is $-1$: $\det(D_a) = -1$.
\vspace{3mm}

{\bf Lemma 1.14.1.} {\it The group $O(n)$ is generated by reflections, that is, 
any $A \in O(n)$ can be expressed by the product of finite number of reflections}: 
$$
        A = D_{a_m} \cdots D_{a_2}D_{a_1}, \quad a_i \in S^{n-1}. $$
{\it In particular, $A \in SO(n)$ can be expressed by the product of even number of reflections.} 
$$
       A = D_{a_{2m}}\cdots D_{a_2}D_{a_1}, \quad a_i \in S^{n-1}.$$

{\bf Proof.} (See Theorem 24 of Yokota [58).
\vspace{2mm}

From now on, the group $SO(8)$ is identified with the group
$$
      SO(\gC) = \{ \alpha \in \Iso_{\sR}(\gC) \,| \, 
                  (\alpha x,\alpha y) = (x,y), \det \alpha = 1 \}.$$

{\bf Theorem 1.14.2.} ({\bf Principle of triality in $SO(8)$}). {\it For any $\alpha_3 \in SO(8)$, there exist $\alpha_1, \alpha_2 \in SO(8)$ such that
$$
      (\alpha_1 x)(\alpha_2 y) = \alpha_3 (xy), \quad x,y \in \gC.$$
Moreover, $\alpha_1,\alpha_2$ are determined uniquely for $\alpha_3$ up to the 
sign, that is, for $\alpha_3$, such $\alpha_1,\alpha_2$ have to be $\alpha_1,\alpha_2$ or $- \alpha_1, - \alpha_2$.}
\vspace{2mm}

{\bf Proof.} Since $\alpha_3 \in SO(8)$ is expressed by the product of even number of reflections (Lemma 1.14.1), it is sufficient to show the existence of $\alpha_1$, $\alpha_2$ only for $\alpha_3 = D_bD_a, a,b \in \gC, |a| = |b| = 1$. Now, since
$$
     D_a x = x - 2(x,a)a = x - (x\overline{a} + a\overline{x})a = -a\overline{x}a, \quad  x \in \gC,$$
we have $\alpha_3x = D_bD_ax = b(\overline{a}x\overline{a})b.$  Define mappings $\alpha_1,\alpha_2 : \gC \to \gC$ respectively by $\alpha_1x = b(\overline{a}x)$, $\alpha_2x=(x\overline{a})b$. We then see that $\alpha_1, \alpha_2 \in SO(8)$ and
$$
       (\alpha_1x)(\alpha_2y)
         = (b(\overline{a}x))((y\overline{a})b) = b(\overline{a}(xy)\overline{a})b
         = \alpha_3(xy), \quad x,y \in \gC.$$
Next, we shall show the uniqueness of $\alpha_1, \alpha_2$ up to the sign. To prove this, it is sufficient to show in the case $\alpha_3 = 1$. Now, let
$$
         (\alpha_1x)(\alpha_2y) = xy, \quad x,y \in \gC. $$
Let $\alpha_11 = p$, then $|p| = 1$ and $p(\alpha_2y) = y$, so $\alpha_2y = \bar{p}y$. Similarly we have $\alpha_1x = x \bar{q}$, where $q = \alpha_21$. Hence $(x\bar{q})(\bar{p}y) = xy$. If we let $x = y = 1$, then $\bar{q} \,\bar{p} = 1$, so $\bar{q} = p$. Therefore we have
$$
              (xp)(\bar{p}y) = xy, \quad x, y \in \gC.$$
Putting $py$ instead of $y$, we have
$$
              (xp)y = x(py), \quad x,y \in \gC.$$
From this, we see that $p$ is a real number. Hence $p=\pm 1$ because $|p|=1$. Thus we have $\alpha_1 = \alpha_2 = 1$ or $\alpha_1 = \alpha_2 = -1$.
\vspace{3mm}

{\bf Lemma 1.14.3.} {\it For $\alpha_1,\alpha_2,\alpha_3 \in O(8)$, the relation}
\begin{eqnarray*}
         (\alpha_1x)(\alpha_2y) \!\! &=& \!\! 
             \overline{\alpha_3(\ov{xy})}, \quad x,y \in \gC
\end{eqnarray*}
{\it implies}
\begin{eqnarray*}
         (\alpha_2x)(\alpha_3y) \!\! &=& \!\! 
             \overline{\alpha_1(\ov{xy})}, \quad x,y \in \gC,\\
         (\alpha_3x)(\alpha_1y) \!\! &=& \!\! 
             \overline{\alpha_2(\ov{xy})}, \quad x,y \in \gC.
\end{eqnarray*}

{\bf Proof.} If $x = 0$ or $y = 0$, then the statement is trivially valid, so we may assume $x, y \not= 0$. Now, multiply $\overline{\alpha_1x}$ from left and $\alpha_3(\ov{xy})$ from right on the 
\vspace{1mm}
relation $(\alpha_1x)(\alpha_2y) = \overline{\alpha_3(\ov{xy})})$. Then we get $|x|^2(\alpha_2y)(\alpha_3(\ov{xy})) = \overline{\alpha_1x}|xy|^2$. Therefore
$$
        (\alpha_2y)(\alpha_3(\ov{xy})) = \overline{\alpha_1x}|y|^2.$$
Putting $\overline{yz}$ instead of $x$, we have $(\alpha_2y)(\alpha_3(\overline{y}yz)) = \overline{\alpha_1(\ov{yz})}|y|^2$, that is,
$$
        (\alpha_2y)(\alpha_3z) = \overline{\alpha_1(\ov{yz})}.$$
The other relation is similarly obtained.
\vspace{3mm}

{\bf Lemma 1.14.4.} {\it If $\alpha_1,\alpha_2,\alpha_3 \in O(8)$ satisfy
$$
        (\alpha_1x)(\alpha_2y) = \alpha_3(xy), \quad x,y \in \gC, $$ 
then $\alpha_1,\alpha_2,\alpha_3 \in SO(8)$.}
\vspace{2mm}

{\bf Proof.} Suppose $\alpha_1 \notin SO(8)$. Since the mapping $\epsilon : \gC \to \gC$ defined by $\epsilon x = \overline{x}$ belongs to $O(8)$,  and since $\det\,\epsilon = -1$, we have $\beta_1 = \epsilon\alpha_1^{-1} \in SO(8)$. Using the principle of triality (Theorem 1.14.2) on the element $\beta_1$ (cf. Lemma 1.14.3), there exist $\beta_2, \beta_3 \in SO(8)$ such that
$$
       (\beta_1(\alpha_1x))(\beta_2(\alpha_2y)) =   
        \beta_3((\alpha_1x)(\alpha_2 y)) = 
        \beta_3(\alpha_3(xy)), \quad x,y \in \gC .$$
Setting $\beta_2 \alpha_2 = \gamma_2$ and $\beta_3\alpha_3 = \gamma_3$, the above relation becomes
$$
       \overline{x}(\gamma_2y) = \gamma_3(xy), \quad x,y \in \gC.$$
Put $x = 1$, then $\gamma_2y = \gamma_3y$, so $\gamma_2 = \gamma_3$. Hence we have
$$
    \displaylines{\hfill
    \overline{x}(\gamma_2y) = \gamma_2(xy), \quad x,y \in \gC.      
    \hfill\mbox{(i)}}$$
Put $\gamma_21 = p$, then $|p| = 1$. Lett $y = 1$ in (i), then $\overline{x}p = \gamma_2x$. Hence (i) becomes
$$
    \displaylines{\hfill
    \overline{x}(\overline{y}p) = (\overline{xy})p, \quad x,y \in \gC.      
    \hfill\mbox{(ii)}}$$
Again let $y = p$, then $\overline{x} = (\overline{xp})p$, and so $\overline{x}\ \overline{p}=\overline{xp}$. Then
$$
          px = xp \quad \mbox{for all}\;\; x \in \gC.$$
Therefore $p \in \R$, hence $p = \pm 1$ since $|p|=1$. Thus (ii) becomes $\overline{x}\ \overline{y}=\overline{xy}$, and so
$$
           xy = yx \quad \mbox{for all}\;\; x,y \in \gC.$$
But this is a contradiction. Therefore $\alpha_1 \in SO(8)$. As for $\alpha_2, \alpha_3$, use Lemma 1.14.3 and the argument above, to show $\alpha_2$, $\alpha_3 \in SO(8)$.
\vspace{2mm}

As a corollary of the principle of triality (Theorem 1.14.2), we have the following proposition.
\vspace{3mm}

{\bf Proposition 1.14.5.} {\it For $a \in \gC$ such that $|a| = 1$, the mapping $\alpha_a : \gC \to \gC$ defined by
$$
                 \alpha_ax = axa^{-1}, \quad x \in \gC $$
belongs to the group $G_2$ if and only if $a^3 = \pm 1$.}
\vspace{2mm}

{\bf Proof.} If we apply Lemma 1.14.3 to the Moufang formula $(\overline{a}x)(y\overline{a}) = \overline{a}(xy)\overline{a}$, then
$$
     (x\overline{a})(aya) = (xy)a, \quad x, y \in \gC. $$
If we replace $x$ by $ax$ and $y$ by $ya$ respectively, then
$$
      (ax\overline{a})(aya^2) = a(xy)a^2, \quad x, y \in \gC. $$
Now, the mapping $\alpha_a$ belongs to $G_2$ if and only if
$$
       (ax\overline{a})(ay\overline{a}) = a(xy)\overline{a}, \quad x, y \in \gC. $$
From the uniqueness of the principle of triality up to sign, we have
$$
            aya^2 = \pm ay\overline{a}. $$
Therefore $a^2 = \pm\overline{a}$, so that $a^3 = \pm 1$. The converse also holds.
\vspace{2mm}

Let $\omega_1 = - \dfrac{1}{2} + \dfrac{\sqrt{3}}{2}e_1 \in \C \subset \gC$. Then ${\omega_1}^3 = 1$, so $\alpha_{\overline{\omega}_1} \in G_2$ by Proposition 
\vspace{0.9mm}
1.14.5. This $\alpha_{\overline{\omega}_1}$ is nothing but $w$ of Section 1.9: $\alpha_{\overline{\omega}_1} = w$, because
$$
   \alpha_{\overline{\omega}_1}(a + \m) = \overline{\omega}_1(a + \m)\omega_1 
     = \overline{\omega}_1a\omega_1 + {\overline{\omega}_1}^2\m = a + \omega_1\m = w(a + \m), $$
for $a + \m \in \C \oplus \C^3 = \gC$.
\vspace{4mm}

{\bf 1.15. Spinor group $Spin(7)$}
\vspace{3mm}

If we use the principle of triality for $\alpha \in SO(7)$, then there exist $\widetilde{\alpha}$, $\alpha' \in SO(8)$ satisfying
$$
        \displaylines{\hfill
        (\alpha x)(\widetilde{\alpha}y) = \alpha'(xy), \quad x, y \in \gC.
        \hfill\mbox{(i)}}  $$
Putting $x = 1$, we have $\widetilde{\alpha}y = \alpha'y$, and so $\widetilde{\alpha} = \alpha'$. Then (i) becomes
$$
        \displaylines{\hfill
    (\alpha x)(\widetilde{\alpha}y) = \widetilde{\alpha}(xy),\quad x, y \in \gC.        \hfill\mbox{(ii)}}$$
Conversely, suppose that $\alpha,\widetilde{\alpha} \in SO(8)$ satisfy (ii). Putting $x = 1$, we have $(\alpha 1)(\widetilde{\alpha}y) = \widetilde{\alpha}y$, and so we have $\alpha 1 = 1$. Hence $\alpha \in SO(7)$.
\vspace{3mm}

{\bf Definitioin.} We define a subgroup $\widetilde{B}_3$ of $SO(8)$ by
$$
     \widetilde{B}_3 = \{\widetilde{\alpha} \in SO(8) \, | \, (\alpha x)(\widetilde{\alpha}y) = \widetilde{\alpha}(xy), x,y \in \gC \, \, \,  \mbox{for some }\; \alpha \in SO(7)\}.$$

$\wti{B}_3$ is a compact group.
\vspace{3mm}

{\bf Theorem 1.15.1.} \qquad \quad $\wti{B}_3/G_2 \simeq S^7, \quad \wti{B}_3 \cap SO(7) = G_2$.
\vspace{1mm}

\noindent {\it In particular, the group $\widetilde{B}_3$ is connected.}
\vspace{2mm}

{\bf Proof.} $S^7 = \{a \in \gC \, | \, |a| = 1 \}$ is a 7 dimensional sphere. Since $\widetilde{B}_3$ is a subgroup of $SO(8)$, the group $\widetilde{B}_3$ acts on $S^7$. We shall show that $\widetilde{B}_3$ acts transitively on $S^7$. To prove this, it is sufficient to show that any $b_0 \in S^7$ can be transformed to $1 \in S^7$ by some $\alpha \in \wti{B}_3$. Now, we choose any element $a_1 \in S^7$ such that $(1,a_1) = 0$. and choose any element $a_2 \in S^7$ such that $(1, a_2) = (a_1, a_2) = 0$. Let $a_3 \in S^7$ be the element determined by $a_3b_0 = a_2(a_1b_0)$. More precisely, if we let $a_3 = (a_2(a_1b_0))\overline{b}_0$, then $a_3 \in S^7$ and satisfies $(1, a_3) = (a_1, a_3) = (a_2, a_3) = 0$. Choose any element $a_4 \in S^7$ such that $(1, a_4) = (a_1, a_4) = (a_2, a_4) = (a_3, a_4) = 0$. Let $a_5, a_6, a_7 \in S^7$ be elements determined by
$$
      a_5b_0 = a_1(a_4b_0), \quad a_6b_0 = a_2(a_4b_0), \quad  a_7b_0 = a_6(a_1b_0).$$
Then, $\{a_0 = 1, a_1, a_2, \cdots, a_7 \}$ forms an orthonormal $\R$-basis of $\gC$. 
To prove this, we need to show $(a_i, a_j) = \delta_{ij}$, $i,j = 0,1,\cdots,7.$  We will only show the following two since the others can be proved in a similar manner.
\begin{eqnarray*}
     (a_2, a_6) \!\!\! &=& \!\!\! (a_2b_0, a_6b_0)=(a_2b_0, a_2(a_4b_0)) = (b_0, a_4b_0) = (1, a_4) = 0,
\vspace{1mm}\\
     (a_3, a_7) \!\!\! &=& \!\!\! (a_3b_0, a_7b_0) = (a_2(a_1b_0), a_6(a_1b_0)) = (a_1(a_2b_0), a_1(a_6b_0))\\
 \!\!\!&=&\!\!\! (a_2b_0, a_6b_0) = (a_2b_0, a_2(a_4b_0)) = (b_0, a_4b_0) = (1, a_4) = 0.
\end{eqnarray*}
Now, since $\{ e_0, e_1,\cdots, e_7 \}$ and $\{ a_0 = 1, a_1, \cdots, a_7 \}$ are both orthonormal $\R$-bases of $\gC$, the $\R$-linear isomorphism $\alpha : \gC \to \gC$  satisfying
$$
            \alpha e_i = a_i, \quad i = 0, 1, \cdots, 7$$
belongs to $O(7)$. Moreover, this $\alpha$ satisfies
$$\displaylines{\hfill
     (\alpha x)((\alpha y)b_0)=(\alpha(xy))b_0, \quad x,y \in \gC.
\hfill\mbox{(i)}}$$
To prove this, it is sufficient to show that
$$
       (\alpha e_i)((\alpha e_j)b_0) = (\alpha(e_ie_j))b_0, \quad i, j = 0, 1, \cdots, 7.$$
Again we need to verify many cases, but here we will only show the following two examples.
\begin{eqnarray*}
     (\alpha e_1)((\alpha e_3)b_0) \!\!\! &=& \!\!\! a_1(a_3b_0) = a_1(a_2(a_1b_0)) = - a_1(a_1(a_2b_0)) = a_2b_0 
\vspace{1mm}\\
     \!\!\! &=& \!\!\! \alpha(e_2)b_0 = \alpha(e_1e_3)b_0, 
\vspace{1mm}\\
     (\alpha e_2)((\alpha e_5)b_0) \!\!\! &=& \!\!\! a_2(a_5b_0) =  a_2(a_1(a_4b_0)) = - a_1(a_2(a_4b_0)) = - a_1(a_6b_0) 
\vspace{1mm}\\
     \!\!\! &=& \!\!\! a_6(a_1b_0) = a_7b_0 = \alpha(e_7)b_0 = \alpha(e_2e_5)b_0.
\end{eqnarray*}
Now, if we put
$$
           b_i = a_ib_0, \quad i = 0, 1, \cdots, 7,$$
then $\{b_0, b_1, \cdots, b_7 \}$ is an orthonormal $\R$-basis of $\gC$. The $\R$-linear isomorphism $\widetilde{\alpha} : \gC \to \gC$ satisfying
$$
         \widetilde{\alpha}e_i=b_i, \quad  i = 0, 1, \cdots, 7 $$
belongs to $O(8)$. Since $\widetilde{\alpha}x = (\alpha x)b_0$, it follows from (i) that
$$ 
\displaylines{\hfill
  (\alpha x)(\widetilde{\alpha}y) = \widetilde{\alpha}(xy), \quad x, y \in \gC.
\hfill\mbox{(ii)}}$$
Since $\alpha \in O(7)$, $\widetilde{\alpha} \in O(8)$ satisfy (ii), 
$\alpha \in SO(7)$ , $\widetilde{\alpha} \in SO(8)$ (Lemma 1.14.4). Hence $\widetilde{\alpha} \in \widetilde{B}_3$ and $\widetilde{\alpha}1 = b_0$, and so $\wti{\alpha}^{-1}b_0 = 1$. This shows the transitivity. The isotropy subgroup of $\widetilde{B}_3$ at $1 \in S^7$ is $G_2$. Indeed, if $\widetilde{\alpha} \in \widetilde{B}_3$ satisfies $\widetilde{\alpha}1 = 1$, then we have $\alpha = \widetilde{\alpha}$, so $\widetilde{\alpha} \in G_2$. Conversely, $\alpha \in G_2$ satisfies $\alpha \in \widetilde{B}_3$ and $\alpha 1 = 1$. Thus we have the homeomorphism $\widetilde{B}_3/G_2 \simeq S^7$.
\vspace{3mm}

{\bf Theorem 1.15.2.} \qquad \qquad \quad $\widetilde{B}_3 \cong Spin(7).$
\vspace{1mm}

\noindent (From now on, we identify these groups).
\vspace{2mm}

{\bf Proof.} Suppose $\alpha \in SO(7)$ and $\widetilde{\alpha} \in \widetilde{B}_3$ satisfy the principle of triality
$$
   (\alpha x)(\widetilde{\alpha} y) = \widetilde{\alpha}(xy), \quad x, y \in \gC.$$
We define a mapping $p : \widetilde{B_3} \to SO(7)$ by $p(\widetilde{\alpha}) = \alpha$. It is not difficult to see that $p$ is a homomorphism. The principle of triality implies that $p$ is onto and $\Ker \, p = \{1, -1 \}$. Next, we shall prove that $p$ is continuous. From Lemma 1.14.3, we have
$$
     \alpha(xy) = (\widetilde{\alpha}x)\overline{(\widetilde{\alpha}\overline{y})}, \quad x,y \in \gC.$$
Consider the matrices of $\alpha$ and $\widetilde{\alpha}$ with respect to the $\R$-basis $\{ e_0, e_1, \cdots, e_7 \}$. Then we can see that each component of matrix $\alpha$ is a polynominal of components of matrix $\widetilde{\alpha}$ (for example, $\alpha e_1 = (\widetilde{\alpha}e_2)(\overline{\widetilde{\alpha}\overline{e}_3})$). Therefore $p$ is continuous. Hence we have the isomorphism
$$
          \widetilde{B}_3/\{1, -1 \} \cong SO(7).$$
Therefore $\widetilde{B}_3$ is isomorphic to $Spin(7)$ as the universal covering group of $SO(7)$.
\vspace{4mm}

{\bf 1.16. Spinor group $Spin(8)$}
\vspace{3mm}
 
{\bf Definition.} We define a subgroup $\widetilde{D}_4$ of $SO(8) \times SO(8) \times SO(8)$ by
$$
   \widetilde{D}_4 = \{(\alpha_1, \alpha_2, \alpha_3) \in SO(8) \times SO(8) \times SO(8) \, | \, (\alpha_1x)(\alpha_2y) = \overline{\alpha_3(\overline{xy})}, x,y \in \gC \}. $$

$\widetilde{D}_4$ is a compact group. 
\vspace{2mm}

Since an element $(\alpha, \widetilde{\alpha}, \kappa\widetilde{\alpha})$ of $\wti{D}_4$ satisfies $(\alpha x)(\widetilde{\alpha}y) = \widetilde{\alpha}(xy)$, $x, y \in \gC$, we see that $\widetilde{D}_4$ contains $Spin(7)$ as a subgroup under the identification

$$
       Spin(7) \ni \widetilde{\alpha} \;\longleftrightarrow \; (\alpha,\widetilde{\alpha},\kappa \widetilde{\alpha}) 
\vspace{3mm}
\in \widetilde{D}_4.$$

{\bf Proposition 1.16.1.} \qquad \qquad $\widetilde{D}_4/Spin(7) \simeq S^7.$
\vspace{1mm}

\noindent {\it In particular, the group $\widetilde{D}_4$ is connected.}
\vspace{2mm}

{\bf Proof.} $S^7 = \{ a \in \gC \, | \, |a| = 1 \}$ is a 7 dimensional sphere.  We define an action of $\widetilde{D}_4$ on $S^7$ by
$$
     (\alpha_1, \alpha_2, \alpha_3)a = \alpha_1a, \quad a \in S^7. $$
This action is transitive. Let $a \in S^7$. Since $SO(8)$ acts transitively on $S^7$, there exists $\alpha_1 \in SO(8)$ such that $\alpha_11 = a$. For $\alpha_1$, choose $\alpha_2$, $\alpha_3 \in SO(8)$ satisfying the principle of triality 
$$
      (\alpha_1x)(\alpha_2y) = \overline{\alpha_3(\overline{xy})}. \quad x, y \in \gC. $$ 
Then $(\alpha_1, \alpha_2, \alpha_3) \in \widetilde{D}_4$ and $(\alpha_1, \alpha_2, \alpha_3)1 = a$, which shows the transitivity. The isotropy subgroup of $\widetilde{D}_4$ at $1 \in S^7$ is $Spin(7)$. Indeed, if $(\alpha_1, \alpha_2, \alpha_3) \in \widetilde{D}_4$ satisfies $(\alpha_1, \alpha_2, \alpha_3)1 = 1$, then $\alpha_11 = 1$. Therefore $\alpha_1 \in SO(7)$, which shows that $(\alpha_1, \alpha_2, \alpha_3) \in Spin(7)$ and vice versa. Thus we have the homeomorphism $\widetilde{D}_4/Spin(7)$ $ \simeq S^7$.
\vspace{3mm}

{\bf Theorem 1.16.2.} \qquad \qquad $\widetilde{D}_4 \cong Spin(8).$
\vspace{1mm}

\noindent (From now on, we identify these groups).
\vspace{2mm}

{\bf Proof.} We define a mapping $p : \widetilde{D}_4 \to SO(8)$ by
$$
           p(\alpha_1, \alpha_2, \alpha_3) = \alpha_1.$$
Evidently, $p$ is a homomorphism. The principle of triality implies that $p$ 
is onto and $\Ker \, p = \{(1,1,1),(1,-1,-1) \}$. Thus we obtain the isomorphism
$$
     \widetilde{D}_4/\{(1,1,1),(1,-1,-1)\} \cong SO(8).$$
Therefore $\widetilde{D}_4$ is isomorphic to $Spin(8)$ as the universal covering group of $SO(8)$.
\vspace{3mm}

{\bf Theorem 1.16.3.} {\it The center $z(Spin(8))$ of the group $Spin(8)$ is 
isomorphic to the group $\Z_2 \times \Z_2$}: 
\begin{eqnarray*}
   z(Spin(8)) \!\!\! &=& \!\!\! \{(1,1,1),(1,-1,-1),(-1,-1,1),(-1,1,-1) \}\\
     \!\!\! &=& \!\!\! \{(1,1,1),(1,-1,-1) \} \times \{(1,1,1),(-1,-1,1) \} \cong \Z_2 \times \Z_2.
\end{eqnarray*}

{\bf Proof.} The proof follows easily from $z(SO(8))=\{1,-1 \}$ and the principle of triality.
\vspace{3mm}

{\bf Theorem 1.16.4.} {\it We define automorphisms $\kappa,\pi,\nu : Spin(8) \to Spin(8)$ respectively by
$$
\begin{array}{l}
    \kappa(\alpha_1, \alpha_2, \alpha_3) = (\kappa\alpha_1, \kappa\alpha_3,
                                            \kappa\alpha_2),
\vspace{1mm}\\
    \pi(\alpha_1, \alpha_2, \alpha_3) = (\kappa\alpha_3, \kappa\alpha_2, \kappa\alpha_1),
\vspace{1mm}\\
    \nu(\alpha_1, \alpha_2, \alpha_3) = (\alpha_2, \alpha_3, \alpha_1),
\end{array} $$
where $\kappa : SO(8) \to SO(8)$ in the right side is defined by $(\kappa\alpha)x = \overline{\alpha\overline{x}}$, $x \in \gC$. Then we have relations
$$
      \kappa^2 = 1, \quad \pi^2 = 1, \quad \nu^3 = 1, \quad \nu = \pi\kappa.$$
The subgroup $\gS_3$ generated by $\kappa,\pi$ in the automorphism group} 
Aut$(Spin(8))$ {\it of $Spin(8)$ is isomorphic to the symmetric group $S_3$. 
Also, we have
\begin{eqnarray*}
   Spin(7) \!\!\! &=& \!\!\! \{ \alpha \in Spin(8) \, | \, \kappa\alpha = 
\alpha \},\\
   G_2 \!\!\! &=&  \!\!\! \{ \alpha \in Spin(8) \,| \, \lambda\alpha = \alpha,
 \lambda \in \gS_3 \}\\
       \!\!\! &=&  \!\!\! \{ \alpha \in Spin(8) \,| \, \pi\alpha = \alpha, 
\nu\alpha = \alpha \} \\
       \!\!\! &=&  \!\!\! \{ \alpha \in Spin(8) \,| \, \nu\alpha = \alpha \}\\
       \!\!\! &=&  \!\!\! \{ \alpha \in Spin(7) \,| \, \pi\alpha = \alpha \}. 
\end{eqnarray*}
Moreover we have the isomorphism
$$
        Spin(8)/\{(1,1,1),(-1,-1,1)\} \cong Spin(8)/\{(1,1,1),(1,-1,-1) \}, $$
that is,}
$$
          SO(8) \cong Ss(8). $$

{\bf Proof.} The group multiplication between $\kappa,\pi,\nu$ is the same as 
that of Theorem 1.3.5. Next, we shall show
$$
      G_2 = \{ \alpha \in Spin(8) \, | \, \nu\alpha = \alpha \}.$$
If $(\alpha_1,\alpha_2,\alpha_3) \in Spin(8)$ satisfies $\nu(\alpha_1, \alpha_2, \alpha_3) = (\alpha_1, \alpha_2, \alpha_3)$, then we have $\alpha_1 = \alpha_2
 = \alpha_3 \; (= \alpha)$, that is,
$$
        \displaylines{\hfill
     (\alpha x)(\alpha y) = \kappa\alpha(xy), \quad x,y \in \gC,
\hfill\mbox{(i)}}$$
Put $a = \alpha 1$, then $|a| = 1$, and put $x = 1$ and $y = 1$ in (i), 
then we have $a(\alpha y) = \kappa\alpha(y)$ and $(\alpha x)a = \kappa\alpha(x)$, respectively. Hence, we get
$$
        a(\alpha x)=(\alpha x)a  \quad \mbox{ for all }\;\; x \in \gC.$$
hence $a \in \R$, and so $a = \pm 1$ from $|a| = 1$. In the case $a = -1$, 
let $x = y = 1$ in (i), then $(-1)(-1) = -1$ which is a contradiction. Hence $a
 = 1$, so that $\kappa\alpha = \alpha$. Therefore we have $\alpha \in G_2$. Finally, the automorphism $\nu^2 : Spin(8) \to Spin(8)$ satisfies
$$
              \nu^2(-1,-1,1)=(1,-1,-1),$$ 
hence $\nu^2$ induces the isomorphism $SO(8) \cong Ss(8)$.

\newpage

\vspace{5mm}

\begin{center}
\large{\bf Exceptional Lie group $F_4$}
\end{center}
\vspace{4mm} 

{\bf 2.1. Exceptional Jordan algebra $\gJ$}
\vspace{3mm}

Let $\gJ = \gJ(3, \gC)$ denote all $3 \times 3$ Hermitian matrices with 
entries in the Cayley algebra $\gC$:
$$
           \gJ  = \{ X \in M(3, \gC) \, | \, X^* = X \}, $$
where $X^* = {}^t\!\ov{X}$. Any element $X \in \gJ$ is of the form
$$
      X = X(\xi, x) = \pmatrix{\xi_1 & x_3 & \ov{x}_2 \cr
                                \ov{x}_3 & \xi_2 & x_1 \cr 
                                x_2 & \ov{x}_1 & \xi_3}, 
          \quad \xi_i \in \R, x_i \in \gC.$$
$\gJ$ is a 27 dimensional $\R$-vector space. In $\gJ$, the multiplication $X 
\circ Y$, called the Jordan multiplication, is defined by
$$
         X \circ Y = \frac{1}{2}(XY + YX).$$
In $\gJ$, we define the trace $\tr(X)$, an inner product $(X,Y)$ and a 
trilinear form $\tr(X,Y,Z)$ respectively by
$$\begin{array}{c}
      \tr(X) = \xi_1 + \xi_2 + \xi_3, \quad X = X(\xi, x), 
\vspace{1mm}\\ 
      (X, Y) = \tr(X \circ Y), \quad \tr(X, Y, Z) = (X, Y \circ Z).
\end{array} $$
Moreover, in $\gJ$, we define a multiplication $X \times Y$, called the 
Freudenthal multiplication, by
$$
       X \times Y = \frac{1}{2}(2X \circ Y-\tr(X)Y - \tr(Y)X + (\tr(X)\tr(Y) - 
(X, Y))E), $$
(where $E$ is the $3 \times 3$ unit matrix) and a trilinear form $(X, Y, Z)$ 
and the determinant $\det X$ respectively by
$$
     (X, Y, Z) = (X, Y \times Z), \quad \det X = \frac{1}{3}(X, X, X).$$
For $X = X(\xi, x)$, $Y = Y(\eta, y)$ and $Z = Z(\zeta, z) \in \gJ$, the 
explicit forms in the terms of their entries are as follows.
\begin{eqnarray*}
      (X, Y) \!\!\! &=& \!\!\! \sum_{i=1}^3(\xi_i\eta_i + 2(x_i,y_i)),\\
      \tr(X, Y, Z) \!\!\! &=& \!\!\!\sum_{i=1}^3(\xi_i\eta_i\zeta_i + R(x_iy_{i+1}z_{i+2} + x_iz_{i+1}y_{i+2}) \\
    &+& \!\!\! \xi_i((y_{i+1}, z_{i+1}) + (y_{i+2}, z_{i+2})) + 
                \eta_i((z_{i+1}, x_{i+1}) + (z_{i+2}, x_{i+2})) \\
    &+& \!\!\! \zeta_i((x_{i+1}, y_{i+1}) + (x_{i+2}, y_{i+2}))), \\
       (X, Y, Z) \!\!\! &=& \!\!\! \sum_{i=1}^3\Big(\frac{1}{2}(\xi_i\eta_{i+1}\zeta_{i+2} + \xi_i\eta_{i+2}\zeta_{i+1}) + R(x_iy_{i+1}z_{i+2} + x_iz_{i+1}y_{i+2}) \\
      &-& \!\!\! (\xi_i(y_i, z_i) + \eta_i(z_i, x_i) + \zeta_i(x_i, y_i))\Big), \\     
  \det X \!\!\! &=& \!\!\! \xi_1\xi_2\xi_3 + 2R(x_1x_2x_3) -\xi_1x_1\overline{x}_1 - \xi_2x_2\overline{x}_2 - \xi_3x_3\overline{x}_3.
\end{eqnarray*}

{\bf Lemma 2.1.1.} {\it The followings hold in} $\gJ$. 
\vspace{2mm}

(1) $\;$ (i) \quad  $X \circ Y = Y \circ X, \quad X \times Y = Y \times X$. \vspace{1mm}

\qquad (ii) \quad $E \circ X = X, \quad E \times X = \dfrac{1}{2}(\tr (X)E - X), \quad E \times E = E$. 
\vspace{1mm}

(2) $\;$ (i) \quad {\it The inner product} $\;(X,Y) \;$ {\it is symmetric and positive definite}.
\vspace{1mm}

    \qquad (ii) \quad $\tr(X,Y,Z) = \tr(Y,Z,X) = \tr(Z,X,Y) = \tr(X,Z,Y) 
= \tr(Y,X,Z)$ 
\vspace{1mm}

     \qquad \qquad $= \tr(Z,Y,X)$.
\vspace{1mm}

    \qquad \qquad {\it The similar statement is also valid for} $\; (X,Y,Z)$.
\vspace{1mm}

    \qquad (iii) \quad $(X, E) = (X, E, E) = \tr(X, E, E) = \tr(X), \quad 
\tr(X, Y, E) = (X, Y)$.
\vspace{1mm}

    \qquad (iv) \quad $\tr(X \times Y) = \dfrac{1}{2}(\tr(X)\tr(Y)-(X, Y))$.
\vspace{1mm}

(3) $\;$ (i) \quad  $(X \times X) \circ X = (\det X)E \;$ (Hamilton-Cayley). 
\vspace{1mm}

    \qquad (ii) \quad $(X \times X) \times (X \times X) = (\det X)X$.
\vspace{2mm}

{\bf Proof.} (1) is evident.
\vspace{1mm}

(2) is clear from the explicit forms of $(X, Y), \tr(X, Y, Z), (X, Y, Z)$ etc.
\vspace{1mm}

(3) Using the following explicit form
$$
  X \times X = \pmatrix{
 \xi_2\xi_3 - x_1\ov{x}_1 & \ov{x_1x_2} - \xi_3x_3 & x_3x_1 - \xi_2\ov{x}_2 
\vspace{1mm}\cr
 x_1x_2 - \xi_3\ov{x}_3 & \xi_3\xi_1 - x_2\ov{x}_2 & \ov{x_2x_3} - \xi_1x_1 
\vspace{1mm}\cr
 \ov{x_3x_1} - \xi_2x_2 & x_2x_3 - \xi_1\ov{x}_1 & \xi_1\xi_2 - x_3\ov{x}_3},
    \quad X = X(\xi,x),$$
each formula is obtained by direct calculations.
\vspace{2mm}

In $\gJ$, we adopt the following notations:
$$
       E_1 = \pmatrix{1 & 0 & 0 \cr 
                      0 & 0 & 0 \cr 
                      0 & 0 & 0 \cr}, \quad
       E_2 = \pmatrix{0 & 0 & 0 \cr 
                      0 & 1 & 0 \cr 
                      0 & 0 & 0 \cr}, \quad
       E_3 = \pmatrix{0 & 0 & 0 \cr 
                      0 & 0 & 0 \cr
                      0 & 0 & 1 \cr}, $$
$$
       F_1(x) = \pmatrix{0 & 0 & 0 \cr 
                         0 & 0 & x \cr 
                         0 & \ov{x} & 0 \cr}, \quad
       F_2(x) = \pmatrix{0 & 0 & \ov{x} \cr 
                         0 & 0 & 0 \cr 
                         x & 0 & 0 \cr}, \quad
       F_3(x) = \pmatrix{0 & x & 0 \cr 
                         \ov{x} & 0 & 0 \cr 
                         0 & 0 & 0 \cr}. $$
The tables of the Jordan multiplication and the Freudenthal multiplication 
among elements above are given as follows.
$$
\left\{\begin{array}{ll}
         E_i \circ E_i = E_i
\vspace{2mm}\\
         E_i \circ F_i(x) = 0
\vspace{2mm}\\
         F_i(x) \circ F_i(y) = (x,y)(E_{i+1} + E_{i+2}),
\end{array}\right. 
\left\{\begin{array}{ll}
       E_i \circ E_j = 0, \quad i \neq j 
\vspace{1mm}\\
       E_i \circ F_j(x) = \dfrac{1}{2}F_j(x), \quad i \neq j 
\vspace{1mm}\\
       F_i(x) \circ F_{i+1}(y) = \dfrac{1}{2}F_{i+2}(\ov{xy}),
\end{array}\right.$$
$$
\left\{\begin{array}{ll}
         E_i \times E_i = 0
\vspace{1mm}\\
         E_i \times F_i(x) = - \dfrac{1}{2}F_i(x)
\vspace{1mm}\\
         F_i(x) \times F_i(y) = -(x,y)E_i,
\end{array}\right. 
\qquad
\left\{\begin{array}{ll}
       E_i \times E_{i+1} = \dfrac{1}{2}E_{i+2} 
\vspace{1mm}\\
       E_i \times F_j(x) = 0, \quad i \neq j 
\vspace{1mm}\\
       F_i(x) \times F_{i+1}(y) = \dfrac{1}{2}F_{i+2}(\ov{xy}),
\end{array}\right. $$
where the indexes are considered as mod 3.
\vspace{4mm}

{\bf 2.2. Compact exceptional Lie group $F_4$}
\vspace{3mm}

{\bf Definition.} The group $F_4$ is defined to be the automorphism group of 
the Jordan algebra $\gJ$:
$$
       F_4 = \{ \alpha \in \Iso_{\sR}(\gJ) 
           \, | \  \alpha(X \circ Y) = \alpha X \circ \alpha Y \}.$$

{\bf Lemma 2.2.1.} (1) {\it For $\alpha \in F_4$, we have $\alpha E = E$.} 
\vspace{1mm}

(2) {\it For $\alpha \in F_4$, we have} $\tr(\alpha X) = \tr(X)$, $X \in \gJ$.
\vspace{2mm}

{\bf Proof.} (1) Applying $\alpha$ on $E \circ X = X$, we have $\alpha E \circ 
\alpha X = \alpha X$.  Let $X = \alpha^{-1}E$, then $\alpha E \circ E = E$, 
that is, $\alpha E = E$.
\vspace{1mm}

(2)  We use the Hamilton-Cayley identity $X \circ (X \times X) = (\det X)E$ 
(Lemma 2.1.1.(3)), that is,
$$
      \displaylines{\hfill
       X \circ (X \circ X) - \tr(X)X^2 + \frac{1}{2}(\tr(X)^2 - \tr(X^2))X = 
(\det X)E.
      \hfill\mbox{(i)}} $$
We put $\alpha X$ in the place of $X$ of (i) and then apply $\alpha^{-1} \in 
F_4$ on the obtained expression. Then
$$
        \displaylines{\hfill
       X \circ (X \circ X) - \tr(\alpha X)X^2 + \frac{1}{2}(\tr(\alpha X)^2 - 
\tr((\alpha X^2))X = (\det\alpha X)E.
        \hfill\mbox{(ii)}}$$
By subtracting (i)$-$(ii), we get
$$
\begin{array}{l}
        (\tr(\alpha X) - \tr(X))X^2 + \displaystyle{\frac{1}{2}}(\tr(X)^2 - 
\tr(\alpha X)^2 + \tr((\alpha X)^2)-\tr (X^2))X 
\vspace{1mm}\\
        \qquad \qquad = (\det X - \det(\alpha X))E.
\end{array} $$
Let $X = F_i(e_j)$, $i = 1, 2, 3, j = 0, 1,\cdots,7$, then
$$
\begin{array}{l}
     \tr(\alpha F_i(e_j))(E_{i+1} + E_{i+2}) + \displaystyle{\frac{1}{2}}(-\tr(\alpha F_i(e_j))^2 + \tr((\alpha F_i(e_j))^2) - 2)F_i(e_j)
\vspace{1mm}\\
      \qquad \qquad = -\det(\alpha F_i(e_j))E.
\end{array} $$
Comparing the entries of both sides of the equation above, we have
$$
        \tr(\alpha F_i(e_j)) = 0 \, \, (= \tr(F_i(e_j))) $$
and $\tr((\alpha F_i(e_j))^2) = 2$, hence
$$
    \tr(\alpha E_i) = \tr(\alpha(E - F_i(1)^2)) = \tr(E) - \tr((\alpha F_i(1))^2) = 3 - 2 = 1 = \tr(E_i),$$
for $i = 1, 2, 3$. Consequently, we have $\tr(\alpha X) = \tr(X)$ for every $X = E_i, F_i(e_j)$ of the $\R$-basis of $\gJ$. Thus the lemma is proved.
\vspace{2mm}

For $\alpha \in \Hom_{\sR}(\gJ)$, we denote by $^t\alpha$ the transpose of 
$\alpha$ with respect to the inner product $(X, Y)$: $(^t\alpha X, Y) = (X, 
\alpha Y)$.
\vspace{3mm}

{\bf Lemma 2.2.2.} {\it For} $\alpha \in \Iso_{\sR}(\gJ)$, {\it the following 
four conditions are equivalent.}
\vspace{1mm}

(1) $\;\; \det\,(\alpha X) = \det \,X,$ \qquad \mbox{\it for all} $X \in \gJ$.
\vspace{1mm}

(2) $\;\; (\alpha X, \alpha Y, \alpha Z) = (X, Y, Z),$ \quad \mbox{\it for 
all} $X, Y, Z \in \gJ$.
\vspace{1mm}

(3) $\;\; \alpha X \times \alpha Y = {}^t\alpha^{-1}(X \times Y)$, \quad 
\mbox{\it for all} $X, Y \in \gJ$.
\vspace{1mm}

(4) $\;\; \alpha X \times \alpha X = {}^t\alpha^{-1}(X \times X)$, \quad 
\mbox{\it for all} $X \in \gJ$.
\vspace{2mm}

{\bf Proof.} (1) $\Rightarrow$ (2) $\det\,(\alpha X) = \det \,X$ implies that 
$(\alpha X, \alpha X, \alpha X) = (X, X, X).$ Putting $\lambda X + \mu Y + 
\nu Z$ in place of $X$ and comparing the coefficient of $\lambda\mu\nu$, we obtain (2). 
\vspace{1mm}

(2) $\Rightarrow$ (1) is evident.
\vspace{1mm}

(4) $\Rightarrow$ (3) Putting $\lambda X + \mu Y$ in place of $X$ and comparing the coefficient of $\lambda\mu$, we obtain (3). 
\vspace{1mm}

(3) $\Rightarrow$ (4) is evident.
\vspace{1mm}

(2) $\Leftrightarrow$ (3) is easily obtained.
\vspace{3mm}
 
{\bf Lemma 2.2.3.} {\it If $\alpha \in \Iso_{\sR}(\gJ)$ satisfies $\det\,(\alpha 
X) = \det\,X$ for all $X \in \gJ$, then we have}
$$
      \det\,(^t\alpha^{-1}X) = \det\,(^t\alpha X) = \det\,X, \quad \mbox{{\it for all}} \, \, \, X \in \gJ.
$$

{\bf Proof.} We have
$$
\begin{array}{l}{}^t\alpha^{-1}(Y \times Y) \times {}^t\alpha^{-1}(Y \times Y) 
 = (\alpha Y \times \alpha Y)\times (\alpha Y \times \alpha Y) \;\; 
\mbox{(Lemma 2.2.2.(4))} 
\vspace{1mm}\\
\hspace{15mm} 
     = (\det\,\alpha Y)\alpha Y \;\; \mbox{(Lemma 2.1.1.(3))} = (\det\,Y)\alpha Y 
= \alpha((\det\,Y)Y) 
\vspace{1mm}\\

\hspace{15mm}
     = \alpha((Y \times Y) \times (Y \times Y)),\;\; Y \in \gJ.
\end{array} $$
Let $Y = X \times X$, $X \in \gJ$ in the above, then 
$$
     ^t\alpha^{-1}((\det\,X)X) \times {}^t\alpha^{-1}((\det\,X)X)
     = \alpha((\det\,X)X \times (\det\,X)X).$$
We now consider the following two cases.
\vspace{1mm}

(1) Case $\det\,X\neq 0$. In this case, we have ${}^t\alpha^{-1}X \times {}^t 
\alpha^{-1}X = \alpha(X \times X)$. Hence,
$$
\begin{array}{l} 
         3\det\,({}^t\alpha^{-1}X) =  
           ({}^t \alpha^{-1}X,{}^t \alpha^{-1}X \times {}^t \alpha^{-1} X) 
\vspace{1mm}\\
\hspace{15mm} = (^t\alpha^{-1}X, \alpha(X \times X)) 
              = (X, X \times X) = 3\det\,X,
\end{array} $$
hence we have $\det\,(^t\alpha^{-1}X) = \det\,X$. Next, if we use $\alpha^{-1}$ 
instead of $\alpha$, we can see also that $\det\,(^t\alpha X) = \det\,X$. 
\vspace{1mm}

(2) Case $\det\,X = 0$. If $\det\,(^t\alpha^{-1}X) \neq 0$, we can use the result 
of (1). If we put $^t\alpha X$ instead of $X$ in $\det\,(^t\alpha^{-1}X) = 
\det X$ of (1), then $\det\,^t\alpha^{-1}(^t\alpha X)) = \det\,(^t\alpha X)$. 
Hence $0 = \det\,X = \det\,^t\alpha X) \neq 0$ which is a contradiction. Thus we
 have $\det\,(^t\alpha^{-1}X) = \det\,(^t\alpha X) = 0$, so $\det\,(^t\alpha^{-1}X) = \det\,(^t\alpha X) = \det\,X$ is also valid.
\vspace{3mm}

{\bf Lemma 2.2.4.} {\it For} $\alpha \in \Iso_{\sR}(\gJ)$, {\it the following 
five conditions are equivalent.}
\vspace{1mm}

(1) \quad $\alpha(X \circ Y) = \alpha X \circ \alpha Y$.
\vspace{1mm}

(2) \quad $\tr(\alpha X,\alpha Y,\alpha Z) = \tr(X,Y,Z), \,(\alpha X,\alpha Y) = 
(X,Y)$.
\vspace{1mm}

(3) \quad $\det\,(\alpha X) = \det\,X, \,(\alpha X,\alpha Y) = (X,Y)$.
\vspace{1mm}

(4) \quad $\det\,(\alpha X) = \det\,X, \,\alpha E = E$.
\vspace{1mm}

(5) \quad $\alpha(X \times Y) = \alpha X \times \alpha Y$.
\vspace{2mm}

{\bf Proof.} $(1) \Rightarrow (2)$  $(\alpha X,\alpha Y) = \tr(\alpha X \circ 
\alpha Y) = \tr(\alpha(X \circ Y)) = \tr(X \circ Y)$ (Lemma 2.2.1) $ = (X,Y)$. 
Also  $\tr(\alpha X, \alpha Y, \alpha Z) =( \alpha X, \alpha Y \circ \alpha Z) 
= (\alpha X, \alpha(Y \circ Z)) = (X, Y \circ Z) = \tr(X, Y, Z)$.
\vspace{1mm}
  
$(2) \Rightarrow (1)$  $(\alpha X \circ \alpha Y, \alpha Z) = \tr(\alpha X, 
\alpha Y, \alpha Z) = \tr(X, Y,Z ) = (X \circ Y, Z) = (\alpha(X \circ Y), 
\alpha Z)$ holds for all $\alpha Z$, so we have $\alpha X \circ \alpha Y = 
\alpha(X \circ Y)$.
\vspace{1mm}

$(2) \Rightarrow (3)$  Since we have already shown $(2) \Rightarrow (1)$, we 
can use $\tr(\alpha X) = \tr(X)$ (Lemma 2.2.1.(2)). Now,
$$
\begin{array}{l}
     3\det(\alpha X) = \tr(\alpha X, \alpha X, \alpha X) - \dfrac{3}{2}\tr(\alpha X)(\alpha X, \alpha X) + \dfrac{1}{2}\tr(\alpha X)^3 
\vspace{1mm}\\
    \quad \quad = \tr(X, X, X) - \dfrac{3}{2}\tr(X)(X, X) + \dfrac{1}{2}\tr(X)^3 = 3\det X.
\end{array} $$

$(3) \Rightarrow (5)$  $(\alpha(X \times Y), \alpha Z) = (X \times Y, Z) = (X, 
Y, Z) = (\alpha X, \alpha Y, \alpha Z)$ \,(Lemma 2.2.2) $ \; = (\alpha X \times 
\alpha Y, \alpha Z)$ holds for all $\alpha Z \in \gJ$, so we have $\alpha X 
\times \alpha Y = \alpha(X \times Y)$.
\vspace{1mm}

$(5) \Rightarrow (4)$ \quad $(\det\,(\alpha X))\alpha X = (\alpha X \times \alpha 
X) \times(\alpha X \times \alpha X)\;\; \mbox{(Lemma 2.1.1.(3))}$
\vspace{1mm}

$\qquad \qquad \qquad \quad = \alpha((X \times X) \times (X \times X)) = (\det\,X)\alpha X\;\; \mbox{(Lemma 2.1.1.(3))}$

\noindent and so we have $\det\,(\alpha X) = \det\,X$. In the relation
$$
     \alpha X \circ \alpha E = \alpha(X \times E) = \dfrac{1}{2}\alpha(\tr(X)E 
-X),$$
if we denote $\alpha E = P$, then
$$
      \alpha X \times P = \dfrac{1}{2}\tr(X)P - \dfrac{1}{2}\alpha X, \quad X 
\in \gJ. $$
Let $X = \alpha^{-1}E_1$ and $P = \rho_1E_1 + \rho_2E_2 + \rho_3E_3 + F_1(p_1) 
+ F_2(p_2) + F_3(p_3)$, then
$$
\begin{array}{l}
     \dfrac{1}{2}(\rho_2E_3 + \rho_3E_2 - F_1(p_1)) \\
     \qquad = \dfrac{1}{2}(\lambda(\rho_1E_1 + \rho_2E_2 + \rho_3E_3 + F_1(p_1) + F_2(p_2) + F_3(p_3)) - E_1),  
\end{array} $$
where $\lambda = \tr(\alpha^{-1}E_1)$. By comparing entries of both sides, we 
have

$$
   0 = \lambda\rho_1 - 1, \quad \rho_3 = \lambda\rho_2, \quad \rho_2 = \lambda\rho_3, \quad - p_1 = \lambda p_1, \quad 0 = \lambda p_2, \quad 0 = \lambda p_3.$$ Consequently we have $p_2 = p_3 = 0$. Similarly, by letting $X = \alpha^{-1}E_2$, we also have $p_1 = 0$. Again put $X = \alpha^{-1}F_1(1)$, then
$$
    - \dfrac{1}{2}\rho_1F_1(1) = \dfrac{1}{2}(\mu(\rho_1E_1 + \rho_2E_2 + 
\rho_3E_3) - F_1(1)), $$
where $\mu = \tr(\alpha^{-1}F_1(1))$. By Comparing entries of $F_1$-parts, we 
see that $\rho_1 = 1$. Similarly $\rho_2 = \rho_3 = 1$. Therefore we have 
$\alpha E = E$.
\vspace{1mm}

$(4) \Rightarrow (2)$  $\tr(\alpha X) = (\alpha X, E, E) = (\alpha X, \alpha E,
 \alpha E) = (X, E, E)$ 
\vspace{0.5mm}
(Lemma 2.2.2.(2)) $ $ = $\tr(X)$. Hence 
$\dfrac{1}{2}(\tr(X)\tr(Y) - (X, Y)) = (X, Y, E) \; \mbox{(Lemma 2.1.1)} 
\vspace{0.7mm}
= (\alpha
 X, \alpha Y, \alpha E)$ $ = (\alpha X, \alpha Y, E) = \dfrac{1}{2}(\tr(\alpha 
X)$ $\tr(\alpha Y) - (\alpha X, \alpha Y)) = \dfrac{1}{2}(\tr(X)\tr(Y) - 
(\alpha X, 
\vspace{0.5mm}
\alpha Y))$. Therefore we obtain 
$$
               (\alpha X, \alpha Y) = (X, Y). $$
Next, using the relation $(X, Y, Z) = \tr(X, Y, Z) - \dfrac{1}{2}\tr(X)(Y, Z) - \displaystyle{\frac{1}{2}}\tr(Y)(Z,X) - \dfrac{1}{2}\tr(Z)(X, Y) + \dfrac{1}{2}\tr(X)\tr(Y)\tr(Z)$ 
\vspace{0.5mm}
and $(\alpha X, \alpha Y, \alpha Z) = (X, Y, Z)$, we obtain 
$$
         \tr(\alpha X, \alpha Y, \alpha Z) = \tr(X, Y, Z). $$

{\bf Theorem 2.2.5.} $F_4$ {\it is a compact Lie group.}
\vspace{2mm}

{\bf Proof.}  $F_4$ is a compact Lie group as a closed subgroup of the 
orthogonal group 
$$
      O(27) = O(\gJ) = \{ \alpha \in \Iso_{\sR}(\gJ) \, | \, (\alpha X, \alpha 
Y) = (X, Y) \}\;\;\mbox{(Lemma 2.2.4.(2))}. $$

The group $F_4$ contains $G_2$ as a subgroup in the following way. For $\alpha 
\in G_2$, we consider a mapping $\wti{\alpha} : \gJ \to \gJ$,
$$
     \wti{\alpha}\pmatrix{\xi_1 & x_3 & \ov{x}_2 \cr
                                \ov{x}_3 & \xi_2 & x_1 \cr 
                                x_2 & \ov{x}_1 & \xi_3} = 
                 \pmatrix{\xi_1 & \alpha x_3 & \ov{\alpha x_2} \cr
                                \ov{\alpha x_3} & \xi_2 & \alpha x_1 \cr 
                                \alpha x_2 & \ov{\alpha x_1} & \xi_3}. $$
Then $\wti{\alpha} \in F_4$. So we identify $\alpha \in G_2$ with $\wti{\alpha}
 \in F_4$\,: $G_2 \subset F_4$.
\vspace{4mm}

{\bf 2.3. Lie algebra $\gf_4$ of $F_4$}
\vspace{3mm}

In order to investigate the Lie algebra $\gf_4$ of the group $F_4$, it will be 
helpful to study the Lie algebra $\gge_{6({}-26)}$ of the group
\begin{eqnarray*}
     E_{6({}-26)} \!\!\! &=& \!\!\! \{\alpha \in \Iso_{\sR}(\gJ) \, | \, \det\,(\alpha X) = \det\,X \} \\
                  \!\!\! &=&\!\!\! \{\alpha \in \Iso_{\sR}(\gJ) \, | \, {}^t\alpha^{-1}(X \times Y) = \alpha X \times \alpha Y \}.
\end{eqnarray*}

{\bf Lemma 2.3.1.} {\it The Lie algebra $\gge_{6(-26)}$ of the group $E_{6(-26)}$ is given by}
\begin{eqnarray*}
      \gge_{6(-26)} \!\!\! &=& \!\!\! \{ \phi \in \Hom_{\sR}(\gJ) \, | \, 
(\phi X, X, X) = 0 \} \\
                    \!\!\! &=& \!\!\! \{ \phi \in \Hom_{\sR}(\gJ) \, | \, 
(\phi X, Y, Z) + (X, \phi Y, Z) + (X, Y, \phi Z) = 0 \} \\
                    \!\!\! &=& \!\!\! \{ \phi \in \Hom_{\sR}(\gJ) \, | \, -{}^t\phi(X \times Y) = \phi X \times Y + X \times \phi Y \}.
\end{eqnarray*}

{\bf Proof.} It is easy to verify that these conditions in $\gge_{6(-26)}$ are 
equivalent (see Lemma 2.2.2). Now, if $\phi \in \Hom_{\sR}(\gJ)$ satisfies 
$((\exp t\phi)X, (\exp t\phi)X, (\exp t\phi)X) = (X, X, X)$ for all $t \in \R$,
 then we have $(\phi X, X, X) = 0$ by putting $t = 0$ after differentiating 
with respect to $t$. Conversely, if $\phi \in \Hom_{\sR}(\gJ)$ satisfies 
$-{}^t\phi(X \times Y) = \phi X \times Y + X \times \phi Y$, then it is easy to
 verify that $\alpha = \exp t\phi$ satisfies ${}^t\alpha^{-1}(X \times Y) = 
\alpha X \times 
\vspace{3mm}
\alpha Y$.

{\bf Theorem 2.3.2.} {\it The Lie algebra $\gf_4$ of the group $F_4$ is given 
by}
\begin{eqnarray*}
      \gf_4 \!\!\! &=& \!\!\! \{ \delta \in \Hom_{\sR}(\gJ) \, | \, \delta(X 
\circ Y) = \delta X \circ Y + X \circ \delta Y \} 
\vspace{1mm}\\
           \!\!\! &=& \!\!\! \Big\{ {\delta \in \Hom_{\sR}(\gJ) \, \left| 
\begin{array}{l}
     \tr(\delta X, Y, Z) + \tr(X, \delta Y, Z) + \tr(X, Y, \delta Z) = 0 \\
     (\delta X, Y) + (X, \delta Y) = 0
\end{array}\right. } \Big\}
\vspace{1mm}\\
           \!\!\! &=& \!\!\! \{ \delta \in \Hom_{\sR}(\gJ) \, | \, (\delta X, 
X, X) = 0, \; (\delta X, Y) + (X, \delta Y) = 0 \} 
\vspace{1mm}\\
           \!\!\! &=& \!\!\! \{ \delta \in \Hom_{\sR}(\gJ) \, | \, (\delta X, 
X, X)=0, \; \delta E = 0 \} 
\vspace{1mm}\\
           \!\!\! &=& \!\!\! \{ \delta \in \Hom_{\sR}(\gJ) \, | \, \delta(X 
\times Y) = \delta X \times Y + X \times \delta Y \}.
\end{eqnarray*}

{\bf Proof.} The theorem follows easily from Lemma 2.2.4.
\vspace{2mm}

We define an $\R$-vector space $\gM^-$ by
$$  
      \gM^- = \{ A \in M(3, \gC) \, | \, A^* = -A \}.$$
For $X, Y \in M(3,\gC)$, we define $[X, Y] \in M(3,\gC)$ by
$$
                  [X, Y] = XY - YX. $$

{\bf Lemma 2.3.3.} \qquad \quad $[\gM^-, \gJ] \subset \gJ, \quad [\gJ, \gJ] 
\subset \gM^-.$
\vspace{2mm}

Since $[\gM^-, \gJ] \subset \gJ$, any element $A \in \gM^-$ induces an 
$\R$-linear mapping $\widetilde{A} : \gJ \to \gJ$ defined by
$$
            \wti{A}X = \frac{1}{2}[A, X], \quad X \in \gJ. $$

{\bf Lemma 2.3.4.} {\it For $X \in \gJ$, there exists $a \in \gC_0$ such that}
$$
                 [X, XX] = aE. $$

{\bf Proof.}  Let $X = \Big(x_{ij} \Big)$, $x_{ij} \in \gC$, $\ov{x}_{ij} = 
x_{ji}$. The $(i,j)$-entry $a_{ij}$ of $[X, XX] = X(XX) - (XX)X$ is given by
$$
    a_{ij} = \sum_{k,l}(x_{ik}(x_{kl}x_{lj}) - (x_{ik}x_{kl})x_{lj}) 
           = - \sum_{k,l}\{x_{ik}, x_{kl}, x_{lj}\}, \quad i, j = 1, 2, 3. $$
Since $x_{ii}$ is real, if the bracket $\{ \;\;, \; \;, \;\; \}$ contains 
$x_{ii}$, then $\{\;\;, \; \;, \;\; \}$ is $0$. If $i \neq j$, then $a_{ij}$ is
a sum of $a(xa) - (ax)a$, $a(\ov{a}x) - (a\ov{a})x$ and so on, so $a_{ij} = 
0$. If $i = j$, then
\begin{eqnarray*}
     a_{11} \!\!\! &=& \!\!\! - \{x_{12}, x_{23}, x_{31}\} - \{x_{13}, x_{32}, 
x_{21}\}, 
\vspace{1mm}\\
     a_{22} \!\!\! &=& \!\!\! - \{x_{21}, x_{13}, x_{32}\} - \{x_{23}, x_{31}, 
x_{12}\}, 
\vspace{1mm}\\
     a_{33} \!\!\! &=& \!\!\! - \{x_{31}, x_{12}, x_{23}\} - \{x_{32}, x_{21}, 
x_{13}\}, 
\end{eqnarray*}
however they are equal, that is, $a_{11} = a_{22} = a_{33} \, (\,  = a)$. 
Hence we have $[X, XX] = aE$. Since $X, XX \in \gJ$, we have $[X, XX] \in \gM^-$ (Lemma 2.3.3). Therefore, $(aE)^* = - aE$ and so $\ov{a}E = - aE$, which imply 
that $\ov{a} = - a$.
\vspace{2mm}

To prove the following Proposition 2.3.6, in $M(3,\gC)$, we define a real 
valued symmetric inner product $(X,Y)$ by
$$ 
        (X, Y) = \frac{1}{2}\tr(XY + Y^*X^*). $$

{\bf Lemma 2.3.5.} {\it The inner product $(X,Y)$ of $M(3,\gC)$ satisfies} 
$$
      (XY, Z) = (YZ, X) = (ZX, Y) = (X, YZ) = (Y, ZX) = (Z, XY). $$

{\bf Proof.} Let $ X = \Big(x_{ij} \Big)$, $Y = \Big(y_{ij} \Big)$, $Z = 
\Big(z_{ij} \Big)$. Then we have
\begin{eqnarray*}
   (XY, Z) \!\!\! &=& \!\!\! R(XY, Z) = \dfrac{1}{2}R(\tr((XY)Z + Z^*(Y^*X^*)))\\           \!\!\! &=& \!\!\! \frac{1}{2}R\Big(\sum_{i,j,k}((x_{ij}y_{jk})z_{ki} + \ov{z}_{ji}(\ov{y}_{kj}\ov{x}_{ik})) \Big)
\vspace{1mm}\\
           \!\!\! &=& \!\!\! \dfrac{1}{2}R\Big(\dsum_{i,j,k}((y_{jk}z_{ki})x_{ij} + \ov{x}_{ik}(\ov{z}_{ji}\ov{y}_{kj})) \Big)
\vspace{1mm}\\
           \!\!\! &=& \!\!\! \dfrac{1}{2}R(\tr((YZ)X + X^*(Z^*Y^*)))=(YZ, X).
\end{eqnarray*}

{\bf Proposition 2.3.6.} {\it For} $A \in \gM^-$, $\tr(A) = 0$, {\it we have 
$\widetilde{A} \in \gf_4$.}
\vspace{2mm}

{\bf Proof.}  From the equivalent conditions in Theorem 2.3.2, it is sufficient
 to show the following two formulas:
$$
\left\{
\begin{array}{l}
     ([A, X], Y) + (X, [A, Y]) = 0, \quad X, Y \in \gJ,
\vspace{1mm}\\
      \tr([A, X], Y, Z) + \tr(X, [A, Y], Z) + \tr(X, Y, [A, Z]) = 0, \quad X, 
Y, Z \in \gJ.
\end{array}\right.$$
Now, the left side of the first formula $= (AX, Y) - (XA, Y) + (X, AY) - (X, 
YA) = 0$ (Lemma 2.3.5). Next, we show that
$$
       ([A, X], XX) = 0, \quad X \in \gJ. $$
Certainly, if $[X, XX] = aE$, $a \in \gC_0$ (Lemma 2.3.4), then
\vspace{2mm}

\qquad \qquad
   $([A, X], XX) = (AX, XX) - (XA, XX)$ 
\vspace{1mm}

\qquad \qquad \qquad \qquad \qquad 
    $= (A, X(XX)) - (A, (XX)X) \;\;\mbox{(Lemma 2.3.5)}$
\vspace{1mm}

\qquad \qquad \qquad \qquad \qquad 
    $= (A, [X, XX]) = (A, aE) = \dfrac{1}{2}\tr(Aa + \ov{a}A^*)$
\vspace{1mm}

\qquad \qquad \qquad \qquad \qquad 
   $= \dfrac{1}{2}\tr(Aa - aA) = \dfrac{1}{2}(\tr(A)a - a\tr(A)) = 0$.
\vspace{2mm}

\noindent Now, putting $\lambda X + \mu X + \nu Z$ in the place of $X$, and 
comparing the coefficient of $\lambda\mu\nu$, we obtain
$$
      ([A, X],YZ + ZY) + ([A, Y], XZ + ZX) + ([A, Z], XY + YX) = 0, $$
which is the required second formula.
\vspace{2mm}

In $\gM^-$, we adopt the following notation:
$$
           A_1(a) = \pmatrix{0 & 0 & 0 \cr
                             0 & 0 & a \cr 
                             0 & -\ov{a} & 0},\quad
           A_2(a) = \pmatrix{0 & 0 & -\ov{a} \cr
                             0 & 0 & 0 \cr
                             a & 0 & 0}, \quad
           A_3(a) = \pmatrix{ 0 & a & 0 \cr
                              -\ov{a} & 0 & 0 \cr
                              0 & 0 & 0}. $$
Then $\widetilde{A}_i(a) \in \gf_4$ (Proposition 2.3.6) and the operation of 
$\wti{A}_i(a)$ on $\gJ$ is given by
$$
\left\{\begin{array}{l}
         \widetilde{A}_i(a)E_{i} = 0 
\vspace{1mm}\\
         \widetilde{A}_i(a)E_{i+1} = - \dfrac{1}{2}F_i(a)
\vspace{1mm}\\
         \widetilde{A}_i(a)E_{i+2} = \dfrac{1}{2}F_i(a),
\end{array} \right.
\quad
\left\{\begin{array}{l}
    \widetilde{A}_i(a)F_i(x) = (a, x)(E_{i+1} - E_{i+2}) 
\vspace{1mm}\\
    \widetilde{A}_i(a)F_{i+1}(x) = \dfrac{1}{2}F_{i+2}(\ov{ax}) 
\vspace{1mm}\\
    \widetilde{A}_i(a)F_{i+2}(x) = - \dfrac{1}{2}F_{i+1}(\ov{xa}).
\end{array} \right. $$

{\bf Proposition 2.3.7.} {\it The Lie subalgebra $\gd_4$ of $\gf_4$}:
$$
      \gd_4 = \{ \delta \in \gf_4 \, | \, \delta E_i = 0, i = 1, 2, 3 \} $$
{\it is isomorphic to the Lie algebra $\gD_4 = \so(8)$}: 
$$
      \gD_4 = \{ D \in \Hom_{\sR}(\gC) \, | \, (Dx, y) + (x, Dy) = 0\} $$
{\it under the correspondence}
$$
        \gD_4 \ni D_1 \longrightarrow \delta \in \gd_4 $$
{\it given by}
$$
      \delta\pmatrix{\xi_1 & x_3 & \ov{x}_2 \cr
                     \ov{x}_3 & \xi_2 & x_1 \cr
                     x_2 & \ov{x}_1 & \xi_3}
          = \pmatrix{0 & D_3x_3 & \ov{D_2x_2} 
\vspace{0.5mm}\cr
                     \ov{D_3x_3} & 0 & D_1x_1 
\vspace{0.5mm}\cr
                     D_2x_2 & \ov{D_1x_1} & 0}, $$
{\it where $D_2, D_3$ are elements of $\gD_4$ which are determined by $D_1$ 
from the principle of triality}\,:
$$
         (D_1x)y + x(D_2y) = \ov{D_3(\ov{xy})}, \quad x, y \in \gC.$$
(From now on, we identify these Lie algebras $\mbox{\es d}_4$ and $\gD_4$).
\vspace{2mm}

{\bf Proof.} We define a mapping $\varphi_* : \gD_4 \to \gd_4$ by $\varphi_*(D_1) = \delta$. We first prove that $\delta \in \gd_4 \subset \gf_4$. Indeed, 
$$
\begin{array}{l}
   (\delta X, X, X) = (\delta X, X \times X)
\vspace{1.5mm}\\
\quad
   = \Big(\pmatrix{0 & D_3x_3 & \ov{D_2x_2} 
\vspace{1mm}\cr
                     \ov{D_3x_3} & 0 & D_1x_1 
\vspace{1mm}\cr
                     D_2x_2 & \ov{D_1x_1} & 0}, 
       \pmatrix{\xi_2\xi_3 - x_1\ov{x}_1 & \ov{x_1x_2} - \xi_3x_3 & x_3x_1 - 
\xi_2\ov{x}_2 
\vspace{1mm}\cr
 x_1x_2 - \xi_3\ov{x}_3 & \xi_3\xi_1 - x_2\ov{x}_2 & \ov{x_2x_3} - \xi_1x_1 
\vspace{1mm}\cr
 \ov{x_3x_1} - \xi_2x_2 & x_2x_3 - \xi_1\ov{x}_1 & \xi_1\xi_2 - x_3\ov{x}_3} 
\Big)
\vspace{1.5mm}\\
\quad
   = 2(D_1x_1, \ov{x_2x_3} - \xi_1x_1) + 2(D_2x_2, \ov{x_3x_1} - \xi_2x_2) + 
2(D_3x_3, \ov{x_1x_2} - \xi_3x_3)
\vspace{1.5mm}\\
\quad
   = 2((D_1x_1, \ov{x_2x_3}) + (D_2x_2, \ov{x_3x_1}) + (D_3x_3, \ov{x_1x_2}))
\vspace{1.5mm}\\
\quad
   = 2((D_1x_1, \ov{x_2x_3}) + (\ov{D_2x_2}, x_3x_1) + (\ov{D_3x_3}, x_1x_2))
\end{array} $$
which is equal to 0, if we use the relation
$$
\begin{array}{l}
    (D_1x_1, \ov{x_2x_3}) = - (x_1, D_1(\ov{x_2x_3})) = - (x_1, \ov{(D_2x_2)x_3
 + x_2(D_3x_3)})
\vspace{1.5mm}\\
\qquad
   = - (x_1, \ov{x_3}\ov{(D_2x_2)} + \ov{(D_3x_3)}\ov{x_2}) = - (x_3x_1, 
\ov{D_2x_2}) - (x_1x_2, \ov{D_3x_3}).
\end{array} $$
Hence $\delta \in \gf_4$. $\delta E_i = 0, i = 1, 2, 3$ are evident, so that 
$\delta \in \gd_4$. Clearly $\varphi_*$ is injective. We shall show that 
$\varphi_*$ is onto. Let $\delta \in \gd_4$. We put
$$
  \gJ_i = \{F_i(x) \, | \, x \in \gC \} = \{X \in \gJ \, | \, 2E_{i+1} \circ X 
= 2E_{i+2} \circ X = X \}. $$
Since $\delta X \in \gJ_i$ for $X \in \gJ_i$, $\delta$ induces $\R$-linear 
mappings $\delta : \gJ_i \to \gJ_i$ and $D_i : \gC \to \gC$ satisfying
$$
      \delta F_i(x) = F_i(D_ix), \quad x \in \gC, $$ 
for $i = 1, 2, 3$. Applying $\delta$ on $F_i(x) \circ F_i(y) = (x, y)(E_{i+1} + E_{i+2})$, we have $F_i(D_ix) \circ F_i(y) + F_i(x) \circ F_i(D_iy) = 0$, and hence we have
$$
            (D_ix, y) + (x, D_iy) = 0, \quad x, y \in \gC. $$
Hence $D_i \in \gD_4, i = 1, 2, 3$. Moreover, by applying $\delta$ on $F_1(x) \circ F_2(y) = \dfrac{1}{2}F_3(\ov{xy})$, we see that
$$
         (D_1x)y + x(D_2y) = \ov{D_3(\ov{xy})}, \quad x, y \in \gC. $$
This shows that $\varphi_*$ is onto. Thus Proposition 2.3.7 is proved.
\vspace{3mm}

{\bf Theorem 2.3.8.} {\it Any element $\delta \in \gf_4$ is uniquely expressed 
by}
$$
     \delta = D + \widetilde{A}, \quad D \in \gd_4, A \in \gM^-, \diag A = 0,$$
{\it where} $\diag A = 0$ {\it means that all diagonal elements} $a_{ii}$ of 
$A$ are $0$. {\it In particular, the dimension of $\gf_4$ is}
$$
             \dim\gf_4 = 28 + 24 = 52.$$ 

{\bf Proof.} Applying $\delta \in \gf_4$ on $E_i \circ E_i = E_i$ and $E_i \circ E_j = 0$, $i \neq j$, we have 
$$
     2\delta E_i \circ E_i = \delta E_i, \quad
     \delta E_i \circ E_j + E_i \circ \delta E_j = 0. $$ 
From these relations, we see that each $\delta E_i$ is of the form
$$
       \delta E_1 = \pmatrix{0 \! & \! -a_3 \! & \! \ov{a}_2 \cr
                             -\ov{a}_3 \! & \! 0 \! & \! 0 \cr 
                             a_2 \! & \! 0 \! & \! 0},
       \delta E_2 = \pmatrix{0 \! & \! a_3 \! & \! 0 \cr 
                             \ov{a}_3 \! & \! 0 \! & \! -a_1 \cr
                             0 \! & \! -\ov{a}_1 \! & \! 0},
       \delta E_3 = \pmatrix{0 \! & \! 0 \! & \! -\ov{a}_2 \cr
                             0 \! & \! 0 \! & \! a_1 \cr
                             -a_2 \! & \! \ov{a}_1 \! & \! 0},$$ 
$a_i \in \gC$. If we construct a matrix $A 
          = 2\pmatrix{0 \! & \! a_3 \! & \! -\ov{a}_2 \cr
                     -\ov{a}_3 \! & \! 0 \! & \! a_1 \cr
                     a_2 \! & \! -\ov{a}_1 \! & \! 0}$ 
using these elements $a_i$, then $A \in \gM^-$ with $\diag A = 0$, and we have 
$$
        \delta E_i = \wti{A}E_i, \quad i = 1, 2, 3. $$
If we put $D = \delta - \wti{A}$, then $DE_i = 0$, $i = 1, 2, 3$, hence $D \in
 \gd_4$. Thus, $\delta$ can be expressed by $\delta = D + \wti{A}$, where $D 
\in \gd_4$, $A \in \gM^-$, $\diag A = 0$. To show the uniqueness of the 
expression, it is sufficient to prove that 
$$
      D + \wti{A} = 0, \; \; D \in \gd_4, A \in \gM^-, \diag A = 0, \;\; 
\mbox{then} \;\; D = 0, A = 0.$$
Certainly, if we apply it on $E_i$, then $\widetilde{A}E_i = 0$, $i = 1, 2, 3$,
 hence $A = 0$ and so $D = 0$. Finally, we have $\dim\gf_4 = 28 + 24 = 52$ from
 the expression above. Thus the theorem is proved.
\vspace{4mm}

{\bf 2.4. Simplicity of ${\gf_4}^C$}
\vspace{3mm}

Let $\gJ^C = \{ X_1 + i X_2 \, | \,  X_1, X_2 \in \gJ \}$ be the 
complexification  of the Jordan algebra $\gJ$.  In the same manner as in $\gJ$,
 in $\gJ^C$ we can also define the multiplications $X \circ Y$, $X \times Y$, 
the inner product $(X, Y)$, the trilinear forms $\tr(X, Y, Z)$, $(X, Y, Z)$ and the determinant $\det X$. They have the same properties as those of $\gJ$. 
$\gJ^C$ is called the complex exceptional Jordan algebra. $\gJ^C$ has two complex conjugations, namely,
$$
      \ov{X_1 + i X_2} = \ov{X}_1 + i\ov{X}_2,\;\; \tau(X_1 + iX_2) = X_1 - 
iX_2, \quad X_i \in \gJ.$$
The complex conjugation $\tau$ of $\gJ^C$ satisfies
$$
        \tau(X \circ Y) = \tau X \circ \tau Y,\;\; \tau(X \times Y) = \tau X 
\times \tau Y, \quad X, Y \in \gJ^C.$$

We define Lie algebras ${\gge_6}^C$ and ${\gf_4}^C$ respectively by
\begin{eqnarray*}
   {\gge_6}^C \!\!\!&=&\!\!\! \{\phi \in \Hom_C(\gJ^C) \, | \, (\phi X, X, X) = 0\}
\vspace{1mm}\\
    \!\!\!&=&\!\!\! \{\phi \in \Hom_C(\gJ^C) \, | \, -{}^t\phi(X \times Y) = \phi X \times Y + X \times \phi Y\},
\vspace{1mm}\\  
  {\gf_4}^C \!\!\!&=&\!\!\! \{\delta \in \Hom_C(\gJ^C) \, | \, \delta(X \circ Y) = \delta X \circ Y + X \circ \delta Y\}
\vspace{1mm}\\
            \!\!\!&=&\!\!\! \{\delta \in \Hom_C(\gJ^C) \, | \, \delta(X \times Y) = \delta X \times Y + X \times \delta Y\}.
\end{eqnarray*}  
These are the complexification of the Lie algebras $\gge_{6(-26)}$ and $f_4$ respectively. Then the properties of $\gge_{6(-26)}$ and $\gf_4$ stated in Section 2.3 also hold for ${\gge_6}^C$ and ${\gf_4}^C$.  Hereafter we will describe ${\gf_4}^ C$, but the statements are also valid for $\gf_4$ using $\gJ$ instead of $\gJ^C$.\vspace{2mm}

For $A \in \gJ^C$, we define a $C$-linear mapping $\wti{A}: \gJ^C \to \gJ^C$ by
$$
        \wti{A}X = A \circ X, \quad X \in \gJ^C. $$ 

{\bf Proposition 2.4.1.} (1) {\it For} $A \in \gJ^C$, $\tr(A) = 0$, {\it we 
have} $\wti{A} \in {\gge_6}^C$.
\vspace{1mm}

(2) {\it For} $A,B \in \gJ^C$, {\it we have} $[\wti{A}, \wti{B}] \in {\gf_4}^C$.
\vspace{2mm}

{\bf Proof.} (1) $(\wti{A}X, X, X) = (A \circ X, X \times X) = (A, X \circ (X \times X))$ 
\vspace{1mm}

\quad \quad \quad \quad
   $= (A, (\det\,X)E) = (\det\,X)(A, E) = (\det\,X)\tr(A) = 0, \quad X \in \gJ^C. 
$\vspace{1mm}

\noindent Hence $\wti{A} \in {\gge_6}^C$.
\vspace{1mm}

(2)  \quad $[\wti{A}, \wti{B}] = \Big[\Big(A - \displaystyle{\frac{1}{3}}\tr(A)E \Big)^{\sim}, \Big(B - \displaystyle{\frac{1}{3}}\tr(B)E \Big)^{\sim}\Big] \in
 {e_6}^C,$ 
\vspace{1mm}

\qquad \,\, 
      $ [\wti{A}, \wti{B}]E = \wti{A}(\wti{B}E) - \wti{B}(\wti{A}E) = A \circ B
 -  B \circ A = 0. $
\vspace{1mm}

\noindent Hence $[\wti{A}, \wti{B}] \in {\gf_4}^C$ (Theorem 2.3.2). 
\vspace{3mm}

{\bf Lemma 2.4.2.} {\it In} ${\gf_4}^C$, {\it we have}
$$
\begin{array}{ll}
    [\wti{E}_i, \wti{E}_j] = 0, & [\wti{E}_i, \wti{F}_i(a)] = 0,  
\vspace{1mm}\\
    {[}\wti{E}_i, \wti{F}_{i+1}(a){]} = -\dfrac{1}{2}\wti{A}_{i+1}(a), &
    [\wti{E}_i, \wti{F}_{i+2}(a)] = \dfrac{1}{2}\wti{A}_{i+2}(a), 
\vspace{1mm}\\
    {[}\wti{F}_i(a), \wti{F}_{i}(b){]} \in  {\gd_4}^C, &
    [\wti{F}_i(a), \wti{F}_{i+1}(b)] = -\dfrac{1}{2}\wti{A}_{i+2}(\ov{ab}), 
\vspace{1mm}\\
    {[}D, \wti{A}_i(a){]} = \wti{A}_i(D_ia), & 
     \mbox{\it where } D = (D_1, D_2, D_3) \in {\gd_4}^C,
\vspace{1mm}\\
    {[}\wti{A}_i(a), \wti{A}_{i}(b){]} \in {\gd_4}^C, &
    [\wti{A}_i(a), \wti{A}_{i+1}(b)] = - \dfrac{1}{2}\wti{A}_{i+2}(\ov{ab}), 
\vspace{1mm}\\
    {[}\wti{F}_1(e_i), \wti{F}_1(e_j){]} = G_{ij}.  &
    {[}\wti{A}_i(a), \wti{A}_{i+2}(b)] = \dfrac{1}{2}\wti{A}_{i+1}(\ov{ba}).   
\end{array} $$

{\bf Lemma 2.4.3.} (1) {\it Any non-zero element} $x \in \gC^C$ {\it can be 
transformed to $1$ by successive actions of ${\gd_4}^C$.}  
\vspace{1mm}

(2) $\gC^C$ {\it is a ${\gd_4}^C$-irreducible $C$-module.}
\vspace{1mm}

(3) ${\gd_4}^C\gC^C = \big\{\dsum_iD_ia_i \, | \, D_i \in {\gd_4}^C, a_i \in 
\gC^C \big\} = \gC^C. $
\vspace{1mm}
                 
{\bf Proof.} (1) Let $0 \neq x = \sum_{i=0}^7x_ie_i$, $x_i \in C$. Suppose $x_i
 \neq 0$ (moreover we may assume $i\neq 0$). Then we have
$$
       x \stackrel{G_{i0}}{\longrightarrow}x_ie_0 - x_0e_i 
         \stackrel{G_{j0} (j \neq i)}{\longrightarrow}x_ie_j 
         \stackrel{{x_j}^{-1}G_{0j}}{\longrightarrow}e_0 = 1.$$ 
The case $x_0 \neq 0$ is reduced to the above by the action of $G_{i0}$.
\vspace{1mm}

(2) Let $W \neq \{0\}$ be a ${\gd_4}^C$-invariant $C$-submodule of $\gC^C$. If 
$e_0 = 1 \in W$, then $G_{i0}e_0 = e_i$. Therefore $W = \gC^C$. Now, since a 
non-zero element $x$ of $W$ can be transformed to $1$ by ${\gd_4}^C$ from (1) 
above, $W$ 
contains $1$, so that $W = \gC^C$, because $G_{i0}e_0 = e_i$. 
\vspace{1mm}

(3) Since ${\gd_4}^C\gC^C$ is a ${\gd_4}^C$-invariant $C$-submodule of $\gC^C$,
 from the irreducibility of $\gC^C$ of (2), we have ${\gd_4}^C\gC^C = \gC^C$. 
\vspace{2mm}

Recall that any element $\delta \in {\gf_4}^C$ is uniquely expressed by
$$
       \delta = D + \wti{A}_1(a_1) + \wti{A}_2(a_2) + \wti{A}_3(a_3), \quad D 
\in {\gd_4}^C, a_i \in \gC^C, $$            
where $A_i(a_i) \in (\gM^{-})^C$ (Theorem 2.3.8).

\vspace{3mm}

{\bf Theorem 2.4.4.} {\it The Lie algebra ${\gf_4}^C$ is simple and so $\gf_4$ 
is also simple.}\vspace{2mm}

{\bf Proof.}  Let denote ${\wti{\gA}_i}^{\;C} = \{\wti{A}_i(a) \, | \, a \in \gC^C \}$ 
and $\wti{\gA}^C = {\wti{\gA}_1}^{\;C} \oplus {\wti{\gA}_2}^{\;C} \oplus {\wti{\gA}_3}^{\;C}$, then
$$
     {\gf_4}^C = {\gd_4}^C \oplus \wti{\gA}_1^{\ C} \oplus \wti{\gA}_2^{\ C} 
\oplus \wti{\gA}_3^{\ C} = {\gd_4}^C \oplus \wti{\gA}^C. $$
Let $p : {\gf_4}^C \to {\gd_4}^C$ and $q : {\gf_4}^C \to \wti{\gA}^C$ be 
projections of ${\gf_4}^C = {\gd_4}^C \oplus \wti{\gA}^C$. Now, let $\ga$ be a 
non-zero ideal of ${\gf_4}^C$. Then $p(\ga)$ is an ideal of ${\gd_4}^C$. 
Indeed, if $D \in p(\ga)$, then there exist $a_i \in \gC^C$, $i = 1, 2, 3$, such that 
$D + \displaystyle{\sum_{i=1}^3}\wti{A}_i(a_i) \in \ga$. For any $D' \in 
{\gd_4}^C$, we have 
$$
         \ga \ni \Big[D', D + \displaystyle{\sum_{i=1}^3}\wti{A}_i(a_i)\Big] 
         = [D',D] + \displaystyle{\sum_{i=1}^3}\wti{A}_i({D_i}'a_i)\; 
\mbox{(Lemma 2.4.2)},$$  
hence $[D', D] \in p(\ga)$. 
\vspace{1mm}

We show that either ${\gd_4}^C \cap \ga \neq \{ 0 \}$ or $\wti{\gA}^C \cap \ga 
\neq \{ 0 \}$. Assume that ${\gd_4}^C \cap \ga = \{ 0 \}$ and $\wti{\gA}^C \cap
 \ga = \{ 0 \}$. The mapping $p|\ga : \ga \to {\gd_4}^C$ is injective because 
$\gA^C \cap \ga = \{ 0 \}$. Since $p(\ga)$ is a non-zero ideal of ${\gd_4}^C$ 
and ${\gd_4}^C$ is simple, we have $p(\ga) = {\gd_4}^C$. Hence $\dim_C\ga = 
\dim_Cp(\ga) = \dim_C{\gd_4}^C = 28$. On the other hand, since ${\gd_4}^C \cap 
\ga =\{ 0 \}$, $q|\ga : \ga \to \wti{\gA}^C$ is also injective, we have 
$\dim_C\ga \le \dim_C\wti{\gA}^C = 8 \times 3 = 24$. This leads to a 
contradiction.
\vspace{1mm}

We now consider the following two cases.
\vspace{1mm}

(1) Case ${\gd_4}^C \cap \ga \neq \{ 0 \}$. From the simplicity of ${\gd_4}^C$,
 we have ${\gd_4}^C \cap \ga = {\gd_4}^C$, hence $\ga \supset {\gd_4}^C$. On 
the other hand, we have
$$
    \ga \supset [\ga, {\gf_4}^C] \supset [{\gd_4}^C, \wti{\gA}_i^{\ C}] = 
\wti{\gA}_i^{\ C}, \quad  i = 1, 2, 3. $$
The last equation follows from $[D, \wti{A}_i(a_i)] = \wti{A}_i(Da_i)$ 
(Lemma 2.4.2) and ${\gd_4}^C\gC^C = \gC^C$ (Lemma 2.4.3.(3)). Hence $\ga 
\supset {\gd_4}^C \oplus \wti{\gA}_1^{\ C} \oplus \wti{\gA}_2^{\ C} \oplus 
\wti{\gA}_3^{\ C} = {\gf_4}^C$.
\vspace{1mm}

(2) Case $\wti{\gA}^C \cap \ga \neq \{ 0 \}$. Choose a non-zero element 
$$
  \wti{A}_1(a_1) + \wti{A}_2(a_2) + \wti{A_3}(a_3) \in \wti{\gA}^C \cap \ga \subset \ga. $$
If $a_1 \neq 0$, under the actions of ${\gd_4}^C$ of Lemma 2.4.3.(1), we have 
$$
  \wti{A}_1(1) + \wti{A}_2(b) + \wti{A}_3(c) \in \ga. $$
If $b = c = 0$, then $\wti{A}_1(1) \in \ga$, hence from Lemma 2.4.2,
$$
    0 \neq [\wti{A}_1(1), \wti{A}_1(e_1)] \in {\gd_4}^C \cap \ga $$
(as for the first inequality, note that $[\wti{A}_1(1), \wti{A}_1(e_1)]F_2(1) = -2F_2(e_1))$, so the case is reduced to the case (1). If $b \neq 0$, then take the Lie bracket with $\wti{A}_1(1)$, then 
$$
    -\wti{A}_3(\ov{b}) + \wti{A}_2(\ov{c}) \in \ga \quad \mbox{and then $\quad \wti{A}_3(1) + \wti{A}_2(c') \in \ga$} $$
 under some actions of ${\gd_4}^C$. If $ c' = 0$, then this can be reduced to the case (1) as in the above. The case $c'\neq 0$ can also be reduced to the case (1) by the same argument as above. Thus we have $\ga = {\gf_4}^C$, which proves the simplicity of ${\gf_4}^C$.
\vspace{3mm}

{\bf Lemma 2.4.5.} (1) {\it For $\delta \in {\gf_4}^C$ and $A, B \in \gJ^C$, we have}
$$ 
          [\delta, [\wti{A}, \wti{B}]\,] 
            = [\wti{\delta A}, \wti{B}] + [\wti{A}, \wti{\delta B}].$$

(2) {\it Any element $\delta \in {\gf_4}^C$ is expressed as} $\delta = 
\displaystyle{\sum_{i}}[\wti{A}_i, \wti{B}_i], A_i, B_i \in \gJ^C$.
\vspace{2mm}

{\bf Proof.} (1) $[\delta, [\wti{A}, \wti{B}]\,]X = \delta[\wti{A}, \wti{B}]X - 
[ \wti{A}, \wti{B}]\delta X 
\vspace{1mm}\\
\begin{array}{l}
     \qquad \qquad = \delta(A \circ (B \circ X) - B \circ (A \circ X)) - (A 
\circ (B \circ \delta X) - B \circ (A \circ \delta X)) 
\vspace{1mm}\\
     \qquad \qquad = \delta A \circ (B \circ X) + A \circ (\delta B \circ X) +
  A \circ (B \circ \delta X) - \delta B \circ (A \circ X) 
\vspace{1mm}\\
     \qquad \qquad \;\; -B \circ (\delta A \circ X) - B \circ (A \circ \delta 
X) - A \circ (B \circ \delta X) + B \circ (A \circ \delta X) 
\vspace{1mm}\\
     \qquad \qquad = \delta A \circ (B \circ X) - B \circ (\delta A \circ X) +
A \circ (\delta B \circ X) - \delta B \circ (A \circ X) 
\vspace{1mm}\\
     \qquad \qquad = [\wti{\delta A}, \wti{B}]X + [\wti{A}, \wti{\delta B}]X, 
\quad X \in \gJ^C.
\end{array}$
\vspace{1mm}

(2) By (1) we see that $\ga = \{ \sum_i[\wti{A}_i, \wti{B}_i] \, | \, A_i,B_i 
\in \gJ^C \}$ is an ideal of ${\gf_4}^C$. From the simplicity of ${\gf_4}^C$, 
we have $\ga = {\gf_4}^C$.
\vspace{3mm}

{\bf Proposition 2.4.6.} (1) ${\gJ_0}^C = \{X \in \gJ^C \, | \, \tr(X) = 0 \}$
 {\it is an ${\gf_4}^C$-irreducible $C$-module.} 
\vspace{0.5mm}

(2) ${\gf_4}^C{\gJ_0}^C = \Big\{ \displaystyle{\sum_i}\delta_i B_i \, | \, \delta_i \in {\gf_4}^C,B_i \in {\gJ_0}^C \Big\} = {\gJ_0}^C$. 
\vspace{2mm}

{\bf Proof.} (1) Evidently ${\gJ_0}^C$ is an ${\gf_4}^C$-$C$-module (Lemma 2.2.1.(2)). Now, let $W$ be a non-zero ${\gf_4}^C$-invariant $C$-submodule of ${\gJ_0}^C$. 
\vspace{1mm}

Case (i) If $W$ contains some $F_i(x) \neq 0$ $(i = 1, 2, 3)$, then $W = 
{\gJ_0}^C$. Indeed, under some actions of ${\gd_4}^C$, we can assume $F_i(1) 
\in W$ (Lemma 2.4.3.(1)). Next, from
$$
     F_i(1) \;\stackrel{\wti{A}_i(1)}{\longrightarrow}\; 
     E_{i+1} - E_{i+2} \stackrel{\wti{A}_{i+2}(1) }{\longrightarrow}\;
     F_{i+2}(1) {\longrightarrow}\;
     \cdots \;{\longrightarrow}\; F_{i+1}(1), $$ 
and $\wti{A}_i(a)F_{i+2}(1) = F_{i+1}(\ov{a})$, etc., we have $E_i - E_{i+1}$,
$F_i(a) \in W$ and hence we have $W = {\gJ_0}^C$. Now, in the general case, 
choose a non-zero element $X$ from W:
$$
       0 \neq X = \xi(E_1 - E_2) + \eta(E_2 - E_3) 
              + F_1(x_1) + F_2(x_2) + F_3(x_3) \in W.$$

Case $x_1 \neq 0$. Under the action of ${\gd_4}^C$ of Lemma 2.4.3.(1), we have 
$$
       F_1(1) + F_2(y_2) + F_3(y_3) \in W. $$ 
If $y_2 = y_3 = 0$, then $F_1(1) \in W$, hence we can reduce to the case (i). If 
$y_2 \neq 0$, choose $a \in \gC^C$, $a \ne 0$, such that $(y_2, a) = 0$, and apply $\wti{A}_2(a)$, then we have $-F_3(\ov{a}) + F_1(\ov{ay_3}) \in W$. Again under the 
action of ${\gd_4}^C$, we have 
$$
          F_3(1) +  F_1(z_1) \in W. $$
If $z_1 = 0$, we can reduce to the case (i). If $z_1 \neq 0$, choose $b \in \gC^C$, $b \ne 0$, such that $(z_1, b) = 0$ and apply $\wti{A}_1(b)$, then $F_2(\ov{b}) \in W$, 
which is again reduced to the case (i). Similarly for $x_2 \neq 0$ or $x_3 
\neq 0$ the arguments above hold. 
\vspace{1mm}

Case $x_1 = x_2 = x_3 = 0$. Applying $\wti{A}_1(1)$ on the non-zero element 
$$
        X = \xi(E_1 - E_2) + \eta(E_2 - E_3) \in W, $$
we have $F_1(\xi - 2\eta) \in W$. If  $\xi - 2\eta \neq 0$, we can reduce to 
the case (i). If $\xi - 2\eta = 0$, applying $\wti{A}_3(1)$, we have $F_3(3\eta) \in W$ which is again the case (i). We have just proved that $W = {\gJ_0}^C$.
\vspace{1mm}

(2) Since ${\gf_4}^C{\gJ_0}^C$ is an ${\gf_4}^C$-invariant $C$-submodule of ${\gJ_0}^C$, we have ${\gf_4}^C{\gJ_0}^C = {\gJ_0}^C$ from (1).
\vspace{4mm}

{\bf 2.5. Killing form of ${\gf_4}^C$}
\vspace{3mm}

{\bf Lemma 2.5.1.} {\it For} $A, B, C, D \in \gJ^C$, {\it we have}
$$
       ([\wti{A}, \wti{B}]C, D) = ([\wti{C}, \wti{D}]A, B). $$

{\bf Proof.} \qquad $([\wti{A}, \wti{B}]C, D) = (A \circ (B \circ C), D) - (B \circ (A \circ C), D)$ 
\vspace{1mm}

$\qquad \qquad \qquad \qquad = (B \circ C,A \circ D) - (A \circ C, B \circ D) $
\vspace{1mm}

$\qquad \qquad \qquad \qquad = (C \circ (D \circ A), B) - (D \circ (C \circ A), B) = ([\wti{C}, \wti{D}]A, B).$
\vspace{3mm}

{\bf Definition.} In ${\gf_4}^C$, we define an inner product $(\delta_1, \delta_2)_4$ by 
$$
          (\delta, [\wti{A}, \wti{B}])_4 = (\delta A, B), \quad 
              \delta \in {\gf_4}^C,\;\;A, B \in \gJ^C.$$
More precisely, we define  
$$
     (\delta_1, \delta_2)_4 = \displaystyle{\sum_{i,j}}([\wti{A}_i, \wti{B}_i]C_j, D_j) = \displaystyle{\sum_{i,j}}([\wti{C}_j, \wti{D}_j]A_i, B_i), $$
for $\delta_1 = \sum_{i}[\wti{A}_i, \wti{B}_i]$, $\delta_2 = \sum_{j}[\wti{C}_j, \wti{D}_j]$, $A_i, B_i, C_j, D_j \in \gJ^C$  (Lemma 2.4.5.(2)). Then, Lemma 2.5.1 shows that the definition of the inner product $(\delta_1, \delta_2)_4$ is 
independent of the choice of the expressions of $\delta_1$, $\delta_2$, and 
that $(\delta_1,\delta_2)_4$ is symmetric. 

{\bf Lemma 2.5.2.} {\it The inner product} $(\delta_1, \delta_2)_4$ {\it of} ${\gf_4}^C$ {\it is} ${\gf_4}^C$-{\it adjoint invariant}: 
$$
       ([\delta, \delta_1], \delta_2)_4 + (\delta_1, [\delta, \delta_2])_4 = 0,
 \quad \delta, \delta_i \in {\gf_4}^C.$$

{\bf Proof.} It is sufficient to show the lemma for $\delta_1 = [\wti{A}, \wti{B}]$, $\delta_2 = [\wti{C}, \wti{D}]$, $A, B, C, D \in \gJ$.
$$
\begin{array}{l}
   ([\delta, [\wti{A}, \wti{B}]\,], [\wti{C}, \wti{D}])_4 + ([\wti{A}, \wti{B}], [\delta, [\wti{C}, \wti{D}]\,])_4  
\vspace{1mm}\\
  \;\;= ([\wti{\delta A}, \wti{B}] + [\wti{A}, \wti{\delta B}], [\wti{C}, \wti{D}])_4 + ([\wti{A}, \wti{B}], [\wti{\delta C}, \wti{D}] + [\wti{C}, \wti{\delta D}])_4 \;\; \mbox{(Lemma 2.4.5.(1))} 
\vspace{1mm}\\
  \;\;= ([\wti{\delta A}, \wti{B}]C, D) + ([\wti{A}, \wti{\delta B}]C, D) + ([\wti{A}, \wti{B}]\delta C, D) + ([\wti{A}, \wti{B}]C, \delta D) 
\vspace{1mm}\\
  \;\;= (\delta A \circ (B \circ C), D) - (B \circ (\delta A \circ C), D) + (A \circ (\delta B \circ C), D) - (\delta B \circ (A \circ C), D) 
\vspace{1mm}\\
\;\;\;\;\;+ (A \circ (B \circ \delta C), D) - (B \circ (A \circ \delta C), D) + (A \circ (B \circ C), \delta D) - (B \circ (A \circ C), \delta D),
\end{array}$$
which is equal to $0$, if we use the relation $(X \circ Y, \delta Z) = -(\delta X \circ Y, Z) - (X \circ \delta Y, Z)$ for the above last two terms.
\vspace{3mm}

{\bf Theorem 2.5.3.} {\it The Killing form $B_4$ of the Lie algebra ${\gf_4}^C$ is given by}
$$
    B_4(\delta_1, \delta_2) = 9(\delta_1, \delta_2)_4 = 3\tr(\delta_1\delta_2), \quad \delta_1, \delta_2 \in {\gf_4}^C.$$

{\bf Proof.} Since ${\gf_4}^C$ is simple (Theorem 2.4.4), there exist $k, k' \in C$ such that
$$
 B_4(\delta_1, \delta_2) = k(\delta_1, \delta_2)_4 = k'\tr(\delta_1\delta_2).$$
To determine these $k, k'$, we put $\delta = \delta_1 = \delta_2 = \wti{A}_1(1)$. Since $\wti{A}_1(1) = - 2[\wti{E}_3, \wti{F}_1(1)]$ (Lemma 2.4.2),  we have
\begin{eqnarray*}
      (\delta, \delta)_4 \!\!\! &=& \!\!\! (\wti{A}_1(1), \wti{A}_1(1))_4 = -2(\wti{A}_1(1), [\wti{E}_3, \wti{F}_1(1)])_4 \\
   \!\!\! &=& \!\!\! -2(\wti{A}_1(1)E_3, F_1(1)) = -(F_1(1), F_1(1)) = - 2.
\end{eqnarray*}
On the other hand, $(\mbox{ad}\delta)^2$ is calculated as follows.
$$
\begin{array}{l}
    [\wti{A}_1(1), [\wti{A}_1(1), G_{i0}]\,] = -[\wti{A}_1(1), \wti{A}_1(e_i)] = -G_{i0}, \;\; i \neq 0 
\vspace{1mm}\\
    {[}\wti{A}_1(1), [\wti{A}_1(1), \wti{A}_1(e_i)]\,{]} = [\wti{A}_1(1), G_{i0}]  = - \wti{A}_1(e_i),\;\; i \neq 0 
\vspace{1mm}\\
    {[}\wti{A}_1(1), [\wti{A}_1(1), \wti{A}_2(e_i)]\,{]} = \displaystyle{\frac{1}{2}}[\wti{A}_1(1), \wti{A}_3(e_i)] = - \displaystyle{\frac{1}{4}}\wti{A}_2(e_i), \vspace{1mm}\\
  {[}\wti{A}_1(1), [\wti{A}_1(1), \wti{A}_3(e_i)]\,{]} = - \displaystyle{\frac{1}{2}}[\wti{A}_1(1), \wti{A}_2(e_i)] = - \displaystyle{\frac{1}{4}}\wti{A}_3(e_i),
\vspace{1mm}\\
   \mbox{the others} \; = 0.
\end{array} $$
Hence we have
$$
      B_4(\delta, \delta) = \tr((\mbox{ad}\wti{A}_1(1))^2) = (-1) \times 7 \times 2 + \Big(-\dfrac{1}{4} \Big) \times 8 \times 2 = -18.$$
Therefore $k = 9$. Next, we will calculate $\tr(\delta\delta)$.

$$
\begin{array}{l}
      \wti{A}_1(1)\wti{A}_1(1)E_2 = - \displaystyle{\frac{1}{2}}\wti{A}_1(1)F_1(1) = - \displaystyle{\frac{1}{2}}(E_2 - E_3),
\vspace{1mm} \\
      \wti{A}_1(1)\wti{A}_1(1)E_3 = \displaystyle{\frac{1}{2}}\wti{A}_1(1)F_1(1)    = \displaystyle{\frac{1}{2}}(E_2 - E_3),
\vspace{1mm}\\
      \wti{A}_1(1)\wti{A}_1(1)F_1(1) = \wti{A}_1(1)(E_2 - E_3) = - F_1(1), 
\vspace{1mm}\\
      \wti{A}_1(1)\wti{A}_1(1)F_2(e_i) = \displaystyle{\frac{1}{2}}\wti{A}_1(1)F    _3(\ov{e_i}) = - \displaystyle{\frac{1}{4}}F_2(e_i), 
\vspace{1mm}\\
      \wti{A}_1(1)\wti{A}_1(1)F_3(e_i) = - \displaystyle{\frac{1}{2}}\wti{A}_1(1)F_3(\ov{e_i}) = - \displaystyle{\frac{1}{4}}F_3(e_i),
\vspace{1mm}\\
  \mbox{the others}\; = 0.
\end{array} $$
Hence we have
$$
      \tr(\delta\delta) = \tr((\wti{A}_1(1))^2) = \Big(-\displaystyle{\frac{1}{2}}\Big) \times 2 - 1 + \Big(-\displaystyle{\frac{1}{4}}\Big) \times 8 \times 2 = - 6.$$
Therefore $k' = 3$.
\vspace{3mm}

{\bf Lemma 2.5.4.} {\it For a non-zero element $A \in {\gJ_0}^C$, there exists $B \in {\gJ_0}^C$ such that $[\wti{A}, \wti{B}] \neq 0$.}
\vspace{2mm}

{\bf Proof.} Assume that $[\wti{A},\wti{B}] = 0$ for all $B \in {\gJ_0}^C$. Then  $0 = (\delta, [\wti{A}, \wti{B}])_4 = - (\delta, [\wti{B}, \wti{A}])_4 = - (\delta B, A)$ for any $\delta \in {\gf_4}^C$. Since ${\gf_4}^C{\gJ_0}^C = {\gJ_0}^C$ (Proposition 2.4.6.(2)), we have $({\gJ_0}^C, A) = 0$, so that $A = 0$.
\vspace{4mm}

{\bf 2.6. Roots of ${\gf_4}^C$}
\vspace{3mm}

Before we obtain the roots of the Lie algebra ${\gf_4}^C$, we recall the roots of the Lie algebra ${\gD_4}^C$: 
$$
   {\gD_4}^C = \{ D \in \Hom_C(\gC^C) \, | \, (Dx, y) + (x, Dy) = 0 \}.$$
If we let $H_k = -iG_{k 4+k}$ for $k = 0, 1, 2, 3$, then
$$
      \gh = \{ H = \sum_{k=0}^3\lambda_iH_k \, | \, \lambda_k \in C \}$$
is a Cartan subalgebra of ${\gD_4}^C$, and roots of ${\gD_4}^C$ relative to $\gh$ are given by
$$
      \pm (\lambda_k - \lambda_l),\;\; \pm (\lambda_k + \lambda_l), \quad 0 \le  k < l \le 3. $$
The root vectors associated with these roots are respectively given by
\begin{eqnarray*}
     \;\;\;\lambda_k - \lambda_l &:& 
          (G_{kl} + G_{4+k 4+l}) - i(G_{k 4+l} + G_{l 4+k}), 
\vspace{1mm}\\
     -\lambda_k + \lambda_l &:& 
          i(G_{k l} + G_{4+k 4+l}) - (G_{k 4+l} + G_{l 4+k}),
\vspace{1mm}\\
          \;\;\;\lambda_k + \lambda_l &:& 
          (G_{kl} - G_{4+k 4+l}) + i(G_{k 4+l} - G_{l 4+k}), 
\vspace{1mm}\\
     -\lambda_k - \lambda_l &:& 
          i(G_{k l} - G_{4+k 4+l}) + (G_{k 4+l} - G_{l 4+k}).
\end{eqnarray*}

\noindent The Lie algebra ${\gD_4}^C$ is contained in ${\gf_4}^C$ as
$$
      {\gD_4}^C \ni D \; \to \; (D, \nu D, \kappa\pi D) \in {\gd_4}^C \subset {\gf_4}^C $$
(Theorem 1.3.6, Proposition 2.3.7). Outer automorphisms $\nu, \pi$ of ${\gD_4}^C$ induce $C$-linear transformations of $\gh$ and the matrices of $\nu$, $\pi$ with respect to the $C$-basis $H_0, H_1, H_2, H_3$ of $\gh$ are respectively given by
$$
       \nu = \pi\kappa = \frac{1}{2}\pmatrix{-1 & -1 & \hfill 1 & -1 \cr
                                           \hfill 1 & \hfill 1 & \hfill 1 & -1 \cr
                                           -1 & \hfill 1 & \hfill 1 & \hfill 1 \cr
                                           \hfill 1 & -1 & \hfill 1 & \hfill 1}, \quad
       \pi = \frac{1}{2}\pmatrix{\hfill 1 & -1 & \hfill 1 & -1 \cr
                                 -1 & \hfill 1 & \hfill 1 & -1 \cr 
                                 \hfill 1 & \hfill 1 & \hfill 1 & \hfill 1 \cr
                                 -1 & -1 & \hfill 1 & \hfill 1}  $$
(Lemma 1.3.1). Note that they are orthogonal matrices.
\vspace{3mm}

{\bf Theorem 2.6.1.} {\it The rank of the Lie algebra ${\gf_4}^C$ is $4$. The roots of ${\gf_4}^C$ relative to some Cartan subalgebra of ${\gf_4}^C$ are given by}
$$
\begin{array}{c}
     \pm (\lambda_k - \lambda_l), \quad \pm (\lambda_k + \lambda_l) \;\;,\;\; 0 \le k < l \le 3, 
\vspace{1mm}\\
       \pm \lambda_0, \quad \pm \lambda_1, \quad \pm \lambda_2, \quad \pm \lambda_3,  
\end{array} $$
\vspace{-2.5mm}
$$
\begin{array}{ll}
   \pm \dfrac{1}{2}(- \lambda_0 - \lambda_1 + \lambda_2 - \lambda_3), & \quad
   \pm \dfrac{1}{2}(\hfill \, \lambda_0 + \lambda_1 + \lambda_2 - \lambda_3), 
\vspace{1mm}\\
   \pm \dfrac{1}{2}(- \lambda_0 + \lambda_1 + \lambda_2 + \lambda_3), & \quad
   \pm \dfrac{1}{2}(\hfill \, \lambda_0 - \lambda_1 + \lambda_2 + \lambda_3), 
\vspace{1mm}\\
   \pm \dfrac{1}{2}(\hfill \, \lambda_0 - \lambda_1 + \lambda_2 - \lambda_3), & \quad    
   \pm \dfrac{1}{2}(- \lambda_0 + \lambda_1 + \lambda_2 - \lambda_3), 
\vspace{1mm}\\
   \pm \dfrac{1}{2}(\hfill \, \lambda_0 + \lambda_1 + \lambda_2 + \lambda_3), & \quad    
   \pm \dfrac{1}{2}(- \lambda_0 - \lambda_1 + \lambda_2 + \lambda_3). 
\end{array} $$

{\bf Proof.}  We use the decomposition in Theorem 2.4.4:
$$
      {\gf_4}^C = {\gd_4}^C \oplus \wti{\gA}_1^{\ C} \oplus \wti{\gA}_2^{\ C} \oplus \wti{\gA}_3^{\ C}.$$
Let $\gh = \Big\{ H = \dsum_{k=0}^3\lambda_kH_k \, | \,  \lambda_k \in C \Big\} \subset {\gd_4}^C \subset {\gf_4}^C$.  Since $\gh$ is a Cartan subalgebra of ${\gd_4}^C$ (it will be also a Cartan subalgebra of ${\gf_4}^C$), the roots of ${\gd_4}^C$:
$$
      \pm (\lambda_k - \lambda_l),\;\; \pm (\lambda_k + \lambda_l), \quad 0 \le k < l \le 3 $$
are also roots of ${\gf_4}^C$. Furthermore, from $[H, \wti{A}_1(a)] = \wti{A}_1(Ha)$ (Lemma 2.4.2), where
\begin{eqnarray*}
        H(e_k + ie_{4+k}) \!\!\! &=& \!\!\! - \sum_{j=0}^3\lambda_kiG_{j4+j}(e_k + ie_{4+k}) \\
 \!\!\! &=& \!\!\! -i\lambda_k(- e_{4+k} + ie_k) = \lambda_k(e_{k} + ie_{4+k}),
\end{eqnarray*}
we see that $\lambda_k$ is a root of ${\gf_4}^C$ and $\wti{A}_1(e_k + ie_{4+k})$ is an associated root vector for $0 \le k \le 3$. Similarly, $-\lambda_k$ for $0 \le k \le 3$ is root of ${\gf_4}^C$ and $\wti{A}_1(e_k - ie_{4+k})$ is an associated root vector. Next, by Lemma 2.4.2 we have $[H, \wti{A}_2(a)] = \wti{A}_2((\nu H)a)$,  where
$$
\begin{array}{l}
      \nu H = \nu\Big(\dsum_{k=0}^3\lambda_kH_k\Big) 
\vspace{1mm}\\
\quad \;\;\;
   = \dfrac{1}{2}(\lambda_0(- H_0 + H_1 - H_2 + H_3) + \lambda_1(- H_0 + H_1 + H_2 - H_3) 
\vspace{1mm}\\
\qquad \;\; 
    + \lambda_2(H_0 + H_1 + H_2 + H_3) + \lambda_3(- H_0 - H_1 + H_2 +H _3)) 
\vspace{1mm}\\
\quad \;\;\;
   = \dfrac{1}{2}(- \lambda_0 - \lambda_1 + \lambda_2 - \lambda_3)H_0 + \dfrac{1}{2}(\lambda_0 + \lambda_1 + \lambda_2 - \lambda_3)H_1 
\vspace{1mm}\\
\qquad \;\; 
   + \dfrac{1}{2}(- \lambda_0 + \lambda_1 + \lambda_2 + \lambda_3)H_2 + \dfrac{1}{2}(\lambda_0 - \lambda_1 + \lambda_2 + \lambda_3)H_3,
\end{array}$$
and so we see that coefficients of $H_0$, $H_1$, $H_2$, $H_3$ are roots of ${\gf_4}^C$ and that $\wti{A}_2(e_k+i + ie_{4+k})$ is associated root vector for $0 \le k \le 3$. The roots above with negative sign are also roots of ${\gf_4}^C$ and $\wti{A}_2(e_k - ie_{4+k})$ are associated root vectors. Finally, from $[H, \wti{A}_3(a)] = \wti{A}_3((\kappa\pi H)a)$ (Lemma 2.4.2), where
$$
\begin{array}{l}
     \kappa\pi H = \dfrac{1}{2}(- \lambda_0 + \lambda_1 - \lambda_2 + \lambda_3)H_0 + \dfrac{1}{2}(- \lambda_0 + \lambda_1 + \lambda_2 - \lambda_3)H_1 
\vspace{1mm}\\
\qquad \;\;\; 
     + \dfrac{1}{2}(\lambda_0 + \lambda_1 + \lambda_2 + \lambda_3)H_2 + \dfrac{1}{2}(- \lambda_0 - \lambda_1 + \lambda_2 + \lambda_3)H_3,
\end{array}$$
we obtain the remainders of roots. 
\vspace{3mm}

{\bf Theorem 2.6.2.} {\it In the root system of Theorem} 2.6.1,
$$
   \alpha_1 = \lambda_0 - \lambda_1, \;\; \alpha_2 = \lambda_1 - \lambda_2,\;\;    \alpha_3 = \lambda_2,\;\; \alpha_4 = \frac{1}{2}(- \lambda_0 - \lambda_1 - \lambda_2 + \lambda_3) $$
{\it is a fundamental root system of the Lie algebra ${\gf_4}^C$ and
$$
        \mu = 2\alpha_1 + 3\alpha_2 + 4\alpha_3 + 2\alpha_4$$
is the highest root. The Dynkin diagram and the extended Dynkin diagram of ${\gf_4}^C$ are respectively given by}
\vspace{-2mm}

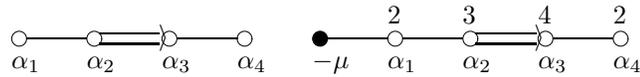
\setlength{\unitlength}{1mm}
\begin{picture}(150,20)
\put(20,10){\circle{2}} \put(19,6){$\alpha_1$}
\put(21,10){\line(1,0){8}}
\put(30,10){\circle{2}} \put(29,6){$\alpha_2$}
\put(30.7,10.7){\line(1,0){8}}
\put(30.7,9.3){\line(1,0){8}}
\put(38.5,9.2){$\rangle$}
\put(40,10){\circle{2}} \put(39,6){$\alpha_3$}
\put(41,10){\line(1,0){8}}
\put(50,10){\circle{2}} \put(49,6){$\alpha_4$}
\put(60,10){\circle*{2}} \put(59,6){$-\mu$} 
\put(61,10){\line(1,0){8}}
\put(70,10){\circle{2}} \put(69,6){$\alpha_1$} \put(69,12){$2$}
\put(71,10){\line(1,0){8}}
\put(80,10){\circle{2}} \put(79,6){$\alpha_2$} \put(79,12){$3$}
\put(80.7,10.7){\line(1,0){8}}
\put(80.7,9.3){\line(1,0){8}}
\put(88.5,9.2){$\rangle$}
\put(90,10){\circle{2}} \put(89,6){$\alpha_3$} \put(89,12){$4$}
\put(91,10){\line(1,0){8}}
\put(100,10){\circle{2}} \put(99,6){$\alpha_4$} \put(99,12){$2$}
\end{picture}

{\bf Proof.}  All positive roots of ${\gf_4}^C$ are expressed by $\alpha_1$, $\alpha_2$, $\alpha_3$, $\alpha_4$ as follows.
$$
\begin{array}{llrrrrrrr}
    \lambda_0 \!\!\! &= \alpha_1 \!\!\!\! &+& \!\!\!\! \alpha_2 \!\!\!\! &+& \!\!\!\! \alpha_3 &&
\vspace{1mm}\\
    \lambda_1 \!\!\! &= && \!\!\!\! \alpha_2 \!\!\!\! &+& \!\!\!\!\alpha_3  &&
\vspace{1mm}\\
    \lambda_3 \!\!\! &= &&  \!\!\!\! & \!\!\!\! & \!\!\!\! \alpha_3 &&
\vspace{1mm}\\
    \lambda_3 \!\!\! &= \alpha_1 \!\!\!\! &+& \!\!\!\! 2\alpha_2 \!\!\!\! &+& \!\!\!\! 3\alpha_3& \!\!\!\! +& \!\!\!\! 2\alpha_4,
\end{array}$$ 
\vspace{-2mm}
$$
\begin{array}{llrrrrrrr}
   \hfill \lambda_0 - \lambda_1 \!\!\! &= \!\!\! &\alpha_1&  
\vspace{1mm}\\
   \hfill\lambda_0 - \lambda_2 \!\!\! &= \!\!\! &\alpha_1& \!\!\! + \!\!\! &\alpha_2& 
\vspace{1mm}\\
  - \lambda_0 + \lambda_3 \!\!\! &= && \!\!\! &\alpha_2& \!\!\! + \!\!\! &2\alpha_3& \!\!\! + \!\!\! &2\alpha_4
\vspace{1mm}\\
   \hfill \lambda_1 - \lambda_2 \!\!\! &= && \!\!\! &\alpha_2& 
\vspace{1mm}\\
   - \lambda_1 + \lambda_3 \!\!\! &= \!\!\! &\alpha_1& \!\!\! + \!\!\! &\alpha_2& \!\!\! + \!\!\! &2\alpha_3& \!\!\! + \!\!\! &2\alpha_4
\vspace{1mm}\\
   - \lambda_2 + \lambda_3 \!\!\! &= \!\!\! &\alpha_1& \!\!\! + \!\!\! &2\alpha_2& \!\!\! + \!\!\! &2\alpha_3& \!\!\! + \!\!\! &2\alpha_4
\vspace{2mm}\\
   \hfill \lambda_0 + \lambda_1 \!\!\! &= \!\!\! &\alpha_1& \!\!\!+ \!\!\! &2\alpha_2& \!\!\! + \!\!\! &2\alpha_3& 
\vspace{1mm}\\
   \hfill \lambda_0 + \lambda_2 \!\!\! &= \!\!\! &\alpha_1& \!\!\! + \!\!\! &\alpha_2& \!\!\! + \!\!\! &2\alpha_3&
\vspace{1mm}\\
   \hfill \lambda_0 + \lambda_3 \!\!\! &= \!\!\! &2\alpha_1& \!\!\! + \!\!\! &3\alpha_2& \!\!\! + \!\!\! &4\alpha_3& \!\!\! + \!\!\! &2\alpha_4
\vspace{1mm}\\
   \hfill \lambda_1 + \lambda_2 \!\!\! &= && \!\!\! &\alpha_2& \!\!\! + \!\!\! &2\alpha_3&
\vspace{1mm}\\
   \hfill \lambda_1 + \lambda_3 \!\!\! &= \!\!\! &\alpha_1& \!\!\! + \!\!\! &3\alpha_2& \!\!\! + \!\!\! &4\alpha_3& \!\!\! + \!\!\! &2\alpha_4
\vspace{1mm}\\
   \hfill \lambda_2 + \lambda_3 \!\!\! &= \!\!\! &\alpha_1& \!\!\! + \!\!\! &2\alpha_2& \!\!\! + \!\!\! &4\alpha_3& \!\!\! + \!\!\! &4\alpha_4,
\end{array}$$
$$
\begin{array}{llrrrrrrrrr}
\displaystyle{\frac{1}{2}}(\hfill \lambda_0 + \lambda_1 + \lambda_2 + \lambda_3) \!\!\! & = \!\!\! &\alpha_1& \!\!\! + \!\!\! &2\alpha_2& \!\!\! + \!\!\! &3\alpha_3& \!\!\! + \!\!\! &\alpha_4&
\vspace{1mm}\\
\displaystyle{\frac{1}{2}}(- \lambda_0 - \lambda_1 - \lambda_2 + \lambda_3) \!\!\! &= && \!\!\! && \!\!\!  && \!\!\! &\alpha_4&
\vspace{1mm}\\
\displaystyle{\frac{1}{2}}(\hfill \lambda_0 + \lambda_1 - \lambda_2 + \lambda_3) \!\!\! &= \!\!\! &\alpha_1& \!\!\! + \!\!\! &2\alpha_2& \!\!\! + \!\!\! &2\alpha_3& \!\!\! + \!\!\! &\alpha_4&
\vspace{1mm}\\
\displaystyle{\frac{1}{2}}(\hfill \lambda_0 - \lambda_1 + \lambda_2 + \lambda_3) \!\!\! &= \!\!\! &\alpha_1& \!\!\! + \!\!\! &\alpha_2& \!\!\! + \!\!\! &2\alpha_3& \!\!\! + \!\!\! &\alpha_4&
\end{array}$$
$$
\begin{array}{llrrrrrrrrr}
\displaystyle{\frac{1}{2}}(- \lambda_0 - \lambda_1 + \lambda_2 + \lambda_3) \!\!\! &= && \!\!\!\! && \!\!\! &\alpha_3& \!\!\! + \!\!\! &\alpha_4&
\vspace{1mm}\\
\displaystyle{\frac{1}{2}}(- \lambda_0 + \lambda_1 - \lambda_2 + \lambda_3) \!\!\! &= && \!\!\! &\alpha_2& \!\!\! + \!\!\!\! &\alpha_3& \!\!\! + \!\!\!\! &\alpha_4&
\vspace{1mm}\\
\displaystyle{\frac{1}{2}}(\hfill \lambda_0 - \lambda_1 - \lambda_2 + \lambda_3) \!\!\! &= \!\!\! &\alpha_1& \!\!\! + \!\!\!\! &\alpha_2& \!\!\! + \!\!\! &\alpha_3& \!\!\! + \!\!\! &\alpha_4&
\vspace{1mm}\\
\displaystyle{\frac{1}{2}}(- \lambda_0 + \lambda_1 + \lambda_2 + \lambda_3) \!\!\! &= && \!\!\! &\alpha_2& \!\!\! + \!\!\! &2\alpha_3& \!\!\! + \!\!\! &\alpha_4.&
\end{array}$$ 
Hence ${\mit\Pi} = \{ \alpha_1, \alpha_2, \alpha_3, \alpha_4 \}$ is a fundamental root system of ${\gf_4}^C$. The real part $\gh_{\sR}$ of $\gh$ is
$$
      \gh_{\sR} = \{H = \dsum_{k=0}^3\lambda_kH_k \, | \, \lambda_k \in \R \}. $$
The Killing form $B_4$ of ${\gf_4}^C$ is $B_4(\delta_1, \delta_2) = 3\tr(\delta_1\delta_2)$ (Theorem 2.5.3), so that
$$
         B_4(H,H') = 18\sum_{k=0}^3\lambda_k{\lambda_k}', \quad
      H = \sum_{k=0}^3\lambda_kH_k, H' = \sum_{k=0}^3{\lambda_k}'H_k \in \gh_{\sR}.$$
Indeed, since 
$$
\begin{array}{ll}
      HE_i = 0, \quad i = 1, 2, 3, 
\vspace{1mm}\\
      HF_1(x) = F_1(Hx), \, \, HF_2(x) = F_2((\nu H)x), \, \, HF_3(x) = F_3((\kappa\pi H)x),
\end{array} $$
we have 

$$
\begin{array}{l}
    B_4(H,H') = \dsum_{k,l=0}^3\lambda_k{\lambda_l}'B_4(H_k, H_l) 
\vspace{1mm}\\
   \qquad 
    = 3\displaystyle{\sum_{k,l}}(\lambda_k{\lambda_l}'(\tr(H_kH_l) + \tr((\nu H_k)(\nu H_l)) + \tr((\kappa\pi H_k)(\kappa\pi H_l)))) 
\\
   \qquad  = 3\displaystyle{\sum_{k,l}}\lambda_k{\lambda_l}'(2 + 2 + 2)\delta_{kl} = 18\sum_{k=0}^3\lambda_k{\lambda_k}'.
\end{array} $$
Now, the canonical elements $H_{\alpha_i} \in \gh_{\sR}$ corresponding to $\alpha_i$ \, ($B_4 (H_\alpha, H) = \alpha(H), H \in \gh_{\sR}$) are determined by 
$$
\begin{array}{ll}
       H_{\alpha_1} = \displaystyle{\frac{1}{18}}(H_0 - H_1), & 
       H_{\alpha_2} = \displaystyle{\frac{1}{18}}(H_1 - H_2), 
\vspace{1mm}\\
       H_{\alpha_3} = \displaystyle{\frac{1}{18}}H_2,&
       H_{\alpha_4} = \displaystyle{\frac{1}{36}}(- H_0 - H_1 - H_2 + H_3).
\end{array} $$
Hence we have

$$
       (\alpha_1, \alpha_1) = B_4(H_{\alpha_1}, H_{\alpha_1}) = 18\frac{1}{18}\frac{1}{18}(1 + 1) = \frac{1}{9}$$
and the other inner products are similarly calculated. Hence, the inner products induced by the Killing form $B_4$ between $\alpha_1$, $\alpha_2$, $\alpha_3$, $\alpha_4$ and $-\mu$ are given by 
\vspace{2mm}

\qquad \qquad
   $(\alpha_1, \alpha_1) = (\alpha_2, \alpha_2) = \displaystyle{\frac{1}{9}}, \quad
  (\alpha_3, \alpha_3) = (\alpha_4, \alpha_4) = \displaystyle{\frac{1}{18}}$, 
\vspace{1mm}

\qquad \qquad
  $(\alpha_1, \alpha_2) = - \displaystyle{\frac{1}{18}}, \quad
  (\alpha_2, \alpha_3) = - \displaystyle{\frac{1}{18}}, \quad
  (\alpha_3, \alpha_4) = - \displaystyle{\frac{1}{36}}$,
\vspace{1mm}

\qquad \qquad
  $(\alpha_1, \alpha_3) = (\alpha_1, \alpha_4) = (\alpha_2, \alpha_4) = 0$,
\vspace{1mm}

\qquad \qquad
  $(-\mu, -\mu) = \displaystyle{\frac{1}{9}}, \quad
  (-\mu, \alpha_1) = - \displaystyle{\frac{1}{18}}, \quad (-\mu, \alpha_i) = 0,\;\; i = 2, 3, 4$,
\vspace{2mm}

\noindent using them, we can draw the Dynkin diagram and the extended Dynkin diagram of ${\gf_4}^C$.                 
\vspace{2mm}

According Borel-Siebenthal theory, the Lie algebra $\gf_4$ has three subalgebras as maximal subalgebras with the maximal rank 4. 
\vspace{1mm}

(1) The first one is a subalgebra of type $B_4$ obtained as the fixed points of an involution $\sigma$ of $\gf_4$.
\vspace{1mm}

(2) The second one is a subalgebra of type $C_1 \oplus C_3$ obtained as the fixed points of an involution $\gamma$ of $\gf_4$.
\vspace{1mm}

(3) The third one is a subalgebra of type $A_2 \oplus A_2$ obtained as the fixed points of an automorphism $w$ of order 3 of $\gf_4$.
\vspace{1mm}

These subalgebras will be realized as subgroups of the group $F_4$ in Theorems 2.9.1, 2.11.2 and 2.12.2, respectively.
\vspace{3mm}

{\bf 2.7. Subgroup $Spin(9)$ of $F_4$}
\vspace{3mm}

{\bf Theorem 2.7.1.} \qquad $\{ \alpha \in F_4 \, | \, \alpha E_i = E_i,\,i = 1, 2, 3 \} \cong Spin(8).$ 
\vspace{1mm}

\noindent (From now on, we identify these groups).
\vspace{2mm}

{\bf Proof.} We recall the group $Spin(8)$ of Theorem 1.16.2 and define a mapping $\varphi : Spin(8) \to D_4 = \{ \alpha \in F_4 \, | \, \alpha E_i = E_i,\,i = 1, 2, 3 \} $ by
$$
         \varphi(\alpha_1, \alpha_2, \alpha_3)
               \pmatrix{\xi_1 & x_3 & \ov{x}_2 \cr 
                        \ov{x}_3 & \xi_2 & x_1 \cr
                        x_2 & \ov{x}_1 & \xi_3} 
             = \pmatrix{\xi_1 & \alpha_3x_3 & \ov{\alpha_2x_2} 
\vspace{0.5mm}\cr 
                        \ov{\alpha_3x_3} & \xi_2 & \alpha_1x_1 
\vspace{0.5mm}\cr
                        \alpha_2x_2 & \ov{\alpha_1x_1} & \xi_3}. $$
We first prove that $\alpha = \varphi(\alpha_1, \alpha_2, \alpha_3) \in D_4$. Indeed, the fact that $\det\,(\alpha X) = \det\,X$ can be seen by observing that
$$
\begin{array}{l}
  R((\alpha_1x_1)(\alpha_2x_2)(\alpha_3x_3)) = R((\ov{\alpha_3(\ov{x_1x_2})})\alpha_3x_3)
\vspace{1mm}\\
 \qquad \quad = (\alpha_3(\ov{x_1x_2}), \alpha_3x_3) = (\ov{x_1x_2}, x_3) = R(x_1x_2x_3),\end{array}$$
which together with $\alpha E = E$, shows that $\alpha \in F_4$ and $\alpha E_i = E_i, i = 1, 2, 3$. Therefore $\alpha \in D_4$. Certainly $\varphi$ is a homomorphism. We shall show that $\varphi$ is onto. Let $\alpha \in D_4$. We put
$$
    \gJ_i = \{ F_i(x) \, | \, x \in \gC \}
         = \{ X \in \gJ \, | \, 2E_{i+1} \circ X = 2E_{i+2} \circ X = X \}, \quad i = 1, 2, 3. $$
Since $\alpha X \in \gJ_i$, $X \in \gJ_i$, $\alpha$ induces $\R$-isomorphisms $\alpha : \gJ_i \to \gJ_i$ and $\alpha_i : \gC \to \gC$ satisfying
$$
     \alpha F_i(x) = F_i(\alpha_i x), \quad x \in \gC, $$ 
for $i = 1, 2, 3$. Applying $\alpha$ on $F_i(x) \circ F_i(y) = (x, y)(E_{i+1} + E_{i+2})$, since the left side becomes $F_i(\alpha_ix) \circ F_i(\alpha_iy) = (\alpha_ix, \alpha_iy)(E_{i+1} + E_{i+2})$, we have 
$$
       (\alpha_ix, \alpha_iy) = (x, y), \quad x,y \in \gC.$$
Hence $\alpha_i \in O(8)$, $i = 1, 2, 3$. Moreover, by applying $\alpha$ on $F_1(x) \circ F_2(y) = \displaystyle{\frac{1}{2}}F_3(\ov{xy})$, we see that

$$
     (\alpha_1x)(\alpha_2y) = \ov{\alpha_3(\ov{xy})}, \quad x,y \in \gC.$$
From Lemma 1.14.4, we have $\alpha_1, \alpha_2, \alpha_3 \in SO(8)$, so that $(\alpha_1, \alpha_2, \alpha_3) \in Spin(8)$ and $\varphi(\alpha_1, \alpha_2, \alpha_3) = \alpha$. Therefore $\varphi$ is onto. Evidently $\Ker\,\varphi = \{ (1,1,1) \}$. Consequently $\varphi$ is an isomorphism.
\vspace{2mm}

We shall study the following subgroup $(F_4)_{E_1}$ of $F_4$: 
$$
         (F_4)_{E_1} = \{ \alpha \in F_4  \, | \, \alpha E_1 = E_1 \}.$$
We define $\R$-vector subspaces $\gJ_{01}$, $\gJ_{23}$ of $\gJ$ respectively by$$
      \gJ_{01} = \{ X \in \gJ  \, | \, E_1 \circ X = 0, \tr(X) = 0 \}
              = \Bigl\{\pmatrix{ 0 & 0 & 0 \cr 
                                 0 & \xi & x \cr 
                                0 & \ov{x} & - \xi} \, \Big| \, \xi \in \R, x \in \gC \Bigl\}, $$
\vspace{-2mm}
$$
      \gJ_{23} = \{ Y \in \gJ  \, | \, 2E_1 \circ Y = Y \}
              = \Bigl\{\pmatrix{0 & y_3 & \ov{y}_2 \cr
                                \ov{y}_3 & 0 & 0 \cr
                                y_2 & 0 & 0} \, \Big| \, y_2,y_3 \in \gC \Bigl\}. $$

{\bf Lemma 2.7.2.} {\it Suppose that we are given an element $A \in \gJ_{01}$ such that $(A, A) = 2$. Choose any element $X_0 \in \gJ_{01}$ such that $(X_0, X_0) = 2$, $(A, X_0) = 0$. Choose any element $Y_0 \in \gJ_{23}$ such that $(Y_0, Y_0) = 2$, $2A \circ Y_0 = - Y_0$. Let
$$ 
                  Z_0 = 2X_0 \circ Y_0. $$
Choose any element $X_1 \in \gJ_{01}$ such that $(X_1, X_1) = 2$, $(A, X_1) = (X_0, X_1) = 0$. Choose any element $X_2 \in \gJ_{01}$ such that $(X_2, X_2) = 2$, $(A, X_2) = (X_0, X_2) = (X_1, X_2) = 0$. Let
$$
        Y_1 = - 2Z_0 \circ X_1, \quad Z_2 = -2X_2 \circ Y_0, \quad
        X_3 = - 2Y_1 \circ Z_2.$$
Finally, choose any element $X_4 \in \gJ_{01}$ such that $(X_4, X_4) = 2$, $(A, X_4) = (X_0, X_4) = (X_1, X_4) = (X_2, X_4) = (X_3, X_4) = 0$. Let
$$
\begin{array}{lll}
   Z_4 = - 2X_4 \circ Y_0, & Y_2 = - 2Z_0 \circ X_2, & Y_3 = - 2Z_0 \circ X_3,
\vspace{1mm}\\
   X_5 = -2Y_1 \circ Z_4, & X_6 = 2Y_2 \circ Z_4, & X_7 = - 2Y_3 \circ Z_4
\end{array} $$
and moreover let
$$
\begin{array}{l}
        Y_i = - 2Z_0 \circ X_i, \quad i = 4, 5, 6, 7, 
\vspace{1mm}\\
        Z_i =  -2X_i \circ Y_0, \quad i= 1, 3, 5, 6, 7.
\end{array}$$
Then, an $\R$-linear mapping $\alpha : \gJ \to \gJ$ satisfying
$$
\begin{array}{lll}
     \alpha E = E, & \alpha E_1 = E_1, & \alpha(E_2 - E_3) = A,  
\vspace{1mm}\\
     \alpha F_1(e_i) = X_i, & \alpha F_2(e_i) = Y_i, & \alpha F_3(e_i) = Z_i, \quad i = 0, 1,\cdots,7 
\end{array}$$
belongs to $(F_4)_{E_1}$.}
\vspace{2mm}

{\bf Proof.} We have to show that the setting of the lemma is appropriate and prove that $\alpha$ satisfies
$$ 
         \alpha(X \circ Y) = \alpha X \circ \alpha Y, \quad X, Y \in \gJ. $$
For this purpose, it is sufficient to show that this holds for generators $X, Y  = E_i$, $F_j(e_k)$ of $\gJ$. We have to check $27^2 = 729$ times if honestly doing so we reduce the number of times after preparing some lemma. However we omit its proof. In details, see Yokota [40].
\vspace{3mm}

{\bf Proposition 2.7.3.} \qquad \qquad $(F_4)_{E_1}/Spin(8) \simeq S^8.$
\vspace{1mm}                                   

\noindent {\it In particular, the group $(F_4)_{E_1}$ is connected.}
\vspace{2mm}

{\bf Proof.} $S^8 = \{ X \in \gJ_{01} \, | \, (X, X) = 2 \}$ is an 8 dimentional sphere. For $\alpha \in (F_4)_{E_1}$ and $X \in S^8$, we have $\alpha X \in S^8$. Hence the group $(F_4)_{E_1}$ acts on $S^8$. This action is transitive. Indeed, for a given $A \in S^8$, by constructing $\alpha \in (F_4)_{E_1}$ of Lemma 2.7.2, we have $\alpha(E_2 - E_3) = A$, and which shows the transitivity. We determine the isotropy subgroup of $(F_4)_{E_1}$ at $E_2 - E_3 \in S^8$. Let $\alpha \in (F_4)_{E_1}$ satisfy $\alpha(E_2 - E_3) = E_2 - E_3$. Since $\alpha \in (F_4)_{E_1}$ satisfies $\alpha E_1 = E_1$ and $\alpha E = E$, it also satisfies $\alpha(E_2 + E_3) = E_2 + E_3$. Therefore we have $\alpha E_2 = E_2$ and $\alpha E_3 = E_3$, so that $\alpha \in Spin(8)$. Conversely, $\alpha \in Spin(8)$ satisfies $\alpha(E_2 - E_3) = E_2 - E_3$. Thus we have the homeomorphism $(F_4)_{E_1}/Spin(8) \simeq S^8$.
\vspace{3mm}

{\bf Remark.} If we know the dimension of the group $(F_4)_{E_1}$, without using Lemma 2.7.2, Proposition 2.7.3 can be simply proved as follows. The group $(F_4)_{E_1}$ acts on $S^8$. The isotropy subgroup of $(F_4)_{E_1}$ at $E_2 - E_3$ is $Spin(8)$ and $\dim((F_4)_{E_1}/Spin(8)) = \dim(F_4)_{E_1} - \dim Spin(8) = 36 - 28 = \dim S^8$. Therefore, we have $(F_4)_{E_1}/Spin(8) \simeq S^8$.
\vspace{3mm}

{\bf Theorem 2.7.4.} \qquad \qquad $(F_4)_{E_1} \cong Spin(9).$
\vspace{1mm}

\noindent (From now on, we identify these groups).
\vspace{2mm}

{\bf Proof.} Let $O(9) = O(\gJ_{01}) = \{ \alpha' \in \Iso_{\sR}(\gJ_{01}) \, | \, (\alpha'X, \alpha'Y) = (X, Y) \}$. Consider the restriction $\alpha' = \alpha|\gJ_{01}$ of $\alpha \in (F_4)_{E_1}$ to $\gJ_{01}$, then $\alpha' \in O(9)$. Hence we can define a homomorphism $p : (F_4)_{E_1} \to O(9)$ by $p(\alpha) = \alpha'$. Since $p$ is continuous and $(F_4)_{E_1}$ is connected (Proposition 2.7.3), the mapping $p$ induces a homomorphism
$$
       p : (F_4)_{E_1} \to SO(9).$$
We shall show $p$ is onto. Let $SO(8)=\{ \alpha' \in SO(9) \, | \, \alpha'(E_2 - E_3) = E_2 - E_3 \}$. The restriction $p'$ of $p : (F_4)_{E_1} \to SO(9)$ to 
$Spin(8) = \{ \alpha \in (F_4)_{E_1} \, | \, \alpha(E_2 - E_3) = E_2 - E_3 \}$ 
coincides with the homomorphism $p : Spin(9) \to SO(8)$ in Theorem 1.16.2. In particular, $p' : Spin(8) \to SO(8)$ is onto. Hence, from the following commutative diagram
\begin{center}
\begin{tabular}{ccccccccc}
$1$ & $\longrightarrow$ & $Spin(8)$ & $\longrightarrow$ & $(F_4)_{E_1}$ & $\longrightarrow$ & $S^8$ & $\longrightarrow$ & $*$\vspace{0.7mm}\\
$$ & $$ & $\downarrow p'$ & $$ & $\downarrow p$ & $$ & $\downarrow =$ & $$ & $$ \vspace{0.7mm}\\
$1$ & $\longrightarrow$ & $SO(8)$ & $\longrightarrow$ & $SO(9)$ & $\longrightarrow$ & $S^8$ & $\longrightarrow$ & $*$\\
\end{tabular}
\end{center}
we see that $p : (F_4)_{E_1} \to SO(9)$ is onto by the five lemma. $\Ker \,p = \{1, \sigma\},$  where $\sigma = \varphi(1, -1, -1)$ (which corresponds to $\sigma$ defined in the following Section 2.9). Indeed, let $\alpha \in \Ker \,p$, then $\alpha$ satisfies $\alpha X = X$ for all $X \in \gJ_{01}$. From $\alpha(E_2 - E_3) = E_2 - E_3$ follows that $\alpha \in Spin(8)$. Let $\alpha = (\alpha_1, \alpha_2, \alpha_3) \in Spin(8)$. Since $\alpha F_1(x) = F_1(x)$, that is, $F_1(\alpha_1x) = F_1(x)$, so $\alpha_1x = x$ for all $x \in \gC$, hence we have $\alpha_1 = 1$. From the principle of triality, we have $\alpha = (1, 1, 1) = 1$ or $\alpha = (1, -1, -1) = \sigma$. Hence we have $\Ker\, p = \{1, \sigma \}$, and so we have the isomorphism 
$$
            (F_4)_{E_1}/\{1,\sigma\} \cong SO(9).$$
Therefore $(F_4)_{E_1}$ is isomorphic to the group $Spin(9)$ as the universal covering group of $SO(9)$.
\vspace{3mm}

{\bf Theorem 2.7.5.} \qquad \qquad \quad $Spin(9)/Spin'(7) \simeq S^{15}$.
\vspace{1mm}

\noindent {\it where $Spin'(7) = \{\wti{\alpha} \in SO(8) \, | \, (\wti{\alpha} x)(\alpha y) = \wti{\alpha}(xy), x, y \in \gC \; \mbox{for some}\;\, \alpha \in SO(7) \}$.}
\vspace{2mm}

{\bf Proof.} Let $S^{15} = \{ Y \in \gJ_{23} \, | \, (Y, Y) = 2 \}$. For 
$\alpha \in Spin(9)$ and $Y \in S^{15}$, we have $\alpha Y \in S^{15}$. Hence 
$Spin(9)$ acts on $S^{15}$. We shall show that this action is transitive. Let 
$Y_0 \in S^{15}$. Choose any element $A \in \gJ_{01}$ such that $(A, A) = 2$, 
$2A \circ Y_0 = - Y_0$. Using these $A$ and $Y_0$, construct $X_i, Y_i, Z_i$ 
and $\alpha$ of Lemma 2.7.2, then $\alpha \in Spin(9)$ and satisfies $\alpha 
F_2(1) = Y_0$,  which shows the transitivity.  We determine the isotropy 
subgroup $Spin(9)_{F_2(1)}$ of $Spin(9)$ at $F_2(1) \in S^{15}$. Let $\alpha 
\in Spin(9)$ satisfies $\alpha F_2(1) = F_2(1)$. Applying $\alpha$ on $F_2(1) 
\circ F_2(1) = E_1 + E_3$ we have $F_2(1) \circ F_2(1) = E_1 + \alpha E_3$, so 
we get $\alpha E_3 = E_3$. Since $\alpha$ always satisfies $\alpha(E_2 + E_3) =
 E_2 + E_3$, we have $\alpha E_2 = E_2$, and so $\alpha \in Spin(8)$. Denote 
$\alpha = (\alpha_1, \alpha_2, \alpha_3)$, then $F_2(\alpha_2(1)) = F_2(1)$ implies $ 
\alpha_21 = 1$, so that $\alpha_2 \in SO(7)$. Hence $\alpha = (\alpha_1, \alpha_2, \alpha_3) \in Spin'(7)$. Conversely $\alpha \in Spin(7)$ satisfies $\alpha F_2(1) = F_2(1)$. Therefore $Spin(9)_{F_2(1)} = Spin'(7)$ is proved. Thus we 
have the homeomorphism $Spin(9)/Spin'(7) \simeq S^{15}$.
\vspace{1mm}

($Spin'(7)$ and $Spin(7)$ are conjugate in the group $O(7)$. Indeed, a mapping $f : Spin'(7) \to Spin(7)$,
$$
     f(\alpha') = \varepsilon\alpha'\varepsilon^{-1}, \quad \mbox{($\epsilon x = \ov{x}, x \in \gC$)} $$
gives the conjugation. 
\vspace{4mm}
In particular, $Spin'(7) \cong Spin(7)$)

{\bf 2.8. Connectedness of $F_4$}
\vspace{3mm}
   
We denote by $(F_4)_0$ the connected component of $F_4$ containing the identity $1$. 
\vspace{3mm}

{\bf Lemma 2.8.1.} {\it For $a \in \gC$, we define a mapping $\beta_1(a) : \gJ \to  \gJ$ by $\beta_1(a)X(\xi, x) = Y(\eta, y)$, where}
\begin{eqnarray*}
&&\left\{\begin{array}{l}
     \eta_1 = \xi_1 
\vspace{1mm}\\
     \eta_2 = \displaystyle{\frac{\xi_2 + \xi_3}{2}} 
            + \displaystyle{\frac{\xi_2 - \xi_3}{2}}\cos2|a|
            + \displaystyle{\frac{(a,x_1)}{|a|}}\sin2|a|
\vspace{1mm}\\
     \eta_3 = \displaystyle{\frac{\xi_2 + \xi_3}{2}}
            - \displaystyle{\frac{\xi_2 - \xi_3}{2}}\cos2|a|
            - \displaystyle{\frac{(a,x_1)}{|a|}}\sin2|a|,
\end{array}\right.\\
&&\left\{\begin{array}{l}
     y_1 = x_1 
            - \displaystyle{\frac{(\xi_2 - \xi_3)a}{2|a|}}\sin2|a|
            - \displaystyle{\frac{2(a, x_1)a}{|a|^2}}\sin^2|a|
\vspace{1mm}\\
     y_2 = x_2\cos|a| - \displaystyle{\frac{\ov{x_3a}}{|a|}}\sin|a|
\vspace{1mm}\\
     y_3 = x_3\cos|a| + \displaystyle{\frac{\ov{ax_2}}{|a|}}\sin|a|
\end{array}\right.
\end{eqnarray*}
\Big({\it if $a = 0$, then $\dfrac{\sin |a|}{|a|}$ means $1$}\Big), {\it then} $\beta_1(a) \in (F_4)_0$. 
\vspace{2mm}

{\bf Proof.} For $A_1(a)=\pmatrix{0 & 0 & 0 \cr 
                                  0 & 0 & a \cr 
                                  0 & -\ov{a} & 0}$, 
we have $\widetilde{A}_1(a) \in \gf_4$ (Proposition 2.3.6) and $\beta_1(a) = \exp \widetilde{A}_1(a)$. Hence $\beta _1(a)\in (F_4)_0$.
\vspace{3mm}

{\bf Proposition 2.8.2.} {\it Any element $X \in \gJ$ can be transformed to a diagonal form by some element $\alpha \in (F_4)_0$}: 
$$
        \alpha X = \pmatrix{\xi_1 & 0 & 0 \cr
                            0 & \xi_2 & 0 \cr
                            0 & 0 & \xi_3},
                            \quad \xi_i \in \R.$$
{\it Moreover, $\xi_1, \xi_2, \xi_3$ are uniquely determined} ({\it up to their permutation}) {\it independent of the choice of} $\alpha \in (F_4)_0$.
\vspace{2mm}

{\bf Proof.} For a given $X \in \gJ$, consider a space $\gX = \{ \alpha X \,| \, \alpha \in (F_4)_0 \}$. Since $(F_4)_0$ is compact (Theorem 2.2.5), $\gX$ is also compact. Let ${\xi_1}^2 + {\xi_2}^2 + {\xi_3}^2$ be the maximal value of all ${\eta_1}^2 + {\eta_2}^2 + {\eta_3}^2$ for $Y = Y(\eta, y) \in \gX$ and let $X_0 = X(\xi,x)$ be an element of $\gX$ which attains its maximal value. Then $X_0$ is of diagonal form. Certainly, suppose $X_0$ is not of diagonal form, for example, the $2 \times 3$ entry $x_1$ of $X_0$ is non-zero: $x_1 \neq 0$. Let $a(t) = \dfrac{x_1}{|x_1|}t$, $t > 0$, and construct $\beta_1(a(t)) \in (F_4)_0$ of 
\vspace{0.5mm}
Lemma 2.8.1. Since $|a(t)| = t$ and $\dfrac{(a(t),x_1)}{|a(t)|} = |x_1|$ for $Y(\eta(t), y(t)) = \beta_1(a(t))X_0 \in \gX$, we have
\vspace{2mm}
 
    $\eta_1(t)^2 + \eta_2(t)^2 + \eta_3(t)^2$ 
\vspace{1mm}

\quad
      $= {\xi_1}^2 + \Big(\dfrac{\xi_2 + \xi_3}{2} + \dfrac{\xi_2 - \xi_3}{2}\cos2t + |x_1|\sin2t \Big)^2$ 
\vspace{1mm}

\quad \qquad \quad
    $+ \,\Big(\dfrac{\xi_2 + \xi_3}{2} - \dfrac{\xi_2 - \xi_3}{2}\cos2t - |x_1|\sin2t \Big)^2$ 
\vspace{1mm}

\quad
    $= {\xi_1}^2 + 2\Big(\dfrac{\xi_2 + \xi_3}{2} \Big)^2 + 2\Big(\dfrac{\xi_2 - \xi_3}{2}\cos2t + |x_1|\sin2t \Big)^2$ 
\vspace{1mm}

\quad
   $= {\xi_1}^2 + 2\Big(\dfrac{\xi_2 + \xi_3}{2} \Big)^2 + 2\Big(\Big(\dfrac{\xi_2 - \xi_3}{2}\Big)^2 + |x_1|^2 \Big)\sin^2(2t + t_0) \;\mbox{(for some}\; \; t_0 \in \R)$ 
\vspace{1mm}

\quad
    $\leq {\xi_1}^2 + 2\Big(\dfrac{\xi_2 + \xi_3}{2} \Big)^2 + 2\Big(\Big(\dfrac{\xi_2 - \xi_3}{2}\Big)^2 + |x_1|^2 \Big)$ 
\vspace{1mm}

\quad       
    $= {\xi_1}^2 + {\xi_2}^2 + {\xi_3}^2 + 2|x_1|^2$

\vspace{2mm}

\noindent which is the maximal value and attains at some $t > 0$. This contradicts the maximum of ${\xi_1}^2 + {\xi_2}^2 + {\xi_3}^2$. Hence $x_1 = 0$. $x_2 = x_3 = 0$ can be similarly proved by constructing $\beta_2(a)$, $\beta_3(a) \in (F_4)_0$  analogous to $\beta_1(a)$ of Lemma 2.8.1. Hence $X_0$ is of diagonal form. We now give the proof of the latter half of the proposition. If $X \in \gJ$ is transformed to a diagonal form 
        $\alpha X = \pmatrix{\xi_1 & 0 & 0 \cr
                      0 & \xi_2 & 0 \cr 
                      0 & 0 & \xi_3}$ 
by $\alpha \in (F_4)_0$, then 
\vspace{2mm}

\qquad \qquad
  $\tr(X) = \tr(\alpha X) \; \mbox{(Lemma 2.2.1.(2))} \; = \sum_{i=1}^3\xi_i$, 
\vspace{1mm}

\qquad \qquad
   $(X, X) = (\alpha X, \alpha X) \; \mbox{(Lemma 2.2.4)} \; = \sum_{i=1}^3 \xi_i^{\ 2}$,
\vspace{1mm}

\qquad \qquad
   $\tr(X, X, X) = \tr(\alpha X, \alpha X, \alpha X)\; \mbox{(Lemma 2.2.4)}\; = \sum_{i=1}^3 \xi_i^{\ 3}$.
\vspace{2mm}

\noindent Hence, $\xi_1, \xi_2, \xi_3$ are uniquely determined (up to order) as the solutions of the following simultaneous equation:

$$
\left\{\begin{array}{l}
       \xi_1 + \xi_2 + \xi_3 = \tr(X) 
\vspace{1mm}\\
       {\xi_1}^2 + {\xi_2}^2 + {\xi_3}^2 = (X, X) 
\vspace{1mm}\\
       {\xi_1}^3 + {\xi_2}^3 + {\xi_3}^3 
= \tr(X, X, X).
\end{array}\right.$$
\vspace{2mm}

The space $\gC P_2$, called the Cayley projective plane, is defined as
$$
     \gC P_2 = \{ X \in \gJ \, | \, X^2 = X, \tr(X) = 1 \}.$$

{\bf Theorem 2.8.3.} \qquad \qquad $F_4/Spin(9) \simeq \gC P_2.$
\vspace{1mm}

\noindent {\it In particular, the group $F_4$ is connected.}
\vspace{2mm}

{\bf Proof.} For $\alpha \in F_4$ and $X \in \gC P_2$, we have $\alpha X \in \gC P_2$. Hence $F_4$ acts on $\gC P_2$. We shall prove that the group $(F_4)_0$ acts transitively on $\gC P_2$. To prove this, it is sufficient to show that any element $X \in \gC P_2$ can be transformed to $E_1 \in \gC P_2$ by some $\alpha \in (F_4)_0$. Now, $X \in \gC P_2 \subset \gJ$ can be transformed to a diagonal form by $\alpha \in (F_4)_0$:
$$
        \alpha X = \pmatrix{\xi_1 & 0 & 0 \cr 
                            0 & \xi_2 & 0 \cr
                            0 & 0 & \xi_3}, \quad \xi_i \in \R $$
(Proposition 2.8.2).  From the condition $X \circ X = X$, we have $\alpha X \circ \alpha X = \alpha X$, that is,
$$
         \pmatrix{{\xi_1}^2 & 0 & 0 \cr
                  0 & {\xi_2}^2 & 0 \cr
                  0 & 0 & {\xi_3}^2 }
        = \pmatrix{\xi_1 & 0 & 0 \cr
                   0 & \xi_2 & 0 \cr
                   0 & 0 & \xi_3}.$$
Hence ${\xi_i}^2 = \xi_i$, so that $\xi_i = 1$ or $\xi_i = 0$, $i = 1, 2, 3$. 
Next, from $\tr(\alpha X) = \tr(X) = 1$, we have $\xi_i = 1$, $\xi_{i+1} = 
\xi_{i+2} = 0$ for some $i$, that is, $\alpha X = E_i$. Moreover, $E_2$, $E_3$ are 
transformed to $E_1$ respectively by $(F_4)_0$. Certainly, if we define a mapping
 $\beta: 
\vspace{0.5mm}
\gJ \to \gJ$, 
$\beta X = TXT^{-1}$, where $T = \pmatrix{0 & 1 & 0 \cr 
                                          1 & 0 & 0 \cr
                                          0 & 0 & -1} \in SO(3)$, 
then $\beta \in (F_4)_0$ and $\beta E_2 
\vspace{0.5mm}
= E_1$. Hence $\beta\alpha X = E_1$. In the case $\alpha X = E_3$, the 
situation is similar to the above. Therefore the transitivity is proved. Since
 we have $\gC P_2 = (F_4)_0E_1$, $\gC P_2$ is connected. Now, the group $F_4$ 
acts transitively on $\gC P_2$ and the isotropy subgroup of $F_4$ at $E_1 \in 
\gC P_2$ is $Spin(9)$ (Theorem 2.7.4). Thus we have the homeomorphism 
$F_4/Spin(9)\simeq \gC P_2$. Finally, the connectedness of $F_4$ follows from 
the connectedness of $\gC P_2$ and $Spin(9)$.
\vspace{4mm}

{\bf 2.9. Involution $\sigma$ and subgroup $Spin(9)$ of $F_4$}
\vspace{3mm}

{\bf Definition.} We define an $\R$-linear transformation $\sigma$ of $\gJ$ by
$$
      \sigma\pmatrix{\xi_1 & x_3 & \ov{x}_2 \cr
                     \ov{x}_3 & \xi_2 & x_1 \cr
                     x_2 & \ov{x}_1 & \xi_3} =
            \pmatrix{\xi_1 & -x_3 & -\ov{x}_2 \cr
                     -\ov{x}_3 & \xi_2 & x_1 \cr
                     -x_2 & \ov{x}_1 & \xi_3}.$$
Then $\sigma \in F_4$ and $\sigma^2 = 1$. Observe that this $\sigma$ is the same as $\sigma = (1, -1, -1) \in Spin(8) \subset Spin(9) \subset F_4$ of Theorem 2.7.4.
\vspace{2mm}

We shall study the following subgroup $(F_4)^{\sigma}$ of $F_4$:
$$
    (F_4)^{\sigma} = \{ \alpha \in F_4 \, | \, \sigma\alpha = \alpha\sigma \}.$$
For this end, we consider $\R$-vector subspaces $\gJ_{\sigma}$ and $\gJ_{-\sigma}$ of $\gJ$, which are eigenspaces of $\sigma$, respectively by 
\begin{eqnarray*}
  \gJ_{\sigma} \!\!\! &=& \!\!\! \{ X \in \gJ \, | \, \sigma X = X \} 
       = \Big\{\pmatrix{\xi_1 & 0 & 0 \cr
                        0 & \xi_2 & x_1 \cr
                        0 & \ov{x}_1 & \xi_3} \,
                   \Big| \, \, \xi_i \in \R, x_1 \in \gC \Big\} \\
        \!\!\! &=& \!\!\! \{ X \in \gJ \, | \, E_1 \circ X = 0 \} \oplus \mbox{\es E}_1 \;\; (\mbox{where {\es E}}_1 = \{ \xi E_1 \, | \, \xi \in \R \}) \\
        \!\!\! &=& \!\!\! \gJ(2,\gC) \oplus \mbox{\es E}_1,\\
  \gJ_{-\sigma} \!\!\! &=& \!\!\! \{ X \in \gJ \, | \, \sigma X = - X \}
      = \Big\{\pmatrix{0 & x_3 & \ov{x}_2 \cr
                       \ov{x}_3 & 0 & 0 \cr
                       x_2 & 0 & 0} \, \Big| \, x_i \in \gC \Big\} \\
       \!\!\! &=& \!\!\! \{ X \in \gJ \, | \, 2E_1 \circ X = X \}.
\end{eqnarray*}
Then, $\gJ = \gJ_{\sigma} \oplus \gJ_{-\sigma}$ and $\gJ_{\sigma}$, $\gJ_{-\sigma}$ are invariant under the action of $(F_4)^{\sigma}$.
\vspace{3mm}

{\bf Theorem 2.9.1.} \qquad \quad $(F_4)^{\sigma} = (F_4)_{E_1} = Spin(9).$
\vspace{2mm}

{\bf Proof.} We shall show that for $\alpha \in (F_4)^{\sigma}$ we have $\alpha E_1 = E_1$. Let $\alpha \in (F_4)^{\sigma}$. We first show that
$$
          \alpha E_2, \alpha E_3 \in \gJ(2, \gC). $$
Certainly we have
\begin{eqnarray*}     
     \alpha E_2 \!\!\! &=& \!\!\! \alpha(-F_2(1) \times F_2(1)) = -\alpha F_2(1) \times \alpha F_2(1) \\
    \!\!\! &=& \!\!\! -(F_2(x_2) + F_3(x_3)) \times (F_2(x_2) + F_3(x_3)) \;\;(\mbox{for some } x_2, x_3 \in \gC) \\
    \!\!\! &=& \!\!\! (x_2, x_2)E_2 + (x_3, x_3)E_3 - F_1(\ov{x_2x_3})\in \gJ(2,\gC).
\end{eqnarray*}
Similarly $\alpha E_3 \in \gJ(2, \gC)$. So $\alpha E_1 = \alpha(E - E_2 - E_3) = E - \alpha E_2 - \alpha E_3$ is of the form
$$
      \alpha E_1 = E_1 + \xi_2E_2 + \xi_3E_3 + F_1(x_1). $$ 
Then we have
$$
   1 = (E_1, E_1) = (\alpha E_1, \alpha E_1) = 1 + {\xi_2}^2 + {\xi_3}^2 + 2|x_1|^2, $$  
so that $\xi_2 = \xi_3 = x_1 = 0$, therefore $\alpha E_1 = E_1$, that is, 
$\alpha \in (F_4)_{E_1}$. Conversely let $\alpha \in (F_4)_{E_1}$, that is, 
$\alpha \in F_4$ satisfies $\alpha E_1 = E_1$. Noting that $\gJ = \gJ_{\sigma} 
\oplus \gJ_{-\sigma}$ and $\gJ_{\sigma} = \{ X \in \gJ \, | \, E_1 \circ X = 0 
\} \oplus \mbox{\es E}_1$, $\gJ_{-\sigma} = \{ X \in \gJ \, | \, 2E_1 \circ X =
 X \}$ are invariant under $\alpha$, we have
\begin{eqnarray*}
    \alpha\sigma X \!\!\! &=& \!\!\! \alpha\sigma(X_1 + X_2) \qquad X_1 \in \gJ_{\sigma},\;\; X_2 \in \gJ_{-\sigma} \\
    \!\!\! &=& \!\!\! \alpha(X_1 - X_2) = \alpha X_1 - \alpha X_2 = \sigma(\alpha X_1) + \sigma(\alpha X_2) \\
    \!\!\! &=& \!\!\! \sigma\alpha(X_1 + X_2) = \sigma\alpha X, \quad X \in \gJ.
\end{eqnarray*}
Thus $\alpha\sigma = \sigma\alpha$, and so $\alpha \in (F_4)^{\sigma}$. Therefore we have shown that $(F_4)^{\sigma} = (F_4)_{E_1} \cong Spin(9)$ (Theorem 2.7.4).
\vspace{4mm}

{\bf 2.10. Center $z(F_4)$ of $F_4$}
\vspace{3mm}

{\bf Theorem 2.10.1.} {\it The center $z(F_4)$ of the group $F_4$ is trivial}: 
$$
               z(F_4) = \{ 1 \}. $$

{\bf Proof.}  Let $\alpha \in z(F_4)$. From the commutativity with $\sigma$: $\sigma\alpha = \alpha\sigma$, we have $\alpha \in Spin(9)$ (Theorem 2.9.1) and so $\alpha \in z(Spin(9))$. Since the center $z(Spin(9))$ of $Spin(9)$ 
is the group of order 2, we have
$$
       \alpha = 1 \quad \mbox{ or } \quad \alpha = \sigma. $$
However, $\sigma \notin z(F_4)$ (Theorem 2.9.1), which implies that $\alpha = 1$.
\vspace{2mm}

According to a general theory of compact Lie groups, it is known that the center of the simply connected simple Lie group of type $F_4$ is trivial. Hence $F_4$ has to be simply connected. Thus we have the following Theorem.
\vspace{3mm}

{\bf Theorem 2.10.2.} $F_4 = \{ \alpha \in \Iso_{\sR}(\gJ) \, | \ \alpha(X \circ Y) = \alpha X \circ \alpha Y \}$ {\it is a simply connected compact Lie group of type $F_4$}.
\vspace{3mm}

{\bf Remark.} If we know that the space $\gC P_2$ is simply connected, the simply connectedness of the group $F_4$ follows from $F_4/Spin(9) \simeq \gC P_2$ of Theorem 2.8.3.
\vspace{4mm}

{\bf 2.11. Involution $\gamma$ and subgroup $(Sp(1) \times Sp(3))/\Z_2$ of $F_4$}
\vspace{3mm}

We define an involutive $\R$-linear transformation $\gamma$ of $\gJ$ by
$$
       \gamma\pmatrix{\xi_1 & x_3 & \ov{x}_2 \cr
                      \ov{x}_3 & \xi_2 & x_1 \cr 
                      x_2 & \ov{x}_1 & \xi_3} 
           = \pmatrix{\xi_1 & \gamma x_3 & \ov{\gamma x_2} \cr
                      \ov{\gamma x_3} & \xi_2 & \gamma x_1 \cr 
                      \gamma x_2 & \ov{\gamma x_1} & \xi_3}. $$
This $\gamma$ is the same as $\gamma \in G_2 \subset F_4$.
\vspace{2mm}

We shall study the following subgroup $(F_4)^{\gamma}$ of $F_4$:
$$
  (F_4)^{\gamma} = \{ \alpha \in F_4 \, | \, \gamma\alpha = \alpha\gamma \}.$$

Any element $X \in \gJ$ is expressed by
$$
      X = \pmatrix{\xi_1 & x_3 & \ov{x}_2 \cr
                   \ov{x}_3 & \xi_2 & x_1 \cr 
                   x_2 & \ov{x}_1 & \xi_3} 
        = \pmatrix{\xi_1 & m_3 & \ov{m}_2 \cr
                   \ov{m}_3 & \xi_2 & m_1 \cr 
                   m_2 & \ov{m}_1 & \xi_3} 
        + \pmatrix{0 & a_3e_4 & -a_2e_4 \cr
                   -a_3e_4 & 0 & a_1e_4 \cr 
                   a_2e_4 & -a_1e_4 & 0}, $$
where $x_i = m_i + a_ie_4 \in \H \oplus \H e_4 = \gC$. We associate such element $X \in \gJ$ with the element of $\gJ(3, \H) \oplus \H^3$,
$$
        \pmatrix{\xi_1 & m_3 & \ov{m}_2 \cr
                 \ov{m}_3 & \xi_2 & m_1 \cr 
                 m_2 & \ov{m}_1 & \xi_3} + (a_1, a_2, a_3). $$
In $\gJ(3, \H) \oplus \H^3$, we define a multiplication $\times$, an inner product $(\; \; , \; \;)$ and an $\R$-linear transformation $\gamma$ respectively by
\begin{eqnarray*}
    (M + \a) \times (N + \b) \!\!\! &=& \!\!\! \Big(M \times N - \dfrac{1}{2}        (\a^*\b + \b^*\a)\Big) - \dfrac{1}{2}(\a N + \b M), \\
        (M + \a, N + \b) \!\!\! &=& \!\!\! (M, N) + 2(\a, \b), \\
         \gamma(M + \a) \!\!\! &=& \!\!\! M - \a.
\end{eqnarray*}
(In $\gJ(3, \H)$ the multiplication $M \times N$ and the inner product $(M, N)$ are analogously defined as in $\gJ$, and in $\H^3$ the inner product $(\a, \b)$ is defined naturally by $\dfrac{1}{2}(\a^*\b + \b^*\a))$. Since these operations correspond to their respective operations 
in $\gJ$, hereafter, we identify $\gJ(3, \H) \oplus \H^3$ with $\gJ$, that is,
$$
                 \gJ(3, \H) \oplus \H^3 = \gJ. $$ 

The group $F_{4,{\sH}}$ is defined to be the automorphism group of the Jordan 
algebra $\gJ_{\sH} = \gJ(3,\H)$  (in which multiplications $\circ$ and $\times$ are analogously defined as in $\gJ$): 
\begin{eqnarray*}
      F_{4,{\sH}} \!\!\! &=& \!\!\! \{ \alpha \in \Iso_{\sR}(\gJ_{\sH}) \, | \, \alpha(M \circ N) = \alpha M \circ \alpha N \} \\
               \!\!\! &=& \!\!\! \{ \alpha \in \Iso_{\sR}(\gJ_{\sH}) \, | \, \alpha(M \times N) = \alpha M \times \alpha N \}.
\end{eqnarray*}

{\bf Proposition 2.11.1.} \qquad $F_{4,{\sH}} \cong Sp(3)/\Z_2, \; \; \Z_2 = \{ E, -E \}. $
\vspace{2mm}

{\bf Proof.}  We define a mapping $\varphi : Sp(3) \to F_{4,{\sH}}$ by
$$
          \varphi(A)M = AMA^*, \quad M \in \gJ_{\sH}.$$
It is not difficult to see that $\varphi$ is well-defined and is a homomorphism. We shall show that $\varphi$ is onto. For a given $\alpha \in F_{4,{\sH}}$, we consider $\alpha E_i$, $i = 1, 2, 3$. Since $\alpha E_i$ satisfies $\alpha E_i 
\circ \alpha E_i = \alpha E_i$ and $\tr(\alpha E_i) = 1$ (that is, $\alpha E_i 
\in \H \!P_2 = \{M \in \gJ_{\sH} \, | \, M \circ M = M, \tr(M) = 1\}$ (the quaternion projective plane)), there exists $A_i \in Sp(3)$ such that $\alpha E_i =A_iE_i{A_i}^*$. Let $\a_i =  
\pmatrix{a_{i1} \cr 
         a_{i2} \cr
         a_{i3}}$ 
be the $i$-th column 
\vspace{0.5mm}
vector of $A_i$, then we have $\alpha E_i = \pmatrix{
     a_{i1}\ov{a}_{i1} & a_{i1}\ov{a}_{i2} & a_{i1}\ov{a}_{i3} 
\vspace{0.7mm}\cr 
     a_{i2}\ov{a}_{i1} & a_{i2}\ov{a}_{i2} & a_{i2}\ov{a}_{i3} 
\vspace{0.7mm}\cr 
     a_{i3}\ov{a}_{i1} & a_{i3}\ov{a}_{i2} & a_{i3}\ov{a}_{i3}}$. 
Construct  a matrix $A = (\a_1, \a_2, \a_3)$, then we have
$$
          \alpha E_i = AE_iA^*, \quad i = 1, 2, 3.$$
Since $AA^* = A(E_1 + E_2 + E_3)A^* = \alpha E_1 + \alpha E_2 + \alpha E_3 = \alpha E = E$, we have $A \in Sp(3)$. If we let $\beta = \varphi(A)^{-1}\alpha$, then $\beta \in F_{4,{\sH}}$ and satisfies
$$
           \beta E_i = E_i, \quad i = 1, 2, 3.$$
Analogously as in Theorem 2.7.1, $\beta$ induces orthogonal transformations $\beta_1,\beta_2,\beta_3$ $: \H \to \H$ satisfying
$$
\displaylines{\hfill
      (\beta_1m)(\beta_2n) = \ov{\beta_3(\ov{mn})}, \quad m, n \in \H.
\hfill\mbox{(i)}} $$
Put $p = \beta_11$ and $q = \beta_21$, then we have $|p| = |q| = 1$. Letting $m = 1$ and $n = 1$ in (i), we have $p(\beta_2n) = \ov{\beta_3\ov{n}}$ and $(\beta_1m)q = \ov{\beta_3\ov{m}}$, and hence
$$
       \beta_2m = \ov{p}(\beta_1m)q, \quad \beta_3m = \ov{(\beta_1\ov{m})}q. $$
Substitute $\zeta m = \ov{p}(\beta_1m)$ in (i), then we see that $\zeta$ satisfies
$$
       (\zeta m)(\zeta n) = \zeta(mn), \quad m, n \in \H, $$
that is, $\zeta$ is an automorphism of $\H$: $\zeta \in \mbox{Aut}(\H)$. Hence, $\zeta$ is expressed by $\zeta m = rm\ov{r}$, using $r \in Sp(1)$ (Proposition 108 of Yokota [58]). Consequently
$$
    \beta_1m = prm\ov{r}, \;\; \beta_2m = rm\ov{r}q, \;\; \beta_3m = \ov{q}rm \ov{r}\,\ov{p}, \quad m \in \H.$$
Construct a matrix $B = \pmatrix{\ov{q}r & 0 & 0 \cr
                                       0 & pr & 0 \cr
                                       0 & 0 & r}$, then $B \in Sp(3)$ and we have 
$$
          \beta M = BMB^*, \quad M \in \gJ_{\sH}, $$
that is, $\beta = \varphi(B)$. Hence
$$
      \alpha = \varphi(A)\beta = \varphi(A)\varphi(B) = \varphi(AB), \quad AB \in Sp(3).$$ 
Therefore $\varphi$ is onto. $\Ker\,\varphi = \{ E, - E \}$ can be easily obtained. Thus we have the isomorphism $Sp(3)/\Z_2 \cong F_{4,{\sH}}$.
\vspace{3mm}

{\bf Theorem 2.11.2.} \quad $(F_4)^\gamma \cong (Sp(1) \times Sp(3))/\Z_2, \; \Z_2 = \{ (1,E),(-1,-E) \}$.
\vspace{2mm}

{\bf Proof.}  We define a mapping $\varphi : Sp(1) \times Sp(3) \to (F_4)^\gamma$ by$$
        \varphi(p,A)(M + \a) = AMA^* + p\a A^*, \quad M + \a \in \gJ_{\sH} \oplus \H^3 = \gJ. $$                                
We first show that $\varphi(p,A) \in (F_4)^\gamma$. For $p \in Sp(1)$, $A \in Sp(3)$, $M, N \in \gJ_{\sH}$ and $\a, \b \in \H^3$, we have
\begin{eqnarray*}
       AMA^* \times ANA^* \!\!\!&=&\!\!\! A(M \times N)A^*, 
\vspace{1mm}\\
       (p\a A^*)^*(p\b A^*) \!\!\!&=&\!\!\! A\a^*\ov{p}p\b A^* = A(\a^*\b)A^*, 
\vspace{1mm}\\
       (p\a A^*)(ANA^*) \!\!\!&=&\!\!\! p(\a N)A^*,\;\;\mbox{etc.}
\end{eqnarray*}
From which, we see that $\varphi(p, A)$ satisfies
$$
   \varphi(p, A)((M + \a) \times (N + \b)) 
      = \varphi(p, A)(M + \a) \times \varphi(p, A)(N + \b), $$
hence $\varphi(p,A) \in F_4$. Clearly $\gamma\varphi(p,A) = \varphi(p,A)\gamma$, so that $\varphi(p,A) \in (F_4)^\gamma$. Certainly $\varphi$ is a homomorphism. We shall show that $\varphi$ is onto.  Let $\alpha \in (F_4)^\gamma$. Since 
the restriction $\alpha' = \alpha|\gJ_{\sH}$ of $\alpha$ to $\gJ_{\sH} = \{X \in \gJ \,
 | \, \gamma X = X \}$ belongs to $F_{4, {\sH}}$, there exists $A \in Sp(3)$ such 
that 
$$
             \alpha M = AMA^*, \quad M \in \gJ_{\sH} $$
(Proposition 2.11.1). Let $\beta = \varphi(1,A)^{-1}\alpha$, then $\beta|\gJ_{\sH} 
= 1$. Therefore $\beta \in G_2$. Certainly, since $\beta \in Spin(7)$, by 
letting  $\beta = (\beta_1, \beta, \kappa\beta)$, we see that $(\beta_1x)(\beta 
y) = \beta(xy)$. Now let $y = 1$, then we see $\beta_1 = \beta$. Since $\beta 
\in (G_2)^{\gamma}$ and $\beta|\H = 1$, there exists $p \in Sp(1)$ such that
$$
   \beta(m + ae_4) = m + (pa)e_4, \quad m + ae_4 \in \H \oplus \H e_4 =\gC,$$
(Theorem 1.10.1).  Hence we have
$$
    \beta(M + \a) = M + p\a = \varphi(p, E)(M + \a),$$
so that $\beta = \varphi(p, E)$. Hence we have $\alpha = \varphi(1, A)\beta = 
\varphi(1, A)\varphi(p, E) = \varphi(p, A)$. Therefore $\varphi$ is onto. 
$\Ker\,\varphi = \{ (1, E),(-1, -E) \} = \Z_2$ can be easily obtained. Thus we have the
 isomorphism $(Sp(1)\times Sp(3))/\Z_2 \cong (F_4)^\gamma$.
\vspace{2mm}

{\bf Remark.} Since $(F_4)^{\gamma}$ is connected as fixed points subgroup of 
$F_4$ by the involution $\gamma$ of the simply connected Lie group $F_4$, the fact that $\varphi : Sp(1) \times Sp(3) \to F_4$ is onto can be proved as follows. The 
elements
$$
\begin{array}{l}
    G_{ij}, \quad 0 \le i < j \le 3, 4 \le i < j \le 7, 
\vspace{1mm}\\
    \wti{A}_k(e_j), \quad 0 \le j \le 3, k = 1, 2, 3
\end{array} $$
forms an $\R$-basis of $(\gf_4)^{\gamma}$. So $\dim((\gf_4)^{\gamma}) = 6 \times 2 + 4 \times 3 = 24 = 3 + 21 = \dim(\sp(1) \oplus \sp(3))$. Hence $\varphi$ is onto.
\vspace{4mm}

{\bf 2.12. Automorphism $w$ of order 3 and subgroup $(SU(3) \times SU(3))/\Z_3$ of $F_4$}
\vspace{3mm}

We define an $\R$-linear transformation $w$ of $\gJ$ by
$$
       w\pmatrix{\xi_1 & x_3 & \ov{x}_2 \cr
                                \ov{x}_3 & \xi_2 & x_1 \cr 
                                x_2 & \ov{x}_1 & \xi_3} =
        \pmatrix{\xi_1 & \omega x_3 & \ov{\omega x_2} \cr
                                \ov{\omega x_3} & \xi_2 & \omega x_1 \cr 
                                \omega x_2 & \ov{\omega x_1} & \xi_3}. $$
This $w$ is the same as $w \in G_2 \subset F_4$.
\vspace{2mm}

We shall study the following subgroup $(F_4)^w$ of $F_4$:
$$
      (F_4)^w = \{\alpha \in F_4 \, |\, w\alpha = \alpha w \}. $$

We associate an element 
$$
          \pmatrix{\xi_1 & x_3 & \ov{x}_2 \cr
                                \ov{x}_3 & \xi_2 & x_1 \cr 
                                x_2 & \ov{x}_1 & \xi_3}, \quad 
         \xi_i \in \R, \,  x_i = a_i + \m_i \in \C \oplus \C^3 = \gC $$
of $\gJ$ with the element 
$$
          \pmatrix{\xi_1 & a_3 & \ov{a}_2 \cr
                                \ov{a}_3 & \xi_2 & a_1 \cr 
                                a_2 & \ov{a}_1 & \xi_3}
           + \Big(\m_1, \m_2, \m_3 \Big) $$
of $\gJ(3, \C) \oplus M(3, \C)$. In $M(3, \C)$, we define the exterior product $M \times N$ by
$$
    M \times N = \pmatrix{\m_2 \times \n_3 & \m_3 \times \n_1 & \m_1 \times \n_2 \cr 
                                 + & + & + \cr
                          \n_2 \times \m_3 & \n_3 \times \m_1 & \n_1 \times \m_2} \in M(3, \C), $$ 
where $M = (\m_1, \m_2, \m_3), N = (\n_1, \n_2, \n_3) \in M(3, \C)$. Then, for $P, A \in M(3, \C)$ and $M, N \in M(3, \C)$, we have
$$
    PM \times PN = {}^t\wti{P}(M \times N), \quad MA \times NA = (M \times N)\,{}^t\!\wti{A}, $$
where $\wti{P}$ and $\wti{A}$ are the adjoint matrices of $P$ and $A$, respectively. Further, in $M(3, \C)$, we define a real valued symmetric product $(M, N)$ by
$$
    (M, N) = \dfrac{1}{2}\tr(M^*N + N^*M) = \dsum_{i,j}(m_{ij}, n_{ij}), $$
where $M = \Big(m_{ij}\Big), N = \Big(n_{ij}\Big) \in M(3, \C)$.
\vspace{1mm}

In $\gJ(3, \C) \oplus M(3, \C)$, we define the multiplication $\times$, the inner product $(\; \;, \;\;)$ and the $\R$-linear transformation $w$ respectively by
\begin{eqnarray*}
 (X + M) \times (Y + N) \!\!\! &=& \!\!\! (X \times Y - \dfrac{1}{2}(M^*N + N^*M)) 
      - \dfrac{1}{2}(MY + NX + \ov{M \times N}), 
\vspace{1mm}\\
  (X + M, Y + N) \!\!\! &=& \!\!\! (X, Y) + 2(M, N),
\vspace{1mm}\\
     w(X + M) \!\!\! &=& \!\!\! X + \omega_1M, \quad \Big(\omega_1 = - \dfrac{1}{2} + \dfrac{\sqrt{3}}{2}e_1\Big),
\end{eqnarray*}
(in $\gJ(3, \C)$, multiplications $\circ$ and $\times$ are defined as in $\gJ$). Since these operations correspond to their respective operations in $\gJ$, we identify $\gJ(3, \C) \oplus M(3, \C)$ with $\gJ$, that is,
$$
        \gJ(3, \C) \oplus M(3, \C) = \gJ. $$

The group $F_{4,\sC}$ is defined to be the automorphism group of the Jordan algebra $\gJ_{\sC} = \gJ(3, \C)$\,:
\begin{eqnarray*}
     F_{4,\sC} \!\!\! &=& \!\!\! \{\alpha \in \Iso_{\sR}(\gJ_{\sC}) \, | \, \alpha(X \circ Y) = \alpha X \circ \alpha Y \}
\vspace{1mm}\\
        \!\!\! &=& \!\!\! \{\alpha \in \Iso_{\sR}(\gJ_{\sC}) \, | \, \alpha(X \times Y) = \alpha X \times \alpha Y \}.
\end{eqnarray*}

The group $\Z_2 = \{1, \epsilon \}$ acts on the group $SU(3)$ by
$$
             \epsilon A = \ov{A}, \quad A \in SU(3), $$
and the group $SU(3)\cdot\Z_2$ be the semi-direct product of groups $SU(3)$ and $\Z_2$ under this action.
\vspace{3mm}

{\bf Proposition 2.12.1.}  $F_{4,\sC} \cong (SU(3)/\Z_3)\cdot\Z_2, \; \Z_3 = \{E, \omega_1E, {\omega_1}^2E \},\, \omega_1 = - \dfrac{1}{2} + \dfrac{\sqrt{3}}{2}e_1$.
\vspace{2mm}

{\bf Proof.} We define a mapping $\varphi : SU(3)\cdot\Z_2 \to F_{4,\sC}$ by
$$
    \varphi(A, 1)X = AXA^*, \quad \varphi(A, \epsilon)X = A\ov{X}A^*, \quad X \in \gJ_{\sC}. $$
It is not difficult to see that $\varphi$ is well-defined and that it is a 
homomorphism. The proof that $\varphi$ is onto is the same as that of 
Proposition 2.11.1 (use $C, SU(3)$, instead of $\H, Sp(3)$), and we only need to modify a little the last part. Since 
$\zeta \in \mbox{Aut}(\C)$ can be either $\zeta = 1$ or $\zeta = \epsilon$ 
(where $\epsilon x = \ov{x}, x \in \C$), we have
$$
  \beta_1x = px, \beta_2x = xq, \beta_3x = \ov{q}x\ov{p} \quad \mbox{or} \quad 
  \beta_1x = p\ov{x}, \beta_2x = \ov{x}q, \beta_3x = \ov{q}\,\ov{x}\,\ov{p}. $$
Now, choose $r \in \C$ such that $r^{-3} = \ov{q}p$ and construct a matrix $B = 
   \pmatrix{\ov{q}r & 0 & 0 \cr 
            0 & pr & 0 \cr
            0 & 0 & r}$. Then, $B \in SU(3)$ and we have
$$
    \beta X = BXB^* \quad \mbox{or} \quad \beta X = B\ov{X}B^*, \quad X \in \gJ_{\sC}. $$
The remaining proof is again analogous to that of Proposition 2.11.1.
\vspace{3mm}

{\bf Theorem 2.12.2.}  $(F_4)^w \cong (SU(3) \times SU(3))/\Z_3, \;\, \Z_3 =\{(E, E), (\omega_1E, \omega_1E),$ $ ({\omega_1}^2E, {\omega_1}^2E) \}, \,\omega_1 = - \dfrac{1}{2} + \dfrac{\sqrt{3}}{2}e_1.$
\vspace{2mm}

{\bf Proof.} We define a mapping $\varphi : SU(3) \times SU(3) \to (F_4)^w$ by
$$
    \varphi(P, A)(X + M) = AXA^* + PMA^*, \quad X + M \in \gJ(3, \C) \oplus M(3, \C) = \gJ. $$
We first show that $\varphi(P, A) \in (F_4)^w$. For $P, A \in SU(3)$ and $X + M, Y + N \in \gJ(3, \C) \oplus M(3, \C)$, we have
$$
\begin{array}{l}
     AXA^* \times AYA^* = A(X \times Y)A^*, 
\vspace{1mm}\\
     (PMA^*)^*(PNA^*) = AM^*P^*PNA = A(M^*N)A,
\vspace{1mm}\\
     (PMA^*)(AYA^*) = P(MY)A^*,
\vspace{1mm}\\
     \ov{PMA^* \times PNA^*} = \ov{{}^t\wti{P}(M \times N){}^t\wti{A}^*} = 
P\ov{M \times N}A^*, \; \mbox{etc.}
\end{array} $$
Using them, we see that $\varphi(P, A)$ satisfies
$$
   \varphi(P, A)((X + M) \times (Y + N)) 
      = \varphi(P, A)(X + M) \times \varphi(P, A)(Y + N), $$
hence $\varphi(P, A) \in F_4$. Clearly $w\varphi(P, A) = \varphi(P, A)w$, so 
that $\varphi(P, A) \in (F_4)^w$. Certainly $\varphi$ is a homomorphism. We 
shall show that $\varphi$ is onto. Let $\alpha \in (F_4)^w$. Since the 
restriction $\alpha' = \alpha|\gJ_{\sC}$ of $\alpha$ to $\gJ_{\sC} = \{X \in \gJ \, | \, wX = X \}$ belongs to $F_{4,\sC}$, there exists $A \in SU(3)$ such that$$
    \alpha X = AXA^* \quad \mbox{or} \quad \alpha X = A\ov{X}A^*, \quad X \in \gJ_{\sC} $$
(Proposition 2.12.1). In the former case, let $\beta = \varphi(E, A)^{-1}\alpha$, then $\beta|\gJ_{\sC} = 1$, and so $\beta \in G_2$. Moreover $\beta \in (G_2)_{e_1} = (G_2)^w = SU(3)$ (Theorem 1.9.4), and hence, there exists $P \in SU(3)$ such that
$$
    \beta(X + M) = X + PM = \varphi(P, E)(X + M), \quad X + M \in \gJ_C \oplus M(3, \C) = \gJ. $$
Therefore we have
$$
          \alpha = \varphi(E, A)\beta = \varphi(E, A)\varphi(P, E) = \varphi(P, A). $$
In the latter case, consider the mapping $\gamma_1 : \gJ \to \gJ$ given by 
$\gamma_1(X + M) = \ov{X} + \ov{M}$, $X + M \in \gJ$ and remember that 
$\gamma_1 \in G_2 \subset F_4$. Let $\beta = \alpha^{-1}\varphi(E, A)\gamma_1$.
 Then $\beta \in F_4$ and $\beta|\gJ_{\sC} = 1$, which shows that $\beta \in 
(G_2)_{e_1} = (G_2)^w$ (Theorem 1.9.4) $\subset (F_4)^w$. Since $\alpha, 
\varphi(E, A) \in (F_4)^w$, we have $\gamma_1 \in (F_4)^w$, so that $\gamma_1 
\in (G_2)^w$ which is a contradiction (Theorem 1.9.4). Consequently the proof 
of $\varphi$ is onto is completed. The fact that $\Ker\varphi = \{(E, E), (\omega_1E, \omega_1E)$, $({\omega_1}^2E, {\omega_1}^2E) \} = \Z_3$ can be 
easily obtained. Thus we have the isomorphism $(SU(3) \times SU(3))/\Z_3 \cong
 (F_4)^w$.
\vspace{3mm}
 
{\bf Remark 1.} Since $(F_4)^w$ is connected as the fixed points subgroup of 
$F_4$ by an automorphism of order 3 of the simply connected group $F_4$ (Rasevskii [32]), that $\varphi : SU(3) \times SU(3) \to (F_4)^w$ is onto is proved as follows. The 
elements
$$
\begin{array}{l}
    G_{01}, \quad G_{23}, \quad G_{45}, \quad G_{67}, \quad G_{26} + G_{37}, 
\quad - G_{27} + G_{36},
\vspace{1mm}\\
    G_{24} + G_{35}, \quad - G_{25} + G_{34}, \quad G_{46} + G_{57}, \quad - 
G_{47} + G_{56},
\vspace{1mm}\\
    \wti{A}_1(1), \quad \wti{A}_2(1), \quad \wti{A}_3(1), \quad \wti{A}_1(e_1), \quad \wti{A}_2(e_1), \quad \wti{A}_3(e_1) 
\end{array}$$
forms an $\R$-basis of $(\gf_4)^w$. So $\dim(\gf_4)^w = 10 + 6 = 16 = 8 + 8 = \dim(\su(3) + \su(3))$. Hence $\varphi$ is onto.
\vspace{2mm}

{\bf Remark 2.} The group $F_4$ has a subgroup which is isomorphic to the group
 $(SU(3) \times SU(3))/\Z_3)\cdot\Z_2$, which is the semi-direct product of the groups $(SU(3) \times SU(3))/\Z_3$ and $\Z_2$ (the action of $\Z_2 = \{1, 
\gamma_1 \}$ on the group $(SU(3) \times SU(3))$ is $\gamma_1(P, A) = (\ov{P}, 
\ov{A}))$.
\vspace{4mm}

{\bf 2.13. Complex exceptional Lie group ${F_4}^C$}
\vspace{3mm}

{\bf Definition.} The group ${F_4}^C$ is defined to be the automorphism group of the complex Jordan algebra $\gJ^C$:
\begin{eqnarray*}
  {F_4}^C 
    \!\!\! &=& \!\!\!\{ \alpha \in \Iso_C(\gJ^C) \, | \, \alpha(X \circ Y) = \alpha X \circ \alpha Y \} \\     
   \!\!\! &=& \!\!\! \{ \alpha \in \Iso_C(\gJ^C) \, | \, \det\,(\alpha X) = \det\,X, \,(\alpha X, \alpha Y) = (X, Y) \} \\
   \!\!\! &=& \!\!\! \{ \alpha \in \Iso_C(\gJ^C) \, | \, \det\,(\alpha X) = \det\, X, \alpha E = E \} \\
   \!\!\! &=& \!\!\! \{ \alpha \in \Iso_C(\gJ^C) \, | \, \alpha(X \times Y) = \alpha X \times \alpha Y \}.
\end{eqnarray*}

We define a positive definite Hermitian inner product $\langle X, Y \rangle $ in $\gJ^C$ by
$$
         \langle X, Y \rangle = (\tau X, Y). $$
For $\alpha \in \Hom_C(\gJ^C)$, we denote the complex conjugate transpose of  $\alpha$ with respect to $\langle X, Y \rangle $ by $\alpha^{*}$: $\langle \alpha^{*}X, Y \rangle = \langle X, \alpha Y \rangle$.
\vspace{3mm}

{\bf Lemma 2.13.1.} (1) {\it For} $\alpha \in {F_4}^C$, {\it we have} $\alpha^{*} = \tau\alpha^{-1}\tau \in {F_4}^C$.
\vspace{1mm}

(2) {\it For any} $\alpha \in F_4$, {\it its complexificated mapping $\alpha^C
 : \gJ^C \to \gJ^C$ belongs to} ${F_4}^C$: $\alpha^C \in {F_4}^C$. {\it 
Identifying $\alpha$ with $\alpha^C$, we regard $F_4$ as a subgroup of 
${F_4}^C$}: {\it $F_4 \subset {F_4}^C$. For $\alpha \in {F_4}^C$, we have 
$\alpha \in F_4$ if and only if $\tau\alpha = \alpha\tau$, that is, }
$$
            F_4 = \{ \alpha \in {F_4}^C \, | \, \tau\alpha = \alpha\tau \}. $$ 

{\bf Proof.} (1) $\langle \alpha^{*}X, Y \rangle = \langle X, \alpha Y \rangle 
=(\tau X, \alpha Y) = (\alpha^{-1}\tau X, Y) = \langle \tau\alpha^{-1}\tau X, Y \rangle $ for all $X, Y \in \gJ^C$. Hence $\alpha^{*} = \tau\alpha^{-1}\tau \in
 {F_4}^C$. 
\vspace{1mm}

(2) Let $\alpha \in {F_4}^C$ satisfy $\tau\alpha = \alpha\tau$. Then, since $\tau\alpha X = \alpha\tau X = \alpha X$, we have $\alpha X \in \gJ$ for $X \in \gJ$. Hence $\alpha$ induces an $\R$-transformation $\alpha'$ of $\gJ$ and $\alpha' \in F_4$, further we have $\alpha = (\alpha')^C$.
\vspace{3mm}

{\bf Theorem 2.13.2.} {\it The polar decomposition of the Lie group} ${F_4}^C$ {\it is given by}
$$
             {F_4}^C \simeq F_4 \times \R^{52}.$$
{\it In particular}, ${F_4}^C$ {\it is a simply connected complex Lie group of type $F_4$}.
\vspace{2mm}

{\bf Proof.} Evidently ${F_4}^C$ is an algebraic subgroup of $\Iso_C(\gJ^C) = 
GL(27, C)$. If $\alpha \in {F_4}^C$, then $\alpha^{*} \in {F_4}^C$ (Lemma 2.13.1.(1)). Hence from Chevalley's lemma, we have
$$
     {F_4}^C \simeq ({F_4}^C \cap U(\gJ^C)) \times \R^d  = F_4 \times \R^d, $$
where $U(\gJ^C) = \{\alpha \in \Iso_C(\gJ^C) \, | \, \langle \alpha X, \alpha Y \rangle = \langle X, Y \rangle\}$ and $d = \dim{F_4}^C - \dim F_4 = 2 \times 52 - 52 = 52.$ Since $F_4$ is 
simply connected (Theorem 2.10.2), ${F_4}^C$ is also simply connected. The Lie algebra of the group ${F_4}^C$ is ${\gf_4}^C$, so that ${F_4}^C$ is a complex simple Lie group of type $F_4$.
\vspace{4mm}

{\bf 2.14.  Non-compact exceptional Lie groups $F_{4(4)}$ and $F_{4(-20)}$ of type 
\vspace{3mm}
$F_4$}

Consider the following two $\R$-vector spaces 
\begin{eqnarray*}
   \gJ(3, \gC') \!\!\! &=& \!\!\! \{ X \in M(3, \gC') \, | \, X^* = X \}, \\
  \gJ(1, 2, \gC) \!\!\! &=& \!\!\! \{ X \in M(3, \gC) \, | \, I_1X^*I_1 = X \}, 
\end{eqnarray*}
where $I_1 = \diag(-1, 1, 1) \in M(3, \R)$. We define the Jordan multiplication $X \circ Y$ in $\gJ(3, \gC')$ and $\gJ(1, 2, \gC)$ respectively by
$$
          X \circ Y = \frac{1}{2}(XY + YX). $$
Then we have the isomorphisms 
\begin{eqnarray*}
   \gJ(3, \gC') \!\!\! &\cong& \!\!\! \{ X \in \gJ(3, \gC^C) \, | \, \tau\gamma X = X \} = (\gJ(3, \gC^C))_{\tau\gamma}, \\
  \gJ(1, 2, \gC) \!\!\! &\cong& \!\!\! \{ X \in \gJ(3, \gC^C) \, | \, \tau\sigma X = X \} = (\gJ(3, \gC^C))_{\tau\sigma} 
\end{eqnarray*}
as Jordan algebras, therefore we identify them, respectively. We denote by $F_{4(4)}$ and $F_{4(-20)}$ respectively the automorphism groups of the Jordan algebras $\gJ(3, \gC')$ and $\gJ(1, $ $2, \gC)$:
\begin{eqnarray*}
  F_{4(4)} \!\!\! &=& \!\!\! \{\alpha \in \Iso_{\sR}(\gJ(3, \gC')) \, | \, \alpha(X \circ Y) = \alpha X \circ \alpha Y \},
\vspace{1mm}\\ 
   F_{4(-20)} \!\!\! &=& \!\!\! \{\alpha \in \Iso_{\sR}(\gJ(1, 2, \gC)) \, | \, \alpha(X \circ Y) = \alpha X \circ \alpha Y \}. 
\end{eqnarray*}
They can also be defined by
$$
       F_{4(4)} = ({F_4}^C)^{\tau\gamma}, \quad                     
       F_{4(-20)} = ({F_4}^C)^{\tau\sigma}. $$

{\bf Theorem 2.14.1.} {\it The polar decompositions of the Lie groups $F_{4(4)}$ and $F_{4(-20)}$ are respectively given by} 
\begin{eqnarray*}
   F_{4(4)}  \!\!\! &\simeq& \!\!\! (Sp(1) \times Sp(3))/\Z_2 \times \R^{28},\\
   F_{4(-20)}  \!\!\! &\simeq& \!\!\! Spin(9) \times \R^{16}.
\end{eqnarray*}

{\bf Proof.} These are facts corresponding to Theorems 2.11.2 and 2.9.1.
\vspace{3mm}

{\bf Theorem 2.14.2.} {\it The centers $z(F_{4(4)})$ and $z(F_{4(-20)})$ are trivial}\,:
$$
         z(F_{4(4)}) = \{1\}, \quad z(F_{4(-20)}) = \{1\}. $$

\newpage

\vspace{5mm}

\begin{center}
\large{\bf Exceptional Lie group $E_6$}
\end{center}
\vspace{4mm} 

{\bf 3.1. Compact exceptional Lie group $E_6$}
\vspace{3mm}

Let $\gJ^C$ be the complex exceptional Jordan algebra.\\

{\bf Definition.} We define the groups ${E_6}^C$ and $E_6$ respectively by
\begin{eqnarray*}
   {E_6}^C \!\!\! &=& \!\!\! \{ \alpha \in \Iso_C(\gJ^C) \, | \, \det\,(\alpha X) = \det\,X \} \\
           \!\!\! &=& \!\!\! \{ \alpha \in \Iso_C(\gJ^C) \, | \, (\alpha X, \alpha Y, \alpha Z) = (X, Y, Z) \},
\vspace{1mm}\\
   E_6 \!\!\! &=& \!\!\! \{ \alpha \in \Iso_C(\gJ^C) \, | \, \det\,(\alpha X) = \det\,X, \langle \alpha X, \alpha Y \rangle = \langle X, Y \rangle \} \\
       \!\!\! &=& \!\!\! \{ \alpha \in \Iso_C(\gJ^C) \, | \, (\alpha X, \alpha Y, \alpha Z) = (X, Y, Z), \langle \alpha X, \alpha Y  \rangle = \langle X, Y \rangle \} \\
        \!\!\! &=& \!\!\! \{ \alpha \in \Iso_C(\gJ^C) \, | \, \alpha X \times \alpha Y = {}^t\alpha^{-1}(X \times Y), \langle \alpha X, \alpha Y \rangle =  \langle X, Y \rangle \} \\
        \!\!\! &=& \!\!\! \{ \alpha \in \Iso_C(\gJ^C) \, | \, \alpha X \times  \alpha Y = \tau\alpha\tau(X \times Y), \langle \alpha X, \alpha Y \rangle =\ \langle X, Y \rangle \}. 
\end{eqnarray*}
If we define an involutive automorphism $\lambda$ of the group ${E_6}^C$ by
$$
                  \lambda(\alpha) = {}^t\alpha^{-1}, \quad \alpha \in {E_6}^C $$
(Lemma 2.2.3), then the definition of the group $E_6$ can be also given by
$$
        E_6 = \{ \alpha \in {E_6}^C \, | \, \tau\lambda(\alpha)\tau = \alpha \}
            = ({E_6}^C)^{\tau\lambda}. $$

{\bf Theorem 3.1.1.} $\;${\it $E_6$ is a compact Lie group.}
\vspace{2mm}

{\bf Proof.}  $E_6$ is a compact Lie group as a closed subgroup of the unitary group                                                       
$$
      U(27) = U(\gJ^C) = \{ \alpha \in \Iso_C(\gJ^C) \, | \,  \langle \alpha X, \alpha Y \rangle = \langle X, Y \rangle \}. $$

{\bf 3.2. Lie algebra $\gge_6$ of $E_6$}
\vspace{3mm}

Before investigating the Lie algebra $\gge_6$ of the group $E_6$, we will study the Lie algebra ${\gge_6}^C$ of the group ${E_6}^C$.
\vspace{3mm}

{\bf Theorem 3.2.1.} (1) {\it The Lie algebra ${\gge_6}^C$ of the Lie group ${E_6}^C$ is given by}
$$
     {\gge_6}^C = \{ \phi \in \Hom_C (\gJ^C) \, | \, (\phi X, X, X) = 0 \}.$$ 

(2) {\it Any element $\phi \in {\gge_6}^C$ is uniquely expressed by}
$$
        \phi = \delta + \wti{T}, \quad \delta \in {\gf_4}^C,T \in {\gJ_0}^C.$$
{\it In particular, the dimension of} ${\gge_6}^C$ {\it is}
$$
         \dim_C({\gge_6}^C) = 52 + 26 = 78. $$

{\bf Proof.} (1) is proved as similar way to Lemma 2.3.1.
\vspace{1mm}

(2) For $\phi \in {\gge_6}^C$, by letting $T = \phi E$, we obtain $T \in \gJ^C$, and $\tr(T) = 0$. Certainly $\tr(T) = (T, E, E) = (\phi E, E, E) = 0$. If we put $\delta = \phi - \wti{T}$, then $\delta \in {\gge_6}^C$. Moreover $\delta \in {\gf_4}^C$, because $\delta E = \phi E - \wti{T}E = T - T = 0$. Hence we have $\phi  = \delta + \wti{T}$. To prove the uniqueness of the expression, it is sufficient to show that
$$
   \delta + \wti{T} = 0,\;\; \delta \in {\gf_4}^C,T \in {\gJ_0}^C \quad \mbox{implies} \quad \delta = 0,T = 0.$$
Certainly, let apply it on $E$, then we have $T = 0$, so that $\delta = 0$. Finally, we have $\dim_C({\gge_6}^C) = 52 + 26 = 78$ from the expression above. Therefore the theorem is 
\vspace{3mm}
proved. 

{\bf Theorem 3.2.2.}  {\it The Lie bracket $[\phi_1, \phi_2]$ in ${\gge_6}^C$ is given by}
$$
      [\delta_1 + \wti{T}_1, \delta_2 + \wti{T}_2] = 
      ([\delta_1, \delta_2] + [\wti{T}_1, \wti{T}_2]) + 
      (\wti{\delta_1T_2} - \wti{\delta_2T_1}), $$
{\it where} $\phi_i = \delta_i + \wti{T}_i, \delta_i \in {\gf_4}^C, T_i \in {\gJ_0}^C.$  
\vspace{2mm}

{\bf Proof.} It is sufficient to show that $[\delta, \wti{T}] = \wti{\delta T}$ for $\delta \in {\gf_4}^C$, $T \in {\gJ_0}^C$. Now, 
\begin{eqnarray*}
   [\delta, \wti{T}]X \!\!\! &=& \!\!\! \delta(T \circ X) - T \circ \delta X
     = \delta T \circ X + T \circ \delta X - T \circ \delta X  \\
     \!\!\! &=& \!\!\! \delta T \circ X = \wti{\delta T}X, \quad X \in \gJ^C. \;\end{eqnarray*}                                     
    
We shall investigate the Lie algebra $\gge_6$ of the Lie group $E_6$.
\vspace{3mm}

{\bf Lemma 3.2.3.} {\it For} $\phi = \delta + \wti{T} \in {\gge_6}^C$, $\delta \in {\gf_4}^C$, $T \in {\gJ_0}^C$, {\it we have}
$$
           \lambda(\phi) = - ^t\phi = - ^t(\delta + \wti{T}) = \delta - \wti{T}. $$
\noindent {\it In particular,} $-{}^t\phi \in {\gge_6}^C$.
\vspace{2mm}

{\bf Proof.} $(- {}^t\phi X, Y) = - (X, \phi Y) = - (X, \delta Y + \wti{T}Y) = - (X, \delta Y) - (X, T \circ Y)$
\vspace{1mm}
      
\qquad \qquad \quad   $= (\delta X, Y) - (\wti{T}X, Y) = ((\delta - \wti{T})X, Y) \quad X, Y \in \gJ^C.  $
\vspace{1mm}

\noindent Therefore $- {}^t\phi = \delta - \wti{T} \in {\gge_6}^C$.
\vspace{3mm}

{\bf Theorem 3.2.4.} (1) {\it The Lie algebra $\gge_6$ of the Lie group $E_6$ is given by}
$$
     \gge_6 = \{ \phi \in \Hom_C(\gJ^C) \, | \, (\phi X, X, X) = 0, \langle \phi X, Y \rangle + \langle X, \phi Y \rangle = 0 \}. $$

(2) {\it Any element $\phi \in \gge_6$ is uniquely expressed by}
$$
        \phi = \delta + i\wti{T}, \quad \delta \in \gf_4, T \in \gJ_0.$$

{\bf Proof.} (1) The proof is evident (cf. Lemma 2.3.1).
\vspace{1mm}

(2) For $\phi \in {\gge_6}^C$, the condition $\phi \in \gge_6$ is equivalent to $\tau\lambda(\phi)\tau = \phi$. Now, if $\phi$ is of the form $\phi = \delta + \wti{T'}, \delta \in {\gf_4}^C, T' \in {\gJ_0}^C$ (Theorem 3.2.1.(2)), then $\tau\lambda(\phi)\tau = \phi$ is $\tau\delta\tau - \wti{\tau T'} = \delta + \wti{T'}$ (Lemma 3.2.3), that is, $\tau\delta\tau = \delta$ and $\tau T' = -T'$. Hence $\delta \in \gf_4$, and $T'$ is of the form $T' = iT, T \in \gJ_0$.
\vspace{3mm}

{\bf Proposition 3.2.5.} {\it The complexification of the Lie algebra $\gge_6$ is} ${\gge_6}^C.$
\vspace{2mm}

{\bf Proof.} For $\phi \in {\gge_6}^C$, the conjugate transposed mapping $\phi^*$ of $\phi$ with respect to the inner product $\langle X, Y \rangle$ of $\gJ^C$ is $ \phi^* = \tau\,^t\phi\tau \in {\gge_6}^C$, and for $\phi \in {\gge_6}^C$, $\phi$ belongs to $\gge_6$ if and only if $\phi^* = - \phi$.  Now, any element $\phi \in {\gge_6}^C$ can be uniquely expressed as
$$
   \phi = \frac{\phi - \phi^{*}}{2} + i\frac{\phi + \phi^{*}}{2i},
\qquad \frac{\phi - \phi^{*}}{2}, \frac{\phi + \phi^{*}}{2i} \in \gge_6. $$
Hence ${\gge_6}^C$ is the complexification of $\gge_6$. 
\vspace{4mm}

{\bf 3.3. Simplicity of ${\gge_6}^C$}
\vspace{4mm}

{\bf Theorem 3.3.1.} {\it The Lie algebra ${\gge_6}^C$ is simple and so $\gge_6$ is also simple.}
\vspace{2mm}

{\bf Proof.} We use the decomposition of ${\gge_6}^C$ of Theorem 3.2.1.(2):
$$
      {\gge_6}^C = {\gf_4}^C \oplus \wti{\gJ}_0^{\ C}. $$
Let $p : {\gge_6}^C \to {\gf_4}^C$ and $q : {\gge_6}^C \to \wti{\gJ}_0^{\ C}$ be projections of ${\gge_6}^C = {\gf_4}^C \oplus \wti{\gJ}_0^{\ C}$. Now, let $\ga$ be a non-zero ideal of ${\gge_6}^C$. Then $p(\ga)$ is an ideal of ${\gf_4}^C$. Indeed, if $\delta \in p(\ga)$, then there exists $T \in {\gJ_0}^C$ such that $\delta + \wti{T} \in \ga$. For any $\delta_1 \in {\gf_4}^C$, we have 
$$
 \ga \ni [\delta_1, \delta + \wti{T}] = [\delta_1, \delta] + \wti{\delta_1T}\;\mbox{(Theorem 3.2.2)}, $$
hence $[\delta_1, \delta] \in \ga$. 
\vspace{1mm}

We shall show that either ${\gf_4}^C \cap \ga \neq \{0\}$ or ${\wti{\gJ}_0}^{\ C} \cap \ga \neq \{ 0 \}$. Assume that ${\gf_4}^C \cap \ga = \{0\}$ and ${\wti{\gJ}_0}^{\ C} \cap \ga = \{ 0 \}$. Then the mapping $p|\ga : \ga \to {\gf_4}^C$ is injective because $\wti{\gJ}_0^{\ C} \cap \ga = \{0\}$. Since $p(\ga)$ is a non-zero ideal of ${\gf_4}^C$ and ${\gf_4}^C$ is simple, we have $p(\ga) = {\gf_4}^C$. Hence $\dim_C\ga = \dim_Cp(\ga) = \dim_C{\gf_4}^C = 52$. On the other hand, since ${\gf_4}^C \cap \ga = \{ 0 \}$, $q|\ga : \ga \to {\wti{\gJ}_0}^{\ C}$ is also injective, we have $\dim_C\ga \le \dim_C{\wti{\gJ}_0}^{\ C} = \dim_C{\gJ_0}^C = 26$. This leads to a contradiction. 
\vspace{1mm}

We now consider the following two cases.
\vspace{1mm}

(1) Case ${\gf_4}^C \cap \ga \neq \{0\}$. From the simplicity of ${\gf_4}^C$, we have ${\gf_4}^C \cap \ga = {\gf_4}^C$, hence $\ga \supset {\gf_4}^C$. On the other hand, we have
$$
         \ga \supset [\ga, {\gge_6}^C] \supset [{\gf_4}^C, {\wti{\gJ}_0}^{\ C}] =
\widetilde{{\gf_4}^C{\gJ_0}^C}\; \mbox{(Lemma 3.2.2)}\; = {\wti{\gJ}_0}^{\ C} \; \mbox{(Proposition 2.4.6.(2))}. $$ 
Consequently $\ga \supset {\gf_4}^C \oplus \wti{\gJ}_0^{\ C} = {\gge_6}^C$.
\vspace{1mm}

(2) Case $\wti{\gJ}_0^{\ C} \cap \ga \neq \{ 0 \}$. Let $\wti{A}$ ($A \in {\gJ_0}^C$) be a non-zero element of $\wti{\gJ}_0^{\ C} \cap \ga \subset \ga$. Choose $B \in {\gJ_0}^C$ such that $[\wti{A}, \wti{B}] \neq 0$ (Lemma 2.5.4), then $0 \neq [\wti{A} ,\wti{B}] \in {\gf_4}^C \cap \ga$. Hence this case is reduced to the case (1).

\noindent Therefore we have $\ga = {\gge_6}^C$.
\vspace{3mm}

{\bf Proposition 3.3.2.} (1) {\it $\gJ^C$ is a simple Jordan algebra.}
\vspace{1mm}

(2) {\it $\gJ^C$ is an ${\gge_6}^C$-irreducible $C$-module.}
\vspace{1mm}

(3) ${\gge_6}^C{\gJ}^C = \Big\{ \dsum_{i}\phi_iA_i \, | \, \phi_i \in {\gge_6}^C, A_i \in \gJ^C \Big\} = \gJ^C$.
\vspace{1mm}

{\bf Proof.} (1) Let $\ga$ be a non-zero ideal of $\gJ^C$ and $X = X(\xi, x)$ a
non-zero element of $\ga$.
\vspace{1mm}

(i) Case $\xi_1 \neq 0$. From $\xi_1E_1 = (2X \circ E_1 - X) \circ E_1 \in \ga$, we have $E_1 \in \ga$. Next, from $F_2(1) = 2E_1 \circ F_2(1) \in \ga$ and $E_1 + E_3 = F_2(1) \circ F_2(1) \in \ga$, we have $E_3 = (E_1 + E_3) - E_1 \in \ga$. Similarly $E_2 \in \ga$, so that $E = E_1 + E_2 + E_3 \in \ga$. Now, for any $X \in \gJ^C$, $X = E \circ X \in \ga$, and so $\ga = \gJ^C$. In the case $\xi_2 \neq 0$ or $\xi_3 \neq 0$, the statement is also valid.
\vspace{1mm}

(ii) Case $\xi_1 = \xi_2 = \xi_3 = 0$, $x_1 \neq 0$. We have $F_1(x_1) = 4(X \circ E_2) \circ E_3 \in \ga$. Choose $a \in \gC^C$ such that $(x_1, a) = 1$, then $F_1(x_1) \circ F_1(a) = (x_1, a)(E_2 + E_3) = E_2 + E_3 \in \ga$. Hence this case is reduced to the case (i) and so $\ga = \gJ^C$.
\vspace{1mm}

(2) Let $W$ be a non-zero ${\gge_6}^C$-invariant $C$-submodule of $\gJ^C$. For any $A \in \gJ^C$ and $X \in W$, we have
$$
      A \circ X = \Big(A - \frac{1}{3}\tr(A)E \Big)^{\sim}X + \frac{1}{3}\tr(A)X \in W.$$
Hence, $W$ is an ideal of $\gJ^C$. From the simplicity of $\gJ^C$ ((1) above), we have 
\vspace{1mm}
$W = \gJ^C$.

(3) ${\gge_6}^C\gJ^C$ is an ${\gge_6}^C$-invariant $C$-submodule of $\gJ^C$. Hence from the irreducibility of $\gJ^C$ ((2) above), we have ${\gge_6}^C\gJ^C = \gJ^C$.
\vspace{4mm}

{\bf 3.4. Element $A \vee B$ of ${\gge_6}^C$}
\vspace{3mm}

{\bf Definition.} For $A, B \in \gJ^C$, we define an element $A \vee B \in {\gge_6}^C$ by
$$
    A \vee B = [\wti{A}, \wti{B}] + \Big(A \circ B - \frac{1}{3}(A, B)E \Big)^{\sim} $$
(Proposition 2.4.1, Theorem 3.2.1).
\vspace{3mm}

{\bf Lemma 3.4.1.} {\it For $A, B \in \gJ^C$, we have}
$$
    (A \vee B)X = \frac{1}{2}(B, X)A + \frac{1}{6}(A, B)X - 2B \times (A \times X), \quad X \in \gJ^C.$$

{\bf Proof.} Consider Hamilton-Cayley formula $X \circ (X \times X) = (\det\,X)E$, $X \in \gJ^C$, that is, 
$$
    X \circ (X \circ X) - \tr(X)X \circ X + \frac{1}{2}(\tr(X)^2 - (X, X))X = 
\frac{1}{3}(X, X, X)E. $$ 
If we put $\lambda A + \mu B + \nu X$ in place of $X$, then taking the coefficient of $\lambda\mu\nu$, we have
\vspace{2mm}

\qquad \quad 
     $A \circ (B \circ X) + B \circ (X \circ A) + X \circ (A \circ B) - \tr(A)B \circ X - \tr(B)A \circ X$  
\vspace{1mm}

\qquad \qquad
   $- \tr(X)A \circ B + \dfrac{1}{2}(\tr(A)\tr(B)X + \tr(B)\tr(X)A + \tr(X)\tr(A)B)$ 
\vspace{1mm}

\qquad \qquad 
  $- \dfrac{1}{2}((A, B)X + (B, X)A + (X, A)B) = (A, B, X)E$.
\vspace{2mm}

\noindent Therefore, using the above, we have
\vspace{1mm}

$\quad  \dfrac{1}{2}(B, X)A + \dfrac{1}{6}(A, B)X - 2B \times (A \times X) + \frac{1}{3}(A, B)X$ 
\vspace{1mm}

$\quad  = \dfrac{1}{2}(B, X)A + \dfrac{1}{2}(A, B)X - 2B \circ (A \times X) +  \tr(B)A \times X + \tr(A \times X)B$
\vspace{1mm}

$\quad  \;\;-(\tr(B)\tr(A \times X) - (B, A \times X))E$
\vspace{1mm}

$\quad  = \dfrac{1}{2}(B, X)A + \dfrac{1}{2}(A, B)X - 2B \circ (A \circ X) + \tr(A)B \circ X + \tr(X)B \circ A$
\vspace{1mm}

$\quad  \;\; - \tr(A)\tr(X)B + (A, X)B + \tr(B)A \circ X - \dfrac{1}{2}\tr(B)\tr(A)X$ 
\vspace{1mm}

$\quad  \;\; -\dfrac{1}{2}\tr(B)\tr(X)A + \dfrac{1}{2}\tr(B)\tr(A)\tr(X)E - \dfrac{1}{2}\tr(B)(A, X)E$ 
\vspace{1mm}

$\quad  \;\; + \dfrac{1}{2}(\tr(A)\tr(X) - (A, X))B - \dfrac{1}{2}\tr(B)(\tr(A)\tr(X) - (A, X))E + (A, B, X)E$ 
\vspace{1mm}

$\quad  = \tr(A)B \circ X + \tr(B)X \circ A + \tr(X)A \circ B $
\vspace{1mm}

$\quad  \;\; - \dfrac{1}{2}(\tr(A)\tr(B)X + \tr(B)\tr(X)A + \tr(A)\tr(X)B)$
\vspace{1mm}

$\quad  \;\; + \dfrac{1}{2}((A, B)X + (B, X)A + (A, X)B) + (A, B, X)E - 2B \circ (A \circ X)$                                                           
\vspace{1mm}

$\quad  = A \circ (B \circ X) + B \circ (A \circ X) +X \circ (A \circ B) - 2B \circ (A \circ X)$
\vspace{1mm}

$\quad   = A \circ (B \circ X) - B \circ (A \circ X) + (A \circ B) \circ X$
\vspace{1mm}

$\quad   = [\wti{A}, \wti{B}]X + (A \circ B)^{\sim}X.$
\vspace{3mm}

{\bf Lemma 3.4.2.} (1) $E_i\vee E_j=0$, $i \neq j$.
\vspace{1mm}

(2)$\;\;E_1 \vee E_1 = \dfrac{1}{3}(2E_1 - E_2 - E_3)^{\sim},$
{\it more explicitly, we have}
$$
          (E_1 \vee E_1)\pmatrix{\xi_1 & x_3 & \ov{x}_2 \cr
                                 \ov{x}_3 & \xi_2 & x_1 \cr
                                 x_2 & \ov{x}_1 & \xi_3}
      = \dfrac{1}{6}\pmatrix{4\xi_1 & x_3 & \ov{x}_2 \cr
                             \ov{x}_3 & -2\xi_2 & -2x_1 \cr
                             x_2 & -2\ov{x}_1 & -2\xi_3}.$$

For $\phi \in \Hom_C(\gJ^C)$, we often denote $- ^t\phi$ by $\phi'$: 
$$
         (\phi'X, Y) = -(X, \phi Y), \quad X, Y \in \gJ^C. $$
\vspace{-2mm}

{\bf Lemma 3.4.3.} (1) {\it For $\phi \in {\gge_6}^C$, we have}
$$
    \phi(X \times Y) = \phi'X \times Y + X \times \phi'Y, \quad X, Y \in \gJ^C.$$
(2) {\it For $A, B \in \gJ^C$, we have}
$$
      (A \vee B)' = - B \vee A. $$

{\bf Proof.} (1) it is evident from Lemma 2.3.1.
\vspace{1mm}

(2) $(A \vee B)' = \Big([\wti{A}, \wti{B}] + \Big(A \circ B - \dfrac{1}{3}(A, B)E\Big)^{\sim}\Big)' = [\wti{A}, \wti{B}] - \Big(A \circ B - \dfrac{1}{3}(A, B)E \Big)^{\sim}$ $
\mbox{(Lemma 3.2.3)} = - [\wti{B}, \wti{A}] - \Big(A \circ B - \dfrac{1}{3}(A, B)E\Big)^{\sim} = - B \vee A. $
\vspace{3mm}

{\bf Lemma 3.4.4.} (1) {\it For} $\phi \in {\gge_6}^C$ {\it and} $A, B \in \gJ^C$, {\it we have}
$$
      [\phi, A \vee B] = \phi A \vee B + A \vee \phi'B.$$

(2) {\it Any element $\phi \in {\gge_6}^C$ is expressed by} $\phi = \dsum_{i}(A_i \vee B_i), A_i, B_i \in \gJ^C$.
\vspace{2mm}

{\bf Proof.} (1) $[\phi, A \vee B]X = \phi(A \vee B)X - (A \vee B)\phi X $
\vspace{1mm}

\quad    
   $= \phi \Big(\dfrac{1}{2}(B, X)A + \dfrac{1}{6}(A, B)X - 2B \times (A \times X) \Big) - (A \vee B)\phi X \,\mbox{(Lemma 3.4.1)}$
\vspace{1mm}

$\quad     = \dfrac{1}{2}(B, X)\phi A + \dfrac{1}{6}(A, B)\phi X -2  \phi'B \times (A \times X) - 2B \times (\phi A \times X)$
\vspace{1mm}

$\quad  \;\;  -2B \times (A \times \phi X) - \dfrac{1}{2}(B, \phi X)A - \dfrac{1}{6}(A, B)\phi X + 2B \times (A \times \phi X) \, \mbox{(Lemma 3.4.3)}$
\vspace{0.5mm}

$\quad     = \dfrac{1}{2}(B, X)\phi A + \dfrac{1}{6}(\phi A, B)X - 2B \times (\phi A \times X)$
\vspace{1mm}

$\qquad \;\;  + \dfrac{1}{2}(\phi'B, X)A + \dfrac{1}{6}(A, \phi'B) X - 2\phi'B \times (A \times X)$
\vspace{1mm}

$\quad     = (\phi A \vee B)X + (A \vee \phi'B)X, \quad X \in \gJ^C.$
\vspace{1mm}

(2) We can see from (1) that $\ga = \Big\{ \dsum_{i}(A_i \vee B_i) \, | \, A_i, B_i \in \gJ^C \Big\}$ is an ideal of ${\gge_6}^C$. From the simplicity of ${\gge_6}^C$ (Theorem 3.3.1) we have $\ga = {\gge_6}^C$.              
\vspace{4mm}

{\bf 3.5. Killing form of ${\gge_6}^C$}
\vspace{3mm}

{\bf Definition.} We define a symmetric inner product $(\phi_1, \phi_2)_6$ in ${\gge_6}^C$ by
$$ 
       (\phi_1, \phi_2)_6 = (\delta_1, \delta_2)_4 + (T_1, T_2), $$
{\it where} $\phi_i = \delta_i + \wti{T}_i$, $\delta_i \in {\gf_4}^C$, $T_i \in {\gJ_0}^C$.
\vspace{3mm}

{\bf Lemma 3.5.1.} (1) {\it The inner product $(\phi_1,\phi_2)_6$ of ${\gge_6}^C$ is ${\gge_6}^C$-adjoint invariant}\,:
\vspace{-2mm}
$$
      ([\phi, \phi_1], \phi_2)_6 + (\phi_1, [\phi, \phi_2])_6 = 0, \quad \phi, \phi_i \in {\gge_6}^C.$$

(2) {\it For} $\phi \in {\gge_6}^C, A, B \in \gJ^C$, {\it we have}
$$
        (\phi, A \vee B)_6 = (\phi A, B). $$

{\bf Proof.} (1) For $\phi = \delta + \wti{T}, \phi_i = \delta_i + \wti{T}_i, \delta, \delta_i \in {\gf_4}^C, T, T_i \in {\gJ_0}^C$, we have
\vspace{2mm}

$\qquad     ([\phi, \phi_1], \phi_2)_6$
\vspace{1mm}

$\qquad   = ([\delta + \wti{T}, \delta_1 + \wti{T}_1], \delta_2 + \wti{T}_2)_6$
\vspace{1mm}

$\qquad   =(([\delta, \delta_1] + [\wti{T}, \wti{T}_1]) + (\wti{\delta T_1} - \wti{\delta_1T}), \delta_2 + \wti{T}_2)_6 \;\; \mbox{(Theorem 3.2.2)}$
\vspace{1mm}

$\qquad   =([\delta, \delta_1], \delta_2)_4 + ([\wti{T}, \wti{T}_1], \delta_2)_4 + (\delta T_1 - \delta_1T, T_2)$
\vspace{1mm}

$\qquad   = - (\delta_1, [\delta, \delta_2])_4 + (\delta_2 T, T_1) + (\delta T_1, T_2) - (\delta_1T, T_2)\;\;\mbox{(Lemma 2.5.2)}$
\vspace{1mm}

$\qquad  = - (\delta_1, [\delta, \delta_2])_4 - (\delta_1T, T_2) - (T_1, \delta T_2) + (T_1, \delta_2T)$
\vspace{1mm}

$\qquad   = - (\delta _1 + \wti{T}_1, [\delta, \delta_2] + [\wti{T}, \wti{T}_2]) + (\wti{\delta T_2} - \wti{\delta_2 T}))_6$
\vspace{1mm}

$\qquad   = - (\delta_1 + \wti{T}_1, [\delta + \wti{T}, \delta_2 + \wti{T}_2])_6$
\vspace{1mm}

$\qquad   = - (\phi_1, [\phi, \phi_2])_6.$
\vspace{1mm}

(2) For $\phi = \delta + \wti{T}, \delta \in {\gf_4}^C, T \in {\gJ_0}^C$, we have
\begin{eqnarray*}
     (\phi, A \vee B)_6 \!\!\!&=&\!\!\! \Big(\delta + \wti{T}, [\wti{A},\wti{B}] + \Big(A \circ B - \dfrac{1}{3}(A, B)E \Big)^{\sim} \Big)_6
\vspace{1mm}\\
  \!\!\!&=&\!\!\! (\delta, [\wti{A}, \wti{B}])_4 + \Big(T, A \circ B - \dfrac{1}{3}(A, B)E\Big)
\vspace{1mm}\\
  \!\!\!&=&\!\!\! (\delta A, B) + (\wti{T}A, B) = ((\delta + \wti{T})A, B) = (\phi A, B).
\end{eqnarray*}

{\bf Lemma 3.5.2.} {\it In ${\gge_6}^C$, we have}
$$
\begin{array}{ll}
    [(E_i - E_{i+1})^{\sim}, D] = 0, \, D \in {\gd_4}^C, \!\!&\!\!
    [(E_i - E_{i+1})^{\sim}, (E_j - E_{j+1})^{\sim}] = 0,  
\vspace{1mm}\\
    {[}(E_i - E_{i+1})^{\sim}, \wti{A}_i(a){]} = - \displaystyle{\frac{1}{2}}\wti{F}_i(a), \!\!&\!\!
    [(E_i - E_{i+1})^{\sim}, \wti{F}_i(a)] = - \displaystyle{\frac{1}{2}}\wti{A}_i(a),
\vspace{1mm}\\
    {[}(E_i - E_{i+1})^{\sim}, \wti{A}_{i+1}(a){]} = - \displaystyle{\frac{1}{2}}\wti{F}_{i+1}(a), \!\!&\!\!
    [(E_i - E_{i+1})^{\sim}, \wti{F}_{i+1}(a)] = - \displaystyle{\frac{1}{2}}\wti{A}_{i+1}(a),
\vspace{1mm}\\   
 {[}(E_i - E_{i+1})^{\sim}, \wti{A}_{i+2}(a){]} = \wti{F}_{i+2}(a), \!\!&\!\!
  [(E_i - E_{i+1})^{\sim}, \wti{F}_{i+2}(a)] = \wti{A}_{i+2}(a).
\end{array} 
\vspace{2mm}$$

{\bf Theorem 3.5.3.} {\it The Killing form $B_6$ of the Lie algebra ${\gge_6}^C$ is given by}
\begin{eqnarray*}
   B_6(\phi_1, \phi_2) \!\!\! &=& \!\!\! 12(\phi_1, \phi_2)_6 
\vspace{1mm}\\
   \!\!\! &=& \!\!\! 12(\delta_1, \delta_2)_4 + 12(T_1, T_2)
\vspace{1mm}\\
   \!\!\! &=& \!\!\! \frac{4}{3}B_4(\delta_1, \delta_2) + 12(T_1, T_2) 
\vspace{1mm}\\
   \!\!\! &=& \!\!\! 4\tr(\phi_1 \phi_2),
\end{eqnarray*}
{\it where} $\phi_i = \delta_i + \wti{T}_i$, $\delta_i \in {\gf_4}^C$, $T_i \in {\gJ_0}^C$ {\it and $B_4$ is the Killing form of ${\gf_4}^C$}.
\vspace{2mm}

{\bf Proof.} Since ${\gge_6}^C$ is simple (Theorem 3.3.1), there exist $k, k' \in C$ such that
$$
     B_6(\phi_1, \phi_2) = k(\phi_1, \phi_2)_6 = k'\tr(\phi_1\phi_2). $$
To determine these $k, k'$, let $\phi = \phi_1 = \phi_2 = (E_1 - E_2)^{\sim}.$ 
Then we have
$$
      (\phi, \phi)_6 = ((E_1 - E_2)^{\sim}, (E_1 - E_2)^{\sim})_6 = (E_1 - E_2, E_1 - E_2) = 2.$$
On the other hand, $(\mbox{ad}\phi)^2$ is calculated as follows.
\vspace{1mm}

\qquad \qquad \qquad \qquad 
    $[\phi, [\phi, \wti{A}_1(e_i)]\,] = \Big[\phi,- \displaystyle{\frac{1}{2}}\wti{F}_1(e_i)\Big] = \displaystyle{\frac{1}{4}}\wti{A}_1(e_i)$,
\vspace{1mm}

\qquad \qquad \qquad \qquad 
    ${[}\phi, [\phi, \wti{A}_2(e_i) ]\, {]} = \Big[\phi,- \displaystyle{\frac{1}{2}}\wti{F}_2(e_i)\Big] = \displaystyle{\frac{1}{4}}\wti{A}_2(e_i)$,
\vspace{1mm}

\qquad \qquad \qquad \qquad 
    ${[}\phi, [\phi, \wti{A}_3(e_i)]\, {]} = [\phi, \wti{F}_3(e_i)] = \wti{A}_3(e_i)$,
\vspace{1mm}

\qquad \qquad \qquad \qquad 
    ${[}\phi, [\phi, \wti{F}_1(e_i)]\,{]} = \Big[\phi,- \displaystyle{\frac{1}{2}}\wti{A}_1(e_i)\Big] = \displaystyle{\frac{1}{4}}\wti{F}_1(e_i)$,
\vspace{1mm}

\qquad \qquad \qquad \qquad 
    ${[}\phi, [\phi, \wti{F}_2(e_i)]\,{]} = \Big[\phi,- \displaystyle{\frac{1}{2}}\wti{A}_2(e_i)\Big] = \displaystyle{\frac{1}{4}}\wti{F}_2(e_i)$,
\vspace{1mm}

\qquad \qquad \qquad \qquad 
  ${[}\phi, [\phi, \wti{F}_3(e_i)]\,{]} = [\phi, \wti{A}_3(e_i)] = \wti{F}_3(e_i)$,
\vspace{1mm}

\qquad \qquad \qquad \qquad 
  $\mbox{the others} = 0$.
\vspace{1mm}

\noindent Hence
$$
     B_6(\phi, \phi) = \tr((\mbox{ad}\phi)^2) = \Big(\displaystyle{\frac{1}{4}} \times 4 + 1 \times 2 \Big) \times 8 = 24.$$
Therefore $k = 12$. Next, we will calculate $\tr(\phi\phi)$ as follows.
$$
\begin{array}{ll}
   \phi\phi E_1 = \phi E_1 = E_1, &  
   \phi\phi\wti{F}_1(e_i) = - \displaystyle{\frac{1}{2}}\phi\wti{F}_1(e_i) = \displaystyle{\frac{1}{4}}\wti{F}_1(e_i),
\vspace{1mm}\\
   \phi\phi E_2 = - \phi E_2 = E_2, & 
   \phi\phi\wti{F}_2(e_i) = \displaystyle{\frac{1}{2}}\phi\wti{F}_2(e_i) = \displaystyle{\frac{1}{4}}\wti{F}_2(e_i),
\vspace{1mm}\\
   \phi\phi E_3 = \phi 0 = 0, &  
   \phi\phi\wti{F}_3(e_i) = \phi 0 = 0.
\end{array} $$
Hence
$$
       \tr(\phi\phi) = 1 \times 2 + \frac{1}{4} \times 8 \times 2 = 6.$$
Therefore $k' = 4$.
\vspace{3mm}

{\bf Lemma 3.5.4.} {\it The followings hold in ${\gge_6}^C$.}
\vspace{1mm}

(1)$ \; \;  A \vee (A \times A) = 0, \quad A \in \gJ^C$,
\vspace{1mm}

\hspace{5.5mm} $A \vee (B \times C) + B \vee (C \times A) + C \vee (A \times B) = 0, \quad  A, B, C \in \gJ^C$.
\vspace{1mm}

(2) {\it For} $A \in \gJ^C$, $A \neq 0$, {\it there exists} $B \in \gJ^C$ {\it such that} $A \vee B \neq 0$.
\vspace{2mm}

{\bf Proof.} (1) $(\phi,(A \times A) \vee A)_6 = (\phi (A \times A), A)$ 
(Lemma 3.5.1.(2))
\vspace{1mm}

\hspace{8mm} $= 2(\phi'A \times A, A)\; \mbox{(Lemma 3.4.3.(1))}\; = 2(\phi'A, A \times A) = - 2(A, \phi(A \times A))$ 
\vspace{1mm}

\hspace{8mm} $= -2(\phi,(A \times A) \vee A)_6$.\\
Hence $(\phi, (A \times A) \vee A)_6 = 0$ for all $\phi \in {\gge_6}^C$, so that $(A \times A) \vee A = 0$, that is, $A \vee (A \times A) = 0$ (Lemma 3.4.3.(2)). If we put $\lambda A + \mu B + \nu C$ in the place of $A$, then the result follows from the coefficient of $\lambda\mu\nu.$
\vspace{1mm}

(2) Assume that $A \vee B = 0$, that is $B \vee A = 0$ (Lemma 3.4.3.(2)) for all $B \in \gJ^C$. Then for any $\phi \in {\gge_6}^C$, $0 = (\phi, B \vee A)_6 = (\phi B, A)$ (Lemma 3.5.1.(2)). Since ${\gge_6}^C\gJ^C = \gJ^C$ (Proposition 3.3.2.(3)), we have $(\gJ^C, A) = 0$, so that $A = 0$.
\vspace{4mm}

{\bf 3.6. Roots of ${\gge_6}^C$}
\vspace{3mm}

Let
$$
        \gM^r = \{ X \in M(3,\gC) \, | \,  \mbox{all diagonal elements of }\, X \, \mbox{are real}\},$$
and let $(\gM^r)^C$ be its complexification. In $(\gM^r)^C$, we define a multiplication $X \circ Y$ by
$$
        X \circ Y = \frac{1}{2}(XY + Y^*X^*), $$
where $X^* = {}^t \ov{X}$. For $\delta = (\delta_1, \delta_2, \delta_3) \in {\gd_4}^C$ satisfying the principle of triality $(\delta_1x)y + x(\delta_2y) = \ov{\delta_3(\ov{xy})}$, $x, y \in \gC^C$, we define a $C$-linear mapping $\delta : (\gM^r)^C \to (\gM^r)^C$ by
$$
       \delta\pmatrix{\xi_1 & x_{12} & x_{13} \cr
                      x_{21} & \xi_2 & x_{23} \cr
                      x_{31} & x_{32} & \xi_3}
      = \pmatrix{0 & \delta_3x_{12} & \ov{\delta_2\ov{x}_{13}} 
\vspace{0.5mm}\cr 
                 \ov{\delta_3\ov{x}_{21}} & 0 & \delta_1x_{23} 
\vspace{0.5mm}\cr
                 \delta_2 x_{31} & \ov{\delta_1\ov{x}_{32}} & 0}.$$
Observe that this mapping $\delta$ is an extension of $\delta : \gJ^C \to \gJ^C.$ 
\vspace{3mm}

{\bf Lemma 3.6.1.} {\it For $\delta \in {\gd_4}^C$, we have}
$$
        \delta(X \circ Y) = \delta X \circ Y + X \circ \delta Y, \quad X, Y \in (\gM^r)^C. $$

{\bf Proof.} We use the following notations:
\begin{center}\begin{tabular}{lll}
        $\sigma_{23} = \delta_1$, & $\sigma_{31} = \delta_2$, & $\sigma_{12} = \delta_3$,
\vspace{1mm}\\
        $\sigma_{32} = \kappa \delta_1$, & $\sigma_{13} = \kappa\delta_2$, & $\sigma_{21} = \kappa\delta_3$.
\end{tabular}\end{center}
The $(i,i)$-element of $\delta X \circ Y + Y \circ \delta Y$ (note that this is contained in $C$) is equal to
$$
\begin{array}{l}
   R(\sum_{k=1}^3((\sigma_{ik}x_{ik})y_{ki} + \ov{y}_{ki}(\ov{\sigma_{ik}x_{ik}})) + \sum_{k=1}^3(x_{ik}(\sigma_{ki}y_{ki}) + (\ov{\sigma_{ki}y_{ki}})\ov{x}_{ik})) 
\vspace{1.5mm}\\
\qquad \qquad
     = 2R(\sum_{k}((\sigma_{ik}x_{ik})y_{ki} + x_{ik}(\sigma_{ki}y_{ki}))) 
\vspace{1.5mm}\\
\qquad \qquad
     = 2\sum_{k}((\sigma_{ik}x_{ik}, \ov{y}_{ki}) + (x_{ik}, \ov{\sigma_{ki}y_{ki}})) 
\vspace{1.5mm}\\
\qquad \qquad
     = 2\sum_{k}((x_{ik}, - \sigma_{ik}\ov{y}_{ki}) + (x_{ik}, \sigma_{ik}\ov{y}_{ki})) = 0. 
\end{array} $$
The $(i,j)$-element of $\delta X \circ Y + Y \circ \delta Y \;(i \neq j)$ is equal to
$$
\begin{array}{l}
    \sum_{k=1}^3((\sigma_{ik}x_{ik})y_{kj} + \ov{y}_{ki}(\ov{\sigma_{jk}x_{jk}})) + \sum_{k=1}^3(x_{ik}(\sigma_{kj}y_{kj}) + (\ov{\sigma_{ki}y_{ki}})\ov{x}_{jk}) 
\vspace{1.5mm}\\
    = \sum_{k}((\sigma_{ik}x_{ik})y_{kj} + y_{ik}(\sigma_{kj}x_{kj})) + \sum_{k}\ov{((\sigma_{jk}x_{jk})y_{ki} + x_{jk}(\sigma_{ki}y_{ki}))}
\end{array} $$
(If $i, j, k$ are all distinct, from the principle of triality, we obtain
$$
       (\sigma_{ik}x_{ik})y_{kj} + x_{ik}(\sigma_{kj}y_{kj}) =
        \sigma_{ij}(x_{ik}y_{kj}).$$
Even if $i = k \neq j$ or $i \neq k =j $, considering the fact that $\sigma_{ll} = 0$, we see that the formula above is also valid, since $x_{ll} \in C$)
$$
\begin{array}{l}
   = \sum_{k}\sigma_{ij}(x_{ik}y_{kj}) + \sum_{k}\ov{\sigma_{ji}(x_{jk}y_{ki})} \qquad 
\vspace{1.5mm}\\
   = \sum_{k}\sigma_{ij}(x_{ik}y_{kj}) + \sum_{k}\sigma_{ij}(\ov{x_{jk}y_{ki}})    = (i,j)\mbox{-element of}\; \,  \delta(X \circ Y).\qquad
\end{array}$$

For $T \in M(3,\gC^C)$, we define a $C$-linear mapping $\wti{T} : \gJ^C \to \gJ^C$ by$$
     \wti{T}X = \frac{1}{2}(TX + XT^*), \quad \mbox{where} \quad T^* = {}^t\ov{T}.$$

{\bf Proposition 3.6.2.} {\it For $T \in M(3, \gC^C), \tr(T) = 0$, we have $\wti{T} \in {\gge_6}^C$}.
\vspace{2mm}

{\bf Proof.} We decompose $T = T_1 + T_2$, $T_1 = \displaystyle{\frac{T + T^*}{2}}$, $T_2 = \displaystyle{\frac{T - T^*}{2}}$. Then $\wti{T}_1 \in {\gge_6}^C$ (Lemma 2.4.1.(1)) and $\wti{T}_2 \in {\gf_4}^C$ (Proposition 2.3.6) $\subset {\gge_6}^C$. Therefore $\wti{T} = \wti{T}_1 + \wti{T}_2 \in {\gge_6}^C$.
\vspace{3mm}

{\bf Lemma 3.6.3.} (1) {\it For} $\delta \in {\gd_4}^C$ {\it and} $R \in (\gM^r)^C$, $\tr(R)=0$, {\it we have}
$$
          [\delta, \wti{R} ]=\wti{\delta R}.$$

(2) {\it For} $H \in M(3 ,C)$, $\tr(H) = 0$ and $T \in M(3, \gC^C)$, $\tr(T) = 0$, {\it we have}
$$
     [\wti{H}, \wti{T}] = \frac{1}{2}[H, T]^{\sim}.$$

{\bf Proof.} (1) $(\wti{\delta R})X = \delta R \circ X = \delta(R \circ X) - R \circ \delta X$ (Lemma 3.6.1) $ = \delta(\wti{R}X) - \wti{R}(\delta X) = [\delta, \wti{R}]X$, $X \in \gJ^C$. Hence, $\wti{\delta R} = [\delta, \wti{R}]$.
\vspace{1mm}

(2) \quad $[\wti{H}, \wti{T}]X = \wti{H}\wti{T}X - \wti{T}\wti{H}X $
\vspace{1mm}

   $= \displaystyle{\frac{1}{2}}(\wti{H}(TX + XT^*) - \wti{T}(HX + XH^*))$
\vspace{1mm}

    $= \displaystyle{\frac{1}{4}}(HTX + XT^*H^* + HXT^* + TXH^* - THX - XH^*T^* - TXH^* - HXT^*)$
\vspace{1mm}

\noindent (since $H \in M(3, C)$, products of matrices above are associative)
\vspace{1mm}

   $ = \displaystyle{\frac{1}{4}}([H, T]X + X[H, T]^*) =  \displaystyle{\frac{1}{2}}[H, T]^{\sim}X, \quad X \in \gJ^C$.
\vspace{3mm}

{\bf Theorem 3.6.4.} {\it The rank of the Lie algebra ${\gge_6}^C$ is} 6. {\it The roots of ${\gge_6}^C$ relative to some Cartan subalgebra are given by}
$$
\begin{array}{cl}
     \pm(\lambda_k - \lambda_l), \quad \pm(\lambda_k + \lambda_l), & 0 \le k < l \le 3, 
\vspace{1mm}\\
     \pm \lambda_k \pm \displaystyle{\frac{1}{2}}(\mu_2 - \mu_3), & 0 \le k \le 3,
\end{array} $$
\vspace{-2mm}
$$
\begin{array}{l}
   \pm \displaystyle{\frac{1}{2}}(- \lambda_0 - \lambda_1 + \lambda_2 - \lambda_3) \pm \displaystyle{\frac{1}{2}}(\mu_3 - \mu_1), 
\vspace{1mm}\\
   \pm \displaystyle{\frac{1}{2}}(\;\;\; \lambda_0 + \lambda_1 + \lambda_2 - \lambda_3) \pm \displaystyle{\frac{1}{2}}(\mu_3 - \mu_1),
\vspace{1mm}\\
   \pm \displaystyle{\frac{1}{2}}(- \lambda_0 + \lambda_1 + \lambda_2 + \lambda_3) \pm \displaystyle{\frac{1}{2}}(\mu_3 - \mu_1), 
\vspace{1mm}\\
   \pm \displaystyle{\frac{1}{2}}(\;\;\; \lambda_0 - \lambda_1 + \lambda_2 + \lambda_3) \pm \displaystyle{\frac{1}{2}}(\mu_3 - \mu_1), 
\end{array}$$
\vspace{-1mm}
$$
\begin{array}{l}
   \pm \displaystyle{\frac{1}{2}}(\;\;\; \lambda_0 - \lambda_1 + \lambda_2 - \lambda_3) \pm \displaystyle{\frac{1}{2}}( \mu_1 - \mu_2), 
\vspace{1mm}\\
   \pm \displaystyle{\frac{1}{2}}(- \lambda_0 + \lambda_1 + \lambda_2 - \lambda_3) \pm \displaystyle{\frac{1}{2}}(\mu_1 - \mu_2),
\vspace{1mm}\\
   \pm \displaystyle{\frac{1}{2}}(\;\;\; \lambda_0 + \lambda_1 + \lambda_2 + \lambda_3) \pm \displaystyle{\frac{1}{2}}(\mu_1 - \mu_2), 
\vspace{1mm}\\
   \pm \displaystyle{\frac{1}{2}}(- \lambda_0 - \lambda_1 + \lambda_2 + \lambda_3) \pm \displaystyle{\frac{1}{2}}(\mu_1 - \mu_2),
\end{array} $$
{\it with} $\mu_1 + \mu_2 + \mu_3 = 0.$ 
\vspace{2mm}

{\bf Proof.} We use the decomposition of Theorem 3.3.1.(2):
$$
          {\gge_6}^C = {\gf_4}^C \oplus {\wti{\gJ}_0}^{\ C}. $$
Let
$$
      \gh = \Bigg\{ h = h_{\delta} + \wti{H} \in {\gge_6}^C \, \left| \; 
\begin{array}{l}
     h_\delta = \displaystyle{\sum_{k=0}^3}\lambda_kH_k = - \displaystyle{\sum_{k=0}^3}\lambda_kiG_{k 4+k}, \lambda_k \in C \\
     H  = \displaystyle{\sum_{j=1}^3}\mu_jE_j,\mu_j \in C,\mu_1 + \mu_2 + \mu_3 = 0
\end{array}\right.
\Bigg\}, $$
then $\gh$ is an abelian subalgebra of ${\gge_6}^C$ (it will be a Cartan subalgebra of ${\gge_6}^C$). That $\gh$ is abelian is clear from
$$
\begin{array}{l}
      [h_\delta, h_{\delta'}] = 0,
\vspace{1mm}\\
      {[}h_{\delta}, \wti{H}'{]} = \wti{h_\delta H'} \; (\mbox{Lemma 3.6.3.(1)})\; = \wti{0} = 0, 
\vspace{1mm}\\
      {[}\wti{H}, \wti{H}'{]} = \displaystyle{\frac{1}{2}}[H, H']^{\sim} \; 
(\mbox{Lemma 3.6.3.(2)})\; = 0.
\end{array} $$

I $\;$ The roots $\pm \lambda_k \pm \lambda_l$ of ${\gd_4}^C (\subset {\gf_4}^C \subset {\gge_6}^C$) are also roots of ${\gge_6}^C$. Indeed, let $\alpha$ be a root of ${\gd_4}^C$ and $S \in {\delta_4}^C \subset {\gge_6}^C$ be an associated root vector. Then
\begin{eqnarray*}
       [h ,S] \!\!\! &=& \!\!\! [h_\delta + \wti{H}, S] = [h_\delta, S] - [S, \wti{H}]
\vspace{1mm}\\
      \!\!\! &=& \!\!\! \alpha(h_\delta)S - \wti{SH}\;\; \mbox{(Lemma 3.6.3.(1))}\vspace{1mm}\\
      \!\!\! &=& \!\!\! \alpha(h_\delta)S = (\pm \lambda_k \pm \lambda_l)S.
\end{eqnarray*}
Hence $\pm \lambda_k \pm \lambda_l$ are roots of ${\gge_6}^C$.
\vspace{1mm}

II $\;$ We denote $aE_{kl} \in M(3, \gC^C)$ by $F_{kl}(a)$: $F_{kl}(a) = aE_{kl}, a \in \gC^C, k \neq l$. For example, $F_{23}(a) = \pmatrix{0 & 0 & 0 \cr 
                   0 & 0 & a \cr 
                   0 & 0 & 0}$. 
Then, we have
\begin{eqnarray*}
    [h, \wti{F}_{23}(a)] \!\!\! &=& \!\!\! [h_\delta, \wti{F}_{23}(a)] + [\wti{H}, \wti{F}_{23}(a)] \\
    \!\!\! &=& \!\!\! (h_{\delta}F_{23}(a))^{\sim} + \displaystyle{\frac{1}{2}}[H, F_{23}(a)]^{\sim} \;\;\mbox{(Lemma 3.6.3)} \\
    \!\!\! &=& \!\!\! \wti{F}_{23}(h_\delta a) + \displaystyle{\frac{1}{2}}[H, F_{23}(a)]^{\sim} \\
    \!\!\! &=& \!\!\! \Big(\lambda_k + \displaystyle{\frac{1}{2}}(\mu_2 - \mu_3)\Big)\wti{F}_{23}(a),
\end{eqnarray*}
where $a = e_k + ie_{4+k}$ (since $h_{\delta}a = h_{\delta}(e_k + ie_{4+k}) = \lambda_k(e_k + ie_{4+k}) = \lambda_ka$, 
\vspace{0.5mm}
and $[H, F_{23}(a)] = HF_{23}(a) - F_{23}(a)H = (\mu_2 - \mu_3)F_{23}(a)$). Hence $\lambda_k + \displaystyle{\frac{1}{2}}(\mu_2 - \mu_3)$ is a root of ${\gge_6}^C$ and $\wti{F}_{23}(e_k + ie_{4+k})$ is its root vector. Similarly $- \lambda_k + \displaystyle{\frac{1}{2}}(\mu_2 - \mu_3)$ is a root of ${\gge_6}^C$ and $\wti{F}_{23}(e_k - ie_{4+k})$ is its root vector. Again, from
$$
    [h, \wti{F}_{32}(a)] = \Big(- \lambda_k + \displaystyle{\frac{1}{2}}(\mu_3 - \mu_2)\Big)\wti{F}_{32}(a), \quad a = e_k + i e_{4+k}, $$
we see that $ - \lambda_k + \displaystyle{\frac{1}{2}}(\mu_3 - \mu_2)$ is a root of ${\gge_6}^C$. Similarly $\lambda_k + \displaystyle{\frac{1}{2}}(\mu_3 - \mu_2)$ is a root of ${\gge_6}^C$ and $\wti{F}_{32}(e_k-ie_{4+k})$ is its root 
vector. Next, using the relation 
$$
\begin{array}{l}
    [h, \wti{F}_{31}(a)] = (h_\delta F_{31}(a))^{\sim} + \displaystyle{\frac{1}{2}}[\wti{H}, \wti{F}_{31}(a)] = \wti{F}_{31}((\nu h_{\delta})a) + \displaystyle{\frac{1}{2}}(\mu_3 - \mu_1)\wti{F}_{31}(a),
\vspace{1mm}\\
    {[} h, \wti{F}_{12}(a){]} = (h_\delta F_{12}(a))^{\sim} + \displaystyle{\frac{1}{2}}[\wti{H}, \wti{F}_{12}(a)] = \wti{F}_{12}((\kappa\pi h_{\delta})a) + \displaystyle{\frac{1}{2}}(\mu_1 - \mu_2)\wti{F}_{12}(a),
\end{array}$$
where
\begin{eqnarray*}
\nu h_\delta \!\!\! &=& \!\!\! 
\displaystyle{\frac{1}{2}}(- \lambda_0 - \lambda_1 + \lambda_2 - \lambda_3)H_0  + \displaystyle{\frac{1}{2}}(\lambda_0 + \lambda_1 + \lambda_2 - \lambda_3)H_1
\\
\!\!\! &+& \!\!\! \displaystyle{\frac{1}{2}}(-\lambda_0 + \lambda_1 + \lambda_2 + \lambda_3)H_2 + \displaystyle{\frac{1}{2}}(\lambda_0 - \lambda_1 + \lambda_2 + \lambda_3)H_3,
\vspace{1mm}\\
\kappa\pi h_\delta \!\!\! &=& \!\!\! 
\displaystyle{\frac{1}{2}}(- \lambda_0 + \lambda_1 - \lambda_2 + \lambda_3)H_0 
+ \displaystyle{\frac{1}{2}}(- \lambda_0 + \lambda_1 + \lambda_2 - \lambda_3)H_1\\
\!\!\! &+& \!\!\! \displaystyle{\frac{1}{2}}(\lambda_0 + \lambda_1 + \lambda_2 + \lambda_3)H_2 
+ \displaystyle{\frac{1}{2}}(-\lambda_0 - \lambda_1 + \lambda_2 + \lambda_3)H_3
\end{eqnarray*}
etc., we can obtain the remainders of roots.
\vspace{3mm}

{\bf Theorem 3.6.5.} {\it In the root system of Theorem} 3.6.4,
$$
\begin{array}{l}
    \alpha_1 = \lambda_0 - \lambda_1, \quad \alpha_2 = \lambda_1 - \lambda_2, \quad \alpha_3 = \lambda_2 - \lambda_3,
\vspace{1mm}\\
    \alpha_4 = \lambda_3 + \displaystyle{\frac{1}{2}}(\mu_2 - \mu_3), 
\vspace{1mm}\\
    \alpha_5 = \displaystyle{\frac{1}{2}}(- \lambda_0 - \lambda_1 - \lambda_2 + \lambda_3) + \displaystyle{\frac{1}{2}}(\mu_3 - \mu_1),
\vspace{1mm}\\
    \alpha_6 = \displaystyle{\frac{1}{2}}(\lambda_0 + \lambda_1 + \lambda_2 + \lambda_3) + \displaystyle{\frac{1}{2}}(\mu_1 - \mu_2)
\end{array} $$
{\it is a fundamental root system of the Lie algebra ${\gge_6}^C$ and
$$
  \mu = \alpha_1 + 2\alpha_2 + 3\alpha_3 + 2\alpha_4 + 2\alpha_5 + \alpha_6$$
is the highest root. The Dynkin diagram and the extended Dynkin diagram of ${\gge_6}^C$ are respectively given by}
\vspace{-2mm}

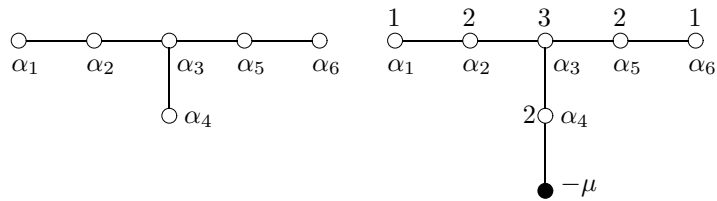
\setlength{\unitlength}{1mm}
\begin{picture}(150,20)
\put(10,10){\circle{2}} \put(9,6){$\alpha_1$}
\put(11,10){\line(1,0){8}}
\put(20,10){\circle{2}} \put(19,6){$\alpha_2$}
\put(21,10){\line(1,0){8}}
\put(30,10){\circle{2}} \put(31,6){$\alpha_3$}
\put(30,9){\line(0,-1){8}}
\put(30,0){\circle{2}} \put(32,-1){$\alpha_4$}
\put(31,10){\line(1,0){8}}
\put(40,10){\circle{2}} \put(39,6){$\alpha_5$}
\put(41,10){\line(1,0){8}}
\put(50,10){\circle{2}} \put(49,6){$\alpha_6$}

\put(60,10){\circle{2}} \put(59,6){$\alpha_1$} \put(59,12){$1$}
\put(61,10){\line(1,0){8}}
\put(70,10){\circle{2}} \put(69,6){$\alpha_2$} \put(69,12){$2$}
\put(71,10){\line(1,0){8}}
\put(80,10){\circle{2}} \put(81,6){$\alpha_3$} \put(79,12){$3$}
\put(80,9){\line(0,-1){8}}
\put(80,0){\circle{2}} \put(82,-1){$\alpha_4$} \put(77,-1){$2$}
\put(80,-1){\line(0,-10){8}}
\put(80,-10){\circle*{2}} \put(82,-10){$-\mu$} 
\put(81,10){\line(1,0){8}}
\put(90,10){\circle{2}} \put(89,6){$\alpha_5$} \put(89,12){$2$}
\put(91,10){\line(1,0){8}}
\put(100,10){\circle{2}} \put(99,6){$\alpha_6$} \put(99,12){$1$}
\end{picture}
\vspace{12mm}

{\bf Proof.} In the following, the notation $n_1n_2 \cdots n_6$ denotes the root  $n_1\alpha_1 + n_2\alpha_2 + \cdots + n_6\alpha_6$. Now, all positive roots of ${\gge_6}^C$ are represented by 
$$
\begin{array}{lllllllll}
      \lambda_0 - \lambda_1 \!\!\! &=& \!\!\! 1 & 0 & 0 & 0 & 0 & 0 
\vspace{1mm}\\
      \lambda_0 - \lambda_2 \!\!\! &=& \!\!\! 1 & 1 & 0 & 0 & 0 & 0 
\vspace{1mm}\\
      \lambda_0 - \lambda_3 \!\!\! &=& \!\!\! 1 & 1 & 1 & 0 & 0 & 0 
\vspace{1mm}\\
      \lambda_1 - \lambda_2 \!\!\! &=& \!\!\! 0 & 1 & 0 & 0 & 0 & 0 
\vspace{1mm}\\
      \lambda_1 - \lambda_3 \!\!\! &=& \!\!\! 0 & 1 & 1 & 0 & 0 & 0 
\vspace{1mm}\\
      \lambda_2 - \lambda_3 \!\!\! &=& \!\!\! 0 & 0 & 1 & 0 & 0 & 0
\end{array}
\quad
\begin{array}{lllllllll}
      \lambda_0 + \lambda_1 \!\!\! &=& \!\!\! 1 & 2 & 2 & 1 & 1 & 1 
\vspace{1mm}\\
      \lambda_0 + \lambda_2 \!\!\! &=& \!\!\! 1 & 1 & 2 & 1 & 1 & 1 
\vspace{1mm}\\
      \lambda_0 + \lambda_3 \!\!\! &=& \!\!\! 1 & 1 & 1 & 1 & 1 & 1 
\vspace{1mm}\\
      \lambda_1 + \lambda_2 \!\!\! &=& \!\!\! 0 & 1 & 2 & 1 & 1 & 1 
\vspace{1mm}\\
      \lambda_1 + \lambda_3 \!\!\! &=& \!\!\! 0 & 1 & 1 & 1 & 1 & 1 
\vspace{1mm}\\
      \lambda_2 + \lambda_3 \!\!\! &=& \!\!\! 0 & 0 & 1 & 1 & 1 & 1
\end{array} $$
$$
\begin{array}{llllllllll}
\lambda_0 + \dfrac{1}{2}(\mu_2 - \mu_3) \!\!\! &=& \!\!\! 1 & 1 & 1 & 1 & 0 & 0 \vspace{1mm}\\
\lambda_1 + \dfrac{1}{2}(\mu_2 - \mu_3) \!\!\! &=& \!\!\! 0 & 1 & 1 & 1 & 0 & 0 \vspace{1mm}\\
\lambda_2 + \dfrac{1}{2}(\mu_2 - \mu_3) \!\!\! &=& \!\!\! 0 & 0 & 1 & 1 & 0 & 0 \vspace{1mm}\\
\lambda_3 + \dfrac{1}{2}(\mu_2 - \mu_3) \!\!\! &=& \!\!\! 0 & 0 & 0 & 1 & 0 & 0 \end{array}$$
\vspace{-1mm}
$$
\begin{array}{llllllllllll}
\lambda_0 - \dfrac{1}{2}(\mu_2 - \mu_3) \!\!\! &=& \!\!\! 1 & 1 & 1 & 0 & 1 & 1 \vspace{1mm}\\
\lambda_1 - \dfrac{1}{2}(\mu_2 - \mu_3) \!\!\! &=& \!\!\! 0 & 1 & 1 & 0 & 1 & 1 \vspace{1mm}\\
\lambda_2 - \dfrac{1}{2}(\mu_2 - \mu_3) \!\!\! &=& \!\!\! 0 & 0 & 1 & 0 & 1 & 1 \vspace{1mm}\\
\lambda_3 - \dfrac{1}{2}(\mu_2 - \mu_3) \!\!\! &=& \!\!\! 0 & 0 & 0 & 0 & 1 & 1
\end{array}$$
$$
\begin{array}{llllllllll}
     \dfrac{1}{2}(- \lambda_0 - \lambda_1 + \lambda_2 - \lambda_3) + \dfrac{1}{2}(\mu_3 - \mu_1) \!\!\! &=& \!\!\! 0 & 0 & 1 & 0 & 1 & 0 
\vspace{1mm}\\
     \dfrac{1}{2}(\;\;\; \lambda_0 + \lambda_1 + \lambda_2 - \lambda_3) + \dfrac{1}{2}(\mu_3 - \mu_1) \!\!\! &=& \!\!\! 1 & 2 & 3 & 1 & 2 & 1 
\vspace{1mm}\\
     \dfrac{1}{2}(- \lambda_0 + \lambda_1 + \lambda_2 + \lambda_3) + \dfrac{1}{2}(\mu_3 - \mu_1) \!\!\! &=& \!\!\! 0 & 1 & 2 & 1 & 2 & 1 
\vspace{1mm}\\
     \dfrac{1}{2}(\;\;\; \lambda_0 - \lambda_1 + \lambda_2 + \lambda_3) + \dfrac{1}{2}(\mu_3 - \mu_1) \!\!\! &=& \!\!\! 1 & 1 & 2 & 1 & 2 & 1 
\end{array}$$
\vspace{-1mm}
$$
\begin{array}{llllllllllll}
     \dfrac{1}{2}(\;\;\; \lambda_0 + \lambda_1 - \lambda_2 + \lambda_3) + \dfrac{1}{2}(\mu_3 - \mu_1) \!\!\! &=& \!\!\! 1 & 2 & 2 & 1 & 2 & 1 
\vspace{1mm}\\
     \dfrac{1}{2}(- \lambda_0 - \lambda_1 - \lambda_2 + \lambda_3) + \dfrac{1}{2}(\mu_3 - \mu_1) \!\!\! &=& \!\!\! 0 & 0 & 0 & 0 & 1 & 0 
\vspace{1mm}\\
     \dfrac{1}{2}(\;\;\; \lambda_0 - \lambda_1 - \lambda_2 - \lambda_3) + \dfrac{1}{2}(\mu_3 - \mu_1) \!\!\! &=& \!\!\! 1 & 1 & 1 & 0 & 1 & 0 
\vspace{1mm}\\
     \dfrac{1}{2}(- \lambda_0 + \lambda_1 - \lambda_2 - \lambda_3) + \dfrac{1}{2}(\mu_3 - \mu_1) \!\!\! &=& \!\!\! 0 & 1 & 1 & 0 & 1 & 0
\end{array} $$
$$
\begin{array}{llllllll}
   \dfrac{1}{2}(\;\;\; \lambda_0 - \lambda_1 + \lambda_2 - \lambda_3) - \dfrac{1}{2}(\mu_1 - \mu_2) \!\!\! &=& \!\!\! 1 & 1 & 2 & 1 & 1 & 0 
\vspace{1mm}\\
   \dfrac{1}{2}(\;\;\; \lambda_0 - \lambda_1 - \lambda_2 + \lambda_3) - \dfrac{1}{2}(\mu_1-\mu_2) \!\!\! &=& \!\!\! 1 & 1 & 1 & 1 & 1 & 0 
\vspace{1mm}\\
   \dfrac{1}{2}(\;\;\; \lambda_0 + \lambda_1 + \lambda_2 + \lambda_3) + \dfrac{1}{2}(\mu_1 - \mu_2) \!\!\! &=& \!\!\! 0 & 0 & 0 & 0 & 0 & 1 
\end{array}$$
$$
\begin{array}{llllllllllll}
   \dfrac{1}{2}(\;\;\; \lambda_0 + \lambda_1 - \lambda_2 - \lambda_3) - \frac{1}{2}(\mu_1 - \mu_2) \!\!\! &=& \!\!\! 1 & 2 & 2 & 1 & 1 & 0 
\vspace{1mm}\\
   \dfrac{1}{2}(- \lambda_0 + \lambda_1 - \lambda_2 + \lambda_3) - \dfrac{1}{2}(\mu_1-\mu_2) \!\!\! &=& \!\!\! 0 & 1 & 1 & 1 & 1 & 0 
\vspace{1mm}\\
   \dfrac{1}{2}(- \lambda_0 + \lambda_1 + \lambda_2 - \lambda_3) - \dfrac{1}{2}(\mu_1 - \mu_2) \!\!\! &=& \!\!\! 0 & 1 & 2 & 1 & 1 & 0 
\vspace{1mm}\\
   \dfrac{1}{2}(\;\;\; \lambda_0 + \lambda_1 + \lambda_2 + \lambda_3) - \dfrac{1}{2}(\mu_1 - \mu_2) \!\!\! &=& \!\!\! 1 & 2 & 3 & 2 & 2 & 1 
\vspace{1mm}\\
   \dfrac{1}{2}(- \lambda_0 - \lambda_1 + \lambda_2 + \lambda_3) - \dfrac{1}{2}(\mu_1 - \mu_2) \!\!\! &=& \!\!\! 0 & 0 & 1 & 1 & 1 & 0.
\end{array} $$
Hence ${\mit\Pi} = \{\alpha_1, \alpha_2, \cdots, \alpha_6 \}$ is a fundamental root system of ${\gge_6}^C$. The real part $\gh_{\sR}$ of $\gh$ is
$$
    \gh_{\sR} = \Big\{ \dsum_{k=0}^3\lambda_kH_k + \Big(\dsum_{j=1}^3\mu_jE_j\Big)^\sim \, | \, \lambda_k, \mu_j \in \R, \mu_1 + \mu_2 + \mu_3 = 0 \Big\}. $$
The Killing form $B_6$ of ${\gge_6}^C$ on $\gh_{\sR}$ is given by 
$$  
        B_6(h, h') = 12\Big(2\dsum_{k=0}^3\lambda_k{\lambda_k}'+
               \dsum_{j=1}^3\mu_j{\mu_j}' \Big)$$
for $h = \dsum_{k=0}^3 \lambda_kH_k + \Big(\dsum_{j=1}^3\mu_jE_j\Big)^\sim$, $h' = \dsum_{k=0}^3{\lambda_k}'H_k + \Big(\dsum_{j=1}^3{\mu_j}'E_j\Big)^\sim \in \gh_{\sR}$. Indeed, from Theorem 3.5.3, we have 
\begin{eqnarray*}
      B_6(h, h') \!\!\! &=& \!\!\! \dfrac{4}{3}B_4\Big(\dsum_{k=0}^3\lambda_kH_k, \dsum_{k=0}^3{\lambda_k}'H_k \Big) + 12\Big(\dsum_{j=1}^3\mu_kE_k, \dsum_{j=1}^3{\mu_j}'E_j \Big) \\    
       \!\!\! &=& \!\!\! \dfrac{4}{3}18\dsum_{k=0}^3\lambda_k{\lambda_k}' + 12\dsum_{j=1}^3\mu_j{\mu_j}'\;\;\mbox{(Theorem 2.6.2)}
\\ 
       \!\!\! &=& \!\!\! 12\Big(2\dsum_{k=0}^3\lambda_k{\lambda_k}'+
               \dsum_{j=1}^3\mu_j{\mu_j}'\Big).
\end{eqnarray*}
Now, the canonical elements $H_{\alpha_i} \in \gh_{\sR}$ corresponding to $\alpha_i$ ($B_6 (H_\alpha, H) = \alpha(H), H \in \gh_{\sR} $) are determined as follows.
\vspace{2mm}

\qquad \quad
       $H_{\alpha_1} = \dfrac{1}{24}(H_0 - H_1), \quad
       H_{\alpha_2} = \dfrac{1}{24}(H_1 - H_2), \quad
       H_{\alpha_3} = \dfrac{1}{24}(H_2 - H_3)$, 
\vspace{1mm}

\qquad \quad
       $H_{\alpha_4} = \dfrac{1}{24}(H_3 + (E_2 - E_3)^{\sim})$,
\vspace{1mm}

\qquad \quad
       $H_{\alpha_5} = \dfrac{1}{48}((- H_0 - H_1 - H_2 + H_3) + 2(E_3 - E_1)^{\sim})$,
\vspace{1mm}

\qquad \quad
       $H_{\alpha_6} = \dfrac{1}{48}((H_0 + H_1 + H_2 + H_3) + 2(E_1 - E_2)^{\sim})$.
\vspace{2mm}

\noindent Therefore, we have
$$
        (\alpha_1, \alpha_1) = B_6(H_{\alpha_1}, H_{\alpha_1}) = 24\frac{1}{24}\frac{1}{24}2 = \frac{1}{12}, $$
and the other inner products are similarly calculated. Hence, the inner product induced by the Killing form $B_6$ between $\alpha_1$, $\alpha_2, \cdots, \alpha_6$ and $- \mu$ are given by
$$
\begin{array}{l}
      (\alpha_i, \alpha_i) = \dfrac{1}{12}, \quad i=1, 2, 3, 4, 5, 6,
\vspace{1mm}\\
      (\alpha_1, \alpha_2) = (\alpha_2, \alpha_3) = (\alpha_3, \alpha_4) = (\alpha_3, \alpha_5) = (\alpha_5, \alpha_6) = - \dfrac{1}{24}, 
\vspace{1mm}\\
       (\alpha_i, \alpha_j) = 0, \quad \mbox{otherwise},
\vspace{1mm}\\
       (-\mu, -\mu) = \dfrac{1}{12}, \;\;\; (- \mu, \alpha_4) = - \dfrac{1}{24}, \;\;\; (- \mu, \alpha_i) = 0, \quad i = 1, 2, 3, 5, 6, 
\end{array} $$                                                                  
using them, we can draw the Dynkin diagram and the extended Dynkin diagram of ${\gge_6}^C$.                            
\vspace{2mm}

According to Borel-Siebenthal theory, the Lie algebra $\gge_6$ has three subalgebras as maximal subalgebras with the maximal rank 6.
\vspace{1mm}

(1) The first one is a subalgebra of type $T \oplus D_5$ which is obtained as the fixed points of an involution $\sigma$ of $\gge_6$.
\vspace{1mm}

(2) The second one is a subalgebra of type $C_1 \oplus A_5$ which is obtained as the fixed points of an involution $\gamma$ of $\gge_6$.
\vspace{1mm}

(3) The third one is a subalgebra of type $A_2 \oplus A_2 \oplus A_2$ which is obtained as the fixed points of an automorphism $w$ of order 3 of $\gge_6$.
\vspace{1mm}

The Lie algebra $\gge_6$ has furthermore two outer involutions $\tau$, $\tau\gamma$. The subalgebra obtained as the fixed points of $\tau$ is type $F_4$ and the subalgebra obtained as the fixed points of $\tau\gamma$ is type $C_4$. 
\vspace{1mm}

These subalgebras will be realized as subgroups of the group $E_6$ in Theorems 3.10.7, 3.11.4, 3.13.5, 3.7.1 and 3.12.2, respectively.
\vspace{4mm}

{\bf 3.7. Involution $\tau$ and subgroup $F_4$ of $E_6$}
\vspace{3mm}

We shall study the following subgroup $(E_6)^{\tau}$ of $E_6$:
\begin{eqnarray*}
  (E_6)^{\tau} \!\!\! &=& \!\!\! \{ \alpha \in E_6 \, | \, \tau\alpha = \alpha\tau \} \\                
     \!\!\! &=& \!\!\! \{ \alpha \in E_6 \, | \, \lambda(\alpha) = \alpha \}  = (E_6)^{\lambda}.
\end{eqnarray*}         
If $\alpha \in E_6$ satisfies $\tau\alpha = \alpha\tau$, then $(\alpha X, \alpha Y) = \langle \tau\alpha X, \alpha Y \rangle = \langle \alpha\tau X, \alpha Y \rangle = \langle \tau X, Y \rangle =(X, Y)$ and vice versa. Further, for $\alpha \in E_6$, the conditions $(\alpha X, \alpha Y) = (X, Y)$ and $\alpha E = E$ are equivalent (Lemma 2.2.4). Hence $(E_6)^\tau$ can be also defined by
\begin{eqnarray*}
  (E_6)^{\tau} \!\!\! &=& \!\!\! \{ \alpha \in E_6 \, | \, (\alpha X, \alpha Y) = (X, Y), X, Y \in \gJ^C\} \\                
     \!\!\! &=& \!\!\! \{ \alpha \in E_6 \, | \, \alpha E = E\} = (E_6)_E.
\end{eqnarray*} 

{\bf Theorem 3.7.1.} \qquad \qquad $(E_6)^{\tau} = (E_6)_E \cong F_4.$
\vspace{1mm}

\noindent (From now on, we identify these groups).
\vspace{2mm}

{\bf Proof.} We shall show $(E_6)^\tau \cong F_4$. Let $\alpha \in (E_6)^{\tau}$. Then, for $X \in \gJ$  we have $\tau\alpha X = \alpha\tau X = \alpha X$, so that $\alpha X \in \gJ$. Hence $\alpha$ induces an $\R$-linear transformation of $\gJ$. Therefore, the restriction $\alpha|\gJ$ of $\alpha$ to $\gJ$ belongs to the group
\begin{eqnarray*}
      F_4 \!\!\! &=& \!\!\! \{ \alpha \in \Iso_{\sR}(\gJ) \, | \, \det(\alpha X) = \det X,(\alpha X, \alpha Y) = (X, Y) \} \\
          \!\!\! &=& \!\!\! \{ \alpha \in \Iso_{\sR}(\gJ) \, | \, \alpha(X \times Y) = \alpha X \times \alpha Y \}. 
\end{eqnarray*}
Conversely, for $\alpha \in F_4$, its complexification $\alpha^C : \gJ^C \to \gJ^C$, $\alpha^C(X_1 + iX_2) = \alpha X_1 + i\alpha X_2$ belongs to $(E_6)^{\tau}$. Therefore the correspondence $F_4 \ni \alpha \to \alpha^C \in (E_6)^{\tau}$ 
gives an isomorphism between $F_4$ and $(E_6)^{\tau}$.
\vspace{4mm}

{\bf 3.8. Connectedness of $E_6$}
\vspace{3mm}

We denote by $(E_6)_0$ the connected component of $E_6$ containing the identity $1$.
\vspace{3mm}

{\bf Lemma 3.8.1.} (1) {\it For} $t \in \R$, {\it if we define a mapping} $\alpha_{12}(t) : \gJ^C \to \gJ^C$ {\it by} 
$$
      \alpha_{12}(t)\pmatrix{\xi_1 & x_3 & \ov{x}_2 \cr
                             \ov{x}_3 & \xi_2 & x_1 \cr
                             x_2 & \ov{x}_1 & \xi_3}
                   = \pmatrix{e^{it}\xi_1 & x_3 & e^{it/2}\ov{x}_2 \cr
                                 \ov{x}_3 & e^{-it}\xi_2 & e^{-it/2}x_1 \cr 
                                 e^{it/2}x_2 & e^{-it/2}\ov{x}_1 & \xi_3}, $$
{\it then} $\alpha_{12}(t) \in (E_6)_0$. {\it Similarly we can define } $\alpha_{13}(t)$, $\alpha_{23}(t) \in (E_6)_0$.
\vspace{1mm}

(2) {\it For} $a \in \gC$, {\it if we define a mapping} $\alpha_1(a) : \gJ^C \to \gJ^C$ {\it by} $\alpha_1(a)X(\xi, x) = Y(\eta, y),$ {\it where}
\begin{eqnarray*}
&&\left\{\begin{array}{l}
     \eta_1 = \xi_1 
\vspace{1mm}\\
     \eta_2 = \displaystyle{\frac{\xi_2 - \xi_3}{2}} 
            + \displaystyle{\frac{\xi_2 + \xi_3}{2}}\cos|a|
            + i\displaystyle{\frac{(a,x_1)}{|a|}}\sin|a|
\vspace{1mm}\\
     \eta_3 = - \displaystyle{\frac{\xi_2 - \xi_3}{2}}
            + \displaystyle{\frac{\xi_2 + \xi_3}{2}}\cos|a|
            + i\displaystyle{\frac{(a,x_1)}{|a|}}\sin|a|,
\end{array}\right.\\
&&\left\{\begin{array}{l}
      y_1 = x_1 
            + i\displaystyle{\frac{(\xi_2 + \xi_3)a}{2|a|}}\sin|a| 
            - \displaystyle{\frac{2(a, x_1)a}{|a|^2}}\sin^2\displaystyle{\frac{|a|}{2}}
\vspace{1mm}\\
     y_2 = x_2\cos\displaystyle{\frac{|a|}{2}} + 
           i\displaystyle{\frac{\ov{x_3a}}{|a|}}\sin\displaystyle{\frac{|a|}{2}}
\vspace{1mm}\\
     y_3 = x_3\cos\displaystyle{\frac{|a|}{2}} + 
          i\displaystyle{\frac{\ov{ax_2}}{|a|}}\sin\displaystyle{\frac{|a|}{2}}
\end{array}\right.
\end{eqnarray*}
\Big({\it if $a = 0$, then $\displaystyle{\frac{\sin|a|}{|a|}}$ means $1$}\Big), {\it then} $\alpha_1(a) \in (E_6)_0$.
\vspace{2mm}

{\bf Proof.} (1) For $E_1 - E_2 \in \gJ_0$, we have $i(E_1 - E_2)^{\sim} \in \gge_6$ (Theorem 3.2.4.(2)) and $\alpha_{12}(t) = \exp it(E_1 - E_2)^{\sim}$. Hence $\alpha_{12}(t) \in (E_6)_0$.
\vspace{1mm}

(2) For $F_1(a) \in \gJ_0$, we have $i\wti{F}_1(a) \in \gge_6$ (Theorem 3.2.4.(2)) and $\alpha_1(a) = \exp i\wti{F}_1(a)$. Hence $\alpha_1(a) \in (E_6)_0$.
\vspace{3mm}

{\bf Proposition 3.8.2.} {\it Any element} $X \in \gJ^C$ {\it can be transformed to a diagonal form by some element} $\alpha \in (E_6)_0$: 

$$
      \alpha X = \pmatrix{\xi_1 & 0 & 0 \cr
                          0 & \xi_2 & 0 \cr
                          0 & 0 & \xi_3},
      \quad \xi_i \in C.$$
{\it Moreover, we can choose $\alpha \in (E_6)_0$ so that two of $\xi_1, \xi_2,\xi_3$ are non-negative real numbers.}
\vspace{2mm}

{\bf Proof.} (Hereafter, we use the notation $|\xi|$ instead of $\sqrt{\xi(\tau \xi)}$ for $\xi \in \R^C = C$). For a given $X \in \gJ^C$, consider a space $\gX = \{ \alpha X \, | \, \alpha \in (E_6)_0 \}$. Since $(E_6)_0$ is compact (Theorem 3.1.1), $\gX$ is also compact. Let $|\xi_1|^2 + |\xi_2|^2 + |\xi_3|^2$ be the maximal value of all $|\eta_1|^2 + |\eta_2|^2 + |\eta_3|^2$ for $Y = Y(\eta,y) \in \gX$ and let $X_0 = X(\xi, x)$ be an element of $\gX$ which attains its maximal value. Then $X_0$ is of diagonal form. Certainly, suppose that $X_0$ is not of diagonal form, for example, the $2 \times 3$ entry $x_1$ of $X_0$ is non-zero: 
$$
        0 \neq x_1 = p + i q, \quad p,q \in \gC. $$
It is sufficient to prove in the case that $\xi_2,\xi_3$ are real numbers (otherwise we can apply some $\alpha_{12}(t_1)$ and $\alpha_{13}(t_2)$ of Lemma 3.8.1.(1)).
\vspace{1mm}

(1) Case $q \neq 0$. Let $a(t) = \dfrac{q}{|q|}t$, $t > 0$ and construct $\alpha_1(a(t)) \in (E_6)_0$ of Lemma 3.8.1.(2). Since $|a(t)|=t$ and $i\dfrac{(a(t),x_1)}{|a(t)|} = i\nu - |q|$, where $\nu = \Big(\dfrac{q}{|q|}, p \Big),$ for $Y(\eta(t), y(t)) = \alpha_1(a(t))X_0 \in \gX$, we have
$$
\begin{array}{l}
  |\eta_1(t)|^2 + |\eta_2(t)|^2 + |\eta_3(t)|^2 
\vspace{1mm}\\
\qquad \quad
   = |\xi_1|^2 + \Big|\dfrac{\xi_2 - \xi_3}{2} + \dfrac{\xi_2 + \xi_3}{2}\cos t -|q|\sin t + i\nu\sin t \Big|^2 
\vspace{1mm}\\
\qquad \qquad 
   + \Big|- \dfrac{\xi_2-\xi_3}{2} + \dfrac{\xi_2+\xi_3}{2}\cos t - |q|\sin t + i\nu\sin t \Big|^2 
\vspace{1mm}\\
\qquad \quad
   = |\xi_1|^2 + 2\Big(\dfrac{\xi_2 - \xi_3}{2} \Big)^2 + 2\Big(\dfrac{\xi_2 + \xi_3}{2}\cos t -|q|\sin t \Big)^2 + 2\nu^2\sin^2 t
\vspace{1mm}\\
\qquad \quad
   = |\xi_1|^2 + 2\Big(\dfrac{\xi_2-\xi_3}{2}\Big)^2 + 2\Big(\Big(\dfrac{\xi_2+\xi_3}{2}\Big)^2 + |q|^2 \Big)\sin^2(t + t_0) + 2\nu^2\sin^2t 
\vspace{1mm}\\
\qquad \quad
   \le |\xi_1|^2 + |\xi_2|^2 + |\xi_3|^2 + 2|q|^2 + 2\nu^2\cos^2t_0  \;\;\;(\mbox{for some }\; t_0 \in \R)
\end{array}$$
which is the maximal value and attains at some $t > 0$. This contradicts the maximum of $|\xi_1|^2 + |\xi_2|^2 + |\xi_3|^2$. Hence $q = 0$.
\vspace{1mm}

(2) Case $p \neq 0$. Let $a(t) = \dfrac{p}{|p|}t$, $ t > 0$ and construct $\beta_1(a(t)) \in F_4 \subset (E_6)_0$ of Lemma 2.8.1. For $Y(\eta(t), y(t)) = \beta_1(a(t))X_0 \in \gX$, the maximal value of $|\eta_1(t)|^2 + |\eta_2(t)|^2 + |\eta_3(t)|^2$ is $|\xi_1|^2 + |\xi_2|^2 + |\xi_3|^2 + 2|p|^2$ (the calculation is the same as Proposition 2.8.2) which contradicts the maximum of $|\xi_1|^2 + |\xi_2|^2 + |\xi_3|^2$. Hence $p = 0$. Consequently, we have $x_1 = 0$. $x_2 = x_3 = 0$ can be similarly proved constructing $\alpha_2(a)$, $\alpha_3(a) \in (E_6)_0$ analogous to $\alpha_1(a)$ of Lemma 3.8.1.(2). Hence $X_0$ is of diagonal form.
\vspace{2mm}

The space $EIV$, called the symmetric space of type $EIV$, is defined by
$$
      EIV = \{ X \in \gJ^C \, | \, \det X = 1,\langle X, X \rangle = 3 \}.$$

{\bf Theorem 3.8.3.} \qquad \qquad $E_6/F_4 \simeq EIV.$
\vspace{1mm}

\noindent {\it In particular, $E_6$ is connected.}
\vspace{2mm}

{\bf Proof.} For $\alpha \in E_6$ and $X \in EIV$, we have $\alpha X \in EIV$. Hence $E_6$ acts on $EIV$. We shall prove that the group $(E_6)_0$ acts transitively on $EIV$. To prove this, it is sufficient to show that any element $X \in EIV$ can be transformed to $E \in EIV$ by some $\alpha \in (E_6)_0$. Now, $X \in EIV \subset \gJ^C$ can be transformed to a diagonal form by $\alpha \in (E_6)_0$:
$$
      \alpha X = \pmatrix{\xi_1 & 0 & 0 \cr 
                          0 & \xi_2 & 0 \cr
                          0 & 0 & \xi_3}, 
                 \quad \xi_1 \in C,\xi_2 \ge 0,\xi_3 \ge 0$$
(Proposition 3.8.2). From the condition $X \in EIV$, we have
$$
\begin{array}{l}
      \xi_1\xi_2\xi_3 = \det\,(\alpha X) = \det\,X = 1, \;\;(\mbox{hence} \; \; \xi_i > 0, i=1, 2, 3),
\vspace{1mm}\\
      {\xi_1}^2 + {\xi_2}^2 + {\xi_3}^2 = \langle \alpha X, \alpha Y \rangle = \langle X, Y \rangle = 3.
\end{array} $$
This implies that $\xi_1 = \xi_2 = \xi_3 = 1$. Certainly, from $0 \le (\xi_2 - \xi_3)^2 = {\xi_2}^2 + {\xi_3}^2 
\vspace{0.5mm}
- 2\xi_2\xi_3 = 3 - {\xi_1}^2 - \displaystyle{\frac{2}{\xi_1}} = - \frac{{\xi_1}^3 - 3\xi_1 + 2}{\xi_1} = -\dfrac{(\xi_1 - 1)^2(\xi_1 + 2)}{\xi_1} \le 0$, we have 
\vspace{0.5mm}
$\xi_1 = 1$. Similarly $\xi_2 = \xi_3 = 1$ are obtained. Hence $\alpha X = E,$ which shows the transitivity of $(E_6)_0$. Since we have $EIV = (E_6)_0E$, $EIV$ is connected. Now, the group $E_6$ acts transitively on $EIV$ and the isotropy subgroup of $E_6$ at $E \in EIV$ is $F_4$ (Theorem 3.7.1). Thus we have the homeomorphism $E_6/F_4 \simeq EIV$. Finally, the connectedness of $E_6$ follows from the connectedness of $EIV$ and $F_4$.
\vspace{4mm}

{\bf 3.9. Center $z(E_6)$ of $E_6$}
\vspace{3mm}

{\bf Theorem 3.9.1.}  {\it The center} $z(E_6)$ {\it of the group} $E_6$ {\it is isomorphic to the cyclic group of order} 3:
$$
        z(E_6) = \{ 1, \omega 1, \omega^2 1 \}, \quad 
        \omega = - \frac{1}{2} + \frac{\sqrt{3}}{2}i \in C. $$

{\bf Proof.} Let $\alpha \in z(E_6)$. From the commutativity with $\beta \in F_4 \subset E_6$, we have $\beta\alpha E = \alpha\beta E = \alpha E$. Let denote $\alpha E = Y = Y(\eta, y) \in \gJ^C$, then we have 
$$
        \beta Y = Y, \quad \mbox{for all} \;\; \beta \in F_4.$$
We choose $\beta \in F_4$ such that
$$
         \beta X = TXT^{-1}, \quad X \in \gJ^C,$$
where $T = \pmatrix{1 & 0 & 0 \cr
                    0 & -1 & 0 \cr
                    0 & 0 & -1}$, 
          $\pmatrix{-1 & 0 & 0 \cr
                    0 & 1 & 0 \cr
                    0 & 0 & -1}$ and 
          $\pmatrix{0 & 0 & 1 \cr
                    1 & 0 & 0 \cr 
                    0 & 1 & 0} \in SO(3). $ 
Then we have $y_1 = y_2 = y_3 = 0$ and $\eta_1 = \eta_2 = \eta_3 \;( = \omega )$, that is,
$$
        \alpha E = Y = \omega E, \quad \omega \in C,$$
and $\omega^3 = \det(\alpha E) = \det E = 1$. Since it is easy to verify that $\omega 1 \in z(E_6)$, we have $\omega^{-1}\alpha \in z(E_6)$ and $\omega^{-1}\alpha E = E$, and so $\omega^{-1}\alpha \in z(F_4)$ (Theorem 3.7.1). Since $z(F_4) = \{1\}$ (Theorem 2.10.1), we have $\omega^{-1}\alpha = 1$, that is, $\alpha = \omega 1$. This completes the proof.
\vspace{2mm}

According to the general theory of compact Lie groups, it is known that the center of the simply connected compact simple Lie group of type $E_6$ is the cyclic group of order 3.  Hence the group $E_6$ has to be simply connected. Thus we have the following theorem.
\vspace{3mm}

{\bf Theorem 3.9.2.} $E_6 = \{\alpha \in \Iso_C(\gJ^C) \, | \, \det\,(\alpha X) = \det\,X, \langle \alpha X, \alpha Y \rangle = \langle X, Y \rangle \}$ {\it is a simply connected compact Lie group of type $E_6$.}
\vspace{4mm}

{\bf 3.10. Involution $\sigma$ and subgroup $(U(1) \times Spin(10))/\Z_4$ of $E_6$}
\vspace{3mm}

Let the $C$-linear mapping $\sigma : \gJ^C \to \gJ^C$ be the complexification of $\sigma \in F_4$ of Section 2.9. Then $\sigma \in E_6$ and $\sigma^2 = 1.$
\vspace{2mm}

We shall study the following subgroup $(E_6)^{\sigma}$ of $E_6$:
$$
      (E_6)^{\sigma} = \{ \alpha \in E_6 \, | \, \sigma\alpha = \alpha\sigma \}.$$
To this end, we consider the $C$-subspaces $(\gJ^C)_{\sigma}$ and $(\gJ^C)_{-\sigma}$ of $\gJ^C$, which are the eigenspaces of $\sigma$:
\begin{eqnarray*}
     (\gJ^C)_{\sigma} \!\!\! &=& \!\!\! \{ X \in \gJ^C \, | \, \sigma X = X \}
\vspace{1mm}\\
         \!\!\! &=& \!\!\! \{ X \in \gJ^C \, | \, 4E_1 \times (E_1 \times X) = X  \} \oplus {\gE_1}^C,
\vspace{1mm}\\
     (\gJ^C)_{-\sigma} \!\!\! &=& \!\!\! \{ X \in \gJ^C \, | \, \sigma X = -X \}
\vspace{1mm}\\
      \!\!\! &=& \!\!\! \{ X \in \gJ^C \, | \, E_1 \times X = 0, \langle E_1, X \rangle = 0 \},
\end{eqnarray*}
where ${\gE_1}^C = \{ \xi E_1 \, | \, \xi \in C \}.$  Then $\gJ^C = (\gJ^C)_{\sigma} \oplus (\gJ^C)_{-\sigma}$ (which is the complexification of $\gJ = \gJ_{\sigma} \oplus \gJ_{-\sigma}$ in Section 2.9), and $(\gJ^C)_{\sigma}$, $(\gJ^C)_{-\sigma}$ are invariant under the action of $(E_6)^{\sigma}$.
\vspace{3mm}

{\bf Lemma 3.10.1.} {\it For $\alpha \in (E_6)^{\sigma}$, there exists $\xi \in C$ such that }
$$
         \alpha E_1 = \xi E_1, \quad (\tau\xi)\xi = 1.$$

{\bf Proof.} By the analogous proof to that of Theorem 2.9.1, we see that 
$$ 
        \alpha E_2, \; \alpha E_3 \in \gJ(2, \gC^C). $$
Indeed, we have 
\begin{eqnarray*}
     \alpha E_2 \!\!\! &=& \!\!\! \alpha(- F_2(1) \times F_2(1)) = - \tau\alpha\tau F_2(1) \times \tau\alpha\tau F_2(1) 
\vspace{1mm}\\
  \!\!\! &=& \!\!\! - (F_2(x_2) + F_3(x_3)) \times (F_2(x_2) + F_3(x_3)) \;\; (\mbox{for some} \;\; x_2, x_3 \in \gC^C) 
\vspace{1mm}\\
  \!\!\! &=& \!\!\! (x_2, x_2)E_2 + (x_3, x_3)E_3  - F_1(\ov{x_2x_3}) \in \gJ(2, \gC^C).
\end{eqnarray*}
Next, we shall show
$$ 
                \alpha E_1 \not\in \gJ(2, \gC^C).$$
Suppose that $\alpha E_1 \in \gJ(2, \gC^C)$. Then $\alpha E = \alpha E_1 + \alpha E_2 + \alpha E_3 \in \gJ(2, \gC^C)$, so we can put $\alpha E = \xi_2E_2 + \xi_3 E_3 + F_1(x_1)$, $\xi_2, \xi_3 \in C$, $x_1 \in \gC^C$. Hence
$$
\begin{array}{l}
      \xi_2E_2 + \xi_3E_3+F_1(x_1) = \alpha E = \alpha(E \times E) = \tau\alpha\tau E \times \tau\alpha\tau E 
\vspace{1mm}\\
      \qquad = \tau(\xi_2E_2 + \xi_3E_3 + F_1(x_1)) \times \tau(\xi_2 E_2 + \xi_3 E_3 + F_1(x_1))
\vspace{1mm}\\
      \qquad = (\tau \xi_2\tau \xi_3 - (\tau x_1, \tau x_1))E_1.
\end{array} $$
This implies that $\xi_2 = \xi_3 = x_1 = 0$. Hence $\alpha E = 0$, which is a contradiction. Therefore $\alpha E_1$ is of the form
$$
      \alpha E_1 = \xi E_1 + \xi_2E_2 + \xi_3E_3 + F_1(x_1), \quad \xi \neq 0.$$From $\alpha E_1 \times \alpha E_1 = \tau\alpha\tau(E_1 \times E_1) = 0$, we 
have 
\begin{eqnarray*}
     0 \!\!\! &=& \!\!\! (\xi E_1 + \xi_2E_2 + \xi_3E_3 + F_1(x_1)) \times 
                         (\xi E_1 + \xi_2E_2 + \xi_3E_3 + F_1(x_1))             \vspace{1mm}\\
     \!\!\! &=& \!\!\! (\xi_2\xi_3 - (x_1, x_1))E_1 + \xi \xi_3E_2 + \xi \xi_2E_3 - \xi F_1(x_1).
\end{eqnarray*}
This implies that $\xi_2 = \xi_3 = x_1 = 0$ and so
$$
            \alpha E_1 = \xi E_1, \quad \xi \neq 0.$$
Finally, from
$$
      1 = \langle E_1, E_1 \rangle = \langle \alpha E_1, \alpha E_1 \rangle = \langle \xi E_1, \xi E_1 \rangle = (\tau\xi)\xi\langle E_1, E_1 \rangle = (\tau\xi)\xi, $$
we have $(\tau\xi)\xi = 1$.
\vspace{3mm}

In order to investigate the group $(E_6)^{\sigma}$, we consider the following subgroup $(E_6)_{E_1}$ of $E_6$\,:
$$
          (E_6)_{E_1} = \{ \alpha \in E_6 \, | \, \alpha E_1 = E_1 \}.$$

{\bf Lemma 3.10.2.} $(E_6)_{E_1}$ {\it is a subgroup of} $(E_6)^{\sigma}$: 
$(E_6)_{E_1} \subset (E_6)^{\sigma}$.
\vspace{2mm}

{\bf Proof.} Since $(\gJ^C)_{\sigma} = \{ X \in \gJ^C \, | \, 4E_1 \times (E_1 \times X) = X \} \oplus {\gE_1}^C$ and $(\gJ^C)_{-\sigma} = \{ X \in \gJ^C \, | \, E_1 \times X = 0, \langle E_1, X \rangle = 0 \}$, these spaces are seen to be invariant under the action of $(E_6)_{E_1}$. Hence for $\alpha \in (E_6)_{E_1},$ we have $\sigma\alpha = \alpha\sigma$ (the proof is the same as that of Theorem 2.9.1). Thus we have $\alpha \in (E_6)^{\sigma}.$
\vspace{2mm}

We define a 10 dimensional $\R$-vector space $V^{10}$ by
$$
      V^{10} = \{ X \in \gJ^C \, | \, 2E_1 \times X = - \tau X \} 
             = \Big\{ \pmatrix{0 & 0 & 0 \cr
                                0 & \xi & x \cr
                                0 & \ov{x} & - \tau\xi} \, \Big| \, 
        \xi \in C, x \in \gC \Big\}. $$

{\bf Proposition 3.10.3.} \qquad \qquad $(E_6)_{E_1}/Spin(9) \simeq S^9.$
\vspace{1mm}

\noindent {\it In particular, the group $(E_6)_{E_1}$ is connected.}
\vspace{2mm}

{\bf Proof.} $S^9 = \{ X \in V^{10} \, | \,  \langle X, X \rangle = 2 \}$ is a 9 dimensional sphere. For $\alpha \in (E_6)_{E_1}$ and $X \in S^9$, we have $\alpha X \in S^9$. Indeed,
$$
\begin{array}{c}
     2E_1 \times \alpha X = 2\alpha E_1 \times \alpha X = 2\tau\alpha\tau(E_1 \times X) = \tau\alpha\tau(- \tau X) = - \tau(\alpha X), 
\vspace{1mm}\\
       \langle \alpha X, \alpha X \rangle = \langle X, X \rangle = 2.
\end{array} $$
Hence the group $(E_6)_{E_1}$ acts on $S^9$. We shall prove that the action is transitive. To prove this, it is sufficient to show that any element $X \in S^9$ can be transformed to $i(E_2 + E_3) \in S^9$ by some $\alpha \in (E_6)_{E_1}$. Now, for a given $X \in S^9$, we can choose $\alpha_{23}(t_0)$ of Lemma 3.8.1.(1) such that
$$
      \alpha_{23}(t_0)X \in S^8 = \{ X \in V^9 \, | \, \langle X, X \rangle = 2 \}$$
where $V^9 = \{ X \in V^{10} \, |\, \tau X = X \}$. (Note that $\alpha_{23}(t_0) \in (E_6)_{E_1}$ because $\alpha_{23}(t_0)E_1 = E_1$). Since the group $Spin(9) =(F_4)_{E_1} \subset (E_6)_{E_1}$ acts transitively on $S^8$ (Proposition 2.7.3), there exists $\beta \in Spin(9)$ such that 
$$
           \beta\alpha_{23}(t_0)X = E_2 - E_3 \in S^8. $$
By applying $\alpha_{23}(\pi/2) \in (E_6)_{E_1}$ of Lemma 3.8.1.(1), we get
$$
       \alpha_{23}(\pi/2)\beta\alpha_{23}(t_0)X = i(E_2 + E_3). $$
This shows the transitivity. The isotropy subgroup of $(E_6)_{E_1}$ at $i(E_2 + E_3) \in S^9$ is $Spin(9)$. Indeed, if $\alpha \in (E_6)_{E_1}$ satisfies $\alpha(i(E_2 + E_3)) = i(E_2 + E_3)$, then $\alpha E = \alpha E_1 + \alpha(E_2 + E_3) = E_1 + (E_2 + E_3) = E$, so that $\alpha \in F_4$ and hence $\alpha \in (F_4)_{E_1} = Spin(9)$. Conversely $\alpha \in Spin(9)$ satisfies $\alpha(i(E_2 + E_3)) = i(E_2 + E_3)$. Thus we have the homeomorphism $(E_6)_{E_1}/Spin(9) \simeq S^9$. 
\vspace{3mm}

{\bf Theorem 3.10.4.} \qquad \qquad $(E_6)_{E_1} \cong Spin(10). $
\vspace{1mm}
 
\noindent (From now on, we identify these groups).
\vspace{2mm}

{\bf Proof.} Analogously to Theorem 2.7.4, we can define a homomorphism 
$$
      p : (E_6)_{E_1} \to SO(10) = SO(V^{10}) $$ 
by $p(\alpha) = \alpha|V^{10}$. The restriction $p'$ of $p : (E_6)_{E_1} \to SO(10)$ to $(F_4)_{E_1}$ coincides with the homomorphism $p' : Spin(9) \to SO(9)$ of Theorem 2.7.4. In particular, $p' : Spin(9) \to SO(9)$ is onto. Hence, from the following commutative diagram
\begin{center}
\begin{tabular}{ccccccccc}
     $1$ & $\longrightarrow$ & $Spin(9)$ & $\longrightarrow$ & $(E_6)_{E_1}$
         & $\longrightarrow$ & $S^9$ & $\longrightarrow$ & $*$ \vspace{0.7mm}\\
     $$ & $$ & $\downarrow p'$ & $$ & $\downarrow p$ & $$ & $\downarrow =$ & $$ & $$\vspace{0.7mm}\\
     $1$ & $\longrightarrow$ & $SO(9)$ & $\longrightarrow$ & $SO(10)$ 
         & $\longrightarrow$ & $S^9$ & $\longrightarrow$ & $*$
\end{tabular}
\end{center}
we see that $p : (E_6)_{E_1} \to SO(10)$ is onto by the five lemma. We also have $\Ker \, p = \{1, \sigma \}$. Indeed, $\alpha \in \Ker \, p$ leaves $E_2 - E_3$ and $i(E_2 + E_3)$ invariant, so  that $\alpha E_i = E_i$ for $i = 1, 2, 3$. Hence we have $\alpha \in Spin(8)$ and so $\alpha \in \Ker \, p'$. Therefore $\alpha = 1$ or $\sigma$ by Theorem 2.7.4. Hence we have the isomorphism
$$ 
         (E_6)_{E_1}/\{1, \sigma \} \cong SO(10).$$
Therefore the group $(E_6)_{E_1}$ is isomorphic to the group $Spin(10)$ as 
the universal covering group of $SO(10)$.
\vspace{2mm}

For $\theta \in C$, $\theta \neq 0$, we define a $C$-linear mapping $\phi(\theta) : \gJ^C \to \gJ^C$ by
$$
       \phi(\theta)\pmatrix{\xi_1 & x_3 & \ov{x}_2 \cr 
                            \ov{x}_3 & \xi_2 & x_1 \cr 
                            x_2 & \ov{x}_1 & \xi_3} 
           = \pmatrix{\theta^4\xi_1 & \theta x_3 & \theta\ov{x}_2 \cr
                      \theta\ov{x}_3 & \theta^{-2}\xi_2 & \theta^{-2}x_1 \cr 
                      \theta x_2 & \theta^{-2}\ov{x}_1 & \theta^{-2}\xi_3}. $$

{\bf Theorem 3.10.5.} {\it The group $E_6$ contains a subgroup}
\begin{eqnarray*}
      U(1) \!\!\! &=& \!\!\! \{ \alpha_{12}(t)\alpha_{13}(t) = \exp it(2E_1 - E_2 - E_3)^{\sim} \, | \, t \in \R \}
\vspace{1mm}\\
      \!\!\! &=& \!\!\! \{ \phi(\theta) \, | \, \theta \in C, (\tau\theta)\theta = 1 \}
\end{eqnarray*}
({\it where $\alpha_{12}(t), \alpha_{13}(t)$ are mappings defined in Lemma} 3.8.1) {\it which is isomorphic to the usual unitary group} $U(1) = \{ \theta \in C \, | \, (\tau\theta)\theta = 1 \}$. 
\vspace{1mm}

Note that $U(1)$ is a subgroup of $(E_6)^{\sigma}$. From now on, we identify these two groups $U(1)$.
\vspace{3mm}

{\bf Lemma 3.10.6.} {\it Two subgroups $U(1)$ and $Spin(10)$ of $(E_6)^{\sigma}$ 
are elementwise commutative}.
\vspace{2mm}

{\bf Proof.} We consider the decomposition $\gJ^C = {\gE_1}^C \oplus \gJ(2, \gC)^C \oplus (\gJ^C)_{-\sigma}$, where
$$
\begin{array}{c}
      {\gE_1}^C = \{ \xi E_1 \, | \, \xi \in C \}, \quad 
     \gJ(2, \gC)^C = \{ X \in \gJ^C \, | \, 4E_1 \times (E_1 \times X) = X \},
\vspace{1mm}\\
      (\gJ^C)_{-\sigma} = \{ X \in \gJ^C \, | \, E_1 \times X = 0, 
                           \langle E_1, X \rangle = 0 \}.
\end{array} $$
The restrictions of $\phi(\theta) \in U(1)$ to these spaces are all constant mappings:
$$
      \phi(\theta)|{\gE_1}^C = \theta^41, \quad
      \phi(\theta)|{\gJ(2, \gC)^C} = \theta^{-2}1, \quad
      \phi(\theta)|(\gJ^C)_{-\sigma} = \theta 1. $$
On the other hand, $\beta \in Spin(10)$ also induces $C$-linear transformations of these spaces. From this, the commutativity of $\phi(\theta)$ and $\beta$: $\phi(\theta)\beta = \beta \phi(\theta)$ follows.
\vspace{3mm}

{\bf Theorem 3.10.7.} $(E_6)^{\sigma} \cong (U(1) \times Spin(10))/\Z_4,$ \, $\Z_4 = \{(1, \phi(1)), (-1, \phi(-1)), $ $(i, \phi(-i)), (-i, \phi(i)) \}.$
\vspace{2mm}

{\bf Proof.} We define a mapping $\varphi : U(1) \times Spin(10) \to (E_6)^{\sigma}$ by
$$
       \varphi(\theta, \beta) = \phi(\theta)\beta. $$
Since $\phi(\theta) \in U(1)$ and $\beta \in Spin(10)$ are commutative (Lemma 3.10.6), we see that $\varphi$ is a homomorphism. We shall show that $\varphi$ is onto. For $\alpha \in (E_6)^{\sigma}$, there exists $\theta \in C$, $(\tau\theta)\theta = 1$ satisfying
$$
      \alpha E_1 = \theta^4E_1 = \phi(\theta)E_1 $$
(Lemma 3.10.1). Let $\beta = \phi(\theta)^{-1}\alpha$, then $\beta E_1 = E_1$, so $\beta \in Spin(10)$ (Theorem 3.10.4). Hence, we have $\alpha =\phi(\theta)\beta = \varphi(\theta, \beta)$ and therefore $\varphi$ is onto. It is easily seen that 
$$
      \Ker\,\varphi = \{ (\theta, \phi(\theta)^{-1}) \, | \,  \theta \in C, \theta^4 = 1 \} = \{ (\theta, \phi(\theta^{-1})) \, | \, \theta = \pm 1, \pm i \} = \Z_4. $$ 
Thus we have the isomorphism $(U(1) \times Spin(10))/\Z_4 \cong (E_6)^{\sigma}$.
\vspace{4mm}

{\bf 3.11. Involution $\gamma$ and subgroup $(Sp(1) \times SU(6))/\Z_2$ of $E_6$}\vspace{3mm}

Let the $C$-linear mapping $\gamma : \gJ^C \to \gJ^C$ be the complexification of $\gamma \in G_2 \subset F_4$. Then $\gamma \in E_6$ and $\gamma^2 = 1$.
\vspace{2mm}

We shall study the following subgroup $(E_6)^{\gamma}$ of $E_6$: 
$$
   (E_6)^{\gamma} = \{ \alpha \in E_6 \, | \, \gamma\alpha = \alpha\gamma \}.$$
As in Section 2.11, we use the decomposition
$$ 
      \gJ^C = \gJ(3, \H)^C \oplus (\H^3)^C,  $$
which is the complexification of the decomposition  $\gJ = \gJ(3, \H) \oplus \H^3$. As usual we denote $\gJ(3, \H)$ and $\{ X \in \gJ(3, \H) \, | \, \tr(X) = 0 \}$ by $\gJ_{\sH}$ and $(\gJ_{\sH})_0$, respectively. 
\vspace{2mm}

We consider the embedding $\C = \{x + ye_1 \, | \, x, y \in \R \} \subset \gC$ and, for an element $a = x + ye_1 \in \C,$  we denote by $a'$ the element $x + yi \in C$. Now, we define an $\R$-linear mapping $k : \H \to M(2, C)$ by
$$
           k(a + be_2) = \pmatrix{a' & b' \cr
                                  - \tau b' & \tau a'}, 
            \quad a, b \in \C.$$
The mapping $k$ is naturally extended to $\R$-linear mappings
$$
  k : M(3, \H) \to M(6,C) \quad \mbox{and} \quad k : \H^3 \to M(2, 6, C).$$ 
We adopt the following notations. 
$$
      J = \pmatrix{J & 0 & 0 \cr
                   0 & J & 0 \cr
                   0 & 0 & J} \in M(6,C), \quad
                    J = \pmatrix{0 & 1 \cr
                                 -1 & 0} $$ 
and $M^* = {}^t\ov{M}$, $M \in M(3,\H)$. Then the following four properties hold.
\vspace{2mm}

$\begin{array}{ll}
(1) \; & k(MN) = k(M)k(N), 
\vspace{1mm}\\
(2)       & \tau\,{}^t(k(M)) = k(M^*),
\vspace{1mm}\\
(3)    & J(k(M)) = (\tau(k(M)))J, 
\vspace{1mm}\\
(4)    & k(\a M) = k(\a)k(M), 
\end{array} 
\qquad M, N \in M(3, \H), \a \in \H^3.$
\vspace{2mm}

The $\R$-linear mappings $k : M(3, \H) \to M(6, C)$ and $k : \H^3 \to M(2, 6, C)$ are extended to $C$-linear mappings $k : M(3, \H)^C \to M(6, C)$ and $k : (\H^3)^C \to M(2, 6, C)$ respectively by
\begin{eqnarray*}
    k(M_1 + iM_2) \!\!\! &=& \!\!\!  k(M_1) + ik(M_2), \quad \, M_1, M_2 \in M(3, \H), 
\vspace{1mm}\\
    k(\a_1 + i\a_2) \!\!\! &=& \!\!\! k(\a_1) + ik(\a_2), \qquad \a_1, \a_2 \in \H^3.
\end{eqnarray*}
It is not difficult to see that they satisfy the four properties $(1) \sim (4)$ above. 
\vspace{2mm}

We define a $C$-vector space $\gS(6, C)$ by
$$ 
          \gS(6, C) = \{S \in M(6, C) \, | \, ^tS = - S \} $$
and a $C$-linear mapping $k_J : \gJ(3, \H)^C \to \gS(6, C)$ by
$$
        k_J(M) = k(M)J. $$
Then $k_J$ is well-defined. Indeed, for $M = M_1 + iM_2 \in \gJ(3, \H)^C$, we have 
\begin{eqnarray*}
    ^t(k_J(M)) \!\!\!&=&\!\!\! {}^t((k(M))J) = - J\,{}^t(k(M)) = - J\,{}^t(k(M_1) + ik(M_2))
\vspace{1mm}\\
    \!\!\!&=&\!\!\! - (\tau\,{}^t(k(M_1)) + i\tau\,^t(k(M_2)))J = - (k({M_1}^*) + ik({M_2}^*))J
\vspace{1mm}\\
    \!\!\!&=&\!\!\! - (k(M_1) + ik(M_2))J = - (k(M))J = - k_J(M).
\end{eqnarray*}

Finally, we define Hermitian inner products $\langle S, T \rangle$ in $\gS(6, C)$ and $\langle P, Q \rangle$ in $M(2, 6, $ $C)$ respectively by
$$
       \langle S, T \rangle = \tr((\tau{}\,^tS)T), \quad 
       \langle P, Q \rangle = \tr((\tau{}\,^tP)Q). $$

{\bf Lemma 3.11.1.} (1) $k : M(3, \H)^C \to M(6, C)$, $k : (\H^3)^C \to M(2, 6, C)$ and $k_J : \gJ(3, \H)^C \to \gS(6, C)$ {\it are $C$-linear isomorphisms.}
\vspace{1mm}

(2) $\;\; \langle k_J(M), k_J(N) \rangle = 2\langle M, N \rangle, \qquad M, N \in \gJ(3, \H)^C$, 
\vspace{1mm}

\qquad  
    $\langle k(\a), k(\b) \rangle = 2\langle \a, \b \rangle, \quad \qquad \qquad \a, \b \in (\H^3)^C$.
\vspace{1mm}

(3) $\;\; \det(k_J(M)) = (\det M)^2,  \qquad \qquad M \in \gJ(3, \H)^C$.
\vspace{2mm}

{\bf Proof.} (1) The mapping $k$ is injective. Indeed, if
$$
      k\big((a + be_2) + i(c + de_2)\big) = 0, \quad a, b, c, d \in \C, $$
then
$$
         \pmatrix{a' & b' \cr
                  - \tau b' & \tau a'} 
       + i\pmatrix{c' & d' \cr
                 - \tau d' & \tau c'} = 0. $$
From which, we have $a' = b' = c' = d' = 0$ and so $a = b = c = d = 0$. Now, $k$ is a $C$-linear isomorphism, because $\dim_C(M(3, \H)^C) = 36 = \dim_C(M(6, C))$. The other mapping can be treated analogously. 
\vspace{1mm}

(2)  For $m_k = a_k + b_ke_2$ and $n_k = c_k + d_ke_2, a_k, b_k, c_k, d_k\in \C, k = 1, 2,$ we have, after some calculations, 
$$
          \langle k(m_1 + in_1), k(m_2 + in_2) \rangle 
       = 2\langle m_1 + in_1, m_2 + in_2 \rangle. $$ 

(3) Note that $\det(k_J(M)) = \det(k(M))$ and $\det M \in C$. Since we know that $\det S$ of a skew-symmetric matrix $S$ is expressed in terms of the square of a polynomial with entries $s_{ij}$ in $S$, we can easily see that
$$
\begin{array}{l}
  \quad  \det\pmatrix{ 0 & s_{12} & s_{13} & s_{14} & s_{15} & s_{16} \cr
                  - s_{12} & 0 & s_{23} & s_{24} & s_{25} & s_{26} \cr
                  - s_{13} & - s_{23} & 0 & s_{34} & s_{35} & s_{36} \cr
                  - s_{14} & - s_{24} & - s_{34} & 0 & s_{45} & s_{46} \cr
                  - s_{15} & - s_{25} & - s_{35} & - s_{45} & 0 & s_{56} \cr
                  - s_{16} & - s_{26} & - s_{36} & - s_{46} & -s_{56} & 0} 
\vspace{1mm}\\
   = (s_{12}s_{34}s_{56} - s_{12}s_{35}s_{46} + s_{12}s_{36}s_{45} 
      - s_{13}s_{24}s_{56} + s_{13}s_{25}s_{46} - s_{13}s_{26}s_{45}
\vspace{1mm}\\
   \;\;\; + s_{14}s_{23}s_{56} - s_{14}s_{25}s_{36} + s_{14}s_{26}s_{35} 
      - s_{15}s_{23}s_{46} + s_{15}s_{24}s_{36} - s_{15}s_{26}s_{35}
\vspace{1mm}\\
   \;\;\; + s_{16}s_{23}s_{45} - s_{16}s_{24}s_{35} + s_{16}s_{25}s_{34})^2.
\end{array} $$
Now, since $k_J(M)$ is skew-symmetric, using the above result, we have
$$
\begin{array}{l}
  \det(k_J(M)) = 
     \det\pmatrix{0 & \xi_1 & - n_3 & m_3 & n_2 & \tau(m_2) \cr
          - \xi_1 & 0 & - \tau(m_3) & - \tau(n_3) & - m_2 & \tau(n_2) \cr
          n_3 & \tau(m_3) & 0 & \xi_2 & - n_1 & m_1 \cr
          - m_3 & \tau(n_3) & - \xi_2 & 0 & - \tau(m_1) & - \tau(n_1) \cr
          - n_2 & m_2 & n_1 & \tau(m_1) & 0 & \xi_3 \cr       
          - \tau(m_2) & - \tau(n_2) & - m_1 & \tau(n_1) & - \xi_3 & 0}
\vspace{2mm}\\
  \;\;\;= (\xi_1\xi_2\xi_3 - \xi_1n_1\tau(n_1) - \xi_1m_1\tau(m_1) - n_3\tau(n_3)\xi_3 -
 n_3m_2\tau(n_1) - n_3\tau(n_2)\tau(m_1) 
\vspace{1mm}\\
 \;\;\; - m_3\tau(m_3)\xi_3 + m_3m_2m_1 - m_3\tau(n_2)n_1 - n_2\tau(m_3)\tau(n_1) - n_2\tau(n_3)m_1 - n_2\tau(n_2)\xi_2 
\vspace{1mm}\\
  \;\;\; + \tau(m_2)\tau(m_3)\tau(m_1) - \tau(m_2)\tau(n_3)n_1 - \tau(m_2)m_2\xi_2)^2.
\end{array} $$
On the other hand, we have
$$
\begin{array}{l}
 \det\pmatrix{\xi_1 & m_3 + n_3e_2 & \tau(m_2 + n_2e_2) \vspace{0.5mm}\cr
               \tau(m_3 + n_3e_2) & \xi_2 & m_1 + n_1e_2 \vspace{0.5mm}\cr
                          m_2 + n_2e_2 & \tau(m_1 + n_1e_2) & \xi_3}
\vspace{1mm}\\
   = \xi_1\xi_2\xi_3 + (m_1 + n_1e_2)(m_2 + n_2e_2)(m_3 + n_3e_2)
\vspace{1mm}\\
   \;\;\;+ \tau((m_1 + n_1e_2)(m_2 + n_2e_2)(m_3 + n_3e_2)) - \dsum_{i=1}^3\xi_i(m_i + n_ie_2)\tau(m_i + n_ie_2)
\vspace{1mm}\\
  = \mbox{the contents of in the parenthesis above.}
\end{array} $$  

We now consider a group $E_{6,{\sH}}$ defined by replacing $\gC$ by $\H$ in the definition of the group $E_6$:

$$
     E_{6,{\sH}} = \{ \alpha \in \Iso_C(({\gJ_{\sH}})^C) \, | \, \det(\alpha M) = \det M,
\langle \alpha M, \alpha N \rangle  = \langle M, N \rangle \}.$$

{\bf Lemma 3.11.2.} {\it $E_{6,{\sH}}$ is a connected group of dimension} 35.
\vspace{2mm}

{\bf Proof.} As in the case $E_6$, the group $E_{6,{\sH}}$ contains a subgroup
$$
            F_{4,{\sH}} = \{ \alpha \in E_{6,{\sH}} \, |\, \alpha E = E \}, $$
which is also defined by
$$
    F_{4,{\sH}} = \{ \alpha \in \Iso_{\sR}(\gJ_{\sH}) \, | \, \alpha(M \circ N) = \alpha M \circ \alpha N \} $$
(see Lemma 2.2.4). Moreover, it is isomorphic to the group $Sp(3)/\Z_2$ (Proposition 2.11.1). The Lie algebra $\gge_{6,{\sH}}$ of the group $E_{6,{\sH}}$ is given by
\begin{eqnarray*}
    \gge_{6,{\sH}} \!\!\! &=& \!\!\! \{ \phi \in \Hom_C(({\gJ_{\sH}})^C) \, | \, (\phi M, M, M)=0, \langle \phi M, N \rangle + \langle M, \phi N \rangle = 0 \}\vspace{1mm}\\
    \!\!\! &=& \!\!\! \{ \delta + i\wti{T} \, | \, \delta \in \gf_{4,{\sH}}, T \in (\gJ_{\sH})_0 \},
\end{eqnarray*}
(see Theorem 3.2.1), so that
$$
       \gge_{6,{\sH}} = \gf_{4,{\sH}} \oplus i(\wti{\gJ}_{\sH})_0. $$
Where $\gf_{4,{\sH}} = \{ \delta \in \gge_{6,{\sH}} \, | \, \delta E = 0 \}$ is the Lie algebra of the group $F_{4,{\sH}}$ and is isomorphic to the Lie algebra $\sp(3)$ by the mapping $\varphi_* : \sp(3) \to \gf_{4,{\sH}}$ given by
$$
      \varphi_*(C)M = CM + MC^* = [C, M], \quad M \in \gJ_{\sH} $$
(Proposition 2.11.1) and we see that $\dim\gge_{6,{\sH}} = 21 + 15 = 35$. As in Theorem 3.8.3, we have a homeomorphism
$$
         E_{6,{\sH}}/F_{4,{\sH}} \simeq EIV_{\sH} = \{ X \in ({\gJ_{\sH}})^C \, | \, \det M = 1, \langle M, M \rangle = 3 \} $$
and so we see that the group $E_{6,{\sH}}$ is connected.
\vspace{3mm}

{\bf Proposition 3.11.3.} \quad \qquad $E_{6,{\sH}} \cong SU(6)/\Z_2, \;\; \Z_2 = \{ E,-E \}$.
\vspace{2mm}

{\bf Proof.} Let $SU(6) = \{ A \in M(6, C) \, | \, (\tau\,^tA)A = E, \det A = 1 \}$ and define a mapping $\varphi : SU(6) \to E_{6,{\sH}}$ by
$$
       \varphi(A)M = {k_J}^{-1}(A(k_J(M))\,^tA), \quad M \in ({\gJ_{\sH}})^C.$$
We first have to prove that $\varphi(A) \in E_{6,{\sH}}$. Indeed, we have $(\det(\varphi(A)M))^2 = \det(k_J(\varphi(A)$ $M))$ (Lemma 3.11.1) $= \det(A(k_J(M))^tA) = \det(k_J(M)) = (\det M)^2$ (Lemma 3.11.1), and so $\det(\varphi(A)M) = \pm \det M.$  On the other hand, since $E_{6,{\sH}}$ is connected (Lemma 3.11.2), the sign of $\det(\varphi(A)M)$ is constant, that is, independent of $A$ (assuming that $\det M \neq 0$). Therefore
$$
           \det\,(\varphi(A)M) = \det\,M. $$
Next, again we have
$$
\begin{array}{l}
     2\langle \varphi(A)M, \varphi(A)N \rangle = \langle k_J(\varphi(A)M), k_J(\varphi(A)N) \rangle \;\,\mbox{(Lemma 3.11.1)}
\vspace{1mm}\\
     \quad = \langle A(k_J(M))\,{}^tA, A(k_J(N))\,{}^tA \rangle = \langle k_J(M), k_J(N) \rangle \; \mbox{(because} \; A \in SU(6)) 
\vspace{1mm}\\
 \quad = 2\langle M, N \rangle,
\end{array} $$
so $\varphi(A) \in E_{6,{\sH}}$. Consequently the mapping $\varphi : SU(6) \to E_{6,{\sH}}$ is well-defined. Evidently $\varphi$ is a homomorphism. That $\Ker\varphi = \{ E, - E \} = \Z_2$ can be easily obtained. Since the group $E_{6,{\sH}}$ is connected and $\dim SU(6) = 35 = \dim E_{6,{\sH}}$ (Lemma 3.11.2), $\varphi$ is onto. Thus we have the isomorphism $SU(6)/\Z_2 \cong E_{6,{\sH}}$. 
\vspace{3mm}

{\bf Theorem 3.11.4.} \quad $(E_6)^{\gamma} \cong (Sp(1) \times SU(6))/\Z_2, \;\; \Z_2 = \{(1,E), (-1,-E) \}.$ 
\vspace{2mm}

{\bf Proof.}  We define a mapping $\varphi : Sp(1)\times SU(6) \to (E_6)^\gamma$ by 
$$
\begin{array}{l}
       \varphi(p, A)(M + \a) = {k_J}^{-1}(A(k_J(M))\,^tA) + p\a k^{-1}(\tau \,^tA), $$
\vspace{1mm}\\
\hspace{4cm}
       M + \a \in ({\gJ_{\sH}})^C \oplus (\H^3)^C = \gJ^C.
\end{array} $$
We first need to prove that 
$\varphi(p, A) \in (E_6)^\gamma$. 
\vspace{2mm}

{\bf Claim 1.} \quad ${}^t\varphi(p, A)^{-1} = \tau\varphi(p, A)\tau.$   
\vspace{2mm}

{\bf Proof.} we have
$$
\begin{array}{l}
    2\langle \tau\,{}^t\varphi(p, A)\tau(M + \a), N + \b \rangle \qquad M + \a, N + \b \in \gJ^C
\vspace{1mm}\\
   \qquad = 2\langle M + \a, \varphi(p, A)(N + \b) \rangle 
\vspace{1mm}\\
   \qquad = 2\langle M + \a, {k_J}^{-1}(A(k_J(N))\,^t\!A) + p\b k^{-1}(\tau\,^t\!A) \rangle
\vspace{1mm}\\
   \qquad = 2\langle  M, {k_J}^{-1}(A(k_J(N))\,^t\!A) \rangle + 4\langle \a, p\b k^{-1}(\tau\,^t\!A) \rangle 
\vspace{1mm}\\
   \qquad = \langle k_JM, A(k_J(N))\,^tA \rangle +2\langle  k(\ov{p}\a), (k\b)\tau\,^tA \rangle 
\vspace{1mm}\\
   \qquad = \langle \tau\,{}^tA(k_J(M))\tau A, k_J(N) \rangle + 2\langle k(\ov{p}\a)A, k\b \rangle
\vspace{1mm}\\
   \qquad = 2\langle {k_J}^{-1}(\tau\,^tA(k_J(M))\tau A), N \rangle + 4\langle \ov{p}\a k^{-1}(A), \b \rangle 
\vspace{1mm}\\
   \qquad = 2\langle {k_J}^{-1}(\tau\,^tA(k_J(M))\tau A + \ov{p}\a k^{-1}(A)), N + \b \rangle 
\vspace{1mm}\\
   \qquad = 2\langle  \varphi(\ov{p}, \tau\,^tA)(M + \a), N + \b \rangle,
\end{array} $$
and so we have $\tau\,^t\varphi(p, A)\tau = \varphi(\ov{p}, \tau\,^tA),$ which implies that ${}^t\varphi(p, A)^{-1} = \vspace{2mm}
\tau\varphi(p, A)\tau.$

{\bf Claim 2.} \quad $\varphi(p, A) \in (E_6)^{\gamma}.$
\vspace{2mm}

{\bf Proof.} We will first show that $ \alpha = \varphi(p, A)$ satisfies
$$
        \alpha X \times \alpha Y = {}^t\alpha^{-1}(X \times Y), \quad X, Y \in \gJ^C. $$
Recall the equality
$$
   (M + \a) \times (N + \b) = \Big(M \times N - \dfrac{1}{2}(\a^*\b + \b^*\a) \Big) - \dfrac{1}{2}(\a N + \b M). $$

\noindent The equality $\alpha M \times \alpha N = {}^t\alpha^{-1}(M \times N)$ is evident from $\det(\alpha M) = \det M$ (Proposition 3.11.3). Now, using the equalities
$$
   {k_J}^{-1}(M) = - k^{-1}(MJ) = - k^{-1}(J(\tau M)), \quad 
   \tau k^{-1}(M) = - k^{-1}(JMJ), $$
we have
\begin{eqnarray*}
   (\alpha\a)^*(\alpha\b) \!\!\! &=& \!\!\! (p\a k^{-1}(\tau\,^tA))^*(p\b k^{-1}(\tau \,{}^tA))
\vspace{1mm}\\
                \!\!\! &=& \!\!\! k^{-1}(A)\a^*\b k^{-1}(\tau\,{}^tA), 
\vspace{1mm}\\
    \tau\alpha\tau(\a^*\b) \!\!\! &=& \!\!\! \tau({k_J}^{-1}(A(k_J\tau(\a^*\b))\,{}^tA))
\vspace{1mm}\\
    \!\!\! &=& \!\!\! - \tau k^{-1}(J\tau(A(k(\tau(\a^*\b)))J\,{}^tA))
\vspace{1mm}\\
    \!\!\! &=& \!\!\! - \tau k^{-1}(J\tau(AJ\tau k(\a^*\b)\,{}^tA))
\vspace{1mm}\\
    \!\!\! &=& \!\!\! - \tau k^{-1}(J\tau AJk(\a^*\b)\tau\,{}^tA)
\vspace{1mm}\\
    \!\!\! &=& \!\!\! - k^{-1}(A\tau k(\a^*\b)J\,{}^tAJ)
    = k^{-1}(A)\a^*\b k^{-1}(\tau\,{}^t\!A),
\vspace{1mm}\\
    (\alpha \a)(\alpha N) \!\!\! &=& \!\!\! (p\a k^{-1}(\tau\,{}^tA)){k_J}^{-1}(A(k_J(N))\,{}^t\!A)
\vspace{1mm}\\
    \!\!\! &=& \!\!\! - p\a k^{-1}(\tau\,{}^tA)k^{-1}(A(k(N))J\,{}^tAJ)
\vspace{1mm}\\
    \!\!\! &=& \!\!\! - p\a Nk^{-1}(J\,{}^tAJ) = p\a N\tau k^{-1}(\tau\,{}^tA), \vspace{1mm}\\
    \tau\alpha\tau(\a N) \!\!\! &=& \!\!\! \tau(p\tau(\a N)k^{-1}(\tau\,{}^tA)) = p\a N\tau k^{-1}(\tau\,{}^tA).
\end{eqnarray*}
We have therefore shown that $\varphi(p, A) \in E_6$,  Clearly, $\gamma\varphi(p, A) = \varphi(p, A)\gamma$, so that $\varphi(p, A) \in (E_6)^\gamma$. 
\vspace{2mm}

We will return to the proof of Theorem 3.11.4. Evidently $\varphi$ is a homomorphism. We shall now show that $\varphi$ is onto. Let $\alpha \in (E_6)^{\gamma}$. Since the restriction $\alpha' = \alpha|({\gJ_{\sH}})^C$ of $\alpha$ to $({\gJ_{\sH}})^C$ belongs to $E_{6, {\sH}}$, there exists $A \in SU(6)$ such that $\alpha' = \varphi(A)$ (Proposition 3.11.3). If we put $\beta = \varphi(1,A)^{-1}\alpha$, 
then $\beta|({\gJ_{\sH}})^C = 1$. Hence, by the same argument as in Theorem 2.11.2, there exists $p \in Sp(1)$ such that $\beta = \varphi(p, E)$, and we obtain
$$
    \alpha = \varphi(1,A)\beta = \varphi(1,A)\varphi(p, E) = \varphi(p, A).$$
Therefore $\varphi$ is onto. $\Ker\,\varphi = \{(1, E),(-1,-E) \} = \Z_2$ can be easily obtained. Thus we have the isomorphism $(Sp(1) \times SU(6))/\Z_2 \cong (E_6)^{\gamma}$.
\vspace{4mm}

{\bf 3.12. Involution $\tau\gamma$ and subgroup $Sp(4)/\Z_2$ of $E_6$}
\vspace{3mm}

We consider an involutive complex conjugate linear transformation $\tau\gamma$ of $\gJ^C$, and we shall study the following subgroup $(E_6)^{\tau\gamma}$ of $E_6$:
\begin{eqnarray*}
      (E_6)^{\tau\gamma} \!\!\! &=& \!\!\! \{ \alpha \in E_6 \, | \, \tau\gamma\alpha = \alpha\tau\gamma \} \\
      \!\!\! &=& \!\!\! \{ \alpha \in E_6 \, | \, \gamma\lambda(\alpha)\gamma =  \alpha \} = (E_6)^{\lambda\gamma}.
\end{eqnarray*}
For this end, we consider $\R$-vector subspaces $(\gJ^C)_{\tau\gamma}$ and $(\gJ^C)_{-\tau\gamma}$ of $\gJ^C$, which are eigenspaces of $\tau\gamma$, respectively by
\begin{eqnarray*}
   (\gJ^C)_{\tau\gamma} \!\!\! &=& \!\!\! \{ X \in \gJ^C \, | \, \tau\gamma X = X \}
\vspace{1mm} \\
     \!\!\! &=& \!\!\! \Bigg\{ \pmatrix{\mu_1 & m_3 & \ov{m}_2 \cr 
                                   \ov{m}_3 & \mu_2 & m_1 \cr
                                   m_2 & \ov{m}_1 & \mu_3} 
                + i\pmatrix{0 & a_3e_4 & -a_2e_4 \cr
                            -a_3e_4 & 0 & a_1e_4 \cr
                            a_2e_4 & -a_1e_4 & 0} \, \Bigg| \, 
            \begin{array}{l}
              \mu_i \in \R \\
              m_i, a_i \in \H
            \end{array} \Bigg\} 
\vspace{1mm}\\
     \!\!\! &=& \!\!\! \{ M + iF(\a e_4) \, | \,  M \in \gJ(3, \H), \a \in \H^3 \}\vspace{1mm}\\
      \!\!\! &=& \!\!\! \gJ_{\sH} \oplus i\H^3, 
\vspace{1mm}\\
(\gJ^C)_{-\tau\gamma} \!\!\! &=& \!\!\! \{ X \in \gJ^C \, | \, \tau\gamma X = - X \} 
\vspace{1mm}\\
     \!\!\! &=& \!\!\! \Bigg\{ i\pmatrix{\mu_1 & m_3 & \ov{m}_2 \cr 
                                   \ov{m}_3 & \mu_2 & m_1 \cr
                                   m_2 & \ov{m}_1 & \mu_3} 
                + \pmatrix{0 & a_3e_4 & -a_2e_4 \cr
                            -a_3e_4 & 0 & a_1e_4 \cr
                            a_2e_4 & -a_1e_4 & 0} \, \Bigg| \, 
            \begin{array}{l}
              \mu_i \in \R \\
              m_i, a_i \in \H
            \end{array} \Bigg\} 
\vspace{1mm}\\
     \!\!\! &=& \!\!\! \{ iM + F(\a e_4) \, | \,  M \in \gJ(3, \H), \a \in \H^3 
\}
\vspace{1mm}\\
      \!\!\! &=& \!\!\! i\gJ_{\sH} \oplus \H^3, 
\end{eqnarray*}              
The spaces $(\gJ^C)_{\tau\gamma}$ and $(\gJ^C)_{-\tau\gamma}$ are invariant under the action of $(E_6)^{\tau\gamma}$ and we have the decomposition of $\gJ^C$:
$$
         \gJ^C = (\gJ^C)_{\tau\gamma} \oplus (\gJ^C)_{-\tau\gamma} = 
                (\gJ^C)_{\tau\gamma} \oplus i(\gJ^C)_{\tau\gamma}. $$
In particular, $\gJ^C$ is the complexification of $(\gJ^C)_{\tau\gamma}$: $\gJ^C = ((\gJ^C)_{\tau\gamma})^C$. 
\vspace{2mm}

In the $\R$-vector space
$$
         \gJ(4, \H) = \{ P \in M(4, \H) \, | \, P^* = P \}, $$ 
we define the Jordan multiplication $P \circ Q$ and an inner product $(P, Q)$ respectively by
$$
       P \circ Q = \frac{1}{2}(PQ + QP), \quad (P, Q) = \tr(P \circ Q).$$
The group $Sp(4)$ acts on $\gJ(4, \H)$ by the mapping $\mu : Sp(4) \times \gJ(4, \H) \to \gJ(4, \H)$, $\mu(A, P) = APA^*$. Then we have
$$
\begin{array}{l}
      A(P \circ Q)A^* = APA^* \circ AQA^*, 
\vspace{1mm}\\
     (APA^*, AQA^*) = (P, Q),
\end{array}
\quad A \in Sp(4), P, Q \in \gJ(4, \H).$$
The quaternion projective space $\H \!P_3$ is defined by
\begin{eqnarray*}
      \H \!P_3 \!\!\! &=& \!\!\! \{ P \in \gJ(4,\H) \, | \,  P^2 = P,\tr(P) = 1 \} \vspace{1mm}\\
     \!\!\! &=& \!\!\! \Big\{ AE_1A^* \, | \, A \in Sp(4), 
                E_1 = \diag(1, 0, 0, 0) \in M(4, \H) \Big\}.
\end{eqnarray*}
Finally, in the complexification $\gJ(4, \H)^C$ of $\gJ(4, \H)$, we extend naturally the Jordan multiplication $P \circ Q$ and the inner product $(P, Q)$ and further define a Hermitian inner product $\langle P, Q \rangle $ by
$$
                   \langle P, Q \rangle  = (\tau P, Q).$$
The action of $Sp(4)$ on $\gJ(4, \H)$ is also naturally extended to $\gJ(4, \H)^C$: 
$$
   A(X_1 + iX_2)A^{*} = AX_1A^{*} + iAX_2A^{*}, \quad A \in Sp(4),X_1, X_2 \in \gJ(4, \H).$$
Then we have
$$
    \langle APA^{*}, AQA^{*} \rangle = \langle P, Q \rangle, \quad P,Q \in \gJ(4, \H)^C.$$
We denote by $\gJ(4, \H)_0$ the space $\{ P \in \gJ(4, \H) \, | \, \tr(P) = 0 \}$ and by ${\gJ(4, \H)_0}^C$ its complexification.
\vspace{3mm}

{\bf Definition.} We define a $C$-linear mapping $g : \gJ^C \to {\gJ(4, \H)_0}^C$ by
$$
        g(M + \a) = \pmatrix{\displaystyle{\frac{1}{2}}\tr(M) & i\a \cr
                             i\a^* & M - \displaystyle{\frac{1}{2}}\tr(M)E},
\quad M + \a \in (\gJ_{\sH})^C \oplus (\H^3)^C = \gJ^C. $$
The restriction of $g$ to $(\gJ^C)_{\tau\gamma}$ is given by
$$
      g(M + i\a) = \pmatrix{\tr(M) & -\a \cr
                            -\a^* & M} 
                 - \displaystyle{\frac{1}{2}}\tr(M)E, \quad M + i\a \in \gJ_{\sH} \oplus i\H^3 = (\gJ^C)_{\tau\gamma}. $$
Note that the mapping $g$ in the above definition is the complexification of this 
\vspace{3mm}
restriction.

{\bf Lemma 3.12.1.} {\it The mapping $g : \gJ^C \to {\gJ(4, \H)_0}^C$ is a $C$-linear isomorphism and satisfies}
$$
\begin{array}{l}
        gX \circ gY = g(\gamma(X \times Y)) + \displaystyle{\frac{1}{4}}(\gamma X, Y)E,
\vspace{1mm}\\
        (gX, gY) = (\gamma X, Y), \;
\end{array}
\quad X, Y \in \gJ^C. $$
{\it Moreover, $g$ is an isometry with respect to the inner product $\langle X,  Y \rangle$}:
$$
     \langle gX, gY \rangle = \langle X, Y \rangle, \quad X, Y \in \gJ^C.$$
{\it The restriction of $g$ to $(\gJ^C)_{\tau\gamma}$ induces an $\R$-linear isomorphism} $g : (\gJ^C)_{\tau\gamma} \to 
\vspace{2mm}
\gJ(4, \H)_0$.

{\bf Proof.} It is not difficult to see that $g$ is well-defined and that it is injective. Since $\dim_C\gJ^C = 27 = \dim_C{\gJ(4, \H)_0}^C$, $g$ is a $C$-linear isomorphism.  Now, for $X = M + \a$, $Y = N + \b \in (\gJ_{\sH})^C \oplus (\H^3)^C = \gJ^C$, we have
$$
\begin{array}{l}
   g(\gamma(X \times Y)) = g(\gamma((M + \a) \times (N + \b)))
\vspace{1mm}\\
     \qquad \; \;  = g((M - \a) \times (N - \b)) 
\vspace{1mm}\\
     \qquad \; \;  = g((M \times N - \dfrac{1}{2}(\a^*\b + \b^*\a)) + \dfrac{1}{2}(\a N + \b M))
\end{array} $$
$$                           
  = \pmatrix{\dfrac{1}{2}\tr(M \times N) - \dfrac{1}{2}(\a, \b) & \dfrac{i}{2}(\a N + \b M) \cr
             \dfrac{i}{2}(\a N + \b M)^* & M \times N - \dfrac{1}{2}(\a^*\b + \b^*\a) - \dfrac{1}{2}(\tr(M \times N) - (\a, \b))E} $$
$$
\begin{array}{l}
\qquad \quad = \cdots 
\vspace{1mm}\\
\qquad \quad = g(M + \a) \circ g(N + \b) - \Big(\displaystyle{\frac{1}{4}}(M, N) - \dfrac{1}{2}(\a, \b) \Big)E 
\vspace{1mm}\\
\quad \qquad = g(M + \a) \circ g(N + \b) - \dfrac{1}{4}(\gamma(M + \a), (N + \b))E 
\vspace{1mm}\\
\qquad \quad = gX \circ gY - \displaystyle{\frac{1}{4}}(\gamma X, Y)E.
\end{array} $$
Thus the first equality $gX \circ gY = g(\gamma(X \times Y)) + \displaystyle{\frac{1}{4}}(\gamma X, Y)E$ is proved. Taking the traces of the both sides, we have
$$
\displaylines{\hfill 
        (gX, gY) = (\gamma X, Y), \quad X, Y \in \gJ^C.
\hfill \mbox{(i)}} $$
It is easily seen that the restriction $g : (\gJ^C)_{\tau\gamma} \to \gJ(4, \H)_0$ of $g$ is an $\R$-isomorphism. Finally, if we note that $(\gamma X, Y) = \langle X, Y \rangle$ for $X, Y \in (\gJ^C)_{\tau\gamma}$, then we can easily show that
$$
   \langle g(X_1 + iX_2), g(Y_1+iY_2) \rangle = 
        \langle X_1 + iX_2, Y_1 + iY_2 \rangle, \quad X_i, Y_i \in (\gJ^C)_{\tau\gamma} $$
from (i). Thus Lemma 3.12.1 is proved.
\vspace{3mm}

{\bf Theorem 3.12.2.} \qquad \quad $(E_6)^{\tau\gamma} \cong Sp(4)/\Z_2, \; \; \Z_2 = \{E, - E \}$.
\vspace{2mm}

{\bf Proof.} We define a mapping $\varphi : Sp(4) \to (E_6)^{\tau\gamma}$ by
$$
       \varphi(A)X = g^{-1}(A(gX)A^{*}), \quad X \in \gJ^C.$$
We first have to prove that $\varphi(A) \in (E_6)^{\tau\gamma}$. Let $Z = \varphi(A)X$ 
and use Lemma 3.12.1, then we have
$$
\begin{array}{l}
 3\det(\varphi(A)X) = 3\det Z = (Z \times Z, Z) = (g(\gamma(Z \times Z)), gZ)
\vspace{1mm}\\
\qquad \; \; = \Big(gZ \circ gZ - \displaystyle{\frac{1}{4}}(\gamma Z, Z)E, gZ\Big) 
\vspace{1mm}
\\

\qquad \; \; = \Big(gZ \circ gZ - \displaystyle{\frac{1}{4}}(gZ, gZ)E, 
\vspace{-1mm}
gZ \Big)
\end{array} $$

$$
   = \Big(A(gX)A^{*} \circ A(gX)A^{*} - \displaystyle{\frac{1}{4}}(A(gX)A^{*}, A(gX)A^{*})E, A(gX)A^{*}\Big) 
\vspace{-1mm}$$

\qquad \qquad \qquad \quad 
   $= (gX \circ gX - \displaystyle{\frac{1}{4}}(gX, gX)E, gX)$ 
\vspace{1mm}

\qquad \qquad \qquad \quad 
    $= (gX \circ gX - \displaystyle{\frac{1}{4}}(\gamma X, X)E, gX)$ 
\vspace{1mm}

\qquad \qquad \qquad \quad 
    $= (g(\gamma(X \times X)), gX) = (X \times X, X) = 3\det X$,
\vspace{1mm}

\noindent and
$$
\begin{array}{l}
   \langle \varphi(A)X, \varphi(A)Y \rangle = \langle g\varphi(A)X, g\varphi(A)Y \rangle 
\vspace{1mm}\\
\qquad \; = \langle A(gX)A^{*}, A(gY)A^{*} \rangle = \langle gX, gY \rangle = \langle X, Y \rangle .
\end{array} $$
Hence $\varphi(A) \in E_6$. To prove $\varphi(A) \in (E_6)^{\tau\gamma}$, namely, $\tau\gamma\varphi(A) = \varphi(A)\tau\gamma$, it is sufficient to show that
$$
      \tau\gamma\varphi(A)\tau\gamma X = \varphi(A)X, \quad X \in (\gJ^C)_{\tau\gamma},$$
since $\gJ^C = ((\gJ^C)_{\tau\gamma})^C$. However this is evident. Indeed, if $X \in (\gJ^C)_{\tau\gamma}$, then $gX \in \gJ(4, \H)_0$, so that $\varphi(A)X \in (\gJ^C)_{\tau\gamma}$. Hence $\tau\gamma\varphi(A)\tau\gamma X = \tau\gamma\varphi(A)X = \varphi(A)X.$  Evidently $\varphi$ is a homomorphism. We shall show that $\varphi$ is onto. For $\alpha \in (E_6)^{\tau\gamma}$, we have
$$
        (g(\alpha E))^2 = g(\alpha E) + \frac{3}{4}E.$$
Certainly,
$$
\begin{array}{l}
    (g(\alpha E))^2 = g(\alpha E) \circ g(\alpha E) 
        = g(\gamma(\alpha E \times \alpha E)) + \displaystyle{\frac{1}{4}}(\gamma\alpha E, \alpha E)E \; \mbox{(Lemma 3.12.1)}
\vspace{1mm}\\ 
   \qquad = g(\gamma\tau\alpha\tau E) + \displaystyle{\frac{1}{4}}\langle \tau\gamma\alpha E, \alpha E \rangle E
        = g(\alpha\tau\gamma E) + \dfrac{1}{4}\langle \alpha\tau\gamma E, \alpha E \rangle E
\vspace{1mm}\\
   \qquad = g(\alpha E) + \dfrac{1}{4}\langle \alpha E, \alpha E \rangle E 
        = g(\alpha E) + \displaystyle{\frac{1}{4}}\langle E, E \rangle E 
        = g(\alpha E) + \displaystyle{\frac{3}{4}}E.
\end{array}$$
Now, let
$$
        P = \frac{1}{4}(2g(\alpha E) + E), $$
then  we have $P \in \gJ(4, \H)$ and $P^2 = P$, $\tr(P) = 1$, that is, $P \in 
\H \!P^3$. Indeed,
\begin{eqnarray*}
     P^2 \!\!\! &=& \!\!\! \frac{1}{16}(4(g(\alpha E))^2 + 4g(\alpha E) + E ) = \frac{1}{4}(2g(\alpha E) + E) = P,
\vspace{1mm}\\
     \tr(P) \!\!\! &=& \!\!\! \frac{1}{4}\tr(2g(\alpha E)) + \frac{1}{4}\tr(E) = 0 + 1 = 1.
\end{eqnarray*}
Hence there exists $A \in Sp(4)$ such that
$$
     P = AE_1A^{*}. $$
Since $gE = 2E_1 - \displaystyle{\frac{1}{2}}E$, we have
\begin{eqnarray*}
     \varphi(A)E \!\!\! &=& \!\!\! g^{-1}(A(gE)A^{*}) = g^{-1}\Big(A\Big(2E_1 - \frac{1}{2}E\Big)A^{*}\Big) = g^{-1}\Big(2AE_1A^{*} - \frac{1}{2}E \Big)
\vspace{1mm}\\
     \!\!\! &=& \!\!\! g^{-1}(g(\alpha E)) = \alpha E.
\end{eqnarray*}
Putting $\beta =\varphi(A)^{-1}\alpha$, we have $\beta E = E$, and so $\beta \in F_4$ (Theorem 3.7.1). Also $\beta$ satisfies $\tau\gamma\beta = \beta\tau\gamma$ and $\tau\beta = \beta\tau$ (Theorem 3.7.1), hence $\gamma\beta = \beta\gamma$, and therefore $\beta \in (F_4)^\gamma$. From Theorem 2.11.2, there exist $p \in Sp(1)$ and $D \in Sp(3)$ such that
$$
      \beta(M + \a) = DMD^* + p\a D^*, \quad M + \a \in \gJ_{\sH} \oplus \H^3 = \gJ. $$Let $B = \pmatrix{p & 0 \cr
                  0 & D}$, then $B \in Sp(4)$ and we have
$$
             \beta = \varphi(B).$$
Certainly, for $M + \a \in (\gJ_{\sH})^C \oplus (\H^3)^C = \gJ^C,$
\begin{eqnarray*}
    \varphi(B)(M + \a) \!\!\! &=& \!\!\! g^{-1}(Bg(M + \a)B^{*}) 
\vspace{1mm}\\
    \!\!\! &=& \!\!\! g^{-1}\Big(\pmatrix{p & 0 \cr
                                          0 & D}
       (\pmatrix{\dfrac{1}{2}\tr(M) & i\a \cr
                 i\a^* & M - \dfrac{1}{2}\tr(M)E}
        \pmatrix{\ov{p}& 0 \cr
                 0 & D^*} \Big) 
\vspace{1mm}\\
   \!\!\! &=& \!\!\! g^{-1}\Big(
    \pmatrix{\dfrac{1}{2}\tr(M) & ip\a D^* \cr
             iD\a^*\ov{p} & DMD^* - \dfrac{1}{2}\tr(M)E} \Big)
\vspace{1mm}\\
    \!\!\! &=& \!\!\! DMD^* + p\a D^* = \beta(M + \a). 
\end{eqnarray*}
Hence we have
$$
    \alpha = \varphi(A)\beta = \varphi(A)\varphi(B) = \varphi(AB), \quad AB \in Sp(4), $$
so that $\varphi$ is onto. $\Ker\varphi = \{ E,- E \} = \Z_2$ can be easily obtained. Thus we have the isomorphism $Sp(4)/\Z_2 \cong (E_6)^{\tau\gamma}$.
\vspace{4mm}

{\bf 3.13. Automorphism $w$ of order 3 and subgroup $(SU(3) \times SU(3) \times SU(3))/\Z_3$ of $E_6$}
\vspace{3mm}

Let the $C$-linear mapping $w : \gJ^C \to \gJ^C$ be the complexification of $w \in G_2 \subset F_4$ of Section 2.12. Then $w \in E_6$ and $w^3 = 1$. 
\vspace{2mm}

We shall study the following subgroup $(E_6)^w$ of $E_6$:
$$
           (E_6)^w = \{\alpha \in E_6 \, | \, w\alpha = \alpha w \}. $$

As in Section 2.12, we identify
$$
             \gJ(3, \C)^C \oplus M(3, \C)^C = \gJ^C. $$
For convenience, we denote $\gJ(3, \C)$ and $\{X \in \gJ(3, \C) \, | \, \tr(X) = 0 \}$ by $\gJ_{\sC}$ and $(\gJ_{\sC})_0$, respectively.
\vspace{2mm}

The group $E_{6,{\sC}}$ is defined to be obtained by replacing $\gJ^C$ with $(\gJ_{\sC})^C$ in the definition of the group $E_6$:
$$
    E_{6,{\sC}} = \{ \alpha \in \Iso_C((\gJ_{\sC})^C) \, | \, \det\,(\alpha X) = \det\,X, \langle \alpha X, \alpha Y \rangle = \langle X, Y \rangle \}. $$
As in the case $E_6$, the group $E_{6,{\sC}}$ contains a subgroup
$$
             F_{4,{\sC}} = \{\alpha \in E_{6, {\sC}} \, | \, \alpha E = E \}, $$
which is also defined by the group $F_{4,{\sC}} = \{\alpha \in \Iso_{\sR}(\gJ_{\sC}) \, | \, \alpha(X \circ Y) = \alpha X \circ \alpha Y \}$, moreover, it is isomorphic to the group $(SU(3)/\Z_3)\cdot\Z_2$ (Proposition 2.12.1). The Lie algebra $\gge_{6,{\sC}}$ of the group $E_{6,{\sC}}$ is
\begin{eqnarray*}
   \gge_{6,{\sC}} \!\!\! &=& \!\!\! \{\phi \in \Hom_C((\gJ_{\sC})^C) \, | \, (\phi X, X, X) = 0, \langle \phi X, Y \rangle + \langle X, \phi Y \rangle = 0 \}
\vspace{1mm}\\
              \!\!\! &=& \!\!\! \{\delta + i\wti{T} \, | \, \delta \in \gf_{4,{\sC}}, T \in (\gJ_{\sC})_0 \}
\end{eqnarray*}
(Theorem 3.2.4), where $\gf_{4,{\sC}} = \{\delta \in \gge_{6,{\sC}} \, | \, \delta E = 0 \}$ is the Lie algebra of the group $F_{4,{\sC}}$. In particular, the dimension of $\gge_{6,{\sC}}$ is
$$
        \dim(\gge_{6,{\sC}}) = 8 + 8 = 16. $$
As in Theorem 3.8.3, we see that the space 
$$
          EIV_{\sC} = \{X \in (\gJ_{\sC})^C \, | \, \det X = 1, \langle X, X \rangle = 3 \}$$
is connected and we have the homeomorphism
$$
             E_{6,{\sC}}/F_{4,{\sC}} \simeq EIV_{\sC}. $$

{\bf Lemma 3.13.1.} $E_{6,{\sC}}$ {\it has at most two connected components} (in reality has two connected components (Proposition 3.13.4)).
\vspace{2mm}

{\bf Proof.} From the exact sequence $\pi_0(F_{4,{\sC}}) \to \pi_0(E_{6,{\sC}}) \to \pi_0(EIV_{\sC})$, that is, $\Z_2 \to \pi_0(E_{6,{\sC}}) \to 0$ (Proposition 2.12.1), we see that $\pi_0(E_{6,{\sC}})$ is 0 or $\Z_2$.
\vspace{2mm}

We define mappings $h : \C \oplus \C \to \C^C$ and $h : M(3, \C) \oplus M(3, \C) \to M(3, \C)^C$ respectively by
$$
\begin{array}{ll}
  \;\,h(a, b) = \dfrac{a + b}{2} + i\dfrac{a - b}{2}e_1 = \iota a + \ov{\iota}b, & {}
\vspace{-2mm}\\
   {} & \quad \iota = \dfrac{1 + ie_1}{2}.
\vspace{-2mm}\\
  h(A, B) = \dfrac{A + B}{2} + i\dfrac{A - B}{2}e_1 = \iota A + \ov{\iota}B, & {}
\end{array}$$

{\bf Lemma 3.13.2.} {\it The mappings $h : \C \oplus \C \to \C^C$ and $h : M(3, \C) \oplus M(3, \C) \to M(3, \C)^C$ satisfy the following four conditions.}
\vspace{1mm}

(1) {\it Both are $\C$-linear isomorphisms, that is, they are injective and satisfy}
\vspace{1mm}

\quad $h(a, b) + h(a', b') = h(a + a', b + b'), \;\; h(ca, cb) = ch(a, b), c \in \C,$ 
\vspace{1mm}

\quad $h(A, B) + h(A', B') = h(A + A', B + B'), \;\; h(cA, cB) = h(A, B), c \in \C.$ 
\vspace{1mm}

(2) $h(a, b)h(a', b') = h(aa', bb'), \;\; h(A, B)h(A', B') = h(AA', BB').$
\vspace{1mm}

(3) $\tau h(a, b) = h(b, a), \;\; \ov{h(a, b)} = h(\ov{b}, \ov{a}). $
\vspace{1mm}

\quad $\tau h(A, B) = h(B, A), \;\; \ov{h(A, B)} = h(\ov{B}, \ov{A}), \;\; h(A, B)^* = h(B^*, A^*).$
\vspace{1mm}

(4) $\det\,(h(A, B)) = h(\det\,A, \det\,B)$.
\vspace{2mm}

{\bf Proof.} It is easy to prove, noting that $\iota^2 = \iota, \ov{\iota}^2 = \ov{\iota}, \iota + \ov{\iota} = 1$. 
\vspace{3mm}

{\bf Lemma 3.13.3.} \qquad \quad $\gge_{6,{\sC}} \cong \su(3) \oplus \su(3)$.
\vspace{2mm}

{\bf Proof.} The mapping $\phi_{\sC}: \su(3) \oplus \su(3) \to \gge_{6, {\sC}}$, 
$$
      \phi_{\sC}(C, D)X = h(C, D)X + Xh(C, D)^*, \quad X \in (\gJ_{\sC})^C $$
gives an isomorphism as Lie algebras. This is the direct consequence of the following Proposition 3.13.4, so we will omit its proof here.
\vspace{2mm}

We define an action of the group $\Z_2 = \{1, \epsilon \}$ on $SU(3) \times SU(3)$ by
$$
              \epsilon(A, B) = (\ov{B}, \ov{A}), $$
and let $(SU(3) \times SU(3))\cdot\Z_2$ be the semi-direct product of the groups $SU(3) \times SU(3)$ and $\Z_2$ under this action.
\vspace{3mm}

{\bf Proposition 3.13.4.} $E_{6,{\sC}} \cong ((SU(3) \times SU(3))/\Z_3)\cdot\Z_2, \; \Z_3 = \{(E, E), (\omega_1E, $ $\omega_1E), ({\omega_1}^2E, {\omega_1}^2E) \}, \omega_1 = - \dfrac{1}{2} + \dfrac{\sqrt{3}}{2}e_1$.
\vspace{2mm}

{\bf Proof.} We define a mapping $\varphi : (SU(3) \times SU(3))\cdot\Z_2 \to E_{6,{\sC}}$ by
$$
\begin{array}{ll}
     \varphi((A, B), 1)X = h(A, B)Xh(A, B)^*, & {} 
\vspace{-1mm}\\ 
     {} & X \in (\gJ_{\sC})^C.
\vspace{-1mm}\\
    \varphi((A, B), \epsilon)X = h(A, B)\ov{X}h(A, B)^*, 
\end{array}$$
First we have to show that $\alpha = \varphi((A, B), 1) \in E_{6, {\sC}}$. Using $\det(h(A, B)) = h(\det A, $ $\det B)$ (Lemma 3.13.2.(4)) $= h(1, 1) = 1$ and $\tau h(A, B)^*h(A, B) = h(A^*, B^*)$ $h(A, B)$ (Lemma 3.13.2.(3)) $ = h(A^*A, B^*B) = h(E, E) = E$, we have
$$
\begin{array}{l} 
   \det\,(\alpha X) = (\det\,(h(A, B)))(\det\,X)(\det\,(h(A, B)^*)) = \det\,X, 
\vspace{1mm}\\
   \langle \alpha X, \alpha Y \rangle = \langle h(A, B)Xh(A, B)^*, h(A, B)Yh(A, B)^* \rangle
\vspace{1mm}\\
\qquad \qquad \;
      =(\tau h(A, B)\tau X\tau h(A, B)^*, h(A, B)Yh(A, B)^*)
\vspace{1mm}\\
\qquad \qquad \;
  = (\tau X, Y) = \langle X, Y \rangle. 
\end{array} $$
Hence $\alpha \in E_{6,{\sC}}$. Since $\varphi((E, E), \epsilon) = \epsilon \in G_{2, {\sC}}\,(= \mbox{Aut}(\C))\, \subset F_{4, {\sC}} \subset E_{6,{\sC}}$, 
we also have $\varphi((A, B), \epsilon) = \varphi((A, B),1)\varphi((E, E), \epsilon) \in E_{6,{\sC}}$. Next, we shall show that $\varphi$ is a homomorphism. Indeed, for instance,
\vspace{2mm}

\qquad \qquad 
  $\varphi((A, B), \epsilon)\varphi((C, D), 1)X$
\vspace{1mm}

\qquad \qquad \qquad 
  $= \varphi((A, B), \epsilon)(h(C, D)Xh(C, D)^*)$
\vspace{1mm}

\qquad \qquad \qquad 
  $= h(A, B)\ov{h(C, D)Xh(C, D)^*}h(A, B)^*$
\vspace{1mm}

\qquad \qquad \qquad 
   $= h(A, B)h(\ov{D}, \ov{C})\ov{X}h(\ov{D}, \ov{C})^*h(A, B)^* \;\;\mbox{(Lemma 3.13.3.(3))}$
\vspace{1mm}

\qquad \qquad \qquad 
   $= h(A\ov{D}, B\ov{C})\ov{X}h(A\ov{D}, B\ov{C})^* = \varphi((A\ov{D}, B\ov{C}), \epsilon)X$
\vspace{1mm}

\qquad \qquad \qquad 
   $= \varphi((A, B)\epsilon(C, D), \epsilon)X, \quad X \in (\gJ_{\sC})^C$.
\vspace{2mm}

\noindent That $\Ker\,\varphi = \{(E, E), \,(\omega_1E, \omega_1E), \,({\omega_1}^2E, {\omega_1}^2E)\} \, \times \,1 = \Z_3 \times \,1$ can be easily obtained. In particular, $\Ker\,\varphi$ is discrete. Hence $\varphi$ induces an injective homomorphism 
$$
     \varphi_* : \su(3) \oplus \su(3) \to \gge_{6,{\sC}}. $$
In particular, $\dim(\su(3) \oplus \su(3)) = \dim(\gge_{6,{\sC}})$, so $\varphi_*$ is an isomorphism ($\varphi_*$ coincides with $\phi_{\sC}$ of Lemma 3.13.3). Hence, $\varphi$ induces the surjection $\varphi : SU(3) \times SU(3) \to (E_{6,{\sC}})_0$ (which denotes the connected component of $E_{6, {\sC}}$ containing the identity 1). However $\epsilon = \varphi((E, E), \epsilon) \not\in (E_{6,{\sC}})_0$. Certainly, for any $A, B \in SU(3)$, 
$$
       h(A, B)Xh(A, B)^* = \ov{X}, \quad X \in (\gJ_{\sC})^C $$
does not hold. Therefore $E_{6,{\sC}}$ has just two connected components (Lemma 3.13.1). Consequently, $\varphi : (SU(3) \times SU(3))\cdot\Z_2 \to E_{6,{\sC}}$ is onto and we have the isomorphism $((SU(3) \times SU(3))/\Z_3)\cdot\Z_2 \cong E_{6,{\sC}}$.
\vspace{3mm}

{\bf Theorem 3.13.5.} $(E_6)^w \cong (SU(3) \times SU(3) \times SU(3))/\Z_3, \; \Z_3 = \{(E, E, E), $ $(\omega_1E, \omega_1E, \omega_1E),({\omega_1}^2E, {\omega_1}^2E, {\omega_1}^2E) \}, \omega_1 = - \dfrac{1}{2} + \dfrac{\sqrt{3}}{2}e_1.$
\vspace{2mm}

{\bf Proof.} We define a mapping $\varphi : SU(3) \times SU(3) \times SU(3) \to (E_6)^w$ by
$$
\begin{array}{l}
    \varphi(P, A, B)(X + M) = h(A, B)Xh(A, B)^* + PM\tau h(A, B)^*,
\vspace{1mm}\\
\qquad \qquad \qquad \qquad
   X + M \in (\gJ_{\sC})^C \oplus M(3, \C)^C = \gJ^C.
\end{array} $$
We first have to show that $\alpha = \varphi(P, A, B) \in (E_6)^w$, that is, $\alpha$ leaves
\begin{eqnarray*}
  (X + M) \times (Y + N) \!\!\! &=& \!\!\! (X \times Y - \dfrac{1}{2}(M^*N + N^*M)) - \dfrac{1}{2}(MY + NX + \ov{M \times N}),
\vspace{1mm}\\
    \langle X + M, Y + N \rangle \!\!\! &=& \!\!\! \langle X, Y \rangle + 2\langle M, N \rangle
\end{eqnarray*}
invariant, namely, $\alpha$ satisfies
\begin{eqnarray*}
   \tau\alpha\tau((X + M) \times (Y + N)) \!\!\!&=&\!\!\! \alpha(X + M) \times \alpha(Y + N),
\vspace{1mm}\\
   \langle \alpha(X + M), \alpha(Y + N) \rangle \!\!\!&=&\!\!\! \langle X + M, Y + M \rangle.
\end{eqnarray*}
Now, for $X + M, Y + N \in (\gJ_{\sC})^C \oplus M(3, \C)^C = \gJ^C$, the relations
$$
     \tau\alpha\tau(X \times Y) = \alpha X \times \alpha Y, \quad \langle \alpha X, \alpha Y \rangle = \langle X, Y \rangle $$
are already shown in Proposition 3.13.4. Next, 
$$
\begin{array}{l}
   (PM\tau h(A, B)^*)^*(PN\tau h(A, B)^*) = \tau h(A, B)M^*P^*PN\tau h(A, B)^*
\vspace{1mm}\\
\qquad \qquad \qquad
   = \tau(h(A, B)\tau(M^*N)h(A, B)^*),
\vspace{1mm}\\
   (PM\tau h(A, B)^*)(h(A, B)Yh(A, B)^*) = \tau(P\tau(MY)\tau h(A, B)^*),
\vspace{1mm}\\
   \ov{(PM\tau h(A, B)^*) \times (PN\tau h(A, B)^*)} = \ov{^t\wti{P}(M \times N)\tau\,{}^t(h(A, B)^*)^{\sim}}
\vspace{1mm}\\
\qquad \qquad \qquad
    = P(\ov{M \times N})h(A, B)^* = \tau(P(\ov{M \times N})\tau h(A, B)^*). 
\end{array} $$
Furthermore, we have
$$
\begin{array}{l}
   \langle \alpha X, \alpha M \rangle = 0 = \langle X, M \rangle,
\vspace{1mm}\\
   \langle \alpha M, \alpha N \rangle = (\tau\alpha M, \alpha N) = (P\tau Mh(A, B)^*, PN\tau h(A, B)^*)
\vspace{1mm}\\
\qquad \qquad \;\;\,
   = (\tau M, N) = \langle M, N \rangle, 
\end{array} $$
and so $\alpha \in E_6$. Clearly, $w\alpha = \alpha w$, so that $\alpha \in (E_6)^w$. Evidently $\varphi$ is a homomorphism. We shall now show that $\varphi$ is onto. Let $\alpha \in (E_6)^w$. The restriction $\alpha'$ of $\alpha$ to $(\gJ^C)_w = \{X \in \gJ^C \, | \, wX = X \} = (\gJ_{\sC})^C$ belongs to $E_{6,{\sC}}$: $\alpha' \in E_{6,{\sC}}$. Hence, there exist $A, B \in SU(3)$ such that
$$
    \alpha X = h(A, B)Xh(A, B)^* \quad \mbox{or} \quad \alpha X = h(A, B)\ov{X}h(A, B)^*, \quad X \in (\gJ_{\sC})^C $$
(Proposition 3.13.4). In the former case, let $\beta = \varphi(E, A, B)^{-1}\alpha$, then $\beta|(\gJ_{\sC})^C = 1$. Hence $\beta \in G_2$, moreover $\beta \in (G_2)_{e_1} = (G_2)^w = SU(3)$ (Theorem 1.9.4). Thus there exists $P \in SU(3)$ such that
$$
    \beta(X + M) = X + PM = \varphi(P, E, E)(X + M), \quad X + M \in (\gJ_{\sC})^C \oplus M(3, \C)^C = \gJ^C. $$
Therefore we have
$$
    \alpha = \varphi(E, A, B)\beta = \varphi(E, A, B)\varphi(P, E, E) = \varphi(P, A, B). $$
In the latter case, consider the mapping $\gamma_1 : \gJ^C \to \gJ^C, \gamma_1(X + M) = \ov{X} + \ov{M}, X + M \in (\gJ_{\sC})^C \oplus M(3, \C)^C = \gJ^C$ and remember that $\gamma_1 \in G_2 \subset F_4 \subset E_6$. Let $\beta = \alpha^{-1}\varphi(E, A, B)\gamma_1$, then $\beta \in E_6$ and $\beta|(\gJ_{\sC})^C = 1$. Hence $\beta \in (G_2)_{e_1} = (G_2)^w$ (Theorem 1.9.4) $\subset (E_6)^w$. Since $\alpha$ and $\varphi(E, A, B) \in (E_6)^w$, we have $\gamma_1 \in (E_6)^w$, so that $\gamma_1 \in (G_2)^w$ which is a contradiction (Theorem 1.9.4). Therefor we see that $\varphi$ is onto. That $\Ker\,\varphi = \{(E, E, E), $ $(\omega_1E, \omega_1E, \omega_1E),({\omega_1}^2E, {\omega_1}^2E, {\omega_1}^2E) \} = \Z_3$ can be easily obtained. Thus we have the isomorphism $(SU(3) \times SU(3) \times SU(3))/\Z_3 \cong (E_6)^w$.
\vspace{3mm}

{\bf Remark 1.} Since $(E_6)^w$ is connected as the fixed points subgroup of $E_6$ by an automorphism $w$ of order 3 of the simply connected group $E_6$, the fact that $\varphi : SU(3) \times SU(3) \times SU(3) \to (E_6)^w$ is onto can be proved as follows. The elements
$$
\begin{array}{l}
   G_{01}, \quad G_{23}, \quad \; \, G_{45}, \quad G_{67}, \quad G_{26} + G_{37}, \quad -G_{27} + G_{36}, \vspace{1mm}\\
   G_{24} + G_{35}, \quad -G_{25} + G_{34}, \quad G_{46} + G_{57}, \quad -G_{47} + G_{56}, 
\vspace{1mm}\\
   \wti{A}_1(1), \quad \wti{A}_2(1), \quad \wti{A}_3(1), \quad \wti{A}_1(e_1), \quad \wti{A}_2(e_1), \quad \wti{A}_3(e_1), 
\vspace{1mm}\\
     (E_1 - E_2)^{\sim}, \quad \wti{F}_1(1), \quad \; \wti{F}_2(1), \quad \; \wti{F}_3(1), 
\vspace{1mm}\\
     (E_2 - E_3)^{\sim}, \quad \wti{F}_1(e_1), \quad \wti{F}_2(e_1), \quad \wti{F}_3(e_1) 
\end{array}$$
forms an $\R$-basis of $(\gge_6)^w$. So $\dim((\gge_6)^w) = 16 + 8 = 24 = 8 + 8 + 8 = \dim(\su(3) \oplus \su(3) \oplus \su(3))$, Hence $\varphi$ is onto.
\vspace{2mm}

{\bf Remark 2.} The group $E_6$ has a subgroup which is isomorphic to the group $((SU(3) \times SU(3) \times SU(3))/\Z_3)\cdot\Z_2$ which is the semi-direct product of the groups $(SU(3) \times SU(3) \times SU(3))/\Z_3$ and $\Z_2$ (the action of $\Z_2 = \{1, \gamma_1 \}$ on the group $SU(3) \times SU(3) \times SU(3)$ is $\gamma_1(P, A, B) = (\ov{P}, \ov{B}, \ov{A}))$.
\vspace{4mm}

{\bf 3.14. Complex exceptional Lie group ${E_6}^C$}
\vspace{3mm}

{\bf Theorem 3.14.1.} {\it The polar decomposition of the Lie group ${E_6}^C$ is given by}
$$
      {E_6}^C \simeq E_6 \times \R^{78}.$$
{\it In particular, ${E_6}^C$ is a simply connected complex Lie group of type $E_6$}.
\vspace{2mm}

{\bf Proof.} Evidently ${E_6}^C$ is an algebraic subgroup of $\Iso_C(\gJ^C) = GL(27, C)$. The conjugate transposed mapping $\alpha^*$ of $\alpha \in {E_6}^C$ with respect to the inner product $\langle X, Y \rangle$ is $\alpha^* = \tau^t\alpha\tau \in {E_6}^C$. Hence, from Chevalley's lemma, we have
$$
     {E_6}^C \simeq ({E_6}^C \cap U(\gJ^C)) \times \R^d = E_6 \times \R^{d}, \quad d = 78.$$
Since $E_6$ is simply connected (Theorem 3.9.2), ${E_6}^C$ is also simply connected. The Lie algebra of the group ${E_6}^C$ is ${\gge_6}^C$, so ${E_6}^C$ is a complex simple Lie group of type $E_6$.
\vspace{4mm}

{\bf 3.15. Non-compact exceptional Lie groups $E_{6(6)}, E_{6(2)}, E_{6(-14)}$ and $ E_{6(-26)}$ of type $E_6$}
\vspace{3mm}

We define Hermitian inner products $\langle X, Y \rangle_{\gamma}$ and $\langle X, Y \rangle_{\sigma}$ in $\gJ(3, \gC^C)$ respectively by
$$
       \langle X, Y \rangle_{\gamma} = \langle \gamma X, Y \rangle,
 \quad \langle X, Y \rangle_{\sigma} = \langle \sigma X, Y \rangle. $$
Now, we define groups $E_{6(6)}, E_{6(2)}, E_{6(-14)}$ and $E_{6(-26)}$ respectively by
\begin{eqnarray*}
    E_{6(6)} \!\!\! &=& \!\!\! \{\alpha \in \Iso_{\sR}(\gJ(3, \gC')) \, | \, \det\,(\alpha X) = \det\,X \}, 
\vspace{1mm}\\
    E_{6(2)} \!\!\! &=& \!\!\! \{\alpha \in \Iso_C(\gJ(3, \gC^C)) \, | \, \det\,(\alpha X) = \det\,X,  \langle \alpha X, \alpha Y \rangle_{\gamma} = \langle X, Y \rangle_{\gamma} \}, 
\vspace{1mm}\\
    E_{6(-14)} \!\!\! &=& \!\!\! \{\alpha \in \Iso_C(\gJ(3, \gC^C)) \, | \, \det\,(\alpha X) = \det\, X,  \langle \alpha X, \alpha Y \rangle_{\sigma} = \langle X, Y \rangle_{\sigma} \},
\vspace{1mm}\\
     E_{6(-26)} \!\!\! &=& \!\!\! \{\alpha \in \Iso_{\sR}(\gJ(3, \gC)) \, | \, \det\,(\alpha X) = \det\,X \}. 
\end{eqnarray*}
These groups can also be defined by
$$
         E_{6(6)} \cong ({E_6}^C)^{\tau\gamma}, \, 
         E_{6(2)} \cong ({E_6}^C)^{\tau\lambda\gamma}, \, 
         E_{6(-14)} \cong ({E_6}^C)^{\tau\lambda\sigma}, \,
         E_{6(-26)} \cong ({E_6}^C)^{\tau}. $$

{\bf Theorem 3.15.1.} {\it The polar decompositions of the Lie groups $E_{6(6)}, E_{6(2)}, E_{6(-14)}$ and $E_{6(-26)}$ are respectively given by}
\begin{eqnarray*}
    E_{6(6)} \!\!\! &\simeq& \!\!\! Sp(4)/\Z_2 \times \R^{42}, 
\vspace{1mm}\\
    E_{6(2)} \!\!\! &\simeq& \!\!\! (Sp(1) \times SU(6))/\Z_2 \times \R^{40}, 
\vspace{1mm}\\
  E_{6(-14)} \!\!\! &\simeq& \!\!\! (U(1) \times Spin(10))/\Z_4 \times \R^{32}, \vspace{1mm}\\
    E_{6(-26)} \!\!\! &\simeq& \!\!\! F_4 \times \R^{26}. 
\end{eqnarray*}

{\bf Proof.} These are the facts corresponding to Theorems 3.12.2, 3.11.4, 3.10.7 and 3.7.1.
\vspace{4mm}

{\bf Theorem 3.15.2.} {\it The centers of the groups $E_{6(6)}, E_{6(2)}, E_{6(-14)}$ and $E_{6(-26)}$ are given by}
$$
      z(E_{6(6)}) = \{ 1 \}, \;\; z(E_{6(2)}) = \Z_3, \;\;  z(E_{6(-14)}) = \Z_3, \;\; z(E_{6(-26)}) = \{ 1 \}. $$

\newpage

\vspace{5mm}

\begin{center}
\large{\bf 4. Exceptional Lie group $E_7$}
\end{center}
\vspace{4mm} 

{\bf 4.1. Freudenthal vector space $\gP^C$}
\vspace{3mm}

{\bf Definition.} We define a $C$-vector space $\gP^C$, called the Freudenthal $C$-vector space, by
$$
      \gP^C = \gJ^C \oplus \gJ^C \oplus C \oplus C.$$
$\gP^C$ is a 56 dimensional $C$-vector space. An element $\pmatrix{ X \cr
                                                                     Y \cr
                                                                     \xi \cr
                                                                     \eta}$ 
of $\gP^C$ is often 
\vspace{0.7mm}
denoted by $(X, Y, \xi, \eta)$ or sometimes by $\dot{X} + \d{$Y$} + \dot{\xi} + \d{$\eta$}.$  In $\gP^C$, we define an inner product $(P, Q)$, a positive definite Hermitian inner product $\langle P, Q \rangle$ and a skew-symmetric inner product $\{P, Q\}$ respectively by
\begin{eqnarray*}
      (P, Q) \!\!\! &=& \!\!\! (X, Z) + (Y, W) + \xi\zeta + \eta\omega,
\vspace{1mm}\\
      \langle P, Q \rangle \!\!\! &=& \!\!\! \langle X, Z \rangle + \langle Y, W \rangle + (\tau\xi)\zeta + (\tau\eta)\omega,
\vspace{1mm}\\
      \{ P, Q \} \!\!\! &=& \!\!\! (X, W) - (Z, Y) + \xi\omega - \zeta\eta,
\end{eqnarray*}
where $P = (X, Y, \xi, \eta)$, $Q = (Z, W, \zeta, \omega) \in \gP^C$.
\vspace{2mm}

For $\phi \in {\gge_6}^C$, $A, B \in \gJ^C$, $\nu \in C$, we define a $C$-linear mapping ${\mit\Phi}(\phi, A, B, \nu) : \gP^C \to \gP^C$ by
\begin{eqnarray*}
     {\mit\Phi}(\phi, A, B, \nu)\pmatrix{X \vspace{1mm}\cr 
                                         Y \vspace{1mm}\cr 
                                         \xi \vspace{1mm}\cr \eta}
  \!\!\! &=& \!\!\! \pmatrix{\phi-\dfrac{1}{3}\nu & 2B & 0 & A \cr
                             2A & -{}^t\phi + \dfrac{1}{3}\nu & B & 0 
\vspace{1mm}\cr
                             0 & A & \nu & 0
\vspace{1mm}\cr
                             B & 0 & 0 & -\nu}
                    \pmatrix{X \vspace{1mm}\cr 
                             Y \vspace{1mm}\cr 
                             \xi \vspace{1mm}\cr \eta}
\vspace{1mm}\\
\!\!\! &=& \!\!\! \pmatrix{\phi X - \dfrac{1}{3}\nu X + 2B \times Y + \eta A \cr                           2A \times X - {}^t\!\phi \,Y + \dfrac{1}{3}\nu Y + \xi B 
\vspace{1mm}\cr
                           (A, Y) + \nu\xi 
\vspace{1mm}\cr
                           (B, X) - \nu\eta}.
\end{eqnarray*}

{\bf Definition.} For $P = (X, Y, \xi, \eta)$, $Q = (Z, W, \zeta, \omega) \in \gP^C$, we define a $C$-linear mapping $P \times Q : \gP^C \to \gP^C$ by
$$
        P \times Q = {\mit\Phi}(\phi, A, B, \nu), \quad
\left\{\begin{array}{rcl}
        \phi \!\!\! &=& \!\!\! - \dfrac{1}{2}(X \vee W + Z \vee Y)
\vspace{1mm}\\
        A \!\!\! &=& \!\!\! - \dfrac{1}{4}(2Y \times W - \xi Z - \zeta X)
\vspace{1mm}\\
  B \!\!\! &=& \!\!\! \;\;\; \dfrac{1}{4}(2X \times Z - \eta W - \omega Y)
\vspace{1mm}\\
  \nu \!\!\! &=& \!\!\! \;\;\; \dfrac{1}{8}((X, W) + (Z, Y) - 3(\xi\omega + \zeta\eta)).
\end{array}\right.$$

{\bf Lemma 4.1.1.} {\it For $P, Q, R \in \gP^C$, we have}
\vspace{1mm}

(1)\quad $P \times Q = Q \times P$.
\vspace{1mm}

(2)\quad  $(P \times Q)P - (P \times P)Q + \dfrac{3}{8}\{P, Q \}P = 0$.
\vspace{1mm}

(3)\quad  $(P \times R)Q - (Q \times R)P + \dfrac{1}{8}\{Q, R \}P - \dfrac{1}{8}\{P, R \}Q - \dfrac{1}{4}\{P, Q \}R = 0$.
\vspace{2mm}

{\bf Proof.} (1) is evident.
\vspace{1mm}

(2) For $P = (X, Y, \xi, \eta)$, $Q = (Z, W, \zeta, \omega) \in \gP^C$, we have
$$
\begin{array}{l}
 (P \times Q)P
\vspace{1mm}\\
    = {\mit\Phi}\Big(- \dfrac{1}{2}(X \vee W + Z \vee Y), - \dfrac{1}{4}(2Y \times W - \xi Z - \zeta X), \dfrac{1}{4}(2X \times Z - \eta W - \omega Y),
\vspace{1mm}\\
    \quad \dfrac{1}{8}((X, W) + (Z, Y) - 3(\xi\omega + \zeta\eta))\Big)(X, Y, \xi, \eta)
\vspace{1mm}\\
   = \cdots (\mbox{using the formula}\, \,  X \vee Y \, \, \mbox{of Lemma 3.4.1), etc.)} \cdots 
\vspace{1mm}\\
   = {\mit\Phi}\Big(- X \vee Y, - \dfrac{1}{2}(Y \times Y - \xi X), \dfrac{1}{2}(X \times X - \eta Y), \dfrac{1}{4}((X, Y) - 3\xi\eta)\Big)(Z, W, \zeta, \omega),
\vspace{1mm}\\
   \quad - \dfrac{3}{8}((X, W) - (Z, Y) + \xi\omega - \zeta\eta)(X, Y, \xi, \eta)
\vspace{1mm}\\
   = (P \times P)Q - \dfrac{3}{8}\{P, Q \}P.
\end{array}$$

(3) In the equality (2), put $P + R$ in the place of $P$, then we have
$$
\displaylines{\hfill 
        2(P \times R)Q - (P \times Q)R -(R \times Q)P + 
     \dfrac{3}{8}\{Q, R\}P - \dfrac{3}{8}\{P, Q\}R = 0. 
\hfill \mbox{(i)}} $$
Exchanging $P$ with $Q$, we see that
$$
\displaylines{\hfill 
    2(Q \times R)P - (Q \times P)R - (R \times P)Q +
    \dfrac{3}{8}\{P, R\}Q - \dfrac{3}{8}\{Q, P\}R = 0.
 \hfill \mbox{(ii)}} $$
Taking ((i)$-$(ii))$\div 3$, we have
$$
      (P \times R)Q - (Q \times R)P + \dfrac{1}{8}\{Q, R\}P - \dfrac{1}{8}\{P, R\}Q - \vspace{2mm}
\dfrac{1}{4}\{P, Q\}R = 0.$$

We define a space $\gM^C$, called the complex Freudenthal manifold, by
\begin{eqnarray*}
   \gM^C \!\!\! &=& \!\!\! \{ P \in \gP^C \, | \, P \times P = 0, P \neq 0 \}
\vspace{1mm}\\
    \!\!\! &=& \!\!\! \left\{\begin{array}{l}
         P = (X, Y, \xi, \eta) \in \gP^C 
\vspace{1mm}\\
         P \neq 0
\end{array}\right.
\left|\begin{array}{l}
         X \vee Y = 0, (X, Y) = 3\xi\eta, 
\vspace{1mm}\\
         X \times X = \eta Y, Y \times Y = \xi X
\end{array}\right\} .
\end{eqnarray*}

{\bf Lemma 4.1.2.} {\it The following elements} ({\it assuming $\xi \neq 0, \eta \neq 0$}) {\it of} $\gP^C$
$$
         \pmatrix{X \cr \dfrac{1}{\eta}(X \times X) 
\vspace{1mm}\cr
                  \dfrac{1}{\eta^2}\det\,X 
\vspace{1mm}\cr
                  \eta}, \quad
         \pmatrix{\dfrac{1}{\xi}(Y \times Y) 
\vspace{1mm}\cr 
                  Y 
\vspace{1mm}\cr
                  \xi 
\vspace{1mm}\cr
                  \dfrac{1}{\xi^2}\det\,Y}, \quad
         \dot{1} = \pmatrix{0 
\vspace{1mm}\cr
                  0 
\vspace{1mm}\cr
                  1 \vspace{1mm}\cr 0}, \quad
         \d{1} = \pmatrix{0 
\vspace{1mm}\cr
                  0 
\vspace{1mm}\cr
                  0 \vspace{1mm}\cr 1} $$
{\it belong to} $\gM^C$.
\vspace{2mm}

{\bf Proof.} This is clear from $X \vee (X \times X) = 0$ (Lemma 3.5.4.(1)) and $(X \times X) \times (X \times X) = (\det\,X)X$ (Lemma 2.1.1.(3)).
\vspace{4mm}

{\bf 4.2. Compact exceptional Lie group $E_7$}
\vspace{3mm}

{\bf Definition.} We define the groups ${E_7}^C$ and $E_7$ respectively by 
\begin{eqnarray*}
     {E_7}^C \!\!\! &=& \!\!\! \{ \alpha \in \Iso_C(\gP^C) \, | \, \alpha(P \times Q)\alpha^{-1} = \alpha P \times \alpha Q \},
\vspace{1mm}\\
     E_7 \!\!\! &=& \!\!\! \{ \alpha \in \Iso_C(\gP^C) \, | \, \alpha(P \times Q)\alpha^{-1} = \alpha P \times \alpha Q, \langle \alpha P, \alpha Q \rangle = \langle P, Q \rangle \}.
\end{eqnarray*}

$E_7$ is a subgroup of ${E_7}^C$.
\vspace{3mm}

{\bf Theorem 4.2.1.}$\;$ {\it $E_7$ is a compact Lie group.}
\vspace{2mm}

{\bf Proof.} $E_7$ is a compact Lie group as a closed subgroup of the unitary group 
$$
      U(56) = U(\gP^C) = \{ \alpha \in \Iso_C(\gP^C) \, | \  \langle \alpha P, \alpha Q \rangle = \langle P, Q \rangle \}. $$

{\bf Proposition 4.2.2.} {\it For} $\alpha \in {E_7}^C$ ({\it and so, for} $\alpha \in E_7$), {\it we have}
\vspace{1mm}

(1)\quad $\alpha\gM^C = \gM^C$.
\vspace{1mm}

(2)\quad $\{ \alpha P, \alpha Q \} = \{ P, Q \}, \quad P, Q \in \gP^C$.
\vspace{2mm}

{\bf Proof.} (1) It is sufficient to prove that $\alpha\gM^C \subset \gM^C$. Now, for $P \in \gM^C$, we have $\alpha P \times \alpha P = \alpha(P \times P)\alpha^{-1} = \alpha 0 \alpha^{-1} = 0$. Hence $\alpha P \in \gM^C$.
\vspace{1mm}

$\begin{array}{l}
\mbox{(2)} \quad \{ \alpha P, \alpha Q \}\alpha P = \dfrac{8}{3}((\alpha P \times \alpha P)\alpha Q - (\alpha P \times \alpha Q)\alpha P) \;\; (\mbox{Lemma 4.1.1.(2)})
\vspace{1mm}\\
\qquad \qquad \qquad \qquad = \dfrac{8}{3}(\alpha(P \times P)Q - \alpha(P \times Q)P) = \{ P, Q \}\alpha P.
\end{array}$
\vspace{1mm}

\noindent It follows $\{ \alpha P, \alpha Q \} = \{ P, Q \}$.
\vspace{4mm}

{\bf 4.3. Lie algebra $\gge_7$ of $E_7$}
\vspace{3mm}

Before we investigate the Lie algebra $\gge_7$ of the group $E_7$, we shall first study the Lie algebra ${\gge_7}^C$ of the group ${E_7}^C$.
\vspace{3mm}

{\bf Theorem 4.3.1.} {\it The Lie algebra ${\gge_7}^C$ of the group ${E_7}^C$ is given by}
$$
     {\gge_7}^C = \{ {\mit\Phi}(\phi, A, B, \nu) \in \Hom_C(\gP^C) \, | \, \phi \in {\gge_6}^C,A, B \in \gJ^C, \nu \in C \}.$$
{\it The Lie bracket} $[{\mit\Phi_1}, {\mit\Phi_2}]$ in ${\gge_7}^C$ {\it is given by}
$$
      [{\mit\Phi}(\phi_1, A_1, B_1, \nu_1), {\mit\Phi}(\phi_2, A_2, B_2, \nu_2) ] = {\mit\Phi}(\phi, A, B, \nu), $$
{\it where}

$$
\left\{\begin{array}{rcl}
     \phi \!\!\! &=& \!\!\! [\phi_1, \phi_2] + 2A_1 \vee B_2 - 2A_2 \vee B_1
\vspace{1mm}\\
     A \!\!\! &=& \!\!\! \Big(\phi_1 + \dfrac{2}{3}\nu_1\Big)A_2 - \Big(\phi_2 + \dfrac{2}{3}\nu_2\Big)A_1
\vspace{1mm}\\
     B \!\!\! &=& \!\!\! -\Big({}^t\phi_1 + \dfrac{2}{3}\nu_1\Big)B_2 + \Big({}^t\phi_2 + \dfrac{2}{3}\nu_2\Big)B_1
\vspace{1mm}\\
     \nu \!\!\! &=& \!\!\! (A_1, B_2) - (B_1, A_2).
\end{array}\right.$$
{\it In particular, the dimension of} ${\gge_7}^C$ is
$$
        \dim_C({\gge_7}^C) = 78 + 27 \times 2 + 1 = 133. $$

{\bf Proof.} Before we show that ${\gge_7}^C$ is the Lie algebra of the group ${E_7}^C$, we first check the form of the Lie bracket in ${\gge_7}^C$. For ${\mit\Phi}_i \in {\gge_7}^C$ and $P \in \gP^C$, we have
$$
\begin{array}{l}
     [{\mit\Phi}_1, {\mit\Phi}_2]P = {\mit\Phi}_1{\mit\Phi}_2P - {\mit\Phi}_2{\mit\Phi}_1P 
\vspace{1mm}\\
    \qquad \qquad = \cdots ( \mbox{using the formula of $A \vee B$ (Lemma 3.4.1) etc.)} \cdots 
\vspace{1mm}\\
    \qquad \qquad = {\mit\Phi}P. 
\end{array}$$
This ${\mit\Phi}$ is that of the theorem.
\vspace{1mm}

We now determine the Lie algebra ${\gge_7}^C$ of the group ${E_7}^C$. Since $\alpha \in {E_7}^C$ satisfies
$$
\left\{\begin{array}{ll}
     \alpha P \times \alpha P = 0, & P \in \gM^C 
\vspace{1mm}\\
     \{\alpha P, \alpha Q \} = \{P, Q \}, & P, Q  \in \gP^C \;\;\;
\mbox{(Proposition 4.2.2.(2))},
\end{array}\right.$$
if ${\mit\Phi} \in \Hom_C(\gP^C)$ belongs to ${\gge_7}^C$, then ${\mit\Phi}$ satisfies
$$
\displaylines{\hfill
\left\{\begin{array}{ll}
       {\mit\Phi}P \times P = 0, & P \in \gM^C 
\vspace{1mm}\\
        \{ {\mit\Phi}P, Q \} + \{ P, {\mit\Phi}Q \} = 0, & P, Q  \in \gP^C.
\end{array}\right.
\hfill
\left.\begin{array}{r}
\mbox{(i)}
\vspace{1mm}\\
\mbox{(ii)}
\end{array}\right.} $$
Since ${\mit\Phi} \in {\gge_7}^C$ is a $C$-linear mapping of $\gP^C = \gJ^C \oplus \gJ^C \oplus C \oplus C$, ${\mit\Phi}$ is of the form
$$
      {\mit\Phi} = \pmatrix{g & l & C & A \cr
                            k & h & B & D \cr
                            c & a & \nu & \lambda \cr
                            b & d & \kappa & \mu},\qquad
\begin{array}{l}
      g, h, k, l \in \Hom_C(\gJ^C),\\
      a, b, c, d \in \Hom_C(\gJ^C, C),\\
      A, B, C, D \in \gJ^C, \\
      \nu, \mu, \kappa, \lambda \in C.
\end{array} $$
For $0 \neq r \in C$, we define a $C$-linear isomorphism $f_r : \gP^C \to \gP^C$ by
$$
      f_r (X, Y, \xi, \eta) = (X, rY, r^2\xi, r^{-1}\eta),$$
then it is easy to see that $f_r$ satisfies
$$
      rf_r(P \times Q){f_r}^{-1} = f_rP \times f_rQ, \quad P, Q \in \gP^C.$$
Hence we see that for $\alpha \in {E_7}^C$ we have $f_r\alpha{f_r}^{-1} \in {E_7}^C$. Therefore, for ${\mit\Phi} \in {\gge_7}^C$,  we have

$$
      {\gge_7}^C \ni f_r{\mit\Phi}{f_r}^{-1} = 
        \pmatrix{g & r^{-1}l & r^{-2}C & rA 
\vspace{1mm}\cr
                 rk & h & r^{-1}B & r^2D 
\vspace{1mm}\cr
                 r^2c & ra & \nu & r^3\lambda 
\vspace{1mm}\cr
                 r^{-1}b & r^{-2}d & r^{-3}k & \mu } $$
for any $0 \neq r \in C$. Hence ${\mit\Phi}$ is decomposed as
$$
     {\mit\Phi} = {\mit\Phi}_{-3} + {\mit\Phi}_{-2} + {\mit\Phi}_{-1} + {\mit\Phi}_0 + {\mit\Phi}_1 + {\mit\Phi}_2 + {\mit\Phi}_3, \quad {\mit\Phi}_i \in {\gge_7}^C,$$
where
$$
\begin{array}{l}
     {\mit\Phi}_{-3} = \pmatrix{0 & 0 & 0 & 0 \cr
                                0 & 0 & 0 & 0 \cr
                                0 & 0 & 0 & 0 \cr
                                0 & 0 & \kappa & 0}, \quad
     {\mit\Phi}_{3} = \pmatrix{0 & 0 & 0 & 0 \cr
                              0 & 0 & 0 & 0 \cr
                              0 & 0 & 0 & \lambda \cr
                              0 & 0 & 0 & 0}, 
\vspace{1mm}\\
     {\mit\Phi}_{-2} = \pmatrix{0 & 0 & C & 0 \cr
                                 0 & 0 & 0 & 0 \cr
                                 0 & 0 & 0 & 0 \cr
                                 0 & d & 0 & 0}, \quad
     {\mit\Phi}_{2} = \pmatrix{0 & 0 & 0 & 0 \cr 
                               0 & 0 & 0 & D \cr 
                               c & 0 & 0 & 0 \cr 
                               0 & 0 & 0 & 0}, 
\vspace{1mm}\\
      {\mit\Phi}_{-1} = \pmatrix{0 & l & 0 & 0 \cr
                                 0 & 0 & B & 0 \cr
                                 0 & 0 & 0 & 0 \cr 
                                 b & 0 & 0 & 0}, \quad
      {\mit\Phi}_{1} = \pmatrix{0 & 0 & 0 & A \cr 
                                k & 0 & 0 & 0 \cr
                                0 & a & 0 & 0 \cr
                                0 & 0 & 0 & 0},
\end{array} $$
\vspace{-1mm}
$$
     {\mit\Phi}_{0} = \pmatrix{g & 0 & 0 & 0 \cr
                               0 & h & 0 & 0 \cr
                               0 & 0 & \nu & 0 \cr 
                               0 & 0 & 0 & \mu}.$$
The relation ${\mit\Phi}_{-3}(0, 0, 1, 0) \times (0, 0, 1, 0) = 0$ implies $\kappa = 0$, hence ${\mit\Phi}_{-3} = 0$. Similarly ${\mit\Phi}_3 = 0$. The relation ${\mit\Phi}_{-2}(0, 0, 1, 0) \times (0, 0, 1, 0) = 0$ implies $C = 0$. Moreover, the relation ${\mit\Phi}_{-2}P \times P = 0$ for $P = (Y \times Y, Y, 1, \det\,Y) \in \gM^C$, that is, $(0, 0, 0, d(Y)) \times (Y \times Y, Y, 1, \det\,Y) = 0$ implies $d = 0$. Hence ${\mit\Phi}_{-2} = 0$. Similarly ${\mit\Phi}_2 = 0$. 
\vspace{1mm}
The relation ${\mit\Phi}_{-1}P \times P = 0$ for $P = (Y \times Y, Y, 1, \det\,Y) \in \gM^C$, that is, $(l(Y), B, 0, b(Y \times Y)) \times (Y \times Y, Y,1, \det\,Y) = 0$ implies
$$
         l(Y) = 2B \times Y, \quad Y \in \gJ^C. $$
Next, the relation ${\mit\Phi}_{-1}P \times P = 0$ for $P = (X, X \times X, \det\,X, 1) \in \gM^C$, that is, $(l(X \times X), (\det\,X)B, 0, b(X)) \times (X, X \times X, \det\,X, 1) = 0$, (from the 4-th condition) we have
$$
      2(B \times (X \times X),X \times X) + (\det\,X)(X, B) = 3(\det\,X)b(X), $$
hence, $3(\det\,X)(B, X) = 3(\det\,X)b(X)$, and so we have
$$
              b(X) =  (B, X),  \quad X \in \gJ^C. $$
(Since $b$ is continuous, the above is also valid for $X \in \gJ^C$ such that $\det\,X = 0$). Similarly, using ${\mit\Phi}_1$, we have 

$$
        k(X) = 2A \times X, \quad a(Y)=(A, Y),  \quad X, Y \in \gJ^C. $$ 
The relation ${\mit\Phi}_{0}P \times P=0$ for $P = (X, X \times X, \det\,X, 1) \in \gM^C$, that is, $(g(X), h(X \times X), (\det\,X)\nu, \mu) \times (X, X \times X, \det\,X, 1) = 0$ implies 
$$
\displaylines{\hfill
\begin{array}{c}
      2g(X) \times X = \mu X \times X + h(X \times X), 
\vspace{1mm}\\
    2h(X \times X) \times (X \times X) = (\det\,X)(\nu X + g(X)), 
\vspace{1mm}\\
    (g(X), X \times X) + (h(X \times X), X) = 3(\nu + \mu)\det\,X.
\end{array}
\hfill
\left.\begin{array}{r}
\mbox{(i)} 
\vspace{1mm}\\
\mbox{(ii)} 
\vspace{1mm}\\
\mbox{(iii)}
\end{array}\right.}$$
Putting $\phi = g - \dfrac{1}{3}(\nu + 2\mu)1$, then, using (i) and (iii), we have
$$
\begin{array}{l}
      3(\phi X, X, X) = 3(g(X), X, X) - (\nu + 2\mu)(X, X, X) 
\vspace{1mm}\\
  \quad = (\mu X \times X + h(X \times X), X) + (g(X), X \times X) - 3(\nu + 2\mu)\det\,X 
\vspace{1mm}\\
   \quad = 3\mu\det\,X + 3(\nu + \mu)\det\,X - 3(\nu + 2\mu)\det\,X = 0.
\end{array} $$
Therefore 
$$
                     \phi \in {\gge_6}^C.  $$
Similarly $\psi = h - \dfrac{1}{3}(2\nu + \mu)1 \in {\gge_6}^C$. Furthermore, from (ii), we have
$$
\begin{array}{c}
      2\Big(\psi(X \times X) + \dfrac{1}{3}(2\nu + \mu)(X \times X)\Big) \times (X \times X)\\
\hspace{25mm} = (\det X)\Big(\nu X + \phi X + \dfrac{1}{3}(\nu + 2\mu)X\Big),
\end{array}$$
so $2\psi(X \times X) \times (X \times X) = (\det\,X)\phi X,$ and so (for a while, instead of $-{}^t\psi$, we use $\psi'$ again) $\psi'((X \times X) \times (X \times X)) = (\det\,X)\phi X$ (Lemma 3.4.3.(1)), hence $(\det\,X)\psi'X = (\det\,X)\phi X.$  Therefore, we have $\psi'X = \phi X, X \in \gJ^C,$ (even if for $X \in \gJ^C$ such that $\det\,X = 0$), that is,
$$
                \psi' = \phi. $$
Finally, the relation $\{{\mit\Phi}(0, 0, 1, 0), (0, 0, 0, 1) \} + \{(0, 0, 1, 0), {\mit\Phi}(0, 0, 0, 1) \} = 0$ implies $\nu + \mu = 0$. Thus we see that ${\mit\Phi} \in {\gge_7}^C$ is of the form
$$
\begin{array}{l}
      {\mit\Phi} = \pmatrix{\phi - \dfrac{1}{3}\nu & 2B & 0 & A \cr
                            2A & \phi' + \dfrac{1}{3}\nu & B & 0 
\vspace{1mm}\cr
                            0 & A &  \nu & 0 
\vspace{1mm}\cr
                            B & 0 & 0 & -\nu} 
\vspace{1mm}\\
\;\;\;
        = {\mit\Phi}(\phi, A, B, \nu), \quad \phi \in {\gge_6}^C, A, B \in \gJ^C,\nu \in C.
\end{array} $$
Conversely, we shall show that ${\mit\Phi} = {\mit\Phi}(\phi, A, B, \nu)$, $\phi \in {\gge_6}^C$, $A, B \in \gJ^C$, $\nu \in C$ belongs to ${\gge_7}^C$, that is, $\exp t\mit\Phi \in {E_7}^C$ for all $t \in C$. For this purpose, we prove the following proposition.
\vspace{3mm}

{\bf Proposition 4.3.2.} {\it For} ${\mit\Phi} = {\mit\Phi}(\phi, A, B, \nu)$, $\phi \in {\gge_6}^C$, $A, B \in \gJ^C$, $\nu \in C$ {\it and} $P, Q \in \gP^C$, {\it we have}
$$
   [{\mit\Phi}, P \times Q] = {\mit\Phi}P \times Q + P \times {\mit\Phi} Q.$$

{\bf Proof.} It is sufficient to show that $[{\mit\Phi}, P \times P] = 2{\mit\Phi}P \times P.$  For $P = (X, Y, \xi, \eta) \in \gP^C,$ we have
$$
\begin{array}{l}
  [{\mit\Phi}, P \times P]
\vspace{1mm}\\
  = \Big[{\mit\Phi}(\phi, A, B, \nu), {\mit\Phi}\Big(- X \vee Y, - \dfrac{1}{2}(Y \times Y - \xi X),
 \dfrac{1}{2}(X \times X - \eta Y), \dfrac{1}{4}((X, Y) - 3\xi\eta)\Big)\Big]
\vspace{1mm}\\
 = \cdots \mbox{(using} \;\; \phi(X \times Y) = \phi'X \times Y + X \times \phi'Y \;\; \mbox{(Lemma 3.4.3.(1)),} \;\;[\phi, A \vee B]
\vspace{1mm}\\
\quad = \phi A \vee B + A \vee \phi'B \;\;\mbox{(Lemma 3.4.4.(1)), the formula about}\; A \vee B \;\; \mbox{(Lemma}
\vspace{1mm}\\
\quad \mbox{3.4.1) etc.)} \cdots 
\vspace{1mm}\\
   = 2{\mit\Phi}P \times P.
\end{array} $$

We will now return to the proof of Theorem 4.3.1. For ${\mit\Phi} = {\mit\Phi}(\phi, A, B, \nu)$, $\phi \in {\gge_6}^C$, $A, B \in \gJ^C$, $\nu \in C$ and $t \in C$, we have
$$
\begin{array}{l}
  (\exp t{\mit\Phi})(P \times Q)(\exp t{\mit\Phi})^{-1} 
\vspace{1mm}\\
   \qquad = (\exp t(\mbox{ad}{\mit\Phi}))(P \times Q)\;\;\;\; ((\mbox{ad}{\mit\Phi}){\mit\Phi}_1 = [{\mit\Phi}, {\mit\Phi}_1], \; {\mit\Phi}_1 \in {\gge_7}^C)
\vspace{1mm}\\
    \qquad = \dsum_{n=0}^{\infty}\dfrac{1}{n!}(t(\mbox{ad}{\mit\Phi}))^n(P \times Q)
\vspace{1mm}\\
  \qquad = \dsum_{n=0}^{\infty}\dfrac{t^n}{n!}\Big(\dsum_{k+l=n}\dfrac{n!}{k!\,l!}{\mit\Phi}^kP \times {\mit\Phi}^lQ \Big)\;\;(\mbox{Proposition 4.3.2})
\vspace{1mm}\\
   \qquad = \dsum_{n=0}^{\infty}\Big(\dsum_{k+l=n}\dfrac{t^kt^l}{k!\,l!}{\mit\Phi}^kP \times {\mit\Phi}^lQ \Big)
\vspace{1mm}\\
    \qquad = \Big(\dsum_{k=0}^{\infty}\dfrac{1}{k!}(t{\mit\Phi})^kP \Big) \times \Big(\dsum_{l=0}^{\infty}\dfrac{1}{l!}(t{\mit\Phi})^lQ \Big)
\vspace{1mm}\\
   \qquad = (\exp t{\mit{\mit\Phi}})P \times (\exp t{\mit\Phi})Q.
\end{array} $$
Hence $\exp t{\mit\Phi} \in {E_7}^C$, so that ${\mit\Phi} \in {\gge_7}^C$. Thus the proof of Theorem 4.3.1 is 
\vspace{3mm}
completed.

{\bf Definition.} We define a $C$-linear transformation $\lambda$ of $\gP^C$ by
$$
          \lambda(X, Y, \xi, \eta) = (Y, - X, \eta, - \xi). $$

For $\alpha \in \Hom_C(\gP^C)$, we also denote by $^t\alpha$ the transpose of $\alpha$ with respect to the inner product $(P, Q)$: $(^t\alpha P, Q) = (P, \alpha Q)$.
\vspace{3mm}

{\bf Lemma 4.3.3.} (1) $\lambda \in E_7$ {\it and satisfies} $\lambda^2 = -1$.
\vspace{1mm}

(2) {\it For} $P, Q \in \gP^C$, {\it we have}
$$
     (P, Q) = \{\lambda P, Q \} = - \{P, \lambda Q \}, \quad 
    \langle P, Q \rangle = \{\tau\lambda P, Q \}. $$

(3) {\it For $\alpha \in {E_7}^C$, we have }
$$
         ^t\alpha^{-1} = \lambda\alpha\lambda^{-1}. $$

(4) {\it For} $\alpha \in {E_7}^C$, {\it we have} 
$$
       \alpha \in E_7 \quad \mbox{{\it if and only if}} \quad 
            \tau\lambda\alpha = \alpha\tau\lambda.  $$

(5) {\it For} ${\mit\Phi}(\phi, A, B, \nu) \in {\gge_7}^C$, {\it we have}
$$
        \lambda{\mit\Phi}(\phi, A, B, \nu)\lambda^{-1} = 
         {\mit\Phi}(- {}^t\phi, - B, - A, - \nu). $$

{\bf Proof.} (1) and (2) are evident.
\vspace{1mm}

(3) $(P, \lambda Q) = \{P, Q \} = \{ \alpha P, \alpha Q \} \, \mbox{(Proposition 4.2.2.(2))}\, = (\alpha P, \lambda\alpha Q) \; \mbox{(Lemma} \\ \mbox{4.3.3)}\; = (P, $ ${}^t\alpha\lambda\alpha Q)$ implies $\lambda = {}^t\alpha\lambda\alpha$, hence ${}^t\alpha^{-1} = \lambda\alpha\lambda^{-1}$.
\vspace{1mm}

(4) If $\alpha \in {E_7}^C$ satisfies $\tau\lambda\alpha = \alpha\tau\lambda$, then $\langle \alpha P, \alpha Q \rangle = \{\tau\lambda\alpha P, \alpha Q \} = \{\alpha\tau\lambda P, \alpha Q \} $ $= \{\tau\lambda P, Q \}\, \mbox{(Proposition 4.2.2.(2))}\,  = \langle P, Q \rangle \; \mbox{(Lemma 4.3.3)}$, $P, Q \in \gP^C$.  Hence, $\alpha \in E_7$. The converse implication can also be easily proved.
\vspace{1mm}

(5) is easily checked by direct calculations.
\vspace{3mm}

{\bf Theorem 4.3.4.} {\it The Lie algebra ${\gge_7}$ of the group $E_7$ is given by}
$$
      \gge_7 = \{ {\mit\Phi}(\phi, A, -\tau A, \nu) \, | \, \phi \in \gge_6, A \in \gJ^C, \nu \in i\R \},$$
{\it where}
$$
        {\mit\Phi}(\phi, A, -\tau A, \nu) = 
         \pmatrix{\phi - \dfrac{1}{3}\nu & -2\tau A & 0 & A \cr
                  2A & \tau\phi\tau + \dfrac{1}{3}\nu & -\tau A & 0 
\vspace{1mm}\cr
                  0 & A & \nu & 0 
\vspace{1mm}\cr
                  - \tau A & 0 & 0 & -\nu}.$$
{\it The Lie bracket} $[{\mit\Phi}_1, {\mit\Phi}_2]$ {\it in $\gge_7$ is given by}
$$
       [{\mit\Phi}(\phi_1, A_1, -\tau A_1, \nu_1), {\mit\Phi}(\phi_2, A_2, -\tau A_2, \nu_2)] = {\mit\Phi}(\phi, A, -\tau A, \nu),$$
{\it where}
$$
\left\{\begin{array}{rcl}
       \phi \!\!\! &=& \!\!\! [\phi_1, \phi_2] - 2A_1 \vee \tau A_2 + 2A_2 \vee \tau A_1
\vspace{1mm}\\
       A \!\!\! &=& \!\!\! \Big(\phi_1 + \dfrac{2}{3}\nu_1\Big)A_2 - \Big(\phi_2 + \dfrac{2}{3}\nu_2\Big)A_1 
\vspace{1mm}\\
     \nu \!\!\! &=& \!\!\! \langle A_1, A_2 \rangle - \langle A_2, A_1 \rangle.
\end{array}\right.$$
{\it In particular, the dimension of $\gge_7$ is} 
$$
               \dim(\gge_7) =78 + 54 + 1 = 133. $$

{\bf Proof.} For ${\mit\Phi} \in {\gge_7}^C$,
$$
   {\mit\Phi} \in \gge_7 \quad \mbox{if and only if} \quad \tau\lambda{\mit\Phi}\lambda^{-1}\tau = {\mit\Phi}. $$
Therefore the theorem follows from Theorem 4.3.1, Lemma 4.3.3.(5) and $\tau{\mit\Phi}(\phi, A, $ $B, \nu)\tau = {\mit\Phi}(\tau\phi\tau, \tau A, \tau B, \tau\nu).$
\vspace{3mm}

{\bf Proposition 4.3.5.} {\it The complexification of the Lie algebra $\gge_7$ is ${\gge_7}^C$.}
\vspace{2mm}

{\bf Proof.} For ${\mit\Phi} \in {\gge_7}^C$, the conjugate transposed mapping ${\mit\Phi}^*$ of ${\mit\Phi}$ with respect to the inner product $\langle P, Q \rangle$ of $\gP^C$ is ${\mit\Phi}^* = \tau\lambda{\mit\Phi}\lambda\tau \in {\gge_7}^C$, and for ${\mit\Phi} \in {\gge_7}^C$, ${\mit\Phi}$ belongs to $\gge_7$ if and only if ${\mit\Phi}^* = - {\mit\Phi}$. Now, any element ${\mit\Phi} \in {\gge_7}^C$ is represented as
$$
    {\mit\Phi} = \dfrac{{\mit\Phi} - {\mit\Phi}^*}{2} + i\dfrac{{\mit\Phi} + {\mit\Phi}^*}{2i}, \quad  \dfrac{{\mit\Phi} - {\mit\Phi}^*}{2}, \dfrac{{\mit\Phi} + {\mit\Phi}^*}{2i} \in \gge_7.$$
Hence ${\gge_7}^C$ is the complexification of $\gge_7$.
\vspace{4mm}

{\bf 4.4. Simplicity of ${\gge_7}^C$}
\vspace{3mm}

{\bf Theorem 4.4.1.} {\it The Lie algebra ${\gge_7}^C$ is simple and so $\gge_7$ is also simple.}
\vspace{2mm}

{\bf Proof.} We use the decomposition of ${\gge_7}^C$ of Theorem 4.3.1
$$
      {\gge_7}^C = {\gge_6}^C \oplus \gN^C, $$
where ${\gge_6}^C = \{{\mit\Phi}(\phi, 0, 0, 0) \in {\gge_7}^C \, | \, \phi \in {\gge_6}^C \}$ and $\gN^C = \{{\mit\Phi}(0, A, B, \nu) \in {\gge_7}^C \, | \, A, B \in \gJ^C, \nu \in C \}$.  Let $p : {\gge_7}^C \to {\gge_6}^C$ and $q : {\gge_7}^C \to \gN^C$ be projections of ${\gge_7}^C = {\gge_6}^C \oplus \gN^C$. Now, let $\ga$ be a non-zero ideal of ${\gge_7}^C$. Then $p(\ga)$ is an ideal of ${\gge_6}^C$. Indeed, if $\phi \in p(\ga)$, then there exists ${\mit\Phi}(0, A, B, \nu) \in \gN^C$ such that ${\mit\Phi}(\phi, A, B, \nu) \in \ga$. For any $\phi_1 \in {\gge_6}^C$, we have 
$$
   \ga \ni [{\mit\Phi}(\phi_1, 0, 0, 0), {\mit\Phi}(\phi, A, B, \nu)] = {\mit\Phi}([\phi_1, \phi], \phi_1A, {\phi_1}'B, 0) $$ 
(Theorem 4.3.1), hence $[\phi_1, \phi] \in p(\ga)$.
\vspace{1mm}

We shall show that either ${\gge_6}^C \cap \,\ga \neq \{0\}$ or $\gN^C \cap \ga \neq \{0\}$. Assume that ${\gge_6}^C \,\cap\, \ga = \{0\}$ and $\gN^C \cap \ga = \{0\}$. Then the mapping $p|\ga : \ga \to {\gge_6}^C$ is injective because $\gN^C \cap \ga = \{0\}$. Since $p(\ga)$ is a non-zero ideal of ${\gge_6}^C$ and ${\gge_6}^C$ is simple, we have $p(\ga) = 
{\gge_6}^C$. Hence $\dim_C(\ga) = \dim_C(p(\ga)) = \dim_C({\gge_6}^C) = 78$. On the other hand, since ${\gge_6}^C \cap \ga = \{0\}$, $q|\ga : \ga \to \gN^C$ is also injective. Hence we have $\dim_C(\ga) \le \dim_C(\gN^C) = 27 + 27 + 1 = 55$. This leads to a contradiction. 
\vspace{1mm}

We now consider the following two cases.
\vspace{1mm}

(1) Case ${\gge_6}^C \cap \ga \neq \{0\}$. From the simplicity of ${\gge_6}^C$, we have ${\gge_6}^C \cap \ga = {\gge_6}^C$, hence $\ga \supset {\gge_6}^C$. On the other hand, we have
$$
\begin{array}{l}
   \ga \supset [\ga, {\gge_7}^C] \supset 
     [{\mit\Phi}({\gge_6}^C, 0, 0, 0), {\mit\Phi}(0, \gJ^C, 0, 0)]
\vspace{1mm}\\
  \qquad = {\mit\Phi}(0, {\gge_6}^C\gJ^C, 0, 0) = {\mit\Phi}(0, \gJ^C, 0, 0) \; \mbox{(Proposition 3.3.2.(3))}. 
\end{array}$$
Similarly we have ${\mit\Phi}(0, 0, \gJ^C, 0) \subset \ga$. Moreover, from
$$
     \ga \ni [{\mit\Phi}(0, E_1, 0, 0), {\mit\Phi}(0, 0, E_1, 0)] 
       = {\mit\Phi}(2E_1 \vee E_1, 0, 0, 1),$$
we have ${\mit\Phi}(0, 0, 0, 1) \in \ga$, and so $\ga \supset \gN^C$. Hence $\ga \supset {\gge_6}^C \oplus \gN^C = {\gge_7}^C$.
\vspace{1mm}

(2) Case $\gN^C \cap \ga \neq \{ 0 \}$.  Let ${\mit\Phi}(0, A, B, \nu)$ be a non-zero element of $\gN^C \cap \ga$.
\vspace{1mm}

$\;$ (i) Case ${\mit\Phi}(0, A, B, \nu), A \neq 0$. Choose $B_1 \in \gJ^C$ such that $A \vee B_1 \neq 0$ (Lemma 3.5.4.(2)), and choose $\phi \in {\gge_6}^C$ such that $[A \vee B_1, \phi] \neq 0$ (since ${\gge_6}^C$ is simple, such a $\phi$ exists because the center of ${\gge_6}^C$ is zero). Now, we have
$$
\begin{array}{l}
  \ga \ni \Big[{\mit\Phi}(0, A, B, \nu ), {\mit\Phi}\Big(0, 0, 0, -\dfrac{3}{2}\Big)\Big] = {\mit\Phi}(0, A, -B, 0),
\vspace{1mm}\\
  \ga \ni {[}{\mit\Phi}(0, A,-B, 0), {\mit\Phi}(0, 0, B_1, 0){]} = {\mit\Phi}(2A \vee B_1, 0, 0, (A, B_1)),
\vspace{1mm}\\
  \ga \ni [{\mit\Phi}(2A \vee B_1, 0, 0, (A, B_1)), {\mit\Phi}(\phi, 0, 0, 0)] = {\mit\Phi}(2[A \vee B_1. \phi], 0, 0, 0), $$
\end{array}$$
Hence this case is reduced to the case (1).
\vspace{1mm}

$\;$ (ii) Case ${\mit\Phi}(0, A, B, \nu), B \neq 0$. This case is also reduced to the case (1) in a similar way to (i).
\vspace{1mm}

$\;$ (iii) Case ${\mit\Phi}(0, 0, 0, \nu), \nu \neq 0$. If we choose $ 0 \neq A \in \gJ^C$, then we have
$$
     \ga \ni [{\mit\Phi}(0, 0, 0, \nu), {\mit\Phi}(0, A, 0, 0)] = {\mit\Phi}\Big(0, \frac{2}{3}\nu A, 0, 0\Big). $$
Hence this case is also reduced to the case (1).
\vspace{1mm}

\noindent Thus we have $\ga = {\gge_7}^C$.
\vspace{3mm}

{\bf Proposition 4.4.2.} (1) {\it $\gP^C$ is an ${\gge_7}^C$-irreducible $C$-module.}
\vspace{1mm}

(2) ${\gge_7}^C\gP^C = \Big\{\dsum_{k}{\mit\Phi}_kP_k \, | \, {\mit\Phi}_k \in {\gge_7}^C, P_k \in \gP^C \Big\} = \gP^C.$
\vspace{2mm}

{\bf Proof.} (1) Let $W$ be a non-zero ${\gge_7}^C$-invariant $C$-submodule of $\gP^C$. We first prove that if $(0, 0, 0, 1) \in W$, then we have $W = \gP^C$. Indeed,
$$
\begin{array}{l}
    W \ni {\mit\Phi}(0, X, 0, 0)(0, 0, 0, 1) = (X, 0, 0, 0),
\vspace{1mm}\\
    W \ni {\mit\Phi}(0, E_2, 0, 0)(E_3, 0, 0, 0) = (0, E_1, 0, 0),
\vspace{1mm}\\
    W \ni {\mit\Phi}(0, E_1, 0, 0)(0, E_1, 0, 0) = (0, 0, 1, 0),
\vspace{1mm}\\
    W \ni {\mit\Phi}(0, 0, Y, 0)(0, 0, 1, 0) = (0, Y, 0, 0).
\end{array}$$
Hence we have $W = \gP^C$. Now, let $P=(X, Y, \xi, \eta)$ be a non-zero element of $W$.
$$
\displaylines{\hspace*{5mm}\mbox{(i) Case} \; W \ni P = (X, Y, \xi, \eta), 
X \neq 0. \mbox{ Then we have}\hfill\mbox{(a)}}$$ 
\vspace{-10mm}
$$
\displaylines{\hfill W \ni {\mit\Phi}(0, 0, 0, 3)(X, Y, \xi, \eta) = 
(-X, Y, 3\xi, -3\eta), \hfill\mbox{(b)}}$$
\vspace{-10mm}
$$
\displaylines{\hfill W \ni {\mit\Phi}(0, 0, 0, 3)(-X, Y, 3\xi, -3\eta) = 
(X, Y, 9\xi, 9\eta).\hfill\mbox{(c)}}$$
Taking ((a)$-$(b))$\div 2$, ((a)$-$(c))$\div 8$, we have $(X, 0, -\xi, 2\eta) \in W$, 
$(0, 0, \xi, \eta) \in W$, respectively. Consequently $(X, 0, 0, 3\eta) \in 
W$. Next, if we choose $X_1 \in \gJ^C$ such that $(X_1, X) \neq 0$, from
$$
  W \ni {\mit\Phi}(0, 0, X_1, 0)(X, 0, 0, 3\eta) = (0, 0, 0, (X_1, X)), $$
we have $(0, 0, 0, 1) \in W$. Hence this case is reduced to the first case.
\vspace{1mm}

(ii) Case $P = (0, Y, \xi, \eta), Y \neq 0$. If we choose $B \in \gJ^C$ such that $B \times Y \neq 0$, from
$$
 W \ni {\mit\Phi}(0, 0, B, 0)(0, Y, \xi, \eta) = (2B \times Y, \xi B, 0, 0). $$
Hence this case is reduced to the case (i).
\vspace{1mm}

(iii) Case $P = (0, 0, \xi, \eta)$, $\xi \neq 0$. For $0 \neq B \in \gJ^C$, we have
$$
  W \ni {\mit\Phi}(0, 0, B, 0)(0, 0, \xi, \eta) = (0, \xi B, 0, 0). $$
Hence this case is also reduced to the case (ii).
\vspace{1mm}

\noindent Thus we have $W = \gP^C$.
\vspace{1mm}

(2) Since ${\gge_7}^C\gP^C$ is an ${\gge_7}^C$-invariant $C$-submodule of $\gP^C$, ${\gge_7}^C\gP^C = \gP^C$ follows from the irreducibility of $\gP^C$ (above (1)).
\vspace{3mm}

{\bf Lemma 4.4.3.} {\it Any element ${\mit\Phi} \in {\gge_7}^C$ is expressed by   ${\mit\Phi} = \dsum_{i}(P_i \times Q_i), \;\; P_i, Q_i \in \gP^C$.}
\vspace{2mm}

{\bf Proof.}  Since $[{\mit\Phi}, P \times Q] = {\mit\Phi}P \times Q + P \times  {\mit\Phi}Q$ (Proposition 4.3.2), $\ga = \Big\{ \dsum_{i}(P_i \times Q_i) \, | \, P_i, Q_i \in \gP^C \Big\}$ is an ideal of ${\gge_7}^C$. From the simplicity of ${\gge_7}^C$ (Theorem 4.4.1), we have $\ga = {\gge_7}^C$.
\vspace{4mm}

{\bf 4.5. Killing form of ${\gge_7}^C$}
\vspace{3mm}

{\bf Definition.} We define a symmetric inner product $({\mit\Phi}_1, {\mit\Phi}_2)_7$ in ${\gge_7}^C$ by
$$
   ({\mit\Phi}_1, {\mit\Phi}_2)_7 = -2(\phi_1, \phi_2)_6 - 4(A_1, B_2) - 4(A_2, B_1) - \dfrac{8}{3}\nu_1\nu_2, $$ 
where ${\mit\Phi}_i = {\mit\Phi}(\phi_i, A_i, B_i, \nu_i) \in {\gge_7}^C$.
\vspace{3mm}

{\bf Lemma 4.5.1.} (1) {\it The inner product $({\mit\Phi}_1, {\mit\Phi}_2)_7$ of ${\gge_7}^C$ is ${\gge_7}^C$-adjoint invariant}\,:
$$
  ([{\mit\Phi}, {\mit\Phi}_1], {\mit\Phi}_2)_7 + ({\mit\Phi}_1, [{\mit\Phi}, {\mit\Phi}_2])_7 = 0, \quad {\mit\Phi}, {\mit\Phi}_i \in {\gge_7}^C.$$

(2) {\it For} ${\mit\Phi} \in {\gge_7}^C, P, Q \in \gP^C$, {\it we have}
$$
       ({\mit\Phi}, P \times Q)_7 = \{{\mit\Phi}P, Q \}.  $$

$\begin{array}{l}{\bf Proof.}\; (1) \;
   ([{\mit\Phi}, {\mit\Phi}_1], {\mit\Phi}_2)_7 
\vspace{1mm}\\
 \qquad  = \Biggl({\mit\Phi}\pmatrix{[\phi, \phi_1] + 2A \vee B_1 - 2A_1 \vee B \vspace{1mm}\cr
       \Big(\phi + \dfrac{2}{3}\nu \Big)A_1 - \Big(\phi_1 + \dfrac{2}{3}\nu_1 \Big)A \vspace{1mm}\cr
       -\Big({}^t\phi + \dfrac{2}{3}\nu \Big)B_1 + \Big({}^t\phi_1 + \dfrac{2}{3}\nu_1 \Big)B \vspace{1mm}\cr
       (A, B_1)-(B, A_1)},
        {\mit\Phi}\pmatrix{\phi_2 \vspace{2mm}\cr 
                           A_2 \vspace{2mm}\cr 
                           B_2 \vspace{2mm}\cr \nu_2}\Biggl)_7 
\vspace{1mm}\\
\end{array}$
\vspace{1mm}

$ \qquad \;\,   = \cdots (\mbox{using $([\phi, A \vee B] = \phi A \vee B + A \vee \phi'B$ (Lemma 3.4.4) etc.} ) \cdots $
\vspace{1mm}

$ \qquad \;\, = -({\mit\Phi}_1, [{\mit\Phi}, {\mit\Phi}_2])_7.$
\vspace{1mm}

(2) For $P = (X, Y, \xi, \eta), Q = (Z, W, \zeta, \omega) \in \gP^C$, we have
\vspace{1mm}
   
$({\mit\Phi}, P \times Q)_7 
    = \Biggl({\mit\Phi}\pmatrix{\phi 
                                \vspace{1mm}\cr A 
                                \vspace{1mm}\cr B 
                                \vspace{1mm}\cr \nu}, 
    {\mit\Phi}\pmatrix{- \dfrac{1}{2}(X \vee W + Z \vee Y) \cr
                       - \dfrac{1}{4}(2Y \times W - \xi Z - \zeta X) \cr
                       \dfrac{1}{4}(2X \times Z - \eta W - \omega Y) \cr
            \dfrac{1}{8}((X, W) + (Z, Y) - 3(\xi\omega + \zeta\eta))} \Biggl)_7
$
\vspace{2mm}

$\quad  = (\phi, X \vee W + Z \vee Y)_6 - (A, 2X \times Z - \eta W - \omega Y) + (2Y \times W - \xi Z - \zeta X, B)$
\vspace{1mm}

$\quad \;\;\;- \dfrac{1}{3}\nu((X, W) + (Z, Y) - 3(\xi\omega + \zeta\eta))  $
\vspace{1mm}

$  \quad = (\phi X, W) + (\phi Z, Y) - 2(A, X, Z) + \eta(A, W) + \omega(A, Y) + 2(Y, W, B)$
\vspace{1mm}

$\quad \;\;\;- \xi(Z, B)-\zeta(X, B) - \dfrac{1}{3}\nu(X, W) - \dfrac{1}{3}\nu(Z, Y) + \nu(\xi\omega + \zeta\eta) $
\vspace{1mm}

$  \quad = \Biggl\{\pmatrix{\phi X - \dfrac{1}{3}\nu X + 2B \times Y + \eta A 
\vspace{1mm}\cr
               2A \times X - {}^t\phi \,Y + \dfrac{1}{3}\nu Y + \xi B 
\vspace{1mm}\cr
               (A, Y) + \nu\xi 
\vspace{1mm}\cr
               (B, X) - \nu\eta},
      \pmatrix{Z 
\vspace{1mm}\cr W 
\vspace{1mm}\cr \zeta 
\vspace{1mm}\cr \omega} \Biggl\}= \{\Phi P, Q \}. $
\vspace{3mm}

{\bf Theorem 4.5.2.} {\it The Killing form $B_7$ of the Lie algebra ${\gge_7}^C$  is given by}
\begin{eqnarray*}
B_7 ({\mit\Phi}_1, {\mit\Phi}_2) \!\!\! &=& \!\!\! - 9({\mit\Phi}_1, {\mit\Phi}_2)_7
\vspace{1mm}\\
  \!\!\! &=& \!\!\! 18(\phi_1, \phi_2)_6 + 36(A_1, B_2) + 36(A_2, B_1) + 24\nu_1 \nu_2
\vspace{1mm}\\
  \!\!\! &=& \!\!\! \dfrac{3}{2}B_6(\phi_1, \phi_2) + 36(A_1, B_2) + 36(A_2, B_1) + 24\nu_1 \nu_2
\vspace{1mm}\\
  \!\!\! &=& \!\!\! 3\tr({\mit\Phi}_1 {\mit\Phi}_2),
\end{eqnarray*}
{\it where} ${\mit\Phi}_i = {\mit\Phi}(\phi_i, A_i, B_i, \nu_i) \in {\gge_7}^C$ {\it and $B_6$ is the Killing form of ${\gge_6}^C$}.
\vspace{2mm}

{\bf Proof.} Since ${\gge_7}^C$ is simple (Theorem 4.4.1), there exist $k, k' \in C$ such that
$$
    B_7({\mit\Phi}_1, {\mit\Phi}_2) = k({\mit\Phi}_1, {\mit\Phi}_2)_7 
        = k'\tr({\mit\Phi}_1 {\mit\Phi}_2).$$
To determine these $k, k'$, let ${\mit\Phi}_0 = {\mit\Phi}_1 = {\mit\Phi}_2 = {\mit\Phi}(0, 0, 0, 1)$. Then we have
$$
        ({\mit\Phi}_0, {\mit\Phi}_0)_7 = -\dfrac{8}{3}.$$
On the other hand, we have
$$
    [{\mit\Phi}_0, [{\mit\Phi}_0, {\mit\Phi}(\phi, A, B, \nu)]\,] = \Big[{\mit\Phi}_0, {\mit\Phi}\Big(0, \dfrac{2}{3}A, -\dfrac{2}{3}B, 0 \Big)\Big] = {\mit\Phi}\Big(0, \dfrac{4}{9}A, \dfrac{4}{9}B, 0 \Big). $$
Hence 
$$
      B_7({\mit\Phi}_0, {\mit\Phi}_0) = \tr((\mbox{ad}{\mit\Phi}_0)^2) = \dfrac{4}{9} \times 27 \times 2 = 24. $$
Therefore $k = -9$. Next, from
$$
      {\mit\Phi}_0{\mit\Phi}_0(X, Y, \xi, \eta) = {\mit\Phi}_0\Big(-\dfrac{X}{3}, \dfrac{Y}{3}, \xi, -\eta \Big) = \Big(\dfrac{X}{9}, \dfrac{Y}{9}, \xi, \eta \Big),$$
we have 
$$
         \tr({\mit\Phi}_0^{\ 2}) = \dfrac{1}{9} \times 27 \times 2 + 1 + 1 = 8. $$
Therefore $k' = 3$.
\vspace{3mm}

{\bf Lemma 4.5.3.} {\it For $P \in \gP^C, P \neq 0$, there exists $Q \in \gP^C$ such that $P \times Q \vspace{2mm}
\neq 0$. }

{\bf Proof.} Asumme that $P \times Q = 0$ for all $Q \in \gP^C$. Then for any ${\mit\Phi} \in {\gge_7}^C$, $0 = ({\mit\Phi}, P \times Q)_7 = ({\mit\Phi}, Q \times P)_7 = \{{\mit\Phi}Q, P\}$ (Lemma 4.5.1.(2)). Since ${\gge_7}^C\gP^C = \gP^C$ (Proposition 4.4.2.(2)), we have $\{\gP^C, P\} = 0$, so that $P = 0$.
\vspace{4mm}

{\bf 4.6.  Roots of ${\gge_7}^C$}
\vspace{3mm}

{\bf Theorem 4.6.1} {\it The rank of the Lie algebra ${\gge_7}^C$ is} 7. {\it The roots of ${\gge_7}^C$ relative to some Cartan subalgebra $\gh$ are given by}
$$
\begin{array}{cl}
     \pm(\lambda_k - \lambda_l), \quad \pm(\lambda_k + \lambda_l), & 0 \le k < l \le 3, 
\vspace{1mm}\\
     \pm \lambda_k \pm \dfrac{1}{2}(\mu_2 - \mu_3), & 0 \le k \le 3,
\end{array} $$
\vspace{-2mm}
$$
\begin{array}{l}
   \pm \dfrac{1}{2}(- \lambda_0 - \lambda_1 + \lambda_2 - \lambda_3) \pm \dfrac{1}{2}(\mu_3 - \mu_1), 
\vspace{1mm}\\
   \pm \dfrac{1}{2}(\;\;\; \lambda_0 + \lambda_1 + \lambda_2 - \lambda_3) \pm \dfrac{1}{2}(\mu_3 - \mu_1),  
\vspace{1mm}\\
   \pm \dfrac{1}{2}(- \lambda_0 + \lambda_1 + \lambda_2 + \lambda_3) \pm \dfrac{1}{2}(\mu_3 - \mu_1),   
\vspace{1mm}\\
   \pm \dfrac{1}{2}(\;\;\; \lambda_0 - \lambda_1 + \lambda_2 + \lambda_3) \pm \dfrac{1}{2}(\mu_3 - \mu_1),   
\end{array}$$
\vspace{-1mm}
$$
\begin{array}{l}
   \pm \dfrac{1}{2}(\;\;\; \lambda_0 - \lambda_1 + \lambda_2 - \lambda_3) \pm \dfrac{1}{2}(\mu_1 - \mu_2), 
\vspace{1mm}\\
   \pm \dfrac{1}{2}(- \lambda_0 + \lambda_1 + \lambda_2 - \lambda_3) \pm \dfrac{1}{2}(\mu_1 - \mu_2),  
\vspace{1mm}\\
   \pm \dfrac{1}{2}(\;\;\; \lambda_0 + \lambda_1 + \lambda_2 + \lambda_3) \pm \dfrac{1}{2}(\mu_1 - \mu_2),   
\vspace{1mm}\\
   \pm \dfrac{1}{2}(- \lambda_0 - \lambda_1 + \lambda_2 + \lambda_3) \pm \dfrac{1}{2}(\mu_1 - \mu_2),   
\end{array}$$
\vspace{-1mm}
$$
\begin{array}{c}
     \pm\Big(\mu_j + \dfrac{2}{3}\nu\Big), \quad 0 \le j \le 3, 
\vspace{1mm}\\
     \pm \lambda_k \pm \Big(\dfrac{1}{2}\mu_1 - \dfrac{2}{3}\nu \Big), \quad 0 \le k \le 3,
\end{array} $$
\vspace{-1mm}
$$
\begin{array}{l}
   \pm \dfrac{1}{2}(- \lambda_0 - \lambda_1 + \lambda_2 - \lambda_3) \pm \Big(\dfrac{1}{2}\mu_2 - \dfrac{2}{3}\nu \Big), 
\vspace{1mm}\\
   \pm \dfrac{1}{2}(\;\;\; \lambda_0 + \lambda_1 + \lambda_2 - \lambda_3) \pm \Big(\dfrac{1}{2}\mu_2 - \dfrac{2}{3}\nu \Big), 
\end{array}$$
$$
\begin{array}{l}
   \pm \dfrac{1}{2}(- \lambda_0 + \lambda_1 + \lambda_2 + \lambda_3) \pm \Big(\dfrac{1}{2}\mu_2 - \dfrac{2}{3}\nu \Big),  
\vspace{1mm}\\
   \pm \dfrac{1}{2}(\;\;\; \lambda_0 - \lambda_1 + \lambda_2 + \lambda_3) \pm \Big(\dfrac{1}{2}\mu_2 - \dfrac{2}{3}\nu \Big), 
\end{array} $$
\vspace{-1mm}
$$
\begin{array}{l}
   \pm \dfrac{1}{2}(\;\;\; \lambda_0 - \lambda_1 + \lambda_2 - \lambda_3) \pm \Big(\dfrac{1}{2}\mu_3 - \dfrac{2}{3}\nu \Big), 
\vspace{1mm}\\
   \pm \dfrac{1}{2}(- \lambda_0 + \lambda_1 + \lambda_2 - \lambda_3) \pm \Big(\dfrac{1}{2}\mu_3 - \dfrac{2}{3}\nu \Big), 
\vspace{1mm}\\
   \pm \dfrac{1}{2}(\;\;\; \lambda_0 + \lambda_1 + \lambda_2 + \lambda_3) \pm \Big(\dfrac{1}{2}\mu_3 - \dfrac{2}{3}\nu \Big),  
\vspace{1mm}\\
   \pm \dfrac{1}{2}(- \lambda_0 - \lambda_1 + \lambda_2 + \lambda_3) \pm \Big(\dfrac{1}{2}\mu_3 - \dfrac{2}{3}\nu \Big) 
\end{array} $$
{\it with} $\mu_1 + \mu_2 + \mu_3 = 0.$ 
\vspace{2mm}

{\bf Proof.} We use the decomposition of of ${\gge_7}^C$ in Theorem 4.3.1 
$$
   {\gge_7}^C = {\gge_6}^C \oplus \gJ^C \oplus \gJ^C \oplus C.$$
Let
$$
   \gh = \Big\{{\mit\Phi}\Big(\dsum_{k=0}^3\lambda_kH_k + \Big(\dsum_{j=1}^3\mu_jE_j\Big)^{\sim}, 0, 0, \nu \Big) \in {\gge_7}^C\, \left|
\begin{array}{l}
   \lambda_k, \nu \in C \\
   \mu_j \in C, \mu_1 + \mu_2 + \mu_3 = 0 
\end{array} \right. \Big\} $$
(where $H_k = -iG_{k4+k})$, then $\gh$ is an abelian subalgebra of ${\gge_7}^C$ (it will be a Cartan subalgebra of ${\gge_7}^C$). In the following calculations, we put $h_\delta = \dsum_{k=0}^3\lambda_kH_k, H\vspace{-1mm}
 = \dsum_{j=1}^3\mu_jE_j.$ 

I $\,$ The roots of ${\gge_6}^C$ are also roots of ${\gge_7}^C$. Indeed, let $\alpha$ be a root of ${\gge_6}^C$ and $S \in {\gge_6}^C \subset {\gge_7}^C$ be a root vector belonging to $\alpha$. Then
\begin{eqnarray*}
     [h, S] \!\!\!&=&\!\!\!
     [{\mit\Phi}(h_\delta + \wti{H}, 0, 0, \nu), {\mit\Phi}(S, 0, 0, 0)] 
\vspace{1mm}\\
    \!\!\!&=&\!\!\!
    {\mit\Phi}([h_\delta + \wti{H}, S ], 0, 0, 0) = {\mit\Phi}(\alpha(h_\delta + \wti{H})S, 0, 0, 0) = \alpha(h)S.
\end{eqnarray*}
Hence $\alpha$ is a root of ${\gge_7}^C$.
\vspace{1mm}

II \quad $[{\mit\Phi}(h_{\delta} + \wti{H}, 0, 0, \nu), {\mit\Phi}(0, E_j, 0, 0)]$
\vspace{1mm}

\qquad \qquad = ${\mit\Phi}\Big(0,\Big(h_{\delta} + \wti{H} + \dfrac{2}{3}\nu \Big)E_j, 0, 0 \Big) = \Big(\mu_j + \dfrac{2}{3}\nu \Big){\mit\Phi}(0, E_j, 0, 0)$,
\vspace{1mm}

  \quad \quad $\; [{\mit\Phi}((h_{\delta} + \wti{H}, 0, 0, \nu), {\mit\Phi}(0, 0, E_j, 0)]$ = ${\mit\Phi}\Big(0, 0, \Big((h_{\delta} + \wti{H})' - \dfrac{2}{3}\nu \Big)E_j, 0 \Big)$
\vspace{1mm}

\qquad \qquad = ${\mit\Phi}\Big(0, 0, \Big(h_{\delta} - \wti{H} - \dfrac{2}{3}\nu \Big)E_j, 0 \Big) = \Big(- \mu_j - \dfrac{2}{3}\nu \Big){\mit\Phi}(0, 0, E_j, 0)$.
\vspace{1mm}

\noindent Hence $\pm\Big(\mu_j + \dfrac{2}{3}\nu\Big)$, $0 \le j \le 3$, are roots of ${\gge_7}^C$.
\vspace{1mm}

III  \quad $[{\mit\Phi}(h_{\delta} + \wti{H}, 0, 0, \nu), {\mit\Phi}(0, F_1(a), 0, 0)] \qquad a = e_k \pm ie_{4+k}$ 
\vspace{1mm}

\qquad \qquad $= {\mit\Phi}\Big(0, \Big(h_{\delta} + \wti{H} + \dfrac{2}{3}\nu \Big)F_1(a), 0, 0 \Big)$ 
\vspace{1mm}

\qquad \qquad $= {\mit\Phi}\Big(0, F_1(h_{\delta}a) + \dfrac{1}{2}(\mu_2 + \mu_3)F_1(a) + \dfrac{2}{3}\nu F_1(a), 0, 0 \Big)$ 
\vspace{1mm}

\qquad \qquad $= \Big( \pm \lambda_k - \dfrac{1}{2}\mu_1 + \dfrac{2}{3}\nu \Big){\mit\Phi}(0, F_1(a), 0, 0)$.
\vspace{1mm}

\noindent Hence $\pm \lambda_k - \dfrac{1}{2}\mu_1 + \dfrac{2}{3}\nu$, $0 \le k \le 3$, are roots of ${\gge_7}^C$. The remainders of roots can be similarly found ($\nu h_{\delta}$, $\kappa\pi h_{\delta}$ in Theorem 3.6.4 will be used).
\vspace{3mm}

{\bf Theorem 4.6.2.} {\it In the root system of Theorem} 4.6.1, 
\vspace{2mm}

\qquad \qquad \quad \qquad
       $\alpha_1 = \lambda_0 - \lambda_1, \quad 
       \alpha_2 = \lambda_1 - \lambda_2, \quad
       \alpha_3 = \lambda_2 - \lambda_3$,
\vspace{1mm}

\qquad \qquad \quad \qquad
       $\alpha_4 = \dfrac{1}{2}(- \lambda_0 - \lambda_1 - \lambda_2 + \lambda_3) + \dfrac{1}{2}(\mu_3 - \mu_1)$,
\vspace{1mm}

\qquad \qquad \quad \qquad
       $\alpha_5 = \dfrac{1}{2}(\lambda_0 + \lambda_1 + \lambda_2 + \lambda_3) +  \dfrac{1}{2}(\mu_1 - \mu_2)$,
\vspace{1mm}

\qquad \qquad \qquad \quad
       $\alpha_6 = \mu_2 + \dfrac{3}{2}\nu, \quad
       \alpha_7 = -\mu_3 - \dfrac{3}{2}\nu$
\vspace{2mm}

\noindent {\it is a fundamental root system of the Lie algebra ${\gge_7}^C$ and 
$$
     \mu = \alpha_1 + 2\alpha_2 + 3\alpha_3 + 4\alpha_4 + 3\alpha_5 + 2\alpha_6 + 2\alpha_7 $$
is the highest root. The Dynkin diagram and the extended Dynkin diagram of ${\gge_7}^C$ are respectively given by }
\vspace{-2mm}

\vspace{-2mm}

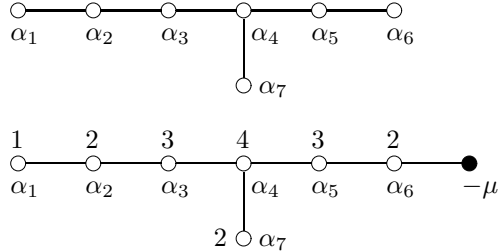
\setlength{\unitlength}{1mm}
\begin{picture}(100,20)
\put(20,10){\circle{2}} \put(19,6){$\alpha_1$}
\put(21,10){\line(1,0){8}}
\put(30,10){\circle{2}} \put(29,6){$\alpha_2$}
\put(31,10){\line(1,0){8}}
\put(40,10){\circle{2}} \put(39,6){$\alpha_3$}
\put(41,10){\line(1,0){8}}
\put(50,10){\circle{2}} \put(51,6){$\alpha_4$}
\put(51,10){\line(1,0){8}}
\put(50,8.8){\line(0,-1){7.8}}
\put(50,0){\circle{2}} \put(52,-1){$\alpha_7$}
\put(51,10){\line(1,0){8}}
\put(60,10){\circle{2}} \put(59,6){$\alpha_5$}
\put(61,10){\line(1,0){8}}
\put(70,10){\circle{2}} \put(69,6){$\alpha_6$}
\end{picture}

\setlength{\unitlength}{1mm}
\begin{picture}(100,20)
\put(20,10){\circle{2}} \put(19,6){$\alpha_1$} \put(19,12){$1$}
\put(21,10){\line(1,0){8}}
\put(30,10){\circle{2}} \put(29,6){$\alpha_2$} \put(29,12){$2$}
\put(31,10){\line(1,0){8}}
\put(40,10){\circle{2}} \put(39,6){$\alpha_3$} \put(39,12){$3$}
\put(41,10){\line(1,0){8}}
\put(50,8.8){\line(0,-1){7.8}}
\put(50,10){\circle{2}} \put(51,6){$\alpha_4$} \put(49,12){$4$}
\put(51,10){\line(1,0){8}}
\put(50,8.8){\line(0,-1){7.8}}
\put(50,0){\circle{2}} \put(52,-1){$\alpha_7$} \put(46,-1){$2$}
\put(51,10){\line(1,0){8}}
\put(60,10){\circle{2}} \put(59,6){$\alpha_5$} \put(59,12){$3$}
\put(61,10){\line(1,0){8}}
\put(70,10){\circle{2}} \put(69,6){$\alpha_6$} \put(69,12){$2$}
\put(71,10){\line(1,0){8}}
\put(80,10){\circle*{2}} \put(79,6){$-\mu$}
\end{picture}
\vspace{4mm}

{\bf Proof.} In the following, the notation $n_1n_2 \cdots n_7$ denotes the root  $n_1\alpha_1 + n_2\alpha_2 + \cdots + n_7\alpha_7$. Now, all positive roots of ${\gge_7}^C$ are represented by 
$$
\begin{array}{lllllllll}
      \lambda_0 - \lambda_1 = 1 & 0 & 0 & 0 & 0 & 0 & 0
\vspace{1mm}\\
      \lambda_0 - \lambda_2 = 1 & 1 & 0 & 0 & 0 & 0 & 0
\vspace{1mm}\\
      \lambda_0 - \lambda_3 = 1 & 1 & 1 & 0 & 0 & 0 & 0
\vspace{1mm}\\
      \lambda_1 - \lambda_2 = 0 & 1 & 0 & 0 & 0 & 0 & 0
\vspace{1mm}\\
      \lambda_1 - \lambda_3 = 0 & 1 & 1 & 0 & 0 & 0 & 0
\vspace{1mm}\\                            
      \lambda_2 - \lambda_3 = 0 & 0 & 1 & 0 & 0 & 0 & 0
\end{array}
\quad
\begin{array}{lllllllll}
      \lambda_0 + \lambda_1 = 1 & 2 & 2 & 2 & 2 & 1 & 1
\vspace{1mm}\\
      \lambda_0 + \lambda_2 = 1 & 1 & 2 & 2 & 2 & 1 & 1
\vspace{1mm}\\
      \lambda_0 + \lambda_3 = 1 & 1 & 1 & 2 & 2 & 1 & 1
\vspace{1mm}\\
      \lambda_1 + \lambda_2 = 0 & 1 & 2 & 2 & 2 & 1 & 1
\vspace{1mm}\\
      \lambda_1 + \lambda_3 = 0 & 1 & 1 & 2 & 2 & 1 & 1 
\vspace{1mm}\\
      \lambda_2 + \lambda_3 = 0 & 0 & 1 & 2 & 2 & 1 & 1
\end{array} $$
\vspace{-1mm}
$$
\begin{array}{llllllllll}
      \lambda_0 + \dfrac{1}{2}(\mu_2 - \mu_3) = 1 & 1 & 1 & 1 & 1 & 1 & 1
\vspace{1mm}\\
      \lambda_1 + \dfrac{1}{2}(\mu_2 - \mu_3) = 0 & 1 & 1 & 1 & 1 & 1 & 1
\vspace{1mm}\\
      \lambda_2 + \dfrac{1}{2}(\mu_2 - \mu_3) = 0 & 0 & 1 & 1 & 1 & 1 & 1
\vspace{1mm}\\
      \lambda_3 + \dfrac{1}{2}(\mu_2 - \mu_3) = 0 & 0 & 0 & 1 & 1 & 1 & 1
\end{array}$$
\vspace{-1mm}
$$
\begin{array}{lllllllllll}
      \lambda_0 - \dfrac{1}{2}(\mu_2 - \mu_3) = 1 & 1 & 1 & 1 & 1 & 0 & 0
\vspace{1mm}\\
      \lambda_1 - \dfrac{1}{2}(\mu_2 - \mu_3) = 0 & 1 & 1 & 1 & 1 & 0 & 0
\vspace{1mm}\\
      \lambda_2 - \dfrac{1}{2}(\mu_2 - \mu_3) = 0 & 0 & 1 & 1 & 1 & 0 & 0
\vspace{1mm}\\
      \lambda_3 - \dfrac{1}{2}(\mu_2 - \mu_3) = 0 & 0 & 0 & 1 & 1 & 0 & 0
\end{array}$$
$$
\begin{array}{llllllllll}
     \dfrac{1}{2}(- \lambda_0 - \lambda_1 + \lambda_2 - \lambda_3) + \dfrac{1}{2}(\mu_3 - \mu_1) = 0 & 0 & 1 & 1 & 0 & 0 & 0 
\vspace{1mm}\\
     \dfrac{1}{2}(\;\;\; \lambda_0 + \lambda_1 + \lambda_2 - \lambda_3) + \dfrac{1}{2}(\mu_3 - \mu_1) = 1 & 2 & 3 & 3 & 2 & 1 & 1
\vspace{1mm}\\
     \dfrac{1}{2}(- \lambda_0 + \lambda_1 + \lambda_2 + \lambda_3) + \dfrac{1}{2}(\mu_3 - \mu_1) = 0 & 1 & 2 & 3 & 2 & 1 & 1
\vspace{1mm}\\
     \dfrac{1}{2}(\;\;\; \lambda_0 - \lambda_1 + \lambda_2 + \lambda_3) + \dfrac{1}{2}(\mu_3 - \mu_1) = 1 & 1 & 2 & 3 & 2 & 1 & 1 
\end{array}$$
\vspace{-1mm}
$$
\begin{array}{lllllllllll}
     \dfrac{1}{2}(\;\;\; \lambda_0 + \lambda_1 - \lambda_2 + \lambda_3) + \dfrac{1}{2}(\mu_3 - \mu_1) = 1 & 2 & 2 & 3 & 2 & 1 & 1 
\vspace{1mm}\\
     \dfrac{1}{2}(- \lambda_0 - \lambda_1 - \lambda_2 + \lambda_3) + \dfrac{1}{2}(\mu_3 - \mu_1) = 0 & 0 & 0 & 1 & 0 & 0 & 0
\vspace{1mm}\\
     \dfrac{1}{2}(\;\;\; \lambda_0 - \lambda_1 - \lambda_2 - \lambda_3) + \dfrac{1}{2}(\mu_3 - \mu_1) = 1 & 1 & 1 & 1 & 0 & 0 & 0 
\vspace{1mm}\\
     \dfrac{1}{2}(- \lambda_0 + \lambda_1 - \lambda_2 - \lambda_3) + \dfrac{1}{2}(\mu_3 - \mu_1) = 0 & 1 & 1 & 1 & 0 & 0 & 0
\end{array} $$
$$
\begin{array}{llllllll}
   \dfrac{1}{2}(\;\;\; \lambda_0 - \lambda_1 + \lambda_2 - \lambda_3) - \dfrac{1}{2}(\mu_1 - \mu_2) = 1 & 1 & 2 & 2 & 1 & 1 & 1 
\vspace{1mm}\\
   \dfrac{1}{2}(\;\;\; \lambda_0 - \lambda_1 - \lambda_2 + \lambda_3) - \dfrac{1}{2}(\mu_1-\mu_2) = 1 & 1 & 1 & 2 & 1 & 1 & 1 
\vspace{1mm}\\
   \dfrac{1}{2}(\;\;\; \lambda_0 + \lambda_1 + \lambda_2 + \lambda_3) + \dfrac{1}{2}(\mu_1 - \mu_2) = 0 & 0 & 0 & 0 & 1 & 0 & 0
\vspace{1mm}\\
   \dfrac{1}{2}(\;\;\; \lambda_0 + \lambda_1 - \lambda_2 - \lambda_3) - \frac{1}{2}(\mu_1 - \mu_2) = 1 & 2 & 2 & 2 & 1 & 1 & 1 
\end{array}$$
\vspace{-1mm}
$$
\begin{array}{lllllllllll}
   \dfrac{1}{2}(- \lambda_0 + \lambda_1 - \lambda_2 + \lambda_3) - \dfrac{1}{2}(\mu_1-\mu_2) = 0 & 1 & 1 & 2 & 1 & 1 & 1 
\vspace{1mm}\\
   \dfrac{1}{2}(- \lambda_0 + \lambda_1 + \lambda_2 - \lambda_3) - \dfrac{1}{2}(\mu_1 - \mu_2) = 0 & 1 & 2 & 2 & 1 & 1 & 1
\vspace{1mm}\\
   \dfrac{1}{2}(\;\;\; \lambda_0 + \lambda_1 + \lambda_2 + \lambda_3) - \dfrac{1}{2}(\mu_1 - \mu_2) = 1 & 2 & 3 & 4 & 3 & 2 & 2 
\vspace{1mm}\\
   \dfrac{1}{2}(- \lambda_0 - \lambda_1 + \lambda_2 + \lambda_3) - \dfrac{1}{2}(\mu_1 - \mu_2) = 0 & 0 & 1 & 2 & 1 & 1 & 1
\end{array} $$
$$
\begin{array}{lllllllllllll}
       - \mu_1 - \dfrac{2}{3}\nu = 1 & 2 & 3 & 4 & 2 & 1 & 2
\vspace{1mm}\\
       \;\;\; \mu_2 + \dfrac{2}{3}\nu = 0 & 0 & 0 & 0 & 0 & 1 & 0
\vspace{1mm}\\
       - \mu_3 - \dfrac{2}{3}\nu = 0 & 0 & 0 & 0 & 0 & 0 & 1
\end{array}$$
$$
\begin{array}{lllllllllll}
      \lambda_0 - \dfrac{1}{2}\mu_1 + \dfrac{2}{3}\nu = 1&1&1&1&1&1&0
\vspace{1mm}\\
      \lambda_1 - \dfrac{1}{2}\mu_1 + \dfrac{2}{3}\nu = 0&1&1&1&1&1&0
\vspace{1mm}\\
      \lambda_2 - \dfrac{1}{2}\mu_1 + \dfrac{2}{3}\nu = 0&0&1&1&1&1&0
\vspace{1mm}\\
      \lambda_3 - \dfrac{1}{2}\mu_1 + \dfrac{2}{3}\nu = 0&0&0&1&1&1&0
\end{array}$$
\vspace{-1mm}
$$
\begin{array}{lllllllllll}
    \lambda_0 + \dfrac{1}{2}\mu_1 - \dfrac{2}{3}\nu = 1 & 1 & 1 & 1 & 1 & 0 & 1
\vspace{1mm}\\
    \lambda_1 + \dfrac{1}{2}\mu_1 - \dfrac{2}{3}\nu = 0 & 1 & 1 & 1 & 1 & 0 & 1
\vspace{1mm}\\
    \lambda_2 + \dfrac{1}{2}\mu_1 - \dfrac{2}{3}\nu = 0 & 0 & 1 & 1 & 1 & 0 & 1
\vspace{1mm}\\
    \lambda_3 + \dfrac{1}{2}\mu_1 - \dfrac{2}{3}\nu = 0 & 0 & 0 & 1 & 1 & 0 & 1
\end{array}$$
$$
\begin{array}{lllllllllll}
   \dfrac{1}{2}(- \lambda_0 - \lambda_1 + \lambda_2 - \lambda_3) + \dfrac{1}{2}\mu_2 - \dfrac{2}{3}\nu = 0 & 0 & 1 & 1 & 0 & 0 & 1
\vspace{1mm}\\
   \dfrac{1}{2}(\;\;\; \lambda_0 + \lambda_1 + \lambda_2 - \lambda_3) + \frac{1}{2}\mu_2 - \dfrac{2}{3}\nu = 1 & 2 & 3 & 3 & 2 & 1 & 2
\vspace{1mm}\\
   \dfrac{1}{2}(- \lambda_0 + \lambda_1 + \lambda_2 + \lambda_3) + \dfrac{1}{2}\mu_2 - \dfrac{2}{3}\nu = 0 & 1 & 2 & 3 & 2 & 1 & 2
\vspace{1mm}\\
   \dfrac{1}{2}(\;\;\; \lambda_0 - \lambda_1 + \lambda_2 + \lambda_3) + \frac{1}{2}\mu_2 - \dfrac{2}{3}\nu = 1 & 1 & 2 & 3 & 2 & 1 & 2
\end{array}$$
\vspace{-1mm}
$$
\begin{array}{lllllllllll}
   \dfrac{1}{2}(\;\;\; \lambda_0 + \lambda_1 - \lambda_2 + \lambda_3) + \dfrac{1}{2}\mu_2 - \dfrac{2}{3}\nu = 1 & 2 & 2 & 3 & 2 & 1 & 2
\vspace{1mm}\\
   \dfrac{1}{2}(- \lambda_0 - \lambda_1 - \lambda_2 + \lambda_3) + \dfrac{1}{2}\mu_2 - \dfrac{2}{3}\nu = 0 & 0 & 0 & 1 & 0 & 0 & 1
\vspace{1mm}\\
   \dfrac{1}{2}(\;\;\; \lambda_0 - \lambda_1 - \lambda_2 - \lambda_3) + \dfrac{1}{2}\mu_2  - \dfrac{2}{3}\nu = 1 & 1 & 1 & 1 & 0 & 0 & 1
\vspace{1mm}\\
   \dfrac{1}{2}(- \lambda_0 + \lambda_1 - \lambda_2 - \lambda_3) + \dfrac{1}{2}\mu_2 - \dfrac{2}{3}\nu = 0 & 1 & 1 & 1 & 0 & 0 & 1
\end{array}$$
$$
\begin{array}{lllllllllll}
   \dfrac{1}{2}(\;\;\; \lambda_0 - \lambda_1 + \lambda_2 - \lambda_3) + \dfrac{1}{2}\mu_3 - \dfrac{2}{3}\nu = 1 & 1 & 2 & 2 & 1 & 0 & 1
\vspace{1mm}\\
   \dfrac{1}{2}(\;\;\; \lambda_0 - \lambda_1 - \lambda_2 + \lambda_3) + \dfrac{1}{2}\mu_3 + \dfrac{2}{3}\nu = 1 & 1 & 1 & 2 & 1 & 0 & 1
\vspace{1mm}\\
   \dfrac{1}{2}(\;\;\; \lambda_0 + \lambda_1 + \lambda_2 + \lambda_3) - \dfrac{1}{2}\mu_3  -\dfrac{2}{3}\nu = 0 & 0 & 0 & 0 & 1 & 1 & 0
\vspace{1mm}\\
   \dfrac{1}{2}(\;\;\; \lambda_0 + \lambda_1 - \lambda_2 - \lambda_3) + \dfrac{1}{2}\mu_3 - \dfrac{2}{3}\nu = 1 & 2 & 2 & 2 & 1 & 0 & 1
\end{array}$$
\vspace{-1mm}
$$
\begin{array}{lllllllllll}
   \dfrac{1}{2}(- \lambda_0 + \lambda_1 - \lambda_2 + \lambda_3) + \dfrac{1}{2}\mu_3 - \dfrac{2}{3}\nu = 0 & 1 & 1 & 2 & 1 & 0 & 1
\vspace{1mm}\\
   \dfrac{1}{2}(- \lambda_0 + \lambda_1 + \lambda_2 - \lambda_3) + \dfrac{1}{2}\mu_3 - \dfrac{2}{3}\nu = 0 & 1 & 2 & 2 & 1 & 0 & 1
\vspace{1mm}\\
   \dfrac{1}{2}(\;\;\; \lambda_0 + \lambda_1 + \lambda_2 + \lambda_3) + \dfrac{1}{2}\mu_3 - \dfrac{2}{3}\nu = 1 & 2 & 3 & 4 & 3 & 1 & 2
\vspace{1mm}\\
   \dfrac{1}{2}(- \lambda_0 - \lambda_1 + \lambda_2 + \lambda_3) + \dfrac{1}{2}\mu_3 - \dfrac{2}{3}\nu = 0 & 0 & 1 & 2 & 1 & 0 & 1.
\end{array}$$
Hence ${\mit\Pi} = \{ \alpha_1, \alpha_2, \cdots, \alpha_7 \}$ is a fundamental root system of ${\gge_7}^C$. The real part $\gh_{\sR}$ of $\gh$ is
$$
   \gh_{\sR} = \Big\{{\mit\Phi}\Big(\dsum_{k=0}^3\lambda_kH_k + \Big(\dsum_{j=1}^3\mu_jE_j\Big)^{\sim}, 0, 0, \nu \Big) \in {\gge_7}^C \, \left|
\begin{array}{l}
   \lambda_k, \nu \in \R \\
   \mu_j \in \R, \mu_1 + \mu_2 + \mu_3 = 0 
\end{array} \right. \Big\} $$
and the Killing form $B_7$ of ${\gge_7}^C$ on $\gh_{\sR}$ is given by
$$
     B_7(h, h') = 6\Big(6\dsum_{k=0}^3\lambda_k{\lambda_k}' + 3\dsum_{j=1}^3\mu_j{\mu_j}' + 4\nu \nu'\Big) $$
for $h = {\mit\Phi}\Big(\dsum_{k=0}^3\lambda_kH_k + \Big(\dsum_{j=1}^3\mu_jE_j\Big)^\sim, 0, 0, \nu \Big)$, $h' = {\mit\Phi}\Big(\dsum_{k=0}^3{\lambda_k}'H_k + \Big(\dsum_{j=1}^3{\mu_j}'E_j\Big)^\sim, 0, $ $0, \nu'\Big) \in \gh_{\sR}$. Indeed, from Theorem 4.5.2, we have
\begin{eqnarray*}
   B_7(h, h') \!\!\!&=&\!\!\!
  \dfrac{3}{2}B_6\Big(\dsum_{k=0}^3\lambda_kH_k + \Big(\dsum_{j=1}^3\mu_jE_j\Big)^\sim,  \dsum_{k=0}^3{\lambda_k}'H_k + \Big(\dsum_{j=1}^3{\mu_j}'E_j\Big)^\sim\Big) + 24\nu\nu'
\vspace{1mm}\\
  \!\!\!&=&\!\!\! \dfrac{3}{2}12\Big(2\dsum_{k=0}^3\lambda_k{\lambda_k}' + \dsum_{j=1}^3\mu_j{\mu_j}'\Big) + 24\nu\nu' \; \mbox{(Theorem 3.6.5)}
\vspace{1mm}\\
 \!\!\!&=&\!\!\! 6\Big(6\dsum_{k=0}^3\lambda_k{\lambda_k}' + 3\dsum_{j=1}^3\mu_j{\mu_j}' + 4\nu\nu'\Big).
\end{eqnarray*}

\noindent Now, the canonical elements $H_{\alpha_i} \in \gh_{\sR}$ corresponding to $\alpha_i$ ($B_7 (H_\alpha, H) = \alpha(H), \; H \in \gh $) are determined as follows.
$$
\begin{array}{l}
    H_{\alpha_1} = \dfrac{1}{36}{\mit\Phi}(H_0 - H_1, 0, 0, 0), 
\vspace{1mm}\\
    H_{\alpha_2} = \dfrac{1}{36}{\mit\Phi}(H_1 - H_2, 0, 0, 0), 
\vspace{1mm}\\
    H_{\alpha_3} = \dfrac{1}{36}{\mit\Phi}(H_2 - H_3, 0, 0, 0),
\vspace{1mm}\\
    H_{\alpha_4} = \dfrac{1}{72}{\mit\Phi}((- H_0 - H_1 - H_2 + H_3) + 2(E_3 - E_1)^{\sim}, 0, 0, 0),
\vspace{1mm}\\
    H_{\alpha_5} = \dfrac{1}{72}{\mit\Phi}((H_0 + H_1 + H_2 + H_3) + 2(E_1 - E_2)^{\sim}, 0, 0, 0),
\vspace{1mm}\\
    H_{\alpha_6} = \dfrac{1}{54}{\mit\Phi}\Big((- E_1 + 2E_2 - E_3)^{\sim}, 0, 0, \dfrac{3}{2} \Big),
\vspace{1mm}\\
    H_{\alpha_7} = \dfrac{1}{54}{\mit\Phi}\Big((E_1 + E_2 - 2E_3)^{\sim}, 0, 0, - \dfrac{3}{2}\Big).
\end{array} $$
Thus we have
$$(\alpha_1, \alpha_1)=B_7(H_{\alpha_1}, H_{\alpha_1})=
36 \frac{1}{36} \frac{1}{36} 2=\frac{1}{18}$$
and the other inner products are similarly calculated. Consequently, the inner product induced by the Killing form $B_7$ between $\alpha_1, \alpha_2, \cdots, \alpha_7$ and $- \mu$ are given by

$$
\begin{array}{l}
      (\alpha_i, \alpha_i) = \dfrac{1}{18}, \quad i = 1, 2, 3, 4, 5, 6, 7,
\vspace{1mm}\\
      (\alpha_1, \alpha_2) = (\alpha_2, \alpha_3) = (\alpha_3, \alpha_4) = (\alpha_4, \alpha_5) = (\alpha_4, \alpha_7) = (\alpha_5, \alpha_6) = - \dfrac{1}{36},
\vspace{1mm}\\
       (\alpha_i, \alpha_j) = 0, \quad \mbox{otherwise},
\vspace{1mm}\\
       (-\mu, -\mu) = \dfrac{1}{18}, \quad (- \mu, \alpha_6) = - \dfrac{1}{36}, \quad (- \mu, \alpha_i) = 0, \;\; i = 1, 2, 3, 5, 7,
\end{array} $$
using them, we can draw the Dynkin diagram and the extended Dynkin diagram of ${\gge_7}^C$.
\vspace{2mm}

According to Borel-Siebenthal theory, the Lie algebra $\gge_7$ has four subalgebras as maximal subalgebras with the maximal rank 7. 
\vspace{1mm}

(1) The first one is a subalgebra of type $T \oplus E_6$ which is obtained as the fixed points of an involution $\iota$ of $\gge_7$.
\vspace{1mm}

(2) The second one is a subalgebra of type $A_1 \oplus D_6$ which is obtained as the fixed points of an involution $\sigma$ of $\gge_7$.
\vspace{1mm}

(3) The third one is a subalgebra of type $A_7$ which is obtained as the fixed points of an involution $\lambda\gamma$ of $\gge_7$.
\vspace{1mm}

(4) The fourth one is a subalgebra of type $A_2 \oplus A_5$ which is obtained as the fixed points of an automorphism $w$ of order 3 of $\gge_7$.
\vspace{1mm}

These subalgebras will be realized as subgroups of the group $E_7$ in Theorems 4.10.2, 4.11.15, 4.12.5 and 4.13.5, respectively.
\vspace{4mm}

{\bf 4.7. Subgroups $E_6$ and $U(1)$ of $E_7$}
\vspace{3mm}

We shall study the following subgroup $(E_7)_{(0, 0, 1, 0)}$ of $E_7$: 
$$
     (E_7)_{(0,0,1,0)} = \{\alpha \in E_7 \, | \, \alpha(0, 0, 1, 0) = (0, 0, 1, 0) \}. $$

{\bf Lemma 4.7.1.} {\it If $\alpha \in E_7$ satisfies $\alpha(0, 0, 1, 0) = (0, 0, 1, 0)$, then} $\alpha$ {\it also satisfies} $\alpha(0, 0, 0, 1) = (0, 0, 0, 1)$, {\it and conversely}.
\vspace{2mm}

{\bf Proof.} If $\alpha \in E_7$ satisfies $\alpha(0, 0, 1, 0) = (0, 0, 1, 0)$, then we have
$$
    \alpha(0, 0, 0, 1) = \alpha(- \tau\lambda(0, 0, 1, 0)) = - \tau\lambda\alpha(0, 0, 1, 0) = - \tau\lambda(0, 0, 1, 0) = (0, 0, 0, 1). $$
The converse can be similarly proved.
\vspace{3mm}

{\bf Theorem 4.7.2.} \qquad \qquad \quad $(E_7)_{(0, 0, 1, 0)} \cong E_6.$
\vspace{2mm}

{\bf Proof.} We associate an element $\alpha \in E_6 = \{ \alpha \in \Iso_C(\gJ^C) \, | \, \det\,(\alpha X) =\det\,X, \langle \alpha X,$ $ \alpha Y \rangle = \langle X, Y \rangle \}$ with the element

$$
      \wti{\alpha} = \pmatrix{\alpha & 0 & 0 & 0 \cr
                              0 & \tau\alpha\tau & 0 & 0 \cr
                              0 & 0 & 1 & 0 \cr
                              0 & 0 & 0 & 1}
            \in (E_7)_{(0, 0, 1, 0)} \subset E_7.$$ 
We first have to prove that $\wti{\alpha} \in E_7$. For $P = (X, Y, \xi, \eta)$, $Q = (Z, W, \zeta, \omega) \in \gP^C$, we have
$$
\begin{array}{l}
   \wti{\alpha}P \times \wti{\alpha}Q = (\alpha X, \tau\alpha\tau Y, \xi, \eta)
 \times (\alpha Z, \tau\alpha\tau W, \zeta, \omega)
\vspace{1mm}\\
    = \cdots \mbox{(using $\alpha(X \vee Y)\alpha^{-1} = \alpha X \vee \tau\alpha\tau Y, (\alpha\phi\alpha^{-1})' = (\tau\alpha\tau)\phi'(\tau\alpha^{-1}\tau)$ etc.)} \cdots 
\vspace{1mm}\\
   = \wti{\alpha}(P \times Q)\wti{\alpha}^{-1}
\end{array} $$
and $\langle \wti{\alpha}P, \wti{\alpha}Q \rangle = \langle P, Q \rangle$ is evident. Hence $\wti{\alpha} \in E_7$, moreover $\wti{\alpha} \in (E_7)_{(0,0,1,0)}$.
\vspace{1mm}

Conversely, suppose that $\alpha \in E_7$ satisfies $\alpha(0, 0, 1, 0) = (0, 0, 1, 0)$ and $\alpha(0, 0, 0, 1) = (0, 0, 0, 1)$ (Lemma 4.7.1). Then $\alpha$ is of the form
$$
     \alpha = \pmatrix{\beta & \epsilon & 0 & 0 \cr
                       \delta & \beta_1 & 0 & 0 \cr
                       0 & 0 & 1 & 0 \cr
                       0 & 0 & 0 & 1}, \quad
                    \beta, \beta_1, \delta, \epsilon \in \Hom_C(\gJ^C).$$
Indeed, the fact that the left bottom parts are $0$ follows from 
$$
\begin{array}{l}
   \langle \alpha\dot{X}, \dot{1} \rangle = \langle \alpha\dot{X}, \alpha\dot{1} \rangle = \langle \dot{X}, \dot{1} \rangle = 0, 
\vspace{1mm}\\ 
   \langle \alpha\dot{X}, \d{1} \rangle = \langle \alpha\dot {X}, \alpha\d{1} \rangle = \langle \dot{X}, \d{1} \rangle = 0.
\end{array} $$
Now, since
$$
    \gM^C \ni \alpha\pmatrix{X \vspace{1mm}
\vspace{1mm}\cr
                            \dfrac{1}{\eta}X \times X 
\vspace{1mm}\cr
                            \dfrac{1}{\eta^2}\det X 
\vspace{1mm}\cr
                            \eta}
     = \pmatrix{\beta X + \dfrac{1}{\eta}\epsilon(X \times X) 
\vspace{1mm}\cr 
                \delta X + \dfrac{1}{\eta}\beta_1(X \times X)
 \vspace{1mm}\cr
                \dfrac{1}{\eta^2}\det\,X 
\vspace{1mm}\cr
                \eta},$$
we can see that
$$
     \Big(\beta X + \dfrac{1}{\eta}\epsilon(X \times X)\Big) \times \Big(\beta X + \dfrac{1}{\eta}\epsilon(X \times X)\Big) = \eta\Big(\delta X + \dfrac{1}{\eta}\beta_1(X \times X)\Big) $$
holds for all $0 \neq \eta \in C$. Comparing the coefficients of $\eta$ of both sides, we have 
\vspace{0.8mm}
$\delta = 0$. Similarly, from $\alpha\Big(\dfrac{1}{\xi}(Y \times Y), Y, \xi, \dfrac{1}{\xi^2} \det\,Y \Big) \in \gM^C$, we have $\epsilon = 0$. Furthermore, since
$$
    \gM^C \ni \alpha(X, X \times X, \det\,X, 1) = (\beta X, \beta_1(X \times X), \det\,X, 1),$$
we have
$$
             \beta X \times \beta X = \beta_1(X \times X), \quad
             (\beta X, \beta_1(X \times X)) = 3\det\,X. $$
and so
$$
        \det\,(\beta X) = \dfrac{1}{3}(\beta X, \beta X \times \beta X) = 
                   \dfrac{1}{3}(\beta X, \beta_1(X \times X)) = \det\,X,$$
which implies that $\beta \in {E_6}^C$. The equality $\langle \alpha\dot{X}, \alpha\dot{Y} \rangle = \langle \dot{X}, \dot{Y} \rangle$ implies $\langle \beta X, \beta Y \rangle = \langle X, Y \rangle$ and therefore $\beta \in E_6$. Moreover from the relation
$$
     \beta_1(X \times X) = \beta X \times \beta X = \tau\beta\tau(X \times X),$$
we obtain $\beta_1 = \tau\beta\tau$. Indeed, putting $X \times X$ instead of $X$, we have
$$
         (\det\,X)\beta_1X = (\det\,X)\tau\beta\tau X.$$
If $\det\,X \neq 0$, then we have $\beta_1X = \tau\beta\tau X$, Since $\beta_1$ and $\tau\beta\tau$ are linear mappings (of course are continuous mappings), we have $\beta_1X = \tau\beta\tau X$ even if $\det\,X = 0$. Thus, the proof of Theorem 4.7.2 is completed.
\vspace{2mm}

For $\theta \in C, \theta \neq 0$, we define a $C$-linear transformation $\varphi_1(\theta) : \gP^C \to \gP^C$ by
$$
        \varphi_1(\theta)(X, Y, \xi, \eta) = (\theta^{-1}X, \theta Y, \theta^3\xi, \theta^{-3}\eta). $$
Then $\varphi_1(\theta) \in {E_7}^C.$
\vspace{3mm}

{\bf Theorem 4.7.3.} {\it The group} $E_7$ {\it contains}
$$
      U(1) = \{ \varphi_1(\theta) \, | \, \theta \in C, (\tau\theta)\theta = 1 \} $$
{\it as a subargroup. This subgroup is isomorphic to the usual unitary group} $U(1) = \{ \theta \in C \, | \, (\tau\theta)\theta = 1 \}$.
\vspace{2mm}

{\bf Proof.} It is easy to check that $\varphi_1(\theta) \in E_7$.
\vspace{4mm}

{\bf 4.8. Connectedness of $E_7$}
\vspace{3mm}

We denote by $(E_7)_0$ the connected component of $E_7$ containing the identity $1$. 
\vspace{3mm}

{\bf Lemma 4.8.1.} {\it For} $a \in C$, {\it if we define a mapping} $\alpha_i(a) : 
\gP^C \to \gP^C$, $i=1, 2, 3$ {\it by}
$$
    \alpha_i(a) = \pmatrix{1 + (\cos|a| - 1)p_i & 2a\dfrac{\sin|a|}{|a|}E_i & 0                   & - \tau a\dfrac{\sin|a|}{|a|}E_i \cr
       - 2\tau a\dfrac{\sin|a|}{|a|}E_i & 1 + (\cos|a| - 1)p_i & a\dfrac{\sin|a|          }{|a|}E_i & 0 \cr
            0 & - \tau a\dfrac{\sin|a|}{|a|}E_i & \cos|a| & 0 \cr
           a\dfrac{\sin|a|}{|a|}E_i & 0 & 0 & \cos|a|} $$
\Big({\it if $a = 0$, then $\dfrac{\sin|a|}{|a|}$ means $1$}\Big), {\it then} $\alpha_i(a) \in (E_7)_0$, {\it where $p_i : \gJ^C \to \gJ^C$ is the $C$-linear mapping defined by}

$$
      p_i\pmatrix{\xi_1 & x_3 & \ov{x}_2 \cr
                  \ov{x}_3 & \xi_2 & x_1 \cr
                  x_2 & \ov{x}_1 & \xi_3}
       = \pmatrix{\xi_1 & \delta_{i3}x_3 & \delta_{i2}\ov{x}_2 \cr
                  \delta_{i3}\ov{x}_3 & \xi_2 & \delta_{i1}x_1 \cr
                  \delta_{i2}x_2 & \delta_{i1}\ov{x}_1 & \xi_3}, $$
{\it where $\delta_{ij}$ is the Kronecker delta symble. The mappings $\alpha_1(a_1), \alpha_2(a_2), \alpha(a_3), (a_i \in C)$ are commutative for each other.}
\vspace{2mm}

{\bf Proof.} For 
$$
    {\mit\Phi}_i(a) = {\mit\Phi}(0, - \tau aE_i, \tau aE_i, 0 ) = 
                       \pmatrix{0 & 2aE_i & 0 & -\tau a E_i \cr
                                -2\tau aE_i & 0 & a E_i & 0 \cr
                                 0 & -\tau aE_i & 0 & 0 \cr
                                 aE_i & 0 & 0 & 0} \in \gge_7$$
(Theorem 4.3.4), we have $\alpha_i(a) = \exp{\mit\Phi}_i(a)$. Hence $\alpha_i(a) \in (E_7)_0$. The relation $[{\mit\Phi}_i(a_i), {\mit\Phi}_j(a_j)] = 0$ shows that $\alpha_i(a_i)$ and $\alpha_j(a_j)$ are commutative.
\vspace{3mm}

{\bf Proposition 4.8.2.} {\it Any element $P \in \gM^C$, $P \neq 0$ can be transformed to a diagonal form by some element $\alpha \in (E_7)_0$}:
$$
        \alpha P = (X, Y, \xi, \eta), \quad  X, Y \; \mbox{\it are diagonal},\; \xi > 0. $$

{\bf Proof.} Let $P = (X, Y, \xi, \eta) \in \gM^C$. We shall first show that $P$  can be transformed to a diagonal form with $\xi \neq 0$.
\vspace{1mm}

(1) Case $P = (X, Y, \xi, \eta)$, $\xi \neq 0$. In this case, $X = \dfrac{1}{\xi}(Y \times Y)$. Choose $\beta \in E_6$ such that $\tau\beta\tau Y$ is diagonal (Proposition 3.8.2), then
$$
           \beta X = \dfrac{1}{\xi}\beta(Y \times Y)
               = \dfrac{1}{\xi}\tau\beta\tau Y \times \tau\beta\tau Y $$
is also diagonal.
\vspace{1mm}

(2) Case $P=(X, Y, 0, \eta)$, $Y \neq 0$. Choose $\beta \in E_6$ so that
$$
       \tau\beta\tau Y = \pmatrix{\eta_1 & 0 & 0 \cr
                                  0 & \eta_2 & 0 \cr
                                  0 & 0 & \eta_3}, \quad \eta_i \in C $$
(Proposition 3.8.2). Since $\tau\beta\tau Y \neq 0$, some $\eta_i$ is non-zero\,: $\eta_i \neq 0$. Applying $\alpha_i(-\pi/2) \in (E_7)_0$ of Lemma 4.8.1 on $\beta P$, we get 
$$
      \alpha_i(-\pi/2)\beta P
        = \pmatrix{1 - p_i & -2E_i & 0 & E_i \cr
                   2E_i & 1 - p_i & -E_i & 0 \cr
                   0 & E_i & 0 & 0 \cr
                   -E_i & 0 & 0 & 0 \cr}
          \pmatrix{\beta X \cr
                   \tau\beta\tau Y \cr
                   0 \cr
                   \eta} = \pmatrix{* \cr
                                    * \cr
                                    \eta_i \cr
                                    *}, \quad \eta_i \neq 0,$$
so that this case is reduced to the case (1).
\vspace{1mm}

(3) Case $P = (X, 0, 0, \eta)$, $X \neq 0$. Choose $\beta \in E_6$ so that $\beta X = \xi_1E_1 + \xi_2E_2 + \xi_3E_3$, $\xi_i \in C$ (Proposition 3.8.2). Since $\beta X \neq 0$, some $\xi_i$ is non-zero\,: $\xi_i \neq 0$. Then
$$
       \alpha_{i+1}(- \pi/2)\beta P = (*, \xi_iE_{i+2} + \xi_{i+2}E_i, 0, *), \quad \xi_i \neq 0, $$
so that this case is reduced to the case (2).
\vspace{1mm}

(4) Case $P = (0, 0, 0, \eta)$, $\eta \neq 0$. Then
$$
         \alpha_1(- \pi/2)P = (\eta E_1, 0, 0, 0), \quad \eta \neq 0,$$
so that this case is also reduced to the case (3).
\vspace{1mm}

\noindent Consequently, any element $P$ can be transformed to a diagonal form with $\xi \neq 0$. Furthermore, by applying some $\phi_1(\theta) \in U(1) \subset (E_7)_0$ of Theorem 4.7.3 on it, then $\xi$ becomes $\xi > 0$. Thus the proof of Proposition 4.8.2 is completed.
\vspace{3mm}

{\bf Remark.} In Proposition 4.8.2, the condition $P \in \gM^C$ does not need. That is, any element $P \in \gP^C$ can be transformed to a diagonal form by some $\alpha \in E_7$. (See Miyasaka, Yasukura and Yokota [23]).
\vspace{2mm}

We define a space $\gM_1$, called the compact Freudenthal manifold, by
$$
   \gM_1 = \{ P \in \gP^C \, | \, P \times P = 0, \langle P, P \rangle = 1 \}.$$

{\bf Theorem 4.8.3.} \qquad \qquad \qquad $E_7/E_6 \simeq \gM_1.$
\vspace{1mm}

\noindent {\it In particular, the group $E_7$ is connected}.
\vspace{2mm}

{\bf Proof.} For $\alpha \in E_7$ and $P \in \gM_1$, we have $\alpha P \in \gM_1$. Hence $E_7$ acts on $\gM_1$. We shall prove that the group $(E_7)_0$ acts transitively on $\gM_1$. To prove this, it is sufficient to show that any element $P \in \gM_1$ can be transformed to $(0, 0, 1, 0) \in \gM_1$ by some $\alpha \in (E_7)_0$. Now, $P \in \gM_1$ can be transformed to a diagonal form by $\alpha \in (E_7)_0$:
$$
  \alpha P=\Big(\dfrac{1}{\xi}\pmatrix{\eta_2\eta_3 & 0 & 0 \cr 
                                       0 & \eta_3\eta_1 & 0 \cr
                                       0 & 0 & \eta_1 \eta_2},
                \pmatrix{\eta_1 & 0 & 0 \cr
                         0 & \eta_2 & 0 \cr 
                         0 & 0 & \eta_3}, 
           \xi, \dfrac{1}{\xi^2}\eta_1\eta_2\eta_3 \Big), \quad \xi > 0 $$
(Proposition 4.8.2). From the condition $\langle \alpha P,\alpha P \rangle = \langle P, P \rangle = 1,$  we have
$$
   \dfrac{1}{\xi^2}(|\eta_2\eta_3|^2 + |\eta_3\eta_1|^2 + |\eta_1\eta_2|^2) +
(|\eta_1|^2 + |\eta_2|^2 + |\eta_3|^2) + \xi^2 + \dfrac{1}{\xi^4}|\eta_1\eta_2\eta_3|^2 = 1,$$ 
that is,
$$
\displaylines{\hfill
      \Big(1 + \dfrac{|\eta_1|^2}{\xi^2}\Big)\Big(1 + \frac{|\eta_2|^2}{\xi^2}\Big)\Big(1 + \dfrac{|\eta_3|^2}{\xi^2}\Big) = \dfrac{1}{\xi^2}.
\hfill\mbox{(i)}}$$
Choose $r_1, r_2, r_3 \in \R$, $0 \le r_i < \dfrac{\pi}{2}$, such that
$$
              \tan r_i = \dfrac{|\eta_i|}{\xi}, \quad i = 1, 2, 3,$$
then (i) becomes
$$
           \xi = \cos r_1\cos r_2\cos r_3.$$
By letting
$$
             a_i = \dfrac{\eta_i}{|\eta_i|}r_i, \quad i = 1, 2, 3 $$
(if $\eta_i = 0$, then $a_i$ means $0$), we have
$$
        r_i=|a_i|, \quad \eta_i = \dfrac{1}{|a_i|}\dfrac{\eta_i}{|\eta_i|}r_i
\dfrac{|\eta_i|}{\xi}\xi = \dfrac{a_i}{|a_i|}\tan r_i\cos r_1\cos r_2\cos r_3.$$Therefore, we see that $\alpha P$ is of the form
$$
     \pmatrix{
      \pmatrix{\begin{array}{l}
    \cos|a_1|a_2\dfrac{\sin|a_2|}{|a_2|}a_3\dfrac{\sin|a_3|}{|a_3|} \quad \qquad  0 \qquad \qquad \qquad \qquad 0 
\vspace{1mm}\cr 
   \qquad \qquad 0 \qquad \quad a_1\dfrac{\sin|a_1|}{|a_1|}\cos|a_2|a_3\dfrac{\sin|a_3|}{|a_3|} \qquad \qquad 0 
\cr 
   \qquad \qquad 0 \qquad \qquad \qquad \qquad \qquad 0 \qquad \quad a_1\dfrac{\sin|a_1|}{|a_1|}a_2\dfrac{\sin|a_2|}{|a_2|}\cos|a_3|
\end{array}}
\vspace{1mm}\cr
   \pmatrix{\begin{array}{l}
     a_1\dfrac{\sin|a_1|}{|a_1|}\cos|a_2|\cos|a_3| \quad \qquad  0 \qquad \qquad \qquad \qquad 0 
\cr 
     \qquad \qquad 0 \quad \qquad \cos|a_1|a_2\dfrac{\sin|a_2|}{|a_2|}\cos|a_3| \qquad \qquad 0 
\cr 
  \qquad \qquad  0 \qquad \qquad \qquad \qquad \quad 0 \qquad \quad \cos|a_1|\cos|a_2|a_3\dfrac{\sin|a_3|}{|a_3|} 
\end{array}} 
\vspace{1mm}\cr
         \cos|a_1|\cos|a_2|\cos|a_3| 
\vspace{2mm}\cr
        a_1\dfrac{\sin|a_1|}{|a_1|}a_2\dfrac{\sin|a_2|}{|a_2|}a_3\dfrac{\sin|a_3|}{|a_3|}} 
$$
$$
  \displaylines{= \alpha_3(a_3)\alpha_2(a_2)\alpha_1(a_1)(0, 0, 1, 0),\hfill}$$
hence we have
$$
      \alpha_1(a_1)^{-1}\alpha_2(a_2)^{-1}\alpha_3(a_3)^{-1}\alpha P = (0, 0, 1, 0).$$
This shows the transitivity of $(E_7)_0$. Since we have $\gM_1=(E_7)_0(0, 0, 1, 0)$, $\gM_1$ is connected. Now, the group $E_7$ acts transitively on $\gM_1$ and the isotropy subgroup of $E_7$ at $(0, 0, 1, 0) \in \gM_1$ is $E_6$ (Theorem 4.7.2). Therefore we have the homeomorphism $E_7/E_6 \simeq \gM_1$. Finally, the connectedness of $E_7$ follows from the connectedness of $\gM_1$ and $E_6$.
\vspace{4mm}

{\bf 4.9. Center $z(E_7)$ of $E_7$}
\vspace{3mm}

{\bf Theorem 4.9.1.} {\it The center $z(E_7)$ of the group $E_7$ is the cyclic 
group of order} $2$:
$$
                    z(E_7) = \{ 1, -1 \}.$$

{\bf Proof.} Let $\alpha \in z(E_7)$. From the commutativity with $\beta \in E_6 \subset E_7$, we have $\beta\alpha(0, 0, 1, 0) = \alpha\beta(0, 0, 1, 0) = \alpha(0, 0, 1, 0)$. If we denote $\alpha(0, 0, 1, 0) = (X, Y, \xi, \eta) \in \gP^C$, then from $(\beta X, \tau\beta\tau Y, \xi, \eta) = (X, Y, \xi, \eta)$, we have
$$
      \beta X = X, \;\; \tau\beta\tau Y = Y  \quad \mbox{for all} \quad \beta \in E_6.$$
Hence $X = Y = 0$. Therefore, $\alpha(0, 0, 1, 0)$ is of the form
$$
            \alpha(0, 0, 1, 0) = (0, 0, \xi, \eta).$$
From the condition $\alpha(0, 0, 1, 0) \in \gM^C$, we have $\xi\eta = 0$. Suppose $\xi = 0$, then $\alpha(0, 0, 1, 0) = (0, 0, 0, \eta), \eta \neq 0$. Also from the commutativity with $\varphi_1(\theta) \in U(1) \subset E_7$ (Theorem 4.7.3), we have
$$
\begin{array}{l}
     (0, 0, 0, \theta^{-3}\eta) = \varphi_1(\theta)(0, 0, 0, \eta) = \varphi_1(\theta)\alpha(0, 0, 1, 0)
\vspace{1mm}\\
  \qquad \quad = \alpha\varphi_1(\theta)(0, 0, 1, 0) = \alpha(0, 0, \theta^3, 0) = (0, 0, 0, \theta^3\eta), 
\end{array}$$
and so $\theta^{-3}\eta = \theta^3\eta$ for all $\theta$. But this is a contradiction. Hence $\xi \neq 0$, $\eta = 0$, that is, $\alpha(0, 0, 1, 0) = (0, 0, \xi, 0)$. Similarly we have $\alpha(0, 0, 0, 1) = (0, 0, 0, \zeta).$ Since $\{ \alpha(0, 0, 1, 0), \alpha(0, 0, 1, 0) \} = 1$, we have $\xi\zeta = 1$, and therefore 
$$
        \alpha(0, 0, 1, 0) = (0, 0, \xi, 0), \quad \alpha(0, 0, 0, 1) = (0, 0, 0, \xi^{-1}). $$ 
Moreover, from the commutativity with $\lambda \in E_7$,
$$
\begin{array}{l}
     (0, 0, 0, -\xi) = \lambda(0, 0, \xi, 0) = \lambda\alpha(0, 0, 1, 0) 
\vspace{1mm}\\
    \qquad = \alpha\lambda(0, 0, 1, 0) = \alpha(0, 0, 0, -1) = (0, 0, 0, - \xi^{-1}). 
\end{array} $$
Hence $\xi = \xi^{-1}$, that is, $\xi = \pm 1$. In the case $\xi =1$, we have $\alpha \in E_6$ (Theorem 4.7.2), so that $\alpha \in z(E_6) = \{1, \omega 1, 
\omega^2 1 \}$ (Theorem 3.9.1), that is, 
$$
         \alpha = \pmatrix{\omega'1 & 0 & 0 & 0 \cr
                           0 & {\omega'}^{-1}1 & 0 & 0 \cr
                           0 & 0 & 1 & 0 \cr
                           0 & 0 & 0 & 1}, \quad \omega' = 1, \omega \;\mbox{or} \; \omega^2. $$
Again from the commutativity with $\lambda$,
$$
\begin{array}{l}
  (0, \omega'X, 0, 0) = - \lambda(\omega'X, 0, 0, 0) = - \lambda\alpha(X, 0, 0, 0)
\vspace{1mm}\\
   \qquad \quad = - \alpha\lambda(X, 0, 0, 0) = \alpha(0, X, 0, 0) = (0, {\omega'}^{-1}X, 0, 0),
\end{array}$$
for all $X \in \gJ^C,$ which shows that $\omega' = {\omega'}^{-1}$, hence $\omega' = 1$. Therefore $\alpha = 1$. In the case $\xi = - 1$, we have $- \alpha \in z(E_6)$, so that by the similar argument as above we have $- \alpha = 1$. Thus we have $z(E_7) = \{ 1, -1 \}$.
\vspace{2mm}

According to a general theory of compact Lie groups, it is known that the center of the simply connected compact simple Lie group of type $E_7$ is the cyclic 
group of order $2$. Hence the group $E_7$ has to be simply connected. Thus we have the following theorem.
\vspace{3mm}

{\bf Theorem 4.9.2.} $E_7 = \{ \alpha \in \Iso_C(\gP^C) \, | \, \alpha(P \times Q)\alpha^{-1} = \alpha P \times \alpha Q, \langle \alpha P, \alpha Q \rangle = \langle P, Q \rangle \}$ {\it is a simply connected compact Lie group of type $E_7$}. 
\vspace{4mm}

{\bf 4.10. Involution $\iota$ and subgroup $(U(1) \times E_6)/\Z_3$ of $E_7$}
\vspace{3mm}

{\bf Definition.} We define a $C$-linear transformation $\iota$ of $\gP^C$ by
$$
     \iota(X, Y, \xi, \eta) = (-iX, iY, -i\xi, i\eta).$$
Then, $\iota = \varphi_1(i) \in U(1) \subset E_7$, $\iota^2 = - 1 \in z(E_7)$ (Theorem 4.9.1) and so $\iota^4 = 1$. 
\vspace{3mm}

{\bf Lemma 4.10.1.} $\iota$ {\it is conjugate to $\lambda$ in $E_7$.}
\vspace{2mm}

{\bf Proof.} For $\delta = \alpha_1\Big(\dfrac{i\pi}{4}\Big)\alpha_2\Big(\dfrac{i\pi}{4}\Big)\alpha_3\Big(\dfrac{i\pi}{4}\Big)$ (Lemma 4.8.1), we have $ \iota = \delta^{-1}\lambda\delta$.
\vspace{3mm}

$\iota$ induces an involutive automorphism $\wti{\iota} : E_7  \to E_7$ by
$$
      \wti{\iota}(\alpha) = \iota\alpha\iota^{-1}, \quad \alpha \in E_7. $$
We shall now study the following subgroup $(E_7)^{\iota}$ of $E_7$:
\begin{eqnarray*}
       (E_7)^{\iota} \!\!\! &=& \!\!\! \{ \alpha \in E_7 \, | \, \iota\alpha = \alpha\iota \}
\vspace{1mm}\\
        \!\!\! &\cong& \!\!\! \{ \alpha \in E_7 \, | \, \lambda\alpha = \alpha\lambda \} = (E_7)^{\lambda}.
\end{eqnarray*}

{\bf Theorem 4.10.2}  $\; (E_7)^{\iota} \cong (U(1) \times E_6)/\Z_3, \;\; \Z_3 = \{ (1, 1), (\omega, \omega 1),(\omega^2, \omega^2 1) \}$,

\noindent $\omega = - \dfrac{1}{2} + \dfrac{\sqrt{3}}{2}i \in C$.
\vspace{2mm}

{\bf Proof.}  We define a mapping $\varphi : U(1) \times E_6 \to (E_7)^{\iota}$ by
$$
       \varphi(\theta, \beta) = \varphi_1(\theta)\beta, \quad
    \varphi_1(\theta) = \pmatrix{\theta^{-1}1 & 0 & 0 & 0 \cr
                                    0 & \theta 1 & 0 & 0 \cr
                                    0 & 0 & \theta^3 & 0 \cr
                                    0 & 0 & 0 & \theta^{-3}},
           \beta = \pmatrix{\beta & 0 & 0 & 0 \cr
                            0 & \tau\beta\tau & 0 & 0 \cr
                            0 & 0 & 1 & 0 \cr
                            0 & 0 & 0 & 1}.$$
Evidently $\varphi(\theta, \beta) \in (E_7)^{\iota}$. Since $\varphi_1(\theta)$ and  $\beta$ are commutaive, $\varphi$ is a homomorphism. We shall prove that $\varphi$ is onto. Let $\alpha \in (E_7)^{\iota}$. From $\iota\alpha = \alpha\iota$, $\alpha$ is seen to be of the form
$$
      \alpha = \pmatrix{\beta & 0 & M & 0 \vspace{0.5mm} \cr
                        0 & \delta & 0 & N \vspace{0.5mm} \cr
                        a & 0 & \mu & 0 \vspace{0.5mm} \cr
                        0 & b & 0 & \nu}, \quad
\left.\begin{array}{l}
     \beta, \delta \in \Hom_C(\gJ^C),
\vspace{0.5mm}\\
     a, b \in \Hom_C(\gJ^C, C),
\vspace{0.5mm}\\
     M, N \in \gJ^C,
\vspace{0.5mm}\\
     \mu, \nu \in \C.
\end{array}\right.$$
The condition $\alpha(0, 0, 1, 0)$, $\alpha(0, 0, 0, 1) \in \gM^C$ implies that
$$
        \mu M = 0, \quad \nu N = 0. $$
We shall first show that $M = N = 0$. Suppose that $M \neq 0$, $\mu = 0$. Then, the condition $\{\alpha(0, 0, 1, 0), \alpha(0, 0, 0, 1) \} = \{(0, 0, 1, 0), (0, 0, 0, 1) \} = 1$ implies that
$$
\displaylines{\hfill
                 (M, N) = 1.
\hfill\mbox{(i)}}$$
Hence we have $N \neq 0$, $\nu = 0$. From 

$$
       \gM^C \ni \alpha\pmatrix{X 
\vspace{1mm}\cr
                               \dfrac{1}{\eta}X \times X 
\vspace{1mm}\cr
                               \dfrac{1}{\eta^2}\det X
\vspace{1mm}\cr                \eta}
       = \pmatrix{\beta X + \dfrac{1}{\eta^2}(\det X)M 
\vspace{1mm}\cr
                  \dfrac{1}{\eta}\delta(X \times X) + \eta N
\vspace{1mm}\cr
                  a(X) 
\vspace{1mm}\cr 
                  \dfrac{1}{\eta}b(X \times X)},$$
we have
$$
\begin{array}{c}
        \Big(\dfrac{1}{\eta}\delta(X \times X) + \eta N \Big) \times \Big(\dfrac{1}{\eta}\delta(X \times X) + \eta N \Big) = a(X)\Big(\beta X +\dfrac{1}{\eta^2}(\det X)M \Big),
\vspace{1mm}\\
    \Big(\beta X + \dfrac{1}{\eta^2}(\det X)M, \dfrac{1}{\eta}\delta(X \times X) + \eta N \Big) = 3a(X)\dfrac{1}{\eta}b(X \times X)
\end{array}$$
hold for all $0 \neq \eta \in C$. Comparing the coefficients of $\eta$, we have
$$\displaylines{\hfill
\left\{\begin{array}{l}
     2\delta(X \times X) \times N = a(X)\beta X
\vspace{1mm}\\
     \delta(X \times X) \times \delta (X \times X) = a(X)(\det\,X)M
\vspace{1mm}\\
     (\beta X, \delta(X \times X)) + \det\,X = 3a(X)b(X \times X)\;\; \mbox{(use (i))}.
\end{array}\right.
\hfill
\left.\begin{array}{r}
\mbox{(ii)}\vspace{1mm}\\
\mbox{(iii)}\vspace{1mm}\\
\mbox{(iv)} 
\end{array}\right.}$$
Therfore, using (i) $\sim$ (iv), we have
$$
\begin{array}{l}
     a(X)\det X = a(X)(\det\,X)(M, N) = (\delta(X \times X) \times \delta(X \times X), N)
\vspace{1mm}\\
   \quad = (\delta(X \times X), \delta(X \times X) \times N) = \dfrac{1}{2}a(X)(\delta(X \times X), \beta X)
\vspace{1mm}\\
   \quad = \dfrac{1}{2}a(X)(3a(X)b(X \times X) - \det\,X).
\end{array}$$
Hence we have $a(X)\det\,X = a(X)^2b(X \times X)$. Furthermore we have
$$\displaylines{\hfill
            \det\,X = a(X)b(X \times X).
\hfill\mbox{(v)}}$$
Indeed, from $\mu = 0$, we deduce that $a \neq 0$. Since $ a : \gJ^C \to C$ is a linear form, the set $\{X \in \gJ^C \, | \, a(X) \neq 0 \}$ is dense in $\gJ^C$ and the correspondence $\det\,X$ and $b(X \times X)$ is continuous with respect to $X$, (v) is also valid for $X \in \gJ^C$ such that $a(X) = 0$. Now, since $a \neq 0$ and $b\neq 0$, (v) contradicts the irreducibility of the determinant $\det\,X$ with respect to the variables of its components. Consequently we have shown that $M = 0$. Similarly we can prove that $N = 0$. Therefore $\alpha(0, 0, 1, 0) = (0, 0, \mu, 0)$, $\alpha(0, 0, 0, 1) = (0, 0, 0, \nu)$. From the condition $\{ \alpha\dot{1}, \alpha\d {1} \} = 1$, $\langle \alpha\dot{1}, \alpha\dot{1} \rangle = 1$, we deduce that$$ 
       \alpha(0, 0, 1, 0) = (0, 0, \mu, 0), \;\; \alpha(0, 0, 0, 1) = (0, 0, 0, \mu^{-1}), \quad \mu \in C, (\tau\mu)\mu = 1.$$
If we choose $\theta \in C$ such that $\theta^3 = \mu$ and let $\beta = \varphi_1(\theta)^{-1}\alpha$, then $\beta(0, 0, 1, 0) = (0, 0, 1, 0)$, $\beta(0, 0, 0, 1) = (0, 0, 0, 1)$. Hence, $\beta \in E_6$ (Theorem 4.7.2) and we have
$$
        \alpha = \varphi_1(\theta)\beta, \quad \theta \in U(1),\beta \in E_6.$$
This shows $\varphi$ is onto. That $\Ker\,\varphi = \{ (1, 1), (\omega, \omega 1),(\omega^2, \omega^2 1) \} = \Z_3$ is easily obtained. Thus we have the isomorphism $(U(1) \times E_6)/\Z_3 \cong (E_7)^{\iota}$.
\vspace{3mm}

{\bf Remark.} $(E_7)^{\iota}$ is connected as a fixed points subgroup under the involution $\iota$ of the simply connected Lie group $E_7$. Hence, to show that $\varphi : U(1) \times E_6 \to (E_6)^{\iota}$ is onto, it is sufficient to show that $\varphi_* : \gu(1) \oplus \gge_6 \to (\gge_7)^{\iota}$ is onto, 
\vspace{4mm}
which is easily shown.

{\bf 4.11. Involution $\sigma$ and subgroup $(SU(2)\times Spin(12))/\Z_2$ of $E_7$}
\vspace{3mm}

We define an involutive $C$-linear transformation $\sigma$ of $\gP^C$ by
$$
     \sigma(X, Y, \xi, \eta) = (\sigma X, \sigma Y, \xi, \eta).$$
which is the extension of the $C$-linear transformation $\sigma$ of $\gJ^C$. This is the same as $\sigma \in F_4$ regarding as $\sigma \in F_4 \subset E_6 \subset E_7.$ 
\vspace{1mm}

We shall now study the following subgroup $(E_7)^{\sigma}$ of $E_7$:
$$
 (E_7)^{\sigma} = \{ \alpha \in E_7  \, | \, \sigma\alpha = \alpha\sigma \}.$$  To this end, we define two $C$-linear mappings $\kappa, \mu : \gP^C \to \gP^C$ by
$$
    \kappa = {\mit\Phi}(-2E_1 \vee E_1, 0, 0, - 1), \quad
    \mu = {\mit\Phi}(0, E_1, E_1, 0). $$
The explicit forms of $\kappa$ and $\mu$ are given by
$$
\begin{array}{l}
      \kappa\pmatrix{X \cr Y \cr
                     \xi \cr \eta}
           = \pmatrix{- \kappa_1X \cr \kappa_1Y \cr
                      - \xi \cr \eta}, 
     \quad \kappa_1X = (E_1, X) - 4E_1 \times (E_1 \times X), 
\vspace{1mm}\\
       \mu\pmatrix{X \cr Y \cr
                   \xi \cr \eta} 
           = \pmatrix{2E_1 \times Y + \eta E_1 \cr
                      2E_1 \times X + \xi E_1 \cr
                      (E_1, Y) \cr
                      (E_1, X)}.
\end{array}$$
More precisely, $\kappa$ and $\mu$ are of the form
$$
\begin{array}{l}
     \kappa(X, Y, \xi, \eta) = 
     \kappa\Big(\pmatrix{\xi_1 & x_3 & \ov{x}_2 \cr
                         \ov{x}_3 & \xi_2 & x_1 \cr
                         x_2 & \ov{x}_1 & \xi_3},
                \pmatrix{\eta_1 & y_3 & \ov{y}_2 \cr
                         \ov{y}_3 & \eta_2 & y_1 \cr
                         y_2 & \ov{y}_1 & \eta_3}, \xi, \eta \Big)
\vspace{1mm}\\
     \qquad \qquad \quad \;\;= \Big(\pmatrix{- \xi_1 & 0 & 0 \cr
                       0 & \xi_2 & x_1 \cr
                       0 & \ov{x}_1 & \xi_3},
               \pmatrix{\eta_1 & 0 & 0 \cr
                        0 & - \eta_2 & - y_1 \cr
                        0 & - \ov{y}_1 & - \eta_3}, - \xi, \eta \Big),
\vspace{1mm}\\ 
             \mu(X, Y, \xi, \eta) = 
             \Big(\pmatrix{\eta & 0 & 0 \cr
                          0 & \eta_3 & -y_1 \cr
                          0 & - \ov{y}_1 & \eta_2}, 
                  \pmatrix{\xi & 0 & 0 \cr
                           0 & \xi_3 & -x_1 \cr
                           0 & - \ov{x}_1 & \xi_2}, \eta_1, \xi_1 \Big).

\end{array} $$

{\bf Lemma 4.11.1.} (1) \qquad
    $  \kappa\mu = \mu\kappa, \quad
\left\{\begin{array}{l} \kappa\sigma = \sigma\kappa \\
                        \mu\sigma = \sigma\mu, 
\end{array}\right. 
\quad                   
\left\{\begin{array}{l}
                  \kappa\lambda = - \lambda\kappa \\
                  \mu\lambda = - \lambda\mu.
\end{array}\right.$
\vspace{1mm}

(2) {\it If $\alpha \in E_7$ satisfies $\kappa\alpha = \alpha\kappa$, then $\alpha$ also satisfies} $\sigma\alpha = \alpha\sigma$.
\vspace{2mm}

{\bf Proof.} (1) These are checked by direct calculations.
\vspace{1mm}

(2) Since $\sigma = \exp \pi i\kappa$, we have $\sigma\alpha = (\exp \pi i\kappa)\alpha = \alpha(\exp \pi i\kappa) = \alpha\sigma.$
\vspace{3mm}

We shall first study the following subgroups $(E_7)^{\kappa,\mu}$ and $((E_7)^{\kappa,\mu})_{(0, E_1, 0, 1)}$ of $E_7$:
\begin{eqnarray*}
    (E_7)^{\kappa,\mu} \!\!\! &=& \!\!\! \{ \alpha \in E_7 \, | \, \kappa\alpha = \alpha\kappa, \mu\alpha = \alpha\mu \},
\vspace{1mm}\\
     ((E_7)^{\kappa,\mu})_{(0, E_1, 0, 1)} \!\!\! &=& \!\!\! \{\alpha \in (E_7)^{\kappa,\mu} \, | \, \alpha(0, E_1, 0, 1) = (0, E_1, 0, 1) \}.
\end{eqnarray*}

{\bf Proposition 4.11.2.} {\it The Lie algebras $(\gge_7)^{\sigma}$, $(\gge_7)^{\kappa,\mu}$, $((\gge_7)^{\kappa,\mu})_{(0, E_1, 0, 1)}$ of the groups $(E_7)^{\sigma}, (E_7)^{\kappa,\mu},$ $ ((E_7)^{\kappa,\mu})_{(0, E_1, 0, 1)}$ are respectively given by}
\vspace{2mm}

$\begin{array}{l}
(1)\;\; (\gge_7)^{\sigma} = \{{\mit\Phi} \in \gge_7 \, | \, \sigma{\mit\Phi} = {\mit\Phi}\sigma \} 
\vspace{1mm}\\
    \qquad \qquad \;= \{{\mit\Phi}(\phi, A, - \tau A, \nu) \in \gge_7 | \, \phi \in (\gge_6)^{\sigma}, A \in (\gJ^C)_{\sigma} \}.
\end{array}$
\vspace{1mm}\\  

$\begin{array}{l}
(2)\;\;  (\gge_7)^{\kappa,\mu} = \{{\mit\Phi} \in \gge_7 \, | \, \kappa{\mit\Phi} = {\mit\Phi}\kappa, \mu{\mit\Phi} = {\mit\Phi}\mu \}
\vspace{1mm}\\
    \qquad \qquad \quad = \Biggl\{ {\mit\Phi}(\phi. A, - \tau A, \nu) \in \gge_7 \left| \, 
\begin{array}{l}
        \phi \in (\gge_6)^{\sigma}, A \in (\gJ^C)_{\sigma}, (E_1, A) = 0, 
\vspace{1mm}\\
        \nu = - \dfrac{3}{2}(\phi E_1, E_1) 
\end{array}\right.        
\Biggl\}.
\end{array}$
\vspace{1mm}\\

$\begin{array}{l}
(3)\;\;  ((\gge_7)^{\kappa,\mu})_{(0, E_1, 0, 1)} = \{ {\mit\Phi} \in (\gge_7)^{\kappa,\mu} \, |\, {\mit\Phi}((0, E_1, 0, 1)) = 0 \}
\vspace{1mm}\\
   \qquad \qquad \quad = \biggl\{ {\mit\Phi}(\phi, A, - \tau A, 0) \in \gge_7 \left| \,\begin{array}{l}
        \phi \in \gge_6, \phi E_1 = 0, 
\vspace{1mm}\\
        A \in \gJ^C, 2E_1 \times A = \tau A 
\end{array} \right.
\biggl\}. 
\end{array}$
\vspace{2mm}

{\bf Proof.} (1) It is not difficult to prove and so is omitted here.
\vspace{1mm}

(2) Suppose that ${\mit\Phi} = {\mit\Phi}(\phi, A, - \tau A, \nu) \in \gge_7$ satisfies $\kappa{\mit\Phi} = {\mit\Phi}\kappa$ and $\mu{\mit\Phi} = {\mit\Phi}\mu.$ Comparing the $\xi$-term of $\kappa{\mit\Phi}P = {\mit\Phi}\kappa P, P = (X, Y, \xi, \eta) \in \gP^C$, we have
$$
           (A, Y) = - (A, \kappa_1Y), \quad Y \in \gJ^C. $$
Let $Y = E_1$, then we have $(A, E_1) = 0$. Next, comparing the $\eta$-term of $\mu{\mit\Phi} = {\mit\Phi}\mu$, we have
$$
\displaylines{\hfill
       (E_1, \phi X) - \dfrac{1}{3}\nu(E_1, X) = - \nu(E_1, X).
\hfill\mbox{(i)}} $$
Since $\phi \in (\gge_6)^{\sigma}$, we can set $\phi E_1 = kE_1, k \in i\R$ (Lemma 3.10.1). Hence let $X = E_1$ in (i), then we have $ k = - \dfrac{2}{3}\nu$. Conversely, if ${\mit\Phi} = {\mit\Phi}(\phi, A, - \tau A, \nu) \in \gge_7$ 
has the 
\vspace{0.7mm}
condition above, then from the following Lemma 4.11.3, we can see that ${\mit\Phi}$ satisfies $\kappa{\mit\Phi} = {\mit\Phi}\kappa$ and $\mu{\mit\Phi} = {\mit\Phi}\mu.$ 
\vspace{3mm}

{\bf Lemma 4.11.3.} {\it In $\gJ^C$, the following hold.}
\vspace{1mm}

(1) {\it For $A \in (\gJ^C)_{\sigma}$, we have} $\kappa_1(A \times X) = \kappa_1A  \times \kappa_1X, \;\; X\in \gJ^C.$
\vspace{1mm}

(2) {\it For $\phi \in (\gge_6)^{\sigma}$, we have} $\kappa_1\phi = \phi\kappa_1.$
\vspace{1mm}

(3) {\it For $A \in (\gJ^C)_{\sigma}, (E_1, A) = 0$, we have $\kappa_1A = - A$ and }
$$
    - 4\tau A \times (E_1 \times X) + (E_1, X)A = 4E_1 \times (A \times X) - \langle A, X \rangle E_1, \quad X \in \gJ^C. $$

We shall now return to the proof of (3) of Proposition 4.11.2. 
\vspace{1mm}

\noindent If ${\mit\Phi} = {\mit\Phi}(\phi, A, - \tau A, \nu) \in (\gge_7)^{\kappa, \mu}$ satisfies ${\mit\Phi}((0, E_1, 0, 1)) = 0$, then
$$
          \nu = - (A, E_1) = - \tau(A, E_1), $$
so that $\nu = \tau\nu$. Together with $\tau\nu = - \nu$, we have $\nu = 0$ and $\phi E_1 = 0$, furthermore, we have $2A \times E_1 = \tau A$ (in this case, $A \in (\gJ^C)_{\sigma}$ and $(E_1, A) = 0$ automatically 
\vspace{2mm}
hold).

For $\nu \in i\R$, we define a $C$-linear mapping $\phi(\nu) : \gJ^C \to \gJ^C$ by
$$
      \phi(\nu) = 2\nu E_1 \vee E_1, $$
that is,
$$
     \phi(\nu)\pmatrix{\xi_1 & x_3 & \ov{x}_2 \cr
                       \ov{x}_3 & \xi_2 & x_1 \cr
                       x_2 & \ov{x}_1 & \xi_3}
         = \dfrac{\nu}{3}\pmatrix{4\xi_1 & x_3 & \ov{x}_2 \cr
                                  \ov{x}_3 & -2 \xi_2 & -2 x_1 \cr
                                  x_2 & -2 \ov{x}_1 & -2\xi_3} $$
(Lemma 3.4.2.(2)). Then $\phi(\nu) \in (\gge_6)^{\sigma}$.
\vspace{3mm}

{\bf Proposition 4.11.4.} (1) \quad $\ga_1 =\{ {\mit\Phi}(\phi (\nu), aE_1,-\tau a E_1, \nu) \, | \, a \in C, \nu \in i\R \}$ 

\noindent {\it is a Lie subalgebra of $(\gge_7)^{\sigma}$ and isomorphic to the Lie algebra $\su(2)$}.
\vspace{1mm}

(2) {\it The Lie algebra $(\gge_7)^{\sigma}$ is isomorphic to the direct sum of Lie algebras $\ga_1$ and $(\gge_7)^{\kappa,\mu}$}:
$$
      (\gge_7)^{\sigma} \cong \ga_1 \oplus (\gge_7)^{\kappa,\mu}.$$

{\bf Proof.} (1) The mapping $\varphi_* : \ga_1 \to \su(2) = \{ D \in M(2, C) \, | \, \tau(^t\!D) = - D \}$ defined by
$$
       \varphi_*({\mit\Phi}(\phi(\nu), aE_1, - \tau aE_1, \nu)) = 
            \pmatrix{\nu & a \cr
                     - \tau a & -\nu} $$
gives an isomorphism as Lie algebras. Indeed, this is clear from
$$
     \Big[\pmatrix{\nu & a \cr
                   - \tau a & -\nu},
          \pmatrix{\rho & b \cr
                   - \tau b & -\rho}\Big]
        = \pmatrix{b(\tau a) - a(\tau b) & 2(b\nu - a\rho) \cr
                   -2\tau(b\nu - a\rho) & a(\tau b) - b(\tau a)},$$
and
$$
\begin{array}{l}
    \Big[{\mit\Phi}(\phi(\nu), aE_1, - \tau aE_1, \nu), {\mit\Phi}(\phi(\rho), bE_1, - \tau bE_1, \rho) \Big]
\vspace{1mm}\\
    = {\mit\Phi}(\phi(b\tau a) - a(\tau b)), 2(b\nu - a\rho)E_1, - 2\tau(b\nu - a\rho)E_1, (\tau a)b - a(\tau b)).
\end{array} $$

(2) Using (1) above and Proposition 4.11.2.(2), the following decomposition of $(\gge_7)^\sigma$,
$$
    (\gge_7)^{\sigma} \ni {\mit\Phi}\pmatrix{\phi \cr A \cr
                                             - \tau A \cr \nu}
        = {\mit\Phi}\pmatrix{\phi(\nu') \cr aE_1 \cr
                             - \tau aE_1 \cr \nu'}
        + {\mit\Phi}\pmatrix{\phi - \phi(\nu') \cr
                             A - aE_1 \cr
                             - \tau A + \tau aE_1 \cr
                             \nu - \nu'}
        \in \ga_1 \oplus (\gge_7)^{\kappa,\mu}, $$
where $\nu' = \dfrac{1}{3}\nu + \dfrac{1}{2}(E_1, \phi E_1)$, $a = (E_1, A)$, gives an isomorphism as Lie algebras.
\vspace{3mm}

{\bf Lemma 4.11.5.} {\it For $a \in C$, we have }
$$
        \alpha_{23}(a) = \alpha_{2}(a)\alpha_{3}(\tau a) \in ((E_7)^{\kappa,\mu})_{(0, E_1, 0, 1)},$$
{\it where $\alpha_i(a) \in E_7, i = 2, 3$ are defined in Lemma} 4.8.1.
\vspace{2mm}

{\bf Proof.} Since ${\mit\Phi}(0, -\tau aE_i, aE_i, 0) \in (\gge_7)^{\kappa,\mu}$, we have $\alpha_i(a) = \exp{\it\Phi}(0, -\tau aE_i, $ $aE_i, 0) \in (E_7)^{\kappa,\mu}$, $i = 2, 3$. Since $\alpha_2(a)$ and $\alpha_3(\tau a)$ are commutative, we have
\vspace{2mm}

\qquad \quad \qquad \qquad
    $\alpha_{23}(a) = \alpha_2(a)\alpha_3(\tau a)$
\vspace{1mm}

\qquad \quad \qquad \qquad \qquad \quad
     $= \exp{\mit\Phi}(0,-\tau aE_2 - aE_3, aE_2 + \tau aE_3, 0)$.
\vspace{2mm}

\noindent Since ${\mit\Phi}(0, - \tau aE_2 - aE_3, aE_2 + \tau aE_3, 0) \in ((\gge_7)^{\kappa,\mu})_{(0, E_1, 0, 1)}$, we have $\alpha_{23}(a) \in \\
\vspace{2mm}
((E_7)^{\kappa,\nu})_{(0, E_1, 0, 1)}$.

We recall the group
\begin{eqnarray*}
     Spin(10) \!\!\! &=& \!\!\! \{ \alpha \in E_6 \, | \, \alpha E_1 = E_1 \}
\vspace{1mm}\\
              \!\!\! &=& \!\!\! \{ \alpha \in E_6 \, | \, \sigma\alpha = \alpha\sigma, \alpha E_1 = E_1 \} \subset E_7
\end{eqnarray*}
(Lemma 3.10.4). The group $Spin(10)$ acts transitively on the 9 dimensional sphere
$$
  S^9 = \Big\{\Big(\pmatrix{0 & 0 & 0 \cr
                                0 & \xi & x \cr
                                0 & \ov{x} & - \tau\xi}, 0, 0, 0 \Big) \, \Big| \, \xi \in C, x \in \gC, \ov{x}x + (\tau\xi)\xi = 1 \Big\}. $$     

{\bf Lemma 4.11.6.} {\it For $\alpha \in ((E_7)^{\kappa,\mu})_{(0, E_1, 0, 1)}$, we have
$$
   \alpha(0, -E_1, 0, 1) = (0, -E_1, 0, 1) \quad \mbox{\it if and only if} \quad     \alpha(0, 0, 1, 0) = (0, 0, 1, 0). $$
In particular, we have}
$$
        \{ \alpha \in ((E_7)^{\kappa,\mu})_{(0, E_1, 0, 1)} \, | \, \alpha(0, -E_1, 0, 1)=(0, -E_1, 0, 1) \} \cong Spin(10).$$

{\bf Proof.} If $\alpha \in (E_7)^{\kappa,\mu}$ satisfies $\alpha(0, E_1, 0, 1) = (0, E_1, 0, 1)$ and $\alpha (0, -E_1, 0, 1) = (0, -E_1, 0, 1)$, then we have $\alpha(0, 0, 0, 1) = (0, 0, 0, 1)$ and $\alpha(0, E_1, 0, 0) = (0, E_1, 0, 0),$ which imply that 
$\alpha(0, 0, 1, 0) = \alpha\mu(0, E_1, 0, 0) = \mu\alpha(0, E_1, 0, 0)
= \mu(0, E_1, 0, 0) = (0, 0, 1, 0)$.
The converse can be  similarly proved. If $\alpha \in E_7$ satisfies $\alpha(0, 0, 1, 0) = (0, 0, 1, 0)$, then $\alpha \in E_6$ (Theorem 4.7.2), and from the condition $\alpha E_1 = E_1$, we obtain $\alpha \in Spin(10)$, (Theorem 3.10.4). The converse also holds.
\vspace{2mm}

We define an 11 dimensional $\R$-vector space $V^{11}$ by
\begin{eqnarray*}
     V^{11} \!\!\!&=&\!\!\! \{ P \in \gP^C \, | \, \kappa P = P, \mu\tau\lambda P = P, P \times (0, E_1, 0, 1, 0) = 0 \}
\vspace{1mm}\\
     \!\!\!&=&\!\!\! \Big\{\Big(\pmatrix{0 & 0 & 0 \cr
                                             0 & \xi & x \cr
                                             0 & \ov{x} & - \tau\xi},
           \pmatrix{\eta & 0 & 0 \cr
                    0 & 0 & 0 \cr
                    0 & 0 & 0}, 0, \tau\eta \Big) \,
\Big| \, x \in \gC, \xi \in C, \eta \in i\R \Big\}
\end{eqnarray*}
with the norm $(P, P)_{\mu}$ given by
$$
     (P, P)_{\mu} = \dfrac{1}{2}\{\mu P, P\} = \ov{x}x + (\tau\xi)\xi + (\tau\eta)\eta. $$

{\bf Proposition 4.11.7.} \qquad \quad $((E_7)^{\kappa,\mu})_{(0, E_1, 0, 1)}/Spin(10) \simeq S^{10}.$
\vspace{1mm}

\noindent {\it In particular, the group $((E_7)^{\kappa,\mu})_{(0, E_1, 0, 1)}$ is connected.}
\vspace{2mm}

{\bf Proof.}  $S^{10} = \{ P \in V^{11} \, | \, (P, P)_{\mu} = 1 \}$ is a 10 dimensional sphere. For $\alpha \in ((E_7)^{\kappa,\mu})_{(0, E_1, 0, 1)}$ and $P \in S^{10}$, we have $\alpha P \in S^{10}$ (Proposition 4.2.2, Lemma 4.3.3). Hence the group $((E_7)^{\kappa,\mu})_{(0, E_1, 0, 1)}$ acts on $S^{10}$. We shall prove that this action is transitive. To prove this, it is sufficient to show that any element $P \in S^{10}$ can be transformed to $(0, - iE_1, 0, i) \in S^{10}$ by some $\alpha \in ((E_7)^{\kappa,\mu})_{(0, E_1, 0, 1)}$. Now, for a given
$$
    P = \Big(\pmatrix{0 & 0 & 0 \cr
                        0 & \xi & x \cr
                        0 & \ov{x} & -\tau\xi},
               \pmatrix{\eta & 0 & 0 \cr
                        0 & 0 & 0 \cr
                        0 & 0 & 0}, 0, \tau\eta \Big) \in S^{10},$$
choose $a \in \R$, $0 \le a < \dfrac{\pi}{4}$, such that
$$
              \tan 2a = \dfrac{2\eta }{\tau\xi - \xi}.$$
\Big(If $\tau\xi - \xi = 0$, then we choose $a = \dfrac{\pi}{4}$\Big). Applying $\alpha_{23}(a)$ of Lemma 4.11.5 on $P$, then the $\eta$-part of $\alpha_{23}(a)P$ becomes
$$
\begin{array}{l}
     2\sin^2 a (E_2, E_3 \times X) + \tau\xi\sin a\cos a - (E_3, Y)\sin a\cos a - \eta\cos^2 a
\vspace{1mm}\\
   \quad \qquad = \eta\sin^2 a + (\tau\xi - \xi)\sin a\cos a - \eta\cos^2a 
\vspace{1mm}\\
   \quad \qquad = \dfrac{1}{2}(\tau\xi - \xi)\sin 2a - \eta\cos 2a = 0.
\end{array}$$
Hence we have
$$
      \alpha_{23}(a)P \in S^9. $$
Since the group $Spin(10)$ acts transitively on $S^9$ (Proposition 3.10.3), there exists $\beta \in Spin(10) = (E_6)_{E_1} \subset ((E_7)^{\kappa,\mu})_{(0, E_1, 0, 1)}$ such that
$$
     \beta\alpha_{23}(a)P = (i(E_2 + E_3), 0, 0, 0).$$
Again, applying $\alpha_{23}(- \pi/4) \in ((E_7)^{\kappa,\mu})_{(0, E_1, 0, 1)}$ of Lemma 4.11.5 on the above, we have
$$
      \alpha_{23}(- \pi/4)\beta\alpha_{23}(a)P = (0, - iE_1, 0, i).$$
This shows the transitivity of $((E_7)^{\kappa,\mu})_{(0, E_1, 0, 1)}$. The isotropy subgroup of \\ 
$((E_7)^{\kappa,\mu})_{(0, E_1, 0, 1)}$ at $(0 -iE_1, 0, i)$ is $Spin(10)$ (Lemma 4.11.6). Thus we have the homeomorphism $((E_7)^{\kappa,\mu})_{(0, E_1, 0, 1)}/Spin(10) \simeq S^{10}$.
\vspace{3mm}

{\bf Theorem 4.11.8.} \qquad \qquad $((E_7)^{\kappa,\mu})_{(0, E_1, 0, 1)} \cong Spin(11).$
\vspace{1mm}

\noindent (From now on, we identify these groups).
\vspace{2mm}

{\bf Proof.} Analogous to Theorem 3.10.4, we can define a homomorphism $p : \\
((E_7)^{\kappa,\mu})_{(0, E_1, 0, 1)} \to SO(11) = SO(V^{11})$  by $p(\alpha) = \alpha | V^{11}$. The restriction $p'$ of $p$ to $(E_6)_{E_1}$ coincides with the homomorphism $p' : Spin(10) \to SO(10) = SO(V^{10})$ (where $V^{10} = \{ P \in V^{11} \, | \, P = (X, 0, 0, 0) \}$). In particular, $p' : Spin(10) \to SO(10)$ is onto. Hence, from the following commutative diagram
\begin{center}
\begin{tabular}{ccccccccc}
$1$ & $\longrightarrow$ & $Spin(10)$ & $\longrightarrow$ & $((E_7)^{\kappa,\mu})_{(0, E_1, 0, 1)}$ & $\longrightarrow$ & $S^{10}$ & $\longrightarrow$ & $*$\vspace{0.7mm}\\
$$ & $$ & $\downarrow p'$ & $$ & $\downarrow p$ & $$ & $\downarrow =$ & $$ & $$\vspace{0.7mm}\\
$1$ & $\longrightarrow$ & $SO(10)$ & $\longrightarrow$ & $SO(11)$ & $\longrightarrow$ & $S^{10}$ & $\longrightarrow$ & $*$\\
\end{tabular}
\end{center}
we see that $p : ((E_7)^{\kappa,\mu})_{(0, E_1, 0, 1)} \to SO(11)$ is onto by the five lemma. Using the five lemme again, we see that $\Ker \,p$ coincides with $\Ker \,p'$. Hence $\Ker \,p = \{1, \sigma\}$ (Theorem 3.10.4). Thus we have the isomorphism 
$$
       ((E_7)^{\kappa, \mu})_{(0, E_1, 0, 1)}/\{1, \sigma \} \cong SO(11).$$
Therefore the group $((E_7)^{\kappa,\mu})_{(0, E_1, 0, 1)}$ is isomorphic to the group $Spin(11)$ as the universal covering group of $SO(11)$.
\vspace{3mm}

{\bf Lemma 4.11.9.} {\it For $t \in \R$, we define a mapping $\alpha(t) : \gP^C \to \gP^C$ by }
$$
\begin{array}{l}
     \alpha(t)\Big(\pmatrix{\xi_1 & x_3 & \ov{x}_2 \cr
                              \ov{x}_3 & \xi_2 & x_1 \cr
                              x_2 & \ov{x}_1 & \xi_3},
                     \pmatrix{\eta_1 & y_3 & \ov{y}_2 \cr
                              \ov{y}_3 & \eta_2 & y_1 \cr
                              y_2 & \ov{y}_1 & \eta_3}, \xi, \eta \Big)
\vspace{1mm}\\
      = \Big(\pmatrix{e^{2it}\xi_1 & e^{it}x_3 & e^{it}\ov{x}_2 \cr
                 e^{it}\ov{x}_3 & \xi_2 & x_1 \cr
                 e^{it}x_2 & \ov{x}_1 & \xi_3},
        \pmatrix{e^{-2it}\eta_1 & e^{-it}y_3 & e^{-it}\ov{y}_2 \cr
                 e^{-it}\ov{y}_3 & \eta_2 & y_1 \cr
                 e^{-it}y_2 & \ov{y}_1 & \eta_3},
                 e^{-2it}\xi, e^{2it}\eta \Big),
\end{array}$$
{\it then} $\alpha(t) \in (E_7)^{\kappa,\mu}$.
\vspace{2mm}

{\bf Proof.}  For $\nu = it \in i\R$, let $\phi(\nu) = 2\nu E_1 \vee E_1 \in (\gge_6)^{\sigma}$. Then, ${\mit\Phi}(\phi(\nu), 0, 0, -2\nu) \in (\gge_7)^{\kappa,\nu}$ (Proposition 4.11.2) and $\alpha(t) = \exp{\mit\Phi}(\phi(\nu), 0, 0, -2 \nu)$. Hence we have $\alpha(t) \in (E_7)^{\kappa,\mu}$.
\vspace{2mm}

We define a 12 dimensional $\R$-vector space $V^{12}$ by
\begin{eqnarray*}
       V^{12} \!\!\!&=& \!\!\!\{ P \in \gP^C  \, | \, \kappa P = P, \mu\tau\lambda P = P  \}
\vspace{1mm}\\
      \!\!\!&=& \!\!\! \Big\{\Big(
            \pmatrix{0 & 0 & 0 \cr
                     0 & \xi & x \cr
                     0 & \ov{x} & -\tau \xi},
            \pmatrix{\eta & 0 & 0 \cr
                     0 & 0 & 0 \cr
                     0 & 0 & 0}, 0, \tau\eta \Big) \, \Big|
            \, x \in \gC, \xi, \eta \in C \Big\}
\end{eqnarray*}
with the norm $(P, P)_{\mu}$ given by
$$
      (P, P)_{\mu} = \dfrac{1}{2}\{ \mu P, P \} = \ov{x}x + (\tau\xi)\xi + (\tau\eta)\eta. $$

{\bf Proposition 4.11.10.} \qquad \qquad $(E_7)^{\kappa,\mu}/Spin(11) \simeq S^{11}.$ 
\vspace{1mm}

\noindent {\it In particular, the group $(E_7)^{\kappa,\mu}$ is connected.}
\vspace{2mm}

{\bf Proof.} $S^{11} = \{ P \in V^{12} \, | \, (P, P)_{\mu} = 1 \}$ is an 11 dimensional sphere. For $\alpha \in (E_7)^{\kappa,\mu}$ and $P \in S^{11}$, we have $\alpha P \in S^{11}$ (Proposition 4.2.2, Lemma 4.3.3). Hence the group $(E_7)^{\kappa,\mu}$ acts on $S^{11}$. We shall prove that this action is transitive. To prove this, it is sufficient to show that any element $P \in S^{11}$ can be transformed to $(0, E_1, 0, 1) \in S^{11}$ by some $\alpha \in (E_7)^{\kappa,\mu}$. Now, for a given
$$
      P = \Big(\pmatrix{0 & 0 & 0 \cr
                          0 & \xi & x \cr
                          0 & \ov{x} & -\tau\xi},
                 \pmatrix{\eta & 0 & 0 \cr
                          0 & 0 & 0 \cr
                          0 & 0 & 0}, 0, \tau\eta \Big) \in S^{11},$$
we choose $t \in \R$ such that $e^{-2it}\eta \in i\R$. Applying $\alpha(t)$ of Lemma 4.11.9 on $P$, we get
$$
         \alpha(t)P \in S^{10}.$$
Since the group $Spin(11)$ acts transitively on $S^{10}$ (Proposition 4.11.7), there exists $\beta \in Spin(11) = ((E_7)^{\kappa,\mu})_{(0, E_1, 0, 1)}$ such that
$$
        \beta\alpha(t)P = (0, - iE_1, 0, i).$$
If we further apply $\alpha(-\pi/4) \in (E_7)^{\kappa,\mu}$ of Lemma 4.11.9 on the above, then we have
$$ 
    \alpha(-\pi/4)\beta\alpha(t)P = (0, E_1, 0, 1).$$
This shows the transitivity of $(E_7)^{\kappa,\mu}$. The isotropy subgroup of $(E_7)^{\kappa,\mu}$ at $(0, E_1, 0, 1)$ is $Spin(11)$ (Theorem 4.11.8). Thus we have the homeomorphism $(E_7)^{\kappa,\mu}/$ $Spin(11)$ $ \simeq S^{11}$.
\vspace{3mm}

{\bf Theorem 4.11.11.} \qquad \qquad $(E_7)^{\kappa,\mu} \cong Spin(12).$
\vspace{1mm}

\noindent (From now on, we identify these groups).
\vspace{1mm}

{\bf Proof.} Analogous to Theorem 4.11.8, we can define a homomorphism 
$$
     p : (E_7)^{\kappa,\mu} \to SO(12) = SO(V^{12}) $$
by $p(\alpha) = \alpha|V^{12}$. The restriction $p'$ of $p$ to $(E_6)^{\kappa,\mu}$  coincides with the homomorphism $p' : Spin(11) \to SO(11)$ of Theorem 4.11.8. In particular, $p' : Spin(11) \to SO(11)$ is onto. Hence from the following commutative diagram
\begin{center}
\begin{tabular}{ccccccccc}
$1$ & $\longrightarrow$ & $Spin(11)$ & $\longrightarrow$ & $(E_7)^{\kappa,\mu}$ & $\longrightarrow$ & $S^{11}$ & $\longrightarrow$ & $*$\vspace{0.7mm}\\
$$ & $$ & $\downarrow p'$ & $$ & $\downarrow p$ & $$ & $\downarrow =$ & $$ & $$\vspace{0.7mm}\\
$1$ & $\longrightarrow$ & $SO(11)$ & $\longrightarrow$ & $SO(12)$ & $\longrightarrow$ & $S^{11}$ & $\longrightarrow$ & $*$\\
\end{tabular}
\end{center}
we see that $p : (E_7)^{\kappa,\mu} \to SO(12)$ is onto by the five lemma. Using the five lemme again we see that $\Ker \,p$ coincides with $\Ker \,p'$. Hence $\Ker \,p = \{1, \sigma\}$ (Theorem 4.11.8). Thus we have the isomorphism
$$
          (E_7)^{\kappa,\mu}/\{1,\sigma \} \cong SO(12).$$
Therefore the group $(E_7)^{\kappa,\mu}$ is isomorphic to the group $Spin(12)$ as the universal covering group of $SO(12)$.
\vspace{3mm}

{\bf Theorem 4.11.12.} {\it The center $z(Spin(12))$ of $Spin(12)$ is}
$$
      z(Spin(12)) = \{ 1, -1, \sigma -\sigma \} \cong \{1, -1 \} \times \{1, \sigma \} \cong \Z_2 \times \Z_2. $$
{\it And we have}
$$
\begin{array}{l}
      Spin(12)/\{1,\sigma \} \cong SO(12), 
\vspace{1mm}\\
      Spin(12)/\{1, -1\} \cong Spin(12)/\{1,-\sigma \} \cong Ss(12).
\end{array} $$      

{\bf Theorem 4.11.13.} {\it The group $(E_7)^{\sigma}$ contains a subgroup
$$
      \varphi_2(SU(2)) = \{ \varphi_2(A) \in E_7  \, | \, A \in SU(2) \}$$
which is isomorphic to the special unitary group $SU(2) = \{ A \in M(2,C) \, | \, {}^t(\tau A)A = E, \det\,A = 1 \}$. Here, for $A \in SU(2)$, a mapping $\varphi_2(A) : \gP^C \to \gP^C$ is defined by}
$$
\begin{array}{l}
      \varphi_2(A)\Big(\pmatrix{\xi_1 & x_3 & \ov{x}_2 \cr
                             \ov{x}_3 & \xi_2 & x_1 \cr 
                             x_2 & \ov{x}_1 & \xi_3},
                    \pmatrix{\eta_1 & y_3 & \ov{y}_2 \cr
                             \ov{y}_3 & \eta_2 & y_1 \cr 
                             y_2 & \ov{y}_1 & \eta_3}, \xi, \eta \Big)
\vspace{1mm}\\
\qquad
      = \Big(\pmatrix{{\xi_1}' & {x_3}' & {\ov{x}_2}' \cr
                        {\ov{x}_3}' & {\xi_2}' & {x_1}' \cr
                         {x_2}' & {\ov{x}_1}' & {\xi_3}'},
               \pmatrix{{\eta_1}' & {y_3}' & {\ov{y}_2}' \cr
                        {\ov{y}_3}' & {\eta_2}' & {y_1}' \cr
                         {y_2}' & {\ov{y}_1}' & {\eta_3}'},
                          \xi', \eta'\Big),$$
\end{array}$$
{\it where}
$$
\begin{array}{c}
     \pmatrix{{\xi_1}' \cr \eta'} = A\pmatrix{\xi_1 \cr \eta}, \;\;
     \pmatrix{\xi' \cr {\eta_1}'} = A\pmatrix{\xi \cr \eta_1}, \;\;
     \pmatrix{{\eta_2}'\cr {\xi_3}'} = A\pmatrix{\eta_2 \cr \xi_3},\;\;
     \pmatrix{{\eta_3}' \cr {\xi_2}'} = A\pmatrix{\eta_3 \cr \xi_2},
\vspace{1mm}\\
     \pmatrix{{x_1}' \cr {y_1}'} = \tau A\pmatrix{x_1 \cr y_1}, \quad
     \pmatrix{{x_2}' \cr {y_2}'} = \pmatrix{x_2 \cr y_2}, \quad
     \pmatrix{{x_3}' \cr {y_3}'} = \pmatrix{x_3 \cr y_3}.
\end{array}$$

{\bf Proof.} The action of ${\mit\Phi}(\phi(\nu), aE_1, - \tau aE_1, \nu) \in \ga_1$ ($\phi(\nu) = 2\nu E_1 \vee E_1, \nu \in i\R, a \in C$) on $\gP^C$ is
$$
     {\mit\Phi}(\phi(\nu), aE_1, -\tau aE_1, \nu)(X, Y, \xi, \eta) =
(X', Y', \xi', \eta') $$
where 
$$
\begin{array}{ll}
     \pmatrix{{\xi_1} '\cr \eta'} = \pmatrix{\nu & a \cr 
                                             - \tau a & -\nu}
     \pmatrix{\xi_1 \cr \eta}, &
     \pmatrix{\xi' \cr {\eta_1}'} = \pmatrix{\nu & a \cr 
                                             - \tau a & -\nu}
     \pmatrix{\xi \cr \eta_1},
\vspace{1mm}\\
     \pmatrix{{\eta_2}' \cr {\xi_3}'} = \pmatrix{\nu & a \cr 
                                             - \tau a & -\nu}
     \pmatrix{\eta_2 \cr \xi_3}, &
     \pmatrix{{\eta_3}' \cr {\xi_2}'} = \pmatrix{\nu & a \cr
                                              - \tau a & - \nu}
     \pmatrix{\eta_3 \cr \xi_2},
\vspace{2mm}\\
     \pmatrix{{x_1}' \cr {y_1}'} = \pmatrix{- \nu & \tau a \cr
                                            -a & \nu}
     \pmatrix{x_1 \cr y_1}, & 
     \pmatrix{{x_2}' \cr {y_2}'} = \pmatrix{{x_3}' \cr {y_3}'} = \pmatrix{0 \cr 0}.
\end{array}$$
Therefore, for $A = \exp\pmatrix{\nu & a \cr
                                -\tau a & -\nu} \in SU(2)$, we have 
$$
     \exp({\mit\Phi}(\phi(\nu), aE_1, -\tau aE_1, \nu)) = \varphi_2(A) \in \varphi_2(SU(2)) \subset (E_7)^\sigma. $$

{\bf Lemma 4.11.14.} {\it The group $(E_7)^{\sigma}$ is connected.}
\vspace{2mm}

{\bf Proof.} The group $(E_7)^{\sigma}$ is the fixed points set obtained by the involutive automorphism $\sigma$ of the simply connected Lie group $E_7$, hence $(E_7)^{\sigma}$ is connected. 
\vspace{2mm}

{\bf Theorem 4.11.15.} $(E_7)^{\sigma} \cong (SU(2) \times Spin(12))/\Z_2, \; \Z_2 = \{ (E, 1), (-E, -\sigma)
\vspace{2mm}
 \}$.

{\bf Proof.}  We define a mapping $\varphi : SU(2) \times Spin(12) \to (E_7)^{\sigma}$ by
$$ 
        \varphi(A, \beta) = \varphi_2(A)\beta .$$
Since the Lie algebras $\ga_1$ and $(\gge_7)^{\kappa, \mu}$ of $SU(2)$ and 
$Spin(12)$ are elementwise commutative (Proposition 4.11.4.(2)), $\varphi_2(A) 
\in SU(2)$ and $\beta \in Spin(12)$ are commutative : $\varphi_2(A)\beta = \beta\varphi_2(A)$. Hence $\varphi$ is a homomorphism. We shall show that $\varphi$ 
is onto. Since the group $(E_7)^{\sigma}$ is connected (Lemma 4.11.14), to prove 
this, it is sufficient to show that its differential mapping $\varphi_{*} : \ga_1 \oplus (\gge_7)^{\kappa,\mu} \to (\gge_7)^{\sigma}$ is onto. However this has 
been already shown in Proposition 4.11.4.(2). $\Ker\,\varphi = \{ (E, 1), (-E, \varphi_2(-E)) \} = \{ (E, 1), (-E, -\sigma) \} = \Z_2$ (Theorem 4.11.13) is easily obtained. Thus we have the isomorphism $(SU(2)\times Spin(12))/\Z_2 \cong (E_7)^{\sigma}$.
\vspace{2mm}

{\bf Remark.} We can give an elementary proof of Lemma 4.11.14 once we have proved the following three claims.
\vspace{1mm}
 
{\bf Claim 1.} {\it Any element $X \in (\gJ^C)_{\sigma}$ can be transformed to a diagonal form by some $\alpha \in (E_6)^{\sigma}$}:
$$
       \alpha X = \pmatrix{\xi_1 & 0 & 0 \cr
                           0 & \xi_2 & 0 \cr
                           0 & 0 & \xi_3}, \quad \xi_i \in C. $$
{\it Moreover we can choose $\alpha \in (E_6)^{\sigma}$ so that $\xi_2 \ge 0$, 
$\xi_3 \ge 0$.}
\vspace{2mm}

{\bf Proof.} Recall that $i(E_1 - E_2)^{\sim}$, $i(E_1 - E_3)^{\sim}$, $i\wti{F}_1(a)$, $\wti{A}_1(a) (a \in \gC) \in (\gge_6)^{\sigma}$, then we can prove analogously as in Proposition 3.8.2.
\vspace{1mm}

{\bf Claim 2.} {\it Any element $P \in (\gM^C)_{\sigma} = \{ P \in \gM^C  \, | \, \sigma P = P \}$ can be transformed to a diagonal form by some $\alpha \in ((E_7)^{\sigma})_0$ }({\it the connected component of $(E_7)^{\sigma}$ containing the identity} 1):
$$
      \alpha P = (X, Y, \xi, \eta), \quad X, Y \;\;\mbox{\it are diagonal}, \;\; \xi > 0.$$

{\bf Proof.} Recall ${\mit\Phi}(0, -\tau aE_i, aE_i, 0) \in (\gge_7)^{\sigma}$, 
$i = 1, 2, 3$, then we can prove analogously as in Proposition 4.8.2.
\vspace{3mm}

{\bf Claim 3.}\qquad \qquad $(E_7)^{\sigma}/(E_6)^{\sigma} \simeq (\gM_1)_{\sigma} = \{ P \in \gM_1  \, | \, \sigma P = P \}.$

\noindent {\it In particular, the group $(E_7)^{\sigma}$ is connected.}
\vspace{2mm}

{\bf Proof.} Remark that $\alpha_i(a)$ of Lemma 4.8.1 belongs to $((E_7)^{\sigma})_0$, then this claim is proved analogously as Theorem 4.8.3. The connectedness of $(E_7)^{\sigma}$ follows from the connectedness of $(E_6)^{\sigma} \cong (U(1)\times Spin(10))/\Z_4$ (Theorem 3.10.7) and \vspace{3mm}
$(\gM_1)_{\sigma}$.
\vspace{4mm}

{\bf 4.12. Involution $\tau\gamma$ and subgroup $SU(8)/\Z_2$ of $E_7$}
\vspace{3mm}

We consider the involutive complex conjugate transformation 
$\tau\gamma$ of $\gP^C$,
$$
      \tau\gamma(X, Y, \xi, \eta) = (\tau\gamma X, \tau\gamma Y, \tau \xi, \tau \eta),$$
where $\gamma$ is the same as $\gamma \in G_2 \subset F_4 \subset E_6 \subset E_7$.
\vspace{2mm}

We shall study the following subgroup $(E_7)^{\tau\gamma}$ of $E_7$:
\begin{eqnarray*}
      (E_7)^{\tau\gamma} \!\!\!&=& \!\!\! \{ \alpha \in E_7 \, | \, \tau\gamma\alpha = \alpha\tau\gamma \}
\vspace{1mm}\\
           \!\!\!&=&\!\!\! \{ \alpha \in E_7 \, | \, \lambda\gamma\alpha = \alpha\lambda\gamma \} = (E_7)^{\lambda\gamma}.
\end{eqnarray*}
To this end, we consider $\R$-vector subspaces $(\gP^C)_{\tau\gamma}$, $(\gP^C)_{-\tau\gamma}$ of $\gP^C$, which are eigenspaces of $\tau\gamma$, respectively by
\begin{eqnarray*}
     (\gP^C)_{\tau\gamma} \!\!\! &=& \!\!\! \{ P \in \gP^C \, | \, \tau\gamma P = P \}
\vspace{1mm}\\
    \!\!\! &=& \!\!\! \{ (X, Y, \xi, \eta) \in \gP^C \, | \, X, Y \in (\gJ^C)_{\tau\gamma}, \xi, \eta \in \R \},
\vspace{1mm}\\
     (\gP^C)_{- \tau\gamma} \!\!\!&=&\!\!\! \{ P \in \gP^C  \, | \, \tau\gamma P = - P \}
\vspace{1mm}\\
   \!\!\! &=& \!\!\! \{ (X, Y, \xi, \eta) \in \gP^C  \, | \, X, Y \in (\gJ^C)_{-\tau\gamma}, \xi, \eta \in i\R \},
\vspace{1mm}\\
   \!\!\! &=& \!\!\! i(\gP^C)_{\tau\gamma}.
\end{eqnarray*}
These spaces $(\gP^C)_{\tau\gamma}$, $(\gP^C)_{-\tau\gamma}$ are invariant under the action of $(E_7)^{\tau\gamma}$ and we have the decomposition
$$
     \gP^C = (\gP^C)_{\tau\gamma} \oplus (\gP^C)_{-\tau\gamma}
           = (\gP^C)_{\tau\gamma} \oplus i(\gP^C)_{\tau\gamma}.$$
In paticular, $\gP^C$ is the complexification of $(\gP^C)_{\tau\gamma}$: $\gP^C = ((\gP^C)_{\tau\gamma})^C$.
\vspace{2mm}

Analogious to Section 3.11, we can define the $\R$-linear mapping $k : M(4, \H) \to M(8, \C)$,
$$
           k\biggl(\Big(a + be_2 \Big)\biggl) = \biggl(\pmatrix{a & b \cr
                                           - \ov{b} & \ov{a}}\biggl), \quad a, b \in \C. $$

{\bf Lemma 4.12.1.} {\it Any element $B \in \su(8)$ is uniquely expressed by}
$$
     B = k(D) + e_1k(T), \quad D \in \sp(4), T \in \gJ(4, \H)_0.$$

{\bf Proof.} For $B \in \su(8)$, let $D_1 = \dfrac{B - J\ov{B}J}{2}, T_1 = \dfrac{B + J\ov{B}J}{2e_1} \in M(8, \C)$, then we have
$$
       B = D_1 + e_1T_1, \quad 
\begin{array}{l}
       {D_1}^* = - D_1, JD_1 = \ov{D}_1J, 
\vspace{1mm}\\
       {T_1}^* = T_1, JT_1 = \ov{T}_1J, \tr(T_1) = 0.
\end{array}$$
Then, $D = k^{-1}(D_1), T = k^{-1}(T_1) \in M(4, \H)$ are the required elements. To prove the uniqueness of the expression, it is sufficient to show that 
$$
        D_1 + e_1T_1 = 0, D_1 \in k(\sp(4)), T_1 \in k(\gJ(4, \H)_0) \quad \mbox{implies} \quad C_1 = T_1 = 0. $$
Certainly, from the condition, we have $JD_1 + e_1JT_1 = 0$, so $\ov{D}_1J + e_1\ov{T}_1J = 0$, and so that $D_1J - e_1T_1J = 0$, that is, $D_1 - e_1T_1 = 0$. Together with the first equation, we have $D_1 = T_1 = 0$.
\vspace{2mm}

After this, we will use the $C$-linear mapping $g : \gJ^C \to \gJ(4, \H)^C$, $g(M + \a) 
\vspace{2mm}
= \pmatrix{\dfrac{1}{2}\tr(M) & i\a \cr
             i\a^* & M - \dfrac{1}{2}\tr(M)E}$, 
the homomorphism $\varphi : Sp(4) \to (E_6)^{\tau\gamma}$, $\varphi(A)X 
\vspace{2mm}
 = g^{-1}(A(gX)A^*)$, $X \in \gJ^C$ and its differential mapping $\varphi_* : \sp(4) \to (\gge_6)^{\tau\gamma}$, $\varphi_{*}(D)X = g^{-1}(D(gX) + (gX)D^*)$, $X \in \gJ^C$ which are defined in Section 3.12.
\vspace{3mm}

{\bf Proposition 4.12.2.} {\it The Lie algebra $(\gge_7)^{\tau\gamma}$ of the group $(E_7)^{\tau\gamma}$ is}
\begin{eqnarray*}
     (\gge_7)^{\tau\gamma} \!\!\!&=&\!\!\! \{ {\mit\Phi} \in \gge_7  \, | \, \tau\gamma{\mit\Phi} = {\mit\Phi}\tau\gamma \}
\vspace{1mm}\\
    \!\!\!&=&\!\!\! \{ {\mit\Phi}(\phi, A, -\gamma A, 0) \in \gge_7  \, | \, \phi \in (\gge_6)^{\tau\gamma}, A \in (\gJ^C)_{\tau\gamma} \}
\vspace{1mm}\\
    \!\!\!&=&\!\!\! \{ {\mit\Phi}(\varphi_*(D), g^{-1}(T), -\gamma g^{-1}(T), 0) \in \gge_7  \, | \, D \in \sp(4), T \in \gJ(4, \H)_0 \}.
\end{eqnarray*}
{\it The Lie bracket} $[{\mit\Phi}_1, {\mit\Phi}_2]$ {\it in $(\gge_7)^{\tau\gamma}$ is given by}
$$
   [{\mit\Phi}(\phi_1, A_1, -\gamma A_1, 0), {\mit\Phi}(\phi_2, A_2, -\gamma A_2, 0)] = {\mit\Phi}(\phi, A, -\gamma A, 0),$$
{\it where}
$$
\left\{\begin{array}{ccl}
       \phi \!\!\!&=&\!\!\! [\phi_1, \phi_2] - 2A_1 \vee \gamma A_2  + 2A_2 \vee \gamma A_1,\vspace{1mm}\\
       A \!\!\!&=&\!\!\! \phi_1A_2 - \phi_2A_1.
\end{array}\right.$$

{\bf Proof.} It is not difficult to verity them.. 
\vspace{2mm}

Analogous to Section 3.11, we define a $\C$-vector space $\gS(8, \C)$ by
$$ 
     \gS(8, \C) = \{ S \in M(8, \C) \, | \ {}^tS = -S \},$$
and a $C$-linear mapping $k_J : \gJ(4, \H)^C \to \gS(8, \C)^C$ by
$$
     k_J(M_1 + iM_2) = k(M_1)J + ik(M_2)J, \quad M_1, M_2 \in \gJ(4, \H),$$ 
where $J = \diag(J, J, J, J) \in M(8, C), J = \pmatrix{0 & 1 \cr
                                                       -1 & 0}$.
\vspace{2mm}

{\bf Definition.} We define a $C$-linear isomorphism $\chi : \gP^C \to \gS(8,\C)^C$ by
$$
       \chi (X, Y, \xi, \eta) = k_J\Big(gX - \dfrac{\xi}{2}E \Big) + e_1k_J\Big(g(\gamma Y) - \dfrac{\eta}{2}E \Big). $$

{\bf Proposition 4.12.3.} \quad \qquad $(\gge_7)^{\tau\gamma} \cong \su(8)$. 
\vspace{1mm}

{\it This isomorphism is given by the mapping} $\varphi_* : \su(8) \to (\gge_7)^{\tau\gamma}$,
$$
    \varphi_*(B)P = \chi^{-1}(B(\chi P) + (\chi P){}^t\!B), \quad P \in \gP^C. $$  
{\bf Proof.} We first prove that for $B \in \su(8)$ we have $\varphi_*(B) \in (\gge_7)^{\tau\gamma}$.
\vspace{1mm}

(i) For $B = k(D)$, $D \in \sp(4)$, we have
\vspace{2mm}

\quad
      $P = (X, Y, \xi, \eta)$
\vspace{1mm}

\qquad
   $\stackrel{\chi}{\longrightarrow} 
    k\Big(gX - \dfrac{\xi}{2}E \Big)J + e_1k\Big(g(\gamma Y) - \dfrac{\eta}{2}E \Big)J$
\vspace{1mm}

\qquad
    $\longrightarrow 
     k(D)k\Big(gX-\dfrac{\xi}{2}E \Big)J + e_1k(D)k\Big(g(\gamma Y) - \dfrac{\eta}{2}E \Big)J$
\vspace{1mm}

\qquad \quad
    $+ k\Big(gX - \dfrac{\xi}{2}E \Big)J\,{}^tk(D) + e_1k\Big(g(\gamma Y) - \dfrac{\eta}{2}E \Big)J\,{}^tk(D)$
\vspace{1mm}

\qquad 
   $= k\Big(D\Big(gX - \dfrac{\xi}{2}E \Big)\Big)J + e_1k\Big(D\Big(g(\gamma Y) - \dfrac{\eta}{2}E \Big)\Big)J$
\vspace{1mm}

\qquad \quad
   $+ k\Big(\Big(gX - \dfrac{\xi}{2}E \Big)D^* \Big)J + e_1k\Big(\Big(g(\gamma Y) - \dfrac{\eta}{2}E \Big)D^* \Big)J$
\vspace{1mm}

\qquad
   $= k(D(gX) + (gX)D^*)J + e_1k(D(g(\gamma Y) + (g(\gamma Y))D^*)J$
\vspace{1mm}

\qquad
   $= k(g(\varphi_*(D)X))J + e_1k(g(\varphi_*(D)(\gamma Y)))J$
\vspace{1mm}

  $(\mbox{recall} \;\;\varphi_{*}(C)X = g^{-1}(C(gX)+(gX)C^{*}))$ \qquad
\vspace{1mm}

\qquad
   $= \chi\pmatrix{\varphi_*(D)X \cr
                   \gamma\varphi_*(D)\gamma Y \cr
                   0 \cr 0}
    = \chi\Big(\pmatrix{\varphi_*(D)X & 0 & 0 & 0 \cr
                        0 & \tau\varphi_*(D)\tau  & 0 & 0 \cr
                        0 & 0 & 0 & 0 \cr 
                        0 & 0 & 0 & 0}
               \pmatrix{X \cr Y \cr
                        \xi \cr \eta}\Big)$
\vspace{1mm}

\qquad
   $= \chi({\mit\Phi}(\varphi_*(D), 0, 0, 0)P)$.
\vspace{2mm}

\noindent Hence, we have $\varphi_*(k(D)) = {\mit\Phi}(\varphi_*(D), 0, 0, 0) \in (\gge_7)^{\tau\gamma}$.
\vspace{1mm}

(ii) For $B = e_1k(T)$, $T \in \gJ(4, \H)_0$  (denote $T = gA$, $A \in (\gJ^C)_{\tau\gamma}$)
\vspace{1mm}

\qquad
   $P = (X, Y, \xi, \eta)$
\vspace{1mm}
 
\qquad
    $\stackrel{\chi}\longrightarrow 
    k\Big(gX - \dfrac{\xi}{2}E \Big)J + e_1k\Big(g(\gamma Y) - \dfrac{\eta}{2}E \Big)J$
\vspace{1mm}

\qquad
  $\longrightarrow 
    e_1k(T)k\Big(gX - \dfrac{\xi}{2}E \Big)J - k(T)k\Big(g(\gamma Y) - \dfrac{\eta}{2}E \Big)J$
\vspace{1mm}

\qquad  \qquad
    $+ e_1k\Big(gX - \dfrac{\xi}{2}E \Big)J\,{}^tk(T) - k\Big(g(\gamma Y) - \dfrac{\eta}{2}E \Big)J{}^tk(T)$ 
\vspace{1mm}

\qquad $\;$
  $= \, k(-Tg(\gamma Y) - g(\gamma Y)T + \eta T)J + e_1k(T(gX) + (gX)T - \xi T)J$
\vspace{1mm}

\qquad $\;$
  $= \, k(-2gA \circ g(\gamma Y) + \eta gA)J + e_1k(2gA \circ gX - \xi gA)J$
\vspace{1mm}

\qquad $\;$
   $= \, k\Big(-2g(\gamma A \times Y) - \dfrac{1}{2}(A, Y)E + \eta gA \Big)J$
\vspace{1mm}

\qquad \qquad
    $+ e_1k\Big(2g(\gamma A \times \gamma X) + \dfrac{1}{2}(\gamma A, X)E - \xi gA \Big)J \,(\mbox{Lemma 3.12.1})$
\vspace{1mm}

\qquad $\;$
   $= \, \chi\pmatrix{-2\gamma A \times Y + \eta A \cr
                    2A\times X - \xi\gamma A \cr
                    (A, Y) \cr
                    (-\gamma A, X)}
   = \chi\biggl(\pmatrix{0 & -2\gamma A & 0 & A \cr 
                         2A & 0 & -\gamma A & 0 \cr
                         0 & A & 0 & 0 \cr
                         -\gamma A & 0 & 0 & 0}
                \pmatrix{X \cr Y \cr
                         \xi \cr \eta}\Big)$
\vspace{1mm}

\qquad $\;$
   $= \, \chi({\mit\Phi}(0, A, -\gamma A, 0)P)$.
\vspace{2mm}

\noindent Hence we have $\varphi_*(e_1k(T)) = {\mit\Phi}(0, A, -\gamma A, 0) \in (\gge_7)^{\tau\gamma}$.
\vspace{1mm}

\noindent Consequently we see that the mapping $\varphi : \su(8) \to (\gge_7)^{\tau\gamma}$ is well-defined. We can easily check that $\varphi_*$ is onto. Finally, we have to prove that $\varphi_*$ is a homomorphism as Lie algebras. However this follows from the following Theorem 4.12.5, so that we shall omit the proof.
\vspace{3mm}

{\bf Lemma 4.12.4.} {\it The group $(E_7)^{\tau\gamma}$ is connected.}
\vspace{2mm}

{\bf Proof.} The group $(E_7)^{\tau\gamma}$ is the fixed points set by the involution $\tau\gamma$ of the simply connected Lie group $E_7$, hence $(E_7)^{\tau\gamma}$ is connected.
\vspace{3mm}

{\bf Theorem 4.12.5.} \qquad $(E_7)^{\tau\gamma} \cong SU(8)/\Z_2, \;\;\Z_2 =\{ E, -E \}$.
\vspace{2mm}

{\bf Proof.} We define a mapping $\varphi : SU(8) \to (E_7)^{\tau\gamma}$ by
$$
      \varphi(A)P = \chi^{-1}(A(\chi P)\,{}^t\!A), \quad P \in \gP^C. $$
We first prove $\varphi(A) \in (E_7)^{\tau\gamma}$. To prove this, for the 
differential mapping $\varphi_{*} : \su(8) \to (\gge_7)^{\tau\gamma}$ of $\varphi$,
$$ 
     \varphi_{*}(D)P = \chi^{-1}(D(\chi P) + (\chi P){}^tD), \quad P \in \gP^C,$$
it is sufficient to show that $\varphi_{*}$ is well-defined, that is, $\varphi_{*}(D) \in (\gge_7)^{\tau\gamma}$. However this fact is already shown in Proposition 4.12.3. Evidently $\varphi : SU(8) \to (E_7)^{\tau\gamma}$ is a homomorphism.  Since $\varphi_* : \su(8) \to (\gge_7)^{\tau\gamma}$ is onto and $(E_7)^{\tau\gamma}$ is connected (Lemma 4.12.4), $\varphi : SU(8) \to (E_7)^{\tau\gamma}$ is also onto. $\Ker\varphi = \{ E, - E \} = \Z_2$ is easily obtained. Thus we have the isomorphism $SU(8)/\Z_2 \cong (E_7)^{\tau\gamma}$.
\vspace{2mm}

{\bf Remark.} Without using Lemma 4.12.4, the fact that the mapping $\varphi : SU(8) \to (E_7)^{\tau\gamma}$ is onto will be followed from two claims. 
\vspace{1mm}

{\bf Claim 1.} {For $a \in \R$, $\alpha_i (a)$ of Lemma 4.8.1 belongs to $\varphi (SU(8))$}.
\vspace{1mm}

{\bf Proof.} $\alpha_i(a) = \exp({\mit\Phi}(0, -aE_i, aE_i, 0)) \in \exp\varphi_*(\su(8))$ (Proposition 4.12.3) $= \varphi(\exp(\su(8))) \in \varphi(SU(8))$.
\vspace{1mm}

{\bf Claim 2.} {\it Any element $P \in (\gM^C)_{\tau\gamma} = \{ P \in \gM^C  \, | \, \tau\gamma P = P \}$ can be transformed to a diagonal form by some $\alpha \in \varphi(SU(8))$}:
$$
      \alpha P = (X, Y, \xi, \eta), \quad X, Y \;\;\mbox{{\it are real diagonal},} \;\; \xi > 0.$$

{\bf Proof.} Let $P = (X, Y, \xi, \eta) \in (\gM^C)_{\tau\gamma}$. If $\xi \neq 0$. Then
$$
      \tau\gamma Y = Y, \quad X = \dfrac{1}{\xi}(Y \times Y), \quad \tau\xi = \xi, \quad \tau\eta = \eta. $$
Since $Y \in (\gJ^C)_{\tau\gamma}$, we have $\gamma Y \in (\gJ^C)_{\tau\gamma}$ , so that $g(\gamma Y) \in \gJ(4, \H)_0$. Hence, there exists $D \in Sp(4)$ such that 
$$
     D(g(\gamma Y))D^{*}  \quad \mbox{is real diagonal.}$$
Then 
$$
\begin{array}{c}
      \gamma\varphi(D)\gamma Y = g^{-1}(D(g(\gamma Y))D^{*}) \quad \mbox{is real diagonal},
\vspace{1mm}\\
   \varphi(D)X = \varphi(D)\Big(\dfrac{1}{\xi}Y \times Y \Big) = \dfrac{1}{\xi}(\gamma\varphi(D)\gamma Y \times \gamma\varphi(D) \gamma Y) \;\; \mbox{is real diagonal}.
\end{array}$$ 
In the case $\xi = 0$, by the same proof of Proposition 4.8.2, we can choose $\alpha \in \varphi(SU(8))$ so that
$$
        \alpha P = (X, Y, \xi, \eta ), \quad X, Y \;\; \mbox{are real diagonal},\;\; 0 \neq \xi \in \R.$$
If $\xi < 0$, apply $\alpha_1(\pi)$ of Claim 1 on it, then $\xi$ becomes $\xi > 0$. 
\vspace{2mm}

Now, we will return to the proof of the surjection of $\varphi : SU(8) \to (E_7)^{\tau\gamma}$ using Claims 1, 2. For a given $\alpha \in (E_7)^{\tau\gamma}$, consider the element $P = \alpha\dot{1} \in (\gM^C)_{\tau\gamma}$. We first transform $P$ to a diagonal form (Claim 2) by some $\beta \in \varphi(SU(8))$, and we have, in a similar way to Theorem 4.8.3,
$$
     \alpha_1(a_1)^{-1}\alpha_2(a_2)^{-1}\alpha_3(a_3)^{-1}\beta\alpha\dot{1} = \dot{1}, $$ 
where $a_i = \dfrac{\eta_i}{|\eta_i|}r_i$ ($\eta_i$ is a diagonal element of $Y$). Since $\eta_i \in \R$, we have $\alpha_i(a_i) \in \varphi(SU(8))$ (Claim 1). If we put $\wti{\alpha} = \alpha_1(a_1)^{-1}\alpha_2(a_2)^{-1}\alpha_3(a_3)^{-1}\beta\alpha$, then, $\wti{\alpha} \in E_6$ (Theorem 4.7.2) and $\wti{\alpha}$ satisfies $\tau\gamma\wti{\alpha} = \wti{\alpha}\tau\gamma$. Hence $\wti{\alpha} \in (E_6)^{\tau\gamma} = \varphi(Sp(4))$ (Theorem 3.12.2) $\subset \varphi(SU(8))$, Therefore $\alpha = \beta^{-1}\alpha_3(a_3)\alpha_2(a_2)\alpha_1(a_1)\wti{\alpha} \in \varphi(SU(8))$. This shows that $\varphi$ is onto. 
\vspace{4mm} 

{\bf 4.13. Automorphism $w$ of order 3 and subgroup $(SU(3) \times SU(6))/\Z_3$ of 
\vspace{3mm}
$E_7$}

We define a $C$-linear transformation $w$ of order 3 of $\gP^C$ by
$$
       w(X, Y, \xi, \eta) = (wX, wY, \xi, \eta). $$
This $w$ is the same as $w \in G_2 \subset F_4 \subset E_6 \subset E_7$.
\vspace{2mm}

We shall study the following subgroup $(E_7)^w$ of $E_7$:
$$
          (E_7)^w = \{\alpha \in E_7 \, | \, w\alpha = \alpha w \}. $$

We consider the group $E_{7,{\sC}}$ replaced with $\C$ in the place $\gC$ in the definition of the group $E_7$:
$$
    E_{7,{\sC}} = \{\alpha \in \Iso_C((\gP_{\sC})^C) \, | \, \alpha(P \times Q)\alpha^{-1} = \alpha P \times \alpha Q, \langle \alpha P, \alpha Q \rangle = \langle P, Q \rangle \}. $$
As in Section 4.7, the group $E_{7,{\sC}}$ contains a subgroup
$$
             E_{6,{\sC}} = \{\alpha \in E_{7,{\sC}} \, | \, \alpha(0, 0, 1, 0) = (0, 0, 1, 0) \}, $$
which is isomorphic to the group $((SU(3) \times SU(3))/\Z_3)\cdot\Z_2$ (Proposition 3.13.4).   
The Lie algebra $\gge_{7,{\sC}}$ of the group $E_{7,{\sC}}$ is given by
$$
   \gge_{7,{\sC}} = \{{\mit\Phi}(\phi, A, -\tau A, \nu) \, | \, \phi \in \gge_{6,{\sC}}, A \in (\gJ_{\sC})^C, \nu \in i\R \} $$
(Theorem 4.3.4). In particular, the dimension of $\gge_{7,{\sC}}$ is
$$
        \dim\gge_{7,{\sC}} = 16 + 18 + 1 = 35. $$
As in Theorem 4.8.3, we see that the space
$$
       (\gM_{\sC})_1 = \{P \in ({\gM_{\sC}})^C \, | \, P \times P = 0, \langle P, P \rangle = 1 \} $$
is connected and we have the homeomorphism
$$
                   E_{7,{\sC}}/E_{6,{\sC}} \simeq (\gM_{\sC})_1. 
\vspace{2mm}$$

{\bf Lemma 4.13.1.} {\it $E_{7,{\sC}}$ has at most two connected components} (in reality has two connected components).
\vspace{2mm}

{\bf Proof.} From the exact sequence $\pi_0(E_{6,{\sC}}) \to \pi_0(E_{7,{\sC}}) \to \pi_0((\gM_{\sC})_1)$, that is, $\Z_2 \to \pi_0(E_{7.{\sC}}) \to 0$ (Proposition 3.13.4), we see that $\pi_0(E_{7,{\sC}})$ is 0 or $\Z_2$.  
\vspace{2mm}

Let $h' : C \to \C$ be the $\R$-linear isomorphism defined by
$$
            h'(a + bi) = a + be_1, \quad a, b \in \R. $$

Now, let $V, W$ be $C$- and $\C$-vector spaces, respectively. A linear mapping $f : V \to W$ is called a $C$-$\C$-linear mapping if
$$
        f(av) = h'(a)f(v), \quad a \in C, v \in V. $$
Similarly, a $\C$-$C$-linear mapping $g : W \to V$ is defined.
\vspace{3mm}

{\bf Definition.} Let $h' : \C^C \to \C$ be a $C$-$\C$-linear mapping defined by$$
         h'(a + bi) = a + be_1, \quad a, b \in \C. $$

Now, let ${\mit\Lambda}^3(\C^6)$ be the third exterior product of $\C$-vector space $\C^6$ and we define a $C$-$\C$-linear isomorphism $f : (\gP_{\sC})^C \to {\mit\Lambda}^3(\C^6)$ by
$$
   f\Big(\pmatrix{\xi_1 & x_3 & \ov{x}_2 \cr
                  \ov{x}_3 & \xi_2 & x_1 \cr
                  x_2 & \ov{x}_1 & \xi_3},
         \pmatrix{\eta_1 & y_3 & \ov{y}_2 \cr
                  \ov{y}_3 & \eta_2 & y_1 \cr
                  y_2 & \ov{y}_1 & \eta_3}, \xi, \eta \Big)
     = \dsum_{i<j<k}x_{ijk}\e_i \land \e_j \land \e_k $$
$\Big(\{\e_1, \e_2, \cdots, \e_6 \}$ is the canonical basis of $\C^6$ and $x_{ijk} \in \C$ are skew-symmetric tensor: $x_{i'j'k'} = \mbox{sgn}\pmatrix{i & j & k \cr
                                               i' & j' & k'}x_{ijk}\Big)$, where\vspace{2mm}

$\qquad \qquad \qquad \quad
      x_{156} = h'(\xi_1), \quad x_{164} = h'(x_3), \quad x_{145} = h'(\ov{x}_2)$, 
\vspace{1mm}

$\qquad \qquad \qquad \quad
      x_{256} = h'(\ov{x}_3), \quad x_{264} = h'(\xi_2), \quad x_{245} = h'(x_1)$, 
\vspace{1mm}

$\qquad \qquad \qquad \quad
      x_{356} = h'(x_2), \quad x_{364} = h'(\ov{x}_1), \quad x_{345} = h'(\xi_3)$, 
\vspace{1mm}

$\qquad \qquad \qquad \quad
      x_{423} = h'(\eta_1), \quad x_{431} = h'(y_3), \quad x_{412} = h'(\ov{y}_2)$, 
\vspace{1mm}

$\qquad \qquad \qquad \quad
      x_{523} = h'(\ov{y}_3), \quad x_{531} = h'(\eta_2), \quad x_{512} = h'(y_1)$, 
\vspace{1mm}

$\qquad \qquad \qquad \quad
      x_{623} = h'(y_2), \quad x_{631} = h'(\ov{y}_1), \quad x_{612} = h'(\eta_3)$, 
\vspace{1mm}

$\qquad \qquad \qquad \qquad \qquad \qquad \quad \quad \;
        x_{123} = h'(\xi)$,  
\vspace{1mm}

$\qquad \qquad \qquad \qquad \qquad \qquad \quad \quad \;
        x_{456} = h'(\eta)$.                        
\vspace{2mm} 

\noindent The inverse mapping $f^{-1} : {\mit\Lambda}^3(\C^6) \to (\gP_{\sC})^C$ of $f$ is given by
$$
   f^{-1}\Big(\dsum_{i<j<k}x_{ijk}\e_i \land \e_j \land \e_k\Big) = 
\pmatrix{
\begin{array}{c}
\pmatrix{
\begin{array}{ccc}
              h(x_{156}) & h(x_{164}, \ov{x}_{256}) & h(x_{145}, \ov{x}_{356})
\vspace{1mm}\cr
              h(x_{256}, \ov{x}_{164}) & h(x_{264}) & h(x_{245}, \ov{x}_{364})
\vspace{1mm}\cr
              h(x_{356}, \ov{x}_{145}) & h(x_{364}, \ov{x}_{245}) & h(x_{345})
\end{array}}
\vspace{3mm}\cr
\pmatrix{
\begin{array}{ccc}
            h(x_{423}) & h(x_{431}, \ov{x}_{523}) & h(x_{412}, \ov{x}_{623})
\vspace{1mm}\cr
            h(x_{523}, \ov{x}_{431}) & h(x_{531}) & h(x_{512}, \ov{x}_{631})
\vspace{1mm}\cr
            h(x_{623}, \ov{x}_{412}) & h(x_{631}, \ov{x}_{512}) & h(x_{612})
\end{array}}
\vspace{2mm}\cr
             h(x_{123})  
\vspace{1mm}\cr
             h(x_{456}) 
\end{array}},$$
where $ h : \C \oplus \C \to \C^C, h : \C \to C$ are $\C$-$C$-linear mappings defined respectively by
$$
\begin{array}{l}
     h(a, b) = \dfrac{a + b}{2} + i\dfrac{(b - a)e_1}{2}, \quad a, b \in \C,
\vspace{1mm}\\
     h(a + be_1) = a + bi, \quad a, b \in \R.
\end{array} $$
It is easy to see that
$$
      f(h(a)P) = a(fP), \quad a \in \C, P \in (\gP_{\sC})^C. $$

The group $SU(6)$ acts naturally on ${\mit\Lambda}^3(\C^6)$, that is, the action of $A \in SU(6)$ on $\a \land \b \land \c \in {\mit\Lambda}^3(\C^6)$ is defined by
$$
        A(\a \land \b \land \c) =  A\a \land A\b \land A\c. $$
Hence, the action of $D \in \su(6)$ on ${\mit\Lambda}^3(\C^6)$ is given by
$$
       D(\a \land \b \land \c) = D\a \land \b \land \c + \a \land D\b \land \c + \a \land \b \land 
\vspace{2mm}
D\c. $$

{\bf Lemma 4.13.2.} (1) {\it Any element $D \in \su(6)$ is uniquely expressed by}
$$
    D = \pmatrix{B & L \cr
                - L^* & C} + \dfrac{\nu}{3}\pmatrix{E & 0 \cr
                                                    0 & -E}, \quad B, C \in \su(3), L \in M(3, \C), \nu \in e_1\R. $$

(2) {\it The Lie algebra $\gge_{7,{\sC}}$ is isomorphic to the Lie algebra $\su(6)$ as Lie algebras}:
$$
      \gge_{7,{\sC}} \cong \su(6). $$
{\it This isomorphism is given by the mapping ${\varphi_{\sC}} : \su(6) \to \gge_{7,{\sC}}$, 
$$
   {\varphi_{\sC}}\Big(\pmatrix{B & L \cr
                - L^* & C} + \dfrac{\nu}{3}\pmatrix{E & 0 \cr
                                                    0 & -E}\Big) =
   {\mit\Phi}({\phi_{\sC}}(B, C), h(L), - \tau h(L), -i\nu e_1) $$ 
where ${\phi_{\sC}} : \su(3) \oplus \su(3) \to \gge_{6,{\sC}}$ is defined by}
$
    {\phi_{\sC}}(B, C)X = h(B, C)X + Xh(B, C)^*, $ $X \in ({\gJ_{\sC}})^C$
 (Lemma 3.13.3).
\vspace{2mm}

{\bf Proof.} (1) For $D = \pmatrix{B' & L \cr
                                 - L^* & C'} \in \su(6), B', C' \in \gu(3), L \in M(3, \C)$, we let
$$
    \nu = \tr(B') = - \tr(C'), \quad B = B' - \dfrac{\nu}{3}E, \quad C = C' + \dfrac{\nu}{3}E,$$
then we have the result.
\vspace{1mm}

(2) This is the direct consequence of the following Proposition 4.13.3, so we will omit its proof.
\vspace{2mm}

We define the action of the group $\Z_2 = \{1, \epsilon \}$ on the group $SU(6)$ by
$$
    \epsilon A = \ov{(\mbox{Ad}J_3)A}, \quad J_3 = \pmatrix{0 & E \cr
                                                            -E & 0}, $$
that is,
$$
    \epsilon A = \epsilon\pmatrix{A_{11} & A_{12} \cr
                                  A_{21} & A_{22}} = \!
   \ov{\pmatrix{0 & E \cr
                -E & 0}\pmatrix{A_{11} & A_{12} \cr
                                A_{21} & A_{22}}\pmatrix{0 & E \cr
                                                         -E & 0}^{-1}} \!= \!
        \pmatrix{\ov{A}_{22} \!&\! -\ov{A}_{21} \vspace{1mm}\cr
                 -\ov{A}_{12} & \ov{A}_{11}}, $$
where $E, A_{ij} \in M(3, \C)$, and let $SU(6)\cdot\Z_2$ be the semi-direct product of the groups $SU(6)$ and $\Z_2$ under this action.
\vspace{3mm}

{\bf Proposition 4.13.3.} $E_{7,{\sC}} \cong (SU(6)/\Z_3)\cdot\Z_2, \; \Z_3 = \{E, \omega_1E, {\omega_1}^2E \}$, $\omega_1 = - \dfrac{1}{2} + \dfrac{\sqrt{3}}{2}e_1$.
\vspace{2mm}

{\bf Proof.} We define a mapping $\psi : SU(6)\cdot\Z_2 \to E_{7,{\sC}}$ by
$$
     \psi(A, 1)P = f^{-1}(A(fP)), \quad  \psi(A, \epsilon)P = f^{-1}(A(f\ov{P})), \quad P \in (\gP_{\sC})^C. $$       
We first have to show that $\psi(A, 1) \in E_{7,{\sC}}$. To prove this, it is sufficient to show that the differential mapping $\psi_* : \su(6) \to \gge_{7,{\sC}}$ of $\psi$:
$$
         \psi_*(D)P= f^{-1}(D(fP)), \quad P \in (\gP_C)^C $$
coincides with the mapping $\psi_{\sC} : \su(6) \to \gge_{7,{\sC}}$ of Lemma 4.13.2. We put
$$
D = \pmatrix{B & L \cr
                 -L^* & C} + \dfrac{\nu}{3}\pmatrix{E & 0 \cr
                                                    0 & -E} \qquad \qquad \qquad \qquad \qquad \qquad \qquad \qquad $$
$$
\qquad 
   = \pmatrix{\begin{array}{llllll}
                    \hfill b_{11} & \hfill b_{12} & \hfill b_{13} & \hfill l_{11} & \hfill l_{12} & l_{13} \cr           
                  -\ov{b}_{12} & \hfill b_{22} & \hfill b_{23} & \hfill l_{21} & \hfill l_{22} & l_{23} \cr  
                  -\ov{b}_{13} & \hfill b_{23} & \hfill b_{33} & \hfill l_{31} & \hfill l_{32} & l_{33} \cr  
                  -\ov{l}_{11} & -\ov{l}_{21} & -\ov{l}_{31} & \hfill c_{11} & \hfill c_{12} & c_{13} \cr
                  -\ov{l}_{12} & -\ov{l}_{22} & -\ov{l}_{32} & -\ov{c}_{12} & \hfill c_{22} & c_{23} \cr
           -\ov{l}_{13} & -\ov{l}_{23} & -\ov{l}_{33} & -\ov{c}_{13} & -\ov{c}_{23} & c_{33} \cr
                  \end{array}}
        + \dfrac{\nu}{3}\pmatrix{E & 0 \cr
                                 0 & -E} \in \su(6), 
$$
where $\ov{b}_{ii} = - b_{ii}, \ov{c}_{ii} = - c_{ii}, b_{11} + b_{22} + b_{33} = c_{11} + c_{22} + c_{33} = 0, \ov{\nu} = - \nu$. 
\vspace{1mm}

(1) For $P = (0, 0, 1, 0), \psi_*(D)P$ is calculated as follows.
\vspace{2mm}

$\quad
  P = (0, 0, 1, 0)$
\vspace{1mm}

$\quad
   \stackrel{f} {\longrightarrow} \e_1 \land \e_2 \land \e_3$
\vspace{1mm}

$\quad
   \stackrel{D} {\longrightarrow} D\e_1 \land \e_2 \land \e_3 + \e_1 \land D\e_2 \land \e_3 + \e_1 \land \e_2 \land D\e_3$
\vspace{1mm}

$\quad
   = \Big(b_{11} + \dfrac{\nu}{3}\Big)\e_1 \land \e_2 \land \e_3 - \ov{l}_{11}\e_4 \land \e_2 \land \e_3 - \ov{l}_{12}\e_5 \land \e_2 \land \e_3 - \ov{l}_{13}\e_6 \land \e_2 \land \e_3$
\vspace{1mm}

$\quad
    \quad + \Big(b_{22} + \dfrac{\nu}{3}\Big)\e_1 \land \e_2 \land \e_3 - \ov{l}_{21}\e_4 \land \e_3 \land \e_1 - \ov{l}_{22}\e_1 \land \e_5 \land \e_3 - \ov{l}_{23}\e_1 \land \e_6 \land \e_3$
\vspace{1mm}

$\quad
   \quad + \Big(b_{33} + \dfrac{\nu}{3}\Big)\e_1 \land \e_2 \land \e_3 - \ov{l}_{31}\e_4 \land \e_1 \land \e_2 - \ov{l}_{32}\e_5 \land \e_1 \land \e_2 - \ov{l}_{33}\e_6 \land \e_1 \land \e_2$
\vspace{1mm}

$\quad
   \stackrel{f^{-1}} {\longrightarrow} \pmatrix{\begin{array}{c}
\pmatrix{0 & 0 & 0 \cr
         0 & 0 & 0 \cr
         0 & 0 & 0}
\vspace{1mm}\cr
\pmatrix{\begin{array}{ccc}
         -h(\ov{l}_{11}) & -h(\ov{l}_{21}, l_{12}) & -h(\ov{l}_{31}, l_{13}) \vspace{1mm}\cr
         -h(\ov{l}_{12}, l_{21}) & -h(\ov{l}_{22}) & -h(\ov{l}_{32}, l_{23}) \vspace{1mm}\cr
         -h(\ov{l}_{13}, l_{31}) & -h(\ov{l}_{23}, l_{32}) & -h(\ov{l}_{33})
\end{array}} 
\vspace{1mm}\cr
           h(\nu)
\vspace{1mm}\cr
            0
\end{array}}
 = \pmatrix{0 \vspace{3mm}\cr
            -\tau h(L) \vspace{3mm}\cr
            - i\nu e_1 \vspace{3mm}\cr
             0}$ 
\vspace{1mm}

$\quad
 = \pmatrix{\phi_{\sC}(B, C) + \dfrac{i\nu e_1}{3} & -2\tau h(L) & 0 & h(L) 
\vspace{1mm}\cr
            2h(L) & \tau\phi_{\sC}(B, C)\tau - \dfrac{i\nu e_1}{3} & - \tau h(L) & 0
\vspace{1mm}\cr
            0 & h(L) & -i\nu e_1 & 0 
\vspace{1mm}\cr
            -\tau h(L) & 0 & 0 & i\nu e_1}\pmatrix{0 \vspace{3mm}\cr 0 \vspace{3mm}\cr 1 \vspace{3mm}\cr 0}$
\vspace{2mm}

$\quad
   = {\mit\Phi}(\phi_{\sC}(B, C), h(L), - \tau h(L), -i\nu e_1)P$. 
\vspace{1mm}

(2) For $P = (E_1, 0, 0, 0), \psi_*(D)P$ is calculated as follows.
\vspace{2mm}

$\quad
  P = (E_1, 0, 0, 0)$
\vspace{1mm}

$\quad
   \stackrel{f} {\longrightarrow} \e_1 \land \e_5 \land \e_6$
\vspace{1mm}

$\quad
   \stackrel{D} {\longrightarrow} D\e_1 \land \e_5 \land \e_6 + \e_1 \land D\e_5 \land \e_6 + \e_1 \land \e_5 \land D\e_6$
\vspace{1mm}

$\quad
   = \Big(b_{11} + \dfrac{\nu}{3}\Big)\e_1 \land \e_5 \land \e_6 - \ov{b}_{12}\e_2 \land \e_5 \land \e_6 - \ov{b}_{13}\e_3 \land \e_5 \land \e_6 - \ov{l}_{11}\e_4 \land \e_5 \land \e_6$
\vspace{1mm}

$\quad
    \quad + l_{22}\e_1 \land \e_2 \land \e_6 + l_{32}\e_1 \land \e_3 \land \e_6 + c_{12}\e_1 \land \e_4 \land \e_6 + \Big(c_{22} - \dfrac{\nu}{3}\Big)\e_1 \land \e_5 \land \e_6$
\vspace{1mm}

$\quad
    \quad + l_{23}\e_1 \land \e_5 \land \e_2 + l_{33}\e_1 \land \e_5 \land \e_3 + c_{13}\e_1 \land \e_5 \land \e_4 + \Big(c_{33} - \dfrac{\nu}{3}\Big)\e_1 \land \e_5 \land \e_6$
\vspace{2mm}

\vspace{2mm}

$\quad
 \stackrel{f^{-1}} {\longrightarrow} \pmatrix{\begin{array}{c}
\pmatrix{\begin{array}{ccc}
         h\Big(b_{11} - c_{11} - \dfrac{\nu}{3}\Big) & -h(c_{12}, b_{12}) & -h(c_{13}, b_{13}) \vspace{1mm}\cr
         -h(\ov{b}_{12}, \ov{c}_{12}) & 0 & 0 \vspace{1mm}\cr
         -h(\ov{b}_{13}, \ov{c}_{13}) & 0 & 0
\end{array}} 
\vspace{2mm}\\
  \pmatrix{0 & 0 & 0 \vspace{1mm}\cr
         0 & h(l_{33}) & -h(l_{23}, \ov{l}_{32}) \vspace{1mm}\cr
         0 & -h(l_{32}, \ov{l}_{23}) & h(l_{22})}
\vspace{1mm}\cr       
    0
\vspace{1mm}\cr
            -h(\ov{l}_{11})
\end{array}} $
\vspace{2mm}

$\quad  
= \pmatrix{\phi_{\sC}(B, C)E_1 + \dfrac{i\nu e_1}{3}E_1
 \vspace{1mm}\cr
            2h(L) \times E_1 
 \vspace{1mm}\cr
             0 
\vspace{1mm}\cr
        - (\tau h(L), E_1)}
    = {\mit\Phi}(\phi_{\sC}(B, C), h(L), -\tau h(L), - i\nu e_1)
      \pmatrix{E_1 \vspace{1mm}\cr
      0  \vspace{1mm}\cr
      0 \vspace{1mm}\cr  0}.$

\vspace{2mm}

(3) For $P = (F_1(1), 0, 0, 0), \psi_*(D)P$ is calculated as follows.
\vspace{1mm}

  $P = (F_1(1), 0, 0, 0)$
\vspace{1mm}

   $\stackrel{f} {\longrightarrow} \e_2 \land \e_4 \land \e_5 + \e_3 \land \e_6 \land \e_4$ 
\vspace{1mm}

   $\stackrel{D} {\longrightarrow} (D\e_2 \land \e_4 \land \e_5 + \e_2 \land D\e_4 \land \e_5 + \e_2 \land \e_4 \land D\e_5)$ 
\vspace{1mm}

\qquad    
       $+ (D\e_3 \land \e_6 \land \e_4 + \e_3 \land D\e_6 \land \e_4 + \e_3 \land \e_6 \land D\e_4)$ 
\vspace{1mm}

   $= \Big(b_{12}\e_1 \land \e_4 \land \e_5 + \Big(b_{22} + \dfrac{\nu}{3}\Big)\e_2 \land \e_4 \land \e_5 - \ov{b}_{23}\e_3 \land \e_4 \land \e_5 - \ov{l}_{23}\e_6 \land \e_4 \land \e_5$
\vspace{1mm}
    
   $\quad + l_{11}\e_2 \land \e_1 \land \e_5 + l_{31}\e_2 \land \e_3 \land \e_5 + \Big(c_{11} - \dfrac{\nu}{3}\Big)\e_2 \land \e_4 \land \e_5 - \ov{c}_{13}\e_2 \land \e_6 \land \e_5$
\vspace{1mm}

\quad 
    $+ l_{12}\e_2 \land \e_4 \land \e_1 + l_{32}\e_2 \land \e_4 \land \e_3 + \Big(c_{22} - \dfrac{\nu}{3}\Big)\e_3 \land \e_6 \land \e_4 - \ov{c}_{23}\e_2 \land \e_4 \land \e_6\Big)$
\vspace{1mm}

\quad
   $+ \Big(b_{13}\e_1 \land \e_6 \land \e_4 + b_{23}\e_2 \land \e_6 \land \e_4
 + \Big(b_{33} + \dfrac{\nu}{3}\Big)\e_3 \land \e_6 \land \e_4 - \ov{l}_{32}\e_5 \land \e_6 \land \e_4$
\vspace{1mm}

\quad
   $ + l_{13}\e_3 \land \e_1 \land \e_4 + l_{23}\e_3 \land \e_2 \land \e_4 + c_{23}\e_3 \land \e_5 \land \e_4 + \Big(c_{33} - \dfrac{\nu}{3}\Big)\e_3 \land \e_6 \land \e_4$
\vspace{1mm}

\quad 
   $+ l_{11}\e_3 \land \e_6 \land \e_1 + l_{21}\e_3 \land \e_6 \land \e_2 + \Big(c_{11} - \dfrac{\nu}{3}\Big)\e_3 \land \e_6 \land \e_4 - \ov{c}_{12}\e_3 \land \e_6 \land \e_5\Big)$    
\vspace{3mm}

 $\stackrel{f^{-1}} {\longrightarrow} \pmatrix{\begin{array}{c}
\pmatrix{\begin{array}{ccc}
         0 & h(b_{13}, c_{13}) & 0 \vspace{1mm}\cr
         * & h(b_{23} + \ov{c}_{23}) & h(b_{22} - c_{33} - \dfrac{\nu}{3}, -b_{33} + c_{33} + \dfrac{\nu}{3}\Big) \vspace{1mm}\cr
         h(\ov{c}_{12}, \ov{b}_{12}) & * & -h(\ov{b}_{23} + c_{23})
\end{array}} 
\vspace{3mm}\\
\pmatrix{-h(l_{23} + l_{32}) & h(l_{13}, \ov{l}_{31}) & * \vspace{1mm}\cr
         * & 0 & -h(l_{11}) \vspace{1mm}\cr
         h(l_{21}, \ov{l}_{12}) & * & 0}
\vspace{1mm}\cr       
    0
\vspace{1mm}\cr
            -h(\ov{l}_{23} + \ov{l}_{32})
\end{array}}
\vspace{2mm}\\
  = \pmatrix{\phi_{\sC}(B, C)F_1(1) + \dfrac{i\nu e_1}{3}F_1(1) 
 \vspace{1mm}\cr
            2h(L) \times F_1(1) 
 \vspace{1mm}\cr
             0 
\vspace{1mm}\cr
        - (h(L), F_1(1))}
    = {\mit\Phi}(\phi_{\sC}(B, C), h(L), -\tau h(L), - i\nu e_1)
      \pmatrix{F_1(1) \vspace{1mm}\cr
      0  \vspace{1mm}\cr
      0 \vspace{1mm}\cr       0}$.
\vspace{2mm}

(4) For the other generator of ${\gP_{\sC}}^C$, that is, $P = (X, 0, 0, 0), (0, X, 0, 0)$ where $X = E_i, F_i(1), F_i(e_1), i = 1, 2, 3$ and $P = (0, 0, 0, 1)$, we have also
$$
   f^{-1}(D(fP)) = {\mit\Phi}(\phi_{\sC}(B, C), h(L), -\tau h(L), -i\nu e_1)P. $$
Thus we see that $\psi(A, 1) \in E_{7,{\sC}}$ for $A \in SU(6)$. Since $\psi(E, \epsilon)P = \ov{P}$, we see $\psi(E, \epsilon) = \epsilon \in G_{2,{\sC}}\, (= \mbox{Aut(\C)})\, \subset F_{4,{\sC}} \subset E_{6,{\sC}} \subset E_{7,{\sC}}$. We shall show that $\psi : SU(6)\cdot\Z_2 \to E_{7,{\sC}}$ is a homomorphism. For this purpose, we first show
$$
     \ov{f^{-1}(A(fP))} = f^{-1}((\ov{(\mbox{Ad}J_3)A}f\ov{P})), \quad A \in SU(6), P \in (\gP_{\sC})^C. $$
Furthermore, to show this, it is sufficient to show that for $D \in \su(6)$ instead of $A \in SU(6)$. Now,
\begin{eqnarray*}
    \ov{f^{-1}(D(fP))} \!\!\! &=& \!\!\! \ov{{\mit\Phi}(\phi_{\sC}(B, C), h(L), - \tau h(L), - i\nu e_1)P}
\vspace{1mm}\\
    \!\!\! &=& \!\!\! {\mit\Phi}(\ov{\phi_{\sC}(B, C)}, \ov{h(L)}, - \tau \ov{h(L)}, -\ov{i\nu e_1})\ov{P}
\vspace{1mm}\\
     \!\!\! &=& \!\!\! {\mit\Phi}(\phi_{\sC}(\ov{C}, \ov{B}), h({\ov{L}}^{\,*}), - \tau h({\ov{L}}^{\,*}), - i\nu e_1)\ov{P}
\vspace{1mm}\\
     \!\!\! &=& \!\!\! f^{-1}(((\ov{\mbox{Ad}J_3)D)}(f\ov{P})).
\end{eqnarray*}
$\psi$ is a homomorphism. Indeed, for example,
$$
\begin{array}{l}
   \psi(A, \epsilon)\psi(B, 1)P = \psi(A, \epsilon)(f^{-1}(B(fP)))
\vspace{1mm}\\
\qquad 
   = f^{-1}(Af\ov{(f^{-1}(B(fP)})) = f^{-1}(Af(f^{-1}(\ov{(\mbox{Ad}J_3)B)}(f\ov{P}))) 
\vspace{1mm}\\
\qquad
  = f^{-1}((A(\epsilon B))(f\ov{P})) = \psi(A(\epsilon B), \epsilon)P.
\end{array}$$
Thus, for $A \in SU(6)$, we have $\psi(A, \epsilon) = \varphi(A, 1)\varphi(E, \epsilon) \in E_{7,{\sC}}$. Since $\psi$ induces a surjection $\psi_* : \su(6) \to \gge_{7,{\sC}}$, $\psi : SU(6) \to (E_{7,{\sC}})_0$ (which is the connected component of $E_{7,{\sC}}$ containing the identity 1) is onto. However $\epsilon = \psi(E, \epsilon) \not\in (E_{7,{\sC}})_0$. Indeed, for any $A \in SU(6)$,
$$
    \dsum(A\a \land A\b \land A\c) = \ov{\a} \land \ov{\b} \land \ov{\c}, \quad \a, \b, \c \in \C^6 $$
does not hold. Therefore $E_{7,{\sC}}$ has just two connected components (see Lemma 4.13.1). Hence $\psi : SU(6)\cdot\Z_2 \to E_{7,{\sC}}$ is onto. $\Ker\,\psi = \{E, \omega_1E, {\omega_1}^2E \} \times 1 = \Z_3 \times 1$ easily obtained. Thus we have the isomorphism $(SU(6)/\Z_3)\cdot\Z_2 \cong E_{7,{\sC}}$.
\vspace{2mm}
    
We identify $(\gP_{\sC})^C \oplus (M(3, \C)^C \oplus M(3, \C)^C)$ with $\gP^C$  (using the identification 
 ${\gJ_{\sC}}^C \oplus M(3, \C)^C$ with $\gJ^C$ in Section 3.13) by
$$
   ((X, Y, \xi, \eta), (M, N)) = (X + M, Y + N, \xi, \eta). $$
Further, we define a $\C$-$C$-linear mapping $\mu : M(6, \C) \to M(3, \C)^C \oplus M(3, \C)^C$ by
$$
\begin{array}{l} 
   \mu\pmatrix{M_{11} & M_{12} \cr
               M_{21} & M_{22}} 
\vspace{1mm}\\
\qquad 
    = \Big(\dfrac{(M_{21} - M_{12})e_1}{2} + i\dfrac{M_{21} + M_{12}}{2},
           \dfrac{(M_{22} + M_{11})e_1}{2} + i\dfrac{M_{22} - M_{11}}{2}\Big), 
\end{array}$$
where $M_{ij} \in M(3, \C)$. The inverse mapping $\mu^{-1} : M(3, \C)^C \oplus M(3, \C)^C \to M(6, \C)$ of $\mu$ is given by
$$
   \mu^{-1}(M_1 + iM_2, N_1 + iN_2) 
   = \pmatrix{- N_2 - N_1e_1 & M_2 + M_1e_1 \cr
              M_2 - M_1e_1 & N_2 - N_1e_1}, \quad M_i, N_i \in M(3, \C). 
\vspace{3mm}
$$

{\bf Lemma 4.13.4.} {\it For $D \in \su(6)$ and $\wti{M} \in M(6, \C)$, we have}
$$
           \mu(\wti{M}D^*) = \psi_*(D)(\mu\wti{M}). 
\vspace{2mm}$$

{\bf Proof.} Let 
$$
\begin{array}{l}
       D = \pmatrix{B & L \cr
                    - L^* & C} + \dfrac{\nu}{3}\pmatrix{E & 0 \cr
                                                        0 & - E} \in \su(6), 
\vspace{1mm}\\
      \wti{M} = \pmatrix{- N_2 - N_1e_1 & M_2 + M_1e_1 \cr
                                   M_2 - M_1e_1 & N_2 - N_1e_1}, 
    \quad M_i, N_i \in M(3, \C), 
\vspace{1mm}\\
    M = M_1 + iM_2, N = N_1 + iN_2.
\end{array}$$
Then we have
\vspace{2mm}

  $\psi_*(D)(\mu\wti{M})$
\vspace{1mm}

\quad
  $= {\mit\Phi}(\phi_{\sC}(B, C), h(L), -\tau h(L), - i\nu e_1)(M, N)$
\vspace{1mm}

\quad
  $= \pmatrix{\phi_{\sC}(B, C) + \dfrac{1}{3}i\nu e_1 & - 2\tau h(L) & 0 & h(L) \cr
       2h(L) & \tau\phi_C(B, C)\tau - \dfrac{1}{3}i\nu e_1 & - \tau h(L) & 0 \cr
       0 & h(L) & -i\nu e_1 & 0 \cr
       -\tau h(L) & 0 & 0 & i\nu e_1}\pmatrix{M \vspace{1.5mm}\cr N \vspace{1.5mm}\cr 0 \vspace{1.5mm}\cr 0}$
\vspace{2mm}

\quad
   $= \pmatrix{\phi_{\sC}(B, C)M + \dfrac{1}{3}i\nu e_1M -2\tau h(L) \times N \cr
      2h(l) \times M + \tau\phi_{\sC}(B, C)\tau N - \dfrac{1}{3}i\nu e_1N  \cr
      (h(L), N) \cr
      - (\tau h(L), M)}$
\vspace{2mm}

\quad 
 $= \pmatrix{- Mh(B, C) + N\tau h(L) + \dfrac{1}{3}i\nu e_1M \vspace{0.5mm}\cr
               - Mh(L) - N\tau h(B, C) - \dfrac{1}{3}i\nu e_1N \cr
               0 \cr 
               0}$ 
\vspace{1mm}

\noindent (using $\phi_{\sC}(B, C)M = M\tau h(B, C)^* = - Mh(B, C)$ and $-2\tau h(L) \times N = N\tau h(L)$ etc.)
$$
\begin{array}{l}
\stackrel{\mu^{-1}}\longrightarrow \cdots \; \mbox{by simple calculations}\; \cdots
\vspace{1mm}\\
   = \pmatrix{- N_2 - N_1e_1 & M_2 + M_1e_1 \cr
                M_2 + M_1e_1 & N_2 + N_1e_1}\Big(\pmatrix{- B & - L \cr
                                                            L^* & - C} - 
         \dfrac{\nu}{3}\pmatrix{E & 0 \cr
                                0 & - E} \Big) = \wti{M}D^*.
\end{array}$$

{\bf Definition.} We define a $C$-$\C$-linear isomorphism $f : \gP^C \to {\mit\Lambda}^3(\C^6) \oplus M(6, \C)$ by
$$
\begin{array}{l}
    f(P_{\sC} + (M + N)) = f(P_{\sC}) + \mu^{-1}(M + N), 
\vspace{1mm}\\
   \qquad \qquad P_{\sC} + (M + N) \in (\gP_C)^C \oplus (M(3, \C)^C \oplus M(3, \C)^C) = \gP^C. 
\end{array} $$

The group $SU(3) \times SU(6)$ acts on ${\mit\Lambda}^3(\C^6) \oplus M(6, \C)$ by
$$
         (Q, A)(\dsum(\a \land \b \land \c) + \wti{M}) = \dsum(A\a \land A\b \land A\c) + Q\wti{M}A^*, $$
where $Q\wti{M}$ means $\pmatrix{Q & 0 \cr
                                 0 & Q}\pmatrix{M_{11} & M_{12} \cr
                                                M_{21} & M_{22}} 
    = \pmatrix{QM_{11} & QM_{12} \cr
               QM_{21} & QM_{22}}, M_{ij} \in M(3, \C)$.
\vspace{3mm}

{\bf Theorem 4.13.5.} $(E_7)^w \cong (SU(3) \times SU(6))/\Z_3, \; \Z_3 = \{(E, E), (\omega_1E, \omega_1E),$ $ ({\omega_1}^2E, {\omega_1}^2E) \}, \omega_1 = - \dfrac{1}{2} + \dfrac{\sqrt{3}}{2}e_1$.
\vspace{2mm}

{\bf Proof.} We defined a mapping $\psi : SU(3) \times SU(6) \to (E_7)^w$ by
$$
      \psi(Q, A)P = f^{-1}((Q, A)(fP)), \quad P\in \gP^C. $$
We first have to prove that $\psi(Q, A) \in E_7$. To prove this, since $\varphi(Q, E) \in (E_6)^w \subset (E_7)^w$, it suffices to show that $\psi(E, A) \in E_7$. Moreover, it is sufficient to show that, for the differential mapping $\psi_* : \su(3) \oplus \su(6) \to \gge_7$ of $\psi$, $\psi_*(0, D)$ coincides with
$$
    {\mit\Phi}(\phi_{\sC}(B, C), h(L),  -\tau h(L), - i\nu e_1) \in \gge_7. $$
However, this is already shown in Proposition 4.13.3 and Lemma 4.13.4. Since $(\gP^C)_w = \{P \in \gP^C \, | \, wP = P \} = (\gP_{\sC})^C$, obviously we have $w\psi(Q, A) = \psi(Q, A)w$, hence $\psi(Q, A) \in (E_7)^w$. Evidently $\psi$ is a homomorphism. We shall show that $\psi$ is onto. Let $\alpha \in (E_7)^w$. Since the restriction $\alpha'$ to $(\gP^C)_w = (\gP_{\sC})^C$ of $\alpha$ belongs to $E_{7,{\sC}}$, there exists $A \in SU(6)$ such that
$$
     \alpha P = f^{-1}(A(fP)) \quad \mbox{or} \quad 
     \alpha P = f^{-1}(A(f\ov{P})), \quad P \in (\gP_{\sC})^C $$
(Proposition 4.13.3). In the former case, let $\beta = \psi(E, A)^{-1}\alpha$, then $\beta|(\gP_{\sC})^C = 1$, hence $\beta \in G_2$. Furthermore, $\beta \in (G_2)^w = SU(3)$ (Theorem 1.9.4), so there exists $Q \in SU(3)$ such that
$$
\begin{array}{l}
     \beta(P_{\sC} + (M + N)) = P_{\sC} + Q(M + N) = P_{\sC} + (QM + QN)
\vspace{1mm}\\
\qquad \quad = \psi(Q, E)(P_{\sC} + (M + N)), \quad P_{\sC} + (M + N) \in \gP^C
\end{array} $$
Hence we have
$$
      \alpha = \psi(E, A)\beta = \psi(E, A)\psi(Q, E) = \psi(Q, A). $$
In this case, this shows that $\psi$ is onto. In the latter case, consider the mapping 
$$
     \gamma_1 : \gP^C \to \gP^C, \gamma_1(P_{\sC} + (M + N)) = \ov{P_{\sC}} + (\ov{M} +\ov{N}), \quad P_{\sC} + (M + N) \in \gP^C. $$
Then, $\gamma_1 \in G_2 \subset F_4 \subset E_6 \subset E_7$. From the same argument as Section 1, we have $\gamma_1\alpha \in (E_7)^w$, hence $\gamma_1 \in (G_2)^w = SU(3)$. However this is a contradiction (Theorem 1.9.4). Therefore that $\psi$ is onto is shown. $\Ker\psi = \{(E, E), (\omega_1E, \omega_1E),$ $ ({\omega_1}^2E, {\omega_1}^2E) \} = \Z_3$ is easily obtained. Thus we have the isomorphism $(SU(3) \times SU(6))/\Z_3 \cong (E_7)^w$.
\vspace{2mm}

{\bf Remark 1.} The group $E_7$ has a subgroup which is isomorphic to the semi-direct product $((SU(3) \times SU(6))/\Z_3) \cdot \Z_2$ (the action of the group $\Z_2 = \{1, \gamma \}$ to the group  $SU(3) \times SU(6)$ is $\gamma(Q, A) = (\ov{Q}, \ov{\mbox{Ad}(J_3)A}))$.

\vspace{2mm}

{\bf Remark 2.} Since $(E_7)^w$ is connected, the fact that $\psi : SU(3) \times SU(6) \to (E_7)^w$ is onto can be proved as follows. The elements
$$
\begin{array}{l}    
    G_{01}, \quad G_{23}, \quad G_{45}, \quad G_{67}, \quad G_{46} + G_{47}, \quad G_{47} - G_{56}, 
\vspace{1mm}\\ 
    G_{24} + G_{35}, \quad \; G_{25} - G_{34}, \quad \, G_{26} + G_{37}, \quad G_{27} - G_{36}, 
\vspace{1mm}\\ 
    \wti{A}_l(1), \quad  \wti{A}_l(e_1), \quad \wti{F}_l(1), \quad \wti{F}_l(e_1), \quad (E_1 - E_2)^{\sim}, \quad (E_2 - E_3)^{\sim}
\vspace{1mm}\\
    \check{F}_l(1), \quad \check{F}_l(e_1), \quad \, \hat{F}_l(1), \quad \hat{F}_l(e_1), \quad \check{E}_l, \quad \hat{E}_l, \quad \1, \quad l=1, 2, 3 
\end{array}$$
forms an $\R$-basis of $(\gge_7)^w$. So, $\dim(\gge_7)^w = 10 + 14 + 6 \times 3 + 1 = 43 = 8 + 35 = \dim(\su(3) \oplus \su(6))$. Hence $\varphi$ is onto.
\vspace{4mm}

{\bf 4.14. Complex exceptional Lie group ${E_7}^C$}
\vspace{3mm}

{\bf Theorem 4.14.1.} {\it The polar decomposition of the Lie group ${E_7}^C$ 
is given by}
$$
       {E_7}^C \simeq E_7 \times \R^{133}.$$
{\it In particular, ${E_7}^C$ is a simply connected complex Lie group of type $E_7$}.
\vspace{2mm}

{\bf Proof.} Evidently ${E_7}^C$ is an algebraic subgroup of $\Iso_C(\gP^C) = GL(78, C)$. If $\alpha \in {E_7}^C$, then, the complex conjugate transpose $\alpha^*$ with repect to the inner product $\langle X, Y \rangle$: $\langle \alpha X, Y \rangle = \langle X, \alpha^*Y \rangle$ is $\alpha^* = \tau\lambda\alpha^{-1}\lambda^{-1}\tau \in {E_7}^C$. Hence, from Chevalley's lemma, we have
$$
     {E_7}^C \simeq ({E_7}^C \cap U(\gP^C)) \times \R^d = E_7 \times \R^d, \quad d = 133. $$
Since $E_7$ is simply connected (Theorem 4.9.2), ${E_7}^C$ is also simply connected. The Lie algebra of the group ${E_7}^C$ is ${\gge_7}^C$, so ${E_7}^C$ is a complex simple Lie group of type $E_7$.
\vspace{4mm}

{\bf 4.15. Non-compact exceptional Lie groups $E_{7(7)}, E_{7(-5)}$ and $E_{7(-25)}$ of type 
\vspace{3mm}
$E_7$}

Let
\begin{eqnarray*}
   \gP \!\!\! &=& \!\!\! \gJ(3, \gC) \oplus \gJ(3, \gC) \oplus \R \oplus \R, 
\vspace{1mm}\\
   \gP' \!\!\! &=& \!\!\! \gJ(3, \gC') \oplus \gJ(3, \gC') \oplus \R \oplus \R.
\end{eqnarray*}
For $P, Q \in \gP$ or $\gP'$, we define an $\R$-linear mapping $P \times Q : \gP  \to \gP$ or $\gP' \to \gP'$ as similar to Section 4.1. And we define a Hermitian inner product $\langle P, Q \rangle_{\sigma}$ in $\gP^C$ by
$$
       \langle P, Q \rangle_{\sigma} = \langle \sigma P, Q \rangle. $$
Now, we define groups $E_{7(7)}, E_{7(-5)}$ and $E_{7(-25)}$ by
\begin{eqnarray*}
    E_{7(7)} \!\!\! &=& \!\!\! \{\alpha \in \Iso_{\sR}(\gP') \, | \, \alpha(P \times Q)\alpha^{-1} = \alpha P \times \alpha Q \}, 
\vspace{1mm}\\
    E_{7(-5)} \!\!\! &=& \!\!\! \{\alpha \in \Iso_C(\gP^C) \, | \, \alpha(P \times Q)\alpha^{-1} = \alpha P \times \alpha Q, \langle \alpha P, \alpha Q \rangle_{\sigma} = \langle P, Q \rangle_{\sigma} \}, 
\vspace{1mm}\\
    E_{7(-25)} \!\!\! &=& \!\!\! \{\alpha \in \Iso_{\sR}(\gP) \, | \, \alpha(P \times Q)\alpha^{-1} = \alpha P \times \alpha Q \}.   
\end{eqnarray*}
These groups can also be defined by
$$
         E_{7(7)} \cong ({E_7}^C)^{\tau\gamma}, \quad 
         E_{7(-5)} \cong ({E_7}^C)^{\tau\lambda\sigma}, \quad 
         E_{7(-25)} \cong ({E_7}^C)^{\tau}. $$

{\bf Theorem 4.15.1.} {\it The polar decompositions of the groups $E_{7(7)}, E_{7(-5)}$ and $E_{7(-25)}$ are respectively given by}
\begin{eqnarray*}
    E_{7(7)} \!\!\! &\simeq& \!\!\! SU(8)/\Z_2 \times \R^{70}, 
\vspace{1mm}\\
 E_{7(-5)} \!\!\! &\simeq& \!\!\! (SU(2) \times Spin(12))/\Z_2 \times \R^{64}, 
\vspace{1mm}\\
    E_{7(-25)} \!\!\! &\simeq& \!\!\! (U(1) \times E_6)/\Z_3 \times \R^{54}. 
\end{eqnarray*}

{\bf Proof.} These are the facts corresponding to Theorems 4.12.5, 4.11.15 and 
\vspace{3mm}
4.10.2.

{\bf Theorem 4.15.2.} {\it The centers of the groups $E_{7(7)}, E_{7(-5)}$ and $E_{7(-25)}$ are the group of order} 2:
$$
      z(E_{7(7)}) = \Z_2, \quad z(E_{7(-5)}) = \Z_2, \quad  z(E_{7(-25)}) = \Z_2. $$

\newpage

\vspace{5mm}

\begin{center}
\large{\bf 5. Exceptional Lie group $E_8$}
\end{center}
\vspace{4mm} 

{\bf 5.1. Lie algebra ${\gge_8}^C$}
\vspace{3mm}

{\bf Theorem 5.1.1.} {\it In a} $133 + 56 \times 2 + 3 = 248$ {\it dimensional $C$-vector space}
$$
     {\gge_8}^C = {\gge_7}^C \oplus \gP^C \oplus \gP^C \oplus C \oplus C \oplus C, $$
{\it if we define a Lie bracket $[R_1, R_2]$ by
$$
   [({\mit\Phi}_1, P_1, Q_1, r_1, s_1, t_1), ({\mit\Phi}_2, P_2, Q_2, r_2, s_2, t_2)] = ({\mit\Phi}, P, Q, r, s, t), $$ 
where}
$$
\left\{\begin{array}{l}
   {\mit\Phi} = [{\mit\Phi}_1, {\mit\Phi}_2] + P_1 \times Q_2 - P_2 \times Q_1 
\vspace{1mm}\\
   P = {\mit\Phi}_1P_2 - {\mit\Phi}_2P_1 + r_1P_2 - r_2P_1 + s_1Q_2 - s_2Q_1
\vspace{1mm}\\              
   Q = {\mit\Phi}_1Q_2 - {\mit\Phi}_2Q_1 - r_1Q_2 + r_2Q_1 + t_1P_2 - t_2P_1
\vspace{1mm}\\              
 r = - \dfrac{1}{8}\{ P_1, Q_2 \} + \dfrac{1}{8}\{ P_2, Q_1 \} + s_1t_2 - s_2t_1
\vspace{1mm}\\
   s = \;\;\,\dfrac{1}{4}\{ P_1, P_2 \} + 2r_1s_2 - 2r_2s_1
\vspace{1mm}\\
   t = - \dfrac{1}{4}\{ Q_1, Q_2 \} - 2r_1t_2 + 2r_2t_1,
\end{array} \right. $$
{\it then ${\gge_8}^C$ is a $C$-Lie algebra.}  
\vspace{2mm}

{\bf Proof.} Among the definition of the Lie algebra, the relations
$$
\begin{array}{l}
        [R_1, R_2 + R_3] = [R_1, R_2] + [R_1, R_3],
\vspace{1mm}\\
        {[}kR_1, R_2] = k[R_1, R_2], \quad k \in C,
\vspace{1mm}\\
        {[}R_1, R_2] = - [R_2, R_1]
\end{array}$$
are evident, and we are left to show the Jacobi identity, which can be proved by direct calculations as follows.
$$
\begin{array}{l}
    [R_1, [R_2, R_3]\,] + [R_2, [R_3, R_1]\,] + [R_3, [R_1, R_2]\,] 
\vspace{1mm}\\
    \quad = \cdots (\mbox{using}\; [{\mit\Phi}, P \times Q] = {\mit\Phi}P \times Q + P \times {\mit\Phi}Q \; \mbox{(Proposition 4.3.2)}, 
\vspace{1mm}\\
    \qquad (P \times R)Q - (Q \times R)P + \dfrac{1}{8}\{Q, R\}P - \dfrac{1}{8}\{P, R\} - \frac{1}{4}\{P, Q\}R  = 0 
\vspace{1mm}\\
    \qquad \mbox{(Lemma 4.1.1.(3)}, \{{\mit\Phi}P, Q\} + \{P, {\mit\Phi}Q\} = 0 \; \mbox{(Proposition 4.2.2.(2)) etc.)} \cdots
\vspace{1mm}\\
    \quad = 0.
\end{array}$$

\vspace{4mm}

{\bf 5.2. Simplicity of ${\gge_8}^C$}
\vspace{3mm}

We use the following notation in ${\gge_8}^C$:
$$
\begin{array}{cc}
       {\mit\Phi} = ({\mit\Phi}, 0, 0, 0, 0, 0),
     & P^- = (0, P, 0, 0, 0, 0), \,\;
\vspace{1mm}\\
       Q_- = (0, 0, Q, 0, 0, 0), \;
     & \;\,r = (0, 0, 0, r, 0, 0),
\vspace{1mm}\\
       s^- = (0, 0, 0, 0, s, 0), \;
      & t_- = (0, 0, 0, 0, 0, t).
\end{array}$$

{\bf Theorem 5.2.1.} {\it The $C$-Lie algebra ${\gge_8}^C$ is simple.}
\vspace{2mm}

{\bf Proof.} We use the decomposition of ${\gge_8}^C$:
$$
      {\gge_8}^C = {\gge_7}^C \oplus \gK^C, $$
where $\gK^C = \gP^C \oplus \gP^C \oplus C \oplus C \oplus C$. Let $p : {\gge_8}^C \to {\gge_7}^C$ and $q : {\gge_8}^C \to \gK^C$ be projections of ${\gge_8}^C = {\gge_7}^C \oplus \gK^C$.  Now, let $\ga$ be a non-zero ideal of ${\gge_8}^C$. Then $p(\ga)$ is an ideal of ${\gge_7}^C$. Indeed, if ${\mit\Phi} \in p(\ga)$, then there exists $(0, P, Q, r, s, t) \in \gK^C$ such that $({\mit\Phi}, P, Q, r, s, t) \in \ga$. For any ${\mit\Phi}_1 \in {\gge_7}^C$, we have 
$$
   \ga \ni [{\mit\Phi}_1, ({\mit\Phi}, P, Q, r, s, t)] = ([{\mit\Phi}_1, {\mit\Phi}], {\mit\Phi}_1P, {\mit\Phi}_1Q, 0, 0, 0), $$
hence $[{\mit\Phi}_1, {\mit\Phi}] \in p(\ga)$.
\vspace{1mm}

We shall show that either ${\gge_7}^C \cap \ga \neq \{0\}$ or $\gK^C \cap \ga \neq \{0\}$. Assume that ${\gge_7}^C \cap \ga = \{0\}$ and $\gK^C \cap \ga = \{0 \}$. Then the mapping $p|\ga : \ga \to {\gge_7}^C$ is injective because $\gK^C \cap \ga = \{0\}$.  Since $p(\ga)$ is a non-zero ideal of ${\gge_7}^C$ and ${\gge_7}^C$ is simple, we have $p(\ga) = {\gge_7}^C$. Hence $\dim_C(\ga) = \dim_C(p(\ga)) = \dim_C({\gge_7}^C) = 133$. On the other hand, since ${\gge_7}^C \cap \ga = \{0\}$, $q|\ga : \ga \to \gK^C$ is also injective. Hence we have $\dim_C(\ga) \le \dim_C(\gK^C) = 56 \times 2 + 3 = 115$. This leads to a contradiction. 
\vspace{1mm}

We now consider the following two cases.
\vspace{1mm}

(1) Case ${\gge_7}^C \cap \ga \neq \{0\}$. From the simplicity of ${\gge_7}^C$, we have ${\gge_7}^C \cap \ga = {\gge_7}^C$, hence $\ga \supset {\gge_7}^C$. On the other hand, we have
$$
\begin{array}{l}
   \ga \ni [{\mit\Phi}(0, 0, 0, 1), (0, 0, 1, 0)^-] = (0, 0, 1, 0)^-,
\vspace{1mm}\\
   \ga \ni [{\mit\Phi}(0, 0, 0, 1), (0, 0, 0, -1)_-] = (0, 0, 0, 1)_-,
\vspace{1mm}\\
   \ga \ni [(0, 0, 1, 0)^-, (0, 0, 0, 4)^-] = \1^-,
\vspace{1mm}\\
   \ga \ni [(0, 0, 0, 1)_-, (0, 0, 4, 0)_-] = \1_-,
\vspace{1mm}\\
   \ga \ni [\1^-, \1_-] = \1,
\vspace{1mm}\\
   \ga \ni [\1^- + \1_-, Q^- + P_-] = P^- + Q_-,
\end{array}$$
Therefore, $\ga \supset {\gge_7}^C \oplus \gK^C = {\gge_8}^C$ which implies $\ga = {\gge_8}^C$.
\vspace{1mm}

(2) Case $\gK^C \cap \ga \neq \{ 0 \}$. Let $R = (0, P, Q, r, s, t)$ be a non-zero element of $\gK^C \cap \ga \subset 
\vspace{1mm}
\ga$.

$\;$(i) Case $R = (0, P, Q, r, s, t), P \neq 0$. We have
$$
\begin{array}{l}
     \ga \ni [\1,[\1_-,[\1, R] \,] \,] = [\1, [\1_-, (0, P, -Q, 0, 2s, -2t)]\,]
\vspace{1mm}\\
     \;\;\; = [\1, (0, 0, P, -2s, 0, 0)] = - (0, 0, P, 0, 0, 0) = - P_-.
\end{array}$$
We choose $P_1 \in \gP^C$ so that $P \times P_1 \neq 0$ (Lemma 4.5.3) and choose ${\mit\Phi} \in {\gge_7}^C$ so that $[{\mit\Phi}, P \times P_1] \neq 0$. (Since ${\gge_7}^C$ simple, the center of ${\gge_7}^C$ consists only of $0$, so such ${\mit\Phi}$ exists). Then we have
$$
      \ga \ni [{\mit\Phi}, [{P_1}^-, P_-]\,] = [{\mit\Phi}, P \times P_1]. $$
Hence this case is reduced to the case (1).
\vspace{1mm}

(ii) Case $R = (0, P, Q, r, s, t), Q \neq 0$. The argument is similar to (i).
\vspace{1mm}

(iii) Case $R = (0, 0, 0, r, s, t), r \neq 0$. For $0 \neq P \in \gP^C$, we have
$$
\begin{array}{l}
  \ga \ni [P^-, [\1^-, [\1_-, R]\,]\,] = [P^-, [\1^-, (0, 0, 0, -s, 0, 2r)]\,]
\vspace{1mm}\\
     \;\;\; = [P^-, (0, 0, 0, 2r, 2s, 0)] = (0, -2rP, 0, 0, 0, 0).
\end{array}$$
Hence this case is reduced to the case (ii) above.
\vspace{1mm}

(iv) Case $R = (0, 0, 0, 0, s, t), s \neq 0$. We have
$$
       \ga \ni [\1_-, R] = (0, 0, 0, -s, 0, 0). $$
Hence this case is reduced to the case (iii) above.
\vspace{1mm}

(v) Case $R= (0, 0, 0, 0, 0, t), t \neq 0$. The argument is similar to (iv).
\vspace{1mm}

\noindent Consequently, we have $\ga = {\gge_8}^C$, which proves the simplicity of ${\gge_8}^C$. 
\vspace{2mm}

For $R \in {\gge_8}^C$ we denote ad$R : {\gge_8}^C \to {\gge_8}^C$, that is, $\mbox{ad}R(R_1) = [R, R_1]$, by ${\mit\Theta}(R) = \mbox{ad}R$. Since ${\gge_8}^C$ is simple (Theorem 5.2.1), we obtain an isomorphism of Lie algebras
$$
     {\gge_8}^C \cong {\mit\Theta}({\gge_8}^C) = \{{\mit\Theta}(R) \, | \, R \in {\gge_8}^C \} $$
by assigning ${\mit\Theta}(R)$ to $R$. Moreover, ${\gge_8}^C$ is isomorphic to the algebra
$$
  \mbox{Der}({\gge_8}^C) = \{ {\mit\Theta} \in \Hom_C({\gge_8}^C) \, | \, {\mit\Theta}[R_1, R_2] = [{\mit\Theta}R_1, R_2] + [R_2, {\mit\Theta}R_2] \}. $$
Hereafter we often denote ${\mit\Theta}(R)$ by $R$ identifying ${\gge_8}^C \cong \mbox{Der}({\gge_8}^C)$. 
\vspace{4mm}

{\bf 5.3. Killing form of ${\gge_8}^C$}
\vspace{3mm}

{\bf Definition.} We define a symmetric inner product $(R_1, R_2)_8$ in ${\gge_8}^C$ by
$$
   (R_1, R_2)_8 = ({\mit\Phi}_1, {\mit\Phi}_2)_7 - \{ Q_1, P_2 \} + \{ P_1, Q_2 \} - 8r_1r_2 - 4t_1s_2 - s_1t_2, $$
where $R_i = ({\mit\Phi}_i, P_i, Q_i, r_i, s_i, t_i) \in {\gge_8}^C$.
\vspace{3mm}

{\bf Lemma 5.3.1.} {\it The inner product $(R_1, R_2)_8$ of ${\gge_8}^C$ is ${\gge_8}^C$-adjoint invariant}\,:
$$
     ([R, R_1], R_2)_8 + (R_1, [R, R_2])_8 = 0, \quad R, R_i \in {\gge_8}^C.$$

{\bf Proof.} $([R, R_1], R_2)_8$

$$
\begin{array}{l}
  = \Biggl(\pmatrix{[{\mit\Phi}, {\mit\Phi}_1] + P \times Q_1 - P_1 \times Q 
\vspace{1mm}\cr
    {\mit\Phi}P_1 - {\mit\Phi}_1P + rP_1 - r_1P + sQ_1 - s_1Q
\vspace{1mm}\cr              
    {\mit\Phi}Q_1 - {\mit\Phi}_1Q - rQ_1 + r_1Q + tP_1 - t_1P
\vspace{1mm}\cr              
     - \dfrac{1}{8}\{ P, Q_1 \} + \dfrac{1}{8}\{ P_1, Q \} + st_1 - s_1t
\vspace{1mm}\cr
      \dfrac{1}{4}\{ P, P_1 \} + 2rs_1 - 2r_1s
\vspace{1mm}\cr
     - \dfrac{1}{4}\{ Q, Q_1 \} - 2rt_1 + 2r_1t},
    \pmatrix{{\mit\Phi}_2 
\vspace{2mm}\cr
              P_2
\vspace{2mm}\cr
              Q_2
\vspace{2mm}\cr
              r_2
\vspace{2mm}\cr
              s_2
\vspace{2mm}\cr
              t_2} \Biggl)_8 
\vspace{1mm}\\
    \quad = \cdots (\mbox{using}\; ([{\mit\Phi}, {\mit\Phi}_1], {\mit\Phi}_2)_7 + ({\mit\Phi}_1, [{\mit\Phi}, {\mit\Phi}_2])_7 = 0 \; \mbox{(Lemma 4.5.1.(1))}, \vspace{1mm}\\
   \quad \qquad({\mit\Phi}, P \times Q)_7 = \{{\mit\Phi}P, Q \} \; \mbox{(Lemma 4.5.1.(2)) etc.)} \cdots 
\vspace{1mm}\\
    \quad = - (R_1, [R, R_2])_8.
\end{array}$$

{\bf Theorem 5.3.2.} {\it The Killing form $B_8$ of the Lie algebra ${\gge_8}^C$  is given by}
$$
\begin{array}{l}
  B_8(R_1, R_2)
\vspace{1mm}\\ 
  = - 15(R_1, R_2)_8
\vspace{1mm}\\
  = - 15({\mit\Phi}_1, {\mit\Phi}_2)_7 + 15\{ Q_1, P_2 \} - 15\{ P_1, Q_2 \} +120r_1r_2 + 60t_1s_2 + 60s_1t_2
\vspace{1mm}\\ 
   = \dfrac{5}{3}B_7({\mit\Phi}_1, {\mit\Phi}_2) + 15\{ Q_1, P_2 \} - 15\{ P_1, Q_2 \} +120r_1r_2 + 60t_1s_2 + 60s_1t_2,
\end{array}$$
{\it where} $R_i = ({\mit\Phi}_i, P_i, Q_i, r_i, s_i, t_i) \in {\gge_8}^C$ {\it and $B_7$ is the Killing form of ${\gge_7}^C$}.
\vspace{2mm}

{\bf Proof.} Since ${\gge_8}^C$ is simple (Theorem 5.2.1), there exist $k \in C$ such that
$$
    B_8(R_1, R_2) = k(R_1, R_2)_8, \quad R_i \in {\gge_8}^C. $$
To determine $k$, let $R_1 = R_2 = (0, 0, 0, 1, 0, 0) = \1$. Then, we have
$$
        (\1, \1)_8 = - 8.$$
On the other hand, since
$$
   [\1, [\1, ({\mit\Phi}, P, Q, r, s, t)]\,] = [\1, (0, P, - Q, 0, 2s, - 2t)] = (0, P, Q, 0, 4s, 4t), $$
we have 
$$
      B_8(\1, \1) = 56 \times 2 + 4 \times 2 = 120. $$
Therefore $k = - 15$. Thus we have $B_8(R_1, R_2) = -15(R_1, R_2)_8$.
\vspace{4mm}

{\bf 5.4. Complex exceptional Lie group ${E_8}^C$}
\vspace{3mm}
                                                                             
{\bf Definition.} The group ${E_8}^C$ is defined to be the automorphism group of the Lie algebra ${\gge_8}^C$:
$$
    {E_8}^C = \{\alpha \in \Iso_C({\gge_8}^C) \, |\, \alpha[R_1, R_2] = [\alpha R_1, \alpha R_2] \}. 
\vspace{2mm}
$$

{\bf Theorem 5.4.1.} {\it The group ${E_8}^C$ is connected.}
\vspace{2mm}

{\bf Proof.} Denote by Inn(${\gge_8}^C$) the subgroup generated by the inner automorphisms \\ $\exp({\mit\Theta}(R))$, $R \in {\gge_8}^C$ in the automorphism group $\mbox{Aut}({\gge_8}^C) = {E_8}^C$ of the Lie algebra ${\gge_8}^C$. It is known that $\mbox{Aut}({\gge_8}^C)/\mbox{Inn}({\gge_8}^C) = \{1\}$ holds for the $C$-algebra of $E_8$ type (see for example, Matsushima [19]), that is,
$$
       \mbox{Aut}({\gge_8}^C) = \mbox{Inn}({\gge_8}^C). $$
Since $\mbox{Inn}({\gge_8}^C)$ is connected, ${E_8}^C$ which is equal to Inn$({\gge_8}^C)$ is also connected. 
\vspace{2mm}

{\bf Remark.} We can also prove the connectedness of the group ${E_8}^C$ from the following fact. For $R \in {\gge_8}^C$, we define a $C$-linear mapping $R \times R : {\gge_8}^C \to {\gge_8}^C$ by
$$
     (R \times R)R_1 = {\mit\Theta}(R)^2R_1 + \dfrac{1}{30}B_8(R, R_1)R, \quad R_1 \in {\gge_8}^C, $$
and we define a space $\gW^C$ by
$$
   \gW^C = \{ R \in {\gge_8}^C \, | \, R \times R = 0, R \neq 0 \}. $$
Then we have
$$
         {E_8}^C/({E_8}^C)_{1_-} \simeq \gW^C, $$
where $({E_8}^C)_{1_-} = \{ \alpha \in {E_8}^C \, | \, \alpha 1_- = 1_- \} = \exp({\mit\Phi}(0, 0, \gP^C, 0, 0, C))E_7$. The connectedness of ${E_8}^C$ follows from the connectedness of $({E_8}^C)_{1_-}$ and $\gW^C$. (See Imai and Yokota [13]).
\vspace{2mm}

The Lie algebra of the group ${E_8}^C$ is Der(${\gge_8}^C$) $\cong {\gge_8}^C$, and therefore we have shown that ${E_8}^C$ is a complex Lie group of type $E_8$, since we will show in Theorem 5.6.2 that ${\gge_8}^C$ is a Lie algebra of type $E_8$. It is known by the general theory of Lie groups, that if a complex Lie group of type $E_8$ is connected, then it is simply connected, and hence we have obtained the following result.
\vspace{3mm}

{\bf Theorem 5.4.2.} {\it ${E_8}^C$ is a simply connected complex Lie group of type $E_8$}.
\vspace{4mm}

{\bf 5.5. Compact exceptional Lie group $E_8$}
\vspace{3mm}

We define $C$-linear transformations $\lambda, \lambda'$ of ${\gge_8}^C$ respectively by
$$
\begin{array}{l}
      \lambda({\mit\Phi}, P, Q, r, s, t) = (\lambda{\mit\Phi}\lambda^{-1}, \lambda P, \lambda Q, r, s, t),
\vspace{1mm}\\
      \lambda'({\mit\Phi}, P, Q, r, s, t) = ({\mit\Phi}, Q, - P, -r, -t, -s),
\end{array}$$
where $\lambda$ in the right hand side is the same as $\lambda \in E_7$ defined in Section 4.3. The mappings $\lambda$ and $\lambda'$ preserve the Lie bracket in ${\gge_8}^C$, that is, $\lambda, \lambda' \in \mbox{Aut}({\gge_8}^C) = {E_8}^C$. We set
$$
        \wti{\lambda} = \lambda\lambda' = \lambda'\lambda. $$
Finally, we denote by $\tau$ the complex cojugation in ${\gge_8}^C$, that is,
$$
      \tau({\mit\Phi}, P, Q, r, s, t) = (\tau{\mit\Phi}\tau, \tau P, \tau Q, \tau r, \tau s, \tau t),$$ 
where $\tau$ in the right hand side is the usual complex cojugation in the 
\vspace{3mm}
complexification.

{\bf Definition.} We define a Hermitian inner product $\langle R_1, R_2 \rangle$ in ${\gge_8}^C$ by
$$
      \langle R_1, R_2 \rangle = - \dfrac{1}{15}B_8(\tau\wti{\lambda}R_1, R_2). $$

{\bf Proposition 5.5.1.} {\it The Hermitian inner product $\langle R_1, R_2 \rangle$ in ${\gge_8}^C$ is positive definite.}
\vspace{2mm}

{\bf Proof.} Let $R_i = ({\mit\Phi}_i, P_i, Q_i, r_i, s_i, t_i) \in {\gge_8}^C, i = 1, 2$. Since
$
   \tau\wti{\lambda}R_1 = (\tau\lambda{\mit\Phi}_1\lambda^{-1}\tau, $ $\tau\lambda Q_1, - \tau\lambda P_1, - \tau r_1, - \tau t_1, - \tau s_1), $
we have, by Theorem 5.3.2, that
$$
\begin{array}{l}
   \langle R_1, R_2 \rangle 
\vspace{1mm}\\
  = (\tau\lambda{\mit\Phi}_1\lambda^{-1}\tau, {\mit\Phi}_2)_7 + \langle P_1, P_2 \rangle + \langle Q_1, Q_2 \rangle + 8(\tau r_1)r_2 + 4(\tau s_1)s_2 + 4(\tau t_1)t_2.
\end{array} $$
Hence, it is sufficient to show that $(\tau\lambda{\mit\Phi}_1\lambda^{-1}\tau, {\mit\Phi}_2)_7$ is positive definite. Let
$
     {\mit\Phi}_i = {\mit\Phi}(\phi_i, A_i, B_i, \nu_i), \,i = 1, 2. $
Since 
$
        \tau\lambda{\mit\Phi}_1\lambda^{-1}\tau = {\mit\Phi}(- \tau\,{}^t\phi_1\tau, - \tau B_1, - \tau A_1, - \tau\nu_1),$
we have
$$
     (\tau\lambda{\mit\Phi}_1\lambda^{-1}\tau, {\mit\Phi}_2)_7 = 2(\tau\,{}^t\phi_1\tau, \phi_2)_6 + 4\langle A_1, A_2 \rangle + 4\langle B_1, B_2 \rangle + \dfrac{8}{3}(\tau\nu_1)\nu_2. $$
Therefore, it is enough to show that $(\tau{}^t\phi_1\tau, \phi_2)_6$ is positive definite. Let $\phi_i = \delta_i + \wti{T}_i \in {\gge_6}^C, \delta_i \in {\gf_4}^C, \wti{T}_i \in {\gJ_0}^C, i = 1, 2$. Since $\tau{}^t\phi_1\tau = - \tau\delta_1\tau + \tau\wti{T}_1$, we have
$$
      (\tau\,{}^t\phi_1\tau, \phi_2)_6 = - (\tau\delta_1\tau, \delta_2)_4 + \langle T_1, T_2 \rangle. $$
Consequently, it is sufficient to show that $-(\tau\delta_1\tau, \delta_2)_4$ is positive definite, which can be seen, however, from the fact that the following set
$$
     \sqrt{2}[\wti{E}_1, \wti{F}_2(e_i)], \quad
     \sqrt{2}[\wti{E}_1, \wti{F}_3(e_i)], \quad
     \sqrt{2}[\wti{E}_3, \wti{F}_1(e_i)], \quad  0 \le i \le 7, $$
$$
      \dfrac{1}{\sqrt{2}}[\wti{F}_1(e_i), \wti{F}_1(e_j)], \quad 0 \le i < j \le 7 $$
forms an orthonormal $C$-basis of ${\gf_4}^C$ with respect to the inner product $-(\tau\delta_1\tau, \delta_2)_4$. Thus the proposition is proved.
\vspace{2mm}

{\bf Definition.} We define a group $E_8$ by
\begin{eqnarray*}
        E_8 \!\!\!&=&\!\!\! \{ \alpha \in {E_8}^C \, | \, \langle \alpha R_1, \alpha R_2 \rangle = \langle R_1, R_2 \rangle \}
\vspace{1mm}\\
      \!\!\!&=&\!\!\! \{ \alpha \in {E_8}^C \, | \, \tau\wti{\lambda}\alpha = \alpha\tau\wti{\lambda}\}.
\end{eqnarray*}

{\bf Theorem 5.5.2.} {\it The group $E_8$ is a compact Lie group.}
\vspace{2mm}

{\bf Proof.} $E_8$ is a compact Lie group as a closed subgroup of the unitary group 
$$
     U(248) = U({\gge_8}^C) = \{\alpha \in \Iso_C({\gge_8}^C) \, | \, \langle \alpha R_1, \alpha R_2 \rangle = \langle R_1, R_2 \rangle \}.$$

{\bf Theorem 5.5.3.} {\it The Lie algebra $\gge_8$ of the group $E_8$ is} 
$$
\begin{array}{l}
   \gge_8 = \{R \in {\gge_8}^C \, | \, \tau\wti{\lambda}R = R\}
\vspace{1mm}\\
\quad
     = \{({\mit\Phi}, P, - \tau\lambda P, r, s, - \tau s) \in {\gge_8}^C \, | \, {\mit\Phi} \in \gge_7, P \in \gP^C, r \in i\R, s \in C \}. 
\end{array}$$ 

{\bf Proof.} For $R = ({\mit\Phi}, P, Q, r, s, t) \in {\gge_8}^C$, since
$$  
      \tau\wti{\lambda}R = (\tau\lambda{\mit\Phi}\lambda^{-1}\tau, \tau\lambda Q, - \tau\lambda P, - \tau r, - \tau t, - \tau s). $$
the condition $\tau\wti{\lambda}R = R$ is equivalent to $\tau\lambda{\mit\Phi} = {\mit\Phi}\lambda\tau, Q = -\tau\lambda P, \tau r = - r, t = - \tau s$, hence we have the theorem.
\vspace{3mm}

{\bf Proposition 5.5.4.} {\it The complexification of the Lie algebra $\gge_8$ is ${\gge_8}^C$. Hence $\gge_8$ is simple.}
\vspace{2mm}

{\bf Proof.} For $R \in {\gge_8}^C$, the conjugate transposed mapping $R^*$ of $R$ with respect to the inner product $\langle R_1, R_2 \rangle$ of ${\gge_8}^C$ is $R^* = \tau\wti{\lambda}R\wti{\lambda}\tau \in {\gge_8}^C$, and for $R \in {\gge_8}^C$, $R$ belongs to $\gge_8$ if and only if $R^* = - R$. Now, any element $R \in {\gge_8}^C$ is represented by
$$
   R = \dfrac{R - R^*}{2} + i\dfrac{R + R^*}{2i},
\quad \dfrac{R - R^*}{2}, \dfrac{R + R^*}{2i} \in \gge_8. $$
Hence ${\gge_8}^C$ is the complexification of $\gge_8$. Since ${\gge_8}^C$ is simple (Theorem 5.2.1), $\gge_8$ is also simple.
\vspace{2mm}

Analogously as in ${\gge_8}^C$, for $R \in \gge_8$, we identify $R$ with ${\mit\Theta}(R)$ and regard $\gge_8 \cong 
\vspace{3mm}
{\mit\Theta}(\gge_8)$.

{\bf Theorem 5.5.5.} {\it The polar decomposition of the Lie group ${E_8}^C$ is given by}
$$
     {E_8}^C \simeq E_8 \times \R^{248}. $$
{\it In particular, the group $E_8$ is simply connected.}
\vspace{2mm}

{\bf Proof.} Evidently ${E_8}^C$ is an algebraic subgroup of $\Iso_C({\gge_8}^C) = GL(248, C)$. For $\alpha \in {E_8}^C$, the conjugate transposed mapping $\alpha^*$ of $\alpha$ with respect to the inner product $\langle R_1, R_2 \rangle$ is $\alpha^* = \tau\wti{\lambda}\alpha^{-1}\wti{\lambda}\tau \in {E_8}^C$. Therefore, by Chevalley's lemma, we have
$$
     {E_8}^C \simeq ({E_8}^C \cap U({\gge_8}^C)) \times \R^d = E_8 \times \R^d, \quad d = 248. $$
Since ${E_8}^C$ is simply connected (Theorem 5.4.2), $E_8$ is also simply connected.
\vspace{4mm}

{\bf 5.6  Roots of ${\gge_8}^C$}
\vspace{3mm}

{\bf Theorem 5.6.1.} {\it The rank of the Lie algebra ${\gge_8}^C$ is} 8. {\it The roots of ${\gge_8}^C$ relative to some Cartan subalgebra are given by}

$$
\begin{array}{cl}
     \pm(\lambda_k - \lambda_l), \quad \pm(\lambda_k + \lambda_l), & 0 \le k < l \le 3, 
\vspace{1mm}\\
     \pm \lambda_k \pm \dfrac{1}{2}(\mu_2 - \mu_3), & 0 \le k \le 3,
\end{array} $$
\vspace{-2mm}
$$
\begin{array}{l}
   \pm \dfrac{1}{2}(- \lambda_0 - \lambda_1 + \lambda_2 - \lambda_3) \pm \dfrac{1}{2}(\mu_3 - \mu_1), 
\vspace{1mm}\\
   \pm \dfrac{1}{2}(- \lambda_0 + \lambda_1 + \lambda_2 - \lambda_3) \pm \dfrac{1}{2}(\mu_3 - \mu_1),
\end{array}$$
$$
\begin{array}{l}
   \pm \dfrac{1}{2}(- \lambda_0 + \lambda_1 + \lambda_2 + \lambda_3) \pm \dfrac{1}{2}(\mu_3 - \mu_1), 
\vspace{1mm}\\
   \pm \dfrac{1}{2}(\;\;\; \lambda_0 - \lambda_1 + \lambda_2 + \lambda_3) \pm \dfrac{1}{2}(\mu_3 - \mu_1), 
\end{array}$$
\vspace{-2mm}
$$
\begin{array}{l}
   \pm \dfrac{1}{2}(\;\;\; \lambda_0 - \lambda_1 + \lambda_2 - \lambda_3) \pm \dfrac{1}{2}( \mu_1 - \mu_2), 
\vspace{1mm}\\
   \pm \dfrac{1}{2}(- \lambda_0 + \lambda_1 + \lambda_2 - \lambda_3) \pm \dfrac{1}{2}(\mu_1 - \mu_2),
\vspace{1mm}\\
   \pm \dfrac{1}{2}(\;\;\; \lambda_0 + \lambda_1 + \lambda_2 + \lambda_3) \pm \dfrac{1}{2}(\mu_1 - \mu_2), 
\vspace{1mm}\\
   \pm \dfrac{1}{2}(- \lambda_0 - \lambda_1 + \lambda_2 + \lambda_3) \pm \dfrac{1}{2}(\mu_1 - \mu_2),
\end{array} $$
$$
\begin{array}{cl}
     \pm \Big(\mu_k + \dfrac{2}{3}\nu \Big),  & 0 \le k < l \le 3, 
\vspace{1mm}\\
\pm \lambda_k \pm \Big(\dfrac{1}{2}\mu_1 - \dfrac{2}{3}\nu \Big), & 0 \le k \le 3,
\end{array} $$
\vspace{-1mm}
$$
\begin{array}{l}
   \pm \dfrac{1}{2}(- \lambda_0 - \lambda_1 + \lambda_2 - \lambda_3) \pm \Big(\dfrac{1}{2}\mu_2 - \dfrac{2}{3}\nu \Big), 
\vspace{1mm}\\
   \pm \dfrac{1}{2}(\;\;\; \lambda_0 + \lambda_1 + \lambda_2 - \lambda_3) \pm \Big(\dfrac{1}{2}\mu_2 - \dfrac{2}{3}\nu \Big), 
\vspace{1mm}\\
   \pm \dfrac{1}{2}(- \lambda_0 + \lambda_1 + \lambda_2 + \lambda_3) \pm \Big(\dfrac{1}{2}\mu_2 - \dfrac{2}{3}\nu \Big),  
\vspace{1mm}\\
   \pm \dfrac{1}{2}(\;\;\; \lambda_0 - \lambda_1 + \lambda_2 + \lambda_3) \pm \Big(\dfrac{1}{2}\mu_2 - \dfrac{2}{3}\nu \Big),  
\end{array}$$
\vspace{-1mm}
$$
\begin{array}{l}
   \pm \dfrac{1}{2}(\;\;\; \lambda_0 - \lambda_1 + \lambda_2 - \lambda_3) \pm \Big(\dfrac{1}{2}\mu_3 - \dfrac{2}{3}\nu \Big), 
\vspace{1mm}\\
   \pm \dfrac{1}{2}(- \lambda_0 + \lambda_1 + \lambda_2 - \lambda_3) \pm \Big(\dfrac{1}{2}\mu_3 - \dfrac{2}{3}\nu \Big), 
\vspace{1mm}\\
   \pm \dfrac{1}{2}(\;\;\; \lambda_0 + \lambda_1 + \lambda_2 + \lambda_3) \pm \Big(\dfrac{1}{2}\mu_3 - \dfrac{2}{3}\nu \Big),  
\vspace{1mm}\\
   \pm \dfrac{1}{2}(- \lambda_0 - \lambda_1 + \lambda_2 + \lambda_3) \pm \Big(\dfrac{1}{2}\mu_3 - \dfrac{2}{3}\nu \Big), 
\end{array} $$
$$
\begin{array}{cl}
     \pm \Big(\mu_j - \dfrac{1}{3}\nu + r \Big),  & 1 \le j \le 3, 
\vspace{1mm}\\
   \pm \lambda_k \pm \Big(\dfrac{1}{2}\mu_1 + \dfrac{1}{3}\nu \Big) \pm r, & 0 \le k \le 3,
\end{array} $$
\vspace{-1mm}
$$
\begin{array}{l}
   \pm \dfrac{1}{2}(\;\;\; \lambda_0 + \lambda_1 - \lambda_2 - \lambda_3) \pm \Big(\dfrac{1}{2}\mu_2 + \dfrac{1}{3}\nu \Big) \pm r,
\vspace{1mm}\\
   \pm \dfrac{1}{2}(\;\;\; \lambda_0 + \lambda_1 - \lambda_2 - \lambda_3) \pm \Big(\dfrac{1}{2}\mu_2 + \dfrac{1}{3}\nu \Big) \pm r,
\end{array}$$
$$
\begin{array}{l}
   \pm \dfrac{1}{2}(\;\;\; \lambda_0 + \lambda_1 + \lambda_2 - \lambda_3) \pm \Big(\dfrac{1}{2}\mu_2 + \dfrac{1}{3}\nu \Big) \pm r, 
\vspace{1mm}\\
   \pm \dfrac{1}{2}(\;\;\; \lambda_0 + \lambda_1 - \lambda_2 + \lambda_3) \pm \Big(\dfrac{1}{2}\mu_2 + \dfrac{1}{3}\nu \Big) \pm r, 
\end{array}$$
\vspace{-1mm}
$$
\begin{array}{l}
   \pm \dfrac{1}{2}(\;\;\; \lambda_0 + \lambda_1 + \lambda_2 + \lambda_3) \pm \Big(\dfrac{1}{2}\mu_3 + \dfrac{1}{3}\nu \Big) \pm r
\vspace{1mm}\\
   \pm \dfrac{1}{2}(\;\;\; \lambda_0 + \lambda_1 - \lambda_2 - \lambda_3) \pm \Big(\dfrac{1}{2}\mu_3 + \dfrac{1}{3}\nu \Big) \pm r,
\vspace{1mm}\\
   \pm \dfrac{1}{2}(\;\;\; \lambda_0 - \lambda_1 + \lambda_2 + \lambda_3) \pm \Big(\dfrac{1}{2}\mu_3 + \dfrac{1}{3}\nu \Big) \pm r, 
\vspace{1mm}\\
   \pm \dfrac{1}{2}(\;\;\; \lambda_0 - \lambda_1 - \lambda_2 + \lambda_3) \pm \Big(\dfrac{1}{2}\mu_3 + \dfrac{2}{3}\nu \Big) \pm r,
\end{array} $$
$$
\begin{array}{cl}
                  \pm 2r,
\vspace{1mm}\\
             \pm \nu \pm r
\end{array} $$
{\it with} $\mu_1 + \mu_2 + \mu_3 = 0.$  
\vspace{2mm}

{\bf Proof.} Let $\gh_7$ be the Cartan subalgebra of ${\gge_7}^C$ given in Theorem 4.6.1, then
$$
  \gh = \{ ({\mit\Phi}(h), 0, 0, r, 0, 0) \, | \, {\mit\Phi}(h) \in \gh_7, r \in C \} $$
is an abelian subalgebra of ${\gge_8}^C$ (it will be a Cartan subalgebra of ${\gge_8}^C$).
\vspace{1mm}

I $\;$ The roots of ${\gge_7}^C$ are also roots of ${\gge_8}^C$. Indeed, we have
$$
  [({\mit\Phi}(h), 0, 0, r, 0, 0), ({\mit\Phi}, 0, 0, 0, 0, 0)] =
         ([{\mit\Phi}(h), {\mit\Phi}], 0, 0, 0, 0, 0). $$

II $\;$ We have
$$ 
      [({\mit\Phi}(h), 0, 0, r, 0, 0), (0, P, 0, 0, 0, 0)] = 
         (0, ({\mit\Phi}(h) + r)P, 0, 0, 0, 0), $$
and using the same notation as in Theorem 4.6.1, we also have
$$
\begin{array}{l}
   ({\mit\Phi}(h_{\delta} + \wti{H}, 0, 0, \nu) + r1)(X, Y, \xi, \eta) 
\vspace{1mm}\\
$$
   = \Big(\Big(h_{\delta} + \wti{H} - \dfrac{1}{3}\nu + r \Big)X, \Big(h_{\delta} - \wti{H} + \dfrac{1}{3}\nu + r \Big)Y, (\nu + r)\xi, (- \nu + r)\eta \Big).  \end{array}$$
By putting $Y =0, \xi = \eta = 0$, we obtain
$$
\begin{array}{lll}
     \mbox{the root}\;\; \mu_k - \dfrac{1}{3}\nu + r & \mbox{by letting}\;\; X = E_k, 
\vspace{1mm}\\
     \mbox{the root}\;\; \pm \lambda_k - \dfrac{1}{2}\mu_1 - \dfrac{1}{3}\nu + r & \mbox{by letting}\;\; X = F_1(a), a = e_k \pm ie_{4+k}. 
\end{array} $$
We can also obtain roots by letting $X = F_2(a), F_3(a)$.

\noindent By putting $X = 0, \xi = \eta = 0$ and further $Y = E_k, Y = F_i(a)$, we can again roots.

\noindent By putting $X = Y = 0, \xi = 1,\eta = 0$, we can again the root $\nu + r$. Similarly, we can obtain the roort $- \nu + r$.

\noindent By using
$$
   [({\mit\Phi}(h), 0, 0, r, 0, 0), (0, 0, Q, 0, 0, 0)] = 
         (0, 0, ({\mit\Phi}(h) - r)Q, 0, 0, 0), $$
we can also obtain roots.
\vspace{1mm}

III $\;$ From the relations
$$
\begin{array}{l}
     [({\mit\Phi}(h), 0, 0, r, 0), (0, 0, 0, 0, 1, 0)] = 
         (0, 0, 0, 0, 2r, 0), 
\vspace{1mm}\\
     {[}({\mit\Phi}(h), 0, 0, r, 0), (0, 0, 0, 0, 0, 1){]} = 
         (0, 0, 0, 0, 0, - 2r),
\end{array}$$
we obtain the roots $2r$ and $-2r$. 
\vspace{3mm}

{\bf Theorem 5.6.2.} {\it In the root system of Theorem} 5.6.1, 
$$
\begin{array}{l}
      \alpha_1 = \dfrac{1}{2}(\lambda_0 - \lambda_1 - \lambda_2 - \lambda_3) + \dfrac{1}{2}(\mu_3 - \mu_1),
\vspace{1mm}\\
      \alpha_2 = \mu_1 - \dfrac{1}{3}\nu - r, \quad
      \alpha_3 = 2r,
\vspace{1mm}\\
      \alpha_4 = \mu_2 - \dfrac{1}{3}\nu - r, \quad
      \alpha_5 = \lambda_3 - \dfrac{1}{2}(\mu_2 - \mu_3), 
\vspace{1mm}\\
      \alpha_6 = \lambda_2 - \lambda_3, \quad
      \alpha_7 = \lambda_1 - \lambda_2, \quad
      \alpha_8 = \nu - r
\end{array}$$
{\it is a fundamental root system of the Lie algebra ${\gge_8}^C$ and 
$$
     \mu = 2\alpha_1 + 4\alpha_2 + 6\alpha_3 + 5\alpha_4 + 4\alpha_5 + 3\alpha_6 + 2\alpha_7 + 3\alpha_8 $$
is the highest root. The Dynkin diagram and the extended Dynkin diagram of ${\gge_8}^C$ are respectively given by }

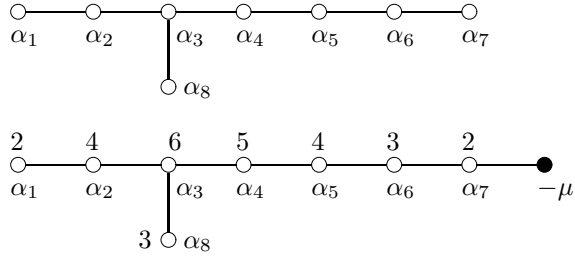
\setlength{\unitlength}{1mm}
\begin{picture}(100,20)
\put(20,10){\circle{2}} \put(19,6){$\alpha_1$}
\put(21,10){\line(1,0){8}}
\put(30,10){\circle{2}} \put(29,6){$\alpha_2$}
\put(31,10){\line(1,0){8}}
\put(40,10){\circle{2}} \put(41,6){$\alpha_3$}
\put(40,8.8){\line(0,-1){7.8}}
\put(40,0){\circle{2}} \put(42,-1){$\alpha_8$}
\put(41,10){\line(1,0){8}}
\put(50,10){\circle{2}} \put(49,6){$\alpha_4$}
\put(51,10){\line(1,0){8}}
\put(60,10){\circle{2}} \put(59,6){$\alpha_5$}
\put(61,10){\line(1,0){8}}
\put(70,10){\circle{2}} \put(69,6){$\alpha_6$}
\put(71,10){\line(1,0){8}}
\put(80,10){\circle{2}} \put(79,6){$\alpha_7$}
\end{picture}

\setlength{\unitlength}{1mm}
\begin{picture}(100,20)
\put(20,10){\circle{2}} \put(19,6){$\alpha_1$} \put(19,12){$2$}
\put(21,10){\line(1,0){8}}
\put(30,10){\circle{2}} \put(29,6){$\alpha_2$} \put(29,12){$4$}
\put(31,10){\line(1,0){8}}
\put(40,10){\circle{2}} \put(41,6){$\alpha_3$} \put(40,12){$6$}
\put(40,8.8){\line(0,-1){7.8}}
\put(40,0){\circle{2}} \put(42,-1){$\alpha_8$} \put(36,-1){$3$}
\put(41,10){\line(1,0){8}}
\put(50,10){\circle{2}} \put(49,6){$\alpha_4$} \put(49,12){$5$}
\put(51,10){\line(1,0){8}}
\put(60,10){\circle{2}} \put(59,6){$\alpha_5$} \put(59,12){$4$}
\put(61,10){\line(1,0){8}}
\put(70,10){\circle{2}} \put(69,6){$\alpha_6$} \put(69,12){$3$}
\put(71,10){\line(1,0){8}}
\put(80,10){\circle{2}} \put(79,6){$\alpha_7$} \put(79,12){$2$}
\put(81,10){\line(1,0){8}}
\put(90,10){\circle*{2}} \put(89,6){$-\mu$}
\end{picture}
\vspace{2mm}

{\bf Proof.} In the following, the notation $n_1n_2 \cdots n_8$ denotes the root  $n_1\alpha_1 + n_2\alpha_2 + \cdots + n_8\alpha_8$. Now, all positive roots of ${\gge_8}^C$ are expressed by 
$$
\begin{array}{lllllllll}
      \lambda_0 - \lambda_1 = 2 & 3 & 4 & 3 & 2 & 1 & 0 & 2
\vspace{1mm}\\
      \lambda_0 - \lambda_2 = 2 & 3 & 4 & 3 & 2 & 1 & 1 & 2
\vspace{1mm}\\
      \lambda_0 - \lambda_3 = 2 & 3 & 4 & 3 & 2 & 2 & 1 & 2
\vspace{1mm}\\
      \lambda_1 - \lambda_2 = 0 & 0 & 0 & 0 & 0 & 0 & 1 & 0
\vspace{1mm}\\
      \lambda_1 - \lambda_3 = 0 & 0 & 0 & 0 & 0 & 1 & 1 & 0
\vspace{1mm}\\                            
      \lambda_2 - \lambda_3 = 0 & 0 & 0 & 0 & 0 & 1 & 0 & 0
\end{array}
\quad
\begin{array}{lllllllll}
      \lambda_0 + \lambda_1 = 2 & 4 & 6 & 5 & 4 & 3 & 2 & 3
\vspace{1mm}\\
      \lambda_0 + \lambda_2 = 2 & 4 & 6 & 5 & 4 & 3 & 1 & 3
\vspace{1mm}\\
      \lambda_0 + \lambda_3 = 2 & 4 & 6 & 5 & 4 & 2 & 1 & 3
\vspace{1mm}\\
      \lambda_1 + \lambda_2 = 0 & 1 & 2 & 2 & 2 & 2 & 1 & 1
\vspace{1mm}\\
      \lambda_1 + \lambda_3 = 0 & 1 & 2 & 2 & 2 & 1 & 1 & 1 
\vspace{1mm}\\
      \lambda_2 + \lambda_3 = 0 & 1 & 2 & 2 & 2 & 1 & 0 & 1
\end{array} $$
$$
\begin{array}{llllllllll}
      \lambda_0 + \dfrac{1}{2}(\mu_2 - \mu_3) = 2 & 4 & 6 & 5 & 3 & 2 & 1 & 3
\vspace{1mm}\\
      \lambda_1 + \dfrac{1}{2}(\mu_2 - \mu_3) = 0 & 1 & 2 & 2 & 1 & 1 & 1 & 1
\end{array}$$
$$
\begin{array}{llllllllll}
      \lambda_2 + \dfrac{1}{2}(\mu_2 - \mu_3) = 0 & 1 & 2 & 2 & 1 & 1 & 0 & 1
\vspace{1mm}\\
      \lambda_3 + \dfrac{1}{2}(\mu_2 - \mu_3) = 0 & 1 & 2 & 2 & 1 & 0 & 0 & 1
\end{array}$$
\vspace{-1mm}
$$
\begin{array}{lllllllllllll}
      \lambda_0 - \dfrac{1}{2}(\mu_2 - \mu_3) = 2 & 3 & 4 & 3 & 3 & 2 & 1 & 2
\vspace{1mm}\\
      \lambda_1 - \dfrac{1}{2}(\mu_2 - \mu_3) = 0 & 0 & 0 & 0 & 1 & 1 & 1 & 0
\vspace{1mm}\\
      \lambda_2 - \dfrac{1}{2}(\mu_2 - \mu_3) = 0 & 0 & 0 & 0 & 1 & 1 & 0 & 0
\vspace{1mm}\\
      \lambda_3 - \dfrac{1}{2}(\mu_2 - \mu_3) = 0 & 0 & 0 & 0 & 1 & 0 & 0 & 0
\end{array}$$
$$
\begin{array}{llllllllll}
     \dfrac{1}{2}(\lambda_0 + \lambda_1 + \lambda_2 - \lambda_3) + \dfrac{1}{2}(\mu_3 - \mu_1) = 1 & 1 & 2 & 2 & 2 & 2 & 1 & 1 
\vspace{1mm}\\
     \dfrac{1}{2}(\lambda_0 + \lambda_1 - \lambda_2 + \lambda_3) + \dfrac{1}{2}(\mu_3 - \mu_1) = 1 & 1 & 2 & 2 & 2 & 1 & 1 & 1
\vspace{1mm}\\
     \dfrac{1}{2}(\lambda_0 - \lambda_1 + \lambda_2 + \lambda_3) + \dfrac{1}{2}(\mu_3 - \mu_1) = 1 & 1 & 2 & 2 & 2 & 1 & 0 & 1
\vspace{1mm}\\
     \dfrac{1}{2}(\lambda_0 - \lambda_1 - \lambda_2 - \lambda_3) + \dfrac{1}{2}(\mu_3 - \mu_1) = 1 & 0 & 0 & 0 & 0 & 0 & 0 & 0 
\end{array}$$
\vspace{-1mm}
$$
\begin{array}{llllllllllll}
     \dfrac{1}{2}(\lambda_0 + \lambda_1 + \lambda_2 - \lambda_3) - \dfrac{1}{2}(\mu_3 - \mu_1) = 1 & 3 & 4 & 3 & 2 & 2 & 1 & 2 
\vspace{1mm}\\
     \dfrac{1}{2}(\lambda_0 + \lambda_1 - \lambda_2 + \lambda_3) - \dfrac{1}{2}(\mu_3 - \mu_1) = 1 & 3 & 4 & 3 & 2 & 1 & 1 & 2
\vspace{1mm}\\
     \dfrac{1}{2}(\lambda_0 - \lambda_1 + \lambda_2 + \lambda_3) - \dfrac{1}{2}(\mu_3 - \mu_1) = 1 & 3 & 4 & 3 & 2 & 1 & 0 & 2 
\vspace{1mm}\\
     \dfrac{1}{2}(\lambda_0 - \lambda_1 - \lambda_2 - \lambda_3) - \dfrac{1}{2}(\mu_3 - \mu_1) = 1 & 2 & 2 & 1 & 0 & 0 & 0 & 1
\end{array} $$
$$
\begin{array}{llllllll}
   \dfrac{1}{2}(\lambda_0 - \lambda_1 + \lambda_2 - \lambda_3) + \dfrac{1}{2}(\mu_1 - \mu_2) = 1 & 2 & 2 & 1 & 1 & 1 & 0 & 1 
\vspace{1mm}\\
   \dfrac{1}{2}(\lambda_0 - \lambda_1 - \lambda_2 + \lambda_3) + \dfrac{1}{2}(\mu_1 - \mu_2) = 1 & 2 & 2 & 1 & 1 & 0 & 0 & 1 
\vspace{1mm}\\
   \dfrac{1}{2}(\lambda_0 + \lambda_1 - \lambda_2 - \lambda_3) + \dfrac{1}{2}(\mu_1 - \mu_2) = 1 & 2 & 2 & 1 & 1 & 1 & 1 & 1
\vspace{1mm}\\
   \dfrac{1}{2}(\lambda_0 + \lambda_1 + \lambda_2 + \lambda_3) + \dfrac{1}{2}(\mu_1 - \mu_2) = 1 & 3 & 4 & 3 & 3 & 2 & 1 & 2 
\end{array}$$
\vspace{-1mm}
$$
\begin{array}{lllllllllllll}
   \dfrac{1}{2}(\lambda_0 - \lambda_1 + \lambda_2 - \lambda_3) - \dfrac{1}{2}(\mu_1 - \mu_2) = 1 & 1 & 2 & 2 & 1 & 1 & 0 & 1 
\vspace{1mm}\\
   \dfrac{1}{2}(\lambda_0 - \lambda_1 - \lambda_2 + \lambda_3) - \dfrac{1}{2}(\mu_1 - \mu_2) = 1 & 1 & 2 & 2 & 1 & 0 & 0 & 1
\vspace{1mm}\\
   \dfrac{1}{2}(\lambda_0 + \lambda_1 - \lambda_2 - \lambda_3) - \dfrac{1}{2}(\mu_1 - \mu_2) = 1 & 1 & 2 & 2 & 1 & 1 & 1 & 1 
\vspace{1mm}\\
   \dfrac{1}{2}(\lambda_0 + \lambda_1 + \lambda_2 + \lambda_3) - \dfrac{1}{2}(\mu_1 - \mu_2) = 1 & 2 & 4 & 4 & 3 & 2 & 1 & 2
\end{array} $$
$$
\begin{array}{lllllllllllll}
       \;\;\;\mu_1 + \dfrac{2}{3}\nu = 0 & 1 & 1 & 0 & 0 & 0 & 0 & 1
\vspace{1mm}\\
       \;\;\;\mu_2 + \dfrac{2}{3}\nu = 0 & 0 & 1 & 1 & 0 & 0 & 0 & 1
\vspace{1mm}\\
       - \mu_3 - \dfrac{2}{3}\nu = 0 & 1 & 1 & 1 & 0 & 0 & 0 & 0
\end{array}$$
$$
\begin{array}{lllllllllll}
      \lambda_0 + \dfrac{1}{2}\mu_1 - \dfrac{2}{3}\nu = 2 & 4 & 5 & 4 & 3 & 2 &    1 & 2
\vspace{1mm}\\
      \lambda_0 - \dfrac{1}{2}\mu_1 + \dfrac{2}{3}\nu = 2 & 3 & 5 & 4 & 3 & 2 &  1 & 3
\vspace{1mm}\\
      \lambda_1 + \dfrac{1}{2}\mu_1 - \dfrac{2}{3}\nu = 0 & 1 & 1 & 1 & 1 & 1 &  1 & 0
\end{array}$$
$$
\begin{array}{lllllllllll}
      \lambda_1 - \dfrac{1}{2}\mu_1 + \dfrac{2}{3}\nu = 0 & 0 & 1 & 1 & 1 & 1 & 1 & 1
\end{array}$$
\vspace{-1mm}
$$
\begin{array}{lllllllllll}
      \lambda_2 + \dfrac{1}{2}\mu_1 - \dfrac{2}{3}\nu = 0 & 1 & 1 & 1 & 1 & 1 & 0 & 0
\vspace{1mm}\\
      \lambda_2 - \dfrac{1}{2}\mu_1 + \dfrac{2}{3}\nu = 0 & 0 & 1 & 1 & 1 & 1 & 0 & 1
\vspace{1mm}\\
      \lambda_3 + \dfrac{1}{2}\mu_1 - \dfrac{2}{3}\nu = 0 & 1 & 1 & 1 & 1 & 0 & 0 & 0
\vspace{1mm}\\
      \lambda_3 - \dfrac{1}{2}\mu_1 + \dfrac{2}{3}\nu = 0 & 0 & 1 & 1 & 1 & 0 & 0 & 1
\end{array}$$
$$
\begin{array}{lllllllllll}
   \dfrac{1}{2}(\lambda_0 + \lambda_1 + \lambda_2 - \lambda_3) + \dfrac{1}{2}\mu_2 - \dfrac{2}{3}\nu = 1 & 2 & 3 & 3 & 2 & 2 & 1 & 1
\vspace{1mm}\\
   \dfrac{1}{2}(\lambda_0 + \lambda_1 - \lambda_2 + \lambda_3) + \frac{1}{2}\mu_2 - \dfrac{2}{3}\nu = 1 & 2 & 3 & 3 & 2 & 1 & 1 & 1
\vspace{1mm}\\
   \dfrac{1}{2}(\lambda_0 - \lambda_1 + \lambda_2 + \lambda_3) + \dfrac{1}{2}\mu_2 - \dfrac{2}{3}\nu = 1 & 2 & 3 & 3 & 2 & 1 & 0 & 1
\vspace{1mm}\\
   \dfrac{1}{2}(\lambda_0 - \lambda_1 - \lambda_2 - \lambda_3) + \frac{1}{2}\mu_2 - \dfrac{2}{3}\nu = 1 & 1 & 1 & 1 & 0 & 0 & 0 & 0
\end{array}$$
\vspace{-1mm}
$$
\begin{array}{lllllllllllll}
   \dfrac{1}{2}(\lambda_0 + \lambda_1 + \lambda_2 - \lambda_3) - \dfrac{1}{2}\mu_2 + \dfrac{2}{3}\nu = 1 & 2 & 3 & 2 & 2 & 2 & 1 & 2
\vspace{1mm}\\
   \dfrac{1}{2}(\lambda_0 + \lambda_1 - \lambda_2 + \lambda_3) - \dfrac{1}{2}\mu_2 + \dfrac{2}{3}\nu = 1 & 2 & 3 & 2 & 2 & 1 & 1 & 2
\vspace{1mm}\\
   \dfrac{1}{2}(\lambda_0 - \lambda_1 + \lambda_2 + \lambda_3) - \dfrac{1}{2}\mu_2 - \dfrac{2}{3}\nu = 1 & 2 & 3 & 2 & 2 & 1 & 0 & 2
\vspace{1mm}\\
   \dfrac{1}{2}(\lambda_0 - \lambda_1 - \lambda_2 - \lambda_3) - \dfrac{1}{2}\mu_2 + \dfrac{2}{3}\nu = 1 & 1 & 1 & 0 & 0 & 0 & 0 & 1
\end{array}$$
$$
\begin{array}{lllllllllll}
   \dfrac{1}{2}(\lambda_0 + \lambda_1 + \lambda_2 + \lambda_3) + \dfrac{1}{2}\mu_3 - \dfrac{2}{3}\nu = 1 & 2 & 3 & 3 & 3 & 2 & 1 & 1
\vspace{1mm}\\
   \dfrac{1}{2}(\lambda_0 + \lambda_1 - \lambda_2 - \lambda_3) + \dfrac{1}{2}\mu_3 - \dfrac{2}{3}\nu = 1 & 1 & 1 & 1 & 1 & 1 & 1 & 0
\vspace{1mm}\\
   \dfrac{1}{2}(\lambda_0 - \lambda_1 + \lambda_2 - \lambda_3) + \dfrac{1}{2}\mu_3 - \dfrac{2}{3}\nu = 1 & 1 & 1 & 1 & 1 & 1 & 0 & 0
\vspace{1mm}\\
   \dfrac{1}{2}(\lambda_0 - \lambda_1 - \lambda_2 + \lambda_3) + \dfrac{1}{2}\mu_3 - \dfrac{2}{3}\nu = 1 & 1 & 1 & 1 & 1 & 0 & 0 & 0
\end{array}$$
\vspace{-1mm}
$$
\begin{array}{llllllllllll}
   \dfrac{1}{2}(\lambda_0 + \lambda_1 + \lambda_2 + \lambda_3) - \dfrac{1}{2}\mu_3 + \dfrac{2}{3}\nu = 1 & 3 & 5 & 4 & 3 & 2 & 1 & 3
\vspace{1mm}\\
   \dfrac{1}{2}(\lambda_0 + \lambda_1 - \lambda_2 - \lambda_3) - \dfrac{1}{2}\mu_3 + \dfrac{2}{3}\nu = 1 & 2 & 3 & 2 & 1 & 1 & 1 & 2
\vspace{1mm}\\
   \dfrac{1}{2}(\lambda_0 - \lambda_1 + \lambda_2 - \lambda_3) - \dfrac{1}{2}\mu_3 + \dfrac{2}{3}\nu = 1 & 2 & 3 & 2 & 1 & 1 & 0 & 2
\vspace{1mm}\\
   \dfrac{1}{2}(\lambda_0 - \lambda_1 - \lambda_2 + \lambda_3) - \dfrac{1}{2}\mu_3 + \dfrac{2}{3}\nu = 1 & 2 & 3 & 2 & 1 & 0 & 0 & 2
\end{array}$$
$$
\begin{array}{lllllllll}
      \;\;\; \mu_1 - \dfrac{1}{3}\nu + r = 0 & 1 & 1 & 0 & 0 & 0 & 0 & 0
\vspace{1mm}\\
      \;\;\; \mu_2 - \dfrac{1}{3}\nu + r = 0 & 0 & 1 & 1 & 0 & 0 & 0 & 0
\vspace{1mm}\\
      - \mu_3 + \dfrac{1}{3}\nu + r = 0 & 1 & 2 & 1 & 0 & 0 & 0 & 1
\end{array}$$
\vspace{-1mm}
$$
\begin{array}{llllllllllll}
      \;\;\; \mu_1 - \dfrac{1}{3}\nu - r = 0 & 1 & 0 & 0 & 0 & 0 & 0 & 0
\vspace{1mm}\\
      \;\;\; \mu_2 - \dfrac{1}{3}\nu + r = 0 & 0 & 0 & 1 & 0 & 0 & 0 & 0
\vspace{1mm}\\
       - \mu_3 - \dfrac{1}{3}\nu - r = 0 & 1 & 1 & 1 & 0 & 0 & 0 & 1
\vspace{1mm}\\
\end{array} $$
$$
\begin{array}{lllllllll}
    \lambda_0 + \dfrac{1}{2}\mu_1 + \dfrac{1}{3}\nu + r = 2 & 4 & 6 & 4 & 3 & 2 & 1& 3
\vspace{1mm}\\
    \lambda_0 + \dfrac{1}{2}\mu_1 + \dfrac{1}{3}\nu - r = 2 & 4 & 5 & 4 & 3 & 2 & 1 & 3
\vspace{1mm}\\
    \lambda_0 - \dfrac{1}{2}\mu_1 - \dfrac{1}{3}\nu + r = 2 & 3 & 5 & 4 & 3 & 2 & 1& 2
\vspace{1mm}\\
    \lambda_0 - \dfrac{1}{2}\mu_1 - \dfrac{1}{3}\nu - r = 2 & 3 & 4 & 4 & 3 & 2 & 1 & 2
\end{array}$$
\vspace{-1mm}
$$
\begin{array}{llllllllllll}
    \lambda_1 + \dfrac{1}{2}\mu_1 + \dfrac{1}{3}\nu + r = 0 & 1 & 2 & 1 & 1 & 1 & 1 & 1
\vspace{1mm}\\
    \lambda_1 + \dfrac{1}{2}\mu_1 + \dfrac{1}{3}\nu - r = 0 & 1 & 1 & 1 & 1 & 1 & 1 & 1
\vspace{1mm}\\
    \lambda_1 - \dfrac{1}{2}\mu_1 - \dfrac{1}{3}\nu + r = 0 & 0 & 1 & 1 & 1 & 1 & 1 & 0
\vspace{1mm}\\
    \lambda_1 - \dfrac{1}{2}\mu_1 - \dfrac{1}{3}\nu - r = 0 & 0 & 0 & 1 & 1 & 1 & 1 & 0
\end{array}$$
$$
\begin{array}{lllllllll}
    \lambda_2 + \dfrac{1}{2}\mu_1 + \dfrac{1}{3}\nu + r = 0 & 1 & 2 & 1 & 1 & 1 & 0 & 1
\vspace{1mm}\\
    \lambda_2 + \dfrac{1}{2}\mu_1 + \dfrac{1}{3}\nu - r = 0 & 1 & 1 & 1 & 1 & 1 & 0 & 1
\vspace{1mm}\\
    \lambda_2 - \dfrac{1}{2}\mu_1 - \dfrac{1}{3}\nu + r = 0 & 0 & 1 & 1 & 1 & 1 & 0 & 0
\vspace{1mm}\\
    \lambda_2 - \dfrac{1}{2}\mu_1 - \dfrac{1}{3}\nu - r = 0 & 0 & 0 & 1 & 1 & 1 & 0 & 0
\end{array}$$
\vspace{-1mm}
$$
\begin{array}{llllllllllll}
    \lambda_3 + \dfrac{1}{2}\mu_1 + \dfrac{1}{3}\nu + r = 0 & 1 & 2 & 1 & 1 & 0 & 0 & 1
\vspace{1mm}\\
    \lambda_3 + \dfrac{1}{2}\mu_1 + \dfrac{1}{3}\nu - r = 0 & 1 & 1 & 1 & 1 & 0 & 0 & 1
\vspace{1mm}\\
    \lambda_3 - \dfrac{1}{2}\mu_1 - \dfrac{1}{3}\nu + r = 0 & 0 & 1 & 1 & 1 & 0 & 0 & 0
\vspace{1mm}\\
    \lambda_3 - \dfrac{1}{2}\mu_1 - \dfrac{1}{3}\nu - r = 0 & 0 & 0 & 1 & 1 & 0 & 0 & 0
\end{array}$$
\vspace{-1mm}
$$
\begin{array}{lllllllll}
   \dfrac{1}{2}(\lambda_0 + \lambda_1 + \lambda_2 - \lambda_3) + \dfrac{1}{2}\mu_2 + \dfrac{1}{3}\nu + r = 1 & 2 & 4 & 3 & 2 & 2 & 1 & 2
\vspace{1mm}\\  
   \dfrac{1}{2}(\lambda_0 + \lambda_1 + \lambda_2 - \lambda_3) + \dfrac{1}{2}\mu_2 + \dfrac{1}{3}\nu - r = 1 & 2 & 3 & 3 & 2 & 2 & 1 & 2
\end{array}$$
$$
\begin{array}{lllllllll}
\vspace{1mm}\\  
   \dfrac{1}{2}(\lambda_0 + \lambda_1 + \lambda_2 - \lambda_3) - \dfrac{1}{2}\mu_2 - \dfrac{1}{3}\nu + r = 1 & 2 & 3 & 2 & 2 & 2 & 1 & 1
\vspace{1mm}\\  
   \dfrac{1}{2}(\lambda_0 + \lambda_1 + \lambda_2 - \lambda_3) - \dfrac{1}{2}\mu_2 - \dfrac{1}{3}\nu - r = 1 & 2 & 2 & 2 & 2 & 2 & 1 & 1
\end{array}$$
\vspace{-1mm}
$$
\begin{array}{llllllllllll}
   \dfrac{1}{2}(\lambda_0 + \lambda_1 - \lambda_2 + \lambda_3) + \dfrac{1}{2}\mu_2 + \dfrac{1}{3}\nu + r = 1 & 2 & 4 & 3 & 2 & 1 & 1 & 2
\vspace{1mm}\\  
   \dfrac{1}{2}(\lambda_0 + \lambda_1 - \lambda_2 + \lambda_3) + \dfrac{1}{2}\mu_2 + \dfrac{1}{3}\nu - r = 1 & 2 & 3 & 3 & 2 & 1 & 1 & 2
\vspace{1mm}\\  
   \dfrac{1}{2}(\lambda_0 + \lambda_1 - \lambda_2 + \lambda_3) - \dfrac{1}{2}\mu_2 - \dfrac{1}{3}\nu + r = 1 & 2 & 3 & 2 & 2 & 1 & 1 & 1
\vspace{1mm}\\  
   \dfrac{1}{2}(\lambda_0 + \lambda_1 - \lambda_2 + \lambda_3) - \dfrac{1}{2}\mu_2 - \dfrac{1}{3}\nu - r = 1 & 2 & 2 & 2 & 2 & 1 & 1 & 1
\end{array}$$
$$
\begin{array}{lllllllll}
   \dfrac{1}{2}(\lambda_0 - \lambda_1 + \lambda_2 + \lambda_3) + \dfrac{1}{2}\mu_2 + \dfrac{1}{3}\nu + r = 1 & 2 & 4 & 3 & 2 & 1 & 0 & 2
\vspace{1mm}\\  
   \dfrac{1}{2}(\lambda_0 - \lambda_1 + \lambda_2 + \lambda_3) + \dfrac{1}{2}\mu_2 + \dfrac{1}{3}\nu - r = 1 & 2 & 3 & 3 & 2 & 1 & 0 & 2
\vspace{1mm}\\  
   \dfrac{1}{2}(\lambda_0 - \lambda_1 + \lambda_2 + \lambda_3) - \dfrac{1}{2}\mu_2 - \dfrac{1}{3}\nu + r = 1 & 2 & 3 & 2 & 2 & 1 & 0 & 1
\vspace{1mm}\\  
   \dfrac{1}{2}(\lambda_0 - \lambda_1 + \lambda_2 + \lambda_3) - \dfrac{1}{2}\mu_2 - \dfrac{1}{3}\nu - r = 1 & 2 & 2 & 2 & 2 & 1 & 0 & 1
\end{array}$$
\vspace{-1mm}
$$
\begin{array}{llllllllllll}  
   \dfrac{1}{2}(\lambda_0 - \lambda_1 - \lambda_2 - \lambda_3) + \dfrac{1}{2}\mu_2 + \dfrac{1}{3}\nu + r = 1 & 1 & 2 & 1 & 0 & 0 & 0 & 1
\vspace{1mm}\\  
   \dfrac{1}{2}(\lambda_0 - \lambda_1 - \lambda_2 - \lambda_3) + \dfrac{1}{2}\mu_2 + \dfrac{1}{3}\nu - r = 1 & 1 & 1 & 1 & 0 & 0 & 0 & 1
\vspace{1mm}\\  
   \dfrac{1}{2}(\lambda_0 - \lambda_1 - \lambda_2 - \lambda_3) - \dfrac{1}{2}\mu_2 - \dfrac{1}{3}\nu + r = 1 & 1 & 1 & 0 & 0 & 0 & 0 & 0
\vspace{1mm}\\  
   \dfrac{1}{2}(\lambda_0 - \lambda_1 - \lambda_2 - \lambda_3) - \dfrac{1}{2}\mu_2 - \dfrac{1}{3}\nu - r = 1 & 1 & 0 & 0 & 0 & 0 & 0 & 0
\end{array}$$  
$$
\begin{array}{lllllllll}
   \dfrac{1}{2}(\lambda_0 + \lambda_1 - \lambda_2 - \lambda_3) + \dfrac{1}{2}\mu_3 + \dfrac{1}{3}\nu + r = 1 & 1 & 2 & 1 & 1 & 1 & 1 & 1
\vspace{1mm}\\  
   \dfrac{1}{2}(\lambda_0 + \lambda_1 - \lambda_2 - \lambda_3) + \dfrac{1}{2}\mu_3 + \dfrac{1}{3}\nu - r = 1 & 1 & 1 & 1 & 1 & 1 & 1 & 1
\vspace{1mm}\\  
   \dfrac{1}{2}(\lambda_0 + \lambda_1 - \lambda_2 - \lambda_3) - \dfrac{1}{2}\mu_3 - \dfrac{1}{3}\nu + r = 1 & 2 & 3 & 2 & 1 & 1 & 1 & 1
\vspace{1mm}\\  
   \dfrac{1}{2}(\lambda_0 + \lambda_1 - \lambda_2 - \lambda_3) - \dfrac{1}{2}\mu_3 - \dfrac{1}{3}\nu - r = 1 & 2 & 2 & 2 & 1 & 1 & 1 & 1   
\end{array}$$
\vspace{-1mm}
$$
\begin{array}{llllllllllll}  
   \dfrac{1}{2}(\lambda_0 - \lambda_1 + \lambda_2 - \lambda_3) + \dfrac{1}{2}\mu_3 + \dfrac{1}{3}\nu + r = 1 & 1 & 2 & 1 & 1 & 1 & 0 & 1
\vspace{1mm}\\  
   \dfrac{1}{2}(\lambda_0 - \lambda_1 + \lambda_2 - \lambda_3) + \dfrac{1}{2}\mu_3 + \dfrac{1}{3}\nu - r = 1& 1 & 1 & 1 & 1 & 1 & 0 & 1
\vspace{1mm}\\  
   \dfrac{1}{2}(\lambda_0 - \lambda_1 + \lambda_2 - \lambda_3) - \dfrac{1}{2}\mu_3 - \dfrac{1}{3}\nu + r = 1 & 2 & 3 & 2 & 1 & 1 & 0 & 1
\vspace{1mm}\\  
   \dfrac{1}{2}(\lambda_0 - \lambda_1 + \lambda_2 - \lambda_3) - \dfrac{1}{2}\mu_3 - \dfrac{1}{3}\nu - r = 1 & 2 & 2 & 2 & 1 & 1 & 0 & 1
\end{array}$$
$$
\begin{array}{lllllllll}
   \dfrac{1}{2}(\lambda_0 + \lambda_1 + \lambda_2 + \lambda_3) + \dfrac{1}{2}\mu_3 + \dfrac{1}{3}\nu + r = 1 & 2 & 4 & 3 & 3 & 2 & 1 & 2
\vspace{1mm}\\  
   \dfrac{1}{2}(\lambda_0 + \lambda_1 + \lambda_2 + \lambda_3) + \dfrac{1}{2}\mu_3 + \dfrac{1}{3}\nu - r = 1 & 2 & 3 & 3 & 3 & 2 & 1 & 2
\end{array}$$
$$
\begin{array}{lllllllll}  
   \dfrac{1}{2}(\lambda_0 + \lambda_1 + \lambda_2 + \lambda_3) - \dfrac{1}{2}\mu_3 - \dfrac{1}{3}\nu + r = 1 & 3 & 5 & 4 & 3 & 2 & 1 & 2
\vspace{1mm}\\  
   \dfrac{1}{2}(\lambda_0 + \lambda_1 + \lambda_2 + \lambda_3) - \dfrac{1}{2}\mu_3 - \dfrac{1}{3}\nu - r = 1 & 3 & 4 & 4 & 3 & 2 & 1 & 2
\end{array}$$
\vspace{-1mm}
$$
\begin{array}{llllllllllll}
   \dfrac{1}{2}(\lambda_0 - \lambda_1 - \lambda_2 + \lambda_3) + \dfrac{1}{2}\mu_3 + \dfrac{1}{3}\nu + r = 1 & 1 & 2 & 1 & 1 & 0 & 0 & 1
\vspace{1mm}\\  
   \dfrac{1}{2}(\lambda_0 - \lambda_1 - \lambda_2 + \lambda_3) + \dfrac{1}{2}\mu_3 + \dfrac{1}{3}\nu - r = 1 & 1 & 1 & 1 & 1 & 0 & 0 & 1
\vspace{1mm}\\  
   \dfrac{1}{2}(\lambda_0 - \lambda_1 - \lambda_2 + \lambda_3) - \dfrac{1}{2}\mu_3 - \dfrac{1}{3}\nu + r = 1 & 2 & 3 & 2 & 1 & 0 & 0 & 1
\vspace{1mm}\\  
   \dfrac{1}{2}(\lambda_0 - \lambda_1 - \lambda_2 + \lambda_3) - \dfrac{1}{2}\mu_3 - \dfrac{1}{3}\nu - r = 1 & 2 & 2 & 2 & 1 & 0 & 0 & 1
\end{array}$$
$$
\begin{array}{lllllllllll}
       \nu - r& = 0 & 0 & 0 & 0 & 0 & 0 & 0 & 1 
\vspace{1mm}\\
       {} 2r& = 0 & 0 & 1 & 0 & 0 & 0 & 0 & 0
\vspace{1mm}\\
       \nu + r& = 0 & 0 & 1 & 0 & 0 & 0 & 0 & 1.
\end{array}$$ 
Hence ${\mit\Pi} = \{ \alpha_1, \alpha_2, \cdots, \alpha_8 \}$ is a fundamental root system of ${\gge_8}^C$. The real part $\gh_{\sR}$ of $\gh$ is
$$
   \gh_{\sR} = \{({\mit\Phi}(h), 0, 0, r, 0) \, | \, {\mit\Phi}(h) \in (\gh_7)_{\sR}, r \in \R\},$$
\big(where ${\mit\Phi}(h) = \{{\mit\Phi}\Big(\dsum_{k=0}^3\lambda_kH_k + \Big(\dsum_{j=1}^3\mu_jE_j\Big)^{\sim}, 0, 0, \nu \Big) \in (\gh_7)_{\sR}$ (Theorem 4.6.2)\big). The Killing form $B_8$ of ${\gge_8}^C$ on $\gh_{\sR}$ is given by
$$
      B_8(\wti{h}, \wti{h}') 
     = 60\dsum_{k=0}^3\lambda_k{\lambda_k}' + 30\dsum_{j=1}^3\mu_j{\mu_j}' + 40\nu\nu' + 120rr', 
\vspace{-3mm}$$
for $\wti{h} = ({\mit\Phi}(h), 0, 0, r, 0)$, $\wti{h}' = ({\mit\Phi}(h'), 0, 0,$ $ r', 0) \in \gh$, where ${\mit\Phi}(h) = {\mit\Phi}\Big(\dsum_{k=0}^3\lambda_kH_k + \Big(\dsum_{j=1}^3\mu_jE_j\Big)^\sim, 0, 0, \nu \Big)$, ${\mit\Phi}(h') = {\mit\Phi}(\dsum_{k=0}^3{\lambda_k}'H_k + \Big(\dsum_{j=1}^3\mu_j'E_j\Big)^\sim, 0, 0, \nu'\Big) \in (\gh_7)_{\sR}$, Indeed,

\begin{eqnarray*}
   B_8(\wti{h}, \wti{h}') \!\!\!&=&\!\!\!
 \dfrac{5}{3}B_7({\mit\Phi}(h), {\mit\Phi}(h')) + 120rr' \;\;\mbox{(Theorem 5.3.2)}
\vspace{1mm}\\
    \!\!\!&=&\!\!\! \dfrac{5}{3}6\Big(6\dsum_{k=0}^3\lambda_k{\lambda_k}' + 3\dsum_{j=1}^3\mu_j{\mu_j}' + 4\nu\nu'\Big) + 120rr' \;\; \mbox{(Theorem 4.6.2)}
\vspace{1mm}\\
     \!\!\!&=&\!\!\! 60\dsum_{k=0}^3\lambda_k{\lambda_k}' + 30\dsum_{j=1}^3\mu_j{\mu_j}' + 40\nu\nu' + 120rr'.
\end{eqnarray*}
Now, the canonical elements $H_{\alpha_i} \in \gh_{\sR}$ associated with $\alpha_i$ ($B_8 (H_\alpha, H) = \alpha(H), \; H \in \gh_{\sR} $) are determined as follows.
\vspace{2mm}

\quad
    $H_{\alpha_1} = \Big({\mit\Phi}\Big(\dfrac{1}{120}(H_0 - H_1 - H_2 - H_3) + 2(E_3 - E_1)^{\sim}, 0, 0, 0 \Big), 0, 0, 0, 0, 0 \Big)$,
\vspace{1mm}    

\quad
  $H_{\alpha_2} = \Big({\mit\Phi}\Big(\dfrac{1}{90}(2E_1 - E_2 - E_3)^{\sim}, 0, 0, - \dfrac{1}{120} \Big), 0, 0, - \dfrac{1}{120}, 0, 0 \Big)$,
\vspace{1mm}

\quad
    $H_{\alpha_3} = \Big(0, 0, 0, \dfrac{1}{60}, 0, 0 \Big)$,
\vspace{1mm}

\quad
    $H_{\alpha_4} = \Big({\mit\Phi}\Big(\dfrac{1}{90}(- E_1 + 2E_2 - E_3)^{\sim}, 0, 0, - \dfrac{1}{120} \Big), 0, 0, - \dfrac{1}{120}, 0, 0 \Big)$, 
\vspace{1mm}

\quad
    $H_{\alpha_5} = \Big({\mit\Phi}\Big(\dfrac{1}{60}(H_3 - (E_2 - E_3)^{\sim}), 0, 0, 0 \Big), 0, 0, 0, 0, 0 \Big)$,
\vspace{1mm}

\quad
    $H_{\alpha_6} = \Big({\mit\Phi}\Big(\dfrac{1}{60}(H_2 - H_3), 0, 0, 0 \Big), 0, 0, 0, 0, 0 \Big)$,
\vspace{1mm}

\quad
    $H_{\alpha_7} = \Big({\mit\Phi}\Big(\dfrac{1}{60}(H_1 - H_2), 0, 0, 0 \Big), 0, 0, 0, 0, 0 \Big)$,
\vspace{1mm}

\quad
    $H_{\alpha_8} = \Big({\mit\Phi}\Big(0, 0, 0, \dfrac{1}{40} \Big), 0, 0, - \dfrac{1}{120}, 0, 0 \Big)$.
\vspace{2mm}

\noindent Thus we have
$$
     (\alpha_1, \alpha_1) = B_7(H_{\alpha_1}, H_{\alpha_1}) = 60\frac{1}{120}\frac{1}{120}4 + 30\dfrac{1}{120}\dfrac{1}{120}8 = \frac{1}{30}, $$
and the other inner products are similarly calculated. Consequently, the inner product induced by the Killing form $B_8$ between $\alpha_1, \alpha_2, \cdots, \alpha_8$ and $- \mu$ are given by
$$
\begin{array}{l}
      (\alpha_i, \alpha_i) = \dfrac{1}{30}, \quad i = 1, 2, 3, 4, 5, 6, 7, 8,
\vspace{1mm}\\
      (\alpha_i, \alpha_{i+1}) = - \dfrac{1}{60}, \quad i = 1, 2, \cdots, 6, \quad 
      (\alpha_3, \alpha_8) = - \dfrac{1}{60},
\vspace{1mm}\\
       (\alpha_i, \alpha_j) = 0, \quad \mbox{otherwise},
\vspace{1mm}\\
       (-\mu, -\mu) = \dfrac{1}{30}, \;\; (- \mu, \alpha_7) = - \dfrac{1}{60}, \quad (- \mu, \alpha_i) = 0, \;\; i = 1, 2, 3, 4, 5, 6, 8,
\end{array} $$
using them, we can draw the Dynkin diagram and the extended Dynkin diagram of ${\gge_8}^C$.
\vspace{3mm}

According to Borel-Siebenthal theory, the Lie algebra $\gge_8$ has five subalgebras as maximal subalgebras with the maximal rank 8. 
\vspace{1mm}

(1) The first one is a subalgebra of type $A_1 \oplus E_7$ which is obtained as the fixed points of an involution $\upsilon$ of $\gge_8$.
\vspace{1mm}

(2) The second one is a subalgebra of type $D_8$ which is obtained as the fixed points of an involution $\wti{\lambda}\gamma$ of $\gge_8$.
\vspace{1mm}

(3) The third one is a subalgebra of type $A_2 \oplus E_6$ which is obtained as the fixed points of an automorphism $w$ of order 3 of $\gge_8$.
\vspace{1mm}

(4) The fourth one is a subalgebra of type $A_8$ which is obtained as the fixed points of an automorphism $w_3$ of order 3 of $\gge_8$.
\vspace{1mm}

(5) The fifth one is a subalgebra of type $A_4 \oplus A_4$ which is obtained as the fixed points of an automorphism $z_5$ of order 5 of $\gge_8$.
\vspace{1mm}

These subalgebras will be realized as subgroups of the group $E_8$ in Theorems 5.7.6, 5.8.7, 5.10.2, 5.11.7 and 5.12.5, respectively. As for Theorems 5.10.2, 5.11.7 and 5.12.5, we refer to Gomyo [9]. 
\vspace{4mm}

{\bf 5.7. Involution $\upsilon$ and subgroup $(SU(2) \times E_7)/\Z_2$ of $E_8$}\vspace{3mm}

We shall first study the following subgroup $({E_8}^C)_{1,1^-,1_-}$ of ${E_8}^C$:
$$
     ({E_8}^C)_{1,1^-,1_-}
       = \{ \alpha \in {E_8}^C \, | \, \alpha\1 = \1, \alpha 1^- = 1^-, \alpha 1_- = 1_- \}. $$

{\bf Proposition 5.7.1.} \qquad \qquad \quad $({E_8}^C)_{1,1^-,1_-} \cong {E_7}^C. $
\vspace{2mm}

{\bf Proof.} For $\beta \in {E_7}^C$, we define a $C$-linear mapping $\wti{\beta} : {\gge_8}^C \to {\gge_8}^C$ by
$$
    \wti{\beta} = \pmatrix{\mbox{Ad}\beta & 0 & 0 & 0 & 0 & 0 \cr
                            0 & \beta & 0 & 0 & 0 & 0 \cr
                            0 & 0 & \beta & 0 & 0 & 0 \cr
                            0 & 0 & 0 & 1 & 0 & 0 \cr
                            0 & 0 & 0 & 0 & 1 & 0 \cr
                            0 & 0 & 0 & 0 & 0 & 1 }, $$
where $(\mbox{Ad}\beta){\mit\Phi} = \beta{\mit\Phi}\beta^{-1}$, ${\mit\Phi} \in {\gge_7}^C$. It is easy to check that $\wti{\beta} \in ({E_8}^C)_{1,1^-,1_-}$. Conversely, if $\alpha \in {E_8}^C$ satisfies $ \alpha \1 = \1, \alpha 1^- = 1^-$ and $\alpha 1_- = 1_-$, then $\alpha$ is of the form$$
  \alpha = \pmatrix{\beta_1 & \beta_{12} & \beta_{13} & 0 & 0 & 0 
\vspace{0.5mm}\cr
                    \beta_{21} & \beta_2 & \beta_{23} & 0 & 0 & 0 
\vspace{0.5mm}\cr        
                    \beta_{31} & \beta_{32} & \beta_3 & 0 & 0 & 0 
\vspace{0.5mm}\cr
                    a_1 & b_1 & c_1 & 1 & 0 & 0 
\vspace{0.5mm}\cr
                    a_2 & b_2 & c_2 & 0 & 1 & 0 
\vspace{0.5mm}\cr
                    a_3 & b_3 & c_3 & 0 & 0 & 1}, \quad
\begin{array}{l}
      \beta_1 \in \Hom_C({\gge_7}^C),
\vspace{0.5mm}\cr
      \beta_2, \beta_3, \beta_{23}, \beta_{32} \in \Hom_C(\gP^C),
\vspace{0.5mm}\cr
      \beta_{21}, \beta_{31} \in \Hom_C({\gge_7}^C, \gP^C),
\vspace{0.5mm}\cr
      \beta_{12}, \beta_{13} \in \Hom_C(\gP^C, {\gge_7}^C),
\vspace{0.5mm}\cr
      a_i \in \Hom_C({\gge_7}^C, C),
\vspace{0.5mm}\cr 
      b_i, c_i \in \Hom_C(\gP^C, C).
\vspace{0.5mm}\cr
\end{array}$$
From the relation $[\alpha{\mit\Phi}, \1] = [\alpha{\mit\Phi}, \alpha\1] = \alpha[{\mit\Phi}, \1] = 0$, that is,
\begin{eqnarray*}
    0 \!\!\! &= [(\beta_1{\mit\Phi}, \beta_{21}{\mit\Phi}, \beta_{31}{\mit\Phi}, a_1{\mit\Phi}, a_2{\mit\Phi}, a_3{\mit\Phi}), (0, 0, 0, 1, 0, 0)]
\vspace{1mm}\\
             &= (0, - \beta_{21}{\mit\Phi}, \beta_{31}{\mit\Phi}, 0, -2a_2{\mit\Phi}, 2a_3{\mit\Phi}),
\end{eqnarray*}
we obtain $\beta_{21} = \beta_{31} = 0$ and $a_2 = a_3 = 0$. Furthermore, from $[\alpha{\mit\Phi}, \1^-] = [\alpha{\mit\Phi}, \alpha 1^-] = \alpha[{\mit\Phi}, 1^-] = 0$, that is, 
$$
   0 = [(\beta_1{\mit\Phi}, 0, 0, a_1{\mit\Phi}, 0, 0), (0, 0, 0, 0, 1, 0)] = (0, 0, 0, 0, 2a_1{\mit\Phi}, 0), $$
we obtain $a_1 = 0$. Using the fact that $[\alpha P^-, \1] = [\alpha P^-, \alpha 1] = \alpha[P^-, \1] = - \alpha P^-$, that is,
$$
\begin{array}{l}
    - (\beta_{12}P, \beta_2P, \beta_{32}P, b_1P, b_2P, b_3P)
\vspace{1mm}\\ 
     \quad = [(\beta_{12}P, \beta_2P, \beta_{32}P, b_1P, b_2P, b_3P), (0, 0, 0, 1, 0, 0)]
\vspace{1mm}\\ 
     \quad = (0, - \beta_2P, \beta_{32}P, 0, -2b_2P, 2b_3P),
\end{array} $$
we obtain $\beta_{12} = \beta_{32} = 0$ and $b_1 = b_2 = b_3 = 0$. Similarly, from $[\alpha Q_-, \1] = [\alpha Q_-, \alpha \1] = \alpha[Q_-, \1] = \alpha Q_-$, we obtain $\beta_{13} = \beta_{23} = 0$ and $c_1 = c_2 = c_3 = 0$. Thus we have seen that $\alpha$ is of the form
$$
      \alpha = \pmatrix{ \beta_1 & 0 & 0 & 0 \cr
                         0 & \beta_2 & 0 & 0 \cr
                         0 & 0 & \beta_3 & 0 \cr
                         0 & 0 & 0 & E}. $$
By applying $\alpha$ on $[(0, P, 0, 0, 0, 0), (0, 0, Q, 0, 0, 0)] = \Big(P \times Q, 0, 0, - \dfrac{1}{8}\{P, Q \}, 0, 0\Big)$, we obtain
$$
\displaylines{\hfill    
     \beta_1(P \times Q) = \beta_2P \times \beta_3Q, \quad \{\beta_2P, \beta_3Q \} = \{P, Q \},
\hfill \mbox{(i)}}$$
since $[(0, \beta_2P, 0, 0, 0, 0), (0, 0, \beta_3Q, 0, 0, 0)] = \Big(\beta_1(P \times Q), 0, 0, - \dfrac{1}{8}\{P, Q \}, 0, 0 \Big)$. Again, by applying $\alpha$ on $[(0, P, 0, 0, 0, 0), (0, Q, 0, 0, 0, 0)] = \Big(0, 0, 0, \dfrac{1}{4}\{P, Q \}, 0, 0 \Big)$, we obtain
$$
\displaylines{\hfill
       \{\beta_2P, \beta_2Q \} = \{ P, Q \}. $$
\hfill \mbox{(ii)}}$$
Further, by applying $\alpha$ on $[({\mit\Phi}, 0, 0, 0, 0, 0), (0, P, 0, 0, 0, 0)] = (0, {\mit\Phi}P, 0, 0, 0, 0)$, we obtain
$$
\displaylines{\hfill
        (\beta_1{\mit\Phi})(\beta_2P) = \beta_2({\mit\Phi}P).
\hfill \mbox{(iii)}}$$
From (i), (ii), we have
$$
   \{\beta_2P, \beta_3Q \} = \{\beta_2P, \beta_2Q \}, \quad P, Q \in \gP^C, $$
hence $\beta_2 = \beta_3$, which we denote by $\beta$. If we put $\beta^{-1}P$ in (iii) in the place of $P$, we obtain 
$$
        \beta_1{\mit\Phi} = \beta{\mit\Phi}\beta^{-1}. $$
Therefore, from (i), we have 
$$
        \beta(P \times Q)\beta^{-1} = \beta P \times \beta Q, $$
which implies that $\beta \in {E_7}^C$. Thus the proof of Theorem 5.4.1 is completed.
\vspace{3mm}

{\bf Definition.} We define a $C$-linear mapping $\upsilon : {\gge_8}^C \to {\gge_8}^C$ by
$$
       \upsilon({\mit\Phi}, P, Q, r, s, t) = ({\mit\Phi}. - P, - Q, r, s, t). $$Then $\upsilon \in E_8$ and $\upsilon^2 = 1$. Note that $\upsilon$ is the central element $-1$ of $E_7$ regarding as an element of $E_8$ (see Theorem 5.7.3).
\vspace{2mm}

We shall study the following subgroup $(E_8)^{\upsilon}$ of $E_8$:
$$
        (E_8)^{\upsilon} = \{ \alpha \in E_8 \, | \, \upsilon\alpha = \alpha\upsilon \}. $$

{\bf Lemma 5.7.2.} {\it If $\alpha \in E_8$ satisfies $\alpha 1_- = 1_-$, then it also satisfies $\alpha 1 = 1$ and} $\alpha 1^- = 1^-$.
\vspace{2mm}

{\bf Proof.} Let $\alpha\1 = ({\mit\Phi}, P, Q, r, s, t)$. From the relation $[\alpha\1, 1_-] = [\alpha\1, \alpha 1_-] = \alpha[\1, 1_-] = -2\alpha 1_- = -21_-$. we have
$$
   -21_- = [({\mit\Phi}, P, Q, r, s, t), 1_-] = (0, 0, - P, s, 0, -2r), $$
from which we obtain $P = 0, s = 0, r = 1$. Further, $\langle \alpha\1, \alpha\1 \rangle  = \langle \1, \1 \rangle = 8$, that is, $\langle {\mit\Phi}, {\mit\Phi} \rangle + \langle Q, Q \rangle + 8 + 4(\tau t)t = 8$, which implies that ${\mit\Phi} = 0, Q = 0, t = 0$. Therefore $\alpha\1 = \1$. By applying $\alpha$ on $[1^-, 1_-] = \1$, we obtain $\alpha 1^- = 1^-$ by a similar method to the above.
\vspace{2mm}

To find the structure of the group $(E_8)^{\upsilon}$, we shall first study the following subgroup $(E_8)_{1_-}$ of $E_8$:
$$
      (E_8)_{1_-} = \{ \alpha \in E_8 \, | \, \alpha 1_- = 1_- \}. $$

{\bf Theorem 5.7.3.} \qquad \qquad \quad $(E_8)_{1_-} \cong E_7$.
\vspace{2mm}

{\bf Proof.} $(E_8)_{1_-} = (E_8)_{1,1^-,1_-}$ (Lemma 5.7.2)
\vspace{1mm}

\qquad \quad $= \{ \alpha \in ({E_8}^C)_{1,1^-,1_-} \, | \, \tau\wti{\lambda}\alpha = \alpha\tau\wti{\lambda} \}$
\vspace{1mm}

\qquad \quad $= \{ \alpha \in {E_7}^C \, | \, \tau\wti{\lambda}\alpha = \alpha\tau\wti{\lambda} \}$ (Proposition 5.7.1)
\vspace{1mm}

\qquad \quad $= \{ \alpha \in {E_7}^C \, | \, \tau\lambda\alpha = \alpha\tau\lambda \}$ (by the correspondence to Proposition 5.7.1)
\vspace{1mm}

\qquad \quad $= E_7$ (Lemma 4.3.3.(4)).
\vspace{3mm}

{\bf Remark.} We define a space $\gW_1$ by
$$
        \gW_1 = \{ R \in {\gge_8}^C \, | \, R \times R = 0, \langle R, R \rangle = 4 \} $$
(see Remark of Theorem 5.4.1), then, we obtain a homeomorphism
$$
            E_8/E_7 \simeq \gW_1. $$
(See Yokota, Imai and Yasukura [53]).
\vspace{3mm}

{\bf Theorem 5.7.4.} {\it The group} $(E_8)^{\upsilon}$ {\it contains a subgroup}
$$
        \varphi_3(SU(2)) = \{ \varphi_3(A) \in E_8 \, | \, A \in SU(2) \} $$
{\it which is isomorphic to the group $SU(2) = \{ A \in M(2, C) \, | \, (\tau\,{}^t\!A)A = E, \det A = 1 \}$, where, for $A = \pmatrix{a & - \tau b \cr
                                                 b & \tau a} \in SU(2)$,
$\varphi_3(A) : {\gge_8}^C\to {\gge_8}^C$ is defined by}

$$
   \varphi_3(A) = \pmatrix{1 & 0 & 0 & 0 & 0 & 0 
\vspace{1mm}\cr
                   0 & a1 & -\tau b1 & 0 & 0 & 0 
\vspace{1mm}\cr
                   0 & b1 & \tau a1 & 0 & 0 & 0 
\vspace{1mm}\cr
                   0 & 0 & 0 & (\tau a)a - (\tau b)b & -(\tau a)b & a(\tau b) 
\vspace{1mm}\cr
                   0 & 0 & 0 & 2a(\tau b) & a^2 & -(\tau b)^2 
\vspace{1mm}\cr
                   0 & 0 & 0 & 2(\tau a)b & -b^2 & (\tau a)^2}. $$
\vspace{2mm}

{\bf Proof.} For $A = \pmatrix{a & -\tau b \cr
                               b & \tau a} =
     \exp\pmatrix{- i\nu & - \tau\rho \cr
                  \rho & i\nu} \in SU(2)$, we have $\phi(A) = \exp({\mit\Theta}(0, 0,$ $ 0, i\nu, \rho, -\tau\rho)) \in (E_8)^{\upsilon}$.
\vspace{3mm}

{\bf Lemma 5.7.5.} {\it The group $(E_8)^{\upsilon}$ is connected}.
\vspace{1mm}

{\bf Proof.} $(E_8)^{\upsilon}$ is connected as a fixed points subgroup of the involutive automorphism $\upsilon$ of the simply connected Lie group $E_8$. 
\vspace{3mm}

{\bf Theorem 5.7.6.}  \quad $(E_8)^{\upsilon} \cong (SU(2) \times E_7)/\Z_2, \; \Z_2 = \{ (E, 1), (-E, -1) \}.$
\vspace{2mm}

{\bf Proof.} We define a mapping $\varphi : SU(2) \times E_7 \to (E_8)^{\upsilon}$ by
$$
             \varphi(A, \beta) = \varphi_3(A)\beta. $$
Since $\varphi_3(A) \in \varphi_3(SU(A))$ and $\beta \in E_7$ commute, $\varphi$ is a homomorphism. Since $(E_8)^{\upsilon}$ is connected (Lemma 5.7.5), to prove that $\varphi$ is onto, it is sufficient to show that the differential mapping $\varphi_* : \su(8) \oplus \gge_7 \to (\gge_8)^{\upsilon}$ of $\varphi$ is onto. But which is not difficult to see. Indeed, we have
\begin{eqnarray*}
    (\gge_8)^{\upsilon} \!\!\! &=& \!\!\! \{{\mit\Theta}(R) \in {\mit\Theta}(\gge_8) \, | \, \upsilon{\mit\Theta}(R) = {\mit\Theta}(R)\upsilon \}
\cong \{ R \in \gge_8 \, | \, \upsilon R = R \}
\vspace{1mm}\\
   \!\!\! &=& \!\!\! \{({\mit\Phi}, 0, 0, r, s, -\tau s) \, | \, {\mit\Phi} \in \gge_7, r \in i\R, s \in C \}.
\end{eqnarray*}
$\Ker\varphi = \{(E, 1), (-E, -1) \} = \Z_2$ is easily obtained. Thus we have the isomorphism $(SU(2) \times E_7)/\Z_2 \cong (E_8)^{\upsilon}$.
\vspace{2mm}

{\bf Remark.} We can prove directly that $\varphi$ is onto without using the connectedness of $(E_8)^{\upsilon}$ (Lemma 5.7.5), (see Imai and Yokota [13]).
\vspace{4mm}

{\bf 5.8. Involution $\wti{\lambda}\gamma$ and subgroup $Ss(16)$ of $E_8$}
\vspace{3mm}

We define a $C$-linear mapping $\wti{\lambda}\gamma : {\gge_8}^C \to {\gge_8}^C$ by
$$
 \wti{\lambda}\gamma({\mit\Phi}, P, Q, , s, t) =
      (\lambda\gamma{\mit\Phi}\gamma\lambda^{-1}, \lambda\gamma Q, - \lambda\gamma P, -r, -t, -s). $$
Then $\wti{\lambda}\gamma \in E_8$ and $(\wti{\lambda}\gamma)^2 = 1$. 
\vspace{2mm}

We shall study the following subgroup $(E_8)^{\wti{\lambda}\gamma}$ of $E_8$:
\begin{eqnarray*}
     (E_8)^{\wti{\lambda}\gamma} \!\!\! &=& \!\!\! \{ \alpha \in E_8 \, | \, \wti{\lambda}\gamma\alpha = \alpha\wti{\lambda}\gamma \}
\vspace{1mm}\\
      \!\!\! &=& \!\!\! \{ \alpha \in E_8 \, | \, \tau\gamma\alpha = \alpha\tau\gamma \} = (E_8)^{\tau\gamma}.
\end{eqnarray*}

We define an $\R$-linear mapping $l : M(8, C) \to M(16, \R)$ by
$$
    l\Big((x_{kl} + iy_{kl}) \Big) = \Big(\pmatrix{x_{kl} & y_{kl} \cr
                                                   -y_{kl} & x_{kl}} \Big), \quad x_{kl},y_{kl} \in \R. $$
Further, we define $I, J \in M(16, \R)$ by
$$
  I = \diag(I, \cdots, I), \;\; I = \pmatrix{1 & 0 \cr
                                             0 & -1},
\quad
  J = \diag(J, \cdots, J), \;\; J = \pmatrix{0 & 1 \cr
                                             -1 & 0}, $$
then $IJ = - JI$, and for $X, Y \in M(8, C)$, we have
\vspace{1mm}

$\begin{array}{l}
\mbox{(1)} \quad l(XY) = l(X)l(Y),
\vspace{1mm}\\

\mbox{(2)} \quad Il(X) = l(\tau X)I, \quad Jl(X) = l(X)J, 
\vspace{1mm}\\

\mbox{(3)} \quad {}^tl(X) = l(\tau {}^t\!X).
\end{array}$
\vspace{3mm}

$\begin{array}{ll}
\mbox{{\bf Lemma 5.8.1.}}
\;\;(1) & l(\gu(8)) = \{ B \in \so(16) \, | \, JB = BJ \},
\vspace{1mm}\\
    & l(\gS(8, C))I = \{ B \in \so(16) \, | \, JB = - BJ \}.
\end{array}$
\vspace{1mm}

(2) {\it Any element $B \in \so(16)$ is uniquely expressed by}
$$
\begin{array}{lll}
  B \!\!\! &= l(D') + l(S)I, & D' \in \gu(8), S \in \gS(8, C)
\vspace{1mm}\\
   \!\!\!&= l(D) + l(S)I + l(icE),  & D \in \su(8), S \in \gS(8, C), c \in \R.
\end{array}$$   

{\bf Proof.} (1) If $D \in \gu(8)$, then we have $Jl(D) = l(D)J$ and $^tl(D) = l(\tau\,^t\!D) = - l(D)$. Conversely, suppose that $B \in \so(16)$ satisfies $JB = BJ$. Let $B = l(D), D \in M(8, C)$. Then, the relation 
$$
        l(-D) = - B = {}^tB = {}^tl(D) = l(\tau \,{}^tD) $$
implies $\tau{}\,^tD = -D$, that, is, $D \in \gu(8)$.  Next, for $S \in \gS(8, C)$, we have
$$
\begin{array}{l}
     Jl(S)I = l(S)JI = - l(S)IJ,
\vspace{1mm}\\
     {}^t(l(S)I) = {}^t\,I^t(l(S)) = Il(\tau\,{}^tS) = - Il(\tau S) = - l(S)I.
\end{array} $$
Conversely, suppose that $B \in \so(16)$ satisfies $JB = - BJ$. Since the element $BI$ satisfies $JBI = BIJ$, we let $BI = l(S), S \in M(8, C)$. Then the relation
$$
    l(-S)I = - B = {}^tB = {}^tI\,{}^t(l(S)) = Il(\tau\,{}^tS) = l({}^tS)I $$
implies $-S = {}^tS$, that is, $S \in \gS(8, C)$.
\vspace{2mm}

(2) Let $B = \dfrac{B - JBJ}{2} + \dfrac{B + JBJ}{2}$ and use (1) above.
\vspace{2mm}

The following Lemmas 5.8.2 and 5.8.3 are properties of the mappings $\chi : (\gP^C)_{\tau\gamma} \to \gS(8, C)$ and $\varphi_* : \su(8) \to (\gge_6)^{\lambda\gamma}$ of Section 4.12, and will be used in the proof of Theorem 5.8.4.
\vspace{3mm}

{\bf Lemma 5.8.2} {\it The Lie isomorphism $\varphi_* : \sp(4) \to (\gge_6)^{\lambda\gamma}$ defined by $(\varphi_*D)X = g^{-1}(D(gX) + (gX)D^*), X \in \gJ^C$ of Theorem} 3.1.2 {\it satisfies}
$$
   \varphi_*([gX_1, gX_2]) = 2(X_1 \vee \gamma X_2 - X_2 \vee \gamma X_1), \quad X_1, X_2 \in \gJ^C. $$

{\bf Proof.} Note that $[gX_1, gX_2] \in \sp(4)$. Now, since
$$
\begin{array}{l}
     g(2(X_1 \vee \gamma X_2)X) \qquad X \in \gJ^C
\vspace{1mm}\\
     = g((\gamma X_2, X)X_1 + \dfrac{1}{3}(X_1, \gamma X_2)X - 4\gamma X_2 \times (X_1 \times X))\; \mbox{(Lemma 3.4.1)}
\vspace{1mm}\\
     =(\gamma X_2, X)gX_1 + \dfrac{1}{3}(X_1, \gamma X_2)gX - 4gX_2 \circ (gX_1 \circ gX) + (\gamma X_1, X)gX_2
\vspace{1mm}\\
\qquad + (\gamma X_2, \gamma X_1, \gamma X)E \;\;\mbox{(Lemma 3.12.1)},
\end{array} $$
we have
\vspace{2mm}

\qquad \qquad 
   $g(2(X_1 \vee \gamma X_2 - X_2 \vee \gamma X_1)X)$
\vspace{1mm}

\qquad \qquad \quad
    $= - 4gX_2 \circ(gX_1 \circ gX) + 4gX_1 \circ (gX_2 \circ gX)$
\vspace{1mm}

\qquad \qquad \quad
   $= - gX_2gX_1gX - gX_2gXgX_1 - gX_1gXgX_2 - gXgX_1gX_2$
\vspace{1mm}

\qquad \qquad \quad 
   \quad $+ gX_1gX_2gX + gX_1gXgX_2 + gX_2gXgX_1 + gXgX_2gX_1$
\vspace{1mm}

\qquad \qquad \quad
     $= [gX_1, gX_2]gX - gX[gX_1, gX_2] = g((\varphi_*[gX_1, gX_2])X)$.
\vspace{2mm}

\noindent Consequently, since $g$ is injective, we have the lemma.
\vspace{3mm}

{\bf Lemma 5.8.3.} {\it For $S, S_1, S_2 \in \gS(8, C)$, we have}
\vspace{1mm}

(1) $\quad \lambda\gamma\chi^{-1}(S) = - \chi^{-1}(iS).$
\vspace{1mm}

(2) $\quad \tr(S_1\tau S_2 - S_2\tau S_1) = 4i\{\chi^{-1}S_1, \chi^{-1}S_2 \}.$
\vspace{1mm}

(3) $\quad \varphi_*\Big((S_1\tau S_2 - S_2\tau S_1) - \dfrac{1}{8}\tr(S_1\tau S_2 - S_2\tau S_1)E \Big)$
\vspace{1mm}

\qquad \qquad \qquad \quad
$= 4(\lambda\gamma\chi^{-1}S_1 \times \chi^{-1}S_2 - \lambda\gamma\chi^{-1}S_2 \times \chi^{-1}S_1).$
\vspace{2mm}

{\bf Proof.} (1) Let $\chi^{-1}S = P = (X, Y, \xi, \eta)$. Then
\begin{eqnarray*}
   \chi\lambda\gamma\chi^{-1}S \!\!\! &=& \!\!\! \chi\lambda\gamma(X, Y, \xi,\eta) = \chi(\gamma Y, - \gamma X, \eta, - \xi)
\vspace{1mm}\\
      \!\!\! &=& \!\!\! \Big(k\Big(g(\gamma Y) - \dfrac{\eta}{2}E \Big) + ik\Big(g\Big(\gamma(- \gamma X) - \dfrac{- \xi}{2}E \Big)\Big)\Big)J
\vspace{1mm}\\
      \!\!\! &=& \!\!\! - i\Big(k\Big(gX - \dfrac{\xi}{2}E \Big) + ik\Big(g(\gamma Y) - \dfrac{\eta}{2}E \Big)\Big)J
\vspace{1mm}\\
      \!\!\! &=& \!\!\! - i\chi(X, Y, \xi, \eta) = - i\chi P = - iS.
\end{eqnarray*}

(2),(3) Let $\chi^{-1}S_i = P_i = (X_i, T_i, \xi_i, \eta_i), i = 1, 2.$  Noting
$$
    S_1\tau S_2 - S_2\tau S_1 - \dfrac{1}{8}(S_1\tau S_2 - S_2\tau S_1) \in \su(8), $$
we have
\vspace{2mm}

\quad
   $S_1\tau S_2 = \chi P_1\tau\chi P_2 = \chi(X_1, Y_1, \xi_1, \eta_1)\tau\chi(X_2, Y_2, \xi_2, \eta_2)$
\vspace{1mm}

\qquad
   $= \Big(k\Big(gX_1 - \dfrac{\xi_1}{2}E \Big) + ik\Big(g(\gamma Y_1) - \dfrac{\eta_1}{2}E \Big)\Big)J\tau\Big(\Big(k\Big(gX_2 - \dfrac{\xi_2}{2}E \Big)$ 
\vspace{1mm}

\qquad \qquad \qquad \qquad \qquad \qquad \qquad \qquad \qquad 
  $+ ik\Big(g(\gamma Y_2) - \dfrac{\eta_2}{2}E \Big)\Big)J\Big)$
\vspace{1mm}

\qquad 
  $= \Big(k\Big(gX_1 - \dfrac{\xi_1}{2}E \Big) + ik\Big(g(\gamma Y_1) - \dfrac{\eta_1}{2}E \Big)\Big)\Big(k\Big(gX_2 - \dfrac{\xi_2}{2}E \Big)$ 
\vspace{1mm}

\qquad \qquad \qquad \qquad \qquad \qquad \qquad \qquad \qquad 
   $- ik\Big(g(\gamma Y_2) - \dfrac{\eta_2}{2}E \Big)\Big)J^2$
\vspace{1mm}

\qquad 
   $= - k\Big(\Big(gX_1 - \dfrac{\xi_1}{2}E \Big)\Big(gX_2 - \dfrac{\xi_2}{2}E \Big) + \Big(g(\gamma Y_1) - \dfrac{\eta_1}{2}E \Big)\Big(g(\gamma Y_2) - \dfrac{\eta_2}{2}E\Big)\Big)$
\vspace{1mm}

\qquad $\;\;\;$
   $ -ik\Big(\Big(g(\gamma Y_1) - \dfrac{\eta_1}{2}E \Big)\Big(gX_2 - \dfrac{\xi_2}{2}E \Big) - \Big(gX_1 - \dfrac{\xi_1}{2}E \Big)\Big(g(\gamma Y_2) - \dfrac{\eta_2}{2}E \Big)\Big)$.
\vspace{2mm}

\noindent Therefore
\vspace{2mm}

\quad
  $S_1\tau S_2 - S_2\tau S_1$
\vspace{1mm}

  \qquad \quad
  $= -k(gX_1gX_2 - gX_2gX_1 + g(\gamma Y_1)g(\gamma Y_2) - g(\gamma Y_2)g(\gamma Y_1))$
\vspace{1mm}

   \qquad \qquad  
   $- ik(g(\gamma Y_1)gX_2 - g(\gamma Y_2)gX_1 - gX_1g(\gamma Y_2) + gX_2g(\gamma Y_1))$
\vspace{1mm}

   \qquad \qquad 
    $- \eta_1gX_2 + \eta_2 gX_1 + \xi_1g(\gamma Y_2) - \xi_2g(\gamma Y_1)) + \dfrac{i}{2}(\xi_1\eta_2 - \eta_2\xi_1)E$
\vspace{1mm}

   \qquad \quad
   $= k(- [gX_1, gX_2] - [g(\gamma Y_1),g(\gamma Y_2)])$
\vspace{1mm}

   \qquad \qquad 
   $+ ik(g(2\gamma X_1 \times Y_2 - 2\gamma X_2 \times Y_1 + \eta_1X_2 - \eta_2X_1 - \xi_1\gamma Y_2 + \xi_2\gamma Y_1))$
\vspace{1mm}

\qquad \qquad 
   $+ \dfrac{i}{2}((X_1, Y_2) - (X_2, Y_1) + \xi_1\eta_2 - \xi_2\eta_1))E \;\; \mbox{(Lemma 3.12.1)}$
\vspace{1mm}\\
\mbox{$\begin{array}{ll}
   \mbox{(denote} &D = - [gX_1, gX_2] - [g(\gamma Y_1), g(\gamma Y_2)] \in \sp(4), 
\vspace{1mm}\\
          &A = 2\gamma X_1 \times Y_2 - 2\gamma X_2 \times Y_1 + \eta_1X_2 - \eta_2X_1 - \xi_1\gamma X_2 + \xi_2\gamma Y_1 \in \gJ^C) 
\end{array}$}
\vspace{1mm}

   \qquad \quad
   $= kD + ik(gA) + \dfrac{i}{2}\{P_1, P_2\}E$.
\vspace{2mm}

\noindent Taking the trace of both sides, we obtain
$$
    \tr(S_1\tau S_2 - S_2\tau S_1) = 4i\{P_1, P_2 \} = 4i\{\chi^{-1}S_1, \chi^{-1}S_2 \}. $$
and the expression above equal to
$$
    S_1\tau S_2 - S_2\tau S_1 - \dfrac{1}{8}\tr(S_1\tau S_2 - S_2\tau S_1) = kD + ik(gA), $$
On the other hand, we have
\vspace{1mm}

$\qquad \qquad   \lambda\gamma P_1 \times P_2 = \pmatrix{\gamma Y_1 \cr
                               - \gamma X_1 \cr
                               \eta_1 \cr
                               - \xi_1} 
                \times \pmatrix{X_2 \cr Y_2 \cr
                                \xi_2 \cr \eta_2}$
\vspace{1mm}

$\qquad \qquad \qquad \qquad = {\mit\Phi}\pmatrix{ - \dfrac{1}{2}(\gamma Y_1 \vee Y_2 - X_2 \vee \gamma X_1)
\vspace{1mm}\cr
         \dfrac{1}{4}(2\gamma X_1 \times Y_2 + \eta_1X_2 + \xi_2\gamma Y_1)      \vspace{1mm}\cr
         \dfrac{1}{4}(2\gamma Y_1 \times X_2 + \xi_1Y_2 + \eta_2\gamma X_1)     \vspace{1mm}\cr
    \dfrac{1}{8}((\gamma Y_1, Y_2) - (X_2, \gamma X_1) - 3 (\eta_1\eta_2 - \xi_2\xi_1))}.$
\vspace{1mm}

\noindent Therefore,

$$
\begin{array}{l}
   \lambda\gamma P_1 \times P_2 - \lambda\gamma P_2 \times P_1
\vspace{1mm}\\
\qquad = {\mit\Phi}\pmatrix{
   \dfrac{1}{2}(- X_1 \vee \gamma X_2 + X_2 \vee \gamma X_1 - \gamma Y_1 \vee Y_2 + \gamma Y_2 \vee Y_1)
\vspace{1mm}\cr
   \dfrac{1}{4}(2\gamma X_1 \times Y_2 - 2\gamma X_2 \times Y_1 + \eta_1X_2 - \eta_2X_1 + \xi_2\gamma Y_1 - \xi_1\gamma Y_2)
\vspace{1mm}\cr
   \dfrac{1}{4}(2X_2 \times \gamma Y_1 - 2X_1 \times \gamma Y_2 + \eta_2\gamma X_1 - \eta_1\gamma X_2 + \xi_1Y_2 - \xi_2Y_1)
\vspace{1mm}\cr
        0}
\vspace{1mm}\\
\qquad = {\mit\Phi}\Big(\dfrac{1}{4}\varphi_*D, \dfrac{1}{4}A, - \dfrac{1}{4}\gamma A, 0 \Big) \;\;\mbox{(Lemma 5.8.2)}
\vspace{1mm}\\
\qquad = \dfrac{1}{4}\varphi_*(kD + ik(gA))
\vspace{1mm}\\
\qquad = \dfrac{1}{4}\varphi_*\Big(S_1\tau S_2 - S_2\tau S_1 - \dfrac{1}{8}\tr(S_1\tau S_2 - S_2\tau S_1)E \Big).
\end{array}$$
\vspace{0.5mm}

{\bf Theorem 5.8.4.} {\it The Lie algebra $(\gge_8)^{\wti{\lambda}\gamma}$ of the Lie group $(E_8)^{\wti{\lambda}\gamma}$ is given by}
\begin{eqnarray*}
     (\gge_8)^{\wti{\lambda}\gamma} \!\!\! &=& \!\!\!
        \{ {\mit\Theta} \in {\mit\Theta}({\gge_8}^C) \, | \, \wti{\lambda}\gamma{\mit\Theta} = {\mit\Theta}\wti{\lambda}\gamma \}
\vspace{1mm}\\
     \!\!\! &=&\!\!\! \{{\mit\Theta}({\mit\Phi}, P, - \lambda\gamma P, 0, s, - s) \, | \, {\mit\Phi} \in (\gge_7)^{\lambda\gamma}, P \in (\gP^C)_{\lambda\gamma}, s \in \R \}
\end{eqnarray*}
{\it and $(\gge_8)^{\wti{\lambda}\gamma}$ is isomorphic to the Lie algebra $\so(16)$ by the mapping $\zeta : \so(16) \to (\gge_8)^{\wti{\lambda}\gamma}$ defined by}
$$
    \zeta(l(D) + l(S)I + l(icE)) = {\mit\Theta}(\varphi_*D, 2\lambda\gamma\chi^{-1}S, 2\chi^{-1}S, 0, 2c, -2c), $$
{\it where} $D \in \su(8), S \in \gS(8, C), c \in \R,\; \mbox{and}\; \varphi_* : \su(8) \to (\gge_7)^{\lambda\gamma}, \chi : (\gP^C)_{tau\gamma} \to \gS(8, C)$ {\it are mappings defined in Section} 4.12.
\vspace{2mm}

{\bf Proof.} It is not difficult to see that the first half of the theorem and that $\zeta$ is onto. We will prove that $\zeta$ preserve the Lie bracket.

\vspace{2mm}
$\begin{array}{ll}
(1) \quad \zeta[l(D_1), l(D_2)] = \zeta l[D_1, D_2]                             \vspace{1mm}\\
    \qquad \qquad = {\mit\Theta}(\varphi_*[D_1, D_2], 0, 0, 0, 0, 0) = {\mit\Theta}([\varphi_*D_1, \varphi_*D_2], 0, 0, 0, 0, 0)
\vspace{1mm}\\
    \qquad \qquad = [{\mit\Theta}(\varphi_*D_1, 0, 0, 0, 0, 0), {\mit\Theta}(\varphi_*D_2, 0, 0, 0, 0, 0)]
\vspace{1mm}\\
    \qquad \qquad = [\zeta l(D_1), \zeta l(D_2)].
\end{array}$                                                                    \vspace{1mm}

$\begin{array}{l}
(2) \quad \zeta[l(D), l(S)I)] = \zeta(l(DS - S\tau D)I) = \zeta(l(DS + S\,{}^tD)I)
\vspace{1mm}\\
\qquad \qquad = {\mit\Theta}(0, 2\lambda\gamma\chi^{-1}(DS + S\,{}^tD), 2\chi^{-1}(DS + S\,{}^tD), 0, 0, 0).
\end{array}$                                                                    \vspace{1mm}\\
\noindent On the other hand,
$$
\begin{array}{l}
   [\zeta l(D), \zeta(l(S)I)] 
\vspace{1mm}\\
\qquad = [{\mit\Theta}(\varphi_*D, 0, 0, 0, 0, 0), {\mit\Theta}(0, 2\lambda\gamma\chi^{-1}S, 2\chi^{-1}S, 0, 0, 0)]
\vspace{1mm}\\
\qquad = {\mit\Theta}(0, 2(\varphi_*D)\lambda\gamma\chi^{-1}S, 2(\varphi_*D)\chi^{-1}S, 0, 0, 0).
\end{array} $$
Since $(\varphi_*D)\lambda\gamma = \lambda\gamma(\varphi_*D)$ and $(\varphi_*D)\chi^{-1}S = \chi^{-1}(DS + S^tD)$, they are equal.
\vspace{2mm}

$\begin{array}{l}
(3) \quad \zeta[l(D), l(icE)] = \zeta[D, icE] = \zeta 0 = 0
\vspace{1mm}\\
\qquad = [{\mit\Theta}(\varphi_*D, 0, 0, 0, 0, 0), {\mit\Theta}(0, 0, 0, 0, 2c, -2c) = [\zeta l(D), \zeta l(icE)].
\end{array}$
\vspace{2mm}

$\;\,$(4) \quad $\zeta[l(S_1)I, l(S_2)I] = \zeta l(S_1\tau S_2 - S_2\tau S_1)$
\vspace{1mm}

\qquad \qquad 
   $= \zeta\Big(l\Big(S_1\tau S_2 - S_2\tau S_1 - \dfrac{1}{8}\tr(S_1\tau S_2 - S_2\tau S_1)E \Big)$
\vspace{1mm}

\qquad \qquad \qquad 
    $+ l\Big(\dfrac{1}{8}\tr(S_1\tau S_2 - S_2\tau S_1)E \Big)\Big)$
\vspace{1mm}

\qquad \qquad  
   $= {\mit\Theta}\Big(\varphi_*\Big(S_1\tau S_2 - S_2\tau S_1 - \dfrac{1}{8}\tr(S_1\tau S_2 - S_2\tau S_1)E \Big), 0, 0, 0,$
\vspace{1mm}

\qquad \qquad \qquad  
    $- \dfrac{i}{4}\tr(S_1\tau S_2 - S_2\tau S_1), \dfrac{i}{4}\tr(S_1\tau S_2 - S_2\tau S_1) \Big)$
\vspace{1mm}                                                                  

\qquad \qquad 
    $= {\mit\Theta}(4(\lambda\gamma\chi^{-1}S_1 \times \chi^{-1}S_2 - \lambda\gamma\chi^{-1}S_2 \times \chi^{-1}S_1), 0, 0, 0,$
\vspace{1mm}

\qquad \qquad \qquad 
       $\{\chi^{-1}S_1, \chi^{-1}S_2\} - \{\chi^{-1}S_1, \chi^{-1}S_2 \})\;\; \mbox{(Lemma 5.8.3)}$.
\vspace{1mm}

\noindent On the other hand,
\vspace{2mm}

   $[\zeta l(S_1)I, \zeta(l(S_2)I]$ \qquad \quad
\vspace{1mm}

\quad 
   $= [{\mit\Theta}(0, 2\lambda\gamma\chi^{-1}S_1, 2\chi^{-1}S_1, 0, 0, 0),
       {\mit\Theta}(0, 2\lambda\gamma\chi^{-1}S_2, 2\chi^{-1}S_2, 0, 0, 0)]$
\vspace{1mm}

\quad  
   $= {\mit\Theta}\Big(2\lambda\gamma\chi^{-1}S_1 \times 2\chi^{-1}S_2 - 2\lambda\gamma\chi^{-1}S_2 \times 2\chi^{-1}S_1, 0, 0,$
\vspace{1mm}

\qquad \qquad    
   $\dfrac{1}{8}(- \{2\lambda\gamma\chi^{-1}S_1, 2\chi^{-1}S_2 \} + \{2\lambda\gamma\chi^{-1}S_2, 2\chi^{-1}S_1 \})$,
\vspace{1mm}

\qquad \qquad 
   $\dfrac{1}{4}\{2\lambda\gamma\chi^{-1}S_1, 2\lambda\gamma\chi^{-1}S_2 \}, - \dfrac{1}{4}\{2\chi^{-1}S_1, 2\chi^{-1}S_2 \} \Big)$,
\vspace{1mm}

\noindent which equals to the above.
\vspace{2mm}

$\begin{array}{l}
(5) \quad \zeta[(icE), l(S)I] = \zeta(2l(icS)I)
\vspace{1mm}\\
\qquad \qquad
= {\mit\Theta}(0, 4\lambda\gamma\chi^{-1}(icS), 4\chi^{-1}(icS), 0, 0, 0)\vspace{1mm}\\
\qquad \qquad
= {\mit\Theta}(0, 4\chi^{-1}(cS), - 4\lambda\gamma\chi^{-1}(cS), 0, 0, 0) \;\;\mbox{(Lemma 5.8.3)}
\vspace{1mm}\\
\qquad 
\qquad= [{\mit\Theta}(0, 0, 0, 0, 2c, - 2c), {\mit\Theta}(0, 2\lambda\gamma\chi^{-1}S, 2\chi^{-1}S, 0, 0, 0)]
\vspace{1mm}\\
\qquad \qquad
= [\zeta l(icE), \zeta(l(S)I)]. 
\end{array}$
\vspace{2mm}

\noindent Finally,
\vspace{2mm}

$\begin{array}{l}
(6) \quad  \zeta[l(ic_1E), l(ic_2E)] = \zeta l[ic_1E, ic_2E]  = \zeta 0 = 0
\vspace{1mm}\\
  \qquad \qquad 
= {\mit\Theta}(0, 0, 0, 0, 2c_1, - 2c_1), {\mit\Theta}(0, 0, 0, 0, 2c_2, - 2c_2))]
\vspace{1mm}\\
   \qquad \qquad
= [\zeta l(ic_1E), \zeta l(ic_2E)].
\end{array}$
\vspace{2mm}

\noindent Thus we have proved Theorem 5.8.4. 
\vspace{2mm}

The group $(E_8)^{\wti{\lambda}\gamma}$ is connected as a fixed point subgroup under the involution $\wti{\lambda}\gamma$ of the simply connected Lie group $E_8$. Therefore, by Theorem 5.8.4, $(E_8)^{\wti{\lambda}\gamma}$ is isomorphic to one of the following groups
$$
     Spin(16), \quad SO(16), \quad Ss(16), \quad SO(16)/\Z_2. $$
Precisely we have $(E_8)^{\wti{\lambda}\gamma} \cong Ss(16),$ below we will give an outline of the proof.
\vspace{2mm}

We will use the Lie algebra
$$
     \gge_{8(8)} = ({\gge_8}^C)^{\tau\gamma} = \{R \in {\gge_8}^C \, | \, \tau\gamma R = R \} $$
(see Section 5.13). Now, consider the eigenspace decomposition of ${\gge_8}^C$ by $\wti{\lambda}\gamma$:
$$
   {\gge_8}^C = ({\gge_8}^C)_{\wti{\lambda}\gamma} \oplus ({\gge_8}^C)_{-\wti{\lambda}\gamma}, $$
$$
\begin{array}{l}
\;\; ({\gge_8}^C)_{\wti{\lambda}\gamma} = \{{\mit\Theta} \in \mbox{Der}({\gge_8}^C) \, | \, \wti{\lambda}\gamma{\mit\Theta} = {\mit\Theta}\wti{\lambda}\gamma \} \cong \{R \in {\gge_8}^C \, | \, \wti{\lambda}\gamma R = R \} = ({\gge_8}^C)^{\wti{\lambda}\gamma}, 
\vspace{1mm}\\
   ({\gge_8}^C)_{-\wti{\lambda}\gamma} = \{{\mit\Theta} \in \mbox{Der}({\gge_8}^C) \, | \, \wti{\lambda}\gamma{\mit\Theta} = - {\mit\Theta}\wti{\lambda}\gamma \} \cong \{R \in {\gge_8}^C \, | \, \wti{\lambda}\gamma R = - R \}, 
\end{array}$$
Since we have
\begin{eqnarray*}
       (({\gge_8}^C)_{\wti{\lambda}\gamma})_{\tau\wti{\lambda}} 
     \!\!\!&=&\!\!\! (({\gge_8}^C)_{\tau\gamma})_{\wti{\lambda}\gamma} 
     = (\gge_{8(8)})_{\wti{\lambda}\gamma} = (\gge_{8(8)})^{\wti{\lambda}\gamma},
\vspace{1mm}\\
       (({\gge_8}^C)_{-\wti{\lambda}\gamma})_{\tau\wti{\lambda}} 
     \!\!\!&=&\!\!\! (({\gge_8}^C)_{\tau\gamma})_{-\wti{\lambda}\gamma} 
     = (\gge_{8(8)})_{-\wti{\lambda}\gamma},
\end{eqnarray*} 
we obtain the following decomposition of $\gge_{8(8)}$:
$$
       \gge_{8(8)} = (\gge_{8(8)})^{\wti{\lambda}\gamma} \oplus (\gge_{8(8)})_{-\wti{\lambda}\gamma}. $$
Since $(\gge_{8(8)})^{\wti{\lambda}\gamma} \cong Ss(16)$ (Theorem 5.8.4), this is the Cartan decompposition of $\gge_{8(8)}$. Since $[(\gge_{8(8)})^{\wti{\lambda}\gamma}, (\gge_{8(8)})_{-\wti{\lambda}\gamma}] \subset (\gge_{8(8)})_{-\wti{\lambda}\gamma}$, we obtain a representation $\varphi$ of $(\gge_{8(8)})^{\wti{\lambda}\gamma}$ to $(\gge_{8(8)})_{-\wti{\lambda}\gamma}$:
$$
  \varphi(R)R_1 = [R, R_1], \quad R \in (\gge_{8(8)})^{\wti{\lambda}\gamma},
                                  R_1 \in (\gge_{8(8)})_{-\wti{\lambda}\gamma}. $$
which is irreducible. (See, for example, (8.5.1) of Goto and Grosshans [11]). Furthermore, the complex representation $\varphi^C$ of $\varphi$ to $(((\gge_{8(8)})_{-\wti{\lambda}\gamma})^C = ({e_8}^C)_{-\wti{\lambda}\gamma}$ is also irreducible, since $(\gge_{8(8)})_{-\wti{\lambda}\gamma}$ is simple (see (8.8.3) of the same book). The following lemma follows from above mentioned results.
\vspace{3mm}

{\bf Lemma 5.8.5.} {\it The representation of the group $(E_8)^{\wti{\lambda}\gamma}$ to $(\gge_{8(8)})_{-\wti{\lambda}\gamma}$ is
\vspace{3mm}
irreducible.}

{\bf Proposition 5.8.6.} {\it The center $z((E_8)^{\wti{\lambda}\gamma})$ of the group $(E_8)^{\wti{\lambda}\gamma}$ is a group of order} 2:
$$
         z((E_8)^{\wti{\lambda}\gamma}) = \{1, \wti{\lambda}\gamma \}. $$

{\bf Proof.} Evidently, $\{1, \wti{\lambda}\gamma \} \subset z((E_8)^{\wti{\lambda}\gamma})$. Conversely, let $\alpha \in z((E_8)^{\wti{\lambda}\gamma})$. Since the representation of $(E_8)^{\wti{\lambda}\gamma}$ to $({\gge_8}^C)_{-\wti{\lambda}\gamma}$ irreducible (Lemma 5.8.5), we see, by using Schur's lemma in the theory of groups, that the action of $\alpha$ on $({\gge_8}^C)_{-\wti{\lambda}\gamma}$ is constant. Therefore, there exists an element $k \in C$ such that
$$
      \alpha R = kR, \quad R \in ({\gge_8}^C)_{-\wti{\lambda}\gamma}. $$
Since the Killing form $B_8(R, R')$ is invariant under $\alpha$: $B_8(\alpha R, \alpha R') = B_8(R, R')$, we have
$$
      k^2B_8(R, R') = B_8(\alpha R, \alpha R') = B_8(R, R'), \quad R, R' \in ({\gge_8}^C)_{-\wti{\lambda}\gamma}, $$
which implies that $k^2 = 1$. By Theorem 5.8.4, we have $({\gge_8}^C)^{\wti{\lambda}\gamma} \cong \so(16, C)$ which is simple, and hence we see that $({\gge_8}^C)^{\wti{\lambda}\gamma}$ is generated by $({\gge_8}^C)^{-\wti{\lambda}\gamma}$:
$$
     ({\gge_8}^C)^{\wti{\lambda}\gamma} = \Big\{\dsum_{k,l}[R_k, R_l] \, | \, R_k, R_l \in  ({\gge_8}^C)_{-\wti{\lambda}\gamma} \Big\}. $$
Consequently, $\alpha$ satisfies $k^21 = 1$ on  $({\gge_8}^C)_{-\wti{\lambda}\gamma}$, that is, the identity mapping. When $k = 1$, we have $\alpha = 1$ ,  when $k = - 1$, we have $ \alpha = \wti{\lambda}\gamma$. Thus we have proved the theorem.
\vspace{2mm}

It follows from Proposition 5.8.6, that $(E_8)^{\wti{\lambda}\gamma}$ is isomorphic to one of the following
$$
         SO(16), \quad Ss(16). $$
There are only two, up to equivalence, complex irreducible representations of the Lie algebra $\so(16)$ of dimension 128. In fact, one can obtain the following table, by calculating the dimension of dominant roots by virtue of Weyl's dimension formula (see (7.5.9) of Goto and Grosshans [11]):
$$
    \matrix{\omega_1 & 2\omega_1 & \omega_2 &  2\omega_2 &  \omega_3 &  \omega_4 &  \omega_5 &  \omega_6 &  \omega_7 &  \omega_8 &  \cdots \cr      
16 &  135 &  120 &  5304 &  560 &  1820 & 4368 &  8008 &  128 &  128 &  \cdots}
$$ 
hence the dominant root of dimension 128 is either $\omega_7$ or $\omega_8$. (Here, $\omega_1, \omega_2, \cdots, \omega_8$ are fundamental weights). On the other hand, $Spin(16)$ has two complex irreducible representation ${{\mit\Delta}_{16}}^+$ and ${{\mit\Delta}_{16}}^-$, called spinor representations. Furthermore, both of ${{\mit\Delta}_{16}}^+$ and ${{\mit\Delta}_{16}}^-$ are not representation of $SO(16)$. Now, by Lemma 5.8.5, $(E_8)^{\wti{\lambda}\gamma}$ has a complex irreducible representation $({\gge_8}^C)^{-\wti{\lambda}\gamma}$ of dimension 128, which implies that $(E_8)^{\wti{\lambda}\gamma}$ is not $SO(16)$. So $(E_8)^{\wti{\lambda}\gamma}$ must be $Ss(16)$. Thus we have proved the following 
\vspace{3mm}
theorem.

{\bf Theorem 5.8.7.} \qquad \qquad \quad $(E_8)^{\wti{\lambda}\gamma} \cong Ss(16).$
\vspace{2mm}

{\bf Remark.} We define an involution $C$-linear transformation $\sigma$ of $\gge_8$ by $$
     \sigma({\mit\Phi}, P, Q, r, s, t) = (\sigma{\mit\Phi}\sigma, \sigma P, \sigma Q, r, s, t). $$
This is the same as $\sigma \in F_4 \subset E_6 \subset E_7 \subset E_8$. We define a subgroup $(E_8)^\sigma$ by
$$
     (E_8)^\sigma = \{\alpha \in E_8 \, | \, \sigma\alpha = \alpha\sigma \}. $$
Then we have
$$
       (E_8)^\sigma \cong Ss(16). $$
Indeed, we can prove that the Lie algebra $(\gge_8)^\sigma$ of the group $(E_8)^\sigma$ is isomorphic to the Lie algebra $\so(16) = \{X \in M(16, \R) \, | \, {}^tX = - X \}$. Hnece as the same the case the group $(E_8)^{\wti{\lambda}\gamma}$, the group $(E_8)^\sigma$ have to the semi-spinor group $Ss(16)$. However the proof of $(E_8)^{\wti{\lambda}\gamma} \cong Ss(16) \cong (E_8)^\sigma$ is not concretely. M. Gomyo [10] find the group $Ss(16)$ explicitly in the group $E_8$ (although the definition of the group $E_8$ is different from $E_8$ in Section 5.5).
\vspace{4mm}

{\bf 5.9. Center $z(E_8)$ of $E_8$}
\vspace{3mm}

{\bf Theorem 5.9.1.} {\it The center $z(E_8)$ of the group $E_8$ is trivial}\,:
$$
                z(E_8) = \{ 1 \}. $$

{\bf Proof.} Let $\alpha \in z(E_8)$. The relation  $\upsilon\alpha = \alpha\upsilon$ implies that $\alpha \in \varphi(SU(2) \times E_7) \cong (SU(2) \times E_7)/\Z_2$ (Theorem 5.7.6), and so $\alpha \in z(\varphi(SU(2) \times E_7))$. Therefore by Theorem 4.9.1 we have
$$
      \alpha = \varphi(E, 1) = 1 \quad \mbox{or} \quad \alpha = \varphi(E, -1) = \upsilon. $$
However $\upsilon \notin z(E_8)$ (Theorem 5.7.6). Hence $\alpha = 1$.

\vspace{4mm}

{\bf 5.10. Automorphism $w$ of order 3 and subgroup $(SU(3) \times E_6)/\Z_3$ of $E_8$}
\vspace{3mm}

In this section (also in the following Sections 5.11 and 5.12), we use the same notation as ${\gge_8}^C, \langle R_1, R_2 \rangle, \tau\wti{\lambda}, w$, even if these are different from those used in the proceeding 
\vspace{2mm}
sections.

We consider a $27 \times 3 = 78$ dimensional $C$-vector space
$$
    (\gJ^C)^3 = \Big\{ \X = \pmatrix{X_1 \cr
                                     X_2 \cr
                                     X_3} \, \Big| \, X_i \in \gJ^C \Big\}.$$
In $(\gJ^C)^3$, we define an inner product $(\X, \Y)$, a Hermitian inner product  $\langle \X, \Y \rangle$, a cross product $\X \times \Y$, an element $\X \cdot \Y$ of $\sl(3, C)$ and an element $\X \vee \Y$ of ${\gge_6}^C$ respectively by
\vspace{2mm}

\qquad \quad
   $(\X, \Y) = (X_1, Y_1) + (X_2, Y_2) + (X_3, Y_3) \in C$,
\vspace{1mm}

\qquad \quad
   $\langle \X, \Y \rangle =  \langle X_1, Y_1 \rangle + \langle X_2, Y_2 \rangle + \langle X_3, Y_3 \rangle \in C$,
\vspace{1mm}

\qquad \quad   
  $\X \times \Y = \pmatrix{X_2 \times Y_3 - Y_2 \times X_3 \cr 
                         X_3 \times Y_1 - Y_3 \times X_1 \cr 
                         X_1 \times Y_2 - Y_1 \times X_2} \in (\gJ^C)^3$,
\vspace{1mm}

\qquad \quad
   $\X \cdot \Y = \pmatrix{(X_1, Y_1) & (X_1, Y_2) & (X_1, Y_3) \cr
                          (X_2, Y_1) & (X_2, Y_2) & (X_2, Y_3) \cr
                          (X_3, Y_1) & (X_3, Y_2) & (X_3, Y_3)} - \dfrac{1}{3}(\X, \Y)E \in \sl(3, C)$, 
\vspace{1mm}

\qquad \quad
   $\X \vee \Y = X_1 \vee Y_1 + X_2 \vee Y_2 + X_3 \vee Y_3 \in {\gge_6}^C$,
\vspace{2mm}

\noindent where $\X = \pmatrix{X_1 \cr
                     X_2 \cr
                     X_3}, \Y = \pmatrix{Y_1 \cr
                                    Y_2 \cr
                                    Y_3} \in (\gJ^C)^3$. Further, for $\phi \in \Hom_C(\gJ^C), D = \big(d_{ij}\big) \in M(3, C)$ and $\X = \pmatrix{X_1 \cr
                     X_2 \cr
                     X_3} \in (\gJ^C)^3$, elements $\phi\X, D\X \in (\gJ^C)^3$ are naturally defined by
$$
   \phi\X = \pmatrix{\phi X_1 \cr
                   \phi X_2 \cr
                   \phi X_3}, \quad
       D\X = \pmatrix{d_{11}X_1 + d_{12}X_2 + d_{13}X_3 \cr  
                d_{21}X_1 + d_{22}X_2 + d_{23}X_3 \cr  
                d_{31}X_1 + d_{32}X_2 + d_{33}X_3}. 
\vspace{3mm}$$

{\bf Theorem 5.10.1.} {\it In an} $8 + 78 + 27 \times 3 + 27 \times 3 = 248$ {\it dimensional $C$-vector space
$$
     {\gge}_8^{\;\,C} = \sl(3, C) \oplus {\gge_6}^C \oplus (\gJ^C)^3 \oplus (\gJ^C)^3, $$
we define a Lie bracket $[R_1, R_2]$ by}
$$
   [(D_1, \phi_1, \X_1, \Y_1), (D_2, \phi_2, \X_2, \Y_2)] = (D, \phi, \X, \Y), $${\it where}
$$
\left\{  
\begin{array}{l}
   D = [D_1, D_2] + \dfrac{1}{4}\X_1 \cdot \Y_2 - \dfrac{1}{4}\X_2 \cdot \Y_1
\vspace{1mm}\\
   \phi = [\phi_1, \phi_2] + \dfrac{1}{2}\X_1 \vee \Y_2 - \dfrac{1}{2}\X_2 \vee \Y_1
\vspace{1mm}\\
   \X = \;\;\;\phi_1\X_2 - \phi_2\X_1 + D_1\X_2 - D_2\X_1 - \Y_1 \times \Y_2
\vspace{1mm}\\
   \Y = - {}^t\phi_1\Y_2 + {}^t\phi_2\Y_1 - {}^tD_1\Y_2 + {}^tD_2\Y_1 + \X_1 \times \X_2,
\end{array} \right. $$
{\it then ${\gge}_8^{\;\,C}$ becomes a $C$-Lie algebra of type $E_8$.}
\vspace{2mm}

{\bf Proof.} Let ${\wti{\gge}_8}^{\;\,C} = {\gge_7}^C \oplus \gP^C \oplus \gP^C \oplus C \oplus C \oplus C$ be the $C$-Lie algebra constructed in Theorem 5.1.1. We define a mapping $f : {\wti{\gge}_8}^{\;\,C} \to {\gge_8}^C$ by
$$
\begin{array}{l}
  f({\mit\Phi}(\phi, A, B, \nu), (X, Y, \xi, \eta), (Z, W, \zeta, \omega), r, s, t))
\vspace{2mm}\\
\qquad
   = \Big(\pmatrix{\dfrac{2}{3}\nu & - \dfrac{1}{2}\xi & \dfrac{1}{2}\zeta 
\vspace{1mm}\cr
         \dfrac{1}{2}\omega & - \dfrac{1}{3}\nu - r & t 
\vspace{1mm}\cr
         \dfrac{1}{2}\eta & s & - \dfrac{1}{3}\nu + r}, \phi, 
         \pmatrix{- 2A \vspace{1mm}\cr
                     Z \vspace{1mm}\cr
                     X},
         \pmatrix{- 2B \vspace{1mm}\cr
                     Y \vspace{1mm}\cr
                     - W} \Big),
\end{array}$$
then we can prove that $f$ is an isomorphism as Lie algebras by straightforward calculations. Thus we have the isomorphism ${\wti{\gge}_8}^{\;\,C} \cong {\gge_8}^C$.
\vspace{2mm}

Using the Killing form of ${\wti{\gge}_8}^{\;\,C}$ which is obtained in Theorem 5.3.2, we see that the Killing form $B_8$ of ${\gge_8}^C$ is given by
$$
   B_8(R_1, R_2) = 60\tr(D_1D_2) + \dfrac{5}{2}B_6(\phi_1,\phi_2) + 15(\X_1, \Y_2)    + 15(\X_2, \Y_1)$$ 
($R_i = (D_i, \phi_i, \X_i, \Y_i) \in {\gge_8}^C$), where $B_6$ is the Killing form of ${\gge_6}^C$. We define a complex conjugate transformation $\tau\wti{\lambda}$ of ${\gge_8}^C$ by
$$
   \tau\wti{\lambda}(D, \phi, \X, \Y) = (-\tau\,{}^t\!D, -\tau{}^t\phi\tau, -\tau\Y, -\tau\X). $$
And we define a Hermitian inner product $\langle R_1, R_2 \rangle$ in ${\gge_8}^C$ by
$$
       \langle R_1, R_2 \rangle = - B_8(R_1, \tau\wti{\lambda}R_2). $$
Then we have
$$
   \langle R_1, R_2 \rangle = 60\tr(D_1(\tau{}^tD_2)) + \dfrac{5}{2}B_6(\phi_1, \tau\,{}^t\!\phi_2\tau) + 15\langle \X_1, \X_2 \rangle + 15\langle \Y_1, \Y_2 \rangle. $$

Now, as in Theorem 5.5.3, we see that
$$
   E_8 = \{\alpha \in \mbox{Aut}({\gge_8}^C) \, | \, \langle \alpha R_1, \alpha R_2 \rangle =  \langle R_1, R_2 \rangle \} $$
is a simply connected compact Lie group of type $E_8$.  
\vspace{2mm}

We define a $C$-linear transformation $w$ of ${\gge_8}^C$ by
$$
   w(D, \phi, \X, \Y) =(D, \phi, \omega\X, \omega^2\Y), $$
where $\omega = - \dfrac{1}{2} + \dfrac{\sqrt{3}}{2}i \in C$. Then $w \in E_8$ and $w^3 = 1$.
\vspace{2mm}

Now, we shall study the following subgroup $(E_8)^w$ of $E_8$:
$$
     (E_8)^w = \{\alpha \in E_8 \, | \, w\alpha = \alpha w \}. $$

{\bf Theorem 5.10.2.} $(E_8)^w \cong (SU(3) \times E_6)/\Z_3, \; \Z_3 = \{(E, 1), (\omega E, \omega^21), (\omega^2E, $ $\omega 1) 
\vspace{2mm}
\}.$

{\bf Proof.} We first define a mapping $\varphi_1 : SU(3) \to (E_8)^w$ by
$$
   \varphi_1(A)(D, \phi, \X, \Y) = (ADA^{-1}, \phi, A\X, {}^t\!A^{-1}\Y). $$
We have to prove that $\varphi_1(A) \in (E_8)^w$. Indeed, since the action of $D_1 = (D_1, 0, 0, 0) \in \su(3) \subset \sl(3, C) \subset {\gge_8}^C$ is given by
$$
   \mbox(\mbox{ad}D_1)(D, \phi, \X, \Y) = ((\mbox{ad}D_1)D, 0, D_1\X, - \,{}^t\!D_1\Y), $$
for $A = \exp D_1$, we have $\varphi_1(A) = \exp(\mbox{ad}(D_1)) \in \mbox{Aut}({\gge_8}^C)$. And from
$$
\begin{array}{l}
   \tr(AD_1\tau{}^t\!A(\tau{}^t(AD_2\tau{}^t\!A))) = \tr(AD_1(\tau{}^t\!D_2)A^{-1}) = \tr(D_1(\tau{}^tD_2)), 
\vspace{1mm}\\
\qquad \qquad \qquad \;\;
   \langle A\X, A\Y \rangle = \langle \X, \Y \rangle,
\end{array}$$ 
we see that $\varphi_1(A) \in E_8$. Evidently, $w\varphi_1(A) = \varphi_1(A)w$, hence, $\varphi_1(A) \in (E_8)^w$. Next, we define a mapping $\varphi_2 : E_6 \to (E_8)^w$ by
$$
   \varphi_2(\alpha)(D, \phi, \X, \Y) = (D, \alpha\phi\alpha^{-1}, \alpha\X, {}^t\alpha^{-1}\Y). $$
We have to prove that $\varphi_2(\alpha) \in (E_8)^w$. Indeed, since the action of an element $\phi' = (0, \phi', 0, 0) \in \gge_6 \subset {\gge_6}^C \subset {\gge_8}^C$ is given by
$$
   (\mbox{ad}\phi')(D, \phi, \X, \Y) = (0, (\mbox{ad}\phi')\phi, \phi'\X, - {}^t\phi'\Y), $$
for $\alpha = \exp\phi'$, we have $\varphi_2(\alpha) = \exp(\mbox{ad}(\phi')) \in \mbox{Aut}({\gge_8}^C)$. And from

\begin{eqnarray*}
  B_6(\alpha\phi_1\alpha^{-1}, \tau{}^t(\alpha\phi_2\alpha^{-1})\tau) \!\!\!&=&\!\!\! B_6(\alpha\phi_1\alpha^{-1}, \alpha\tau{}^t\phi_2\tau\alpha^{-1}) = B_6(\phi_1, \tau{}^t\phi_2\tau), 
\vspace{1mm}\\
   \langle \alpha\X, {}^t\alpha^{-1}\Y \rangle \!\!\!&=&\!\!\! \langle \X, \Y \rangle,
\end{eqnarray*}
we see that $\varphi_2(\alpha) \in E_8$. Evidently, $w\varphi_2(\alpha) = \varphi_2(\alpha)w$, hence, $\varphi_2(\alpha) \in (E_8)^w$. Now, we define a mapping $\varphi : SU(3) \times E_6 \to (E_8)^w$ by
$$
   \varphi(A, \alpha) = \varphi_1(A)\varphi_2(\alpha). $$
Since $\varphi_1(A)$ and $\varphi_2(\alpha)$ commute, $\varphi$ is a homomorphism. It is not difficult to show that $\Ker\varphi = \{(E, 1), (\omega E, \omega^21), (\omega^2E, \omega1) \} = \Z_3$. Finally, since $(E_8)^w$ is connected as the fixed points subgroup by automorphims $w$ of the simply connected group $E_8$ and 
$$
\begin{array}{l}
  (\gge_8)^w = \{ R \in {\gge_8}^C \, | \, wR = R, \tau\wti{\lambda}R = R \} 
\vspace{1mm}\\
\qquad \;\;
   = \{(D, \phi, 0, 0) \in {\gge_8}^C \, | \, D \in \su(3), \phi \in \gge_6 \}    \cong \su(3) \oplus \gge_6, 
\end{array}$$
$\varphi$ is onto. Thus we have the isomorphism $(SU(3) \times E_6)/\Z_3 \cong (E_8)^w$.
\vspace{4mm}

{\bf 5.11. Automorphism $w_3$ of order 3 and subgroup $SU(9)/\Z_3$ of $E_8$}
\vspace{3mm}

In order to construct another $C$-Lie algebra of type $E_8$, we investigate the properties of the exterior $C$-vector space ${\mit\Lambda}^k(C^n)$. Let $\e_1, \cdots, \e_n$ be the canonical $C$-basis of $n$-dimensional $C$-vector space $C^n$ and $(\x, \y)$ the inner product in $C^n$ satisfying $(\e_i, \e_j) = \delta_{ij}$. In ${\mit\Lambda}^k(C^n)$, we define an inner product by
$$
\begin{array}{c}
   (\x_1 \land \cdots \land \x_k, \y_1 \land \cdots \land \y_k) = \det\Big((\x_i, \y_j)\Big), \quad k \geq 1, 
\vspace{1mm}\\
   (a, b) = ab, \quad a, b \in {\mit\Lambda}^0(C^n) = C.
\end{array}$$
Then, $\e_{i_1} \land \cdots \land \e_{i_k}, i_1 < \cdots < i_k$ forms an orthonormal $C$-basis of ${\mit\Lambda}^k(C^n)$. For $\u \in {\mit\Lambda}^k(C^n)$, we define an element $*\u \in {\mit\Lambda}^{n-k}(C^n)$ by
$$
   (*\u, \v) = (\u \land \v, \e_1 \land \cdots \land \e_n), \quad \v \in {\mit\Lambda}^{n-k}(C^n). $$
Then, $*$ induces a $C$-linear isomorphism
$$
       * : {\mit\Lambda}^k(C^n) \to {\mit\Lambda}^{n-k}(C^n) $$
and satisfies the following identity:
$$
    *^2\u = (-1)^{k(n-k)}\u, \quad \u \in {\mit\Lambda}^k(C^n). $$
The group $SL(n, C)$ naturally acts on ${\mit\Lambda}^k(C^n)$ as
$$
    A(\x_1 \land \cdots \land \x_k) = A\x_1 \land \cdots \land A\x_k, \quad A1 = 1. $$
Hence the Lie algebra $\sl(n, C)$ acts on ${\mit\Lambda}^k(C^n)$ as
$$
  D(\x_1 \land \cdots \land \x_k) = \dsum_{j=1}^k\x_1 \land \cdots \land D\x_j \land \cdots \land \x_k, \quad D1 = 0. $$   

{\bf Lemma 5.11.1.} {\it For $A \in SL(n, C), D \in \sl(n, C)$ and $\u, \v \in {\mit\Lambda}^k(C^n)$, we have}
\vspace{1mm}

(1) $\;\; (A\u, {}^t\!A^{-1}\v) = (\u, \v), \quad (D\u, \v) + (\u, -{}^tD\v) = 0$.
\vspace{1mm}

(2) $\;\; *(A\u) = {}^t\!A^{-1}(*\u), \quad *(D\u) = -{}^t\!D^{-1}(*\u)$.
\vspace{2mm}

For $\u, \v \in {\mit\Lambda}^k(C^n) \; (1 \le k \le n)$, we define a $C$-linear mapping $\u \times \v$ of $C^n$ by
$$
  (\u \times \v)\x = *(\v \land *(\u \land \x)) + (-1)^{n-k}\dfrac{n - k}{n}(\u, \v)\x, \quad \x \in C^n. $$
Since $\tr(\u \times \v) = 0$, $\u \times \v$ can be regarded as element of $\sl(n, C)$ with respect to the canonical basis of $C^n$. 
\vspace{3mm}

{\bf Lemma 5.11.2.} {\it For $A \in SL(n, C), D \in \sl(n, C)$ and $\u, \v \in {\mit\Lambda}^k(C^n)$, we have}
\vspace{1mm}

(1) $\; A(\u \times \v)A^{-1} = A\u \times {}^tA^{-1}\v, \;\; [D, \u \times \v] = D\u \times \v + \u \times (-{}^t\!D\v)$.
\vspace{1mm}

(2) $\; {}^t(\u \times \v) = \v \times \u, \quad \tau(\u \times \v) = \tau(\u) \times \tau(\v)$.
\vspace{1mm}

(3) $\;\; \tr(D(\u \times \v)) = (-1)^{n-k}(D\u, \v)$.
\vspace{2mm}

Now, we construct another $C$-Lie algebra ${\gge_8}^C$ of type $E_8$. 
\vspace{3mm}

{\bf Theorem 5.11.3.} {\it In an $80 + 84  + 84 = 248$ dimensional $C$-vector space}
$$
   {\gge_8}^C = \sl(9, C) \oplus {\mit\Lambda}^3(C^9) \oplus {\mit\Lambda}^3(C^9), $$
{\it we define a Lie bracket $[R_1, R_2]$ by}
$$
   [(D_1, \u_1, \v_1), (D_2, \u_2, \v_2)] = (D, \u, \v), $$
{\it where}
$$
\left\{
\begin{array}{l}
   D = [D_1, D_2] + \u_1 \times \v_2 - \u_2 \times \v_1
\vspace{1mm}\\
   \u = D_1\u_2 - D_2\u_1 + *(\v_1 \land \v_2)
\vspace{1mm}\\
   \v = - {}^tD_1\v_2 + {}^tD_2\v_1 - *(\u_1 \land \u_2),
\end{array} \right. $$
{\it then ${\gge_8}^{C}$ becomes a $C$-Lie algebra of type $E_8$.}
\vspace{2mm}

{\bf Proof.} In order to prove the Jacobi identity, we need the following Lemma.\vspace{3mm}

{\bf Lemma 5.11.4.} {\it For $\u, \v, \w \in {\mit\Lambda}^3(C^9)$, we have}
\vspace{1mm}

(1) $\;\; \u \times *(\v \land \w) + \v \times *(\w \land \u) + \w \times *(\u \land \v) = 0,$
\vspace{1mm}

(2) $\;\; (\u \times \w)\v - (\v \times \w)\u + *(*(\u \times \v) \land \w) = 0$.
\vspace{2mm}

{\bf Proof.} Let $\u = \u_1 \land \u_2 \land \u_3, \v = \u_4 \land \u_5 \land \u_6$ and $\w = \u_7 \land \u_8 \land \u_9$. For $\x, \y \in C^9$, we have
\begin{eqnarray*}
  ((\u \times \v)\x, \y) \!\!\! &=& \!\!\! (*(\v \land *(\u \land \x)), \y) + \dfrac{2}{3}(\u, \v)(\x, \y)
\vspace{1mm}\\
   \!\!\! &=& \!\!\! - (\x \land \u, \y \land \v) + \dfrac{2}{3}(\u, \v)(\x, \y)
\vspace{1mm}\\
   \!\!\! &=& \!\!\! (\x \land \u_2 \land \u_3, \v)(\u_1, \y) - (\x \land \u_1 \land \u_3, \v)(\u_2, \y) 
\vspace{1mm}\\
   \!\!\! && \quad + (\x \land \u_1 \land \u_2, \v)(\u_3, \y) - \dfrac{1}{3}(\u, \v)(\x, \y).
\end{eqnarray*}   
Hence 
$$
\displaylines {\hfill
\begin{array}{l}
    (\u \times \v)\x = (\x \land \u_2 \land \u_3, \v)\u_1 + (\u_1 \land \x \land \u_3, \v)\u_2
\vspace{1mm}\\
\qquad \quad \qquad + (\u_1 \land \u_2 \land \x, \v)\u_3 - \dfrac{1}{3}(\u, \v)\x.  
\end{array}
\hfill {}
\vspace{1mm}\\
\mbox{(i)}}$$
Using this identity, 
$$
\begin{array}{l}
  (\u \times *(\v \land \w) + \v \times *(\w \land \u) + \w \times *(\u \land \v))\x
\vspace{1mm}\\
\qquad 
   = \dsum_{j=1}^9(\u_1 \land \cdots \land \u_{j-1} \land \x \land \u_{j+1} \land \cdots \land \u_9, \e_1 \land \cdots \land \e_9)\u_j
\vspace{1mm}\\
\qquad \quad 
   - (\u_1 \land \cdots \land \u_9, \e_1 \land \cdots \land \e_9)\x = \mbox{(ii)}.
\end{array}$$
Denote $\x = \sum_{i=1}^9x_i\e_i, \u_j = \sum_{k=1}^9u_{jk}\e_k$ and $U = \big(u_{jk}\big) \in M(9, C)$. Hence we 
\vspace{-2mm}have

$$
\begin{array}{l}
   (\u_1 \land \cdots \land \u_{j-1} \land \x \land \u_{j+1} \land \cdots \land \u_9, e_1 \land \cdots \land \e_9) = \dsum_{k=1}^9\wti{u}_{jk}x_k,
\vspace{1mm}\\
   (\u_1 \land \cdots \land \u_9, \e_1 \land \cdots \land \e_9) = \det U,
\end{array} $$
where $\wti{u}_{jk}$ is the cofactor of $u_{jk}$ of the matrix $U$. Therefore 
\begin{eqnarray*}
   \mbox{(ii)} \!\!\! &=& \!\!\! \dsum_{j,k}x_k\wti{u}_{jk}\u_j - (\det U)\x = \dsum_{i,j,k}x_k\wti{u}_{jk}u_{ji}\e_i - (\det U)\x
\vspace{1mm}\\
   \!\!\! &=& \!\!\! \dsum_{j,k}x_k(\det U)\delta_{ki}\e_i - (\det U)\x = 0.
\end{eqnarray*}
Thus (1) is proved. Next, let $\u = \u_1 \land \u_2 \land \u_3$ and $\v = \v_1 \land \v_2 \land \v_3$. Using (i), for any $\a \in {\mit\Lambda}^3(C^9)$, we have
\vspace{2mm}

$\qquad
   ((\u \times \w)\v - (\v \times \w)\u, \a)$
\vspace{1mm}

$\qquad  \quad
   = (((\u \times \w)\v_1) \land \v_2 \land \v_3, \a) - (((\u \times \w)\v_2) \land \v_1 \land \v_3, \a)$
\vspace{1mm}

$\qquad \qquad
   + (((\u \times \w)\v_3) \land \v_1 \land \v_2, \a) - (((\v \times \w)\u_1) \land \u_2 \land \u_3, \a)$
\vspace{1mm}

$\qquad \qquad
   + (((\v \times \w)\u_2) \land \u_1 \land \u_3, \a) - (((\v \times \w)\u_3) \land \u_1 \land \u_2, \a)$
\vspace{1mm}

$\qquad \quad
   = - (\u, \w)(\v, \a) + \dsum_{i=1}^3\dsum_{j=1}^3(\u_i \land \u_{i+1} \land \v_j, \w)(\u_{i+2} \land \v_{i+1} \land \v_{j+2}, \a)$
\vspace{1mm}

$\qquad \qquad
   +  (\v, \w)(\u, \a) - \dsum_{i=1}^3\dsum_{j=1}^3(\u_i \land \v_j \land \v_{j+1}, \w)(\u_{i+1} \land \u_{i+2} \land \v_{j+2}, \a)$
\vspace{1mm}

$\qquad \quad
   = - (\u \land \v, \w \land a) = - (*(*(\u \land \v) \land \w), \a)$.
\vspace{1mm}

\noindent Thus (2) is proved. 
\vspace{1mm}

From Lemmas 5.11.1, 5.11.2 and 5.11.4, we can prove that ${\gge_8}^C$ is a $C$-Lie algebra. Therefore Theorem 5.11.3 is proved.
\vspace{2mm}

We will show that the Lie algebra ${\gge_8}^C$ given in Theorem 5.11.3 is also the Lie algebra of type $E_8$. Since we can not give explicit isomorphism between this ${\gge_8}^C$ and ${\gge_8}^C$ of Theorem 5.1.1, we shall show that this ${\gge_8}^C$ is simple.
\vspace{3mm}

{\bf Theorem 5.11.5.} {\it ${\gge_8}^C = \sl(9, C) \oplus {\mit\Lambda}^3(C^9) \oplus {\mit\Lambda}^3(C^9)$ is a simple $C$-Lie algebra of type $E_8$.}
\vspace{2mm}

{\bf Proof.} We use the decomposition
$$
     {\gge_8}^C = \sl(9, C) \oplus \gq, \quad \gq = {\mit\Lambda}^3(C^9) \oplus {\mit\Lambda}^3(C^9). $$
For a subset $I = \{i, j, k \}$ $(i < j < k)$ of $\{1, 2, \cdots, 9 \}$, we put
$$
   \e_I = \e_i \land \e_j \land \e_k \in {\mit\Lambda}^3(C^9). $$
Now, let $\ga$ be a non-zero ideal of $\gg = {\gge_8}^C$. 
\vspace{1mm}

(1) Case $\sl(9, C) \cap \ga = \{0\}$ and $\gq \cap \ga = \{0\}$. Let $p : \gg 
\to \sl(9, C)$ be the projection. If $p(\ga) = 0$, then $\ga$ is contained in 
$\gq$, which contradicts $\gq \cap \ga = \{0\}$. Hence, $p(\ga)$ is a non-zero 
ideal of $\sl(9, C)$, so we have $p(\ga) = \sl(9, C)$. For an element $D = 
\sum_{i=1}^8H_i \in \sl(9, C), H_i = E_{ii} - E_{99}$, there exists an element 
$(\u, \v) = \big(\sum_Iu_I\e_I, \sum_Jv_J\e_J\big) \in \gq$ such that $(D, \u, 
\v) \in \ga$. Since $[(D, 0, 0), (X, \u, \v)] = (0, D\u, -{}^tD\v) \in \gq \cap \ga = \{0\}$, we have
$$
\begin{array}{l}
    0 = D\u = \dsum_Iu_ID\e_I = 3\dsum_{I \not\ni 9}u_I\e_I - 6\dsum_{I \ni 9}u_I\e_I, 
\vspace{1mm}\\
    0 = - {}^tD\v = -3\dsum_{J \not\ni 9}v_J\e_J + 6\dsum_{J \ni 9}v_J\e_J, 
\end{array}$$
i.e., $u_I = 0$ and $v_J = 0$. Then, $0 \not= (D, \u, \v) = (D, 0, 0) \in \sl(9, C) \cap \ga = \{0\}$. This is a contradiction.
\vspace{1mm}

(2) Case $\sl(9, C) \cap \ga \not= \{0\}$. Since $\sl(9, C) \cap \ga$ is a non-zero ideal of $\sl(9, C)$, we have $\sl(9, C) \subset \ga$. For any $\e_i \land \e_j \land \e_k \in {\mit\Lambda}^3(C^9)$, put
$$
   D = \dfrac{1}{3}(E_{ii} + E_{jj} + E_{kk}) - E_{ll}, \quad \mbox{$l \not=i, j, k$}. $$
Since $(D, 0, 0) \in \sl(9, C) \subset \ga$, we see that
$$
\begin{array}{l}
   (0, \e_i \land \e_j \land \e_k, 0) = [(D, 0, 0), (0, \e_i \land \e_j \land \e_k, 0)] \in \ga, 
\vspace{1mm}\\
   (0, 0, \e_i \land \e_j \land \e_k) = [(D, 0, 0), (0, 0, - \e_i \land \e_j \land \e_k)] \in \ga. 
\end{array}$$
It follows that $\gq \subset \ga$. Hence we have $\ga = \gg$. 
\vspace{1mm}

(3) Case $\gq \cap \ga \not= \{0\}$. Let $R = (0, \u, \v)$ be a non-zero element of $\gq \cap \ga$. In the case $\u \not= 0$, we put $\u = \sum_Iu_I\e_I$. Without loss of generality, we may assume that $u_{\{123\}} = 1$. Putting $S_{ij} = (E_{ii} - E_{jj}, 0, 0) \in \gg$ and $T = (0, 0, \e_1 \land \e_2 \land \e_4) \in \gg$, we have
$$
\begin{array}{l}
   0 \not= \ad(T)\ad(S_{37})\ad(S_{27})\ad(S_{17})\ad(S_{36})\ad(S_{25})\ad(S_{14})R
\vspace{1mm}\\
  \;\;\, = (-E_{34}, 0, 0) \in \sl(9, C) \cap \ga. 
\end{array}$$
Then we can reduce this case to the case (2). In case $\v \not= 0$, we ca similarly reduce to the case (2).
\vspace{1mm}

\noindent Thus the simplicity of $\gg$ has been proved. Since the dimension of $\gg$ is 248, we see that $\gg$ is a Lie algebra of type $E_8$.
\vspace{3mm}

{\bf Proposition 5.11.6.} {\it The Killing form $B_8$ of the Lie algebra ${\gge_8}^C = \sl(9, C) \oplus {\mit\Lambda}^3(C^9) \oplus {\mit\Lambda}^3(C^9)$ is given by}
$$
   B_8((D_1, \u_1, \v_1), (D_2, \u_2, \v_2)) = 60(\tr(D_1D_2) + (\u_1, \v_2) + (\u_2, \v_1)). $$

{\bf Proof.} We consider a symmetric bilinear form $B$ of ${\gge_8}^C$:
$$
    B((D_1, \u_1, \v_1), (D_2, \u_2, \v_2)) = \tr(D_1D_2) + (\u_1, v_2) + (\u_2, \v_1). $$
Using Lemmas 5.11.2, 5.11.4, we see that $B$ is ${\gge_8}^C$-adjoint invariant. Since ${\gge_8}^C$ is simple, there exists $k \in C$ such that $B_8(R_1, R_2) = kB(R_1, R_2)$ for all $R_i \in {\gge_8}^C$. To determined $k$, let $R = R_1 = R_2 = (E_{11} - E_{22}, 0, 0) \in {\gge_8}^C$. Then we have
$$
     B_8(R, R) = 120, \quad B(R, R) = 2. $$
Therefore $k = 60$.
\vspace{2mm}

We define a complex-conjugate linear transformation $\tau\wti{\lambda}$ of ${\gge_8}^C$ by
$$
    \tau\wti{\lambda}(D, \u, \v) = (-\tau{}^tD, -\tau\v, -\tau\u). $$
and we define a Hermitian inner product $\langle R_1, R_2 \rangle$ in ${\gge_8}^C$ by
$$
    \langle R_1, R_2 \rangle = -B_8(R_1, \tau\wti{\lambda}R_2). $$
Then we have
$$
    \langle R_1, R_2 \rangle = 60(\tr(D_1\tau{}^tD_2) + (\u_1, \tau\u_2) + (\v_2, \tau\v_1)). $$
As in Theorem 5.5.3, 
$$
     E_8 = \{ \alpha \in \mbox{Aut}({\gge_8}^C) \, | \, \langle \alpha R_1, \alpha R_2 \rangle = \langle R_1, R_2 \rangle \} $$
is a simly connected compact simple Lie group of type $E_8$. 
\vspace{1mm}

We define a $C$-linear transformation $w_3$ of ${\gge_8}^C$ by
$$
   w_3(D, \u, \v) = (D, \omega\u, \omega^2\v), $$
where $\omega = - \dfrac{1}{2} + \dfrac{\sqrt{3}}{2}i \in C$. Then, $w_3 \in E_8$ and ${w_3}^3 = 1$.
\vspace{2mm}

Now, we study the following subgroup $(E_8)^{w_3}$ of $E_8$:
$$
     (E_8)^{w_3} = \{\alpha \in E_8 \, | \, w_3\alpha = \alpha w_3 \}. $$

{\bf Theorem 5.11.7.} \qquad $(E_8)^{w_3} \cong SU(9)/\Z_3, \;\; \Z_3 = \{E, \omega E, \omega^2E \}. $
\vspace{2mm}

{\bf Proof.} We define a mapping $\varphi : SU(9) \to (E_8)^{w_3}$ by
$$
   \varphi(A)(D, \u, \v) = (ADA^{-1}, A\u, {}^tA^{-1}\v). $$
$\varphi$ is well-defined\,: $\varphi(A) \in (E_8)^{w_3}$. Indeed, for $A = \exp X, X \in \su(9)$, we have
$$
\begin{array}{l}
   \exp(\ad(X, 0, 0))(D, \u, \v) = (\exp(\ad(X)D, (\exp X)\u, (\exp(-{}^tX))\v)
\vspace{1mm}\\
\qquad \qquad 
   = (\Ad(\exp X)D, (\exp X)\u, {}^t(\exp X)^{-1}\v)
\vspace{1mm}\\
\qquad \qquad 
   = \varphi(\exp X)(D, \u, \v).
\end{array}$$
Hence $\varphi(A) \in E_8$. Clearly $w_3\varphi(A) = \varphi(A)w_3$. Therefore $\varphi(A) \in (E_8)^{w_3}$. Obviously $\varphi$ is a homomorphism. We shall show that $\varphi$ is onto. Since the Lie algebra $(\gge_8)^{w_3}$ of the group $(E_8)^{w_3}$ is 
$$
  (\gge_8)^{w_3} = \{R \in {\gge_8}^C \, | \, \tau\wti{\lambda}R = R, w_3R = R \} = \{(D, 0, 0) \in {\gge_8}^C \, | \, D \in \su(9) \} \cong \su(9), $$
the differential $\varphi_*$ of $\varphi$ is onto. Since $(E_8)^{w_3}$ is connected, $\varphi$ is also onto. It is not difficult to see that $\ker\,\varphi = \{E, \omega E, \omega^2E \} = \Z_3$. Thus we have the isomorphism $SU(9)/\Z_3 \cong (E_8)^{w_3}$.
\vspace{4mm}

{\bf 5.12. Automorphism $z_5$ of order 5 and subgroup $(SU(5) \times SU(5))/\Z_5$ of $E_8$}
\vspace{3mm}

We shall construct one more $C$-Lie algebra of type $E_8$.
\vspace{3mm}

{\bf Theorem 5.12.1.} {\it In a $48 + 50 \times 4 = 248$ dimensional $C$-vector space}
$$
   {\gge_8}^C = \gg_0 \oplus \gg_1 \oplus \gg_2 \oplus \gg_{-2} \oplus \gg_{-1}, $$
(suffices are considered mod 5) {\it where}
$$
\begin{array}{c}
    \gg_0 = \sl(5, C) \oplus \sl(5, C),
\vspace{1mm}\\
    \gg_1 = C^5 \otimes {\mit\Lambda}^2(C^5) = \gg_{-1}, \quad \gg_2 = {\mit\Lambda}^2(C^5) \otimes C^5 = \gg_{-2}, 
\end{array}$$
{\it we define a Lie bracket $[R_1, R_2]$ as follows.}

$$
\begin{array}{ll}
   [\gg_0, \gg_0] \subset \gg_0 & [(C_1, D_1), (C_2, D_2)] = ([C_1, C_2], [D_1, D_2]),
\vspace{1mm}\\
   {[}\gg_0, \gg_1] \subset \gg_1 & [(C, D), \x \otimes \a] = (C\x) \otimes \a + \x \otimes (D\a),
\vspace{1mm}\\
   {[}\gg_0, \gg_2] \subset \gg_2 & [(C, D), \b \otimes \y] = (C\b) \otimes \y + \b \otimes (- {}\,^tD\a),
\vspace{1mm}\\
   {[}\gg_0, \gg_{-2}] \subset \gg_{-2} & [(C, D), \c \otimes \z] = (-{}^tC\c) \otimes \z + \c \otimes (D\z),
\vspace{1mm}\\
   {[}\gg_0, \gg_{-1}] \subset \gg_{-1} & [(C, D), \w \otimes \dd] = (-{}^tC\w)  \otimes \z + \c \otimes (- {}\,^tD\z),
\vspace{1mm}\\
   {[}\gg_1, \gg_{-1}] \subset \gg_0 & [\x \otimes \a, \w \otimes \dd] = (- (\a, \dd)\x \times \w, (\x, \w)\a \times \dd),
\vspace{1mm}\\
   {[}\gg_2, \gg_{-2}] \subset \gg_0 & [\b \otimes \y, \c \otimes \z] = ((\y, \z)\b \times \c, (\b, \c)\z \times \y),
\vspace{1mm}\\
\end{array}$$
\vspace{-2mm}
$$
\begin{array}{ll}
\begin{array}{l}
   {[}\gg_1, \gg_1] \subset \gg_2
\vspace{1mm}\\
   {[}\gg_{-1}, \gg_{-1}] \subset \gg_{-2} 
\end{array}
   & [\x_1 \otimes \a_1, \x_2 \otimes \a_2] = (\x_1 \land \x_2) \otimes *(\a_1 \land \a_2),
\vspace{1mm}\\
\begin{array}{l}
   {[}\gg_2, \gg_2] \subset \gg_{-1} 
\vspace{1mm}\\
   {[}\gg_{-2}, \gg_{-2}] \subset \gg_1 
\end{array}
   & [\b_1 \otimes \y_1, \b_2 \otimes \y_2] = *(\b_1 \land \b_2) \otimes (\y_1 \land \y_2),
\vspace{1mm}\\
\begin{array}{l}
   {[}\gg_1, \gg_2] \subset \gg_{-2} 
\vspace{1mm}\\
   {[}\gg_{-1}, \gg_{-2}] \subset \gg_2 
\end{array}
   & [\x \otimes \a, \b \otimes \y] = *(\b \land \x) \otimes *(*\a \land \y),
\vspace{1mm}\\
\begin{array}{l}
   {[}\gg_2, \gg_{-1}] \subset \gg_1 
\vspace{1mm}\\
   {[}\gg_{-2}, \gg_1] \subset \gg_{-1} 
\end{array}
   & [\b \otimes \y, \w \otimes \dd] = *(*\b \land \w) \otimes *(\dd \land \y).
\end{array} $$
{\it Then ${\gge_8}^{C}$ becomes a $C$-Lie algebra.}
\vspace{2mm}

{\bf Proof.} In order to prove the Jacobi identity, we need the following Lemma.\vspace{3mm}

{\bf Lemma 5.12.2.} {\it For $\x, \y, \z \in {\mit\Lambda}^1(C^5) = C^5$ and $\a, \b, \c \in {\mit\Lambda}^2(C^5)$, we have}
\vspace{1mm}

(1) $\;\; *\a \land *(\b \land \c) + *\b \land *(\c \land \a) + *\c \land *(\a \land \b) = 0$,
\vspace{1mm}

(2) $\;\; *(\a \land *(*\b \land \x)) + *(\b \land *(*\a \land \x)) + \x \land *(\a \land \b) = 0,$
\vspace{1mm}

(3) $\;\; *(*(\x \land \y) \land \z) = (\x, \z)\y - (\y, \z)\x,$
\vspace{1mm}

(4) $\;\; \x \land *(*\a \land \y) + *(\y \land +(\a \land \x)) - (\x, \y)\a = 0,$
\vspace{1mm}

(5) $\;\; *(\a \land *(\b \land \x)) - *(*\b \land *(*\a \land \x)) - (\a, \b)\x = 0,$
\vspace{1mm}

(6) $\;\; \a \times *(\b \land \x) + \b \times *(\a \land \x) - \x \times *(\a \land \b) = 0,$
\vspace{1mm}

(7) $\;\; *(*\a \land \x) \times \y - *(*\a \land \y) \times \x + \a \times (\x \land \y) = 0,$
\vspace{1mm}

(8) $\;\; (\a \land \b)\c = *(*(\a \land \c) \land \b) - \dfrac{1}{5}(\a, \b)\c - (\b, \c)\a,$
\vspace{1mm}

(9) $\;\; (\x \times \y)\a = - *(\y \land *(\x \land \a)) + \dfrac{3}{5}(\x, \y)\a. $
\vspace{2mm}

{\bf Proof.} (1) Let $\a = \a_1 \land \a_2, \b = \a_3 \land \a_4, \c =  \a_5 \land \a_6$ and $\a_i = \sum_{j=1}^5a_{ij}\e_j$. Since
\vspace{2mm}

$\qquad 
  (*(*\a \land *(\b \land \c)), \x) = (\a, *(\b \land \c) \land \x)$
\vspace{1mm}

$\qquad \qquad \quad
    = (\a_1, *(\b \land \c))(\a_2, \x) - (\a_2, *(\b \land \c))(\a_1, \x) $ 
\vspace{1mm}

$\qquad \qquad \quad
    = (\a_1 \land \b \land \c, \e_1 \land \cdots \land \e_5)(\a_2, \x) - (\a_2 \land \b \land \c, \e_1 \land \cdots \land \e_5)(\a_1, \x)$,
\vspace{2mm}    

\noindent we have

\qquad 
$*(*\a \land *(\b \land \c) + *\b \land *(\c \land \a) + *\c \land *(\a \land \b))$
\vspace{-2mm}
$$
\begin{array}{l}
   = \dsum_{j=1}^5\dsum_{i=1}^6(-1)^i(\a_1 \land \cdots \a_{i-1} \land \a_{i+1} \land \cdots \land \a_6, \e_1 \land \cdots \land \e_5)a_{ij}\e_j
\vspace{1mm}\\
\quad
   = - \dsum_{j=1}^5\det\pmatrix{a_{1j} & a_{11} & \cdots & a_{15} \cr 
                                 a_{2j} & a_{21} & \cdots & a_{25} \vspace{1mm}\cr 
                                        & & \cdots &  \vspace{1mm} \cr
                                 a_{6j} & a_{61} & \cdots & a_{65}}\e_j = 0.
\end{array}$$

(2) Using (1), we have
$$
\begin{array}{l}
   (*(\a \land *(*\b \land \x)) + *(\b \land *(*\a \land \x)) + \x \land *(\a \land \b), \c)
\vspace{1mm}\\ 
\qquad
   = (*\b \land \x, \c \land \a) + (*\a \land \x, \b \land \c) + (*\c \land \x, \a \land \b)
\vspace{1mm}\\
\qquad
   = (*\x, *\b \land *(\c \land \a) + *\a \land *(\b \land \c) + *\c \land *(\a \land \b)) = 0.
\end{array}$$

(3) For any $\v \in {\mit\Lambda}^1(C^5) = C^5$, we have
$$
   (*(*(\x \land \y) \land \z), \v) = (\x \land \y, \z \land \v) = (\x, \z)(\y, \v) - (\y, \z)(\x, \v). $$
This shows (3).
\vspace{1mm}

(4),(5) Let $\a = \a_1 \land \a_2$. Since
$$
\begin{array}{l}
   (\x \land \a, \y \land \b) = (\x, \y)(\a, \b) - (\a_1, \y)(\x \land \a_2, \b) + (\a_2, \y)(\x \land \a_1, \b),
\vspace{1mm}\\
   (*\b \land \x, *\a \land \y) = (\a, \y \land *(*\b \land \x)) = (\a_1, \y)(\x \land \a_2, \b) - (\a_2, \y)(\x \land \a_1, \b), 
\end{array}$$
we have
$$
   (\x \land \a, \y \land \b) + (*\b \land \x, *\a \land \y) = (\x, \y)(\a, \b). $$
Using this identity, we have
$$
\begin{array}{l}
   (\x \land *(*\a \land \y) + *(\y \land *(\a \land \x)) - (\x, \y)\a, \b)
\vspace{1mm}\\
\qquad
   = (*\b \land \x, *\a \land \y) + (\x \land \a, \y \land \b) - (\x, \y)(\a, \b) = 0,
\vspace{1mm}\\
   (*(\a \land *(\b \land \x)) - *(*\b \land *(*\a \land \x)) - (\a, \b)\x, \y)
\vspace{1mm}\\
\qquad
   = (\x \land \b, \y \land \a) + (*\a \land \x, *\b \land \y) - (\x, \y)(\a, \b) = 0.
\end{array}$$

(6) Since 
$$
   ((\x \times \y)\z, \v) = -(\x \land \z, \y \land \v) + \dfrac{4}{5}(\x, \y)(\z, \v) = (\y, \z)(\x, \v) - \dfrac{1}{5}(\x, \y)(\z, \v), $$
we have
$$
\displaylines {\hfill
   (\x \times \y)\z = (\y, \z)\x - \dfrac{1}{5}(\x, \y)\z.
\hfill \mbox{(i)}}$$
For $\v, \w \in {\mit\Lambda}^1(C^5) = C^5$, we have
\vspace{2mm}

\qquad \quad
   $((\a \times *(\b \land \x))\v, \w)$
\vspace{1mm}

\qquad \qquad \quad
   $= (\v \land \a, \w \land *(\b \land \x)) - \dfrac{3}{5}(*(\a \land \b), \x)(\v, \w)$
\vspace{1mm}

\qquad \qquad \quad
   $= (\v, \w)(\a, *(\b \land \x)) - (\a_1, \w)(\w \land \a_2, *(\b \land \x))$
\vspace{1mm}

\qquad \quad \qquad
   $+ (\a_2, \w)(\w \land \a_1, *(\b \land \x)) - \dfrac{3}{5}(*(\a \land \b), \x)(\v, \w)$
\vspace{1mm}

\qquad \qquad \quad
   $=\dfrac{2}{5}(*(\a \land \b), \x)(\v, \w) - (\a_1, \w)(\b \land \x \land \w \land \a_2, \e_1 \land \cdots \land \e_5)$
\vspace{1mm}

\qquad \quad \qquad $\;\;$
   $+(\a_2, \w)(\b \land \x \land \w \land \a_1, \e_1 \land \cdots \land \e_5),$\vspace{2mm}

\qquad \qquad 
   $((\b \times *(\a \land \x))\v, \w)$
\vspace{1mm}

\qquad \qquad \quad
   $= (\x \land \a, \w \land *(\b \land \v)) - \dfrac{3}{5}(*(\a \land \b), \x)(\v, \w)$
\vspace{1mm}

\qquad \qquad \quad 
   $= (\v, \w)(*(\a \land \b), \v) + (\a_1, \w)(\b \land \x \land \w \land \a_2, \e_1 \land \cdots \land \e_5)$
\vspace{1mm}

\qquad \qquad \quad
   $- (\a_2, \w)(\b \land \x \land \w \land \a_1, \e_1 \land \cdots \land \e_5) - \dfrac{3}{5}(*(\a \land \b), \x)(\v, \w).$
\vspace{2mm}

\noindent Using (i), we have
$$
\begin{array}{l}
   (\a \times *(\b \land \x))\v + (\b \times *(\a \land \x))\v = (\x, \w)*(\a \land \b) - \dfrac{1}{5}(*(\a \land \b), \x)\v
\vspace{1mm}\\ 
\qquad \qquad \qquad
   = (\x \times *(\a \land \b))\v.
\end{array}$$

(7) We have 
$$
\begin{array}{l}
   ((\a \times (\x \land \y))\v, \w) = (\x \land \y \land \w, \v \land \a) - \dfrac{3}{5}(\a, \x \land \y)(\v, \w)
\vspace{1mm}\\
\qquad
   = (\x, \v)(\y \land \w, \a) - (\y, \v)(\x \land \w, \a) + \dfrac{2}{5}(\a, \x \land \y)(\v, \w)
\vspace{1mm}\\
\qquad
   = (\x, \v)(*(*\a \land \y), \w) - (\y, \v)(*(*\a \land \x), \w) + \dfrac{2}{5}(\a, \x \land \y)(\v, \w).
\end{array}$$
On the other hand, using (i), we have
$$
\begin{array}{l} 
   (*(*\a \land \x) \times \y)\v = (\y, \v)*(*\a \land \x) - \dfrac{1}{5}(*(*\a \land \x), \v)\y,
\vspace{1mm}\\
\qquad \qquad 
   = (\y, \v)*(*\a \land \x) - \dfrac{1}{5}(\a, \x \land \v)\y,
\vspace{1mm}\\
   -(*(*\a \land \x) \times \y)\v = - (\x, \v)*(*\a \land \y) - \dfrac{1}{5}(\a, \x \land \v)\y.
\end{array}$$
Hence (7) is proved.
\vspace{1mm}

(8) Let $\a = \a_1 \land \a_2$ and $\c = \c_1 \land \c_2$. Since
$$
\begin{array}{l}
   ((\a \times \b)\v, \w) = (\a \land \v, \b \land \w) - \dfrac{3}{5}(\a, \b)(\v, \w)
\vspace{1mm}\\
\qquad
   = - (\a_1, \w)(\x \land \a_2, \b) + (\a_2, \w)(\x \land \a_1, \b) + \dfrac{2}{5}(\a, \b)(\v, \w),
\end{array}$$
we have
$$
\begin{array}{l}
   (\a \times \b)\c = - (\a_1 \land \c_1, \b)\a_1 \land \c_2 + (\a_1 \land \c_2, \b)\a_2 \land \c_1
\vspace{1mm}\\
\qquad
   +(\a_2 \land \c_1, \b)\a_1 \land \c_2 - (\a_2 \land \c_2, \b)\a_1 \land \c_1 + \dfrac{4}{5}(\a, \b)\c.
\end{array}$$
On the other hand, for $\dd \in {\mit\Lambda}^2(C^5)$, we have
$$
\begin{array}{l}
   (*(*(\a \land \c) \land \b), \dd) = (\a \land \c, \b \land \dd)
\vspace{1mm}\\
\qquad
   = (\a, \b)(\c, \dd) - (\a_1 \land \c_1, \b)(\a_2 \land \c_2, \dd) + (\a_1 \land \c_2, \b)(\a_2 \land \c_1, \dd)
\vspace{1mm}\\
\qquad \quad
   + (\a_2 \land \c_1, \b)(\a_1 \land \c_2, \dd) - (\a_2 \land \c_2, \b)(\a_1 \land \c_1, \dd) + (\c, \b)(\a, \dd).
\end{array}$$
Hence (8) is proved.
\vspace{1mm}

(9) Using (i), we have
$$
   (\x \times \y)\a = (\y, \a_1)\x \land \a_2 - (\y, \a_2)\x \land \a_1 - \dfrac{2}{5}(\x, \y)\a. $$
On the other hand, we have
$$
\begin{array}{l}
   (-*(\y \land *(\x \land \a)), \b) = - (\x \land \a, \y \land \b)
\vspace{1mm}\\
\qquad
   =(\y, \a_1)(\x \land \a_2, \b) - (\y, \a_2)(\x \land \a_1, \b) - (\x, \y)(\a, \b).
\end{array}$$
Hence (9) is proved.
\vspace{1mm}

From Lemmas 5.11.1, 5.11.2 and 5.12.2, we can prove that ${\gge_8}^C$ be comes a $C$-Lie algebra. Furthermore we have the following theorem.
\vspace{3mm}

{\bf Theorem 5.12.3.} {\it The Lie algebra ${\gge_8}^C = \gg_0 \oplus \gg_1 \oplus \gg_2 \oplus \gg_{-2} \oplus \gg_{-1}$ is a $C$-simple Lie algebra of type $E_8$.}
\vspace{2mm}

{\bf Proof.} We shall show that ${\gge_8}^C$ is simple. For this end, we use the decomposition $\gg = {\gge_8}^C = \gg_{01} \oplus \gg_{12} \oplus \gq$, where
\begin{eqnarray*}
   \gg_{01} \!\!\! &=& \!\!\! \{(C, 0) \in \gg_0 \, | \, C \in \sl(5, C) \} \cong \sl(5, C), 
\vspace{1mm}\\
   \gg_{02} \!\!\! &=& \!\!\! \{(0, D) \in \gg_0 \, | \, D \in \sl(5, C) \} \cong \sl(5, C), 
\vspace{1mm}\\
   \gq \!\!\! &=& \!\!\! \gg_1 \oplus \gg_2 \oplus \gg_{-2} \oplus \gg_{-1}.
\end{eqnarray*}
Now, let $\ga$ be a non-zero ideal of $\gg$. There are three cases to be considered.
\vspace{1mm} 

(1) Case $\gg_{01} \cap \ga = \{0\}, \gg_{02} \cap \ga = \{0\}$ and $\gq \cap \ga = \{0\}$. Let $p_i : \gg \to \gg_{0i}$ ($i = 1, 2$) denote the projection. If $p_1(\ga) = \{0\}$ and $p_2(\ga) = \{0\}$, then $\ga$ is contained in $\gq$, which contradicts $\gq \cap \ga = \{0\}$. Hence, without loss of generality, we may assume that $p_1(\ga) = \gg_{01}$, because $\gg_{01}$ is a simple Lie algebra. For $C = \sum_{i=1}^4H_i \in \sl(3, C)$ where $H_i = E_{ii} - E_{55}$, there exists $(D, g_1, g_2, g_{-2}, g_{-1}) \in \gg_{01} \oplus \gq$ such that $(C, D, g_1, g_2, g_{-2}, g_{-1}) \in \ga$. Since
$$
\begin{array}{l}
   [(C, 0), (C, D, g_1, g_2, g_{-2}, g_{-1})]
\vspace{1mm}\\ 
\qquad
   = (0, 0, [C, g_1], [C, g_2], [C, g_{-2}], [C, g_{-1}]) \in \gq \cap \ga = \{0\},
\end{array}$$
we have $[C, g_i] = 0$ ($i = 1, 2, -2, -1$). Since any eigenvalue of ad$X$ is not 0, we have $g_i = 0$. Then we have $(C, D) \in \gg_0 \cap \ga$. Since
$$
   [(C, D), (E_{45}, 0)] = (5E_{45}, 0) \in \gg_{01} \cap \ga, $$
we have $\gg_{01} \cap \ga \not= \{0\}$. This is a contradiction.
\vspace{1mm}

(2) Case $\gg_{01} \cap \ga \not= \{0\}$ or $\gg_{02} \cap \ga \not= \{0\}$ . We may assume that $\gg_{01} \cap \ga \not= \{0\}$. Since $\gg_{01}$ is simple, we have $\gg_{01} \subset \ga$. Since $[\gg_{01}, \gg_i] = \gg_i$ ($i = 1, 2, -2, -1$), we have $\gq \subset \a$. Since
$$
   \ga \supset [\gg_1, \gg_{-1}] \ni [\e_1 \otimes (\e_1 \land \e_2), \e_1 \otimes (\e_1 \land \e_3)] = (0, -E_{23}), $$
we have $\gg_{02} \cap \ga \not= \{0\}$. It follows that $\gg_{02} \subset \ga$. Hence we have $\ga = \gg$.
\vspace{1mm}

(3) Case $\gq \cap \ga \not= \{0\}$. Let $R = (g_1, g_2, g_{-2}, g_{-1})$ ($g_i \in \gg_i$) be a non-zero element of $\gq \cap \ga$. In the case $g_1 \not= 0$, we put $g_1 = \sum_{i,j<k}g_{ijk}\e_i \otimes (\e_j \land \e_k)$. Without loss of generality, we may assume that $g_{112} \not= 0$. Putting $S_{ijkl} = (E_{ii} - E_{jj}, E_{kk} - E_{ll}) \in \gg_0$ and $T = \e_2 \otimes \e_1 \land \e_2 \in \gg_{-1}$, we have
$$
   \ad(T)\ad(S_{1523})\ad(S_{1415})\ad(S_{1314})\ad(S_{1213})R = (-E_{12}, 0) \in \gg_{01} \cap \ga. $$
Then we can reduce this case to the case (2). In the case $g_i \ne 0$ ($i = 2, -2, -1$), we can similarly reduce to the case (2).
\vspace{1mm}

Thus the simplicity of $\gg$ has been proved. Since the dimension of $\gg$ is 248, we see that $\gg$ is a $C$-Lie algebra of type of $E_8$.
\vspace{3mm}

{\bf Proposition 5.12.4.} {\it The Killing form $B_8$ of the Lie algebra ${\gge_8}^C = \sl(5, C) \oplus \sl(5, C) \oplus \gg_1 \oplus \gg_2 \oplus \gg_{-2} \oplus \gg_{-1}$ is given by}
$$
\begin{array}{l}
    B_8(R_1, R_2) = 60(\tr(C_1C_2) + \tr(D_1D_2) - (\x_1, \w_2)(\a_1, \dd_2) - (\x_2, \w_1)(\a_2, \dd_1)
\vspace{1mm}\\
\qquad \qquad
   - (\y_1, \z_2)(\b_1, \c_2) - (\y_2, \z_1)(\b_2, \c_1)),
\end{array} $$
where $R_i = (C_i, D_i, \x_i \otimes \a_i, \b_i \otimes \y_i, \c_i \otimes \z_i, \w_i \otimes \dd_i) \in {\gge_8}^C$.
\vspace{2mm}

{\bf Proof.} We consider a symmetric bilinear form $B$ of ${\gge_8}^C$:
$$
\begin{array}{l}
    B(R_1, R_2) = \tr(C_1C_2) + \tr(D_1D_2) - (\x_1, \w_2)(\a_1, \dd_2) - (\x_2, \w_1)(\a_2, \dd_1)
\vspace{1mm}\\
\qquad \qquad
   - (\y_1, \z_2)(\b_1, \c_2) - (\y_2, \z_1)(\b_2, \c_1).
\end{array} $$
Using Lemmas 5.12.1, 5.12.2, we see that $B$ is ${\gge_8}^C$-adjoint invariant. Since ${\gge_8}^C$ is simple, there exists $k \in C$ such that $B_8(R_1, R_2) = kB(R_1, R_2)$ for all $R_i \in {\gge_8}^C$. To determined $k$, let $R = R_1 = R_2 = (E_{11} - E_{22}, 0, 0, 0, 0, 0) \in {\gge_8}^C$. Then we have
$$
     B_8(R, R) = 120, \quad B(R, R) = 2. $$
Therefore $k = 60$.
\vspace{2mm}

We define a complex-conjugate linear transformation $\tau\wti{\lambda}$ of ${\gge_8}^C$ by
$$
\begin{array}{l}
    \tau\wti{\lambda}(C, D, \x \otimes \a, \b \otimes \y, \c \otimes \z, \w \otimes \dd)
\vspace{1mm}\\
\qquad \quad
  = (- \tau{}^tC, - \tau{}^tD, \tau \w \otimes \tau\dd, \tau\c \otimes \tau\z, \tau\b \otimes \tau\y, \tau\x \otimes \tau\a).
\end{array}  $$
Further we define a positive dfinite Hermitian inner product $\langle R_1, R_2 \rangle$ in ${\gge_8}^C$ by
$$
    \langle R_1, R_2 \rangle = -B_8(R_1, \tau\wti{\lambda}R_2). $$
Then, we have
$$
\begin{array}{l}
    \langle R_1, R_2 \rangle = 60(\tr(C_1\tau{}^tC_2) + \tr(D_1\tau{}^tD_2) + (\x_1, \tau\x_2)(\a_1, \tau\a_2) + (\y_1, \tau\y_2)(\b_1, \tau\b_2)
\vspace{1mm}\\
\qquad \qquad
   + (\z_1, \tau\z_2)(\c_1, \tau\c_2) + (\w_1, \tau\w_2)(\dd_1, \tau\dd_2)).
\end{array} $$
As in Theorem 5.5.3, 
$$
     E_8 = \{ \alpha \in \mbox{Aut}({\gge_8}^C) \, | \, \langle \alpha R_1, \alpha R_2 \rangle = \langle R_1, R_2 \rangle \} $$
is a simply connected compact simple Lie group of type $E_8$. 
\vspace{1mm}

Let $\zeta = \exp(2\pi i/5) \in C$ and we define a $C$-linear transformation $z_5$ of ${\gge_8}^C$ by
$$
   z_5(C, D, g_1, g_2, g_{-2}, g_{-1}) = (C, D, \zeta(g_1), \zeta^2(g_2), \zeta^3(g_{-2}), \zeta^4(g_{-1})). $$
Then $z_5 \in E_8$ and ${z_5}^5 = 1$.
\vspace{2mm}

Now, we study the following subgroup $(E_8)^{z_5}$ of $E_8$:
$$
     (E_8)^{z_5} = \{\alpha \in E_8 \, | \, z_5\alpha = \alpha z_5 \}. $$

{\bf Theorem 5.12.5.} $(E_8)^{z_5} \cong (SU(5) \times SU(5))/\Z_5, \Z_5 = \{(E, E), (\zeta E, \zeta^2E), $ $(\zeta^2E, \zeta^4E), (\zeta^3E, \zeta E), (\zeta^4E, \zeta^3E) \}, \zeta = \exp(2\pi i/5). $
\vspace{2mm}

{\bf Proof.}@We define mappings $\varphi_1, \varphi_2 : SU(5) \to E_8$ respectively by
$$
\begin{array}{l}
  \varphi_1(A)(C, D, \x \otimes \a, \b \otimes \y, \c \otimes \z, \w \otimes \dd)
\vspace{1mm}\\
\qquad
   = (ACA^{-1}, D, (A\x) \otimes \a, (A\b) \otimes \y, ({}^tA^{-1}\c) \otimes \z, ({}^tA^{-1}\w) \otimes \dd). 
\vspace{1mm}\\
 \varphi_2(B)(C, D, \x \otimes \a, \b \otimes \y, \c \otimes \z, \w \otimes \dd)
\vspace{1mm}\\
\qquad
   = (C, BDB^{-1}, \x \otimes (B\a), \b \otimes ({}^tB^{-1}\y), \c \otimes (B\z), \w \otimes ({}^tB^{-1}\dd)). 
\end{array}$$
$\varphi_1$ and $\varphi_2$ are well-defined: $\varphi_1(A), \varphi_2(B) \in E_8$. Indeed, for $Z \in \su(5)$, we have $(Z, 0) \in \gg_0$ and
$$
\begin{array}{l}
   \exp(\ad(Z, 0))(C, D, \x \otimes \a, \b \otimes \y, \c \otimes \z, \w \otimes \dd)
\vspace{1mm}\\
\qquad
   = (\exp(\ad(Z))C, D, ((\exp Z)\x) \otimes \a, 
\vspace{1mm}\\
\qquad \quad
   ((\exp Z)\b) \otimes \y, ((\exp(-{}^tZ))\c \otimes \z, ((\exp(-{}^tZ))\w) \otimes \dd)
\vspace{1mm}\\
\qquad
  = (\Ad(\exp Z)C, D, ((\exp Z)\x) \otimes \a,
\vspace{1mm}\\
\qquad \quad
   ((\exp Z)\b) \otimes \y, ({}^t(\exp Z)^{-1}\c) \otimes \z, ({}^t(\exp Z)^{-1}\w) \otimes \dd)
\vspace{1mm}\\
\qquad
   = \varphi_1(\exp Z)(C, D, \x \otimes \a, \b \otimes \y, \c \otimes \z, \w \otimes \dd). 
\end{array}$$
Hence $\varphi_1(A) \in \mbox{Aut}({\gge_8}^C) = {E_8}^C$. Using Lemma 5.11.1, we have
$$
   \langle \varphi(A)R_1, \varphi_1(A)R_2 \rangle = \langle R_1, R_2 \rangle. $$ Therefore $\varphi_1(A) \in E_8$. Similarly $\varphi_2(B) \in E_8$. 
\vspace{1mm}

Now, we define a mapping $\varphi : SU(5) \times SU(5) \to E_8$ by
$$
   \varphi(A, B) = \varphi_1(A)\varphi_2(B). $$
Since $\varphi_1(A)$ and $\varphi_2(B)$ commute, $\varphi$ is a homomorphism. We shall show that $\varphi$ is onto. Since $(\gge_8)^{z_5} = \su(5) \oplus \su(5)$, the differential $\varphi_*$ is onto. It is not difficult to see that 
$$
   \Ker\,\varphi = \{(E, E), (\zeta E, \zeta^2E), (\zeta^2E, \zeta^4E), (\zeta^3E, \zeta E), (\zeta^4E, \zeta^3E) \} = \Z_5.$$
 Further, since $(E_8)^{z_5}$ is connected, $\varphi$ is onto. Thus we have the isomorphism $(SU(5) \times SU(5))/\Z_5 \cong (E_8)^{z_5}$.

\vspace{4mm}

{\bf 5.13. Non-compact exceptional Lie groups $E_{8(8)}$ and $E_{8(-24)}$ of type $E_8$}
\vspace{3mm}

Let 
\begin{eqnarray*}
   \gge_{8(8)} \!\!\! &=& \!\!\! \gge_{7(7)} \oplus \gP' \oplus \gP' \oplus \R \oplus \R \oplus \R, \quad (\mbox{where} \; \gge_{7(7)} = ({\gge_7}^C)^{\tau\gamma}),
\vspace{1mm}\\
   \gge_{8(-24)} \!\!\! &=& \!\!\! \gge_{7(-25)} \oplus \gP \oplus \gP \oplus \R \oplus \R \oplus \R, \quad (\mbox{where} \; \gge_{7(-25)} = ({\gge_7}^C)^{\tau}).
\end{eqnarray*}
For $R_1, R_2 \in \gge_{8(8)}$ or $\gge_{8(-24)}$, we define a Lie bracket $[R_1, R_2]$ as similar to ${\gge_8}^C$ of Section 5.1. Now, we define groups $E_{8(8)}$ and $E_{8(-24)}$ by
\begin{eqnarray*}
    E_{8(8)} \!\!\! &=& \!\!\! \{\alpha \in \Iso_R(\gge_{8(8)}) \, | \, \alpha[R_1, R_2] = [\alpha R_1, \alpha R_2] \}, 
\vspace{1mm}\\
    E_{8(-24)} \!\!\! &=& \!\!\! \{\alpha \in \Iso_R(\gge_{8(-24)}) \, | \, \alpha[R_1, R_2] = [\alpha R_1, \alpha R_2] \}.
\end{eqnarray*}
These groups can also be defined by
$$
         E_{8(8)} \cong ({E_8}^C)^{\tau\gamma}, \quad 
         E_{8(-24)} \cong ({E_8}^C)^{\tau}. $$

{\bf Theorem 5.13.1.} {\it The polar decompositions of the Lie groups $E_{8(8)}$ and $E_{8(-24)}$ are respectively given by}

\begin{eqnarray*}
    E_{8(8)} \!\!\! &\simeq& \!\!\! Ss(16) \times \R^{128}, 
\vspace{1mm}\\
    E_{8(-24)} \!\!\! &\simeq& \!\!\! (SU(2) \times E_7)/\Z_2 \times \R^{112}. 
\end{eqnarray*}

{\bf Proof.} These are the facts corresponding to Theorems 5.8.7 and 5.7.6.
\vspace{3mm}

{\bf Theorem 5.13.2.} {\it The centers of the groups $E_{8(8)}$ and $E_{8(-24)}$ are trivial}:
$$
      z(E_{8(8)}) = \{1\}, \quad z(E_{8(-24)}) = \{1\}. $$

$*$ has some errors
\vspace{1mm}

$^{\dagger}$ $\overline{{\;}^{\;}\;\;}$ and $\tau$ are written by the same notaion $\overline{{\;}^{\;}\;\;}$.

\end{document}